\documentclass[10pt]{amsart}
\usepackage[cp866]{inputenc}
\usepackage[english,russian]{babel}
\usepackage{amsmath}
\usepackage{amssymb}
\usepackage{amsfonts}
\usepackage{color}

\begin{document}
\renewcommand{\refname}{References}
\thispagestyle{empty}

\title{Perfect colorings of $Z^2$: nine colors}
\author{{D. S. Krotov}}%
\address{Denis Stanislavovich Krotov
\newline\hphantom{iii} Sobolev Institute of Mathematics,
\newline\hphantom{iii} prosp. Akademika Koptyuga, 4,
\newline\hphantom{iii} 630090, Novosibirsk, Russia}%
\email{krotov@math.nsc.ru}%

\thanks{\sc Krotov, D. S.,
Perfect colorings of $Z^2$: nine colors}%

\maketitle {\small
\begin{quote}
\noindent{\sc Abstract. } We list all perfect colorings of $Z^2$ by $9$ or less colors \medskip

 \end{quote}
}

A vertex coloring $f : V \stackrel{on}{\to} \{1,...,k\}$ of a graph $G=(V,E)$ is called \emph{perfect}
with parameters $(a_{ij})_{i,j=1}^k$
iff for each $i$ and $j$ from $1$ to $k$
every color-$i$ vertex has exactly
$a_{ij}$ color-$j$ neighbors.
We consider the graph of $Z^2$, whose vertices are integer pairs,
two vertices $(x,y)$ and $(x',y')$ being adjacent if and only if $|x-x'|+|y-y'|=1$.
This graph has degree $4$, and so $\sum_{i=1}^k a_{ij} = 4$ for any admissible parameter matrix
$(a_{ij})_{i,j=1}^k$.
The main goal of this note is to present the results of an exhaustive computer search,
listing all admissible parameters for $k \leq 9$.
The cases $k=2$ and $k=3$ were solved (without the use of computer)
by Axenovich \cite{Axe2003} and Puzynina \cite{Puz2005},
respectively.
In the Appendix, we list all the matrices with corresponding colorings
(for convenience, we display graph vertices as colored cells).
Some matrices admit an infinite number of colorings;
such colorings have diagonals colored by two alternating colors
(as proved in \cite{Puz2004}, one the colorings is periodic in the two main directions);
we show only one example, which is not necessarily periodical or canonical
in some other sense (indeed, it is lexicographically first in the order of cells
used by the coloring algorithm).
The others admit a finite number of colorings
(which are periodic in the two main directions \cite{Puz2004});
in this case, we list all possible colorings, up to graph isomorphisms.
There are several matrices admitting more than one but less than infinity number
of non-isomorphic colorings ($k\leq 9$). Here are the numbers of these matrices:
2-7, 3-8, 3-9, 3-10, 3-17, 3-18, 3-20, 4-14, 4-30, 4-32, 4-34, 4-38, 4-39, 4-41, 4-46,
5-36, 5-43, 6-36, 6-37, 6-58, 6-67, 6-80, 6-83, 6-89, 8-123, 8-125, 8-150, 9-117.


\section*{Appendix}

\def\1{\colorbox[rgb]{1.0,0.7,0.7}{\texttt{\footnotesize 1}}}
\def\2{\colorbox[rgb]{0.7,1.0,0.7}{\texttt{\footnotesize 2}}}
\def\3{\colorbox[rgb]{0.7,0.7,1.0}{\texttt{\footnotesize 3}}}
\def\4{\colorbox[rgb]{1.0,1.0,0.5}{\texttt{\footnotesize 4}}}
\def\5{\colorbox[rgb]{0.7,1.0,1.0}{\texttt{\footnotesize 5}}}
\def\6{\colorbox[rgb]{0.95,0.7,1.0}{\texttt{\footnotesize 6}}}
\def\7{\colorbox[rgb]{1.0,0.85,0.6}{\texttt{\footnotesize 7}}}
\def\8{\colorbox[rgb]{0.9,0.9,0.7}{\texttt{\footnotesize 8}}}
\def\9{\colorbox[rgb]{0.7,0.7,0.7}{\texttt{\footnotesize 9}}}
\def\A{\colorbox[rgb]{1.0,1.0,1.0}{\texttt{\footnotesize A}}}
\def\ennn{\\[-0.9ex]}
\def\eeee{\\[-0.6ex]}
\def\baaa{}
\def\eaaa{}

\def\bbbb{\begin{tabular}{|cccccccccc|}}
\def\ebbb{\end{tabular}\ }
\def\baaa{\par}
\def\eaaa{\ \ }



\baaa
2-1
\eaaa
\bbbb
0&4\\
1&3\\
\ebbb
\parbox{7cm}{ 
\1\2\2\2\2\1\2\2\2\2\1\2\2\2\2\1\eeee
\2\2\2\1\2\2\2\2\1\2\2\2\2\1\2\2\eeee
\2\1\2\2\2\2\1\2\2\2\2\1\2\2\2\2\eeee
\2\2\2\2\1\2\2\2\2\1\2\2\2\2\1\2\eeee
\2\2\1\2\2\2\2\1\2\2\2\2\1\2\2\2\eeee
\1\2\2\2\2\1\2\2\2\2\1\2\2\2\2\1\eeee
\2\2\2\1\2\2\2\2\1\2\2\2\2\1\2\2\eeee
\2\1\2\2\2\2\1\2\2\2\2\1\2\2\2\2\eeee
\2\2\2\2\1\2\2\2\2\1\2\2\2\2\1\2\eeee
\2\2\1\2\2\2\2\1\2\2\2\2\1\2\2\2\eeee
\1\2\2\2\2\1\2\2\2\2\1\2\2\2\2\1\eeee
\2\2\2\1\2\2\2\2\1\2\2\2\2\1\2\2\eeee
} 

\baaa
2-2
\eaaa
\bbbb
0&4\\
2&2\\
\ebbb
\parbox{7cm}{ 
\1\2\2\1\2\2\1\2\2\1\2\2\1\2\2\1\eeee
\2\1\2\2\1\2\2\1\2\2\1\2\2\1\2\2\eeee
\2\2\1\2\2\1\2\2\1\2\2\1\2\2\1\2\eeee
\1\2\2\1\2\2\1\2\2\1\2\2\1\2\2\1\eeee
\2\1\2\2\1\2\2\1\2\2\1\2\2\1\2\2\eeee
\2\2\1\2\2\1\2\2\1\2\2\1\2\2\1\2\eeee
\1\2\2\1\2\2\1\2\2\1\2\2\1\2\2\1\eeee
\2\1\2\2\1\2\2\1\2\2\1\2\2\1\2\2\eeee
\2\2\1\2\2\1\2\2\1\2\2\1\2\2\1\2\eeee
\1\2\2\1\2\2\1\2\2\1\2\2\1\2\2\1\eeee
\2\1\2\2\1\2\2\1\2\2\1\2\2\1\2\2\eeee
\2\2\1\2\2\1\2\2\1\2\2\1\2\2\1\2\eeee
} 

\baaa
2-3
\eaaa
\bbbb
0&4\\
4&0\\
\ebbb
\parbox{7cm}{ 
\1\2\1\2\1\2\1\2\1\2\1\2\1\2\1\2\eeee
\2\1\2\1\2\1\2\1\2\1\2\1\2\1\2\1\eeee
\1\2\1\2\1\2\1\2\1\2\1\2\1\2\1\2\eeee
\2\1\2\1\2\1\2\1\2\1\2\1\2\1\2\1\eeee
\1\2\1\2\1\2\1\2\1\2\1\2\1\2\1\2\eeee
\2\1\2\1\2\1\2\1\2\1\2\1\2\1\2\1\eeee
\1\2\1\2\1\2\1\2\1\2\1\2\1\2\1\2\eeee
\2\1\2\1\2\1\2\1\2\1\2\1\2\1\2\1\eeee
\1\2\1\2\1\2\1\2\1\2\1\2\1\2\1\2\eeee
\2\1\2\1\2\1\2\1\2\1\2\1\2\1\2\1\eeee
\1\2\1\2\1\2\1\2\1\2\1\2\1\2\1\2\eeee
\2\1\2\1\2\1\2\1\2\1\2\1\2\1\2\1\eeee
} 

\baaa
2-4
\eaaa
\bbbb
1&3\\
1&3\\
\ebbb
\parbox{7cm}{ 
\1\2\2\2\1\2\2\2\1\2\2\2\1\2\2\2\eeee
\1\2\2\2\1\2\2\2\1\2\2\2\1\2\2\2\eeee
\2\2\1\2\2\2\1\2\2\2\1\2\2\2\1\2\eeee
\2\2\1\2\2\2\1\2\2\2\1\2\2\2\1\2\eeee
\1\2\2\2\1\2\2\2\1\2\2\2\1\2\2\2\eeee
\1\2\2\2\1\2\2\2\1\2\2\2\1\2\2\2\eeee
\2\2\1\2\2\2\1\2\2\2\1\2\2\2\1\2\eeee
\2\2\1\2\2\2\1\2\2\2\1\2\2\2\1\2\eeee
\1\2\2\2\1\2\2\2\1\2\2\2\1\2\2\2\eeee
\1\2\2\2\1\2\2\2\1\2\2\2\1\2\2\2\eeee
\2\2\1\2\2\2\1\2\2\2\1\2\2\2\1\2\eeee
\2\2\1\2\2\2\1\2\2\2\1\2\2\2\1\2\eeee
} 

\baaa
2-5
\eaaa
\bbbb
1&3\\
2&2\\
\ebbb
\parbox{7cm}{ 
\2\2\1\1\2\2\2\1\1\2\2\2\1\1\2\2\eeee
\1\2\2\2\1\1\2\2\2\1\1\2\2\2\1\1\eeee
\2\1\1\2\2\2\1\1\2\2\2\1\1\2\2\2\eeee
\2\2\2\1\1\2\2\2\1\1\2\2\2\1\1\2\eeee
\1\1\2\2\2\1\1\2\2\2\1\1\2\2\2\1\eeee
\2\2\1\1\2\2\2\1\1\2\2\2\1\1\2\2\eeee
\1\2\2\2\1\1\2\2\2\1\1\2\2\2\1\1\eeee
\2\1\1\2\2\2\1\1\2\2\2\1\1\2\2\2\eeee
\2\2\2\1\1\2\2\2\1\1\2\2\2\1\1\2\eeee
\1\1\2\2\2\1\1\2\2\2\1\1\2\2\2\1\eeee
\2\2\1\1\2\2\2\1\1\2\2\2\1\1\2\2\eeee
\1\2\2\2\1\1\2\2\2\1\1\2\2\2\1\1\eeee
} 

\baaa
2-6
\eaaa
\bbbb
1&3\\
3&1\\
\ebbb
\parbox{7cm}{ 
\1\2\1\2\1\2\1\2\1\2\1\2\1\2\1\2\eeee
\1\2\1\2\1\2\1\2\1\2\1\2\1\2\1\2\eeee
\2\1\2\1\2\1\2\1\2\1\2\1\2\1\2\1\eeee
\2\1\2\1\2\1\2\1\2\1\2\1\2\1\2\1\eeee
\1\2\1\2\1\2\1\2\1\2\1\2\1\2\1\2\eeee
\1\2\1\2\1\2\1\2\1\2\1\2\1\2\1\2\eeee
\2\1\2\1\2\1\2\1\2\1\2\1\2\1\2\1\eeee
\2\1\2\1\2\1\2\1\2\1\2\1\2\1\2\1\eeee
\1\2\1\2\1\2\1\2\1\2\1\2\1\2\1\2\eeee
\1\2\1\2\1\2\1\2\1\2\1\2\1\2\1\2\eeee
\2\1\2\1\2\1\2\1\2\1\2\1\2\1\2\1\eeee
\2\1\2\1\2\1\2\1\2\1\2\1\2\1\2\1\eeee
} 

\baaa
2-7
\eaaa
\bbbb
2&2\\
1&3\\
\ebbb
\parbox{7cm}{ 
\1\1\2\2\1\1\2\2\1\1\2\2\1\1\2\2\eeee
\1\1\2\2\1\1\2\2\1\1\2\2\1\1\2\2\eeee
\2\2\2\2\2\2\2\2\2\2\2\2\2\2\2\2\eeee
\2\2\1\1\2\2\1\1\2\2\1\1\2\2\1\1\eeee
\2\2\1\1\2\2\1\1\2\2\1\1\2\2\1\1\eeee
\2\2\2\2\2\2\2\2\2\2\2\2\2\2\2\2\eeee
\1\1\2\2\1\1\2\2\1\1\2\2\1\1\2\2\eeee
\1\1\2\2\1\1\2\2\1\1\2\2\1\1\2\2\eeee
\2\2\2\2\2\2\2\2\2\2\2\2\2\2\2\2\eeee
\2\2\1\1\2\2\1\1\2\2\1\1\2\2\1\1\eeee
\2\2\1\1\2\2\1\1\2\2\1\1\2\2\1\1\eeee
\2\2\2\2\2\2\2\2\2\2\2\2\2\2\2\2\eeee
} 
\baaa
\phantom{0-}\#
\eaaa
\mbox{}\phantom{\bbbb
0&0\\
\ebbb}
\parbox{7cm}{ 
\1\2\2\1\2\2\1\2\2\1\2\2\1\2\2\1\eeee
\1\2\2\1\2\2\1\2\2\1\2\2\1\2\2\1\eeee
\1\2\2\1\2\2\1\2\2\1\2\2\1\2\2\1\eeee
\1\2\2\1\2\2\1\2\2\1\2\2\1\2\2\1\eeee
\1\2\2\1\2\2\1\2\2\1\2\2\1\2\2\1\eeee
\1\2\2\1\2\2\1\2\2\1\2\2\1\2\2\1\eeee
\1\2\2\1\2\2\1\2\2\1\2\2\1\2\2\1\eeee
\1\2\2\1\2\2\1\2\2\1\2\2\1\2\2\1\eeee
\1\2\2\1\2\2\1\2\2\1\2\2\1\2\2\1\eeee
\1\2\2\1\2\2\1\2\2\1\2\2\1\2\2\1\eeee
\1\2\2\1\2\2\1\2\2\1\2\2\1\2\2\1\eeee
\1\2\2\1\2\2\1\2\2\1\2\2\1\2\2\1\eeee
} 

\baaa
2-8
\eaaa
\bbbb
2&2\\
2&2\\
\ebbb
\parbox{7cm}{ 
\1\1\2\1\1\1\2\1\1\1\2\1\1\1\2\1\eeee
\1\1\2\2\2\1\2\2\2\1\2\2\2\1\2\2\eeee
\2\2\1\1\2\1\1\1\2\1\1\1\2\1\1\1\eeee
\1\2\1\1\2\2\2\1\2\2\2\1\2\2\2\1\eeee
\1\2\2\2\1\1\2\1\1\1\2\1\1\1\2\1\eeee
\1\1\1\2\1\1\2\2\2\1\2\2\2\1\2\2\eeee
\2\2\1\2\2\2\1\1\2\1\1\1\2\1\1\1\eeee
\1\2\1\1\1\2\1\1\2\2\2\1\2\2\2\1\eeee
\1\2\2\2\1\2\2\2\1\1\2\1\1\1\2\1\eeee
\1\1\1\2\1\1\1\2\1\1\2\2\2\1\2\2\eeee
\2\2\1\2\2\2\1\2\2\2\1\1\2\1\1\1\eeee
\1\2\1\1\1\2\1\1\1\2\1\1\2\2\2\1\eeee
} 

\baaa
2-9
\eaaa
\bbbb
3&1\\
1&3\\
\ebbb
\parbox{7cm}{ 
\1\1\2\2\1\1\2\2\1\1\2\2\1\1\2\2\eeee
\1\1\2\2\1\1\2\2\1\1\2\2\1\1\2\2\eeee
\1\1\2\2\1\1\2\2\1\1\2\2\1\1\2\2\eeee
\1\1\2\2\1\1\2\2\1\1\2\2\1\1\2\2\eeee
\1\1\2\2\1\1\2\2\1\1\2\2\1\1\2\2\eeee
\1\1\2\2\1\1\2\2\1\1\2\2\1\1\2\2\eeee
\1\1\2\2\1\1\2\2\1\1\2\2\1\1\2\2\eeee
\1\1\2\2\1\1\2\2\1\1\2\2\1\1\2\2\eeee
\1\1\2\2\1\1\2\2\1\1\2\2\1\1\2\2\eeee
\1\1\2\2\1\1\2\2\1\1\2\2\1\1\2\2\eeee
\1\1\2\2\1\1\2\2\1\1\2\2\1\1\2\2\eeee
\1\1\2\2\1\1\2\2\1\1\2\2\1\1\2\2\eeee
} 



\baaa
3-1
\eaaa
\bbbb
0&0&4\\
0&0&4\\
1&1&2\\
\ebbb
\parbox{7cm}{ 
\1\3\3\1\3\3\1\3\3\1\3\3\1\3\3\1\eeee
\3\2\3\3\2\3\3\2\3\3\2\3\3\2\3\3\eeee
\3\3\1\3\3\1\3\3\1\3\3\1\3\3\1\3\eeee
\1\3\3\2\3\3\2\3\3\2\3\3\2\3\3\2\eeee
\3\2\3\3\1\3\3\1\3\3\1\3\3\1\3\3\eeee
\3\3\1\3\3\2\3\3\2\3\3\2\3\3\2\3\eeee
\1\3\3\2\3\3\1\3\3\1\3\3\1\3\3\1\eeee
\3\2\3\3\1\3\3\2\3\3\2\3\3\2\3\3\eeee
\3\3\1\3\3\2\3\3\1\3\3\1\3\3\1\3\eeee
\1\3\3\2\3\3\1\3\3\2\3\3\2\3\3\2\eeee
\3\2\3\3\1\3\3\2\3\3\1\3\3\1\3\3\eeee
\3\3\1\3\3\2\3\3\1\3\3\2\3\3\2\3\eeee
} 

\baaa
3-2
\eaaa
\bbbb
0&0&4\\
0&0&4\\
1&3&0\\
\ebbb
\parbox{7cm}{ 
\1\3\2\3\1\3\2\3\1\3\2\3\1\3\2\3\eeee
\3\2\3\2\3\2\3\2\3\2\3\2\3\2\3\2\eeee
\2\3\1\3\2\3\1\3\2\3\1\3\2\3\1\3\eeee
\3\2\3\2\3\2\3\2\3\2\3\2\3\2\3\2\eeee
\1\3\2\3\1\3\2\3\1\3\2\3\1\3\2\3\eeee
\3\2\3\2\3\2\3\2\3\2\3\2\3\2\3\2\eeee
\2\3\1\3\2\3\1\3\2\3\1\3\2\3\1\3\eeee
\3\2\3\2\3\2\3\2\3\2\3\2\3\2\3\2\eeee
\1\3\2\3\1\3\2\3\1\3\2\3\1\3\2\3\eeee
\3\2\3\2\3\2\3\2\3\2\3\2\3\2\3\2\eeee
\2\3\1\3\2\3\1\3\2\3\1\3\2\3\1\3\eeee
\3\2\3\2\3\2\3\2\3\2\3\2\3\2\3\2\eeee
} 

\baaa
3-3
\eaaa
\bbbb
0&0&4\\
0&0&4\\
2&2&0\\
\ebbb
\parbox{7cm}{ 
\1\3\2\3\1\3\2\3\1\3\2\3\1\3\2\3\eeee
\3\1\3\2\3\1\3\2\3\1\3\2\3\1\3\2\eeee
\2\3\1\3\2\3\1\3\2\3\1\3\2\3\1\3\eeee
\3\2\3\1\3\2\3\1\3\2\3\1\3\2\3\1\eeee
\1\3\2\3\1\3\2\3\1\3\2\3\1\3\2\3\eeee
\3\1\3\2\3\1\3\2\3\1\3\2\3\1\3\2\eeee
\2\3\1\3\2\3\1\3\2\3\1\3\2\3\1\3\eeee
\3\2\3\1\3\2\3\1\3\2\3\1\3\2\3\1\eeee
\1\3\2\3\1\3\2\3\1\3\2\3\1\3\2\3\eeee
\3\1\3\2\3\1\3\2\3\1\3\2\3\1\3\2\eeee
\2\3\1\3\2\3\1\3\2\3\1\3\2\3\1\3\eeee
\3\2\3\1\3\2\3\1\3\2\3\1\3\2\3\1\eeee
} 

\baaa
3-4
\eaaa
\bbbb
0&0&4\\
0&2&2\\
1&2&1\\
\ebbb
\parbox{7cm}{ 
\1\3\3\1\3\3\1\3\3\1\3\3\1\3\3\1\eeee
\3\2\2\3\2\2\3\2\2\3\2\2\3\2\2\3\eeee
\3\2\2\3\2\2\3\2\2\3\2\2\3\2\2\3\eeee
\1\3\3\1\3\3\1\3\3\1\3\3\1\3\3\1\eeee
\3\2\2\3\2\2\3\2\2\3\2\2\3\2\2\3\eeee
\3\2\2\3\2\2\3\2\2\3\2\2\3\2\2\3\eeee
\1\3\3\1\3\3\1\3\3\1\3\3\1\3\3\1\eeee
\3\2\2\3\2\2\3\2\2\3\2\2\3\2\2\3\eeee
\3\2\2\3\2\2\3\2\2\3\2\2\3\2\2\3\eeee
\1\3\3\1\3\3\1\3\3\1\3\3\1\3\3\1\eeee
\3\2\2\3\2\2\3\2\2\3\2\2\3\2\2\3\eeee
\3\2\2\3\2\2\3\2\2\3\2\2\3\2\2\3\eeee
} 

\baaa
3-5
\eaaa
\bbbb
0&0&4\\
0&2&2\\
2&2&0\\
\ebbb
\parbox{7cm}{ 
\1\3\2\2\3\1\3\2\2\3\1\3\2\2\3\1\eeee
\3\1\3\2\2\3\1\3\2\2\3\1\3\2\2\3\eeee
\2\3\1\3\2\2\3\1\3\2\2\3\1\3\2\2\eeee
\2\2\3\1\3\2\2\3\1\3\2\2\3\1\3\2\eeee
\3\2\2\3\1\3\2\2\3\1\3\2\2\3\1\3\eeee
\1\3\2\2\3\1\3\2\2\3\1\3\2\2\3\1\eeee
\3\1\3\2\2\3\1\3\2\2\3\1\3\2\2\3\eeee
\2\3\1\3\2\2\3\1\3\2\2\3\1\3\2\2\eeee
\2\2\3\1\3\2\2\3\1\3\2\2\3\1\3\2\eeee
\3\2\2\3\1\3\2\2\3\1\3\2\2\3\1\3\eeee
\1\3\2\2\3\1\3\2\2\3\1\3\2\2\3\1\eeee
\3\1\3\2\2\3\1\3\2\2\3\1\3\2\2\3\eeee
} 

\baaa
3-6
\eaaa
\bbbb
0&1&3\\
1&0&3\\
1&1&2\\
\ebbb
\parbox{7cm}{ 
\1\3\2\3\3\1\3\2\3\3\1\3\2\3\3\1\eeee
\2\3\3\1\3\2\3\3\1\3\2\3\3\1\3\2\eeee
\3\1\3\2\3\3\1\3\2\3\3\1\3\2\3\3\eeee
\3\2\3\3\1\3\2\3\3\1\3\2\3\3\1\3\eeee
\3\3\1\3\2\3\3\1\3\2\3\3\1\3\2\3\eeee
\1\3\2\3\3\1\3\2\3\3\1\3\2\3\3\1\eeee
\2\3\3\1\3\2\3\3\1\3\2\3\3\1\3\2\eeee
\3\1\3\2\3\3\1\3\2\3\3\1\3\2\3\3\eeee
\3\2\3\3\1\3\2\3\3\1\3\2\3\3\1\3\eeee
\3\3\1\3\2\3\3\1\3\2\3\3\1\3\2\3\eeee
\1\3\2\3\3\1\3\2\3\3\1\3\2\3\3\1\eeee
\2\3\3\1\3\2\3\3\1\3\2\3\3\1\3\2\eeee
} 

\baaa
3-7
\eaaa
\bbbb
0&1&3\\
1&2&1\\
3&1&0\\
\ebbb
\parbox{7cm}{ 
\1\3\1\3\1\3\1\3\1\3\1\3\1\3\1\3\eeee
\2\2\2\2\2\2\2\2\2\2\2\2\2\2\2\2\eeee
\3\1\3\1\3\1\3\1\3\1\3\1\3\1\3\1\eeee
\1\3\1\3\1\3\1\3\1\3\1\3\1\3\1\3\eeee
\2\2\2\2\2\2\2\2\2\2\2\2\2\2\2\2\eeee
\3\1\3\1\3\1\3\1\3\1\3\1\3\1\3\1\eeee
\1\3\1\3\1\3\1\3\1\3\1\3\1\3\1\3\eeee
\2\2\2\2\2\2\2\2\2\2\2\2\2\2\2\2\eeee
\3\1\3\1\3\1\3\1\3\1\3\1\3\1\3\1\eeee
\1\3\1\3\1\3\1\3\1\3\1\3\1\3\1\3\eeee
\2\2\2\2\2\2\2\2\2\2\2\2\2\2\2\2\eeee
\3\1\3\1\3\1\3\1\3\1\3\1\3\1\3\1\eeee
} 

\baaa
3-8
\eaaa
\bbbb
0&2&2\\
1&1&2\\
1&2&1\\
\ebbb
\parbox{7cm}{ 
\1\2\3\2\3\1\2\3\2\3\1\2\3\2\3\1\eeee
\2\3\1\2\3\2\3\1\2\3\2\3\1\2\3\2\eeee
\2\3\2\3\1\2\3\2\3\1\2\3\2\3\1\2\eeee
\3\1\2\3\2\3\1\2\3\2\3\1\2\3\2\3\eeee
\3\2\3\1\2\3\2\3\1\2\3\2\3\1\2\3\eeee
\1\2\3\2\3\1\2\3\2\3\1\2\3\2\3\1\eeee
\2\3\1\2\3\2\3\1\2\3\2\3\1\2\3\2\eeee
\2\3\2\3\1\2\3\2\3\1\2\3\2\3\1\2\eeee
\3\1\2\3\2\3\1\2\3\2\3\1\2\3\2\3\eeee
\3\2\3\1\2\3\2\3\1\2\3\2\3\1\2\3\eeee
\1\2\3\2\3\1\2\3\2\3\1\2\3\2\3\1\eeee
\2\3\1\2\3\2\3\1\2\3\2\3\1\2\3\2\eeee
} 
\baaa
\phantom{0-}\#
\eaaa
\mbox{}\phantom{\bbbb
0&0&0\\
\ebbb}
\parbox{7cm}{ 
\1\3\2\2\3\1\3\2\2\3\1\3\2\2\3\1\eeee
\2\2\3\1\3\2\2\3\1\3\2\2\3\1\3\2\eeee
\3\1\3\2\2\3\1\3\2\2\3\1\3\2\2\3\eeee
\3\2\2\3\1\3\2\2\3\1\3\2\2\3\1\3\eeee
\2\3\1\3\2\2\3\1\3\2\2\3\1\3\2\2\eeee
\1\3\2\2\3\1\3\2\2\3\1\3\2\2\3\1\eeee
\2\2\3\1\3\2\2\3\1\3\2\2\3\1\3\2\eeee
\3\1\3\2\2\3\1\3\2\2\3\1\3\2\2\3\eeee
\3\2\2\3\1\3\2\2\3\1\3\2\2\3\1\3\eeee
\2\3\1\3\2\2\3\1\3\2\2\3\1\3\2\2\eeee
\1\3\2\2\3\1\3\2\2\3\1\3\2\2\3\1\eeee
\2\2\3\1\3\2\2\3\1\3\2\2\3\1\3\2\eeee
} 

\baaa
3-9
\eaaa
\bbbb
0&2&2\\
1&2&1\\
1&1&2\\
\ebbb
\parbox{7cm}{ 
\1\2\3\3\2\1\3\2\2\3\1\2\3\3\2\1\eeee
\2\2\3\1\2\3\3\2\1\3\2\2\3\1\2\3\eeee
\2\1\3\2\2\3\1\2\3\3\2\1\3\2\2\3\eeee
\2\3\3\2\1\3\2\2\3\1\2\3\3\2\1\3\eeee
\2\3\1\2\3\3\2\1\3\2\2\3\1\2\3\3\eeee
\1\3\2\2\3\1\2\3\3\2\1\3\2\2\3\1\eeee
\3\3\2\1\3\2\2\3\1\2\3\3\2\1\3\2\eeee
\3\1\2\3\3\2\1\3\2\2\3\1\2\3\3\2\eeee
\3\2\2\3\1\2\3\3\2\1\3\2\2\3\1\2\eeee
\3\2\1\3\2\2\3\1\2\3\3\2\1\3\2\2\eeee
\1\2\3\3\2\1\3\2\2\3\1\2\3\3\2\1\eeee
\2\2\3\1\2\3\3\2\1\3\2\2\3\1\2\3\eeee
} 
\baaa
\phantom{0-}\#
\eaaa
\mbox{}\phantom{\bbbb
0&0&0\\
\ebbb}
\parbox{7cm}{ 
\1\2\2\3\3\1\3\3\2\2\1\2\2\3\3\1\eeee
\3\2\2\1\2\2\3\3\1\3\3\2\2\1\2\2\eeee
\3\1\3\3\2\2\1\2\2\3\3\1\3\3\2\2\eeee
\2\2\3\3\1\3\3\2\2\1\2\2\3\3\1\3\eeee
\2\2\1\2\2\3\3\1\3\3\2\2\1\2\2\3\eeee
\1\3\3\2\2\1\2\2\3\3\1\3\3\2\2\1\eeee
\2\3\3\1\3\3\2\2\1\2\2\3\3\1\3\3\eeee
\2\1\2\2\3\3\1\3\3\2\2\1\2\2\3\3\eeee
\3\3\2\2\1\2\2\3\3\1\3\3\2\2\1\2\eeee
\3\3\1\3\3\2\2\1\2\2\3\3\1\3\3\2\eeee
\1\2\2\3\3\1\3\3\2\2\1\2\2\3\3\1\eeee
\3\2\2\1\2\2\3\3\1\3\3\2\2\1\2\2\eeee
} 

\baaa
3-10
\eaaa
\bbbb
0&2&2\\
1&2&1\\
2&2&0\\
\ebbb
\parbox{7cm}{ 
\1\2\2\3\1\2\2\3\1\2\2\3\1\2\2\3\eeee
\2\3\1\2\2\3\1\2\2\3\1\2\2\3\1\2\eeee
\2\1\3\2\2\1\3\2\2\1\3\2\2\1\3\2\eeee
\3\2\2\1\3\2\2\1\3\2\2\1\3\2\2\1\eeee
\1\2\2\3\1\2\2\3\1\2\2\3\1\2\2\3\eeee
\2\3\1\2\2\3\1\2\2\3\1\2\2\3\1\2\eeee
\2\1\3\2\2\1\3\2\2\1\3\2\2\1\3\2\eeee
\3\2\2\1\3\2\2\1\3\2\2\1\3\2\2\1\eeee
\1\2\2\3\1\2\2\3\1\2\2\3\1\2\2\3\eeee
\2\3\1\2\2\3\1\2\2\3\1\2\2\3\1\2\eeee
\2\1\3\2\2\1\3\2\2\1\3\2\2\1\3\2\eeee
\3\2\2\1\3\2\2\1\3\2\2\1\3\2\2\1\eeee
} 
\baaa
\phantom{0-0}\#
\eaaa
\mbox{}\phantom{\bbbb
0&0&0\\
\ebbb}
\parbox{7cm}{ 
\1\2\3\2\1\2\3\2\1\2\3\2\1\2\3\2\eeee
\3\2\1\2\3\2\1\2\3\2\1\2\3\2\1\2\eeee
\1\2\3\2\1\2\3\2\1\2\3\2\1\2\3\2\eeee
\3\2\1\2\3\2\1\2\3\2\1\2\3\2\1\2\eeee
\1\2\3\2\1\2\3\2\1\2\3\2\1\2\3\2\eeee
\3\2\1\2\3\2\1\2\3\2\1\2\3\2\1\2\eeee
\1\2\3\2\1\2\3\2\1\2\3\2\1\2\3\2\eeee
\3\2\1\2\3\2\1\2\3\2\1\2\3\2\1\2\eeee
\1\2\3\2\1\2\3\2\1\2\3\2\1\2\3\2\eeee
\3\2\1\2\3\2\1\2\3\2\1\2\3\2\1\2\eeee
\1\2\3\2\1\2\3\2\1\2\3\2\1\2\3\2\eeee
\3\2\1\2\3\2\1\2\3\2\1\2\3\2\1\2\eeee
} 

\baaa
3-11
\eaaa
\bbbb
0&2&2\\
2&0&2\\
2&2&0\\
\ebbb
\parbox{7cm}{ 
\1\2\3\1\2\3\1\2\3\1\2\3\1\2\3\1\eeee
\2\3\1\2\3\1\2\3\1\2\3\1\2\3\1\2\eeee
\3\1\2\3\1\2\3\1\2\3\1\2\3\1\2\3\eeee
\1\2\3\1\2\3\1\2\3\1\2\3\1\2\3\1\eeee
\2\3\1\2\3\1\2\3\1\2\3\1\2\3\1\2\eeee
\3\1\2\3\1\2\3\1\2\3\1\2\3\1\2\3\eeee
\1\2\3\1\2\3\1\2\3\1\2\3\1\2\3\1\eeee
\2\3\1\2\3\1\2\3\1\2\3\1\2\3\1\2\eeee
\3\1\2\3\1\2\3\1\2\3\1\2\3\1\2\3\eeee
\1\2\3\1\2\3\1\2\3\1\2\3\1\2\3\1\eeee
\2\3\1\2\3\1\2\3\1\2\3\1\2\3\1\2\eeee
\3\1\2\3\1\2\3\1\2\3\1\2\3\1\2\3\eeee
} 

\baaa
3-12
\eaaa
\bbbb
0&2&2\\
2&1&1\\
2&1&1\\
\ebbb
\parbox{7cm}{ 
\1\2\2\1\2\2\1\2\2\1\2\2\1\2\2\1\eeee
\2\1\3\3\1\3\3\1\3\3\1\3\3\1\3\3\eeee
\2\3\1\2\2\1\2\2\1\2\2\1\2\2\1\2\eeee
\1\3\2\1\3\3\1\3\3\1\3\3\1\3\3\1\eeee
\2\1\2\3\1\2\2\1\2\2\1\2\2\1\2\2\eeee
\2\3\1\3\2\1\3\3\1\3\3\1\3\3\1\3\eeee
\1\3\2\1\2\3\1\2\2\1\2\2\1\2\2\1\eeee
\2\1\2\3\1\3\2\1\3\3\1\3\3\1\3\3\eeee
\2\3\1\3\2\1\2\3\1\2\2\1\2\2\1\2\eeee
\1\3\2\1\2\3\1\3\2\1\3\3\1\3\3\1\eeee
\2\1\2\3\1\3\2\1\2\3\1\2\2\1\2\2\eeee
\2\3\1\3\2\1\2\3\1\3\2\1\3\3\1\3\eeee
} 

\baaa
3-13
\eaaa
\bbbb
0&2&2\\
2&2&0\\
2&0&2\\
\ebbb
\parbox{7cm}{ 
\1\2\2\1\3\3\1\2\2\1\3\3\1\2\2\1\eeee
\2\2\1\3\3\1\2\2\1\3\3\1\2\2\1\3\eeee
\2\1\3\3\1\2\2\1\3\3\1\2\2\1\3\3\eeee
\1\3\3\1\2\2\1\3\3\1\2\2\1\3\3\1\eeee
\3\3\1\2\2\1\3\3\1\2\2\1\3\3\1\2\eeee
\3\1\2\2\1\3\3\1\2\2\1\3\3\1\2\2\eeee
\1\2\2\1\3\3\1\2\2\1\3\3\1\2\2\1\eeee
\2\2\1\3\3\1\2\2\1\3\3\1\2\2\1\3\eeee
\2\1\3\3\1\2\2\1\3\3\1\2\2\1\3\3\eeee
\1\3\3\1\2\2\1\3\3\1\2\2\1\3\3\1\eeee
\3\3\1\2\2\1\3\3\1\2\2\1\3\3\1\2\eeee
\3\1\2\2\1\3\3\1\2\2\1\3\3\1\2\2\eeee
} 

\baaa
3-14
\eaaa
\bbbb
1&0&3\\
0&1&3\\
1&1&2\\
\ebbb
\parbox{7cm}{ 
\1\3\2\3\3\2\3\1\3\3\1\3\2\3\3\2\eeee
\1\3\3\1\3\2\3\3\2\3\1\3\3\1\3\2\eeee
\3\2\3\1\3\3\1\3\2\3\3\2\3\1\3\3\eeee
\3\2\3\3\2\3\1\3\3\1\3\2\3\3\2\3\eeee
\3\3\1\3\2\3\3\2\3\1\3\3\1\3\2\3\eeee
\2\3\1\3\3\1\3\2\3\3\2\3\1\3\3\1\eeee
\2\3\3\2\3\1\3\3\1\3\2\3\3\2\3\1\eeee
\3\1\3\2\3\3\2\3\1\3\3\1\3\2\3\3\eeee
\3\1\3\3\1\3\2\3\3\2\3\1\3\3\1\3\eeee
\3\3\2\3\1\3\3\1\3\2\3\3\2\3\1\3\eeee
\1\3\2\3\3\2\3\1\3\3\1\3\2\3\3\2\eeee
\1\3\3\1\3\2\3\3\2\3\1\3\3\1\3\2\eeee
} 

\baaa
3-15
\eaaa
\bbbb
1&0&3\\
0&1&3\\
1&2&1\\
\ebbb
\parbox{7cm}{ 
\1\3\2\3\2\3\1\3\2\3\2\3\1\3\2\3\eeee
\1\3\2\3\2\3\1\3\2\3\2\3\1\3\2\3\eeee
\3\2\3\1\3\2\3\2\3\1\3\2\3\2\3\1\eeee
\3\2\3\1\3\2\3\2\3\1\3\2\3\2\3\1\eeee
\1\3\2\3\2\3\1\3\2\3\2\3\1\3\2\3\eeee
\1\3\2\3\2\3\1\3\2\3\2\3\1\3\2\3\eeee
\3\2\3\1\3\2\3\2\3\1\3\2\3\2\3\1\eeee
\3\2\3\1\3\2\3\2\3\1\3\2\3\2\3\1\eeee
\1\3\2\3\2\3\1\3\2\3\2\3\1\3\2\3\eeee
\1\3\2\3\2\3\1\3\2\3\2\3\1\3\2\3\eeee
\3\2\3\1\3\2\3\2\3\1\3\2\3\2\3\1\eeee
\3\2\3\1\3\2\3\2\3\1\3\2\3\2\3\1\eeee
} 

\baaa
3-16
\eaaa
\bbbb
1&1&2\\
1&1&2\\
1&1&2\\
\ebbb
\parbox{7cm}{ 
\1\2\3\1\1\2\3\1\1\2\3\1\1\2\3\1\eeee
\1\2\3\3\3\2\3\3\3\2\3\3\3\2\3\3\eeee
\3\3\1\2\3\1\1\2\3\1\1\2\3\1\1\2\eeee
\2\3\1\2\3\3\3\2\3\3\3\2\3\3\3\2\eeee
\1\3\3\3\1\2\3\1\1\2\3\1\1\2\3\1\eeee
\1\2\2\3\1\2\3\3\3\2\3\3\3\2\3\3\eeee
\3\3\1\3\3\3\1\2\3\1\1\2\3\1\1\2\eeee
\2\3\1\2\2\3\1\2\3\3\3\2\3\3\3\2\eeee
\1\3\3\3\1\3\3\3\1\2\3\1\1\2\3\1\eeee
\1\2\2\3\1\2\2\3\1\2\3\3\3\2\3\3\eeee
\3\3\1\3\3\3\1\3\3\3\1\2\3\1\1\2\eeee
\2\3\1\2\2\3\1\2\2\3\1\2\3\3\3\2\eeee
} 

\baaa
3-17
\eaaa
\bbbb
1&1&2\\
1&2&1\\
2&1&1\\
\ebbb
\parbox{7cm}{ 
\1\2\3\1\2\3\1\2\3\1\2\3\1\2\3\1\eeee
\1\2\3\1\2\3\1\2\3\1\2\3\1\2\3\1\eeee
\3\2\1\3\2\1\3\2\1\3\2\1\3\2\1\3\eeee
\3\2\1\3\2\1\3\2\1\3\2\1\3\2\1\3\eeee
\1\2\3\1\2\3\1\2\3\1\2\3\1\2\3\1\eeee
\1\2\3\1\2\3\1\2\3\1\2\3\1\2\3\1\eeee
\3\2\1\3\2\1\3\2\1\3\2\1\3\2\1\3\eeee
\3\2\1\3\2\1\3\2\1\3\2\1\3\2\1\3\eeee
\1\2\3\1\2\3\1\2\3\1\2\3\1\2\3\1\eeee
\1\2\3\1\2\3\1\2\3\1\2\3\1\2\3\1\eeee
\3\2\1\3\2\1\3\2\1\3\2\1\3\2\1\3\eeee
\3\2\1\3\2\1\3\2\1\3\2\1\3\2\1\3\eeee
} 
\baaa
\phantom{0-0}\#
\eaaa
\mbox{}\phantom{\bbbb
0&0&0\\
\ebbb}
\parbox{7cm}{ 
\1\2\2\3\1\2\2\3\1\2\2\3\1\2\2\1\eeee
\1\2\2\3\1\2\2\3\1\2\2\3\1\2\2\3\eeee
\3\3\1\1\3\3\1\1\3\3\1\1\3\3\1\1\eeee
\2\1\3\2\2\1\3\2\2\1\3\2\2\1\3\2\eeee
\2\1\3\2\2\1\3\2\2\1\3\2\2\1\3\2\eeee
\3\3\1\1\3\3\1\1\3\3\1\1\3\3\1\1\eeee
\1\2\2\3\1\2\2\3\1\2\2\3\1\2\2\3\eeee
\1\2\2\3\1\2\2\3\1\2\2\3\1\2\2\3\eeee
\3\3\1\1\3\3\1\1\3\3\1\1\3\3\1\1\eeee
\2\1\3\2\2\1\3\2\2\1\3\2\2\1\3\2\eeee
\2\1\3\2\2\1\3\2\2\1\3\2\2\1\3\2\eeee
\3\3\1\1\3\3\1\1\3\3\1\1\3\3\1\1\eeee
} 
\baaa
\phantom{0-0}\#
\eaaa
\mbox{}\phantom{\bbbb
0&0&0\\
\ebbb}
\parbox{7cm}{ 
\1\2\2\3\1\2\2\3\1\2\2\3\1\2\2\3\eeee
\1\3\1\3\1\3\1\3\1\3\1\3\1\3\1\3\eeee
\2\3\1\2\2\3\1\2\2\3\1\2\2\3\1\2\eeee
\2\1\3\2\2\1\3\2\2\1\3\2\2\1\3\2\eeee
\3\1\3\1\3\1\3\1\3\1\3\1\3\1\3\1\eeee
\3\2\2\1\3\2\2\1\3\2\2\1\3\2\2\1\eeee
\1\2\2\3\1\2\2\3\1\2\2\3\1\2\2\3\eeee
\1\3\1\3\1\3\1\3\1\3\1\3\1\3\1\3\eeee
\2\3\1\2\2\3\1\2\2\3\1\2\2\3\1\2\eeee
\2\1\3\2\2\1\3\2\2\1\3\2\2\1\3\2\eeee
\3\1\3\1\3\1\3\1\3\1\3\1\3\1\3\1\eeee
\3\2\2\1\3\2\2\1\3\2\2\1\3\2\2\1\eeee
} 
\baaa
\phantom{0-0}\#
\eaaa
\mbox{}\phantom{\bbbb
0&0&0\\
\ebbb}
\parbox{7cm}{ 
\1\3\1\3\1\3\1\3\1\3\1\3\1\3\1\3\eeee
\1\3\1\3\1\3\1\3\1\3\1\3\1\3\1\3\eeee
\2\2\2\2\2\2\2\2\2\2\2\2\2\2\2\2\eeee
\3\1\3\1\3\1\3\1\3\1\3\1\3\1\3\1\eeee
\3\1\3\1\3\1\3\1\3\1\3\1\3\1\3\1\eeee
\2\2\2\2\2\2\2\2\2\2\2\2\2\2\2\2\eeee
\1\3\1\3\1\3\1\3\1\3\1\3\1\3\1\3\eeee
\1\3\1\3\1\3\1\3\1\3\1\3\1\3\1\3\eeee
\2\2\2\2\2\2\2\2\2\2\2\2\2\2\2\2\eeee
\3\1\3\1\3\1\3\1\3\1\3\1\3\1\3\1\eeee
\3\1\3\1\3\1\3\1\3\1\3\1\3\1\3\1\eeee
\2\2\2\2\2\2\2\2\2\2\2\2\2\2\2\2\eeee
} 

\baaa
3-18
\eaaa
\bbbb
2&0&2\\
0&2&2\\
1&1&2\\
\ebbb
\parbox{7cm}{ 
\1\1\3\3\1\1\3\3\1\1\3\3\1\1\3\3\eeee
\1\1\3\3\1\1\3\3\1\1\3\3\1\1\3\3\eeee
\3\3\2\2\3\3\2\2\3\3\2\2\3\3\2\2\eeee
\3\3\2\2\3\3\2\2\3\3\2\2\3\3\2\2\eeee
\1\1\3\3\1\1\3\3\1\1\3\3\1\1\3\3\eeee
\1\1\3\3\1\1\3\3\1\1\3\3\1\1\3\3\eeee
\3\3\2\2\3\3\2\2\3\3\2\2\3\3\2\2\eeee
\3\3\2\2\3\3\2\2\3\3\2\2\3\3\2\2\eeee
\1\1\3\3\1\1\3\3\1\1\3\3\1\1\3\3\eeee
\1\1\3\3\1\1\3\3\1\1\3\3\1\1\3\3\eeee
\3\3\2\2\3\3\2\2\3\3\2\2\3\3\2\2\eeee
\3\3\2\2\3\3\2\2\3\3\2\2\3\3\2\2\eeee
} 
\baaa
\phantom{0-0}\#
\eaaa
\mbox{}\phantom{\bbbb
0&0&0\\
\ebbb}
\parbox{7cm}{ 
\1\3\2\3\1\3\2\3\1\3\2\3\1\3\2\3\eeee
\1\3\2\3\1\3\2\3\1\3\2\3\1\3\2\3\eeee
\1\3\2\3\1\3\2\3\1\3\2\3\1\3\2\3\eeee
\1\3\2\3\1\3\2\3\1\3\2\3\1\3\2\3\eeee
\1\3\2\3\1\3\2\3\1\3\2\3\1\3\2\3\eeee
\1\3\2\3\1\3\2\3\1\3\2\3\1\3\2\3\eeee
\1\3\2\3\1\3\2\3\1\3\2\3\1\3\2\3\eeee
\1\3\2\3\1\3\2\3\1\3\2\3\1\3\2\3\eeee
\1\3\2\3\1\3\2\3\1\3\2\3\1\3\2\3\eeee
\1\3\2\3\1\3\2\3\1\3\2\3\1\3\2\3\eeee
\1\3\2\3\1\3\2\3\1\3\2\3\1\3\2\3\eeee
\1\3\2\3\1\3\2\3\1\3\2\3\1\3\2\3\eeee
} 

\baaa
3-19
\eaaa
\bbbb
2&0&2\\
0&3&1\\
1&1&2\\
\ebbb
\parbox{7cm}{ 
\1\3\2\2\3\1\3\2\2\3\1\3\2\2\3\1\eeee
\1\3\2\2\3\1\3\2\2\3\1\3\2\2\3\1\eeee
\1\3\2\2\3\1\3\2\2\3\1\3\2\2\3\1\eeee
\1\3\2\2\3\1\3\2\2\3\1\3\2\2\3\1\eeee
\1\3\2\2\3\1\3\2\2\3\1\3\2\2\3\1\eeee
\1\3\2\2\3\1\3\2\2\3\1\3\2\2\3\1\eeee
\1\3\2\2\3\1\3\2\2\3\1\3\2\2\3\1\eeee
\1\3\2\2\3\1\3\2\2\3\1\3\2\2\3\1\eeee
\1\3\2\2\3\1\3\2\2\3\1\3\2\2\3\1\eeee
\1\3\2\2\3\1\3\2\2\3\1\3\2\2\3\1\eeee
\1\3\2\2\3\1\3\2\2\3\1\3\2\2\3\1\eeee
\1\3\2\2\3\1\3\2\2\3\1\3\2\2\3\1\eeee
} 

\baaa
3-20
\eaaa
\bbbb
2&1&1\\
1&2&1\\
1&1&2\\
\ebbb
\parbox{7cm}{ 
\1\1\2\2\1\1\2\2\1\1\2\2\1\1\2\2\eeee
\1\1\3\3\1\1\3\3\1\1\3\3\1\1\3\3\eeee
\2\2\3\3\2\2\3\3\2\2\3\3\2\2\3\3\eeee
\2\2\1\1\2\2\1\1\2\2\1\1\2\2\1\1\eeee
\3\3\1\1\3\3\1\1\3\3\1\1\3\3\1\1\eeee
\3\3\2\2\3\3\2\2\3\3\2\2\3\3\2\2\eeee
\1\1\2\2\1\1\2\2\1\1\2\2\1\1\2\2\eeee
\1\1\3\3\1\1\3\3\1\1\3\3\1\1\3\3\eeee
\2\2\3\3\2\2\3\3\2\2\3\3\2\2\3\3\eeee
\2\2\1\1\2\2\1\1\2\2\1\1\2\2\1\1\eeee
\3\3\1\1\3\3\1\1\3\3\1\1\3\3\1\1\eeee
\3\3\2\2\3\3\2\2\3\3\2\2\3\3\2\2\eeee
} 
\baaa
\phantom{0-0}\#
\eaaa
\mbox{}\phantom{\bbbb
0&0&0\\
\ebbb}
\parbox{7cm}{ 
\1\1\3\3\1\1\3\3\1\1\3\3\1\1\3\3\eeee
\1\1\3\3\1\1\3\3\1\1\3\3\1\1\3\3\eeee
\2\2\2\2\2\2\2\2\2\2\2\2\2\2\2\2\eeee
\3\3\1\1\3\3\1\1\3\3\1\1\3\3\1\1\eeee
\3\3\1\1\3\3\1\1\3\3\1\1\3\3\1\1\eeee
\2\2\2\2\2\2\2\2\2\2\2\2\2\2\2\2\eeee
\1\1\3\3\1\1\3\3\1\1\3\3\1\1\3\3\eeee
\1\1\3\3\1\1\3\3\1\1\3\3\1\1\3\3\eeee
\2\2\2\2\2\2\2\2\2\2\2\2\2\2\2\2\eeee
\3\3\1\1\3\3\1\1\3\3\1\1\3\3\1\1\eeee
\3\3\1\1\3\3\1\1\3\3\1\1\3\3\1\1\eeee
\2\2\2\2\2\2\2\2\2\2\2\2\2\2\2\2\eeee
} 
\baaa
\phantom{0-0}\#
\eaaa
\mbox{}\phantom{\bbbb
0&0&0\\
\ebbb}
\parbox{7cm}{ 
\1\2\3\1\2\3\1\2\3\1\2\3\1\2\3\1\eeee
\1\2\3\1\2\3\1\2\3\1\2\3\1\2\3\1\eeee
\1\2\3\1\2\3\1\2\3\1\2\3\1\2\3\1\eeee
\1\2\3\1\2\3\1\2\3\1\2\3\1\2\3\1\eeee
\1\2\3\1\2\3\1\2\3\1\2\3\1\2\3\1\eeee
\1\2\3\1\2\3\1\2\3\1\2\3\1\2\3\1\eeee
\1\2\3\1\2\3\1\2\3\1\2\3\1\2\3\1\eeee
\1\2\3\1\2\3\1\2\3\1\2\3\1\2\3\1\eeee
\1\2\3\1\2\3\1\2\3\1\2\3\1\2\3\1\eeee
\1\2\3\1\2\3\1\2\3\1\2\3\1\2\3\1\eeee
\1\2\3\1\2\3\1\2\3\1\2\3\1\2\3\1\eeee
\1\2\3\1\2\3\1\2\3\1\2\3\1\2\3\1\eeee
} 

\baaa
3-21
\eaaa
\bbbb
2&1&1\\
1&3&0\\
1&0&3\\
\ebbb
\parbox{7cm}{ 
\1\2\2\1\3\3\1\2\2\1\3\3\1\2\2\1\eeee
\1\2\2\1\3\3\1\2\2\1\3\3\1\2\2\1\eeee
\1\2\2\1\3\3\1\2\2\1\3\3\1\2\2\1\eeee
\1\2\2\1\3\3\1\2\2\1\3\3\1\2\2\1\eeee
\1\2\2\1\3\3\1\2\2\1\3\3\1\2\2\1\eeee
\1\2\2\1\3\3\1\2\2\1\3\3\1\2\2\1\eeee
\1\2\2\1\3\3\1\2\2\1\3\3\1\2\2\1\eeee
\1\2\2\1\3\3\1\2\2\1\3\3\1\2\2\1\eeee
\1\2\2\1\3\3\1\2\2\1\3\3\1\2\2\1\eeee
\1\2\2\1\3\3\1\2\2\1\3\3\1\2\2\1\eeee
\1\2\2\1\3\3\1\2\2\1\3\3\1\2\2\1\eeee
\1\2\2\1\3\3\1\2\2\1\3\3\1\2\2\1\eeee
} 


\baaa
4-1
\eaaa
\bbbb
0&0&0&4\\
0&0&0&4\\
0&0&0&4\\
1&1&2&0\\
\ebbb
\parbox{7cm}{ 
\1\4\3\4\1\4\3\4\1\4\3\4\1\4\3\4\eeee
\4\2\4\3\4\2\4\3\4\2\4\3\4\2\4\3\eeee
\3\4\1\4\3\4\1\4\3\4\1\4\3\4\1\4\eeee
\4\3\4\2\4\3\4\2\4\3\4\2\4\3\4\2\eeee
\1\4\3\4\1\4\3\4\1\4\3\4\1\4\3\4\eeee
\4\2\4\3\4\2\4\3\4\2\4\3\4\2\4\3\eeee
\3\4\1\4\3\4\1\4\3\4\1\4\3\4\1\4\eeee
\4\3\4\2\4\3\4\2\4\3\4\2\4\3\4\2\eeee
\1\4\3\4\1\4\3\4\1\4\3\4\1\4\3\4\eeee
\4\2\4\3\4\2\4\3\4\2\4\3\4\2\4\3\eeee
\3\4\1\4\3\4\1\4\3\4\1\4\3\4\1\4\eeee
\4\3\4\2\4\3\4\2\4\3\4\2\4\3\4\2\eeee
} 

\baaa
4-2
\eaaa
\bbbb
0&0&0&4\\
0&0&0&4\\
0&0&2&2\\
1&1&2&0\\
\ebbb
\parbox{7cm}{ 
\1\4\3\3\4\1\4\3\3\4\1\4\3\3\4\1\eeee
\4\2\4\3\3\4\2\4\3\3\4\2\4\3\3\4\eeee
\3\4\1\4\3\3\4\1\4\3\3\4\1\4\3\3\eeee
\3\3\4\2\4\3\3\4\2\4\3\3\4\2\4\3\eeee
\4\3\3\4\1\4\3\3\4\1\4\3\3\4\1\4\eeee
\1\4\3\3\4\2\4\3\3\4\2\4\3\3\4\2\eeee
\4\2\4\3\3\4\1\4\3\3\4\1\4\3\3\4\eeee
\3\4\1\4\3\3\4\2\4\3\3\4\2\4\3\3\eeee
\3\3\4\2\4\3\3\4\1\4\3\3\4\1\4\3\eeee
\4\3\3\4\1\4\3\3\4\2\4\3\3\4\2\4\eeee
\1\4\3\3\4\2\4\3\3\4\1\4\3\3\4\1\eeee
\4\2\4\3\3\4\1\4\3\3\4\2\4\3\3\4\eeee
} 

\baaa
4-3
\eaaa
\bbbb
0&0&0&4\\
0&0&1&3\\
0&4&0&0\\
1&3&0&0\\
\ebbb
\parbox{7cm}{ 
\1\4\2\4\2\3\2\4\2\4\1\4\2\4\2\3\eeee
\4\2\4\1\4\2\4\2\3\2\4\2\4\1\4\2\eeee
\2\3\2\4\2\4\1\4\2\4\2\3\2\4\2\4\eeee
\4\2\4\2\3\2\4\2\4\1\4\2\4\2\3\2\eeee
\2\4\1\4\2\4\2\3\2\4\2\4\1\4\2\4\eeee
\3\2\4\2\4\1\4\2\4\2\3\2\4\2\4\1\eeee
\2\4\2\3\2\4\2\4\1\4\2\4\2\3\2\4\eeee
\4\1\4\2\4\2\3\2\4\2\4\1\4\2\4\2\eeee
\2\4\2\4\1\4\2\4\2\3\2\4\2\4\1\4\eeee
\4\2\3\2\4\2\4\1\4\2\4\2\3\2\4\2\eeee
\1\4\2\4\2\3\2\4\2\4\1\4\2\4\2\3\eeee
\4\2\4\1\4\2\4\2\3\2\4\2\4\1\4\2\eeee
} 

\baaa
4-4
\eaaa
\bbbb
0&0&0&4\\
0&0&2&2\\
0&2&0&2\\
1&1&1&1\\
\ebbb
\parbox{7cm}{ 
\1\4\4\1\4\4\1\4\4\1\4\4\1\4\4\1\eeee
\4\2\3\4\2\3\4\2\3\4\2\3\4\2\3\4\eeee
\4\3\2\4\3\2\4\3\2\4\3\2\4\3\2\4\eeee
\1\4\4\1\4\4\1\4\4\1\4\4\1\4\4\1\eeee
\4\2\3\4\2\3\4\2\3\4\2\3\4\2\3\4\eeee
\4\3\2\4\3\2\4\3\2\4\3\2\4\3\2\4\eeee
\1\4\4\1\4\4\1\4\4\1\4\4\1\4\4\1\eeee
\4\2\3\4\2\3\4\2\3\4\2\3\4\2\3\4\eeee
\4\3\2\4\3\2\4\3\2\4\3\2\4\3\2\4\eeee
\1\4\4\1\4\4\1\4\4\1\4\4\1\4\4\1\eeee
\4\2\3\4\2\3\4\2\3\4\2\3\4\2\3\4\eeee
\4\3\2\4\3\2\4\3\2\4\3\2\4\3\2\4\eeee
} 

\baaa
4-5
\eaaa
\bbbb
0&0&0&4\\
0&0&2&2\\
0&2&2&0\\
2&2&0&0\\
\ebbb
\parbox{7cm}{ 
\1\4\2\3\3\2\4\1\4\2\3\3\2\4\1\4\eeee
\4\1\4\2\3\3\2\4\1\4\2\3\3\2\4\1\eeee
\2\4\1\4\2\3\3\2\4\1\4\2\3\3\2\4\eeee
\3\2\4\1\4\2\3\3\2\4\1\4\2\3\3\2\eeee
\3\3\2\4\1\4\2\3\3\2\4\1\4\2\3\3\eeee
\2\3\3\2\4\1\4\2\3\3\2\4\1\4\2\3\eeee
\4\2\3\3\2\4\1\4\2\3\3\2\4\1\4\2\eeee
\1\4\2\3\3\2\4\1\4\2\3\3\2\4\1\4\eeee
\4\1\4\2\3\3\2\4\1\4\2\3\3\2\4\1\eeee
\2\4\1\4\2\3\3\2\4\1\4\2\3\3\2\4\eeee
\3\2\4\1\4\2\3\3\2\4\1\4\2\3\3\2\eeee
\3\3\2\4\1\4\2\3\3\2\4\1\4\2\3\3\eeee
} 

\baaa
4-6
\eaaa
\bbbb
0&0&0&4\\
0&0&2&2\\
0&4&0&0\\
2&2&0&0\\
\ebbb
\parbox{7cm}{ 
\1\4\2\3\2\4\1\4\2\3\2\4\1\4\2\3\eeee
\4\1\4\2\3\2\4\1\4\2\3\2\4\1\4\2\eeee
\2\4\1\4\2\3\2\4\1\4\2\3\2\4\1\4\eeee
\3\2\4\1\4\2\3\2\4\1\4\2\3\2\4\1\eeee
\2\3\2\4\1\4\2\3\2\4\1\4\2\3\2\4\eeee
\4\2\3\2\4\1\4\2\3\2\4\1\4\2\3\2\eeee
\1\4\2\3\2\4\1\4\2\3\2\4\1\4\2\3\eeee
\4\1\4\2\3\2\4\1\4\2\3\2\4\1\4\2\eeee
\2\4\1\4\2\3\2\4\1\4\2\3\2\4\1\4\eeee
\3\2\4\1\4\2\3\2\4\1\4\2\3\2\4\1\eeee
\2\3\2\4\1\4\2\3\2\4\1\4\2\3\2\4\eeee
\4\2\3\2\4\1\4\2\3\2\4\1\4\2\3\2\eeee
} 

\baaa
4-7
\eaaa
\bbbb
0&0&0&4\\
0&1&1&2\\
0&1&1&2\\
1&1&1&1\\
\ebbb
\parbox{7cm}{ 
\1\4\4\1\4\4\1\4\4\1\4\4\1\4\4\1\eeee
\4\2\3\4\2\3\4\2\3\4\2\3\4\2\3\4\eeee
\4\2\3\4\2\3\4\2\3\4\2\3\4\2\3\4\eeee
\1\4\4\1\4\4\1\4\4\1\4\4\1\4\4\1\eeee
\4\3\2\4\3\2\4\3\2\4\3\2\4\3\2\4\eeee
\4\3\2\4\3\2\4\3\2\4\3\2\4\3\2\4\eeee
\1\4\4\1\4\4\1\4\4\1\4\4\1\4\4\1\eeee
\4\2\3\4\2\3\4\2\3\4\2\3\4\2\3\4\eeee
\4\2\3\4\2\3\4\2\3\4\2\3\4\2\3\4\eeee
\1\4\4\1\4\4\1\4\4\1\4\4\1\4\4\1\eeee
\4\3\2\4\3\2\4\3\2\4\3\2\4\3\2\4\eeee
\4\3\2\4\3\2\4\3\2\4\3\2\4\3\2\4\eeee
} 

\baaa
4-8
\eaaa
\bbbb
0&0&0&4\\
0&1&1&2\\
0&1&1&2\\
2&1&1&0\\
\ebbb
\parbox{7cm}{ 
\1\4\2\2\4\1\4\2\2\4\1\4\2\2\4\1\eeee
\4\1\4\3\3\4\1\4\3\3\4\1\4\3\3\4\eeee
\2\4\1\4\2\2\4\1\4\2\2\4\1\4\2\2\eeee
\2\3\4\1\4\3\3\4\1\4\3\3\4\1\4\3\eeee
\4\3\2\4\1\4\2\2\4\1\4\2\2\4\1\4\eeee
\1\4\2\3\4\1\4\3\3\4\1\4\3\3\4\1\eeee
\4\1\4\3\2\4\1\4\2\2\4\1\4\2\2\4\eeee
\2\4\1\4\2\3\4\1\4\3\3\4\1\4\3\3\eeee
\2\3\4\1\4\3\2\4\1\4\2\2\4\1\4\2\eeee
\4\3\2\4\1\4\2\3\4\1\4\3\3\4\1\4\eeee
\1\4\2\3\4\1\4\3\2\4\1\4\2\2\4\1\eeee
\4\1\4\3\2\4\1\4\2\3\4\1\4\3\3\4\eeee
} 

\baaa
4-9
\eaaa
\bbbb
0&0&0&4\\
0&1&1&2\\
0&1&2&1\\
1&2&1&0\\
\ebbb
\parbox{7cm}{ 
\1\4\3\3\2\4\2\2\4\2\3\3\4\1\4\3\eeee
\4\2\3\3\4\1\4\3\3\2\4\2\2\4\2\3\eeee
\3\2\4\2\2\4\2\3\3\4\1\4\3\3\2\4\eeee
\3\4\1\4\3\3\2\4\2\2\4\2\3\3\4\1\eeee
\2\2\4\2\3\3\4\1\4\3\3\2\4\2\2\4\eeee
\4\3\3\2\4\2\2\4\2\3\3\4\1\4\3\3\eeee
\2\3\3\4\1\4\3\3\2\4\2\2\4\2\3\3\eeee
\2\4\2\2\4\2\3\3\4\1\4\3\3\2\4\2\eeee
\4\1\4\3\3\2\4\2\2\4\2\3\3\4\1\4\eeee
\2\4\2\3\3\4\1\4\3\3\2\4\2\2\4\2\eeee
\3\3\2\4\2\2\4\2\3\3\4\1\4\3\3\2\eeee
\3\3\4\1\4\3\3\2\4\2\2\4\2\3\3\4\eeee
} 

\baaa
4-10
\eaaa
\bbbb
0&0&0&4\\
0&2&0&2\\
0&0&2&2\\
1&1&1&1\\
\ebbb
\parbox{7cm}{ 
\1\4\4\1\4\4\1\4\4\1\4\4\1\4\4\1\eeee
\4\2\2\4\3\3\4\2\2\4\3\3\4\2\2\4\eeee
\4\2\2\4\3\3\4\2\2\4\3\3\4\2\2\4\eeee
\1\4\4\1\4\4\1\4\4\1\4\4\1\4\4\1\eeee
\4\3\3\4\2\2\4\3\3\4\2\2\4\3\3\4\eeee
\4\3\3\4\2\2\4\3\3\4\2\2\4\3\3\4\eeee
\1\4\4\1\4\4\1\4\4\1\4\4\1\4\4\1\eeee
\4\2\2\4\3\3\4\2\2\4\3\3\4\2\2\4\eeee
\4\2\2\4\3\3\4\2\2\4\3\3\4\2\2\4\eeee
\1\4\4\1\4\4\1\4\4\1\4\4\1\4\4\1\eeee
\4\3\3\4\2\2\4\3\3\4\2\2\4\3\3\4\eeee
\4\3\3\4\2\2\4\3\3\4\2\2\4\3\3\4\eeee
} 

\baaa
4-11
\eaaa
\bbbb
0&0&1&3\\
0&0&1&3\\
1&3&0&0\\
1&3&0&0\\
\ebbb
\parbox{7cm}{ 
\1\4\2\3\1\4\2\3\1\4\2\3\1\4\2\3\eeee
\3\2\4\2\4\2\4\2\4\2\4\2\4\2\4\2\eeee
\2\4\1\4\2\3\1\4\2\3\1\4\2\3\1\4\eeee
\4\2\3\2\4\2\4\2\4\2\4\2\4\2\4\2\eeee
\1\4\2\4\1\4\2\3\1\4\2\3\1\4\2\3\eeee
\3\2\4\2\3\2\4\2\4\2\4\2\4\2\4\2\eeee
\2\4\1\4\2\4\1\4\2\3\1\4\2\3\1\4\eeee
\4\2\3\2\4\2\3\2\4\2\4\2\4\2\4\2\eeee
\1\4\2\4\1\4\2\4\1\4\2\3\1\4\2\3\eeee
\3\2\4\2\3\2\4\2\3\2\4\2\4\2\4\2\eeee
\2\4\1\4\2\4\1\4\2\4\1\4\2\3\1\4\eeee
\4\2\3\2\4\2\3\2\4\2\3\2\4\2\4\2\eeee
} 

\baaa
4-12
\eaaa
\bbbb
0&0&1&3\\
0&0&1&3\\
2&2&0&0\\
2&2&0&0\\
\ebbb
\parbox{7cm}{ 
\1\4\2\3\1\4\2\3\1\4\2\3\1\4\2\3\eeee
\3\1\4\2\4\1\4\2\4\1\4\2\4\1\4\2\eeee
\2\4\1\4\2\3\1\4\2\3\1\4\2\3\1\4\eeee
\4\2\3\1\4\2\4\1\4\2\4\1\4\2\4\1\eeee
\1\4\2\4\1\4\2\3\1\4\2\3\1\4\2\3\eeee
\3\1\4\2\3\1\4\2\4\1\4\2\4\1\4\2\eeee
\2\4\1\4\2\4\1\4\2\3\1\4\2\3\1\4\eeee
\4\2\3\1\4\2\3\1\4\2\4\1\4\2\4\1\eeee
\1\4\2\4\1\4\2\4\1\4\2\3\1\4\2\3\eeee
\3\1\4\2\3\1\4\2\3\1\4\2\4\1\4\2\eeee
\2\4\1\4\2\4\1\4\2\4\1\4\2\3\1\4\eeee
\4\2\3\1\4\2\3\1\4\2\3\1\4\2\4\1\eeee
} 

\baaa
4-13
\eaaa
\bbbb
0&0&1&3\\
0&0&2&2\\
1&3&0&0\\
2&2&0&0\\
\ebbb
\parbox{7cm}{ 
\1\4\1\3\2\4\2\3\1\4\2\4\1\3\2\4\eeee
\3\2\4\2\3\1\4\2\4\1\3\2\4\2\3\1\eeee
\2\3\1\4\2\4\1\3\2\4\2\3\1\4\2\4\eeee
\4\2\4\1\3\2\4\2\3\1\4\2\4\1\3\2\eeee
\1\3\2\4\2\3\1\4\2\4\1\3\2\4\2\3\eeee
\4\2\3\1\4\2\4\1\3\2\4\2\3\1\4\2\eeee
\1\4\2\4\1\3\2\4\2\3\1\4\2\4\1\3\eeee
\4\1\3\2\4\2\3\1\4\2\4\1\3\2\4\2\eeee
\2\4\2\3\1\4\2\4\1\3\2\4\2\3\1\4\eeee
\3\1\4\2\4\1\3\2\4\2\3\1\4\2\4\1\eeee
\2\4\1\3\2\4\2\3\1\4\2\4\1\3\2\4\eeee
\3\2\4\2\3\1\4\2\4\1\3\2\4\2\3\1\eeee
} 

\baaa
4-14
\eaaa
\bbbb
0&0&1&3\\
0&0&2&2\\
2&2&0&0\\
3&1&0&0\\
\ebbb
\parbox{7cm}{ 
\2\4\1\3\1\4\2\4\1\3\1\4\2\4\1\3\eeee
\3\1\4\2\4\1\3\1\4\2\4\1\3\1\4\2\eeee
\2\4\1\3\1\4\2\4\1\3\1\4\2\4\1\3\eeee
\3\1\4\2\4\1\3\1\4\2\4\1\3\1\4\2\eeee
\2\4\1\3\1\4\2\4\1\3\1\4\2\4\1\3\eeee
\3\1\4\2\4\1\3\1\4\2\4\1\3\1\4\2\eeee
\2\4\1\3\1\4\2\4\1\3\1\4\2\4\1\3\eeee
\3\1\4\2\4\1\3\1\4\2\4\1\3\1\4\2\eeee
\2\4\1\3\1\4\2\4\1\3\1\4\2\4\1\3\eeee
\3\1\4\2\4\1\3\1\4\2\4\1\3\1\4\2\eeee
\2\4\1\3\1\4\2\4\1\3\1\4\2\4\1\3\eeee
\3\1\4\2\4\1\3\1\4\2\4\1\3\1\4\2\eeee
} 
\baaa
\phantom{0-0}\#
\eaaa
\mbox{}\phantom{\bbbb
0&0&0&0\\
\ebbb}
\parbox{7cm}{ 
\1\4\1\4\1\4\1\4\1\4\1\4\1\4\1\4\eeee
\3\1\4\2\3\1\4\2\3\1\4\2\3\1\4\2\eeee
\2\4\1\3\2\4\1\3\2\4\1\3\2\4\1\3\eeee
\4\1\4\1\4\1\4\1\4\1\4\1\4\1\4\1\eeee
\1\3\2\4\1\3\2\4\1\3\2\4\1\3\2\4\eeee
\4\2\3\1\4\2\3\1\4\2\3\1\4\2\3\1\eeee
\1\4\1\4\1\4\1\4\1\4\1\4\1\4\1\4\eeee
\3\1\4\2\3\1\4\2\3\1\4\2\3\1\4\2\eeee
\2\4\1\3\2\4\1\3\2\4\1\3\2\4\1\3\eeee
\4\1\4\1\4\1\4\1\4\1\4\1\4\1\4\1\eeee
\1\3\2\4\1\3\2\4\1\3\2\4\1\3\2\4\eeee
\4\2\3\1\4\2\3\1\4\2\3\1\4\2\3\1\eeee
} 

\baaa
4-15
\eaaa
\bbbb
0&0&1&3\\
0&0&3&1\\
1&3&0&0\\
3&1&0&0\\
\ebbb
\parbox{7cm}{ 
\1\4\1\4\1\4\1\4\1\4\1\4\1\4\1\4\eeee
\3\2\3\2\3\2\3\2\3\2\3\2\3\2\3\2\eeee
\2\3\2\3\2\3\2\3\2\3\2\3\2\3\2\3\eeee
\4\1\4\1\4\1\4\1\4\1\4\1\4\1\4\1\eeee
\1\4\1\4\1\4\1\4\1\4\1\4\1\4\1\4\eeee
\3\2\3\2\3\2\3\2\3\2\3\2\3\2\3\2\eeee
\2\3\2\3\2\3\2\3\2\3\2\3\2\3\2\3\eeee
\4\1\4\1\4\1\4\1\4\1\4\1\4\1\4\1\eeee
\1\4\1\4\1\4\1\4\1\4\1\4\1\4\1\4\eeee
\3\2\3\2\3\2\3\2\3\2\3\2\3\2\3\2\eeee
\2\3\2\3\2\3\2\3\2\3\2\3\2\3\2\3\eeee
\4\1\4\1\4\1\4\1\4\1\4\1\4\1\4\1\eeee
} 

\baaa
4-16
\eaaa
\bbbb
0&0&1&3\\
0&1&0&3\\
1&0&0&3\\
1&1&1&1\\
\ebbb
\parbox{7cm}{ 
\1\4\2\4\3\4\1\4\2\4\3\4\1\4\2\4\eeee
\3\4\2\4\1\4\3\4\2\4\1\4\3\4\2\4\eeee
\4\1\4\3\4\2\4\1\4\3\4\2\4\1\4\3\eeee
\4\3\4\1\4\2\4\3\4\1\4\2\4\3\4\1\eeee
\1\4\2\4\3\4\1\4\2\4\3\4\1\4\2\4\eeee
\3\4\2\4\1\4\3\4\2\4\1\4\3\4\2\4\eeee
\4\1\4\3\4\2\4\1\4\3\4\2\4\1\4\3\eeee
\4\3\4\1\4\2\4\3\4\1\4\2\4\3\4\1\eeee
\1\4\2\4\3\4\1\4\2\4\3\4\1\4\2\4\eeee
\3\4\2\4\1\4\3\4\2\4\1\4\3\4\2\4\eeee
\4\1\4\3\4\2\4\1\4\3\4\2\4\1\4\3\eeee
\4\3\4\1\4\2\4\3\4\1\4\2\4\3\4\1\eeee
} 

\baaa
4-17
\eaaa
\bbbb
0&0&1&3\\
0&1&2&1\\
1&2&1&0\\
3&1&0&0\\
\ebbb
\parbox{7cm}{ 
\1\4\1\4\1\4\1\4\1\4\1\4\1\4\1\4\eeee
\3\2\3\2\3\2\3\2\3\2\3\2\3\2\3\2\eeee
\3\2\3\2\3\2\3\2\3\2\3\2\3\2\3\2\eeee
\1\4\1\4\1\4\1\4\1\4\1\4\1\4\1\4\eeee
\4\1\4\1\4\1\4\1\4\1\4\1\4\1\4\1\eeee
\2\3\2\3\2\3\2\3\2\3\2\3\2\3\2\3\eeee
\2\3\2\3\2\3\2\3\2\3\2\3\2\3\2\3\eeee
\4\1\4\1\4\1\4\1\4\1\4\1\4\1\4\1\eeee
\1\4\1\4\1\4\1\4\1\4\1\4\1\4\1\4\eeee
\3\2\3\2\3\2\3\2\3\2\3\2\3\2\3\2\eeee
\3\2\3\2\3\2\3\2\3\2\3\2\3\2\3\2\eeee
\1\4\1\4\1\4\1\4\1\4\1\4\1\4\1\4\eeee
} 

\baaa
4-18
\eaaa
\bbbb
0&0&1&3\\
0&2&1&1\\
1&3&0&0\\
1&1&0&2\\
\ebbb
\parbox{7cm}{ 
\1\4\4\4\1\4\4\4\1\4\4\4\1\4\4\4\eeee
\3\2\2\2\3\2\2\2\3\2\2\2\3\2\2\2\eeee
\2\2\3\2\2\2\3\2\2\2\3\2\2\2\3\2\eeee
\4\4\1\4\4\4\1\4\4\4\1\4\4\4\1\4\eeee
\1\4\4\4\1\4\4\4\1\4\4\4\1\4\4\4\eeee
\3\2\2\2\3\2\2\2\3\2\2\2\3\2\2\2\eeee
\2\2\3\2\2\2\3\2\2\2\3\2\2\2\3\2\eeee
\4\4\1\4\4\4\1\4\4\4\1\4\4\4\1\4\eeee
\1\4\4\4\1\4\4\4\1\4\4\4\1\4\4\4\eeee
\3\2\2\2\3\2\2\2\3\2\2\2\3\2\2\2\eeee
\2\2\3\2\2\2\3\2\2\2\3\2\2\2\3\2\eeee
\4\4\1\4\4\4\1\4\4\4\1\4\4\4\1\4\eeee
} 

\baaa
4-19
\eaaa
\bbbb
0&0&2&2\\
0&0&2&2\\
1&1&0&2\\
1&1&2&0\\
\ebbb
\parbox{7cm}{ 
\1\3\4\1\3\4\2\3\4\2\3\4\1\3\4\1\eeee
\3\4\2\3\4\1\3\4\1\3\4\2\3\4\2\3\eeee
\4\1\3\4\2\3\4\2\3\4\1\3\4\1\3\4\eeee
\2\3\4\1\3\4\1\3\4\2\3\4\2\3\4\1\eeee
\3\4\2\3\4\2\3\4\1\3\4\1\3\4\2\3\eeee
\4\1\3\4\1\3\4\2\3\4\2\3\4\1\3\4\eeee
\2\3\4\2\3\4\1\3\4\1\3\4\2\3\4\2\eeee
\3\4\1\3\4\2\3\4\2\3\4\1\3\4\1\3\eeee
\4\2\3\4\1\3\4\1\3\4\2\3\4\2\3\4\eeee
\1\3\4\2\3\4\2\3\4\1\3\4\1\3\4\2\eeee
\3\4\1\3\4\1\3\4\2\3\4\2\3\4\1\3\eeee
\4\2\3\4\2\3\4\1\3\4\1\3\4\2\3\4\eeee
} 

\baaa
4-20
\eaaa
\bbbb
0&0&2&2\\
0&0&2&2\\
1&1&1&1\\
1&1&1&1\\
\ebbb
\parbox{7cm}{ 
\1\3\3\1\3\3\1\3\3\1\3\3\1\3\3\1\eeee
\3\2\4\4\2\4\4\2\4\4\2\4\4\2\4\4\eeee
\3\4\1\3\3\1\3\3\1\3\3\1\3\3\1\3\eeee
\1\4\3\2\4\4\2\4\4\2\4\4\2\4\4\2\eeee
\3\2\3\4\1\3\3\1\3\3\1\3\3\1\3\3\eeee
\3\4\1\4\3\2\4\4\2\4\4\2\4\4\2\4\eeee
\1\4\3\2\3\4\1\3\3\1\3\3\1\3\3\1\eeee
\3\2\3\4\1\4\3\2\4\4\2\4\4\2\4\4\eeee
\3\4\1\4\3\2\3\4\1\3\3\1\3\3\1\3\eeee
\1\4\3\2\3\4\1\4\3\2\4\4\2\4\4\2\eeee
\3\2\3\4\1\4\3\2\3\4\1\3\3\1\3\3\eeee
\3\4\1\4\3\2\3\4\1\4\3\2\4\4\2\4\eeee
} 

\baaa
4-21
\eaaa
\bbbb
0&0&2&2\\
0&0&2&2\\
1&1&2&0\\
1&1&0&2\\
\ebbb
\parbox{7cm}{ 
\1\3\3\1\4\4\2\3\3\2\4\4\1\3\3\1\eeee
\3\3\2\4\4\1\3\3\1\4\4\2\3\3\2\4\eeee
\3\1\4\4\2\3\3\2\4\4\1\3\3\1\4\4\eeee
\2\4\4\1\3\3\1\4\4\2\3\3\2\4\4\1\eeee
\4\4\2\3\3\2\4\4\1\3\3\1\4\4\2\3\eeee
\4\1\3\3\1\4\4\2\3\3\2\4\4\1\3\3\eeee
\2\3\3\2\4\4\1\3\3\1\4\4\2\3\3\2\eeee
\3\3\1\4\4\2\3\3\2\4\4\1\3\3\1\4\eeee
\3\2\4\4\1\3\3\1\4\4\2\3\3\2\4\4\eeee
\1\4\4\2\3\3\2\4\4\1\3\3\1\4\4\2\eeee
\4\4\1\3\3\1\4\4\2\3\3\2\4\4\1\3\eeee
\4\2\3\3\2\4\4\1\3\3\1\4\4\2\3\3\eeee
} 

\baaa
4-22
\eaaa
\bbbb
0&0&2&2\\
0&0&2&2\\
1&1&2&0\\
2&2&0&0\\
\ebbb
\parbox{7cm}{ 
\1\3\3\1\4\2\3\3\1\4\2\3\3\2\4\1\eeee
\3\3\2\4\1\3\3\2\4\1\3\3\1\4\2\3\eeee
\3\1\4\2\3\3\1\4\2\3\3\2\4\1\3\3\eeee
\2\4\1\3\3\2\4\1\3\3\1\4\2\3\3\1\eeee
\4\2\3\3\1\4\2\3\3\2\4\1\3\3\2\4\eeee
\1\3\3\2\4\1\3\3\1\4\2\3\3\1\4\2\eeee
\3\3\1\4\2\3\3\2\4\1\3\3\2\4\1\3\eeee
\3\2\4\1\3\3\1\4\2\3\3\1\4\2\3\3\eeee
\1\4\2\3\3\2\4\1\3\3\2\4\1\3\3\1\eeee
\4\1\3\3\1\4\2\3\3\1\4\2\3\3\2\4\eeee
\2\3\3\2\4\1\3\3\2\4\1\3\3\1\4\2\eeee
\3\3\1\4\2\3\3\1\4\2\3\3\2\4\1\3\eeee
} 

\baaa
4-23
\eaaa
\bbbb
0&0&2&2\\
0&0&2&2\\
2&2&0&0\\
2&2&0&0\\
\ebbb
\parbox{7cm}{ 
\1\3\2\3\1\3\2\3\1\3\2\3\1\3\2\3\eeee
\3\1\4\2\4\1\4\2\4\1\4\2\4\1\4\2\eeee
\2\4\1\3\2\3\1\3\2\3\1\3\2\3\1\3\eeee
\3\2\3\1\4\2\4\1\4\2\4\1\4\2\4\1\eeee
\1\4\2\4\1\3\2\3\1\3\2\3\1\3\2\3\eeee
\3\1\3\2\3\1\4\2\4\1\4\2\4\1\4\2\eeee
\2\4\1\4\2\4\1\3\2\3\1\3\2\3\1\3\eeee
\3\2\3\1\3\2\3\1\4\2\4\1\4\2\4\1\eeee
\1\4\2\4\1\4\2\4\1\3\2\3\1\3\2\3\eeee
\3\1\3\2\3\1\3\2\3\1\4\2\4\1\4\2\eeee
\2\4\1\4\2\4\1\4\2\4\1\3\2\3\1\3\eeee
\3\2\3\1\3\2\3\1\3\2\3\1\4\2\4\1\eeee
} 

\baaa
4-24
\eaaa
\bbbb
0&0&2&2\\
0&1&1&2\\
1&1&1&1\\
1&2&1&0\\
\ebbb
\parbox{7cm}{ 
\1\4\2\3\3\2\4\1\4\2\3\3\2\4\1\4\eeee
\3\2\4\1\4\2\3\3\2\4\1\4\2\3\3\2\eeee
\4\2\3\3\2\4\1\4\2\3\3\2\4\1\4\2\eeee
\2\4\1\4\2\3\3\2\4\1\4\2\3\3\2\4\eeee
\2\3\3\2\4\1\4\2\3\3\2\4\1\4\2\3\eeee
\4\1\4\2\3\3\2\4\1\4\2\3\3\2\4\1\eeee
\3\3\2\4\1\4\2\3\3\2\4\1\4\2\3\3\eeee
\1\4\2\3\3\2\4\1\4\2\3\3\2\4\1\4\eeee
\3\2\4\1\4\2\3\3\2\4\1\4\2\3\3\2\eeee
\4\2\3\3\2\4\1\4\2\3\3\2\4\1\4\2\eeee
\2\4\1\4\2\3\3\2\4\1\4\2\3\3\2\4\eeee
\2\3\3\2\4\1\4\2\3\3\2\4\1\4\2\3\eeee
} 

\baaa
4-25
\eaaa
\bbbb
0&0&2&2\\
0&1&1&2\\
2&2&0&0\\
1&2&0&1\\
\ebbb
\parbox{7cm}{ 
\1\4\4\1\4\4\1\4\4\1\4\4\1\4\4\1\eeee
\3\2\2\3\2\2\3\2\2\3\2\2\3\2\2\3\eeee
\1\4\4\1\4\4\1\4\4\1\4\4\1\4\4\1\eeee
\3\2\2\3\2\2\3\2\2\3\2\2\3\2\2\3\eeee
\1\4\4\1\4\4\1\4\4\1\4\4\1\4\4\1\eeee
\3\2\2\3\2\2\3\2\2\3\2\2\3\2\2\3\eeee
\1\4\4\1\4\4\1\4\4\1\4\4\1\4\4\1\eeee
\3\2\2\3\2\2\3\2\2\3\2\2\3\2\2\3\eeee
\1\4\4\1\4\4\1\4\4\1\4\4\1\4\4\1\eeee
\3\2\2\3\2\2\3\2\2\3\2\2\3\2\2\3\eeee
\1\4\4\1\4\4\1\4\4\1\4\4\1\4\4\1\eeee
\3\2\2\3\2\2\3\2\2\3\2\2\3\2\2\3\eeee
} 

\baaa
4-26
\eaaa
\bbbb
0&0&2&2\\
0&2&0&2\\
2&0&2&0\\
2&2&0&0\\
\ebbb
\parbox{7cm}{ 
\1\3\3\1\4\2\2\4\1\3\3\1\4\2\2\4\eeee
\3\3\1\4\2\2\4\1\3\3\1\4\2\2\4\1\eeee
\3\1\4\2\2\4\1\3\3\1\4\2\2\4\1\3\eeee
\1\4\2\2\4\1\3\3\1\4\2\2\4\1\3\3\eeee
\4\2\2\4\1\3\3\1\4\2\2\4\1\3\3\1\eeee
\2\2\4\1\3\3\1\4\2\2\4\1\3\3\1\4\eeee
\2\4\1\3\3\1\4\2\2\4\1\3\3\1\4\2\eeee
\4\1\3\3\1\4\2\2\4\1\3\3\1\4\2\2\eeee
\1\3\3\1\4\2\2\4\1\3\3\1\4\2\2\4\eeee
\3\3\1\4\2\2\4\1\3\3\1\4\2\2\4\1\eeee
\3\1\4\2\2\4\1\3\3\1\4\2\2\4\1\3\eeee
\1\4\2\2\4\1\3\3\1\4\2\2\4\1\3\3\eeee
} 

\baaa
4-27
\eaaa
\bbbb
0&0&2&2\\
0&2&1&1\\
1&2&0&1\\
1&2&1&0\\
\ebbb
\parbox{7cm}{ 
\1\4\3\1\3\4\1\4\3\1\3\4\1\4\3\1\eeee
\3\2\2\4\2\2\3\2\2\4\2\2\3\2\2\4\eeee
\4\2\2\3\2\2\4\2\2\3\2\2\4\2\2\3\eeee
\1\3\4\1\4\3\1\3\4\1\4\3\1\3\4\1\eeee
\4\2\2\3\2\2\4\2\2\3\2\2\4\2\2\3\eeee
\3\2\2\4\2\2\3\2\2\4\2\2\3\2\2\4\eeee
\1\4\3\1\3\4\1\4\3\1\3\4\1\4\3\1\eeee
\3\2\2\4\2\2\3\2\2\4\2\2\3\2\2\4\eeee
\4\2\2\3\2\2\4\2\2\3\2\2\4\2\2\3\eeee
\1\3\4\1\4\3\1\3\4\1\4\3\1\3\4\1\eeee
\4\2\2\3\2\2\4\2\2\3\2\2\4\2\2\3\eeee
\3\2\2\4\2\2\3\2\2\4\2\2\3\2\2\4\eeee
} 

\baaa
4-28
\eaaa
\bbbb
0&0&2&2\\
0&2&1&1\\
1&2&1&0\\
1&2&0&1\\
\ebbb
\parbox{7cm}{ 
\1\4\4\1\4\4\1\4\4\1\4\4\1\4\4\1\eeee
\3\2\2\3\2\2\3\2\2\3\2\2\3\2\2\3\eeee
\3\2\2\3\2\2\3\2\2\3\2\2\3\2\2\3\eeee
\1\4\4\1\4\4\1\4\4\1\4\4\1\4\4\1\eeee
\3\2\2\3\2\2\3\2\2\3\2\2\3\2\2\3\eeee
\3\2\2\3\2\2\3\2\2\3\2\2\3\2\2\3\eeee
\1\4\4\1\4\4\1\4\4\1\4\4\1\4\4\1\eeee
\3\2\2\3\2\2\3\2\2\3\2\2\3\2\2\3\eeee
\3\2\2\3\2\2\3\2\2\3\2\2\3\2\2\3\eeee
\1\4\4\1\4\4\1\4\4\1\4\4\1\4\4\1\eeee
\3\2\2\3\2\2\3\2\2\3\2\2\3\2\2\3\eeee
\3\2\2\3\2\2\3\2\2\3\2\2\3\2\2\3\eeee
} 

\baaa
4-29
\eaaa
\bbbb
0&1&1&2\\
1&0&1&2\\
1&1&0&2\\
1&1&1&1\\
\ebbb
\parbox{7cm}{ 
\1\3\4\2\4\1\3\4\2\4\1\3\4\2\4\1\eeee
\2\4\1\3\4\2\4\1\3\4\2\4\1\3\4\2\eeee
\3\4\2\4\1\3\4\2\4\1\3\4\2\4\1\3\eeee
\4\1\3\4\2\4\1\3\4\2\4\1\3\4\2\4\eeee
\4\2\4\1\3\4\2\4\1\3\4\2\4\1\3\4\eeee
\1\3\4\2\4\1\3\4\2\4\1\3\4\2\4\1\eeee
\2\4\1\3\4\2\4\1\3\4\2\4\1\3\4\2\eeee
\3\4\2\4\1\3\4\2\4\1\3\4\2\4\1\3\eeee
\4\1\3\4\2\4\1\3\4\2\4\1\3\4\2\4\eeee
\4\2\4\1\3\4\2\4\1\3\4\2\4\1\3\4\eeee
\1\3\4\2\4\1\3\4\2\4\1\3\4\2\4\1\eeee
\2\4\1\3\4\2\4\1\3\4\2\4\1\3\4\2\eeee
} 

\baaa
4-30
\eaaa
\bbbb
0&1&1&2\\
1&0&2&1\\
1&2&0&1\\
2&1&1&0\\
\ebbb
\parbox{7cm}{ 
\1\3\2\4\1\3\2\4\1\3\2\4\1\3\2\4\eeee
\2\4\1\3\2\4\1\3\2\4\1\3\2\4\1\3\eeee
\3\1\4\2\3\1\4\2\3\1\4\2\3\1\4\2\eeee
\4\2\3\1\4\2\3\1\4\2\3\1\4\2\3\1\eeee
\1\3\2\4\1\3\2\4\1\3\2\4\1\3\2\4\eeee
\2\4\1\3\2\4\1\3\2\4\1\3\2\4\1\3\eeee
\3\1\4\2\3\1\4\2\3\1\4\2\3\1\4\2\eeee
\4\2\3\1\4\2\3\1\4\2\3\1\4\2\3\1\eeee
\1\3\2\4\1\3\2\4\1\3\2\4\1\3\2\4\eeee
\2\4\1\3\2\4\1\3\2\4\1\3\2\4\1\3\eeee
\3\1\4\2\3\1\4\2\3\1\4\2\3\1\4\2\eeee
\4\2\3\1\4\2\3\1\4\2\3\1\4\2\3\1\eeee
} 
\baaa
\phantom{0-0}\#
\eaaa
\mbox{}\phantom{\bbbb
0&0&0&0\\
\ebbb}
\parbox{7cm}{ 
\1\4\1\4\1\4\1\4\1\4\1\4\1\4\1\4\eeee
\2\3\2\3\2\3\2\3\2\3\2\3\2\3\2\3\eeee
\4\1\4\1\4\1\4\1\4\1\4\1\4\1\4\1\eeee
\3\2\3\2\3\2\3\2\3\2\3\2\3\2\3\2\eeee
\1\4\1\4\1\4\1\4\1\4\1\4\1\4\1\4\eeee
\2\3\2\3\2\3\2\3\2\3\2\3\2\3\2\3\eeee
\4\1\4\1\4\1\4\1\4\1\4\1\4\1\4\1\eeee
\3\2\3\2\3\2\3\2\3\2\3\2\3\2\3\2\eeee
\1\4\1\4\1\4\1\4\1\4\1\4\1\4\1\4\eeee
\2\3\2\3\2\3\2\3\2\3\2\3\2\3\2\3\eeee
\4\1\4\1\4\1\4\1\4\1\4\1\4\1\4\1\eeee
\3\2\3\2\3\2\3\2\3\2\3\2\3\2\3\2\eeee
} 

\baaa
4-31
\eaaa
\bbbb
0&1&1&2\\
1&0&2&1\\
1&2&1&0\\
2&1&0&1\\
\ebbb
\parbox{7cm}{ 
\1\3\2\4\1\3\2\4\1\3\2\4\1\3\2\4\eeee
\2\3\1\4\2\3\1\4\2\3\1\4\2\3\1\4\eeee
\3\2\4\1\3\2\4\1\3\2\4\1\3\2\4\1\eeee
\3\1\4\2\3\1\4\2\3\1\4\2\3\1\4\2\eeee
\2\4\1\3\2\4\1\3\2\4\1\3\2\4\1\3\eeee
\1\4\2\3\1\4\2\3\1\4\2\3\1\4\2\3\eeee
\4\1\3\2\4\1\3\2\4\1\3\2\4\1\3\2\eeee
\4\2\3\1\4\2\3\1\4\2\3\1\4\2\3\1\eeee
\1\3\2\4\1\3\2\4\1\3\2\4\1\3\2\4\eeee
\2\3\1\4\2\3\1\4\2\3\1\4\2\3\1\4\eeee
\3\2\4\1\3\2\4\1\3\2\4\1\3\2\4\1\eeee
\3\1\4\2\3\1\4\2\3\1\4\2\3\1\4\2\eeee
} 

\baaa
4-32
\eaaa
\bbbb
0&1&1&2\\
1&1&0&2\\
1&0&1&2\\
1&1&1&1\\
\ebbb
\parbox{7cm}{ 
\1\3\4\2\4\1\2\4\3\4\1\3\4\2\4\1\eeee
\2\4\1\2\4\3\4\1\3\4\2\4\1\2\4\3\eeee
\2\4\3\4\1\3\4\2\4\1\2\4\3\4\1\3\eeee
\4\1\3\4\2\4\1\2\4\3\4\1\3\4\2\4\eeee
\4\2\4\1\2\4\3\4\1\3\4\2\4\1\2\4\eeee
\1\2\4\3\4\1\3\4\2\4\1\2\4\3\4\1\eeee
\3\4\1\3\4\2\4\1\2\4\3\4\1\3\4\2\eeee
\3\4\2\4\1\2\4\3\4\1\3\4\2\4\1\2\eeee
\4\1\2\4\3\4\1\3\4\2\4\1\2\4\3\4\eeee
\4\3\4\1\3\4\2\4\1\2\4\3\4\1\3\4\eeee
\1\3\4\2\4\1\2\4\3\4\1\3\4\2\4\1\eeee
\2\4\1\2\4\3\4\1\3\4\2\4\1\2\4\3\eeee
} 
\baaa
\phantom{0-0}\#
\eaaa
\mbox{}\phantom{\bbbb
0&0&0&0\\
\ebbb}
\parbox{7cm}{ 
\1\4\3\3\4\1\4\2\2\4\1\4\3\3\4\1\eeee
\2\2\4\1\4\3\3\4\1\4\2\2\4\1\4\3\eeee
\4\1\4\2\2\4\1\4\3\3\4\1\4\2\2\4\eeee
\4\3\3\4\1\4\2\2\4\1\4\3\3\4\1\4\eeee
\2\4\1\4\3\3\4\1\4\2\2\4\1\4\3\3\eeee
\1\4\2\2\4\1\4\3\3\4\1\4\2\2\4\1\eeee
\3\3\4\1\4\2\2\4\1\4\3\3\4\1\4\2\eeee
\4\1\4\3\3\4\1\4\2\2\4\1\4\3\3\4\eeee
\4\2\2\4\1\4\3\3\4\1\4\2\2\4\1\4\eeee
\3\4\1\4\2\2\4\1\4\3\3\4\1\4\2\2\eeee
\1\4\3\3\4\1\4\2\2\4\1\4\3\3\4\1\eeee
\2\2\4\1\4\3\3\4\1\4\2\2\4\1\4\3\eeee
} 

\baaa
4-33
\eaaa
\bbbb
0&1&1&2\\
1&1&1&1\\
1&1&1&1\\
2&1&1&0\\
\ebbb
\parbox{7cm}{ 
\1\3\2\4\1\3\2\4\1\3\2\4\1\3\2\4\eeee
\2\4\1\3\2\4\1\3\2\4\1\3\2\4\1\3\eeee
\2\1\4\3\2\1\4\3\2\1\4\3\2\1\4\3\eeee
\4\3\2\1\4\3\2\1\4\3\2\1\4\3\2\1\eeee
\1\3\2\4\1\3\2\4\1\3\2\4\1\3\2\4\eeee
\2\4\1\3\2\4\1\3\2\4\1\3\2\4\1\3\eeee
\2\1\4\3\2\1\4\3\2\1\4\3\2\1\4\3\eeee
\4\3\2\1\4\3\2\1\4\3\2\1\4\3\2\1\eeee
\1\3\2\4\1\3\2\4\1\3\2\4\1\3\2\4\eeee
\2\4\1\3\2\4\1\3\2\4\1\3\2\4\1\3\eeee
\2\1\4\3\2\1\4\3\2\1\4\3\2\1\4\3\eeee
\4\3\2\1\4\3\2\1\4\3\2\1\4\3\2\1\eeee
} 

\baaa
4-34
\eaaa
\bbbb
0&1&1&2\\
1&2&0&1\\
1&0&2&1\\
2&1&1&0\\
\ebbb
\parbox{7cm}{ 
\1\3\3\4\1\3\3\4\1\3\3\4\1\3\3\4\eeee
\2\4\1\2\2\4\1\2\2\4\1\2\2\4\1\2\eeee
\2\1\4\2\2\1\4\2\2\1\4\2\2\1\4\2\eeee
\4\3\3\1\4\3\3\1\4\3\3\1\4\3\3\1\eeee
\1\3\3\4\1\3\3\4\1\3\3\4\1\3\3\4\eeee
\2\4\1\2\2\4\1\2\2\4\1\2\2\4\1\2\eeee
\2\1\4\2\2\1\4\2\2\1\4\2\2\1\4\2\eeee
\4\3\3\1\4\3\3\1\4\3\3\1\4\3\3\1\eeee
\1\3\3\4\1\3\3\4\1\3\3\4\1\3\3\4\eeee
\2\4\1\2\2\4\1\2\2\4\1\2\2\4\1\2\eeee
\2\1\4\2\2\1\4\2\2\1\4\2\2\1\4\2\eeee
\4\3\3\1\4\3\3\1\4\3\3\1\4\3\3\1\eeee
} 
\baaa
\phantom{0-0}\#
\eaaa
\mbox{}\phantom{\bbbb
0&0&0&0\\
\ebbb}
\parbox{7cm}{ 
\1\4\1\4\1\4\1\4\1\4\1\4\1\4\1\4\eeee
\2\2\2\2\2\2\2\2\2\2\2\2\2\2\2\2\eeee
\4\1\4\1\4\1\4\1\4\1\4\1\4\1\4\1\eeee
\3\3\3\3\3\3\3\3\3\3\3\3\3\3\3\3\eeee
\1\4\1\4\1\4\1\4\1\4\1\4\1\4\1\4\eeee
\2\2\2\2\2\2\2\2\2\2\2\2\2\2\2\2\eeee
\4\1\4\1\4\1\4\1\4\1\4\1\4\1\4\1\eeee
\3\3\3\3\3\3\3\3\3\3\3\3\3\3\3\3\eeee
\1\4\1\4\1\4\1\4\1\4\1\4\1\4\1\4\eeee
\2\2\2\2\2\2\2\2\2\2\2\2\2\2\2\2\eeee
\4\1\4\1\4\1\4\1\4\1\4\1\4\1\4\1\eeee
\3\3\3\3\3\3\3\3\3\3\3\3\3\3\3\3\eeee
} 

\baaa
4-35
\eaaa
\bbbb
0&1&1&2\\
2&1&1&0\\
2&1&1&0\\
2&0&0&2\\
\ebbb
\parbox{7cm}{ 
\1\3\3\1\4\4\1\2\2\1\4\4\1\3\3\1\eeee
\2\2\1\4\4\1\3\3\1\4\4\1\2\2\1\4\eeee
\3\1\4\4\1\2\2\1\4\4\1\3\3\1\4\4\eeee
\1\4\4\1\3\3\1\4\4\1\2\2\1\4\4\1\eeee
\4\4\1\2\2\1\4\4\1\3\3\1\4\4\1\2\eeee
\4\1\3\3\1\4\4\1\2\2\1\4\4\1\3\3\eeee
\1\2\2\1\4\4\1\3\3\1\4\4\1\2\2\1\eeee
\3\3\1\4\4\1\2\2\1\4\4\1\3\3\1\4\eeee
\2\1\4\4\1\3\3\1\4\4\1\2\2\1\4\4\eeee
\1\4\4\1\2\2\1\4\4\1\3\3\1\4\4\1\eeee
\4\4\1\3\3\1\4\4\1\2\2\1\4\4\1\3\eeee
\4\1\2\2\1\4\4\1\3\3\1\4\4\1\2\2\eeee
} 

\baaa
4-36
\eaaa
\bbbb
1&0&0&3\\
0&1&0&3\\
0&0&1&3\\
1&1&1&1\\
\ebbb
\parbox{7cm}{ 
\1\4\3\4\2\4\1\4\3\4\2\4\1\4\3\4\eeee
\1\4\3\4\2\4\1\4\3\4\2\4\1\4\3\4\eeee
\4\2\4\1\4\3\4\2\4\1\4\3\4\2\4\1\eeee
\4\2\4\1\4\3\4\2\4\1\4\3\4\2\4\1\eeee
\1\4\3\4\2\4\1\4\3\4\2\4\1\4\3\4\eeee
\1\4\3\4\2\4\1\4\3\4\2\4\1\4\3\4\eeee
\4\2\4\1\4\3\4\2\4\1\4\3\4\2\4\1\eeee
\4\2\4\1\4\3\4\2\4\1\4\3\4\2\4\1\eeee
\1\4\3\4\2\4\1\4\3\4\2\4\1\4\3\4\eeee
\1\4\3\4\2\4\1\4\3\4\2\4\1\4\3\4\eeee
\4\2\4\1\4\3\4\2\4\1\4\3\4\2\4\1\eeee
\4\2\4\1\4\3\4\2\4\1\4\3\4\2\4\1\eeee
} 

\baaa
4-37
\eaaa
\bbbb
1&0&0&3\\
0&1&1&2\\
0&3&1&0\\
1&2&0&1\\
\ebbb
\parbox{7cm}{ 
\1\4\2\4\1\4\2\4\1\4\2\4\1\4\2\4\eeee
\1\4\2\4\1\4\2\4\1\4\2\4\1\4\2\4\eeee
\4\2\3\2\4\2\3\2\4\2\3\2\4\2\3\2\eeee
\4\2\3\2\4\2\3\2\4\2\3\2\4\2\3\2\eeee
\1\4\2\4\1\4\2\4\1\4\2\4\1\4\2\4\eeee
\1\4\2\4\1\4\2\4\1\4\2\4\1\4\2\4\eeee
\4\2\3\2\4\2\3\2\4\2\3\2\4\2\3\2\eeee
\4\2\3\2\4\2\3\2\4\2\3\2\4\2\3\2\eeee
\1\4\2\4\1\4\2\4\1\4\2\4\1\4\2\4\eeee
\1\4\2\4\1\4\2\4\1\4\2\4\1\4\2\4\eeee
\4\2\3\2\4\2\3\2\4\2\3\2\4\2\3\2\eeee
\4\2\3\2\4\2\3\2\4\2\3\2\4\2\3\2\eeee
} 

\baaa
4-38
\eaaa
\bbbb
1&0&1&2\\
0&1&1&2\\
1&2&1&0\\
1&2&0&1\\
\ebbb
\parbox{7cm}{ 
\1\3\2\4\2\4\1\3\2\4\2\4\1\3\2\4\eeee
\1\3\2\4\2\4\1\3\2\4\2\4\1\3\2\4\eeee
\4\2\4\1\3\2\4\2\4\1\3\2\4\2\4\1\eeee
\4\2\4\1\3\2\4\2\4\1\3\2\4\2\4\1\eeee
\1\3\2\4\2\4\1\3\2\4\2\4\1\3\2\4\eeee
\1\3\2\4\2\4\1\3\2\4\2\4\1\3\2\4\eeee
\4\2\4\1\3\2\4\2\4\1\3\2\4\2\4\1\eeee
\4\2\4\1\3\2\4\2\4\1\3\2\4\2\4\1\eeee
\1\3\2\4\2\4\1\3\2\4\2\4\1\3\2\4\eeee
\1\3\2\4\2\4\1\3\2\4\2\4\1\3\2\4\eeee
\4\2\4\1\3\2\4\2\4\1\3\2\4\2\4\1\eeee
\4\2\4\1\3\2\4\2\4\1\3\2\4\2\4\1\eeee
} 
\baaa
\phantom{0-0}\#
\eaaa
\mbox{}\phantom{\bbbb
0&0&0&0\\
\ebbb}
\parbox{7cm}{ 
\1\4\2\3\2\4\1\4\2\3\2\4\1\4\2\3\eeee
\1\4\2\3\2\4\1\4\2\3\2\4\1\4\2\3\eeee
\3\2\4\1\4\2\3\2\4\1\4\2\3\2\4\1\eeee
\3\2\4\1\4\2\3\2\4\1\4\2\3\2\4\1\eeee
\1\4\2\3\2\4\1\4\2\3\2\4\1\4\2\3\eeee
\1\4\2\3\2\4\1\4\2\3\2\4\1\4\2\3\eeee
\3\2\4\1\4\2\3\2\4\1\4\2\3\2\4\1\eeee
\3\2\4\1\4\2\3\2\4\1\4\2\3\2\4\1\eeee
\1\4\2\3\2\4\1\4\2\3\2\4\1\4\2\3\eeee
\1\4\2\3\2\4\1\4\2\3\2\4\1\4\2\3\eeee
\3\2\4\1\4\2\3\2\4\1\4\2\3\2\4\1\eeee
\3\2\4\1\4\2\3\2\4\1\4\2\3\2\4\1\eeee
} 

\baaa
4-39
\eaaa
\bbbb
1&0&1&2\\
0&1&2&1\\
1&2&1&0\\
2&1&0&1\\
\ebbb
\parbox{7cm}{ 
\1\3\2\4\1\3\2\4\1\3\2\4\1\3\2\4\eeee
\1\3\2\4\1\3\2\4\1\3\2\4\1\3\2\4\eeee
\4\2\3\1\4\2\3\1\4\2\3\1\4\2\3\1\eeee
\4\2\3\1\4\2\3\1\4\2\3\1\4\2\3\1\eeee
\1\3\2\4\1\3\2\4\1\3\2\4\1\3\2\4\eeee
\1\3\2\4\1\3\2\4\1\3\2\4\1\3\2\4\eeee
\4\2\3\1\4\2\3\1\4\2\3\1\4\2\3\1\eeee
\4\2\3\1\4\2\3\1\4\2\3\1\4\2\3\1\eeee
\1\3\2\4\1\3\2\4\1\3\2\4\1\3\2\4\eeee
\1\3\2\4\1\3\2\4\1\3\2\4\1\3\2\4\eeee
\4\2\3\1\4\2\3\1\4\2\3\1\4\2\3\1\eeee
\4\2\3\1\4\2\3\1\4\2\3\1\4\2\3\1\eeee
} 
\baaa
\phantom{0-0}\#
\eaaa
\mbox{}\phantom{\bbbb
0&0&0&0\\
\ebbb}
\parbox{7cm}{ 
\1\4\1\4\1\4\1\4\1\4\1\4\1\4\1\4\eeee
\1\4\1\4\1\4\1\4\1\4\1\4\1\4\1\4\eeee
\3\2\3\2\3\2\3\2\3\2\3\2\3\2\3\2\eeee
\3\2\3\2\3\2\3\2\3\2\3\2\3\2\3\2\eeee
\1\4\1\4\1\4\1\4\1\4\1\4\1\4\1\4\eeee
\1\4\1\4\1\4\1\4\1\4\1\4\1\4\1\4\eeee
\3\2\3\2\3\2\3\2\3\2\3\2\3\2\3\2\eeee
\3\2\3\2\3\2\3\2\3\2\3\2\3\2\3\2\eeee
\1\4\1\4\1\4\1\4\1\4\1\4\1\4\1\4\eeee
\1\4\1\4\1\4\1\4\1\4\1\4\1\4\1\4\eeee
\3\2\3\2\3\2\3\2\3\2\3\2\3\2\3\2\eeee
\3\2\3\2\3\2\3\2\3\2\3\2\3\2\3\2\eeee
} 

\baaa
4-40
\eaaa
\bbbb
1&0&1&2\\
0&2&1&1\\
1&1&2&0\\
2&1&0&1\\
\ebbb
\parbox{7cm}{ 
\1\3\3\1\4\2\2\4\1\3\3\1\4\2\2\4\eeee
\1\3\3\1\4\2\2\4\1\3\3\1\4\2\2\4\eeee
\4\2\2\4\1\3\3\1\4\2\2\4\1\3\3\1\eeee
\4\2\2\4\1\3\3\1\4\2\2\4\1\3\3\1\eeee
\1\3\3\1\4\2\2\4\1\3\3\1\4\2\2\4\eeee
\1\3\3\1\4\2\2\4\1\3\3\1\4\2\2\4\eeee
\4\2\2\4\1\3\3\1\4\2\2\4\1\3\3\1\eeee
\4\2\2\4\1\3\3\1\4\2\2\4\1\3\3\1\eeee
\1\3\3\1\4\2\2\4\1\3\3\1\4\2\2\4\eeee
\1\3\3\1\4\2\2\4\1\3\3\1\4\2\2\4\eeee
\4\2\2\4\1\3\3\1\4\2\2\4\1\3\3\1\eeee
\4\2\2\4\1\3\3\1\4\2\2\4\1\3\3\1\eeee
} 

\baaa
4-41
\eaaa
\bbbb
1&0&1&2\\
0&2&1&1\\
1&2&1&0\\
1&1&0&2\\
\ebbb
\parbox{7cm}{ 
\1\3\2\4\1\3\2\4\1\3\2\4\1\3\2\4\eeee
\1\3\2\4\1\3\2\4\1\3\2\4\1\3\2\4\eeee
\4\2\2\4\4\2\2\4\4\2\2\4\4\2\2\4\eeee
\4\2\3\1\4\2\3\1\4\2\3\1\4\2\3\1\eeee
\4\2\3\1\4\2\3\1\4\2\3\1\4\2\3\1\eeee
\4\2\2\4\4\2\2\4\4\2\2\4\4\2\2\4\eeee
\1\3\2\4\1\3\2\4\1\3\2\4\1\3\2\4\eeee
\1\3\2\4\1\3\2\4\1\3\2\4\1\3\2\4\eeee
\4\2\2\4\4\2\2\4\4\2\2\4\4\2\2\4\eeee
\4\2\3\1\4\2\3\1\4\2\3\1\4\2\3\1\eeee
\4\2\3\1\4\2\3\1\4\2\3\1\4\2\3\1\eeee
\4\2\2\4\4\2\2\4\4\2\2\4\4\2\2\4\eeee
} 
\baaa
\phantom{0-0}\#
\eaaa
\mbox{}\phantom{\bbbb
0&0&0&0\\
\ebbb}
\parbox{7cm}{ 
\1\4\4\1\4\4\1\4\4\1\4\4\1\4\4\1\eeee
\1\4\4\1\4\4\1\4\4\1\4\4\1\4\4\1\eeee
\3\2\2\3\2\2\3\2\2\3\2\2\3\2\2\3\eeee
\3\2\2\3\2\2\3\2\2\3\2\2\3\2\2\3\eeee
\1\4\4\1\4\4\1\4\4\1\4\4\1\4\4\1\eeee
\1\4\4\1\4\4\1\4\4\1\4\4\1\4\4\1\eeee
\3\2\2\3\2\2\3\2\2\3\2\2\3\2\2\3\eeee
\3\2\2\3\2\2\3\2\2\3\2\2\3\2\2\3\eeee
\1\4\4\1\4\4\1\4\4\1\4\4\1\4\4\1\eeee
\1\4\4\1\4\4\1\4\4\1\4\4\1\4\4\1\eeee
\3\2\2\3\2\2\3\2\2\3\2\2\3\2\2\3\eeee
\3\2\2\3\2\2\3\2\2\3\2\2\3\2\2\3\eeee
} 

\baaa
4-42
\eaaa
\bbbb
1&1&1&1\\
1&1&1&1\\
1&1&1&1\\
1&1&1&1\\
\ebbb
\parbox{7cm}{ 
\1\2\3\1\1\2\3\1\1\2\3\1\1\2\3\1\eeee
\1\2\4\4\3\2\4\4\3\2\4\4\3\2\4\4\eeee
\3\3\1\2\3\1\1\2\3\1\1\2\3\1\1\2\eeee
\2\4\1\2\4\4\3\2\4\4\3\2\4\4\3\2\eeee
\1\4\3\3\1\2\3\1\1\2\3\1\1\2\3\1\eeee
\1\2\2\4\1\2\4\4\3\2\4\4\3\2\4\4\eeee
\3\3\1\4\3\3\1\2\3\1\1\2\3\1\1\2\eeee
\2\4\1\2\2\4\1\2\4\4\3\2\4\4\3\2\eeee
\1\4\3\3\1\4\3\3\1\2\3\1\1\2\3\1\eeee
\1\2\2\4\1\2\2\4\1\2\4\4\3\2\4\4\eeee
\3\3\1\4\3\3\1\4\3\3\1\2\3\1\1\2\eeee
\2\4\1\2\2\4\1\2\2\4\1\2\4\4\3\2\eeee
} 

\baaa
4-43
\eaaa
\bbbb
1&1&1&1\\
1&1&1&1\\
1&1&2&0\\
1&1&0&2\\
\ebbb
\parbox{7cm}{ 
\1\2\4\4\1\2\4\4\1\2\4\4\1\2\4\4\eeee
\1\2\4\4\1\2\4\4\1\2\4\4\1\2\4\4\eeee
\3\3\1\2\3\3\1\2\3\3\1\2\3\3\1\2\eeee
\3\3\1\2\3\3\1\2\3\3\1\2\3\3\1\2\eeee
\1\2\4\4\1\2\4\4\1\2\4\4\1\2\4\4\eeee
\1\2\4\4\1\2\4\4\1\2\4\4\1\2\4\4\eeee
\3\3\1\2\3\3\1\2\3\3\1\2\3\3\1\2\eeee
\3\3\1\2\3\3\1\2\3\3\1\2\3\3\1\2\eeee
\1\2\4\4\1\2\4\4\1\2\4\4\1\2\4\4\eeee
\1\2\4\4\1\2\4\4\1\2\4\4\1\2\4\4\eeee
\3\3\1\2\3\3\1\2\3\3\1\2\3\3\1\2\eeee
\3\3\1\2\3\3\1\2\3\3\1\2\3\3\1\2\eeee
} 

\baaa
4-44
\eaaa
\bbbb
2&0&0&2\\
0&2&1&1\\
0&1&3&0\\
1&1&0&2\\
\ebbb
\parbox{7cm}{ 
\1\4\2\3\3\2\4\1\4\2\3\3\2\4\1\4\eeee
\1\4\2\3\3\2\4\1\4\2\3\3\2\4\1\4\eeee
\1\4\2\3\3\2\4\1\4\2\3\3\2\4\1\4\eeee
\1\4\2\3\3\2\4\1\4\2\3\3\2\4\1\4\eeee
\1\4\2\3\3\2\4\1\4\2\3\3\2\4\1\4\eeee
\1\4\2\3\3\2\4\1\4\2\3\3\2\4\1\4\eeee
\1\4\2\3\3\2\4\1\4\2\3\3\2\4\1\4\eeee
\1\4\2\3\3\2\4\1\4\2\3\3\2\4\1\4\eeee
\1\4\2\3\3\2\4\1\4\2\3\3\2\4\1\4\eeee
\1\4\2\3\3\2\4\1\4\2\3\3\2\4\1\4\eeee
\1\4\2\3\3\2\4\1\4\2\3\3\2\4\1\4\eeee
\1\4\2\3\3\2\4\1\4\2\3\3\2\4\1\4\eeee
} 

\baaa
4-45
\eaaa
\bbbb
2&0&0&2\\
0&2&1&1\\
0&2&2&0\\
1&1&0&2\\
\ebbb
\parbox{7cm}{ 
\1\4\2\3\2\4\1\4\2\3\2\4\1\4\2\3\eeee
\1\4\2\3\2\4\1\4\2\3\2\4\1\4\2\3\eeee
\1\4\2\3\2\4\1\4\2\3\2\4\1\4\2\3\eeee
\1\4\2\3\2\4\1\4\2\3\2\4\1\4\2\3\eeee
\1\4\2\3\2\4\1\4\2\3\2\4\1\4\2\3\eeee
\1\4\2\3\2\4\1\4\2\3\2\4\1\4\2\3\eeee
\1\4\2\3\2\4\1\4\2\3\2\4\1\4\2\3\eeee
\1\4\2\3\2\4\1\4\2\3\2\4\1\4\2\3\eeee
\1\4\2\3\2\4\1\4\2\3\2\4\1\4\2\3\eeee
\1\4\2\3\2\4\1\4\2\3\2\4\1\4\2\3\eeee
\1\4\2\3\2\4\1\4\2\3\2\4\1\4\2\3\eeee
\1\4\2\3\2\4\1\4\2\3\2\4\1\4\2\3\eeee
} 

\baaa
4-46
\eaaa
\bbbb
2&0&1&1\\
0&2&1&1\\
1&1&2&0\\
1&1&0&2\\
\ebbb
\parbox{7cm}{ 
\1\1\4\4\1\1\4\4\1\1\4\4\1\1\4\4\eeee
\1\1\4\4\1\1\4\4\1\1\4\4\1\1\4\4\eeee
\3\3\2\2\3\3\2\2\3\3\2\2\3\3\2\2\eeee
\3\3\2\2\3\3\2\2\3\3\2\2\3\3\2\2\eeee
\1\1\4\4\1\1\4\4\1\1\4\4\1\1\4\4\eeee
\1\1\4\4\1\1\4\4\1\1\4\4\1\1\4\4\eeee
\3\3\2\2\3\3\2\2\3\3\2\2\3\3\2\2\eeee
\3\3\2\2\3\3\2\2\3\3\2\2\3\3\2\2\eeee
\1\1\4\4\1\1\4\4\1\1\4\4\1\1\4\4\eeee
\1\1\4\4\1\1\4\4\1\1\4\4\1\1\4\4\eeee
\3\3\2\2\3\3\2\2\3\3\2\2\3\3\2\2\eeee
\3\3\2\2\3\3\2\2\3\3\2\2\3\3\2\2\eeee
} 
\baaa
\phantom{0-0}\#
\eaaa
\mbox{}\phantom{\bbbb
0&0&0&0\\
\ebbb}
\parbox{7cm}{ 
\1\3\2\4\1\3\2\4\1\3\2\4\1\3\2\4\eeee
\1\3\2\4\1\3\2\4\1\3\2\4\1\3\2\4\eeee
\1\3\2\4\1\3\2\4\1\3\2\4\1\3\2\4\eeee
\1\3\2\4\1\3\2\4\1\3\2\4\1\3\2\4\eeee
\1\3\2\4\1\3\2\4\1\3\2\4\1\3\2\4\eeee
\1\3\2\4\1\3\2\4\1\3\2\4\1\3\2\4\eeee
\1\3\2\4\1\3\2\4\1\3\2\4\1\3\2\4\eeee
\1\3\2\4\1\3\2\4\1\3\2\4\1\3\2\4\eeee
\1\3\2\4\1\3\2\4\1\3\2\4\1\3\2\4\eeee
\1\3\2\4\1\3\2\4\1\3\2\4\1\3\2\4\eeee
\1\3\2\4\1\3\2\4\1\3\2\4\1\3\2\4\eeee
\1\3\2\4\1\3\2\4\1\3\2\4\1\3\2\4\eeee
} 

\baaa
4-47
\eaaa
\bbbb
2&0&1&1\\
0&3&0&1\\
1&0&3&0\\
1&1&0&2\\
\ebbb
\parbox{7cm}{ 
\1\3\3\1\4\2\2\4\1\3\3\1\4\2\2\4\eeee
\1\3\3\1\4\2\2\4\1\3\3\1\4\2\2\4\eeee
\1\3\3\1\4\2\2\4\1\3\3\1\4\2\2\4\eeee
\1\3\3\1\4\2\2\4\1\3\3\1\4\2\2\4\eeee
\1\3\3\1\4\2\2\4\1\3\3\1\4\2\2\4\eeee
\1\3\3\1\4\2\2\4\1\3\3\1\4\2\2\4\eeee
\1\3\3\1\4\2\2\4\1\3\3\1\4\2\2\4\eeee
\1\3\3\1\4\2\2\4\1\3\3\1\4\2\2\4\eeee
\1\3\3\1\4\2\2\4\1\3\3\1\4\2\2\4\eeee
\1\3\3\1\4\2\2\4\1\3\3\1\4\2\2\4\eeee
\1\3\3\1\4\2\2\4\1\3\3\1\4\2\2\4\eeee
\1\3\3\1\4\2\2\4\1\3\3\1\4\2\2\4\eeee
} 

\baaa
5-1
\eaaa
\bbbb
0&0&0&0&4\\
0&0&0&0&4\\
0&0&0&0&4\\
0&0&0&0&4\\
1&1&1&1&0\\
\ebbb
\parbox{7cm}{ 
\1\5\3\5\1\5\3\5\1\5\3\5\1\5\3\5\eeee
\5\2\5\4\5\2\5\4\5\2\5\4\5\2\5\4\eeee
\3\5\1\5\3\5\1\5\3\5\1\5\3\5\1\5\eeee
\5\4\5\2\5\4\5\2\5\4\5\2\5\4\5\2\eeee
\1\5\3\5\1\5\3\5\1\5\3\5\1\5\3\5\eeee
\5\2\5\4\5\2\5\4\5\2\5\4\5\2\5\4\eeee
\3\5\1\5\3\5\1\5\3\5\1\5\3\5\1\5\eeee
\5\4\5\2\5\4\5\2\5\4\5\2\5\4\5\2\eeee
\1\5\3\5\1\5\3\5\1\5\3\5\1\5\3\5\eeee
\5\2\5\4\5\2\5\4\5\2\5\4\5\2\5\4\eeee
\3\5\1\5\3\5\1\5\3\5\1\5\3\5\1\5\eeee
\5\4\5\2\5\4\5\2\5\4\5\2\5\4\5\2\eeee
} 

\baaa
5-2
\eaaa
\bbbb
0&0&0&0&4\\
0&0&0&0&4\\
0&0&0&2&2\\
0&0&2&2&0\\
1&1&2&0&0\\
\ebbb
\parbox{7cm}{ 
\1\5\3\4\4\3\5\1\5\3\4\4\3\5\1\5\eeee
\5\2\5\3\4\4\3\5\2\5\3\4\4\3\5\2\eeee
\3\5\1\5\3\4\4\3\5\1\5\3\4\4\3\5\eeee
\4\3\5\2\5\3\4\4\3\5\2\5\3\4\4\3\eeee
\4\4\3\5\1\5\3\4\4\3\5\1\5\3\4\4\eeee
\3\4\4\3\5\2\5\3\4\4\3\5\2\5\3\4\eeee
\5\3\4\4\3\5\1\5\3\4\4\3\5\1\5\3\eeee
\1\5\3\4\4\3\5\2\5\3\4\4\3\5\2\5\eeee
\5\2\5\3\4\4\3\5\1\5\3\4\4\3\5\1\eeee
\3\5\1\5\3\4\4\3\5\2\5\3\4\4\3\5\eeee
\4\3\5\2\5\3\4\4\3\5\1\5\3\4\4\3\eeee
\4\4\3\5\1\5\3\4\4\3\5\2\5\3\4\4\eeee
} 

\baaa
5-3
\eaaa
\bbbb
0&0&0&0&4\\
0&0&0&0&4\\
0&0&0&2&2\\
0&0&4&0&0\\
1&1&2&0&0\\
\ebbb
\parbox{7cm}{ 
\1\5\3\4\3\5\1\5\3\4\3\5\1\5\3\4\eeee
\5\2\5\3\4\3\5\2\5\3\4\3\5\2\5\3\eeee
\3\5\1\5\3\4\3\5\1\5\3\4\3\5\1\5\eeee
\4\3\5\2\5\3\4\3\5\2\5\3\4\3\5\2\eeee
\3\4\3\5\1\5\3\4\3\5\1\5\3\4\3\5\eeee
\5\3\4\3\5\2\5\3\4\3\5\2\5\3\4\3\eeee
\1\5\3\4\3\5\1\5\3\4\3\5\1\5\3\4\eeee
\5\2\5\3\4\3\5\2\5\3\4\3\5\2\5\3\eeee
\3\5\1\5\3\4\3\5\1\5\3\4\3\5\1\5\eeee
\4\3\5\2\5\3\4\3\5\2\5\3\4\3\5\2\eeee
\3\4\3\5\1\5\3\4\3\5\1\5\3\4\3\5\eeee
\5\3\4\3\5\2\5\3\4\3\5\2\5\3\4\3\eeee
} 

\baaa
5-4
\eaaa
\bbbb
0&0&0&0&4\\
0&0&0&0&4\\
0&0&1&1&2\\
0&0&1&1&2\\
1&1&1&1&0\\
\ebbb
\parbox{7cm}{ 
\1\5\3\3\5\1\5\3\3\5\1\5\3\3\5\1\eeee
\5\2\5\4\4\5\2\5\4\4\5\2\5\4\4\5\eeee
\3\5\1\5\3\3\5\1\5\3\3\5\1\5\3\3\eeee
\3\4\5\2\5\4\4\5\2\5\4\4\5\2\5\4\eeee
\5\4\3\5\1\5\3\3\5\1\5\3\3\5\1\5\eeee
\1\5\3\4\5\2\5\4\4\5\2\5\4\4\5\2\eeee
\5\2\5\4\3\5\1\5\3\3\5\1\5\3\3\5\eeee
\3\5\1\5\3\4\5\2\5\4\4\5\2\5\4\4\eeee
\3\4\5\2\5\4\3\5\1\5\3\3\5\1\5\3\eeee
\5\4\3\5\1\5\3\4\5\2\5\4\4\5\2\5\eeee
\1\5\3\4\5\2\5\4\3\5\1\5\3\3\5\1\eeee
\5\2\5\4\3\5\1\5\3\4\5\2\5\4\4\5\eeee
} 

\baaa
5-5
\eaaa
\bbbb
0&0&0&0&4\\
0&0&0&2&2\\
0&0&0&2&2\\
0&1&1&2&0\\
2&1&1&0&0\\
\ebbb
\parbox{7cm}{ 
\1\5\2\4\4\2\5\1\5\2\4\4\2\5\1\5\eeee
\5\1\5\3\4\4\3\5\1\5\3\4\4\3\5\1\eeee
\2\5\1\5\2\4\4\2\5\1\5\2\4\4\2\5\eeee
\4\3\5\1\5\3\4\4\3\5\1\5\3\4\4\3\eeee
\4\4\2\5\1\5\2\4\4\2\5\1\5\2\4\4\eeee
\2\4\4\3\5\1\5\3\4\4\3\5\1\5\3\4\eeee
\5\3\4\4\2\5\1\5\2\4\4\2\5\1\5\2\eeee
\1\5\2\4\4\3\5\1\5\3\4\4\3\5\1\5\eeee
\5\1\5\3\4\4\2\5\1\5\2\4\4\2\5\1\eeee
\2\5\1\5\2\4\4\3\5\1\5\3\4\4\3\5\eeee
\4\3\5\1\5\3\4\4\2\5\1\5\2\4\4\2\eeee
\4\4\2\5\1\5\2\4\4\3\5\1\5\3\4\4\eeee
} 

\baaa
5-6
\eaaa
\bbbb
0&0&0&0&4\\
0&0&0&2&2\\
0&0&0&2&2\\
0&2&2&0&0\\
2&1&1&0&0\\
\ebbb
\parbox{7cm}{ 
\1\5\2\4\2\5\1\5\2\4\2\5\1\5\2\4\eeee
\5\1\5\3\4\3\5\1\5\3\4\3\5\1\5\3\eeee
\2\5\1\5\2\4\2\5\1\5\2\4\2\5\1\5\eeee
\4\3\5\1\5\3\4\3\5\1\5\3\4\3\5\1\eeee
\2\4\2\5\1\5\2\4\2\5\1\5\2\4\2\5\eeee
\5\3\4\3\5\1\5\3\4\3\5\1\5\3\4\3\eeee
\1\5\2\4\2\5\1\5\2\4\2\5\1\5\2\4\eeee
\5\1\5\3\4\3\5\1\5\3\4\3\5\1\5\3\eeee
\2\5\1\5\2\4\2\5\1\5\2\4\2\5\1\5\eeee
\4\3\5\1\5\3\4\3\5\1\5\3\4\3\5\1\eeee
\2\4\2\5\1\5\2\4\2\5\1\5\2\4\2\5\eeee
\5\3\4\3\5\1\5\3\4\3\5\1\5\3\4\3\eeee
} 

\baaa
5-7
\eaaa
\bbbb
0&0&0&0&4\\
0&0&0&2&2\\
0&0&0&4&0\\
0&2&2&0&0\\
2&2&0&0&0\\
\ebbb
\parbox{7cm}{ 
\1\5\2\4\3\4\2\5\1\5\2\4\3\4\2\5\eeee
\5\1\5\2\4\3\4\2\5\1\5\2\4\3\4\2\eeee
\2\5\1\5\2\4\3\4\2\5\1\5\2\4\3\4\eeee
\4\2\5\1\5\2\4\3\4\2\5\1\5\2\4\3\eeee
\3\4\2\5\1\5\2\4\3\4\2\5\1\5\2\4\eeee
\4\3\4\2\5\1\5\2\4\3\4\2\5\1\5\2\eeee
\2\4\3\4\2\5\1\5\2\4\3\4\2\5\1\5\eeee
\5\2\4\3\4\2\5\1\5\2\4\3\4\2\5\1\eeee
\1\5\2\4\3\4\2\5\1\5\2\4\3\4\2\5\eeee
\5\1\5\2\4\3\4\2\5\1\5\2\4\3\4\2\eeee
\2\5\1\5\2\4\3\4\2\5\1\5\2\4\3\4\eeee
\4\2\5\1\5\2\4\3\4\2\5\1\5\2\4\3\eeee
} 

\baaa
5-8
\eaaa
\bbbb
0&0&0&0&4\\
0&0&0&2&2\\
0&0&0&4&0\\
0&3&1&0&0\\
1&3&0&0&0\\
\ebbb
\parbox{7cm}{ 
\1\5\2\5\1\5\2\5\1\5\2\5\1\5\2\5\eeee
\5\2\4\2\5\2\4\2\5\2\4\2\5\2\4\2\eeee
\2\4\3\4\2\4\3\4\2\4\3\4\2\4\3\4\eeee
\5\2\4\2\5\2\4\2\5\2\4\2\5\2\4\2\eeee
\1\5\2\5\1\5\2\5\1\5\2\5\1\5\2\5\eeee
\5\2\4\2\5\2\4\2\5\2\4\2\5\2\4\2\eeee
\2\4\3\4\2\4\3\4\2\4\3\4\2\4\3\4\eeee
\5\2\4\2\5\2\4\2\5\2\4\2\5\2\4\2\eeee
\1\5\2\5\1\5\2\5\1\5\2\5\1\5\2\5\eeee
\5\2\4\2\5\2\4\2\5\2\4\2\5\2\4\2\eeee
\2\4\3\4\2\4\3\4\2\4\3\4\2\4\3\4\eeee
\5\2\4\2\5\2\4\2\5\2\4\2\5\2\4\2\eeee
} 

\baaa
5-9
\eaaa
\bbbb
0&0&0&0&4\\
0&0&0&2&2\\
0&0&2&2&0\\
0&2&2&0&0\\
2&2&0&0&0\\
\ebbb
\parbox{7cm}{ 
\1\5\2\4\3\3\4\2\5\1\5\2\4\3\3\4\eeee
\5\1\5\2\4\3\3\4\2\5\1\5\2\4\3\3\eeee
\2\5\1\5\2\4\3\3\4\2\5\1\5\2\4\3\eeee
\4\2\5\1\5\2\4\3\3\4\2\5\1\5\2\4\eeee
\3\4\2\5\1\5\2\4\3\3\4\2\5\1\5\2\eeee
\3\3\4\2\5\1\5\2\4\3\3\4\2\5\1\5\eeee
\4\3\3\4\2\5\1\5\2\4\3\3\4\2\5\1\eeee
\2\4\3\3\4\2\5\1\5\2\4\3\3\4\2\5\eeee
\5\2\4\3\3\4\2\5\1\5\2\4\3\3\4\2\eeee
\1\5\2\4\3\3\4\2\5\1\5\2\4\3\3\4\eeee
\5\1\5\2\4\3\3\4\2\5\1\5\2\4\3\3\eeee
\2\5\1\5\2\4\3\3\4\2\5\1\5\2\4\3\eeee
} 

\baaa
5-10
\eaaa
\bbbb
0&0&0&0&4\\
0&0&1&1&2\\
0&1&1&1&1\\
0&1&1&2&0\\
1&2&1&0&0\\
\ebbb
\parbox{7cm}{ 
\1\5\3\2\5\2\4\4\3\3\4\4\2\5\2\3\eeee
\5\2\4\4\3\3\4\4\2\5\2\3\5\1\5\3\eeee
\3\3\4\4\2\5\2\3\5\1\5\3\2\5\2\4\eeee
\2\5\2\3\5\1\5\3\2\5\2\4\4\3\3\4\eeee
\5\1\5\3\2\5\2\4\4\3\3\4\4\2\5\2\eeee
\2\5\2\4\4\3\3\4\4\2\5\2\3\5\1\5\eeee
\4\3\3\4\4\2\5\2\3\5\1\5\3\2\5\2\eeee
\4\2\5\2\3\5\1\5\3\2\5\2\4\4\3\3\eeee
\3\5\1\5\3\2\5\2\4\4\3\3\4\4\2\5\eeee
\3\2\5\2\4\4\3\3\4\4\2\5\2\3\5\1\eeee
\4\4\3\3\4\4\2\5\2\3\5\1\5\3\2\5\eeee
\4\4\2\5\2\3\5\1\5\3\2\5\2\4\4\3\eeee
} 

\baaa
5-11
\eaaa
\bbbb
0&0&0&0&4\\
0&0&1&1&2\\
0&2&1&1&0\\
0&2&1&1&0\\
2&2&0&0&0\\
\ebbb
\parbox{7cm}{ 
\1\5\2\3\3\2\5\1\5\2\3\3\2\5\1\5\eeee
\5\1\5\2\4\4\2\5\1\5\2\4\4\2\5\1\eeee
\2\5\1\5\2\3\3\2\5\1\5\2\3\3\2\5\eeee
\3\2\5\1\5\2\4\4\2\5\1\5\2\4\4\2\eeee
\3\4\2\5\1\5\2\3\3\2\5\1\5\2\3\3\eeee
\2\4\3\2\5\1\5\2\4\4\2\5\1\5\2\4\eeee
\5\2\3\4\2\5\1\5\2\3\3\2\5\1\5\2\eeee
\1\5\2\4\3\2\5\1\5\2\4\4\2\5\1\5\eeee
\5\1\5\2\3\4\2\5\1\5\2\3\3\2\5\1\eeee
\2\5\1\5\2\4\3\2\5\1\5\2\4\4\2\5\eeee
\3\2\5\1\5\2\3\4\2\5\1\5\2\3\3\2\eeee
\3\4\2\5\1\5\2\4\3\2\5\1\5\2\4\4\eeee
} 

\baaa
5-12
\eaaa
\bbbb
0&0&0&1&3\\
0&0&0&1&3\\
0&0&0&1&3\\
1&1&2&0&0\\
1&1&2&0&0\\
\ebbb
\parbox{7cm}{ 
\1\5\3\4\1\5\3\4\1\5\3\4\1\5\3\4\eeee
\4\2\5\3\5\2\5\3\5\2\5\3\5\2\5\3\eeee
\3\5\1\5\3\4\1\5\3\4\1\5\3\4\1\5\eeee
\5\3\4\2\5\3\5\2\5\3\5\2\5\3\5\2\eeee
\1\5\3\5\1\5\3\4\1\5\3\4\1\5\3\4\eeee
\4\2\5\3\4\2\5\3\5\2\5\3\5\2\5\3\eeee
\3\5\1\5\3\5\1\5\3\4\1\5\3\4\1\5\eeee
\5\3\4\2\5\3\4\2\5\3\5\2\5\3\5\2\eeee
\1\5\3\5\1\5\3\5\1\5\3\4\1\5\3\4\eeee
\4\2\5\3\4\2\5\3\4\2\5\3\5\2\5\3\eeee
\3\5\1\5\3\5\1\5\3\5\1\5\3\4\1\5\eeee
\5\3\4\2\5\3\4\2\5\3\4\2\5\3\5\2\eeee
} 

\baaa
5-13
\eaaa
\bbbb
0&0&0&1&3\\
0&0&1&2&1\\
0&1&2&1&0\\
1&2&1&0&0\\
3&1&0&0&0\\
\ebbb
\parbox{7cm}{ 
\1\5\1\5\1\5\1\5\1\5\1\5\1\5\1\5\eeee
\4\2\4\2\4\2\4\2\4\2\4\2\4\2\4\2\eeee
\3\3\3\3\3\3\3\3\3\3\3\3\3\3\3\3\eeee
\2\4\2\4\2\4\2\4\2\4\2\4\2\4\2\4\eeee
\5\1\5\1\5\1\5\1\5\1\5\1\5\1\5\1\eeee
\1\5\1\5\1\5\1\5\1\5\1\5\1\5\1\5\eeee
\4\2\4\2\4\2\4\2\4\2\4\2\4\2\4\2\eeee
\3\3\3\3\3\3\3\3\3\3\3\3\3\3\3\3\eeee
\2\4\2\4\2\4\2\4\2\4\2\4\2\4\2\4\eeee
\5\1\5\1\5\1\5\1\5\1\5\1\5\1\5\1\eeee
\1\5\1\5\1\5\1\5\1\5\1\5\1\5\1\5\eeee
\4\2\4\2\4\2\4\2\4\2\4\2\4\2\4\2\eeee
} 

\baaa
5-14
\eaaa
\bbbb
0&0&0&2&2\\
0&0&0&2&2\\
0&0&0&2&2\\
1&1&2&0&0\\
1&1&2&0&0\\
\ebbb
\parbox{7cm}{ 
\1\4\3\4\1\4\3\4\1\4\3\4\1\4\3\4\eeee
\4\2\5\3\5\2\5\3\5\2\5\3\5\2\5\3\eeee
\3\5\1\4\3\4\1\4\3\4\1\4\3\4\1\4\eeee
\4\3\4\2\5\3\5\2\5\3\5\2\5\3\5\2\eeee
\1\5\3\5\1\4\3\4\1\4\3\4\1\4\3\4\eeee
\4\2\4\3\4\2\5\3\5\2\5\3\5\2\5\3\eeee
\3\5\1\5\3\5\1\4\3\4\1\4\3\4\1\4\eeee
\4\3\4\2\4\3\4\2\5\3\5\2\5\3\5\2\eeee
\1\5\3\5\1\5\3\5\1\4\3\4\1\4\3\4\eeee
\4\2\4\3\4\2\4\3\4\2\5\3\5\2\5\3\eeee
\3\5\1\5\3\5\1\5\3\5\1\4\3\4\1\4\eeee
\4\3\4\2\4\3\4\2\4\3\4\2\5\3\5\2\eeee
} 

\baaa
5-15
\eaaa
\bbbb
0&0&0&2&2\\
0&0&0&2&2\\
0&0&2&0&2\\
1&1&0&2&0\\
1&1&2&0&0\\
\ebbb
\parbox{7cm}{ 
\1\4\4\1\5\3\3\5\1\4\4\1\5\3\3\5\eeee
\4\4\2\5\3\3\5\2\4\4\2\5\3\3\5\2\eeee
\4\1\5\3\3\5\1\4\4\1\5\3\3\5\1\4\eeee
\2\5\3\3\5\2\4\4\2\5\3\3\5\2\4\4\eeee
\5\3\3\5\1\4\4\1\5\3\3\5\1\4\4\1\eeee
\3\3\5\2\4\4\2\5\3\3\5\2\4\4\2\5\eeee
\3\5\1\4\4\1\5\3\3\5\1\4\4\1\5\3\eeee
\5\2\4\4\2\5\3\3\5\2\4\4\2\5\3\3\eeee
\1\4\4\1\5\3\3\5\1\4\4\1\5\3\3\5\eeee
\4\4\2\5\3\3\5\2\4\4\2\5\3\3\5\2\eeee
\4\1\5\3\3\5\1\4\4\1\5\3\3\5\1\4\eeee
\2\5\3\3\5\2\4\4\2\5\3\3\5\2\4\4\eeee
} 

\baaa
5-16
\eaaa
\bbbb
0&0&0&2&2\\
0&0&0&2&2\\
0&0&2&0&2\\
2&2&0&0&0\\
1&1&2&0&0\\
\ebbb
\parbox{7cm}{ 
\1\4\2\5\3\3\5\1\4\2\5\3\3\5\2\4\eeee
\4\1\5\3\3\5\2\4\1\5\3\3\5\1\4\2\eeee
\2\5\3\3\5\1\4\2\5\3\3\5\2\4\1\5\eeee
\5\3\3\5\2\4\1\5\3\3\5\1\4\2\5\3\eeee
\3\3\5\1\4\2\5\3\3\5\2\4\1\5\3\3\eeee
\3\5\2\4\1\5\3\3\5\1\4\2\5\3\3\5\eeee
\5\1\4\2\5\3\3\5\2\4\1\5\3\3\5\2\eeee
\2\4\1\5\3\3\5\1\4\2\5\3\3\5\1\4\eeee
\4\2\5\3\3\5\2\4\1\5\3\3\5\2\4\1\eeee
\1\5\3\3\5\1\4\2\5\3\3\5\1\4\2\5\eeee
\5\3\3\5\2\4\1\5\3\3\5\2\4\1\5\3\eeee
\3\3\5\1\4\2\5\3\3\5\1\4\2\5\3\3\eeee
} 

\baaa
5-17
\eaaa
\bbbb
0&0&0&2&2\\
0&0&0&2&2\\
0&0&2&1&1\\
1&1&2&0&0\\
1&1&2&0&0\\
\ebbb
\parbox{7cm}{ 
\1\4\3\3\4\1\4\3\3\4\1\4\3\3\4\1\eeee
\4\2\5\3\3\5\2\5\3\3\5\2\5\3\3\5\eeee
\3\5\1\4\3\3\4\1\4\3\3\4\1\4\3\3\eeee
\3\3\4\2\5\3\3\5\2\5\3\3\5\2\5\3\eeee
\4\3\3\5\1\4\3\3\4\1\4\3\3\4\1\4\eeee
\1\5\3\3\4\2\5\3\3\5\2\5\3\3\5\2\eeee
\4\2\4\3\3\5\1\4\3\3\4\1\4\3\3\4\eeee
\3\5\1\5\3\3\4\2\5\3\3\5\2\5\3\3\eeee
\3\3\4\2\4\3\3\5\1\4\3\3\4\1\4\3\eeee
\4\3\3\5\1\5\3\3\4\2\5\3\3\5\2\5\eeee
\1\5\3\3\4\2\4\3\3\5\1\4\3\3\4\1\eeee
\4\2\4\3\3\5\1\5\3\3\4\2\5\3\3\5\eeee
} 

\baaa
5-18
\eaaa
\bbbb
0&0&0&2&2\\
0&0&1&1&2\\
0&1&0&1&2\\
2&1&1&0&0\\
1&1&1&0&1\\
\ebbb
\parbox{7cm}{ 
\1\5\5\1\5\5\1\5\5\1\5\5\1\5\5\1\eeee
\4\2\3\4\2\3\4\2\3\4\2\3\4\2\3\4\eeee
\1\5\5\1\5\5\1\5\5\1\5\5\1\5\5\1\eeee
\4\3\2\4\3\2\4\3\2\4\3\2\4\3\2\4\eeee
\1\5\5\1\5\5\1\5\5\1\5\5\1\5\5\1\eeee
\4\2\3\4\2\3\4\2\3\4\2\3\4\2\3\4\eeee
\1\5\5\1\5\5\1\5\5\1\5\5\1\5\5\1\eeee
\4\3\2\4\3\2\4\3\2\4\3\2\4\3\2\4\eeee
\1\5\5\1\5\5\1\5\5\1\5\5\1\5\5\1\eeee
\4\2\3\4\2\3\4\2\3\4\2\3\4\2\3\4\eeee
\1\5\5\1\5\5\1\5\5\1\5\5\1\5\5\1\eeee
\4\3\2\4\3\2\4\3\2\4\3\2\4\3\2\4\eeee
} 

\baaa
5-19
\eaaa
\bbbb
0&0&0&2&2\\
0&0&1&1&2\\
0&1&1&1&1\\
1&1&1&1&0\\
1&2&1&0&0\\
\ebbb
\parbox{7cm}{ 
\1\5\2\4\3\3\4\2\5\1\5\2\4\3\3\4\eeee
\4\2\5\1\5\2\4\3\3\4\2\5\1\5\2\4\eeee
\4\3\3\4\2\5\1\5\2\4\3\3\4\2\5\1\eeee
\1\5\2\4\3\3\4\2\5\1\5\2\4\3\3\4\eeee
\4\2\5\1\5\2\4\3\3\4\2\5\1\5\2\4\eeee
\4\3\3\4\2\5\1\5\2\4\3\3\4\2\5\1\eeee
\1\5\2\4\3\3\4\2\5\1\5\2\4\3\3\4\eeee
\4\2\5\1\5\2\4\3\3\4\2\5\1\5\2\4\eeee
\4\3\3\4\2\5\1\5\2\4\3\3\4\2\5\1\eeee
\1\5\2\4\3\3\4\2\5\1\5\2\4\3\3\4\eeee
\4\2\5\1\5\2\4\3\3\4\2\5\1\5\2\4\eeee
\4\3\3\4\2\5\1\5\2\4\3\3\4\2\5\1\eeee
} 

\baaa
5-20
\eaaa
\bbbb
0&0&0&2&2\\
0&0&1&1&2\\
0&1&2&0&1\\
2&2&0&0&0\\
1&2&1&0&0\\
\ebbb
\parbox{7cm}{ 
\1\5\3\2\4\2\3\5\1\5\3\2\4\2\3\5\eeee
\4\2\3\5\1\5\3\2\4\2\3\5\1\5\3\2\eeee
\1\5\3\2\4\2\3\5\1\5\3\2\4\2\3\5\eeee
\4\2\3\5\1\5\3\2\4\2\3\5\1\5\3\2\eeee
\1\5\3\2\4\2\3\5\1\5\3\2\4\2\3\5\eeee
\4\2\3\5\1\5\3\2\4\2\3\5\1\5\3\2\eeee
\1\5\3\2\4\2\3\5\1\5\3\2\4\2\3\5\eeee
\4\2\3\5\1\5\3\2\4\2\3\5\1\5\3\2\eeee
\1\5\3\2\4\2\3\5\1\5\3\2\4\2\3\5\eeee
\4\2\3\5\1\5\3\2\4\2\3\5\1\5\3\2\eeee
\1\5\3\2\4\2\3\5\1\5\3\2\4\2\3\5\eeee
\4\2\3\5\1\5\3\2\4\2\3\5\1\5\3\2\eeee
} 

\baaa
5-21
\eaaa
\bbbb
0&0&0&2&2\\
0&0&2&0&2\\
0&1&1&1&1\\
1&0&1&1&1\\
1&1&1&1&0\\
\ebbb
\parbox{7cm}{ 
\1\5\2\5\1\5\2\5\1\5\2\5\1\5\2\5\eeee
\4\3\3\4\4\3\3\4\4\3\3\4\4\3\3\4\eeee
\5\2\5\1\5\2\5\1\5\2\5\1\5\2\5\1\eeee
\3\3\4\4\3\3\4\4\3\3\4\4\3\3\4\4\eeee
\2\5\1\5\2\5\1\5\2\5\1\5\2\5\1\5\eeee
\3\4\4\3\3\4\4\3\3\4\4\3\3\4\4\3\eeee
\5\1\5\2\5\1\5\2\5\1\5\2\5\1\5\2\eeee
\4\4\3\3\4\4\3\3\4\4\3\3\4\4\3\3\eeee
\1\5\2\5\1\5\2\5\1\5\2\5\1\5\2\5\eeee
\4\3\3\4\4\3\3\4\4\3\3\4\4\3\3\4\eeee
\5\2\5\1\5\2\5\1\5\2\5\1\5\2\5\1\eeee
\3\3\4\4\3\3\4\4\3\3\4\4\3\3\4\4\eeee
} 

\baaa
5-22
\eaaa
\bbbb
0&0&0&2&2\\
0&0&2&0&2\\
0&2&0&2&0\\
2&0&2&0&0\\
2&2&0&0&0\\
\ebbb
\parbox{7cm}{ 
\1\4\3\2\5\1\4\3\2\5\1\4\3\2\5\1\eeee
\4\3\2\5\1\4\3\2\5\1\4\3\2\5\1\4\eeee
\3\2\5\1\4\3\2\5\1\4\3\2\5\1\4\3\eeee
\2\5\1\4\3\2\5\1\4\3\2\5\1\4\3\2\eeee
\5\1\4\3\2\5\1\4\3\2\5\1\4\3\2\5\eeee
\1\4\3\2\5\1\4\3\2\5\1\4\3\2\5\1\eeee
\4\3\2\5\1\4\3\2\5\1\4\3\2\5\1\4\eeee
\3\2\5\1\4\3\2\5\1\4\3\2\5\1\4\3\eeee
\2\5\1\4\3\2\5\1\4\3\2\5\1\4\3\2\eeee
\5\1\4\3\2\5\1\4\3\2\5\1\4\3\2\5\eeee
\1\4\3\2\5\1\4\3\2\5\1\4\3\2\5\1\eeee
\4\3\2\5\1\4\3\2\5\1\4\3\2\5\1\4\eeee
} 

\baaa
5-23
\eaaa
\bbbb
0&0&0&2&2\\
0&0&2&0&2\\
0&2&2&0&0\\
2&0&0&2&0\\
2&2&0&0&0\\
\ebbb
\parbox{7cm}{ 
\1\4\4\1\5\2\3\3\2\5\1\4\4\1\5\2\eeee
\4\4\1\5\2\3\3\2\5\1\4\4\1\5\2\3\eeee
\4\1\5\2\3\3\2\5\1\4\4\1\5\2\3\3\eeee
\1\5\2\3\3\2\5\1\4\4\1\5\2\3\3\2\eeee
\5\2\3\3\2\5\1\4\4\1\5\2\3\3\2\5\eeee
\2\3\3\2\5\1\4\4\1\5\2\3\3\2\5\1\eeee
\3\3\2\5\1\4\4\1\5\2\3\3\2\5\1\4\eeee
\3\2\5\1\4\4\1\5\2\3\3\2\5\1\4\4\eeee
\2\5\1\4\4\1\5\2\3\3\2\5\1\4\4\1\eeee
\5\1\4\4\1\5\2\3\3\2\5\1\4\4\1\5\eeee
\1\4\4\1\5\2\3\3\2\5\1\4\4\1\5\2\eeee
\4\4\1\5\2\3\3\2\5\1\4\4\1\5\2\3\eeee
} 

\baaa
5-24
\eaaa
\bbbb
0&0&0&2&2\\
0&0&2&1&1\\
0&2&0&1&1\\
1&1&1&0&1\\
1&1&1&1&0\\
\ebbb
\parbox{7cm}{ 
\1\5\4\1\4\5\1\5\4\1\4\5\1\5\4\1\eeee
\4\2\3\5\2\3\4\2\3\5\2\3\4\2\3\5\eeee
\5\3\2\4\3\2\5\3\2\4\3\2\5\3\2\4\eeee
\1\4\5\1\5\4\1\4\5\1\5\4\1\4\5\1\eeee
\5\2\3\4\2\3\5\2\3\4\2\3\5\2\3\4\eeee
\4\3\2\5\3\2\4\3\2\5\3\2\4\3\2\5\eeee
\1\5\4\1\4\5\1\5\4\1\4\5\1\5\4\1\eeee
\4\2\3\5\2\3\4\2\3\5\2\3\4\2\3\5\eeee
\5\3\2\4\3\2\5\3\2\4\3\2\5\3\2\4\eeee
\1\4\5\1\5\4\1\4\5\1\5\4\1\4\5\1\eeee
\5\2\3\4\2\3\5\2\3\4\2\3\5\2\3\4\eeee
\4\3\2\5\3\2\4\3\2\5\3\2\4\3\2\5\eeee
} 

\baaa
5-25
\eaaa
\bbbb
0&0&0&2&2\\
0&0&2&1&1\\
0&2&0&1&1\\
1&1&1&1&0\\
1&1&1&0&1\\
\ebbb
\parbox{7cm}{ 
\1\5\5\1\5\5\1\5\5\1\5\5\1\5\5\1\eeee
\4\2\3\4\2\3\4\2\3\4\2\3\4\2\3\4\eeee
\4\3\2\4\3\2\4\3\2\4\3\2\4\3\2\4\eeee
\1\5\5\1\5\5\1\5\5\1\5\5\1\5\5\1\eeee
\4\2\3\4\2\3\4\2\3\4\2\3\4\2\3\4\eeee
\4\3\2\4\3\2\4\3\2\4\3\2\4\3\2\4\eeee
\1\5\5\1\5\5\1\5\5\1\5\5\1\5\5\1\eeee
\4\2\3\4\2\3\4\2\3\4\2\3\4\2\3\4\eeee
\4\3\2\4\3\2\4\3\2\4\3\2\4\3\2\4\eeee
\1\5\5\1\5\5\1\5\5\1\5\5\1\5\5\1\eeee
\4\2\3\4\2\3\4\2\3\4\2\3\4\2\3\4\eeee
\4\3\2\4\3\2\4\3\2\4\3\2\4\3\2\4\eeee
} 

\baaa
5-26
\eaaa
\bbbb
0&0&0&2&2\\
0&0&2&1&1\\
0&2&1&0&1\\
1&1&0&1&1\\
1&1&1&1&0\\
\ebbb
\parbox{7cm}{ 
\1\5\3\2\4\4\2\3\5\1\5\3\2\4\4\2\eeee
\4\2\3\5\1\5\3\2\4\4\2\3\5\1\5\3\eeee
\5\3\2\4\4\2\3\5\1\5\3\2\4\4\2\3\eeee
\2\3\5\1\5\3\2\4\4\2\3\5\1\5\3\2\eeee
\3\2\4\4\2\3\5\1\5\3\2\4\4\2\3\5\eeee
\3\5\1\5\3\2\4\4\2\3\5\1\5\3\2\4\eeee
\2\4\4\2\3\5\1\5\3\2\4\4\2\3\5\1\eeee
\5\1\5\3\2\4\4\2\3\5\1\5\3\2\4\4\eeee
\4\4\2\3\5\1\5\3\2\4\4\2\3\5\1\5\eeee
\1\5\3\2\4\4\2\3\5\1\5\3\2\4\4\2\eeee
\4\2\3\5\1\5\3\2\4\4\2\3\5\1\5\3\eeee
\5\3\2\4\4\2\3\5\1\5\3\2\4\4\2\3\eeee
} 

\baaa
5-27
\eaaa
\bbbb
0&0&0&2&2\\
0&1&0&1&2\\
0&0&1&1&2\\
2&1&1&0&0\\
1&1&1&0&1\\
\ebbb
\parbox{7cm}{ 
\1\5\5\1\5\5\1\5\5\1\5\5\1\5\5\1\eeee
\4\2\2\4\3\3\4\2\2\4\3\3\4\2\2\4\eeee
\1\5\5\1\5\5\1\5\5\1\5\5\1\5\5\1\eeee
\4\3\3\4\2\2\4\3\3\4\2\2\4\3\3\4\eeee
\1\5\5\1\5\5\1\5\5\1\5\5\1\5\5\1\eeee
\4\2\2\4\3\3\4\2\2\4\3\3\4\2\2\4\eeee
\1\5\5\1\5\5\1\5\5\1\5\5\1\5\5\1\eeee
\4\3\3\4\2\2\4\3\3\4\2\2\4\3\3\4\eeee
\1\5\5\1\5\5\1\5\5\1\5\5\1\5\5\1\eeee
\4\2\2\4\3\3\4\2\2\4\3\3\4\2\2\4\eeee
\1\5\5\1\5\5\1\5\5\1\5\5\1\5\5\1\eeee
\4\3\3\4\2\2\4\3\3\4\2\2\4\3\3\4\eeee
} 

\baaa
5-28
\eaaa
\bbbb
0&0&0&2&2\\
0&1&1&0&2\\
0&1&1&0&2\\
2&0&0&2&0\\
2&1&1&0&0\\
\ebbb
\parbox{7cm}{ 
\1\4\4\1\5\3\3\5\1\4\4\1\5\3\3\5\eeee
\4\4\1\5\2\2\5\1\4\4\1\5\2\2\5\1\eeee
\4\1\5\3\3\5\1\4\4\1\5\3\3\5\1\4\eeee
\1\5\2\2\5\1\4\4\1\5\2\2\5\1\4\4\eeee
\5\3\3\5\1\4\4\1\5\3\3\5\1\4\4\1\eeee
\2\2\5\1\4\4\1\5\2\2\5\1\4\4\1\5\eeee
\3\5\1\4\4\1\5\3\3\5\1\4\4\1\5\3\eeee
\5\1\4\4\1\5\2\2\5\1\4\4\1\5\2\2\eeee
\1\4\4\1\5\3\3\5\1\4\4\1\5\3\3\5\eeee
\4\4\1\5\2\2\5\1\4\4\1\5\2\2\5\1\eeee
\4\1\5\3\3\5\1\4\4\1\5\3\3\5\1\4\eeee
\1\5\2\2\5\1\4\4\1\5\2\2\5\1\4\4\eeee
} 

\baaa
5-29
\eaaa
\bbbb
0&0&0&2&2\\
0&1&1&1&1\\
0&1&1&1&1\\
1&1&1&0&1\\
1&1&1&1&0\\
\ebbb
\parbox{7cm}{ 
\1\5\4\1\4\5\1\5\4\1\4\5\1\5\4\1\eeee
\4\2\3\5\2\3\4\2\3\5\2\3\4\2\3\5\eeee
\5\2\3\4\2\3\5\2\3\4\2\3\5\2\3\4\eeee
\1\4\5\1\5\4\1\4\5\1\5\4\1\4\5\1\eeee
\5\3\2\4\3\2\5\3\2\4\3\2\5\3\2\4\eeee
\4\3\2\5\3\2\4\3\2\5\3\2\4\3\2\5\eeee
\1\5\4\1\4\5\1\5\4\1\4\5\1\5\4\1\eeee
\4\2\3\5\2\3\4\2\3\5\2\3\4\2\3\5\eeee
\5\2\3\4\2\3\5\2\3\4\2\3\5\2\3\4\eeee
\1\4\5\1\5\4\1\4\5\1\5\4\1\4\5\1\eeee
\5\3\2\4\3\2\5\3\2\4\3\2\5\3\2\4\eeee
\4\3\2\5\3\2\4\3\2\5\3\2\4\3\2\5\eeee
} 

\baaa
5-30
\eaaa
\bbbb
0&0&0&2&2\\
0&1&1&1&1\\
0&1&1&1&1\\
1&1&1&1&0\\
1&1&1&0&1\\
\ebbb
\parbox{7cm}{ 
\1\5\5\1\5\5\1\5\5\1\5\5\1\5\5\1\eeee
\4\2\3\4\2\3\4\2\3\4\2\3\4\2\3\4\eeee
\4\2\3\4\2\3\4\2\3\4\2\3\4\2\3\4\eeee
\1\5\5\1\5\5\1\5\5\1\5\5\1\5\5\1\eeee
\4\3\2\4\3\2\4\3\2\4\3\2\4\3\2\4\eeee
\4\3\2\4\3\2\4\3\2\4\3\2\4\3\2\4\eeee
\1\5\5\1\5\5\1\5\5\1\5\5\1\5\5\1\eeee
\4\2\3\4\2\3\4\2\3\4\2\3\4\2\3\4\eeee
\4\2\3\4\2\3\4\2\3\4\2\3\4\2\3\4\eeee
\1\5\5\1\5\5\1\5\5\1\5\5\1\5\5\1\eeee
\4\3\2\4\3\2\4\3\2\4\3\2\4\3\2\4\eeee
\4\3\2\4\3\2\4\3\2\4\3\2\4\3\2\4\eeee
} 

\baaa
5-31
\eaaa
\bbbb
0&0&0&2&2\\
0&1&1&1&1\\
0&1&1&1&1\\
2&1&1&0&0\\
2&1&1&0&0\\
\ebbb
\parbox{7cm}{ 
\1\4\2\2\4\1\4\2\2\4\1\4\2\2\4\1\eeee
\4\1\5\3\3\5\1\5\3\3\5\1\5\3\3\5\eeee
\2\5\1\4\2\2\4\1\4\2\2\4\1\4\2\2\eeee
\2\3\4\1\5\3\3\5\1\5\3\3\5\1\5\3\eeee
\4\3\2\5\1\4\2\2\4\1\4\2\2\4\1\4\eeee
\1\5\2\3\4\1\5\3\3\5\1\5\3\3\5\1\eeee
\4\1\4\3\2\5\1\4\2\2\4\1\4\2\2\4\eeee
\2\5\1\5\2\3\4\1\5\3\3\5\1\5\3\3\eeee
\2\3\4\1\4\3\2\5\1\4\2\2\4\1\4\2\eeee
\4\3\2\5\1\5\2\3\4\1\5\3\3\5\1\5\eeee
\1\5\2\3\4\1\4\3\2\5\1\4\2\2\4\1\eeee
\4\1\4\3\2\5\1\5\2\3\4\1\5\3\3\5\eeee
} 

\baaa
5-32
\eaaa
\bbbb
0&0&0&2&2\\
0&2&0&1&1\\
0&0&2&1&1\\
1&1&1&0&1\\
1&1&1&1&0\\
\ebbb
\parbox{7cm}{ 
\1\5\4\1\4\5\1\5\4\1\4\5\1\5\4\1\eeee
\4\2\2\5\3\3\4\2\2\5\3\3\4\2\2\5\eeee
\5\2\2\4\3\3\5\2\2\4\3\3\5\2\2\4\eeee
\1\4\5\1\5\4\1\4\5\1\5\4\1\4\5\1\eeee
\5\3\3\4\2\2\5\3\3\4\2\2\5\3\3\4\eeee
\4\3\3\5\2\2\4\3\3\5\2\2\4\3\3\5\eeee
\1\5\4\1\4\5\1\5\4\1\4\5\1\5\4\1\eeee
\4\2\2\5\3\3\4\2\2\5\3\3\4\2\2\5\eeee
\5\2\2\4\3\3\5\2\2\4\3\3\5\2\2\4\eeee
\1\4\5\1\5\4\1\4\5\1\5\4\1\4\5\1\eeee
\5\3\3\4\2\2\5\3\3\4\2\2\5\3\3\4\eeee
\4\3\3\5\2\2\4\3\3\5\2\2\4\3\3\5\eeee
} 

\baaa
5-33
\eaaa
\bbbb
0&0&0&2&2\\
0&2&0&1&1\\
0&0&2&1&1\\
1&1&1&1&0\\
1&1&1&0&1\\
\ebbb
\parbox{7cm}{ 
\1\5\5\1\5\5\1\5\5\1\5\5\1\5\5\1\eeee
\4\2\2\4\3\3\4\2\2\4\3\3\4\2\2\4\eeee
\4\2\2\4\3\3\4\2\2\4\3\3\4\2\2\4\eeee
\1\5\5\1\5\5\1\5\5\1\5\5\1\5\5\1\eeee
\4\3\3\4\2\2\4\3\3\4\2\2\4\3\3\4\eeee
\4\3\3\4\2\2\4\3\3\4\2\2\4\3\3\4\eeee
\1\5\5\1\5\5\1\5\5\1\5\5\1\5\5\1\eeee
\4\2\2\4\3\3\4\2\2\4\3\3\4\2\2\4\eeee
\4\2\2\4\3\3\4\2\2\4\3\3\4\2\2\4\eeee
\1\5\5\1\5\5\1\5\5\1\5\5\1\5\5\1\eeee
\4\3\3\4\2\2\4\3\3\4\2\2\4\3\3\4\eeee
\4\3\3\4\2\2\4\3\3\4\2\2\4\3\3\4\eeee
} 

\baaa
5-34
\eaaa
\bbbb
0&0&1&1&2\\
0&0&1&1&2\\
1&1&0&0&2\\
1&1&0&0&2\\
1&1&1&1&0\\
\ebbb
\parbox{7cm}{ 
\1\4\5\1\3\5\2\3\5\2\4\5\1\4\5\1\eeee
\3\5\2\4\5\1\4\5\1\3\5\2\3\5\2\4\eeee
\5\1\3\5\2\3\5\2\4\5\1\4\5\1\3\5\eeee
\2\4\5\1\4\5\1\3\5\2\3\5\2\4\5\1\eeee
\3\5\2\3\5\2\4\5\1\4\5\1\3\5\2\3\eeee
\5\1\4\5\1\3\5\2\3\5\2\4\5\1\4\5\eeee
\2\3\5\2\4\5\1\4\5\1\3\5\2\3\5\2\eeee
\4\5\1\3\5\2\3\5\2\4\5\1\4\5\1\3\eeee
\5\2\4\5\1\4\5\1\3\5\2\3\5\2\4\5\eeee
\1\3\5\2\3\5\2\4\5\1\4\5\1\3\5\2\eeee
\4\5\1\4\5\1\3\5\2\3\5\2\4\5\1\4\eeee
\5\2\3\5\2\4\5\1\4\5\1\3\5\2\3\5\eeee
} 

\baaa
5-35
\eaaa
\bbbb
0&0&1&1&2\\
0&0&1&1&2\\
1&1&1&1&0\\
1&1&1&1&0\\
1&1&0&0&2\\
\ebbb
\parbox{7cm}{ 
\1\4\4\1\5\5\2\3\3\2\5\5\1\4\4\1\eeee
\3\3\2\5\5\1\4\4\1\5\5\2\3\3\2\5\eeee
\4\1\5\5\2\3\3\2\5\5\1\4\4\1\5\5\eeee
\2\5\5\1\4\4\1\5\5\2\3\3\2\5\5\1\eeee
\5\5\2\3\3\2\5\5\1\4\4\1\5\5\2\3\eeee
\5\1\4\4\1\5\5\2\3\3\2\5\5\1\4\4\eeee
\2\3\3\2\5\5\1\4\4\1\5\5\2\3\3\2\eeee
\4\4\1\5\5\2\3\3\2\5\5\1\4\4\1\5\eeee
\3\2\5\5\1\4\4\1\5\5\2\3\3\2\5\5\eeee
\1\5\5\2\3\3\2\5\5\1\4\4\1\5\5\2\eeee
\5\5\1\4\4\1\5\5\2\3\3\2\5\5\1\4\eeee
\5\2\3\3\2\5\5\1\4\4\1\5\5\2\3\3\eeee
} 

\baaa
5-36
\eaaa
\bbbb
0&0&1&1&2\\
0&1&0&1&2\\
1&0&0&1&2\\
1&1&1&1&0\\
1&1&1&0&1\\
\ebbb
\parbox{7cm}{ 
\1\4\2\5\3\5\1\4\2\5\3\5\1\4\2\5\eeee
\3\4\2\5\1\5\3\4\2\5\1\5\3\4\2\5\eeee
\5\1\5\3\4\2\5\1\5\3\4\2\5\1\5\3\eeee
\5\3\5\1\4\2\5\3\5\1\4\2\5\3\5\1\eeee
\1\4\2\5\3\5\1\4\2\5\3\5\1\4\2\5\eeee
\3\4\2\5\1\5\3\4\2\5\1\5\3\4\2\5\eeee
\5\1\5\3\4\2\5\1\5\3\4\2\5\1\5\3\eeee
\5\3\5\1\4\2\5\3\5\1\4\2\5\3\5\1\eeee
\1\4\2\5\3\5\1\4\2\5\3\5\1\4\2\5\eeee
\3\4\2\5\1\5\3\4\2\5\1\5\3\4\2\5\eeee
\5\1\5\3\4\2\5\1\5\3\4\2\5\1\5\3\eeee
\5\3\5\1\4\2\5\3\5\1\4\2\5\3\5\1\eeee
} 
\baaa
\phantom{0-0}\#
\eaaa
\mbox{}\phantom{\bbbb
0&0&0&0&0\\
\ebbb}
\parbox{7cm}{ 
\1\4\3\5\2\5\1\4\3\5\2\5\1\4\3\5\eeee
\3\4\1\5\2\5\3\4\1\5\2\5\3\4\1\5\eeee
\5\2\5\3\4\1\5\2\5\3\4\1\5\2\5\3\eeee
\5\2\5\1\4\3\5\2\5\1\4\3\5\2\5\1\eeee
\1\4\3\5\2\5\1\4\3\5\2\5\1\4\3\5\eeee
\3\4\1\5\2\5\3\4\1\5\2\5\3\4\1\5\eeee
\5\2\5\3\4\1\5\2\5\3\4\1\5\2\5\3\eeee
\5\2\5\1\4\3\5\2\5\1\4\3\5\2\5\1\eeee
\1\4\3\5\2\5\1\4\3\5\2\5\1\4\3\5\eeee
\3\4\1\5\2\5\3\4\1\5\2\5\3\4\1\5\eeee
\5\2\5\3\4\1\5\2\5\3\4\1\5\2\5\3\eeee
\5\2\5\1\4\3\5\2\5\1\4\3\5\2\5\1\eeee
} 

\baaa
5-37
\eaaa
\bbbb
0&0&1&1&2\\
0&1&0&1&2\\
1&0&1&0&2\\
1&1&0&0&2\\
1&1&1&1&0\\
\ebbb
\parbox{7cm}{ 
\1\4\5\3\1\5\3\3\5\1\3\5\4\1\5\2\eeee
\3\5\1\3\5\4\1\5\2\4\5\2\2\5\4\2\eeee
\5\2\4\5\2\2\5\4\2\5\1\4\5\3\1\5\eeee
\4\2\5\1\4\5\3\1\5\3\3\5\1\3\5\4\eeee
\1\5\3\3\5\1\3\5\4\1\5\2\4\5\2\2\eeee
\5\4\1\5\2\4\5\2\2\5\4\2\5\1\4\5\eeee
\2\2\5\4\2\5\1\4\5\3\1\5\3\3\5\1\eeee
\4\5\3\1\5\3\3\5\1\3\5\4\1\5\2\4\eeee
\5\1\3\5\4\1\5\2\4\5\2\2\5\4\2\5\eeee
\2\4\5\2\2\5\4\2\5\1\4\5\3\1\5\3\eeee
\2\5\1\4\5\3\1\5\3\3\5\1\3\5\4\1\eeee
\5\3\3\5\1\3\5\4\1\5\2\4\5\2\2\5\eeee
} 

\baaa
5-38
\eaaa
\bbbb
0&0&1&1&2\\
0&1&0&2&1\\
1&0&2&0&1\\
1&2&0&1&0\\
2&1&1&0&0\\
\ebbb
\parbox{7cm}{ 
\1\5\1\5\1\5\1\5\1\5\1\5\1\5\1\5\eeee
\3\3\3\3\3\3\3\3\3\3\3\3\3\3\3\3\eeee
\5\1\5\1\5\1\5\1\5\1\5\1\5\1\5\1\eeee
\2\4\2\4\2\4\2\4\2\4\2\4\2\4\2\4\eeee
\2\4\2\4\2\4\2\4\2\4\2\4\2\4\2\4\eeee
\5\1\5\1\5\1\5\1\5\1\5\1\5\1\5\1\eeee
\3\3\3\3\3\3\3\3\3\3\3\3\3\3\3\3\eeee
\1\5\1\5\1\5\1\5\1\5\1\5\1\5\1\5\eeee
\4\2\4\2\4\2\4\2\4\2\4\2\4\2\4\2\eeee
\4\2\4\2\4\2\4\2\4\2\4\2\4\2\4\2\eeee
\1\5\1\5\1\5\1\5\1\5\1\5\1\5\1\5\eeee
\3\3\3\3\3\3\3\3\3\3\3\3\3\3\3\3\eeee
} 

\baaa
5-39
\eaaa
\bbbb
0&0&1&1&2\\
0&1&1&2&0\\
1&1&0&1&1\\
1&2&1&0&0\\
2&0&1&0&1\\
\ebbb
\parbox{7cm}{ 
\1\4\2\3\5\1\4\2\3\5\1\4\2\3\5\1\eeee
\3\2\4\1\5\3\2\4\1\5\3\2\4\1\5\3\eeee
\4\2\3\5\1\4\2\3\5\1\4\2\3\5\1\4\eeee
\2\4\1\5\3\2\4\1\5\3\2\4\1\5\3\2\eeee
\2\3\5\1\4\2\3\5\1\4\2\3\5\1\4\2\eeee
\4\1\5\3\2\4\1\5\3\2\4\1\5\3\2\4\eeee
\3\5\1\4\2\3\5\1\4\2\3\5\1\4\2\3\eeee
\1\5\3\2\4\1\5\3\2\4\1\5\3\2\4\1\eeee
\5\1\4\2\3\5\1\4\2\3\5\1\4\2\3\5\eeee
\5\3\2\4\1\5\3\2\4\1\5\3\2\4\1\5\eeee
\1\4\2\3\5\1\4\2\3\5\1\4\2\3\5\1\eeee
\3\2\4\1\5\3\2\4\1\5\3\2\4\1\5\3\eeee
} 

\baaa
5-40
\eaaa
\bbbb
0&0&1&1&2\\
0&2&0&1&1\\
1&0&2&1&0\\
1&1&1&0&1\\
2&1&0&1&0\\
\ebbb
\parbox{7cm}{ 
\1\4\5\1\3\3\4\2\2\5\1\4\5\1\3\3\eeee
\3\3\1\5\4\1\5\2\2\4\3\3\1\5\4\1\eeee
\3\3\4\2\2\5\1\4\5\1\3\3\4\2\2\5\eeee
\4\1\5\2\2\4\3\3\1\5\4\1\5\2\2\4\eeee
\2\5\1\4\5\1\3\3\4\2\2\5\1\4\5\1\eeee
\2\4\3\3\1\5\4\1\5\2\2\4\3\3\1\5\eeee
\5\1\3\3\4\2\2\5\1\4\5\1\3\3\4\2\eeee
\1\5\4\1\5\2\2\4\3\3\1\5\4\1\5\2\eeee
\4\2\2\5\1\4\5\1\3\3\4\2\2\5\1\4\eeee
\5\2\2\4\3\3\1\5\4\1\5\2\2\4\3\3\eeee
\1\4\5\1\3\3\4\2\2\5\1\4\5\1\3\3\eeee
\3\3\1\5\4\1\5\2\2\4\3\3\1\5\4\1\eeee
} 

\baaa
5-41
\eaaa
\bbbb
0&1&1&1&1\\
1&0&1&1&1\\
1&1&0&1&1\\
1&1&1&0&1\\
1&1&1&1&0\\
\ebbb
\parbox{7cm}{ 
\1\3\5\2\4\1\3\5\2\4\1\3\5\2\4\1\eeee
\2\4\1\3\5\2\4\1\3\5\2\4\1\3\5\2\eeee
\3\5\2\4\1\3\5\2\4\1\3\5\2\4\1\3\eeee
\4\1\3\5\2\4\1\3\5\2\4\1\3\5\2\4\eeee
\5\2\4\1\3\5\2\4\1\3\5\2\4\1\3\5\eeee
\1\3\5\2\4\1\3\5\2\4\1\3\5\2\4\1\eeee
\2\4\1\3\5\2\4\1\3\5\2\4\1\3\5\2\eeee
\3\5\2\4\1\3\5\2\4\1\3\5\2\4\1\3\eeee
\4\1\3\5\2\4\1\3\5\2\4\1\3\5\2\4\eeee
\5\2\4\1\3\5\2\4\1\3\5\2\4\1\3\5\eeee
\1\3\5\2\4\1\3\5\2\4\1\3\5\2\4\1\eeee
\2\4\1\3\5\2\4\1\3\5\2\4\1\3\5\2\eeee
} 

\baaa
5-42
\eaaa
\bbbb
0&1&1&1&1\\
1&0&1&1&1\\
1&1&0&1&1\\
1&1&1&1&0\\
1&1&1&0&1\\
\ebbb
\parbox{7cm}{ 
\1\3\5\2\4\1\3\4\2\5\1\3\5\2\4\1\eeee
\2\4\1\3\4\2\5\1\3\5\2\4\1\3\4\2\eeee
\3\4\2\5\1\3\5\2\4\1\3\4\2\5\1\3\eeee
\5\1\3\5\2\4\1\3\4\2\5\1\3\5\2\4\eeee
\5\2\4\1\3\4\2\5\1\3\5\2\4\1\3\4\eeee
\1\3\4\2\5\1\3\5\2\4\1\3\4\2\5\1\eeee
\2\5\1\3\5\2\4\1\3\4\2\5\1\3\5\2\eeee
\3\5\2\4\1\3\4\2\5\1\3\5\2\4\1\3\eeee
\4\1\3\4\2\5\1\3\5\2\4\1\3\4\2\5\eeee
\4\2\5\1\3\5\2\4\1\3\4\2\5\1\3\5\eeee
\1\3\5\2\4\1\3\4\2\5\1\3\5\2\4\1\eeee
\2\4\1\3\4\2\5\1\3\5\2\4\1\3\4\2\eeee
} 

\baaa
5-43
\eaaa
\bbbb
0&1&1&1&1\\
1&1&0&1&1\\
1&0&1&1&1\\
1&1&1&1&0\\
1&1&1&0&1\\
\ebbb
\parbox{7cm}{ 
\1\3\5\2\4\1\2\4\3\5\1\3\5\2\4\1\eeee
\2\4\1\2\4\3\5\1\3\5\2\4\1\2\4\3\eeee
\2\4\3\5\1\3\5\2\4\1\2\4\3\5\1\3\eeee
\5\1\3\5\2\4\1\2\4\3\5\1\3\5\2\4\eeee
\5\2\4\1\2\4\3\5\1\3\5\2\4\1\2\4\eeee
\1\2\4\3\5\1\3\5\2\4\1\2\4\3\5\1\eeee
\3\5\1\3\5\2\4\1\2\4\3\5\1\3\5\2\eeee
\3\5\2\4\1\2\4\3\5\1\3\5\2\4\1\2\eeee
\4\1\2\4\3\5\1\3\5\2\4\1\2\4\3\5\eeee
\4\3\5\1\3\5\2\4\1\2\4\3\5\1\3\5\eeee
\1\3\5\2\4\1\2\4\3\5\1\3\5\2\4\1\eeee
\2\4\1\2\4\3\5\1\3\5\2\4\1\2\4\3\eeee
} 
\baaa
\phantom{0-0}\#
\eaaa
\mbox{}\phantom{\bbbb
0&0&0&0&0\\
\ebbb}
\parbox{7cm}{ 
\1\4\3\3\4\1\5\2\2\5\1\4\3\3\4\1\eeee
\2\2\5\1\4\3\3\4\1\5\2\2\5\1\4\3\eeee
\4\1\5\2\2\5\1\4\3\3\4\1\5\2\2\5\eeee
\4\3\3\4\1\5\2\2\5\1\4\3\3\4\1\5\eeee
\2\5\1\4\3\3\4\1\5\2\2\5\1\4\3\3\eeee
\1\5\2\2\5\1\4\3\3\4\1\5\2\2\5\1\eeee
\3\3\4\1\5\2\2\5\1\4\3\3\4\1\5\2\eeee
\5\1\4\3\3\4\1\5\2\2\5\1\4\3\3\4\eeee
\5\2\2\5\1\4\3\3\4\1\5\2\2\5\1\4\eeee
\3\4\1\5\2\2\5\1\4\3\3\4\1\5\2\2\eeee
\1\4\3\3\4\1\5\2\2\5\1\4\3\3\4\1\eeee
\2\2\5\1\4\3\3\4\1\5\2\2\5\1\4\3\eeee
} 

\baaa
5-44
\eaaa
\bbbb
0&1&1&1&1\\
2&1&0&0&1\\
2&0&1&1&0\\
2&0&1&1&0\\
2&1&0&0&1\\
\ebbb
\parbox{7cm}{ 
\1\3\3\1\2\2\1\3\3\1\2\2\1\3\3\1\eeee
\2\1\4\4\1\5\5\1\4\4\1\5\5\1\4\4\eeee
\2\5\1\3\3\1\2\2\1\3\3\1\2\2\1\3\eeee
\1\5\2\1\4\4\1\5\5\1\4\4\1\5\5\1\eeee
\3\1\2\5\1\3\3\1\2\2\1\3\3\1\2\2\eeee
\3\4\1\5\2\1\4\4\1\5\5\1\4\4\1\5\eeee
\1\4\3\1\2\5\1\3\3\1\2\2\1\3\3\1\eeee
\2\1\3\4\1\5\2\1\4\4\1\5\5\1\4\4\eeee
\2\5\1\4\3\1\2\5\1\3\3\1\2\2\1\3\eeee
\1\5\2\1\3\4\1\5\2\1\4\4\1\5\5\1\eeee
\3\1\2\5\1\4\3\1\2\5\1\3\3\1\2\2\eeee
\3\4\1\5\2\1\3\4\1\5\2\1\4\4\1\5\eeee
} 

\baaa
5-45
\eaaa
\bbbb
1&0&0&1&2\\
0&1&0&1&2\\
0&0&1&1&2\\
1&1&1&1&0\\
1&1&1&0&1\\
\ebbb
\parbox{7cm}{ 
\1\4\3\5\2\5\1\4\3\5\2\5\1\4\3\5\eeee
\1\4\3\5\2\5\1\4\3\5\2\5\1\4\3\5\eeee
\5\2\5\1\4\3\5\2\5\1\4\3\5\2\5\1\eeee
\5\2\5\1\4\3\5\2\5\1\4\3\5\2\5\1\eeee
\1\4\3\5\2\5\1\4\3\5\2\5\1\4\3\5\eeee
\1\4\3\5\2\5\1\4\3\5\2\5\1\4\3\5\eeee
\5\2\5\1\4\3\5\2\5\1\4\3\5\2\5\1\eeee
\5\2\5\1\4\3\5\2\5\1\4\3\5\2\5\1\eeee
\1\4\3\5\2\5\1\4\3\5\2\5\1\4\3\5\eeee
\1\4\3\5\2\5\1\4\3\5\2\5\1\4\3\5\eeee
\5\2\5\1\4\3\5\2\5\1\4\3\5\2\5\1\eeee
\5\2\5\1\4\3\5\2\5\1\4\3\5\2\5\1\eeee
} 

\baaa
5-46
\eaaa
\bbbb
1&0&0&1&2\\
0&1&1&1&1\\
0&1&2&0&1\\
1&2&0&1&0\\
1&1&1&0&1\\
\ebbb
\parbox{7cm}{ 
\1\5\3\2\4\2\3\5\1\5\3\2\4\2\3\5\eeee
\1\5\3\2\4\2\3\5\1\5\3\2\4\2\3\5\eeee
\4\2\3\5\1\5\3\2\4\2\3\5\1\5\3\2\eeee
\4\2\3\5\1\5\3\2\4\2\3\5\1\5\3\2\eeee
\1\5\3\2\4\2\3\5\1\5\3\2\4\2\3\5\eeee
\1\5\3\2\4\2\3\5\1\5\3\2\4\2\3\5\eeee
\4\2\3\5\1\5\3\2\4\2\3\5\1\5\3\2\eeee
\4\2\3\5\1\5\3\2\4\2\3\5\1\5\3\2\eeee
\1\5\3\2\4\2\3\5\1\5\3\2\4\2\3\5\eeee
\1\5\3\2\4\2\3\5\1\5\3\2\4\2\3\5\eeee
\4\2\3\5\1\5\3\2\4\2\3\5\1\5\3\2\eeee
\4\2\3\5\1\5\3\2\4\2\3\5\1\5\3\2\eeee
} 

\baaa
5-47
\eaaa
\bbbb
1&0&0&1&2\\
0&1&1&1&1\\
0&1&2&1&0\\
1&1&1&1&0\\
2&1&0&0&1\\
\ebbb
\parbox{7cm}{ 
\1\4\3\2\5\1\4\3\2\5\1\4\3\2\5\1\eeee
\1\4\3\2\5\1\4\3\2\5\1\4\3\2\5\1\eeee
\5\2\3\4\1\5\2\3\4\1\5\2\3\4\1\5\eeee
\5\2\3\4\1\5\2\3\4\1\5\2\3\4\1\5\eeee
\1\4\3\2\5\1\4\3\2\5\1\4\3\2\5\1\eeee
\1\4\3\2\5\1\4\3\2\5\1\4\3\2\5\1\eeee
\5\2\3\4\1\5\2\3\4\1\5\2\3\4\1\5\eeee
\5\2\3\4\1\5\2\3\4\1\5\2\3\4\1\5\eeee
\1\4\3\2\5\1\4\3\2\5\1\4\3\2\5\1\eeee
\1\4\3\2\5\1\4\3\2\5\1\4\3\2\5\1\eeee
\5\2\3\4\1\5\2\3\4\1\5\2\3\4\1\5\eeee
\5\2\3\4\1\5\2\3\4\1\5\2\3\4\1\5\eeee
} 

\baaa
5-48
\eaaa
\bbbb
1&0&1&1&1\\
0&2&0&1&1\\
1&0&2&0&1\\
1&1&0&2&0\\
1&1&1&0&1\\
\ebbb
\parbox{7cm}{ 
\1\3\5\2\2\5\3\1\4\4\1\3\5\2\2\5\eeee
\1\3\5\2\2\5\3\1\4\4\1\3\5\2\2\5\eeee
\5\3\1\4\4\1\3\5\2\2\5\3\1\4\4\1\eeee
\5\3\1\4\4\1\3\5\2\2\5\3\1\4\4\1\eeee
\1\3\5\2\2\5\3\1\4\4\1\3\5\2\2\5\eeee
\1\3\5\2\2\5\3\1\4\4\1\3\5\2\2\5\eeee
\5\3\1\4\4\1\3\5\2\2\5\3\1\4\4\1\eeee
\5\3\1\4\4\1\3\5\2\2\5\3\1\4\4\1\eeee
\1\3\5\2\2\5\3\1\4\4\1\3\5\2\2\5\eeee
\1\3\5\2\2\5\3\1\4\4\1\3\5\2\2\5\eeee
\5\3\1\4\4\1\3\5\2\2\5\3\1\4\4\1\eeee
\5\3\1\4\4\1\3\5\2\2\5\3\1\4\4\1\eeee
} 

\baaa
5-49
\eaaa
\bbbb
2&0&0&0&2\\
0&2&0&1&1\\
0&0&2&2&0\\
0&1&1&2&0\\
1&1&0&0&2\\
\ebbb
\parbox{7cm}{ 
\1\5\2\4\3\4\2\5\1\5\2\4\3\4\2\5\eeee
\1\5\2\4\3\4\2\5\1\5\2\4\3\4\2\5\eeee
\1\5\2\4\3\4\2\5\1\5\2\4\3\4\2\5\eeee
\1\5\2\4\3\4\2\5\1\5\2\4\3\4\2\5\eeee
\1\5\2\4\3\4\2\5\1\5\2\4\3\4\2\5\eeee
\1\5\2\4\3\4\2\5\1\5\2\4\3\4\2\5\eeee
\1\5\2\4\3\4\2\5\1\5\2\4\3\4\2\5\eeee
\1\5\2\4\3\4\2\5\1\5\2\4\3\4\2\5\eeee
\1\5\2\4\3\4\2\5\1\5\2\4\3\4\2\5\eeee
\1\5\2\4\3\4\2\5\1\5\2\4\3\4\2\5\eeee
\1\5\2\4\3\4\2\5\1\5\2\4\3\4\2\5\eeee
\1\5\2\4\3\4\2\5\1\5\2\4\3\4\2\5\eeee
} 

\baaa
5-50
\eaaa
\bbbb
2&0&0&0&2\\
0&2&0&1&1\\
0&0&3&1&0\\
0&1&1&2&0\\
1&1&0&0&2\\
\ebbb
\parbox{7cm}{ 
\1\5\2\4\3\3\4\2\5\1\5\2\4\3\3\4\eeee
\1\5\2\4\3\3\4\2\5\1\5\2\4\3\3\4\eeee
\1\5\2\4\3\3\4\2\5\1\5\2\4\3\3\4\eeee
\1\5\2\4\3\3\4\2\5\1\5\2\4\3\3\4\eeee
\1\5\2\4\3\3\4\2\5\1\5\2\4\3\3\4\eeee
\1\5\2\4\3\3\4\2\5\1\5\2\4\3\3\4\eeee
\1\5\2\4\3\3\4\2\5\1\5\2\4\3\3\4\eeee
\1\5\2\4\3\3\4\2\5\1\5\2\4\3\3\4\eeee
\1\5\2\4\3\3\4\2\5\1\5\2\4\3\3\4\eeee
\1\5\2\4\3\3\4\2\5\1\5\2\4\3\3\4\eeee
\1\5\2\4\3\3\4\2\5\1\5\2\4\3\3\4\eeee
\1\5\2\4\3\3\4\2\5\1\5\2\4\3\3\4\eeee
} 

\baaa
5-51
\eaaa
\bbbb
2&0&0&1&1\\
0&2&1&0&1\\
0&1&2&1&0\\
1&0&1&2&0\\
1&1&0&0&2\\
\ebbb
\parbox{7cm}{ 
\1\4\3\2\5\1\4\3\2\5\1\4\3\2\5\1\eeee
\1\4\3\2\5\1\4\3\2\5\1\4\3\2\5\1\eeee
\1\4\3\2\5\1\4\3\2\5\1\4\3\2\5\1\eeee
\1\4\3\2\5\1\4\3\2\5\1\4\3\2\5\1\eeee
\1\4\3\2\5\1\4\3\2\5\1\4\3\2\5\1\eeee
\1\4\3\2\5\1\4\3\2\5\1\4\3\2\5\1\eeee
\1\4\3\2\5\1\4\3\2\5\1\4\3\2\5\1\eeee
\1\4\3\2\5\1\4\3\2\5\1\4\3\2\5\1\eeee
\1\4\3\2\5\1\4\3\2\5\1\4\3\2\5\1\eeee
\1\4\3\2\5\1\4\3\2\5\1\4\3\2\5\1\eeee
\1\4\3\2\5\1\4\3\2\5\1\4\3\2\5\1\eeee
\1\4\3\2\5\1\4\3\2\5\1\4\3\2\5\1\eeee
} 

\baaa
5-52
\eaaa
\bbbb
2&0&0&1&1\\
0&2&1&0&1\\
0&1&3&0&0\\
1&0&0&3&0\\
1&1&0&0&2\\
\ebbb
\parbox{7cm}{ 
\1\4\4\1\5\2\3\3\2\5\1\4\4\1\5\2\eeee
\1\4\4\1\5\2\3\3\2\5\1\4\4\1\5\2\eeee
\1\4\4\1\5\2\3\3\2\5\1\4\4\1\5\2\eeee
\1\4\4\1\5\2\3\3\2\5\1\4\4\1\5\2\eeee
\1\4\4\1\5\2\3\3\2\5\1\4\4\1\5\2\eeee
\1\4\4\1\5\2\3\3\2\5\1\4\4\1\5\2\eeee
\1\4\4\1\5\2\3\3\2\5\1\4\4\1\5\2\eeee
\1\4\4\1\5\2\3\3\2\5\1\4\4\1\5\2\eeee
\1\4\4\1\5\2\3\3\2\5\1\4\4\1\5\2\eeee
\1\4\4\1\5\2\3\3\2\5\1\4\4\1\5\2\eeee
\1\4\4\1\5\2\3\3\2\5\1\4\4\1\5\2\eeee
\1\4\4\1\5\2\3\3\2\5\1\4\4\1\5\2\eeee
} 


\baaa
6-1
\eaaa
\bbbb
0&0&0&0&0&4\\
0&0&0&0&0&4\\
0&0&0&0&2&2\\
0&0&0&0&2&2\\
0&0&1&1&2&0\\
1&1&1&1&0&0\\
\ebbb
\parbox{7cm}{ 
\1\6\3\5\5\3\6\1\6\3\5\5\3\6\1\6\eeee
\6\2\6\4\5\5\4\6\2\6\4\5\5\4\6\2\eeee
\3\6\1\6\3\5\5\3\6\1\6\3\5\5\3\6\eeee
\5\4\6\2\6\4\5\5\4\6\2\6\4\5\5\4\eeee
\5\5\3\6\1\6\3\5\5\3\6\1\6\3\5\5\eeee
\3\5\5\4\6\2\6\4\5\5\4\6\2\6\4\5\eeee
\6\4\5\5\3\6\1\6\3\5\5\3\6\1\6\3\eeee
\1\6\3\5\5\4\6\2\6\4\5\5\4\6\2\6\eeee
\6\2\6\4\5\5\3\6\1\6\3\5\5\3\6\1\eeee
\3\6\1\6\3\5\5\4\6\2\6\4\5\5\4\6\eeee
\5\4\6\2\6\4\5\5\3\6\1\6\3\5\5\3\eeee
\5\5\3\6\1\6\3\5\5\4\6\2\6\4\5\5\eeee
} 

\baaa
6-2
\eaaa
\bbbb
0&0&0&0&0&4\\
0&0&0&0&0&4\\
0&0&0&0&2&2\\
0&0&0&0&2&2\\
0&0&2&2&0&0\\
1&1&1&1&0&0\\
\ebbb
\parbox{7cm}{ 
\1\6\3\5\3\6\1\6\3\5\3\6\1\6\3\5\eeee
\6\2\6\4\5\4\6\2\6\4\5\4\6\2\6\4\eeee
\3\6\1\6\3\5\3\6\1\6\3\5\3\6\1\6\eeee
\5\4\6\2\6\4\5\4\6\2\6\4\5\4\6\2\eeee
\3\5\3\6\1\6\3\5\3\6\1\6\3\5\3\6\eeee
\6\4\5\4\6\2\6\4\5\4\6\2\6\4\5\4\eeee
\1\6\3\5\3\6\1\6\3\5\3\6\1\6\3\5\eeee
\6\2\6\4\5\4\6\2\6\4\5\4\6\2\6\4\eeee
\3\6\1\6\3\5\3\6\1\6\3\5\3\6\1\6\eeee
\5\4\6\2\6\4\5\4\6\2\6\4\5\4\6\2\eeee
\3\5\3\6\1\6\3\5\3\6\1\6\3\5\3\6\eeee
\6\4\5\4\6\2\6\4\5\4\6\2\6\4\5\4\eeee
} 

\baaa
6-3
\eaaa
\bbbb
0&0&0&0&0&4\\
0&0&0&0&0&4\\
0&0&0&0&2&2\\
0&0&0&0&4&0\\
0&0&2&2&0&0\\
1&1&2&0&0&0\\
\ebbb
\parbox{7cm}{ 
\1\6\3\5\4\5\3\6\1\6\3\5\4\5\3\6\eeee
\6\2\6\3\5\4\5\3\6\2\6\3\5\4\5\3\eeee
\3\6\1\6\3\5\4\5\3\6\1\6\3\5\4\5\eeee
\5\3\6\2\6\3\5\4\5\3\6\2\6\3\5\4\eeee
\4\5\3\6\1\6\3\5\4\5\3\6\1\6\3\5\eeee
\5\4\5\3\6\2\6\3\5\4\5\3\6\2\6\3\eeee
\3\5\4\5\3\6\1\6\3\5\4\5\3\6\1\6\eeee
\6\3\5\4\5\3\6\2\6\3\5\4\5\3\6\2\eeee
\1\6\3\5\4\5\3\6\1\6\3\5\4\5\3\6\eeee
\6\2\6\3\5\4\5\3\6\2\6\3\5\4\5\3\eeee
\3\6\1\6\3\5\4\5\3\6\1\6\3\5\4\5\eeee
\5\3\6\2\6\3\5\4\5\3\6\2\6\3\5\4\eeee
} 

\baaa
6-4
\eaaa
\bbbb
0&0&0&0&0&4\\
0&0&0&0&0&4\\
0&0&0&0&2&2\\
0&0&0&2&2&0\\
0&0&2&2&0&0\\
1&1&2&0&0&0\\
\ebbb
\parbox{7cm}{ 
\1\6\3\5\4\4\5\3\6\1\6\3\5\4\4\5\eeee
\6\2\6\3\5\4\4\5\3\6\2\6\3\5\4\4\eeee
\3\6\1\6\3\5\4\4\5\3\6\1\6\3\5\4\eeee
\5\3\6\2\6\3\5\4\4\5\3\6\2\6\3\5\eeee
\4\5\3\6\1\6\3\5\4\4\5\3\6\1\6\3\eeee
\4\4\5\3\6\2\6\3\5\4\4\5\3\6\2\6\eeee
\5\4\4\5\3\6\1\6\3\5\4\4\5\3\6\1\eeee
\3\5\4\4\5\3\6\2\6\3\5\4\4\5\3\6\eeee
\6\3\5\4\4\5\3\6\1\6\3\5\4\4\5\3\eeee
\1\6\3\5\4\4\5\3\6\2\6\3\5\4\4\5\eeee
\6\2\6\3\5\4\4\5\3\6\1\6\3\5\4\4\eeee
\3\6\1\6\3\5\4\4\5\3\6\2\6\3\5\4\eeee
} 

\baaa
6-5
\eaaa
\bbbb
0&0&0&0&0&4\\
0&0&0&0&0&4\\
0&0&0&0&4&0\\
0&0&0&0&4&0\\
0&0&1&1&0&2\\
1&1&0&0&2&0\\
\ebbb
\parbox{7cm}{ 
\1\6\5\3\5\6\1\6\5\3\5\6\1\6\5\3\eeee
\6\2\6\5\4\5\6\2\6\5\4\5\6\2\6\5\eeee
\5\6\1\6\5\3\5\6\1\6\5\3\5\6\1\6\eeee
\3\5\6\2\6\5\4\5\6\2\6\5\4\5\6\2\eeee
\5\4\5\6\1\6\5\3\5\6\1\6\5\3\5\6\eeee
\6\5\3\5\6\2\6\5\4\5\6\2\6\5\4\5\eeee
\1\6\5\4\5\6\1\6\5\3\5\6\1\6\5\3\eeee
\6\2\6\5\3\5\6\2\6\5\4\5\6\2\6\5\eeee
\5\6\1\6\5\4\5\6\1\6\5\3\5\6\1\6\eeee
\3\5\6\2\6\5\3\5\6\2\6\5\4\5\6\2\eeee
\5\4\5\6\1\6\5\4\5\6\1\6\5\3\5\6\eeee
\6\5\3\5\6\2\6\5\3\5\6\2\6\5\4\5\eeee
} 

\baaa
6-6
\eaaa
\bbbb
0&0&0&0&0&4\\
0&0&0&0&0&4\\
0&0&0&1&1&2\\
0&0&2&1&1&0\\
0&0&2&1&1&0\\
1&1&2&0&0&0\\
\ebbb
\parbox{7cm}{ 
\1\6\3\4\4\3\6\1\6\3\4\4\3\6\1\6\eeee
\6\2\6\3\5\5\3\6\2\6\3\5\5\3\6\2\eeee
\3\6\1\6\3\4\4\3\6\1\6\3\4\4\3\6\eeee
\4\3\6\2\6\3\5\5\3\6\2\6\3\5\5\3\eeee
\4\5\3\6\1\6\3\4\4\3\6\1\6\3\4\4\eeee
\3\5\4\3\6\2\6\3\5\5\3\6\2\6\3\5\eeee
\6\3\4\5\3\6\1\6\3\4\4\3\6\1\6\3\eeee
\1\6\3\5\4\3\6\2\6\3\5\5\3\6\2\6\eeee
\6\2\6\3\4\5\3\6\1\6\3\4\4\3\6\1\eeee
\3\6\1\6\3\5\4\3\6\2\6\3\5\5\3\6\eeee
\4\3\6\2\6\3\4\5\3\6\1\6\3\4\4\3\eeee
\4\5\3\6\1\6\3\5\4\3\6\2\6\3\5\5\eeee
} 

\baaa
6-7
\eaaa
\bbbb
0&0&0&0&0&4\\
0&0&0&0&2&2\\
0&0&0&0&2&2\\
0&0&0&0&4&0\\
0&1&1&2&0&0\\
2&1&1&0&0&0\\
\ebbb
\parbox{7cm}{ 
\1\6\2\5\4\5\2\6\1\6\2\5\4\5\2\6\eeee
\6\1\6\3\5\4\5\3\6\1\6\3\5\4\5\3\eeee
\2\6\1\6\2\5\4\5\2\6\1\6\2\5\4\5\eeee
\5\3\6\1\6\3\5\4\5\3\6\1\6\3\5\4\eeee
\4\5\2\6\1\6\2\5\4\5\2\6\1\6\2\5\eeee
\5\4\5\3\6\1\6\3\5\4\5\3\6\1\6\3\eeee
\2\5\4\5\2\6\1\6\2\5\4\5\2\6\1\6\eeee
\6\3\5\4\5\3\6\1\6\3\5\4\5\3\6\1\eeee
\1\6\2\5\4\5\2\6\1\6\2\5\4\5\2\6\eeee
\6\1\6\3\5\4\5\3\6\1\6\3\5\4\5\3\eeee
\2\6\1\6\2\5\4\5\2\6\1\6\2\5\4\5\eeee
\5\3\6\1\6\3\5\4\5\3\6\1\6\3\5\4\eeee
} 

\baaa
6-8
\eaaa
\bbbb
0&0&0&0&0&4\\
0&0&0&0&2&2\\
0&0&0&0&2&2\\
0&0&0&0&4&0\\
0&1&2&1&0&0\\
1&1&2&0&0&0\\
\ebbb
\parbox{7cm}{ 
\1\6\3\6\1\6\3\6\1\6\3\6\1\6\3\6\eeee
\6\2\5\3\6\2\5\3\6\2\5\3\6\2\5\3\eeee
\3\5\4\5\3\5\4\5\3\5\4\5\3\5\4\5\eeee
\6\3\5\2\6\3\5\2\6\3\5\2\6\3\5\2\eeee
\1\6\3\6\1\6\3\6\1\6\3\6\1\6\3\6\eeee
\6\2\5\3\6\2\5\3\6\2\5\3\6\2\5\3\eeee
\3\5\4\5\3\5\4\5\3\5\4\5\3\5\4\5\eeee
\6\3\5\2\6\3\5\2\6\3\5\2\6\3\5\2\eeee
\1\6\3\6\1\6\3\6\1\6\3\6\1\6\3\6\eeee
\6\2\5\3\6\2\5\3\6\2\5\3\6\2\5\3\eeee
\3\5\4\5\3\5\4\5\3\5\4\5\3\5\4\5\eeee
\6\3\5\2\6\3\5\2\6\3\5\2\6\3\5\2\eeee
} 

\baaa
6-9
\eaaa
\bbbb
0&0&0&0&0&4\\
0&0&0&0&2&2\\
0&0&0&0&2&2\\
0&0&0&2&2&0\\
0&1&1&2&0&0\\
2&1&1&0&0&0\\
\ebbb
\parbox{7cm}{ 
\1\6\2\5\4\4\5\2\6\1\6\2\5\4\4\5\eeee
\6\1\6\3\5\4\4\5\3\6\1\6\3\5\4\4\eeee
\2\6\1\6\2\5\4\4\5\2\6\1\6\2\5\4\eeee
\5\3\6\1\6\3\5\4\4\5\3\6\1\6\3\5\eeee
\4\5\2\6\1\6\2\5\4\4\5\2\6\1\6\2\eeee
\4\4\5\3\6\1\6\3\5\4\4\5\3\6\1\6\eeee
\5\4\4\5\2\6\1\6\2\5\4\4\5\2\6\1\eeee
\2\5\4\4\5\3\6\1\6\3\5\4\4\5\3\6\eeee
\6\3\5\4\4\5\2\6\1\6\2\5\4\4\5\2\eeee
\1\6\2\5\4\4\5\3\6\1\6\3\5\4\4\5\eeee
\6\1\6\3\5\4\4\5\2\6\1\6\2\5\4\4\eeee
\2\6\1\6\2\5\4\4\5\3\6\1\6\3\5\4\eeee
} 

\baaa
6-10
\eaaa
\bbbb
0&0&0&0&0&4\\
0&0&0&0&2&2\\
0&0&0&1&2&1\\
0&0&4&0&0&0\\
0&2&2&0&0&0\\
1&2&1&0&0&0\\
\ebbb
\parbox{7cm}{ 
\1\6\3\4\3\6\1\6\3\4\3\6\1\6\3\4\eeee
\6\2\5\3\5\2\6\2\5\3\5\2\6\2\5\3\eeee
\3\5\2\6\2\5\3\5\2\6\2\5\3\5\2\6\eeee
\4\3\6\1\6\3\4\3\6\1\6\3\4\3\6\1\eeee
\3\5\2\6\2\5\3\5\2\6\2\5\3\5\2\6\eeee
\6\2\5\3\5\2\6\2\5\3\5\2\6\2\5\3\eeee
\1\6\3\4\3\6\1\6\3\4\3\6\1\6\3\4\eeee
\6\2\5\3\5\2\6\2\5\3\5\2\6\2\5\3\eeee
\3\5\2\6\2\5\3\5\2\6\2\5\3\5\2\6\eeee
\4\3\6\1\6\3\4\3\6\1\6\3\4\3\6\1\eeee
\3\5\2\6\2\5\3\5\2\6\2\5\3\5\2\6\eeee
\6\2\5\3\5\2\6\2\5\3\5\2\6\2\5\3\eeee
} 

\baaa
6-11
\eaaa
\bbbb
0&0&0&0&0&4\\
0&0&0&0&2&2\\
0&0&0&2&0&2\\
0&0&1&1&1&1\\
0&2&0&2&0&0\\
1&1&1&1&0&0\\
\ebbb
\parbox{7cm}{ 
\1\6\4\3\6\2\5\4\4\5\2\6\3\4\6\1\eeee
\6\2\5\4\4\5\2\6\3\4\6\1\6\4\3\6\eeee
\4\5\2\6\3\4\6\1\6\4\3\6\2\5\4\4\eeee
\3\4\6\1\6\4\3\6\2\5\4\4\5\2\6\3\eeee
\6\4\3\6\2\5\4\4\5\2\6\3\4\6\1\6\eeee
\2\5\4\4\5\2\6\3\4\6\1\6\4\3\6\2\eeee
\5\2\6\3\4\6\1\6\4\3\6\2\5\4\4\5\eeee
\4\6\1\6\4\3\6\2\5\4\4\5\2\6\3\4\eeee
\4\3\6\2\5\4\4\5\2\6\3\4\6\1\6\4\eeee
\5\4\4\5\2\6\3\4\6\1\6\4\3\6\2\5\eeee
\2\6\3\4\6\1\6\4\3\6\2\5\4\4\5\2\eeee
\6\1\6\4\3\6\2\5\4\4\5\2\6\3\4\6\eeee
} 

\baaa
6-12
\eaaa
\bbbb
0&0&0&0&0&4\\
0&0&0&0&2&2\\
0&0&0&2&2&0\\
0&0&2&2&0&0\\
0&2&2&0&0&0\\
2&2&0&0&0&0\\
\ebbb
\parbox{7cm}{ 
\1\6\2\5\3\4\4\3\5\2\6\1\6\2\5\3\eeee
\6\1\6\2\5\3\4\4\3\5\2\6\1\6\2\5\eeee
\2\6\1\6\2\5\3\4\4\3\5\2\6\1\6\2\eeee
\5\2\6\1\6\2\5\3\4\4\3\5\2\6\1\6\eeee
\3\5\2\6\1\6\2\5\3\4\4\3\5\2\6\1\eeee
\4\3\5\2\6\1\6\2\5\3\4\4\3\5\2\6\eeee
\4\4\3\5\2\6\1\6\2\5\3\4\4\3\5\2\eeee
\3\4\4\3\5\2\6\1\6\2\5\3\4\4\3\5\eeee
\5\3\4\4\3\5\2\6\1\6\2\5\3\4\4\3\eeee
\2\5\3\4\4\3\5\2\6\1\6\2\5\3\4\4\eeee
\6\2\5\3\4\4\3\5\2\6\1\6\2\5\3\4\eeee
\1\6\2\5\3\4\4\3\5\2\6\1\6\2\5\3\eeee
} 

\baaa
6-13
\eaaa
\bbbb
0&0&0&0&0&4\\
0&0&0&0&2&2\\
0&0&0&2&2&0\\
0&0&4&0&0&0\\
0&2&2&0&0&0\\
2&2&0&0&0&0\\
\ebbb
\parbox{7cm}{ 
\1\6\2\5\3\4\3\5\2\6\1\6\2\5\3\4\eeee
\6\1\6\2\5\3\4\3\5\2\6\1\6\2\5\3\eeee
\2\6\1\6\2\5\3\4\3\5\2\6\1\6\2\5\eeee
\5\2\6\1\6\2\5\3\4\3\5\2\6\1\6\2\eeee
\3\5\2\6\1\6\2\5\3\4\3\5\2\6\1\6\eeee
\4\3\5\2\6\1\6\2\5\3\4\3\5\2\6\1\eeee
\3\4\3\5\2\6\1\6\2\5\3\4\3\5\2\6\eeee
\5\3\4\3\5\2\6\1\6\2\5\3\4\3\5\2\eeee
\2\5\3\4\3\5\2\6\1\6\2\5\3\4\3\5\eeee
\6\2\5\3\4\3\5\2\6\1\6\2\5\3\4\3\eeee
\1\6\2\5\3\4\3\5\2\6\1\6\2\5\3\4\eeee
\6\1\6\2\5\3\4\3\5\2\6\1\6\2\5\3\eeee
} 

\baaa
6-14
\eaaa
\bbbb
0&0&0&0&0&4\\
0&0&0&0&2&2\\
0&0&1&0&2&1\\
0&0&0&2&2&0\\
0&1&1&1&1&0\\
1&2&1&0&0&0\\
\ebbb
\parbox{7cm}{ 
\1\6\3\3\6\1\6\3\3\6\1\6\3\3\6\1\eeee
\6\2\5\5\2\6\2\5\5\2\6\2\5\5\2\6\eeee
\3\5\4\4\5\3\5\4\4\5\3\5\4\4\5\3\eeee
\3\5\4\4\5\3\5\4\4\5\3\5\4\4\5\3\eeee
\6\2\5\5\2\6\2\5\5\2\6\2\5\5\2\6\eeee
\1\6\3\3\6\1\6\3\3\6\1\6\3\3\6\1\eeee
\6\2\5\5\2\6\2\5\5\2\6\2\5\5\2\6\eeee
\3\5\4\4\5\3\5\4\4\5\3\5\4\4\5\3\eeee
\3\5\4\4\5\3\5\4\4\5\3\5\4\4\5\3\eeee
\6\2\5\5\2\6\2\5\5\2\6\2\5\5\2\6\eeee
\1\6\3\3\6\1\6\3\3\6\1\6\3\3\6\1\eeee
\6\2\5\5\2\6\2\5\5\2\6\2\5\5\2\6\eeee
} 

\baaa
6-15
\eaaa
\bbbb
0&0&0&0&0&4\\
0&0&0&0&2&2\\
0&0&1&1&2&0\\
0&0&1&1&2&0\\
0&2&1&1&0&0\\
2&2&0&0&0&0\\
\ebbb
\parbox{7cm}{ 
\1\6\2\5\3\3\5\2\6\1\6\2\5\3\3\5\eeee
\6\1\6\2\5\4\4\5\2\6\1\6\2\5\4\4\eeee
\2\6\1\6\2\5\3\3\5\2\6\1\6\2\5\3\eeee
\5\2\6\1\6\2\5\4\4\5\2\6\1\6\2\5\eeee
\3\5\2\6\1\6\2\5\3\3\5\2\6\1\6\2\eeee
\3\4\5\2\6\1\6\2\5\4\4\5\2\6\1\6\eeee
\5\4\3\5\2\6\1\6\2\5\3\3\5\2\6\1\eeee
\2\5\3\4\5\2\6\1\6\2\5\4\4\5\2\6\eeee
\6\2\5\4\3\5\2\6\1\6\2\5\3\3\5\2\eeee
\1\6\2\5\3\4\5\2\6\1\6\2\5\4\4\5\eeee
\6\1\6\2\5\4\3\5\2\6\1\6\2\5\3\3\eeee
\2\6\1\6\2\5\3\4\5\2\6\1\6\2\5\4\eeee
} 

\baaa
6-16
\eaaa
\bbbb
0&0&0&0&0&4\\
0&0&0&1&1&2\\
0&0&0&1&1&2\\
0&1&1&1&1&0\\
0&1&1&1&1&0\\
2&1&1&0&0&0\\
\ebbb
\parbox{7cm}{ 
\1\6\2\4\4\2\6\1\6\2\4\4\2\6\1\6\eeee
\6\1\6\3\5\5\3\6\1\6\3\5\5\3\6\1\eeee
\2\6\1\6\2\4\4\2\6\1\6\2\4\4\2\6\eeee
\4\3\6\1\6\3\5\5\3\6\1\6\3\5\5\3\eeee
\4\5\2\6\1\6\2\4\4\2\6\1\6\2\4\4\eeee
\2\5\4\3\6\1\6\3\5\5\3\6\1\6\3\5\eeee
\6\3\4\5\2\6\1\6\2\4\4\2\6\1\6\2\eeee
\1\6\2\5\4\3\6\1\6\3\5\5\3\6\1\6\eeee
\6\1\6\3\4\5\2\6\1\6\2\4\4\2\6\1\eeee
\2\6\1\6\2\5\4\3\6\1\6\3\5\5\3\6\eeee
\4\3\6\1\6\3\4\5\2\6\1\6\2\4\4\2\eeee
\4\5\2\6\1\6\2\5\4\3\6\1\6\3\5\5\eeee
} 

\baaa
6-17
\eaaa
\bbbb
0&0&0&0&0&4\\
0&0&0&1&1&2\\
0&0&0&1&1&2\\
0&2&2&0&0&0\\
0&2&2&0&0&0\\
2&1&1&0&0&0\\
\ebbb
\parbox{7cm}{ 
\1\6\2\4\2\6\1\6\2\4\2\6\1\6\2\4\eeee
\6\1\6\3\5\3\6\1\6\3\5\3\6\1\6\3\eeee
\2\6\1\6\2\4\2\6\1\6\2\4\2\6\1\6\eeee
\4\3\6\1\6\3\5\3\6\1\6\3\5\3\6\1\eeee
\2\5\2\6\1\6\2\4\2\6\1\6\2\4\2\6\eeee
\6\3\4\3\6\1\6\3\5\3\6\1\6\3\5\3\eeee
\1\6\2\5\2\6\1\6\2\4\2\6\1\6\2\4\eeee
\6\1\6\3\4\3\6\1\6\3\5\3\6\1\6\3\eeee
\2\6\1\6\2\5\2\6\1\6\2\4\2\6\1\6\eeee
\4\3\6\1\6\3\4\3\6\1\6\3\5\3\6\1\eeee
\2\5\2\6\1\6\2\5\2\6\1\6\2\4\2\6\eeee
\6\3\4\3\6\1\6\3\4\3\6\1\6\3\5\3\eeee
} 

\baaa
6-18
\eaaa
\bbbb
0&0&0&0&0&4\\
0&0&0&1&1&2\\
0&0&0&2&2&0\\
0&2&2&0&0&0\\
0&2&2&0&0&0\\
2&2&0&0&0&0\\
\ebbb
\parbox{7cm}{ 
\1\6\2\4\3\4\2\6\1\6\2\4\3\4\2\6\eeee
\6\1\6\2\5\3\5\2\6\1\6\2\5\3\5\2\eeee
\2\6\1\6\2\4\3\4\2\6\1\6\2\4\3\4\eeee
\4\2\6\1\6\2\5\3\5\2\6\1\6\2\5\3\eeee
\3\5\2\6\1\6\2\4\3\4\2\6\1\6\2\4\eeee
\4\3\4\2\6\1\6\2\5\3\5\2\6\1\6\2\eeee
\2\5\3\5\2\6\1\6\2\4\3\4\2\6\1\6\eeee
\6\2\4\3\4\2\6\1\6\2\5\3\5\2\6\1\eeee
\1\6\2\5\3\5\2\6\1\6\2\4\3\4\2\6\eeee
\6\1\6\2\4\3\4\2\6\1\6\2\5\3\5\2\eeee
\2\6\1\6\2\5\3\5\2\6\1\6\2\4\3\4\eeee
\4\2\6\1\6\2\4\3\4\2\6\1\6\2\5\3\eeee
} 

\baaa
6-19
\eaaa
\bbbb
0&0&0&0&0&4\\
0&0&0&1&1&2\\
0&0&0&2&2&0\\
0&3&1&0&0&0\\
0&3&1&0&0&0\\
1&3&0&0&0&0\\
\ebbb
\parbox{7cm}{ 
\1\6\2\6\1\6\2\6\1\6\2\6\1\6\2\6\eeee
\6\2\5\2\6\2\4\2\6\2\5\2\6\2\4\2\eeee
\2\4\3\4\2\5\3\5\2\4\3\4\2\5\3\5\eeee
\6\2\5\2\6\2\4\2\6\2\5\2\6\2\4\2\eeee
\1\6\2\6\1\6\2\6\1\6\2\6\1\6\2\6\eeee
\6\2\4\2\6\2\5\2\6\2\4\2\6\2\5\2\eeee
\2\5\3\5\2\4\3\4\2\5\3\5\2\4\3\4\eeee
\6\2\4\2\6\2\5\2\6\2\4\2\6\2\5\2\eeee
\1\6\2\6\1\6\2\6\1\6\2\6\1\6\2\6\eeee
\6\2\5\2\6\2\4\2\6\2\5\2\6\2\4\2\eeee
\2\4\3\4\2\5\3\5\2\4\3\4\2\5\3\5\eeee
\6\2\5\2\6\2\4\2\6\2\5\2\6\2\4\2\eeee
} 

\baaa
6-20
\eaaa
\bbbb
0&0&0&0&0&4\\
0&0&0&1&1&2\\
0&0&2&1&1&0\\
0&2&2&0&0&0\\
0&2&2&0&0&0\\
2&2&0&0&0&0\\
\ebbb
\parbox{7cm}{ 
\1\6\2\4\3\3\4\2\6\1\6\2\4\3\3\4\eeee
\6\1\6\2\5\3\3\5\2\6\1\6\2\5\3\3\eeee
\2\6\1\6\2\4\3\3\4\2\6\1\6\2\4\3\eeee
\4\2\6\1\6\2\5\3\3\5\2\6\1\6\2\5\eeee
\3\5\2\6\1\6\2\4\3\3\4\2\6\1\6\2\eeee
\3\3\4\2\6\1\6\2\5\3\3\5\2\6\1\6\eeee
\4\3\3\5\2\6\1\6\2\4\3\3\4\2\6\1\eeee
\2\5\3\3\4\2\6\1\6\2\5\3\3\5\2\6\eeee
\6\2\4\3\3\5\2\6\1\6\2\4\3\3\4\2\eeee
\1\6\2\5\3\3\4\2\6\1\6\2\5\3\3\5\eeee
\6\1\6\2\4\3\3\5\2\6\1\6\2\4\3\3\eeee
\2\6\1\6\2\5\3\3\4\2\6\1\6\2\5\3\eeee
} 

\baaa
6-21
\eaaa
\bbbb
0&0&0&0&1&3\\
0&0&0&0&1&3\\
0&0&0&0&1&3\\
0&0&0&0&1&3\\
1&1&1&1&0&0\\
1&1&1&1&0&0\\
\ebbb
\parbox{7cm}{ 
\1\6\3\5\1\6\3\5\1\6\3\5\1\6\3\5\eeee
\5\2\6\4\6\2\6\4\6\2\6\4\6\2\6\4\eeee
\3\6\1\6\3\5\1\6\3\5\1\6\3\5\1\6\eeee
\6\4\5\2\6\4\6\2\6\4\6\2\6\4\6\2\eeee
\1\6\3\6\1\6\3\5\1\6\3\5\1\6\3\5\eeee
\5\2\6\4\5\2\6\4\6\2\6\4\6\2\6\4\eeee
\3\6\1\6\3\6\1\6\3\5\1\6\3\5\1\6\eeee
\6\4\5\2\6\4\5\2\6\4\6\2\6\4\6\2\eeee
\1\6\3\6\1\6\3\6\1\6\3\5\1\6\3\5\eeee
\5\2\6\4\5\2\6\4\5\2\6\4\6\2\6\4\eeee
\3\6\1\6\3\6\1\6\3\6\1\6\3\5\1\6\eeee
\6\4\5\2\6\4\5\2\6\4\5\2\6\4\6\2\eeee
} 

\baaa
6-22
\eaaa
\bbbb
0&0&0&0&1&3\\
0&0&0&1&0&3\\
0&0&0&1&1&2\\
0&1&3&0&0&0\\
1&0&3&0&0&0\\
1&1&2&0&0&0\\
\ebbb
\parbox{7cm}{ 
\1\6\2\6\3\4\3\5\3\6\1\6\2\6\3\4\eeee
\5\3\6\1\6\2\6\3\4\3\5\3\6\1\6\2\eeee
\3\4\3\5\3\6\1\6\2\6\3\4\3\5\3\6\eeee
\6\2\6\3\4\3\5\3\6\1\6\2\6\3\4\3\eeee
\3\6\1\6\2\6\3\4\3\5\3\6\1\6\2\6\eeee
\4\3\5\3\6\1\6\2\6\3\4\3\5\3\6\1\eeee
\2\6\3\4\3\5\3\6\1\6\2\6\3\4\3\5\eeee
\6\1\6\2\6\3\4\3\5\3\6\1\6\2\6\3\eeee
\3\5\3\6\1\6\2\6\3\4\3\5\3\6\1\6\eeee
\6\3\4\3\5\3\6\1\6\2\6\3\4\3\5\3\eeee
\1\6\2\6\3\4\3\5\3\6\1\6\2\6\3\4\eeee
\5\3\6\1\6\2\6\3\4\3\5\3\6\1\6\2\eeee
} 

\baaa
6-23
\eaaa
\bbbb
0&0&0&0&1&3\\
0&0&0&1&0&3\\
0&0&0&1&2&1\\
0&1&3&0&0&0\\
1&0&3&0&0&0\\
2&1&1&0&0&0\\
\ebbb
\parbox{7cm}{ 
\1\6\1\6\2\6\1\6\1\6\2\6\1\6\1\6\eeee
\5\3\5\3\4\3\5\3\5\3\4\3\5\3\5\3\eeee
\3\4\3\5\3\5\3\4\3\5\3\5\3\4\3\5\eeee
\6\2\6\1\6\1\6\2\6\1\6\1\6\2\6\1\eeee
\1\6\1\6\2\6\1\6\1\6\2\6\1\6\1\6\eeee
\5\3\5\3\4\3\5\3\5\3\4\3\5\3\5\3\eeee
\3\4\3\5\3\5\3\4\3\5\3\5\3\4\3\5\eeee
\6\2\6\1\6\1\6\2\6\1\6\1\6\2\6\1\eeee
\1\6\1\6\2\6\1\6\1\6\2\6\1\6\1\6\eeee
\5\3\5\3\4\3\5\3\5\3\4\3\5\3\5\3\eeee
\3\4\3\5\3\5\3\4\3\5\3\5\3\4\3\5\eeee
\6\2\6\1\6\1\6\2\6\1\6\1\6\2\6\1\eeee
} 

\baaa
6-24
\eaaa
\bbbb
0&0&0&0&1&3\\
0&0&0&1&2&1\\
0&0&0&2&2&0\\
0&2&2&0&0&0\\
1&2&1&0&0&0\\
3&1&0&0&0&0\\
\ebbb
\parbox{7cm}{ 
\1\6\1\6\1\6\1\6\1\6\1\6\1\6\1\6\eeee
\5\2\5\2\5\2\5\2\5\2\5\2\5\2\5\2\eeee
\3\4\3\4\3\4\3\4\3\4\3\4\3\4\3\4\eeee
\5\2\5\2\5\2\5\2\5\2\5\2\5\2\5\2\eeee
\1\6\1\6\1\6\1\6\1\6\1\6\1\6\1\6\eeee
\6\1\6\1\6\1\6\1\6\1\6\1\6\1\6\1\eeee
\2\5\2\5\2\5\2\5\2\5\2\5\2\5\2\5\eeee
\4\3\4\3\4\3\4\3\4\3\4\3\4\3\4\3\eeee
\2\5\2\5\2\5\2\5\2\5\2\5\2\5\2\5\eeee
\6\1\6\1\6\1\6\1\6\1\6\1\6\1\6\1\eeee
\1\6\1\6\1\6\1\6\1\6\1\6\1\6\1\6\eeee
\5\2\5\2\5\2\5\2\5\2\5\2\5\2\5\2\eeee
} 

\baaa
6-25
\eaaa
\bbbb
0&0&0&0&1&3\\
0&0&0&1&2&1\\
0&0&0&3&1&0\\
0&1&3&0&0&0\\
1&2&1&0&0&0\\
3&1&0&0&0&0\\
\ebbb
\parbox{7cm}{ 
\1\6\1\6\1\6\1\6\1\6\1\6\1\6\1\6\eeee
\5\2\5\2\5\2\5\2\5\2\5\2\5\2\5\2\eeee
\3\4\3\4\3\4\3\4\3\4\3\4\3\4\3\4\eeee
\4\3\4\3\4\3\4\3\4\3\4\3\4\3\4\3\eeee
\2\5\2\5\2\5\2\5\2\5\2\5\2\5\2\5\eeee
\6\1\6\1\6\1\6\1\6\1\6\1\6\1\6\1\eeee
\1\6\1\6\1\6\1\6\1\6\1\6\1\6\1\6\eeee
\5\2\5\2\5\2\5\2\5\2\5\2\5\2\5\2\eeee
\3\4\3\4\3\4\3\4\3\4\3\4\3\4\3\4\eeee
\4\3\4\3\4\3\4\3\4\3\4\3\4\3\4\3\eeee
\2\5\2\5\2\5\2\5\2\5\2\5\2\5\2\5\eeee
\6\1\6\1\6\1\6\1\6\1\6\1\6\1\6\1\eeee
} 

\baaa
6-26
\eaaa
\bbbb
0&0&0&0&1&3\\
0&0&0&1&2&1\\
0&0&1&2&1&0\\
0&1&2&1&0&0\\
1&2&1&0&0&0\\
3&1&0&0&0&0\\
\ebbb
\parbox{7cm}{ 
\1\6\1\6\1\6\1\6\1\6\1\6\1\6\1\6\eeee
\5\2\5\2\5\2\5\2\5\2\5\2\5\2\5\2\eeee
\3\4\3\4\3\4\3\4\3\4\3\4\3\4\3\4\eeee
\3\4\3\4\3\4\3\4\3\4\3\4\3\4\3\4\eeee
\5\2\5\2\5\2\5\2\5\2\5\2\5\2\5\2\eeee
\1\6\1\6\1\6\1\6\1\6\1\6\1\6\1\6\eeee
\6\1\6\1\6\1\6\1\6\1\6\1\6\1\6\1\eeee
\2\5\2\5\2\5\2\5\2\5\2\5\2\5\2\5\eeee
\4\3\4\3\4\3\4\3\4\3\4\3\4\3\4\3\eeee
\4\3\4\3\4\3\4\3\4\3\4\3\4\3\4\3\eeee
\2\5\2\5\2\5\2\5\2\5\2\5\2\5\2\5\eeee
\6\1\6\1\6\1\6\1\6\1\6\1\6\1\6\1\eeee
} 

\baaa
6-27
\eaaa
\bbbb
0&0&0&0&2&2\\
0&0&0&0&2&2\\
0&0&0&0&2&2\\
0&0&0&0&2&2\\
1&1&1&1&0&0\\
1&1&1&1&0&0\\
\ebbb
\parbox{7cm}{ 
\1\5\3\5\1\5\3\5\1\5\3\5\1\5\3\5\eeee
\5\2\6\4\6\2\6\4\6\2\6\4\6\2\6\4\eeee
\3\6\1\5\3\5\1\5\3\5\1\5\3\5\1\5\eeee
\5\4\5\2\6\4\6\2\6\4\6\2\6\4\6\2\eeee
\1\6\3\6\1\5\3\5\1\5\3\5\1\5\3\5\eeee
\5\2\5\4\5\2\6\4\6\2\6\4\6\2\6\4\eeee
\3\6\1\6\3\6\1\5\3\5\1\5\3\5\1\5\eeee
\5\4\5\2\5\4\5\2\6\4\6\2\6\4\6\2\eeee
\1\6\3\6\1\6\3\6\1\5\3\5\1\5\3\5\eeee
\5\2\5\4\5\2\5\4\5\2\6\4\6\2\6\4\eeee
\3\6\1\6\3\6\1\6\3\6\1\5\3\5\1\5\eeee
\5\4\5\2\5\4\5\2\5\4\5\2\6\4\6\2\eeee
} 

\baaa
6-28
\eaaa
\bbbb
0&0&0&0&2&2\\
0&0&0&0&2&2\\
0&0&0&2&0&2\\
0&0&2&0&2&0\\
1&1&0&2&0&0\\
1&1&2&0&0&0\\
\ebbb
\parbox{7cm}{ 
\1\5\4\3\6\2\5\4\3\6\2\5\4\3\6\1\eeee
\5\4\3\6\1\5\4\3\6\1\5\4\3\6\2\5\eeee
\4\3\6\2\5\4\3\6\2\5\4\3\6\1\5\4\eeee
\3\6\1\5\4\3\6\1\5\4\3\6\2\5\4\3\eeee
\6\2\5\4\3\6\2\5\4\3\6\1\5\4\3\6\eeee
\1\5\4\3\6\1\5\4\3\6\2\5\4\3\6\2\eeee
\5\4\3\6\2\5\4\3\6\1\5\4\3\6\1\5\eeee
\4\3\6\1\5\4\3\6\2\5\4\3\6\2\5\4\eeee
\3\6\2\5\4\3\6\1\5\4\3\6\1\5\4\3\eeee
\6\1\5\4\3\6\2\5\4\3\6\2\5\4\3\6\eeee
\2\5\4\3\6\1\5\4\3\6\1\5\4\3\6\2\eeee
\5\4\3\6\2\5\4\3\6\2\5\4\3\6\1\5\eeee
} 

\baaa
6-29
\eaaa
\bbbb
0&0&0&0&2&2\\
0&0&0&0&2&2\\
0&0&0&2&0&2\\
0&0&2&2&0&0\\
1&1&0&0&2&0\\
1&1&2&0&0&0\\
\ebbb
\parbox{7cm}{ 
\1\5\5\1\6\3\4\4\3\6\2\5\5\2\6\3\eeee
\5\5\2\6\3\4\4\3\6\1\5\5\1\6\3\4\eeee
\5\1\6\3\4\4\3\6\2\5\5\2\6\3\4\4\eeee
\2\6\3\4\4\3\6\1\5\5\1\6\3\4\4\3\eeee
\6\3\4\4\3\6\2\5\5\2\6\3\4\4\3\6\eeee
\3\4\4\3\6\1\5\5\1\6\3\4\4\3\6\2\eeee
\4\4\3\6\2\5\5\2\6\3\4\4\3\6\1\5\eeee
\4\3\6\1\5\5\1\6\3\4\4\3\6\2\5\5\eeee
\3\6\2\5\5\2\6\3\4\4\3\6\1\5\5\1\eeee
\6\1\5\5\1\6\3\4\4\3\6\2\5\5\2\6\eeee
\2\5\5\2\6\3\4\4\3\6\1\5\5\1\6\3\eeee
\5\5\1\6\3\4\4\3\6\2\5\5\2\6\3\4\eeee
} 

\baaa
6-30
\eaaa
\bbbb
0&0&0&0&2&2\\
0&0&0&0&2&2\\
0&0&0&2&0&2\\
0&0&2&2&0&0\\
2&2&0&0&0&0\\
1&1&2&0&0&0\\
\ebbb
\parbox{7cm}{ 
\1\5\2\6\3\4\4\3\6\2\5\1\6\3\4\4\eeee
\5\1\6\3\4\4\3\6\1\5\2\6\3\4\4\3\eeee
\2\6\3\4\4\3\6\2\5\1\6\3\4\4\3\6\eeee
\6\3\4\4\3\6\1\5\2\6\3\4\4\3\6\1\eeee
\3\4\4\3\6\2\5\1\6\3\4\4\3\6\2\5\eeee
\4\4\3\6\1\5\2\6\3\4\4\3\6\1\5\2\eeee
\4\3\6\2\5\1\6\3\4\4\3\6\2\5\1\6\eeee
\3\6\1\5\2\6\3\4\4\3\6\1\5\2\6\3\eeee
\6\2\5\1\6\3\4\4\3\6\2\5\1\6\3\4\eeee
\1\5\2\6\3\4\4\3\6\1\5\2\6\3\4\4\eeee
\5\1\6\3\4\4\3\6\2\5\1\6\3\4\4\3\eeee
\2\6\3\4\4\3\6\1\5\2\6\3\4\4\3\6\eeee
} 

\baaa
6-31
\eaaa
\bbbb
0&0&0&0&2&2\\
0&0&0&0&2&2\\
0&0&0&2&1&1\\
0&0&2&2&0&0\\
1&1&2&0&0&0\\
1&1&2&0&0&0\\
\ebbb
\parbox{7cm}{ 
\1\5\3\4\4\3\5\1\5\3\4\4\3\5\1\5\eeee
\5\2\6\3\4\4\3\6\2\6\3\4\4\3\6\2\eeee
\3\6\1\5\3\4\4\3\5\1\5\3\4\4\3\5\eeee
\4\3\5\2\6\3\4\4\3\6\2\6\3\4\4\3\eeee
\4\4\3\6\1\5\3\4\4\3\5\1\5\3\4\4\eeee
\3\4\4\3\5\2\6\3\4\4\3\6\2\6\3\4\eeee
\5\3\4\4\3\6\1\5\3\4\4\3\5\1\5\3\eeee
\1\6\3\4\4\3\5\2\6\3\4\4\3\6\2\6\eeee
\5\2\5\3\4\4\3\6\1\5\3\4\4\3\5\1\eeee
\3\6\1\6\3\4\4\3\5\2\6\3\4\4\3\6\eeee
\4\3\5\2\5\3\4\4\3\6\1\5\3\4\4\3\eeee
\4\4\3\6\1\6\3\4\4\3\5\2\6\3\4\4\eeee
} 

\baaa
6-32
\eaaa
\bbbb
0&0&0&0&2&2\\
0&0&0&0&2&2\\
0&0&1&1&0&2\\
0&0&1&1&0&2\\
1&1&0&0&2&0\\
1&1&1&1&0&0\\
\ebbb
\parbox{7cm}{ 
\1\5\5\1\6\4\4\6\1\5\5\1\6\4\4\6\eeee
\5\5\2\6\3\3\6\2\5\5\2\6\3\3\6\2\eeee
\5\1\6\4\4\6\1\5\5\1\6\4\4\6\1\5\eeee
\2\6\3\3\6\2\5\5\2\6\3\3\6\2\5\5\eeee
\6\4\4\6\1\5\5\1\6\4\4\6\1\5\5\1\eeee
\3\3\6\2\5\5\2\6\3\3\6\2\5\5\2\6\eeee
\4\6\1\5\5\1\6\4\4\6\1\5\5\1\6\4\eeee
\6\2\5\5\2\6\3\3\6\2\5\5\2\6\3\3\eeee
\1\5\5\1\6\4\4\6\1\5\5\1\6\4\4\6\eeee
\5\5\2\6\3\3\6\2\5\5\2\6\3\3\6\2\eeee
\5\1\6\4\4\6\1\5\5\1\6\4\4\6\1\5\eeee
\2\6\3\3\6\2\5\5\2\6\3\3\6\2\5\5\eeee
} 

\baaa
6-33
\eaaa
\bbbb
0&0&0&0&2&2\\
0&0&0&0&2&2\\
0&0&1&1&0&2\\
0&0&1&1&0&2\\
2&2&0&0&0&0\\
1&1&1&1&0&0\\
\ebbb
\parbox{7cm}{ 
\1\5\2\6\3\4\6\1\5\2\6\3\3\6\2\5\eeee
\5\1\6\4\3\6\2\5\1\6\4\4\6\1\5\2\eeee
\2\6\3\4\6\1\5\2\6\3\3\6\2\5\1\6\eeee
\6\4\3\6\2\5\1\6\4\4\6\1\5\2\6\3\eeee
\3\4\6\1\5\2\6\3\3\6\2\5\1\6\4\3\eeee
\3\6\2\5\1\6\4\4\6\1\5\2\6\3\4\6\eeee
\6\1\5\2\6\3\3\6\2\5\1\6\4\3\6\2\eeee
\2\5\1\6\4\4\6\1\5\2\6\3\4\6\1\5\eeee
\5\2\6\3\3\6\2\5\1\6\4\3\6\2\5\1\eeee
\1\6\4\4\6\1\5\2\6\3\4\6\1\5\2\6\eeee
\6\3\3\6\2\5\1\6\4\3\6\2\5\1\6\4\eeee
\4\4\6\1\5\2\6\3\4\6\1\5\2\6\3\3\eeee
} 

\baaa
6-34
\eaaa
\bbbb
0&0&0&0&2&2\\
0&0&0&0&2&2\\
0&0&1&1&1&1\\
0&0&1&1&1&1\\
1&1&1&1&0&0\\
1&1&1&1&0&0\\
\ebbb
\parbox{7cm}{ 
\1\5\3\3\5\1\5\3\3\5\1\5\3\3\5\1\eeee
\5\2\6\4\4\6\2\6\4\4\6\2\6\4\4\6\eeee
\3\6\1\5\3\3\5\1\5\3\3\5\1\5\3\3\eeee
\3\4\5\2\6\4\4\6\2\6\4\4\6\2\6\4\eeee
\5\4\3\6\1\5\3\3\5\1\5\3\3\5\1\5\eeee
\1\6\3\4\5\2\6\4\4\6\2\6\4\4\6\2\eeee
\5\2\5\4\3\6\1\5\3\3\5\1\5\3\3\5\eeee
\3\6\1\6\3\4\5\2\6\4\4\6\2\6\4\4\eeee
\3\4\5\2\5\4\3\6\1\5\3\3\5\1\5\3\eeee
\5\4\3\6\1\6\3\4\5\2\6\4\4\6\2\6\eeee
\1\6\3\4\5\2\5\4\3\6\1\5\3\3\5\1\eeee
\5\2\5\4\3\6\1\6\3\4\5\2\6\4\4\6\eeee
} 

\baaa
6-35
\eaaa
\bbbb
0&0&0&0&2&2\\
0&0&0&0&2&2\\
0&0&2&0&0&2\\
0&0&0&2&2&0\\
1&1&0&2&0&0\\
1&1&2&0&0&0\\
\ebbb
\parbox{7cm}{ 
\1\5\4\4\5\2\6\3\3\6\2\5\4\4\5\1\eeee
\5\4\4\5\1\6\3\3\6\1\5\4\4\5\2\6\eeee
\4\4\5\2\6\3\3\6\2\5\4\4\5\1\6\3\eeee
\4\5\1\6\3\3\6\1\5\4\4\5\2\6\3\3\eeee
\5\2\6\3\3\6\2\5\4\4\5\1\6\3\3\6\eeee
\1\6\3\3\6\1\5\4\4\5\2\6\3\3\6\2\eeee
\6\3\3\6\2\5\4\4\5\1\6\3\3\6\1\5\eeee
\3\3\6\1\5\4\4\5\2\6\3\3\6\2\5\4\eeee
\3\6\2\5\4\4\5\1\6\3\3\6\1\5\4\4\eeee
\6\1\5\4\4\5\2\6\3\3\6\2\5\4\4\5\eeee
\2\5\4\4\5\1\6\3\3\6\1\5\4\4\5\2\eeee
\5\4\4\5\2\6\3\3\6\2\5\4\4\5\1\6\eeee
} 

\baaa
6-36
\eaaa
\bbbb
0&0&0&0&2&2\\
0&0&0&1&1&2\\
0&0&0&1&2&1\\
0&2&2&0&0&0\\
1&1&2&0&0&0\\
1&2&1&0&0&0\\
\ebbb
\parbox{7cm}{ 
\1\5\3\5\3\4\2\6\2\6\1\5\3\5\3\4\eeee
\5\3\4\2\6\2\6\1\5\3\5\3\4\2\6\2\eeee
\2\6\2\6\1\5\3\5\3\4\2\6\2\6\1\5\eeee
\6\1\5\3\5\3\4\2\6\2\6\1\5\3\5\3\eeee
\3\5\3\4\2\6\2\6\1\5\3\5\3\4\2\6\eeee
\4\2\6\2\6\1\5\3\5\3\4\2\6\2\6\1\eeee
\2\6\1\5\3\5\3\4\2\6\2\6\1\5\3\5\eeee
\5\3\5\3\4\2\6\2\6\1\5\3\5\3\4\2\eeee
\3\4\2\6\2\6\1\5\3\5\3\4\2\6\2\6\eeee
\6\2\6\1\5\3\5\3\4\2\6\2\6\1\5\3\eeee
\1\5\3\5\3\4\2\6\2\6\1\5\3\5\3\4\eeee
\5\3\4\2\6\2\6\1\5\3\5\3\4\2\6\2\eeee
} 
\baaa
\phantom{0-0}\#
\eaaa
\mbox{}\phantom{\bbbb
0&0&0&0&0&0\\
\ebbb}
\parbox{7cm}{ 
\1\6\2\5\3\4\3\5\2\6\1\6\2\5\3\4\eeee
\5\2\6\1\6\2\5\3\4\3\5\2\6\1\6\2\eeee
\3\4\3\5\2\6\1\6\2\5\3\4\3\5\2\6\eeee
\6\2\5\3\4\3\5\2\6\1\6\2\5\3\4\3\eeee
\2\6\1\6\2\5\3\4\3\5\2\6\1\6\2\5\eeee
\4\3\5\2\6\1\6\2\5\3\4\3\5\2\6\1\eeee
\2\5\3\4\3\5\2\6\1\6\2\5\3\4\3\5\eeee
\6\1\6\2\5\3\4\3\5\2\6\1\6\2\5\3\eeee
\3\5\2\6\1\6\2\5\3\4\3\5\2\6\1\6\eeee
\5\3\4\3\5\2\6\1\6\2\5\3\4\3\5\2\eeee
\1\6\2\5\3\4\3\5\2\6\1\6\2\5\3\4\eeee
\5\2\6\1\6\2\5\3\4\3\5\2\6\1\6\2\eeee
} 

\baaa
6-37
\eaaa
\bbbb
0&0&0&0&2&2\\
0&0&0&1&1&2\\
0&0&0&2&0&2\\
0&2&2&0&0&0\\
2&2&0&0&0&0\\
1&2&1&0&0&0\\
\ebbb
\parbox{7cm}{ 
\1\5\2\6\1\5\2\6\1\5\2\6\1\5\2\6\eeee
\5\1\6\2\5\1\6\2\5\1\6\2\5\1\6\2\eeee
\2\6\3\4\2\6\3\4\2\6\3\4\2\6\3\4\eeee
\6\2\4\3\6\2\4\3\6\2\4\3\6\2\4\3\eeee
\1\5\2\6\1\5\2\6\1\5\2\6\1\5\2\6\eeee
\5\1\6\2\5\1\6\2\5\1\6\2\5\1\6\2\eeee
\2\6\3\4\2\6\3\4\2\6\3\4\2\6\3\4\eeee
\6\2\4\3\6\2\4\3\6\2\4\3\6\2\4\3\eeee
\1\5\2\6\1\5\2\6\1\5\2\6\1\5\2\6\eeee
\5\1\6\2\5\1\6\2\5\1\6\2\5\1\6\2\eeee
\2\6\3\4\2\6\3\4\2\6\3\4\2\6\3\4\eeee
\6\2\4\3\6\2\4\3\6\2\4\3\6\2\4\3\eeee
} 
\baaa
\phantom{0-0}\#
\eaaa
\mbox{}\phantom{\bbbb
0&0&0&0&0&0\\
\ebbb}
\parbox{7cm}{ 
\1\6\3\6\1\6\3\6\1\6\3\6\1\6\3\6\eeee
\5\2\4\2\5\2\4\2\5\2\4\2\5\2\4\2\eeee
\1\6\3\6\1\6\3\6\1\6\3\6\1\6\3\6\eeee
\5\2\4\2\5\2\4\2\5\2\4\2\5\2\4\2\eeee
\1\6\3\6\1\6\3\6\1\6\3\6\1\6\3\6\eeee
\5\2\4\2\5\2\4\2\5\2\4\2\5\2\4\2\eeee
\1\6\3\6\1\6\3\6\1\6\3\6\1\6\3\6\eeee
\5\2\4\2\5\2\4\2\5\2\4\2\5\2\4\2\eeee
\1\6\3\6\1\6\3\6\1\6\3\6\1\6\3\6\eeee
\5\2\4\2\5\2\4\2\5\2\4\2\5\2\4\2\eeee
\1\6\3\6\1\6\3\6\1\6\3\6\1\6\3\6\eeee
\5\2\4\2\5\2\4\2\5\2\4\2\5\2\4\2\eeee
} 

\baaa
6-38
\eaaa
\bbbb
0&0&0&0&2&2\\
0&0&0&1&1&2\\
0&0&1&1&1&1\\
0&1&1&1&1&0\\
1&1&1&1&0&0\\
1&2&1&0&0&0\\
\ebbb
\parbox{7cm}{ 
\1\6\2\5\3\4\4\3\5\2\6\1\6\2\5\3\eeee
\5\2\6\1\6\2\5\3\4\4\3\5\2\6\1\6\eeee
\4\4\3\5\2\6\1\6\2\5\3\4\4\3\5\2\eeee
\2\5\3\4\4\3\5\2\6\1\6\2\5\3\4\4\eeee
\6\1\6\2\5\3\4\4\3\5\2\6\1\6\2\5\eeee
\3\5\2\6\1\6\2\5\3\4\4\3\5\2\6\1\eeee
\3\4\4\3\5\2\6\1\6\2\5\3\4\4\3\5\eeee
\6\2\5\3\4\4\3\5\2\6\1\6\2\5\3\4\eeee
\2\6\1\6\2\5\3\4\4\3\5\2\6\1\6\2\eeee
\4\3\5\2\6\1\6\2\5\3\4\4\3\5\2\6\eeee
\5\3\4\4\3\5\2\6\1\6\2\5\3\4\4\3\eeee
\1\6\2\5\3\4\4\3\5\2\6\1\6\2\5\3\eeee
} 

\baaa
6-39
\eaaa
\bbbb
0&0&0&0&2&2\\
0&0&0&1&1&2\\
0&0&1&2&0&1\\
0&1&2&1&0&0\\
2&2&0&0&0&0\\
1&2&1&0&0&0\\
\ebbb
\parbox{7cm}{ 
\1\6\3\3\6\1\6\3\3\6\1\6\3\3\6\1\eeee
\5\2\4\4\2\5\2\4\4\2\5\2\4\4\2\5\eeee
\1\6\3\3\6\1\6\3\3\6\1\6\3\3\6\1\eeee
\5\2\4\4\2\5\2\4\4\2\5\2\4\4\2\5\eeee
\1\6\3\3\6\1\6\3\3\6\1\6\3\3\6\1\eeee
\5\2\4\4\2\5\2\4\4\2\5\2\4\4\2\5\eeee
\1\6\3\3\6\1\6\3\3\6\1\6\3\3\6\1\eeee
\5\2\4\4\2\5\2\4\4\2\5\2\4\4\2\5\eeee
\1\6\3\3\6\1\6\3\3\6\1\6\3\3\6\1\eeee
\5\2\4\4\2\5\2\4\4\2\5\2\4\4\2\5\eeee
\1\6\3\3\6\1\6\3\3\6\1\6\3\3\6\1\eeee
\5\2\4\4\2\5\2\4\4\2\5\2\4\4\2\5\eeee
} 

\baaa
6-40
\eaaa
\bbbb
0&0&0&0&2&2\\
0&0&0&2&0&2\\
0&0&0&2&2&0\\
0&1&1&0&1&1\\
1&0&1&1&0&1\\
1&1&0&1&1&0\\
\ebbb
\parbox{7cm}{ 
\1\5\6\1\5\6\1\5\6\1\5\6\1\5\6\1\eeee
\5\3\4\5\3\4\5\3\4\5\3\4\5\3\4\5\eeee
\6\4\2\6\4\2\6\4\2\6\4\2\6\4\2\6\eeee
\1\5\6\1\5\6\1\5\6\1\5\6\1\5\6\1\eeee
\5\3\4\5\3\4\5\3\4\5\3\4\5\3\4\5\eeee
\6\4\2\6\4\2\6\4\2\6\4\2\6\4\2\6\eeee
\1\5\6\1\5\6\1\5\6\1\5\6\1\5\6\1\eeee
\5\3\4\5\3\4\5\3\4\5\3\4\5\3\4\5\eeee
\6\4\2\6\4\2\6\4\2\6\4\2\6\4\2\6\eeee
\1\5\6\1\5\6\1\5\6\1\5\6\1\5\6\1\eeee
\5\3\4\5\3\4\5\3\4\5\3\4\5\3\4\5\eeee
\6\4\2\6\4\2\6\4\2\6\4\2\6\4\2\6\eeee
} 

\baaa
6-41
\eaaa
\bbbb
0&0&0&0&2&2\\
0&0&0&2&0&2\\
0&0&0&2&2&0\\
0&1&1&0&1&1\\
1&0&1&1&1&0\\
1&1&0&1&0&1\\
\ebbb
\parbox{7cm}{ 
\1\5\5\1\6\6\1\5\5\1\6\6\1\5\5\1\eeee
\5\3\4\6\2\4\5\3\4\6\2\4\5\3\4\6\eeee
\5\4\2\6\4\3\5\4\2\6\4\3\5\4\2\6\eeee
\1\6\6\1\5\5\1\6\6\1\5\5\1\6\6\1\eeee
\6\2\4\5\3\4\6\2\4\5\3\4\6\2\4\5\eeee
\6\4\3\5\4\2\6\4\3\5\4\2\6\4\3\5\eeee
\1\5\5\1\6\6\1\5\5\1\6\6\1\5\5\1\eeee
\5\3\4\6\2\4\5\3\4\6\2\4\5\3\4\6\eeee
\5\4\2\6\4\3\5\4\2\6\4\3\5\4\2\6\eeee
\1\6\6\1\5\5\1\6\6\1\5\5\1\6\6\1\eeee
\6\2\4\5\3\4\6\2\4\5\3\4\6\2\4\5\eeee
\6\4\3\5\4\2\6\4\3\5\4\2\6\4\3\5\eeee
} 

\baaa
6-42
\eaaa
\bbbb
0&0&0&0&2&2\\
0&0&0&2&0&2\\
0&0&0&2&2&0\\
0&2&2&0&0&0\\
2&0&2&0&0&0\\
2&2&0&0&0&0\\
\ebbb
\parbox{7cm}{ 
\1\5\3\4\2\6\1\5\3\4\2\6\1\5\3\4\eeee
\5\3\4\2\6\1\5\3\4\2\6\1\5\3\4\2\eeee
\3\4\2\6\1\5\3\4\2\6\1\5\3\4\2\6\eeee
\4\2\6\1\5\3\4\2\6\1\5\3\4\2\6\1\eeee
\2\6\1\5\3\4\2\6\1\5\3\4\2\6\1\5\eeee
\6\1\5\3\4\2\6\1\5\3\4\2\6\1\5\3\eeee
\1\5\3\4\2\6\1\5\3\4\2\6\1\5\3\4\eeee
\5\3\4\2\6\1\5\3\4\2\6\1\5\3\4\2\eeee
\3\4\2\6\1\5\3\4\2\6\1\5\3\4\2\6\eeee
\4\2\6\1\5\3\4\2\6\1\5\3\4\2\6\1\eeee
\2\6\1\5\3\4\2\6\1\5\3\4\2\6\1\5\eeee
\6\1\5\3\4\2\6\1\5\3\4\2\6\1\5\3\eeee
} 

\baaa
6-43
\eaaa
\bbbb
0&0&0&0&2&2\\
0&0&0&2&0&2\\
0&0&1&1&1&1\\
0&1&2&0&1&0\\
1&0&2&1&0&0\\
1&1&2&0&0&0\\
\ebbb
\parbox{7cm}{ 
\1\6\2\6\1\6\2\6\1\6\2\6\1\6\2\6\eeee
\5\3\4\3\5\3\4\3\5\3\4\3\5\3\4\3\eeee
\4\3\5\3\4\3\5\3\4\3\5\3\4\3\5\3\eeee
\2\6\1\6\2\6\1\6\2\6\1\6\2\6\1\6\eeee
\4\3\5\3\4\3\5\3\4\3\5\3\4\3\5\3\eeee
\5\3\4\3\5\3\4\3\5\3\4\3\5\3\4\3\eeee
\1\6\2\6\1\6\2\6\1\6\2\6\1\6\2\6\eeee
\5\3\4\3\5\3\4\3\5\3\4\3\5\3\4\3\eeee
\4\3\5\3\4\3\5\3\4\3\5\3\4\3\5\3\eeee
\2\6\1\6\2\6\1\6\2\6\1\6\2\6\1\6\eeee
\4\3\5\3\4\3\5\3\4\3\5\3\4\3\5\3\eeee
\5\3\4\3\5\3\4\3\5\3\4\3\5\3\4\3\eeee
} 

\baaa
6-44
\eaaa
\bbbb
0&0&0&0&2&2\\
0&0&0&2&0&2\\
0&0&1&1&1&1\\
0&1&2&1&0&0\\
1&0&2&0&1&0\\
1&1&2&0&0&0\\
\ebbb
\parbox{7cm}{ 
\1\6\2\6\1\6\2\6\1\6\2\6\1\6\2\6\eeee
\5\3\4\3\5\3\4\3\5\3\4\3\5\3\4\3\eeee
\5\3\4\3\5\3\4\3\5\3\4\3\5\3\4\3\eeee
\1\6\2\6\1\6\2\6\1\6\2\6\1\6\2\6\eeee
\5\3\4\3\5\3\4\3\5\3\4\3\5\3\4\3\eeee
\5\3\4\3\5\3\4\3\5\3\4\3\5\3\4\3\eeee
\1\6\2\6\1\6\2\6\1\6\2\6\1\6\2\6\eeee
\5\3\4\3\5\3\4\3\5\3\4\3\5\3\4\3\eeee
\5\3\4\3\5\3\4\3\5\3\4\3\5\3\4\3\eeee
\1\6\2\6\1\6\2\6\1\6\2\6\1\6\2\6\eeee
\5\3\4\3\5\3\4\3\5\3\4\3\5\3\4\3\eeee
\5\3\4\3\5\3\4\3\5\3\4\3\5\3\4\3\eeee
} 

\baaa
6-45
\eaaa
\bbbb
0&0&0&0&2&2\\
0&0&0&2&0&2\\
0&0&2&0&1&1\\
0&2&0&0&2&0\\
1&0&1&1&0&1\\
1&1&1&0&1&0\\
\ebbb
\parbox{7cm}{ 
\1\5\6\1\5\6\1\5\6\1\5\6\1\5\6\1\eeee
\5\4\2\6\3\3\5\4\2\6\3\3\5\4\2\6\eeee
\6\2\4\5\3\3\6\2\4\5\3\3\6\2\4\5\eeee
\1\6\5\1\6\5\1\6\5\1\6\5\1\6\5\1\eeee
\5\3\3\6\2\4\5\3\3\6\2\4\5\3\3\6\eeee
\6\3\3\5\4\2\6\3\3\5\4\2\6\3\3\5\eeee
\1\5\6\1\5\6\1\5\6\1\5\6\1\5\6\1\eeee
\5\4\2\6\3\3\5\4\2\6\3\3\5\4\2\6\eeee
\6\2\4\5\3\3\6\2\4\5\3\3\6\2\4\5\eeee
\1\6\5\1\6\5\1\6\5\1\6\5\1\6\5\1\eeee
\5\3\3\6\2\4\5\3\3\6\2\4\5\3\3\6\eeee
\6\3\3\5\4\2\6\3\3\5\4\2\6\3\3\5\eeee
} 

\baaa
6-46
\eaaa
\bbbb
0&0&0&0&2&2\\
0&0&0&2&0&2\\
0&0&2&0&2&0\\
0&2&0&2&0&0\\
2&0&2&0&0&0\\
2&2&0&0&0&0\\
\ebbb
\parbox{7cm}{ 
\1\5\3\3\5\1\6\2\4\4\2\6\1\5\3\3\eeee
\5\3\3\5\1\6\2\4\4\2\6\1\5\3\3\5\eeee
\3\3\5\1\6\2\4\4\2\6\1\5\3\3\5\1\eeee
\3\5\1\6\2\4\4\2\6\1\5\3\3\5\1\6\eeee
\5\1\6\2\4\4\2\6\1\5\3\3\5\1\6\2\eeee
\1\6\2\4\4\2\6\1\5\3\3\5\1\6\2\4\eeee
\6\2\4\4\2\6\1\5\3\3\5\1\6\2\4\4\eeee
\2\4\4\2\6\1\5\3\3\5\1\6\2\4\4\2\eeee
\4\4\2\6\1\5\3\3\5\1\6\2\4\4\2\6\eeee
\4\2\6\1\5\3\3\5\1\6\2\4\4\2\6\1\eeee
\2\6\1\5\3\3\5\1\6\2\4\4\2\6\1\5\eeee
\6\1\5\3\3\5\1\6\2\4\4\2\6\1\5\3\eeee
} 

\baaa
6-47
\eaaa
\bbbb
0&0&0&0&2&2\\
0&0&0&2&1&1\\
0&0&0&2&1&1\\
0&1&1&2&0&0\\
2&1&1&0&0&0\\
2&1&1&0&0&0\\
\ebbb
\parbox{7cm}{ 
\1\5\2\4\4\2\5\1\5\2\4\4\2\5\1\5\eeee
\5\1\6\3\4\4\3\6\1\6\3\4\4\3\6\1\eeee
\2\6\1\5\2\4\4\2\5\1\5\2\4\4\2\5\eeee
\4\3\5\1\6\3\4\4\3\6\1\6\3\4\4\3\eeee
\4\4\2\6\1\5\2\4\4\2\5\1\5\2\4\4\eeee
\2\4\4\3\5\1\6\3\4\4\3\6\1\6\3\4\eeee
\5\3\4\4\2\6\1\5\2\4\4\2\5\1\5\2\eeee
\1\6\2\4\4\3\5\1\6\3\4\4\3\6\1\6\eeee
\5\1\5\3\4\4\2\6\1\5\2\4\4\2\5\1\eeee
\2\6\1\6\2\4\4\3\5\1\6\3\4\4\3\6\eeee
\4\3\5\1\5\3\4\4\2\6\1\5\2\4\4\2\eeee
\4\4\2\6\1\6\2\4\4\3\5\1\6\3\4\4\eeee
} 

\baaa
6-48
\eaaa
\bbbb
0&0&0&0&2&2\\
0&0&0&2&1&1\\
0&0&0&2&1&1\\
0&2&2&0&0&0\\
2&1&1&0&0&0\\
2&1&1&0&0&0\\
\ebbb
\parbox{7cm}{ 
\1\5\2\4\2\5\1\5\2\4\2\5\1\5\2\4\eeee
\5\1\6\3\4\3\6\1\6\3\4\3\6\1\6\3\eeee
\2\6\1\5\2\4\2\5\1\5\2\4\2\5\1\5\eeee
\4\3\5\1\6\3\4\3\6\1\6\3\4\3\6\1\eeee
\2\4\2\6\1\5\2\4\2\5\1\5\2\4\2\5\eeee
\5\3\4\3\5\1\6\3\4\3\6\1\6\3\4\3\eeee
\1\6\2\4\2\6\1\5\2\4\2\5\1\5\2\4\eeee
\5\1\5\3\4\3\5\1\6\3\4\3\6\1\6\3\eeee
\2\6\1\6\2\4\2\6\1\5\2\4\2\5\1\5\eeee
\4\3\5\1\5\3\4\3\5\1\6\3\4\3\6\1\eeee
\2\4\2\6\1\6\2\4\2\6\1\5\2\4\2\5\eeee
\5\3\4\3\5\1\5\3\4\3\5\1\6\3\4\3\eeee
} 

\baaa
6-49
\eaaa
\bbbb
0&0&0&0&2&2\\
0&0&1&1&0&2\\
0&2&1&1&0&0\\
0&2&1&1&0&0\\
2&0&0&0&2&0\\
2&2&0&0&0&0\\
\ebbb
\parbox{7cm}{ 
\1\5\5\1\6\2\4\3\2\6\1\5\5\1\6\2\eeee
\5\5\1\6\2\3\4\2\6\1\5\5\1\6\2\4\eeee
\5\1\6\2\4\3\2\6\1\5\5\1\6\2\3\4\eeee
\1\6\2\3\4\2\6\1\5\5\1\6\2\4\3\2\eeee
\6\2\4\3\2\6\1\5\5\1\6\2\3\4\2\6\eeee
\2\3\4\2\6\1\5\5\1\6\2\4\3\2\6\1\eeee
\4\3\2\6\1\5\5\1\6\2\3\4\2\6\1\5\eeee
\4\2\6\1\5\5\1\6\2\4\3\2\6\1\5\5\eeee
\2\6\1\5\5\1\6\2\3\4\2\6\1\5\5\1\eeee
\6\1\5\5\1\6\2\4\3\2\6\1\5\5\1\6\eeee
\1\5\5\1\6\2\3\4\2\6\1\5\5\1\6\2\eeee
\5\5\1\6\2\4\3\2\6\1\5\5\1\6\2\3\eeee
} 

\baaa
6-50
\eaaa
\bbbb
0&0&0&0&2&2\\
0&0&1&1&1&1\\
0&1&0&2&0&1\\
0&1&2&1&0&0\\
1&1&0&0&1&1\\
1&1&1&0&1&0\\
\ebbb
\parbox{7cm}{ 
\1\6\3\4\2\5\5\2\4\3\6\1\6\3\4\2\eeee
\5\2\4\3\6\1\6\3\4\2\5\5\2\4\3\6\eeee
\6\3\4\2\5\5\2\4\3\6\1\6\3\4\2\5\eeee
\2\4\3\6\1\6\3\4\2\5\5\2\4\3\6\1\eeee
\3\4\2\5\5\2\4\3\6\1\6\3\4\2\5\5\eeee
\4\3\6\1\6\3\4\2\5\5\2\4\3\6\1\6\eeee
\4\2\5\5\2\4\3\6\1\6\3\4\2\5\5\2\eeee
\3\6\1\6\3\4\2\5\5\2\4\3\6\1\6\3\eeee
\2\5\5\2\4\3\6\1\6\3\4\2\5\5\2\4\eeee
\6\1\6\3\4\2\5\5\2\4\3\6\1\6\3\4\eeee
\5\5\2\4\3\6\1\6\3\4\2\5\5\2\4\3\eeee
\1\6\3\4\2\5\5\2\4\3\6\1\6\3\4\2\eeee
} 

\baaa
6-51
\eaaa
\bbbb
0&0&0&0&2&2\\
0&0&1&1&1&1\\
0&1&1&1&0&1\\
0&2&2&0&0&0\\
1&1&0&0&1&1\\
1&1&1&0&1&0\\
\ebbb
\parbox{7cm}{ 
\1\6\3\3\6\1\6\3\3\6\1\6\3\3\6\1\eeee
\5\2\4\2\5\5\2\4\2\5\5\2\4\2\5\5\eeee
\6\3\3\6\1\6\3\3\6\1\6\3\3\6\1\6\eeee
\2\4\2\5\5\2\4\2\5\5\2\4\2\5\5\2\eeee
\3\3\6\1\6\3\3\6\1\6\3\3\6\1\6\3\eeee
\4\2\5\5\2\4\2\5\5\2\4\2\5\5\2\4\eeee
\3\6\1\6\3\3\6\1\6\3\3\6\1\6\3\3\eeee
\2\5\5\2\4\2\5\5\2\4\2\5\5\2\4\2\eeee
\6\1\6\3\3\6\1\6\3\3\6\1\6\3\3\6\eeee
\5\5\2\4\2\5\5\2\4\2\5\5\2\4\2\5\eeee
\1\6\3\3\6\1\6\3\3\6\1\6\3\3\6\1\eeee
\5\2\4\2\5\5\2\4\2\5\5\2\4\2\5\5\eeee
} 

\baaa
6-52
\eaaa
\bbbb
0&0&0&0&2&2\\
0&1&0&1&1&1\\
0&0&1&2&0&1\\
0&1&1&2&0&0\\
1&2&0&0&1&0\\
1&2&1&0&0&0\\
\ebbb
\parbox{7cm}{ 
\1\6\3\3\6\1\6\3\3\6\1\6\3\3\6\1\eeee
\5\2\4\4\2\5\2\4\4\2\5\2\4\4\2\5\eeee
\5\2\4\4\2\5\2\4\4\2\5\2\4\4\2\5\eeee
\1\6\3\3\6\1\6\3\3\6\1\6\3\3\6\1\eeee
\5\2\4\4\2\5\2\4\4\2\5\2\4\4\2\5\eeee
\5\2\4\4\2\5\2\4\4\2\5\2\4\4\2\5\eeee
\1\6\3\3\6\1\6\3\3\6\1\6\3\3\6\1\eeee
\5\2\4\4\2\5\2\4\4\2\5\2\4\4\2\5\eeee
\5\2\4\4\2\5\2\4\4\2\5\2\4\4\2\5\eeee
\1\6\3\3\6\1\6\3\3\6\1\6\3\3\6\1\eeee
\5\2\4\4\2\5\2\4\4\2\5\2\4\4\2\5\eeee
\5\2\4\4\2\5\2\4\4\2\5\2\4\4\2\5\eeee
} 

\baaa
6-53
\eaaa
\bbbb
0&0&0&0&2&2\\
0&1&0&1&1&1\\
0&0&2&0&0&2\\
0&2&0&2&0&0\\
2&2&0&0&0&0\\
1&1&1&0&0&1\\
\ebbb
\parbox{7cm}{ 
\1\5\2\2\5\1\6\6\1\5\2\2\5\1\6\6\eeee
\5\1\6\6\1\5\2\2\5\1\6\6\1\5\2\2\eeee
\2\6\3\3\6\2\4\4\2\6\3\3\6\2\4\4\eeee
\2\6\3\3\6\2\4\4\2\6\3\3\6\2\4\4\eeee
\5\1\6\6\1\5\2\2\5\1\6\6\1\5\2\2\eeee
\1\5\2\2\5\1\6\6\1\5\2\2\5\1\6\6\eeee
\6\2\4\4\2\6\3\3\6\2\4\4\2\6\3\3\eeee
\6\2\4\4\2\6\3\3\6\2\4\4\2\6\3\3\eeee
\1\5\2\2\5\1\6\6\1\5\2\2\5\1\6\6\eeee
\5\1\6\6\1\5\2\2\5\1\6\6\1\5\2\2\eeee
\2\6\3\3\6\2\4\4\2\6\3\3\6\2\4\4\eeee
\2\6\3\3\6\2\4\4\2\6\3\3\6\2\4\4\eeee
} 

\baaa
6-54
\eaaa
\bbbb
0&0&0&0&2&2\\
0&2&0&0&0&2\\
0&0&2&0&1&1\\
0&0&0&2&2&0\\
1&0&1&1&1&0\\
1&1&1&0&0&1\\
\ebbb
\parbox{7cm}{ 
\1\5\5\1\6\6\1\5\5\1\6\6\1\5\5\1\eeee
\5\4\4\5\3\3\5\4\4\5\3\3\5\4\4\5\eeee
\5\4\4\5\3\3\5\4\4\5\3\3\5\4\4\5\eeee
\1\5\5\1\6\6\1\5\5\1\6\6\1\5\5\1\eeee
\6\3\3\6\2\2\6\3\3\6\2\2\6\3\3\6\eeee
\6\3\3\6\2\2\6\3\3\6\2\2\6\3\3\6\eeee
\1\5\5\1\6\6\1\5\5\1\6\6\1\5\5\1\eeee
\5\4\4\5\3\3\5\4\4\5\3\3\5\4\4\5\eeee
\5\4\4\5\3\3\5\4\4\5\3\3\5\4\4\5\eeee
\1\5\5\1\6\6\1\5\5\1\6\6\1\5\5\1\eeee
\6\3\3\6\2\2\6\3\3\6\2\2\6\3\3\6\eeee
\6\3\3\6\2\2\6\3\3\6\2\2\6\3\3\6\eeee
} 

\baaa
6-55
\eaaa
\bbbb
0&0&0&1&1&2\\
0&0&0&1&1&2\\
0&0&0&1&1&2\\
1&1&2&0&0&0\\
1&1&2&0&0&0\\
1&1&2&0&0&0\\
\ebbb
\parbox{7cm}{ 
\1\5\3\4\1\5\3\4\1\5\3\4\1\5\3\4\eeee
\4\2\6\3\6\2\6\3\6\2\6\3\6\2\6\3\eeee
\3\6\1\5\3\4\1\5\3\4\1\5\3\4\1\5\eeee
\5\3\4\2\6\3\6\2\6\3\6\2\6\3\6\2\eeee
\1\6\3\6\1\5\3\4\1\5\3\4\1\5\3\4\eeee
\4\2\5\3\4\2\6\3\6\2\6\3\6\2\6\3\eeee
\3\6\1\6\3\6\1\5\3\4\1\5\3\4\1\5\eeee
\5\3\4\2\5\3\4\2\6\3\6\2\6\3\6\2\eeee
\1\6\3\6\1\6\3\6\1\5\3\4\1\5\3\4\eeee
\4\2\5\3\4\2\5\3\4\2\6\3\6\2\6\3\eeee
\3\6\1\6\3\6\1\6\3\6\1\5\3\4\1\5\eeee
\5\3\4\2\5\3\4\2\5\3\4\2\6\3\6\2\eeee
} 

\baaa
6-56
\eaaa
\bbbb
0&0&0&1&1&2\\
0&0&0&1&1&2\\
0&0&2&0&0&2\\
1&1&0&1&1&0\\
1&1&0&1&1&0\\
1&1&2&0&0&0\\
\ebbb
\parbox{7cm}{ 
\1\5\5\1\6\3\3\6\1\5\5\1\6\3\3\6\eeee
\4\4\2\6\3\3\6\2\4\4\2\6\3\3\6\2\eeee
\5\1\6\3\3\6\1\5\5\1\6\3\3\6\1\5\eeee
\2\6\3\3\6\2\4\4\2\6\3\3\6\2\4\4\eeee
\6\3\3\6\1\5\5\1\6\3\3\6\1\5\5\1\eeee
\3\3\6\2\4\4\2\6\3\3\6\2\4\4\2\6\eeee
\3\6\1\5\5\1\6\3\3\6\1\5\5\1\6\3\eeee
\6\2\4\4\2\6\3\3\6\2\4\4\2\6\3\3\eeee
\1\5\5\1\6\3\3\6\1\5\5\1\6\3\3\6\eeee
\4\4\2\6\3\3\6\2\4\4\2\6\3\3\6\2\eeee
\5\1\6\3\3\6\1\5\5\1\6\3\3\6\1\5\eeee
\2\6\3\3\6\2\4\4\2\6\3\3\6\2\4\4\eeee
} 

\baaa
6-57
\eaaa
\bbbb
0&0&0&1&1&2\\
0&0&0&1&1&2\\
0&0&2&1&1&0\\
1&1&2&0&0&0\\
1&1&2&0&0&0\\
1&1&0&0&0&2\\
\ebbb
\parbox{7cm}{ 
\1\5\3\3\4\2\6\6\1\5\3\3\4\2\6\6\eeee
\4\3\3\5\1\6\6\2\4\3\3\5\1\6\6\2\eeee
\3\3\4\2\6\6\1\5\3\3\4\2\6\6\1\5\eeee
\3\5\1\6\6\2\4\3\3\5\1\6\6\2\4\3\eeee
\4\2\6\6\1\5\3\3\4\2\6\6\1\5\3\3\eeee
\1\6\6\2\4\3\3\5\1\6\6\2\4\3\3\5\eeee
\6\6\1\5\3\3\4\2\6\6\1\5\3\3\4\2\eeee
\6\2\4\3\3\5\1\6\6\2\4\3\3\5\1\6\eeee
\1\5\3\3\4\2\6\6\1\5\3\3\4\2\6\6\eeee
\4\3\3\5\1\6\6\2\4\3\3\5\1\6\6\2\eeee
\3\3\4\2\6\6\1\5\3\3\4\2\6\6\1\5\eeee
\3\5\1\6\6\2\4\3\3\5\1\6\6\2\4\3\eeee
} 

\baaa
6-58
\eaaa
\bbbb
0&0&0&1&1&2\\
0&0&0&1&2&1\\
0&0&0&2&1&1\\
1&1&2&0&0&0\\
1&2&1&0&0&0\\
2&1&1&0&0&0\\
\ebbb
\parbox{7cm}{ 
\1\5\3\4\2\6\1\5\3\4\2\6\1\5\3\4\eeee
\4\2\6\1\5\3\4\2\6\1\5\3\4\2\6\1\eeee
\3\5\1\6\2\4\3\5\1\6\2\4\3\5\1\6\eeee
\6\2\4\3\5\1\6\2\4\3\5\1\6\2\4\3\eeee
\1\5\3\4\2\6\1\5\3\4\2\6\1\5\3\4\eeee
\4\2\6\1\5\3\4\2\6\1\5\3\4\2\6\1\eeee
\3\5\1\6\2\4\3\5\1\6\2\4\3\5\1\6\eeee
\6\2\4\3\5\1\6\2\4\3\5\1\6\2\4\3\eeee
\1\5\3\4\2\6\1\5\3\4\2\6\1\5\3\4\eeee
\4\2\6\1\5\3\4\2\6\1\5\3\4\2\6\1\eeee
\3\5\1\6\2\4\3\5\1\6\2\4\3\5\1\6\eeee
\6\2\4\3\5\1\6\2\4\3\5\1\6\2\4\3\eeee
} 
\baaa
\phantom{0-0}\#
\eaaa
\mbox{}\phantom{\bbbb
0&0&0&0&0&0\\
\ebbb}
\parbox{7cm}{ 
\1\5\2\6\1\5\2\6\1\5\2\6\1\5\2\4\eeee
\4\2\5\3\4\2\5\3\4\2\5\3\4\2\5\3\eeee
\3\6\1\4\3\6\1\4\3\6\1\4\3\6\1\4\eeee
\5\1\6\2\5\1\6\2\5\1\6\2\5\1\6\2\eeee
\2\4\3\5\2\4\3\5\2\4\3\5\2\4\3\5\eeee
\6\3\4\1\6\3\4\1\6\3\4\1\6\3\4\1\eeee
\1\5\2\6\1\5\2\6\1\5\2\6\1\5\2\6\eeee
\4\2\5\3\4\2\5\3\4\2\5\3\4\2\5\3\eeee
\3\6\1\4\3\6\1\4\3\6\1\4\3\6\1\4\eeee
\5\1\6\2\5\1\6\2\5\1\6\2\5\1\6\2\eeee
\2\4\3\5\2\4\3\5\2\4\3\5\2\4\3\5\eeee
\6\3\4\1\6\3\4\1\6\3\4\1\6\3\4\1\eeee
} 
\baaa
\phantom{0-0}\#
\eaaa
\mbox{}\phantom{\bbbb
0&0&0&0&0&0\\
\ebbb}
\parbox{7cm}{ 
\1\6\1\6\1\6\1\6\1\6\1\6\1\6\1\6\eeee
\4\3\4\3\4\3\4\3\4\3\4\3\4\3\4\3\eeee
\2\5\2\5\2\5\2\5\2\5\2\5\2\5\2\5\eeee
\6\1\6\1\6\1\6\1\6\1\6\1\6\1\6\1\eeee
\3\4\3\4\3\4\3\4\3\4\3\4\3\4\3\4\eeee
\5\2\5\2\5\2\5\2\5\2\5\2\5\2\5\2\eeee
\1\6\1\6\1\6\1\6\1\6\1\6\1\6\1\6\eeee
\4\3\4\3\4\3\4\3\4\3\4\3\4\3\4\3\eeee
\2\5\2\5\2\5\2\5\2\5\2\5\2\5\2\5\eeee
\6\1\6\1\6\1\6\1\6\1\6\1\6\1\6\1\eeee
\3\4\3\4\3\4\3\4\3\4\3\4\3\4\3\4\eeee
\5\2\5\2\5\2\5\2\5\2\5\2\5\2\5\2\eeee
} 

\baaa
6-59
\eaaa
\bbbb
0&0&0&1&1&2\\
0&0&1&0&1&2\\
0&1&0&2&1&0\\
1&0&1&1&0&1\\
2&1&1&0&0&0\\
2&1&0&1&0&0\\
\ebbb
\parbox{7cm}{ 
\1\5\1\6\2\6\1\5\1\6\2\6\1\5\1\6\eeee
\4\3\4\4\3\4\4\3\4\4\3\4\4\3\4\4\eeee
\6\2\6\1\5\1\6\2\6\1\5\1\6\2\6\1\eeee
\1\5\1\6\2\6\1\5\1\6\2\6\1\5\1\6\eeee
\4\3\4\4\3\4\4\3\4\4\3\4\4\3\4\4\eeee
\6\2\6\1\5\1\6\2\6\1\5\1\6\2\6\1\eeee
\1\5\1\6\2\6\1\5\1\6\2\6\1\5\1\6\eeee
\4\3\4\4\3\4\4\3\4\4\3\4\4\3\4\4\eeee
\6\2\6\1\5\1\6\2\6\1\5\1\6\2\6\1\eeee
\1\5\1\6\2\6\1\5\1\6\2\6\1\5\1\6\eeee
\4\3\4\4\3\4\4\3\4\4\3\4\4\3\4\4\eeee
\6\2\6\1\5\1\6\2\6\1\5\1\6\2\6\1\eeee
} 

\baaa
6-60
\eaaa
\bbbb
0&0&0&1&1&2\\
0&0&1&0&2&1\\
0&1&0&2&0&1\\
1&0&2&0&1&0\\
1&2&0&1&0&0\\
2&1&1&0&0&0\\
\ebbb
\parbox{7cm}{ 
\1\6\1\6\1\6\1\6\1\6\1\6\1\6\1\6\eeee
\4\3\4\3\4\3\4\3\4\3\4\3\4\3\4\3\eeee
\5\2\5\2\5\2\5\2\5\2\5\2\5\2\5\2\eeee
\1\6\1\6\1\6\1\6\1\6\1\6\1\6\1\6\eeee
\4\3\4\3\4\3\4\3\4\3\4\3\4\3\4\3\eeee
\5\2\5\2\5\2\5\2\5\2\5\2\5\2\5\2\eeee
\1\6\1\6\1\6\1\6\1\6\1\6\1\6\1\6\eeee
\4\3\4\3\4\3\4\3\4\3\4\3\4\3\4\3\eeee
\5\2\5\2\5\2\5\2\5\2\5\2\5\2\5\2\eeee
\1\6\1\6\1\6\1\6\1\6\1\6\1\6\1\6\eeee
\4\3\4\3\4\3\4\3\4\3\4\3\4\3\4\3\eeee
\5\2\5\2\5\2\5\2\5\2\5\2\5\2\5\2\eeee
} 

\baaa
6-61
\eaaa
\bbbb
0&0&0&1&1&2\\
0&0&1&0&2&1\\
0&1&1&1&0&1\\
1&0&1&1&1&0\\
1&2&0&1&0&0\\
2&1&1&0&0&0\\
\ebbb
\parbox{7cm}{ 
\1\5\2\6\1\5\2\6\1\5\2\6\1\5\2\6\eeee
\4\4\3\3\4\4\3\3\4\4\3\3\4\4\3\3\eeee
\5\1\6\2\5\1\6\2\5\1\6\2\5\1\6\2\eeee
\2\6\1\5\2\6\1\5\2\6\1\5\2\6\1\5\eeee
\3\3\4\4\3\3\4\4\3\3\4\4\3\3\4\4\eeee
\6\2\5\1\6\2\5\1\6\2\5\1\6\2\5\1\eeee
\1\5\2\6\1\5\2\6\1\5\2\6\1\5\2\6\eeee
\4\4\3\3\4\4\3\3\4\4\3\3\4\4\3\3\eeee
\5\1\6\2\5\1\6\2\5\1\6\2\5\1\6\2\eeee
\2\6\1\5\2\6\1\5\2\6\1\5\2\6\1\5\eeee
\3\3\4\4\3\3\4\4\3\3\4\4\3\3\4\4\eeee
\6\2\5\1\6\2\5\1\6\2\5\1\6\2\5\1\eeee
} 

\baaa
6-62
\eaaa
\bbbb
0&0&0&1&1&2\\
0&0&1&0&2&1\\
0&1&1&1&1&0\\
1&0&1&1&0&1\\
1&2&1&0&0&0\\
2&1&0&1&0&0\\
\ebbb
\parbox{7cm}{ 
\1\5\2\6\1\5\2\6\1\5\2\6\1\5\2\6\eeee
\4\3\3\4\4\3\3\4\4\3\3\4\4\3\3\4\eeee
\6\2\5\1\6\2\5\1\6\2\5\1\6\2\5\1\eeee
\1\5\2\6\1\5\2\6\1\5\2\6\1\5\2\6\eeee
\4\3\3\4\4\3\3\4\4\3\3\4\4\3\3\4\eeee
\6\2\5\1\6\2\5\1\6\2\5\1\6\2\5\1\eeee
\1\5\2\6\1\5\2\6\1\5\2\6\1\5\2\6\eeee
\4\3\3\4\4\3\3\4\4\3\3\4\4\3\3\4\eeee
\6\2\5\1\6\2\5\1\6\2\5\1\6\2\5\1\eeee
\1\5\2\6\1\5\2\6\1\5\2\6\1\5\2\6\eeee
\4\3\3\4\4\3\3\4\4\3\3\4\4\3\3\4\eeee
\6\2\5\1\6\2\5\1\6\2\5\1\6\2\5\1\eeee
} 

\baaa
6-63
\eaaa
\bbbb
0&0&0&1&1&2\\
0&0&1&0&2&1\\
0&1&2&0&1&0\\
1&0&0&2&0&1\\
1&2&1&0&0&0\\
2&1&0&1&0&0\\
\ebbb
\parbox{7cm}{ 
\1\6\1\6\1\6\1\6\1\6\1\6\1\6\1\6\eeee
\4\4\4\4\4\4\4\4\4\4\4\4\4\4\4\4\eeee
\6\1\6\1\6\1\6\1\6\1\6\1\6\1\6\1\eeee
\2\5\2\5\2\5\2\5\2\5\2\5\2\5\2\5\eeee
\3\3\3\3\3\3\3\3\3\3\3\3\3\3\3\3\eeee
\5\2\5\2\5\2\5\2\5\2\5\2\5\2\5\2\eeee
\1\6\1\6\1\6\1\6\1\6\1\6\1\6\1\6\eeee
\4\4\4\4\4\4\4\4\4\4\4\4\4\4\4\4\eeee
\6\1\6\1\6\1\6\1\6\1\6\1\6\1\6\1\eeee
\2\5\2\5\2\5\2\5\2\5\2\5\2\5\2\5\eeee
\3\3\3\3\3\3\3\3\3\3\3\3\3\3\3\3\eeee
\5\2\5\2\5\2\5\2\5\2\5\2\5\2\5\2\eeee
} 

\baaa
6-64
\eaaa
\bbbb
0&0&0&1&1&2\\
0&0&1&1&1&1\\
0&1&0&1&1&1\\
1&1&1&0&1&0\\
1&1&1&1&0&0\\
2&1&1&0&0&0\\
\ebbb
\parbox{7cm}{ 
\1\6\1\6\1\6\1\6\1\6\1\6\1\6\1\6\eeee
\4\2\5\3\4\2\5\3\4\2\5\3\4\2\5\3\eeee
\5\3\4\2\5\3\4\2\5\3\4\2\5\3\4\2\eeee
\1\6\1\6\1\6\1\6\1\6\1\6\1\6\1\6\eeee
\4\2\5\3\4\2\5\3\4\2\5\3\4\2\5\3\eeee
\5\3\4\2\5\3\4\2\5\3\4\2\5\3\4\2\eeee
\1\6\1\6\1\6\1\6\1\6\1\6\1\6\1\6\eeee
\4\2\5\3\4\2\5\3\4\2\5\3\4\2\5\3\eeee
\5\3\4\2\5\3\4\2\5\3\4\2\5\3\4\2\eeee
\1\6\1\6\1\6\1\6\1\6\1\6\1\6\1\6\eeee
\4\2\5\3\4\2\5\3\4\2\5\3\4\2\5\3\eeee
\5\3\4\2\5\3\4\2\5\3\4\2\5\3\4\2\eeee
} 

\baaa
6-65
\eaaa
\bbbb
0&0&0&1&1&2\\
0&0&1&1&1&1\\
0&1&0&1&1&1\\
1&1&1&1&0&0\\
1&1&1&0&1&0\\
2&1&1&0&0&0\\
\ebbb
\parbox{7cm}{ 
\1\6\1\6\1\6\1\6\1\6\1\6\1\6\1\6\eeee
\4\2\5\3\4\2\5\3\4\2\5\3\4\2\5\3\eeee
\4\3\5\2\4\3\5\2\4\3\5\2\4\3\5\2\eeee
\1\6\1\6\1\6\1\6\1\6\1\6\1\6\1\6\eeee
\5\2\4\3\5\2\4\3\5\2\4\3\5\2\4\3\eeee
\5\3\4\2\5\3\4\2\5\3\4\2\5\3\4\2\eeee
\1\6\1\6\1\6\1\6\1\6\1\6\1\6\1\6\eeee
\4\2\5\3\4\2\5\3\4\2\5\3\4\2\5\3\eeee
\4\3\5\2\4\3\5\2\4\3\5\2\4\3\5\2\eeee
\1\6\1\6\1\6\1\6\1\6\1\6\1\6\1\6\eeee
\5\2\4\3\5\2\4\3\5\2\4\3\5\2\4\3\eeee
\5\3\4\2\5\3\4\2\5\3\4\2\5\3\4\2\eeee
} 

\baaa
6-66
\eaaa
\bbbb
0&0&0&1&1&2\\
0&0&1&1&2&0\\
0&1&1&0&1&1\\
1&1&0&1&0&1\\
1&2&1&0&0&0\\
2&0&1&1&0&0\\
\ebbb
\parbox{7cm}{ 
\1\5\2\4\4\2\5\1\6\3\3\6\1\5\2\4\eeee
\4\2\5\1\6\3\3\6\1\5\2\4\4\2\5\1\eeee
\6\3\3\6\1\5\2\4\4\2\5\1\6\3\3\6\eeee
\1\5\2\4\4\2\5\1\6\3\3\6\1\5\2\4\eeee
\4\2\5\1\6\3\3\6\1\5\2\4\4\2\5\1\eeee
\6\3\3\6\1\5\2\4\4\2\5\1\6\3\3\6\eeee
\1\5\2\4\4\2\5\1\6\3\3\6\1\5\2\4\eeee
\4\2\5\1\6\3\3\6\1\5\2\4\4\2\5\1\eeee
\6\3\3\6\1\5\2\4\4\2\5\1\6\3\3\6\eeee
\1\5\2\4\4\2\5\1\6\3\3\6\1\5\2\4\eeee
\4\2\5\1\6\3\3\6\1\5\2\4\4\2\5\1\eeee
\6\3\3\6\1\5\2\4\4\2\5\1\6\3\3\6\eeee
} 

\baaa
6-67
\eaaa
\bbbb
0&0&0&1&1&2\\
0&0&2&1&1&0\\
0&1&1&0&1&1\\
1&1&0&0&0&2\\
1&1&2&0&0&0\\
1&0&1&1&0&1\\
\ebbb
\parbox{7cm}{ 
\1\5\3\6\1\5\3\6\1\5\3\6\1\5\3\6\eeee
\4\2\3\6\4\2\3\6\4\2\3\6\4\2\3\6\eeee
\6\3\5\1\6\3\5\1\6\3\5\1\6\3\5\1\eeee
\6\3\2\4\6\3\2\4\6\3\2\4\6\3\2\4\eeee
\1\5\3\6\1\5\3\6\1\5\3\6\1\5\3\6\eeee
\4\2\3\6\4\2\3\6\4\2\3\6\4\2\3\6\eeee
\6\3\5\1\6\3\5\1\6\3\5\1\6\3\5\1\eeee
\6\3\2\4\6\3\2\4\6\3\2\4\6\3\2\4\eeee
\1\5\3\6\1\5\3\6\1\5\3\6\1\5\3\6\eeee
\4\2\3\6\4\2\3\6\4\2\3\6\4\2\3\6\eeee
\6\3\5\1\6\3\5\1\6\3\5\1\6\3\5\1\eeee
\6\3\2\4\6\3\2\4\6\3\2\4\6\3\2\4\eeee
} 
\baaa
\phantom{0-0}\#
\eaaa
\mbox{}\phantom{\bbbb
0&0&0&0&0&0\\
\ebbb}
\parbox{7cm}{ 
\1\6\4\6\1\6\4\6\1\6\4\6\1\6\4\6\eeee
\4\6\1\6\4\6\1\6\4\6\1\6\4\6\1\6\eeee
\2\3\5\3\2\3\5\3\2\3\5\3\2\3\5\3\eeee
\5\3\2\3\5\3\2\3\5\3\2\3\5\3\2\3\eeee
\1\6\4\6\1\6\4\6\1\6\4\6\1\6\4\6\eeee
\4\6\1\6\4\6\1\6\4\6\1\6\4\6\1\6\eeee
\2\3\5\3\2\3\5\3\2\3\5\3\2\3\5\3\eeee
\5\3\2\3\5\3\2\3\5\3\2\3\5\3\2\3\eeee
\1\6\4\6\1\6\4\6\1\6\4\6\1\6\4\6\eeee
\4\6\1\6\4\6\1\6\4\6\1\6\4\6\1\6\eeee
\2\3\5\3\2\3\5\3\2\3\5\3\2\3\5\3\eeee
\5\3\2\3\5\3\2\3\5\3\2\3\5\3\2\3\eeee
} 

\baaa
6-68
\eaaa
\bbbb
0&0&0&1&1&2\\
0&0&2&1&1&0\\
0&2&1&0&1&0\\
1&1&0&0&1&1\\
1&1&1&1&0&0\\
2&0&0&1&0&1\\
\ebbb
\parbox{7cm}{ 
\1\5\3\2\4\6\1\5\3\2\4\6\1\5\3\2\eeee
\4\2\3\5\1\6\4\2\3\5\1\6\4\2\3\5\eeee
\5\3\2\4\6\1\5\3\2\4\6\1\5\3\2\4\eeee
\2\3\5\1\6\4\2\3\5\1\6\4\2\3\5\1\eeee
\3\2\4\6\1\5\3\2\4\6\1\5\3\2\4\6\eeee
\3\5\1\6\4\2\3\5\1\6\4\2\3\5\1\6\eeee
\2\4\6\1\5\3\2\4\6\1\5\3\2\4\6\1\eeee
\5\1\6\4\2\3\5\1\6\4\2\3\5\1\6\4\eeee
\4\6\1\5\3\2\4\6\1\5\3\2\4\6\1\5\eeee
\1\6\4\2\3\5\1\6\4\2\3\5\1\6\4\2\eeee
\6\1\5\3\2\4\6\1\5\3\2\4\6\1\5\3\eeee
\6\4\2\3\5\1\6\4\2\3\5\1\6\4\2\3\eeee
} 

\baaa
6-69
\eaaa
\bbbb
0&0&0&1&1&2\\
0&1&0&0&2&1\\
0&0&1&2&0&1\\
1&0&2&1&0&0\\
1&2&0&0&1&0\\
2&1&1&0&0&0\\
\ebbb
\parbox{7cm}{ 
\1\6\1\6\1\6\1\6\1\6\1\6\1\6\1\6\eeee
\4\3\4\3\4\3\4\3\4\3\4\3\4\3\4\3\eeee
\4\3\4\3\4\3\4\3\4\3\4\3\4\3\4\3\eeee
\1\6\1\6\1\6\1\6\1\6\1\6\1\6\1\6\eeee
\5\2\5\2\5\2\5\2\5\2\5\2\5\2\5\2\eeee
\5\2\5\2\5\2\5\2\5\2\5\2\5\2\5\2\eeee
\1\6\1\6\1\6\1\6\1\6\1\6\1\6\1\6\eeee
\4\3\4\3\4\3\4\3\4\3\4\3\4\3\4\3\eeee
\4\3\4\3\4\3\4\3\4\3\4\3\4\3\4\3\eeee
\1\6\1\6\1\6\1\6\1\6\1\6\1\6\1\6\eeee
\5\2\5\2\5\2\5\2\5\2\5\2\5\2\5\2\eeee
\5\2\5\2\5\2\5\2\5\2\5\2\5\2\5\2\eeee
} 

\baaa
6-70
\eaaa
\bbbb
0&0&0&1&1&2\\
0&1&0&1&1&1\\
0&0&1&1&1&1\\
1&1&1&1&0&0\\
1&1&1&0&1&0\\
2&1&1&0&0&0\\
\ebbb
\parbox{7cm}{ 
\1\6\1\6\1\6\1\6\1\6\1\6\1\6\1\6\eeee
\4\2\5\3\4\2\5\3\4\2\5\3\4\2\5\3\eeee
\4\2\5\3\4\2\5\3\4\2\5\3\4\2\5\3\eeee
\1\6\1\6\1\6\1\6\1\6\1\6\1\6\1\6\eeee
\5\3\4\2\5\3\4\2\5\3\4\2\5\3\4\2\eeee
\5\3\4\2\5\3\4\2\5\3\4\2\5\3\4\2\eeee
\1\6\1\6\1\6\1\6\1\6\1\6\1\6\1\6\eeee
\4\2\5\3\4\2\5\3\4\2\5\3\4\2\5\3\eeee
\4\2\5\3\4\2\5\3\4\2\5\3\4\2\5\3\eeee
\1\6\1\6\1\6\1\6\1\6\1\6\1\6\1\6\eeee
\5\3\4\2\5\3\4\2\5\3\4\2\5\3\4\2\eeee
\5\3\4\2\5\3\4\2\5\3\4\2\5\3\4\2\eeee
} 

\baaa
6-71
\eaaa
\bbbb
0&0&0&1&1&2\\
0&1&1&0&0&2\\
0&1&1&0&0&2\\
2&0&0&1&1&0\\
2&0&0&1&1&0\\
2&1&1&0&0&0\\
\ebbb
\parbox{7cm}{ 
\1\5\5\1\6\3\3\6\1\5\5\1\6\3\3\6\eeee
\4\4\1\6\2\2\6\1\4\4\1\6\2\2\6\1\eeee
\5\1\6\3\3\6\1\5\5\1\6\3\3\6\1\5\eeee
\1\6\2\2\6\1\4\4\1\6\2\2\6\1\4\4\eeee
\6\3\3\6\1\5\5\1\6\3\3\6\1\5\5\1\eeee
\2\2\6\1\4\4\1\6\2\2\6\1\4\4\1\6\eeee
\3\6\1\5\5\1\6\3\3\6\1\5\5\1\6\3\eeee
\6\1\4\4\1\6\2\2\6\1\4\4\1\6\2\2\eeee
\1\5\5\1\6\3\3\6\1\5\5\1\6\3\3\6\eeee
\4\4\1\6\2\2\6\1\4\4\1\6\2\2\6\1\eeee
\5\1\6\3\3\6\1\5\5\1\6\3\3\6\1\5\eeee
\1\6\2\2\6\1\4\4\1\6\2\2\6\1\4\4\eeee
} 

\baaa
6-72
\eaaa
\bbbb
0&0&0&1&1&2\\
0&1&1&0&1&1\\
0&1&1&1&0&1\\
1&0&2&0&1&0\\
1&2&0&1&0&0\\
1&1&1&0&0&1\\
\ebbb
\parbox{7cm}{ 
\1\6\6\1\6\6\1\6\6\1\6\6\1\6\6\1\eeee
\4\3\3\4\3\3\4\3\3\4\3\3\4\3\3\4\eeee
\5\2\2\5\2\2\5\2\2\5\2\2\5\2\2\5\eeee
\1\6\6\1\6\6\1\6\6\1\6\6\1\6\6\1\eeee
\4\3\3\4\3\3\4\3\3\4\3\3\4\3\3\4\eeee
\5\2\2\5\2\2\5\2\2\5\2\2\5\2\2\5\eeee
\1\6\6\1\6\6\1\6\6\1\6\6\1\6\6\1\eeee
\4\3\3\4\3\3\4\3\3\4\3\3\4\3\3\4\eeee
\5\2\2\5\2\2\5\2\2\5\2\2\5\2\2\5\eeee
\1\6\6\1\6\6\1\6\6\1\6\6\1\6\6\1\eeee
\4\3\3\4\3\3\4\3\3\4\3\3\4\3\3\4\eeee
\5\2\2\5\2\2\5\2\2\5\2\2\5\2\2\5\eeee
} 

\baaa
6-73
\eaaa
\bbbb
0&0&0&1&1&2\\
0&1&1&0&1&1\\
0&1&1&1&0&1\\
1&0&2&1&0&0\\
1&2&0&0&1&0\\
1&1&1&0&0&1\\
\ebbb
\parbox{7cm}{ 
\1\6\6\1\6\6\1\6\6\1\6\6\1\6\6\1\eeee
\4\3\2\5\2\3\4\3\2\5\2\3\4\3\2\5\eeee
\4\3\2\5\2\3\4\3\2\5\2\3\4\3\2\5\eeee
\1\6\6\1\6\6\1\6\6\1\6\6\1\6\6\1\eeee
\5\2\3\4\3\2\5\2\3\4\3\2\5\2\3\4\eeee
\5\2\3\4\3\2\5\2\3\4\3\2\5\2\3\4\eeee
\1\6\6\1\6\6\1\6\6\1\6\6\1\6\6\1\eeee
\4\3\2\5\2\3\4\3\2\5\2\3\4\3\2\5\eeee
\4\3\2\5\2\3\4\3\2\5\2\3\4\3\2\5\eeee
\1\6\6\1\6\6\1\6\6\1\6\6\1\6\6\1\eeee
\5\2\3\4\3\2\5\2\3\4\3\2\5\2\3\4\eeee
\5\2\3\4\3\2\5\2\3\4\3\2\5\2\3\4\eeee
} 

\baaa
6-74
\eaaa
\bbbb
0&0&0&1&1&2\\
0&1&1&0&1&1\\
0&1&1&1&1&0\\
1&0&1&1&0&1\\
1&1&1&0&1&0\\
2&1&0&1&0&0\\
\ebbb
\parbox{7cm}{ 
\1\5\5\1\6\2\2\6\1\5\5\1\6\2\2\6\eeee
\4\3\3\4\4\3\3\4\4\3\3\4\4\3\3\4\eeee
\6\2\2\6\1\5\5\1\6\2\2\6\1\5\5\1\eeee
\1\5\5\1\6\2\2\6\1\5\5\1\6\2\2\6\eeee
\4\3\3\4\4\3\3\4\4\3\3\4\4\3\3\4\eeee
\6\2\2\6\1\5\5\1\6\2\2\6\1\5\5\1\eeee
\1\5\5\1\6\2\2\6\1\5\5\1\6\2\2\6\eeee
\4\3\3\4\4\3\3\4\4\3\3\4\4\3\3\4\eeee
\6\2\2\6\1\5\5\1\6\2\2\6\1\5\5\1\eeee
\1\5\5\1\6\2\2\6\1\5\5\1\6\2\2\6\eeee
\4\3\3\4\4\3\3\4\4\3\3\4\4\3\3\4\eeee
\6\2\2\6\1\5\5\1\6\2\2\6\1\5\5\1\eeee
} 

\baaa
6-75
\eaaa
\bbbb
0&0&0&1&1&2\\
0&1&1&0&1&1\\
0&2&1&1&0&0\\
1&0&1&0&0&2\\
1&2&0&0&1&0\\
1&1&0&1&0&1\\
\ebbb
\parbox{7cm}{ 
\1\6\4\6\1\6\4\6\1\6\4\6\1\6\4\6\eeee
\4\6\1\6\4\6\1\6\4\6\1\6\4\6\1\6\eeee
\3\2\5\2\3\2\5\2\3\2\5\2\3\2\5\2\eeee
\3\2\5\2\3\2\5\2\3\2\5\2\3\2\5\2\eeee
\4\6\1\6\4\6\1\6\4\6\1\6\4\6\1\6\eeee
\1\6\4\6\1\6\4\6\1\6\4\6\1\6\4\6\eeee
\5\2\3\2\5\2\3\2\5\2\3\2\5\2\3\2\eeee
\5\2\3\2\5\2\3\2\5\2\3\2\5\2\3\2\eeee
\1\6\4\6\1\6\4\6\1\6\4\6\1\6\4\6\eeee
\4\6\1\6\4\6\1\6\4\6\1\6\4\6\1\6\eeee
\3\2\5\2\3\2\5\2\3\2\5\2\3\2\5\2\eeee
\3\2\5\2\3\2\5\2\3\2\5\2\3\2\5\2\eeee
} 

\baaa
6-76
\eaaa
\bbbb
0&0&0&1&1&2\\
0&2&0&0&1&1\\
0&0&2&1&0&1\\
1&0&2&1&0&0\\
1&2&0&0&1&0\\
1&1&1&0&0&1\\
\ebbb
\parbox{7cm}{ 
\1\6\6\1\6\6\1\6\6\1\6\6\1\6\6\1\eeee
\4\3\3\4\3\3\4\3\3\4\3\3\4\3\3\4\eeee
\4\3\3\4\3\3\4\3\3\4\3\3\4\3\3\4\eeee
\1\6\6\1\6\6\1\6\6\1\6\6\1\6\6\1\eeee
\5\2\2\5\2\2\5\2\2\5\2\2\5\2\2\5\eeee
\5\2\2\5\2\2\5\2\2\5\2\2\5\2\2\5\eeee
\1\6\6\1\6\6\1\6\6\1\6\6\1\6\6\1\eeee
\4\3\3\4\3\3\4\3\3\4\3\3\4\3\3\4\eeee
\4\3\3\4\3\3\4\3\3\4\3\3\4\3\3\4\eeee
\1\6\6\1\6\6\1\6\6\1\6\6\1\6\6\1\eeee
\5\2\2\5\2\2\5\2\2\5\2\2\5\2\2\5\eeee
\5\2\2\5\2\2\5\2\2\5\2\2\5\2\2\5\eeee
} 

\baaa
6-77
\eaaa
\bbbb
0&0&1&1&1&1\\
0&0&1&1&1&1\\
1&1&0&0&1&1\\
1&1&0&0&1&1\\
1&1&1&1&0&0\\
1&1&1&1&0&0\\
\ebbb
\parbox{7cm}{ 
\1\4\6\1\3\6\2\3\5\2\4\5\1\4\6\1\eeee
\3\5\2\4\5\1\4\6\1\3\6\2\3\5\2\4\eeee
\6\1\3\6\2\3\5\2\4\5\1\4\6\1\3\6\eeee
\2\4\5\1\4\6\1\3\6\2\3\5\2\4\5\1\eeee
\3\6\2\3\5\2\4\5\1\4\6\1\3\6\2\3\eeee
\5\1\4\6\1\3\6\2\3\5\2\4\5\1\4\6\eeee
\2\3\5\2\4\5\1\4\6\1\3\6\2\3\5\2\eeee
\4\6\1\3\6\2\3\5\2\4\5\1\4\6\1\3\eeee
\5\2\4\5\1\4\6\1\3\6\2\3\5\2\4\5\eeee
\1\3\6\2\3\5\2\4\5\1\4\6\1\3\6\2\eeee
\4\5\1\4\6\1\3\6\2\3\5\2\4\5\1\4\eeee
\6\2\3\5\2\4\5\1\4\6\1\3\6\2\3\5\eeee
} 

\baaa
6-78
\eaaa
\bbbb
0&0&1&1&1&1\\
0&0&1&1&1&1\\
1&1&0&0&1&1\\
1&1&0&1&0&1\\
1&1&1&0&1&0\\
1&1&1&1&0&0\\
\ebbb
\parbox{7cm}{ 
\1\4\4\1\6\4\2\3\6\2\5\3\1\5\5\1\eeee
\3\6\2\5\3\1\5\5\1\3\5\2\6\3\2\4\eeee
\5\1\3\5\2\6\3\2\4\6\1\4\4\1\6\4\eeee
\2\4\6\1\4\4\1\6\4\2\3\6\2\5\3\1\eeee
\6\4\2\3\6\2\5\3\1\5\5\1\3\5\2\6\eeee
\3\1\5\5\1\3\5\2\6\3\2\4\6\1\4\4\eeee
\2\6\3\2\4\6\1\4\4\1\6\4\2\3\6\2\eeee
\4\4\1\6\4\2\3\6\2\5\3\1\5\5\1\3\eeee
\6\2\5\3\1\5\5\1\3\5\2\6\3\2\4\6\eeee
\1\3\5\2\6\3\2\4\6\1\4\4\1\6\4\2\eeee
\4\6\1\4\4\1\6\4\2\3\6\2\5\3\1\5\eeee
\4\2\3\6\2\5\3\1\5\5\1\3\5\2\6\3\eeee
} 

\baaa
6-79
\eaaa
\bbbb
0&0&1&1&1&1\\
0&0&1&1&1&1\\
1&1&1&0&0&1\\
1&1&0&1&1&0\\
1&1&0&1&1&0\\
1&1&1&0&0&1\\
\ebbb
\parbox{7cm}{ 
\1\4\4\1\3\3\1\4\4\1\3\3\1\4\4\1\eeee
\3\2\5\5\2\6\6\2\5\5\2\6\6\2\5\5\eeee
\3\6\1\4\4\1\3\3\1\4\4\1\3\3\1\4\eeee
\1\6\3\2\5\5\2\6\6\2\5\5\2\6\6\2\eeee
\4\2\3\6\1\4\4\1\3\3\1\4\4\1\3\3\eeee
\4\5\1\6\3\2\5\5\2\6\6\2\5\5\2\6\eeee
\1\5\4\2\3\6\1\4\4\1\3\3\1\4\4\1\eeee
\3\2\4\5\1\6\3\2\5\5\2\6\6\2\5\5\eeee
\3\6\1\5\4\2\3\6\1\4\4\1\3\3\1\4\eeee
\1\6\3\2\4\5\1\6\3\2\5\5\2\6\6\2\eeee
\4\2\3\6\1\5\4\2\3\6\1\4\4\1\3\3\eeee
\4\5\1\6\3\2\4\5\1\6\3\2\5\5\2\6\eeee
} 

\baaa
6-80
\eaaa
\bbbb
0&0&1&1&1&1\\
0&1&0&1&1&1\\
1&0&0&1&1&1\\
1&1&1&0&0&1\\
1&1&1&0&1&0\\
1&1&1&1&0&0\\
\ebbb
\parbox{7cm}{ 
\1\4\2\6\3\5\1\4\2\6\3\5\1\4\2\6\eeee
\3\6\2\4\1\5\3\6\2\4\1\5\3\6\2\4\eeee
\4\1\5\3\6\2\4\1\5\3\6\2\4\1\5\3\eeee
\6\3\5\1\4\2\6\3\5\1\4\2\6\3\5\1\eeee
\1\4\2\6\3\5\1\4\2\6\3\5\1\4\2\6\eeee
\3\6\2\4\1\5\3\6\2\4\1\5\3\6\2\4\eeee
\4\1\5\3\6\2\4\1\5\3\6\2\4\1\5\3\eeee
\6\3\5\1\4\2\6\3\5\1\4\2\6\3\5\1\eeee
\1\4\2\6\3\5\1\4\2\6\3\5\1\4\2\6\eeee
\3\6\2\4\1\5\3\6\2\4\1\5\3\6\2\4\eeee
\4\1\5\3\6\2\4\1\5\3\6\2\4\1\5\3\eeee
\6\3\5\1\4\2\6\3\5\1\4\2\6\3\5\1\eeee
} 
\baaa
\phantom{0-0}\#
\eaaa
\mbox{}\phantom{\bbbb
0&0&0&0&0&0\\
\ebbb}
\parbox{7cm}{ 
\1\4\3\6\2\5\1\4\3\6\2\5\1\4\3\6\eeee
\3\6\1\4\2\5\3\6\1\4\2\5\3\6\1\4\eeee
\4\2\5\3\6\1\4\2\5\3\6\1\4\2\5\3\eeee
\6\2\5\1\4\3\6\2\5\1\4\3\6\2\5\1\eeee
\1\4\3\6\2\5\1\4\3\6\2\5\1\4\3\6\eeee
\3\6\1\4\2\5\3\6\1\4\2\5\3\6\1\4\eeee
\4\2\5\3\6\1\4\2\5\3\6\1\4\2\5\3\eeee
\6\2\5\1\4\3\6\2\5\1\4\3\6\2\5\1\eeee
\1\4\3\6\2\5\1\4\3\6\2\5\1\4\3\6\eeee
\3\6\1\4\2\5\3\6\1\4\2\5\3\6\1\4\eeee
\4\2\5\3\6\1\4\2\5\3\6\1\4\2\5\3\eeee
\6\2\5\1\4\3\6\2\5\1\4\3\6\2\5\1\eeee
} 

\baaa
6-81
\eaaa
\bbbb
0&0&1&1&1&1\\
0&1&0&1&1&1\\
1&0&0&1&1&1\\
1&1&1&1&0&0\\
1&1&1&0&1&0\\
1&1&1&0&0&1\\
\ebbb
\parbox{7cm}{ 
\1\4\2\5\3\6\1\4\2\5\3\6\1\4\2\5\eeee
\3\4\2\5\1\6\3\4\2\5\1\6\3\4\2\5\eeee
\5\1\6\3\4\2\5\1\6\3\4\2\5\1\6\3\eeee
\5\3\6\1\4\2\5\3\6\1\4\2\5\3\6\1\eeee
\1\4\2\5\3\6\1\4\2\5\3\6\1\4\2\5\eeee
\3\4\2\5\1\6\3\4\2\5\1\6\3\4\2\5\eeee
\5\1\6\3\4\2\5\1\6\3\4\2\5\1\6\3\eeee
\5\3\6\1\4\2\5\3\6\1\4\2\5\3\6\1\eeee
\1\4\2\5\3\6\1\4\2\5\3\6\1\4\2\5\eeee
\3\4\2\5\1\6\3\4\2\5\1\6\3\4\2\5\eeee
\5\1\6\3\4\2\5\1\6\3\4\2\5\1\6\3\eeee
\5\3\6\1\4\2\5\3\6\1\4\2\5\3\6\1\eeee
} 

\baaa
6-82
\eaaa
\bbbb
0&0&1&1&1&1\\
0&1&0&1&1&1\\
1&0&1&0&1&1\\
1&1&0&1&0&1\\
1&1&1&0&1&0\\
1&1&1&1&0&0\\
\ebbb
\parbox{7cm}{ 
\1\5\3\3\5\1\6\2\4\4\2\6\1\5\3\3\eeee
\3\5\1\6\2\4\4\2\6\1\5\3\3\5\1\6\eeee
\6\2\4\4\2\6\1\5\3\3\5\1\6\2\4\4\eeee
\4\2\6\1\5\3\3\5\1\6\2\4\4\2\6\1\eeee
\1\5\3\3\5\1\6\2\4\4\2\6\1\5\3\3\eeee
\3\5\1\6\2\4\4\2\6\1\5\3\3\5\1\6\eeee
\6\2\4\4\2\6\1\5\3\3\5\1\6\2\4\4\eeee
\4\2\6\1\5\3\3\5\1\6\2\4\4\2\6\1\eeee
\1\5\3\3\5\1\6\2\4\4\2\6\1\5\3\3\eeee
\3\5\1\6\2\4\4\2\6\1\5\3\3\5\1\6\eeee
\6\2\4\4\2\6\1\5\3\3\5\1\6\2\4\4\eeee
\4\2\6\1\5\3\3\5\1\6\2\4\4\2\6\1\eeee
} 

\baaa
6-83
\eaaa
\bbbb
1&0&0&0&1&2\\
0&1&0&1&0&2\\
0&0&1&1&1&1\\
0&1&2&1&0&0\\
1&0&2&0&1&0\\
1&1&1&0&0&1\\
\ebbb
\parbox{7cm}{ 
\1\5\3\6\1\5\3\6\1\5\3\6\1\5\3\6\eeee
\1\5\3\6\1\5\3\6\1\5\3\6\1\5\3\6\eeee
\6\3\4\2\6\3\4\2\6\3\4\2\6\3\4\2\eeee
\6\3\4\2\6\3\4\2\6\3\4\2\6\3\4\2\eeee
\1\5\3\6\1\5\3\6\1\5\3\6\1\5\3\6\eeee
\1\5\3\6\1\5\3\6\1\5\3\6\1\5\3\6\eeee
\6\3\4\2\6\3\4\2\6\3\4\2\6\3\4\2\eeee
\6\3\4\2\6\3\4\2\6\3\4\2\6\3\4\2\eeee
\1\5\3\6\1\5\3\6\1\5\3\6\1\5\3\6\eeee
\1\5\3\6\1\5\3\6\1\5\3\6\1\5\3\6\eeee
\6\3\4\2\6\3\4\2\6\3\4\2\6\3\4\2\eeee
\6\3\4\2\6\3\4\2\6\3\4\2\6\3\4\2\eeee
} 
\baaa
\phantom{0-0}\#
\eaaa
\mbox{}\phantom{\bbbb
0&0&0&0&0&0\\
\ebbb}
\parbox{7cm}{ 
\1\6\2\6\1\6\2\6\1\6\2\6\1\6\2\6\eeee
\1\6\2\6\1\6\2\6\1\6\2\6\1\6\2\6\eeee
\5\3\4\3\5\3\4\3\5\3\4\3\5\3\4\3\eeee
\5\3\4\3\5\3\4\3\5\3\4\3\5\3\4\3\eeee
\1\6\2\6\1\6\2\6\1\6\2\6\1\6\2\6\eeee
\1\6\2\6\1\6\2\6\1\6\2\6\1\6\2\6\eeee
\5\3\4\3\5\3\4\3\5\3\4\3\5\3\4\3\eeee
\5\3\4\3\5\3\4\3\5\3\4\3\5\3\4\3\eeee
\1\6\2\6\1\6\2\6\1\6\2\6\1\6\2\6\eeee
\1\6\2\6\1\6\2\6\1\6\2\6\1\6\2\6\eeee
\5\3\4\3\5\3\4\3\5\3\4\3\5\3\4\3\eeee
\5\3\4\3\5\3\4\3\5\3\4\3\5\3\4\3\eeee
} 

\baaa
6-84
\eaaa
\bbbb
1&0&0&0&1&2\\
0&1&0&1&1&1\\
0&0&1&1&2&0\\
0&2&1&1&0&0\\
1&1&1&0&1&0\\
2&1&0&0&0&1\\
\ebbb
\parbox{7cm}{ 
\1\5\3\5\1\6\2\4\2\6\1\5\3\5\1\6\eeee
\1\5\3\5\1\6\2\4\2\6\1\5\3\5\1\6\eeee
\6\2\4\2\6\1\5\3\5\1\6\2\4\2\6\1\eeee
\6\2\4\2\6\1\5\3\5\1\6\2\4\2\6\1\eeee
\1\5\3\5\1\6\2\4\2\6\1\5\3\5\1\6\eeee
\1\5\3\5\1\6\2\4\2\6\1\5\3\5\1\6\eeee
\6\2\4\2\6\1\5\3\5\1\6\2\4\2\6\1\eeee
\6\2\4\2\6\1\5\3\5\1\6\2\4\2\6\1\eeee
\1\5\3\5\1\6\2\4\2\6\1\5\3\5\1\6\eeee
\1\5\3\5\1\6\2\4\2\6\1\5\3\5\1\6\eeee
\6\2\4\2\6\1\5\3\5\1\6\2\4\2\6\1\eeee
\6\2\4\2\6\1\5\3\5\1\6\2\4\2\6\1\eeee
} 

\baaa
6-85
\eaaa
\bbbb
1&0&0&0&1&2\\
0&1&0&1&1&1\\
0&0&1&2&1&0\\
0&1&2&1&0&0\\
1&1&1&0&1&0\\
2&1&0&0&0&1\\
\ebbb
\parbox{7cm}{ 
\1\5\3\4\2\6\1\5\3\4\2\6\1\5\3\4\eeee
\1\5\3\4\2\6\1\5\3\4\2\6\1\5\3\4\eeee
\6\2\4\3\5\1\6\2\4\3\5\1\6\2\4\3\eeee
\6\2\4\3\5\1\6\2\4\3\5\1\6\2\4\3\eeee
\1\5\3\4\2\6\1\5\3\4\2\6\1\5\3\4\eeee
\1\5\3\4\2\6\1\5\3\4\2\6\1\5\3\4\eeee
\6\2\4\3\5\1\6\2\4\3\5\1\6\2\4\3\eeee
\6\2\4\3\5\1\6\2\4\3\5\1\6\2\4\3\eeee
\1\5\3\4\2\6\1\5\3\4\2\6\1\5\3\4\eeee
\1\5\3\4\2\6\1\5\3\4\2\6\1\5\3\4\eeee
\6\2\4\3\5\1\6\2\4\3\5\1\6\2\4\3\eeee
\6\2\4\3\5\1\6\2\4\3\5\1\6\2\4\3\eeee
} 

\baaa
6-86
\eaaa
\bbbb
1&0&0&0&1&2\\
0&1&0&1&1&1\\
0&0&2&1&0&1\\
0&1&1&2&0&0\\
1&2&0&0&1&0\\
1&1&1&0&0&1\\
\ebbb
\parbox{7cm}{ 
\1\6\3\3\6\1\6\3\3\6\1\6\3\3\6\1\eeee
\1\6\3\3\6\1\6\3\3\6\1\6\3\3\6\1\eeee
\5\2\4\4\2\5\2\4\4\2\5\2\4\4\2\5\eeee
\5\2\4\4\2\5\2\4\4\2\5\2\4\4\2\5\eeee
\1\6\3\3\6\1\6\3\3\6\1\6\3\3\6\1\eeee
\1\6\3\3\6\1\6\3\3\6\1\6\3\3\6\1\eeee
\5\2\4\4\2\5\2\4\4\2\5\2\4\4\2\5\eeee
\5\2\4\4\2\5\2\4\4\2\5\2\4\4\2\5\eeee
\1\6\3\3\6\1\6\3\3\6\1\6\3\3\6\1\eeee
\1\6\3\3\6\1\6\3\3\6\1\6\3\3\6\1\eeee
\5\2\4\4\2\5\2\4\4\2\5\2\4\4\2\5\eeee
\5\2\4\4\2\5\2\4\4\2\5\2\4\4\2\5\eeee
} 

\baaa
6-87
\eaaa
\bbbb
1&0&0&0&1&2\\
0&1&0&1&1&1\\
0&0&2&1&1&0\\
0&1&1&2&0&0\\
1&1&1&0&1&0\\
2&1&0&0&0&1\\
\ebbb
\parbox{7cm}{ 
\1\5\3\3\5\1\6\2\4\4\2\6\1\5\3\3\eeee
\1\5\3\3\5\1\6\2\4\4\2\6\1\5\3\3\eeee
\6\2\4\4\2\6\1\5\3\3\5\1\6\2\4\4\eeee
\6\2\4\4\2\6\1\5\3\3\5\1\6\2\4\4\eeee
\1\5\3\3\5\1\6\2\4\4\2\6\1\5\3\3\eeee
\1\5\3\3\5\1\6\2\4\4\2\6\1\5\3\3\eeee
\6\2\4\4\2\6\1\5\3\3\5\1\6\2\4\4\eeee
\6\2\4\4\2\6\1\5\3\3\5\1\6\2\4\4\eeee
\1\5\3\3\5\1\6\2\4\4\2\6\1\5\3\3\eeee
\1\5\3\3\5\1\6\2\4\4\2\6\1\5\3\3\eeee
\6\2\4\4\2\6\1\5\3\3\5\1\6\2\4\4\eeee
\6\2\4\4\2\6\1\5\3\3\5\1\6\2\4\4\eeee
} 

\baaa
6-88
\eaaa
\bbbb
1&0&0&1&1&1\\
0&1&0&1&1&1\\
0&0&1&1&1&1\\
1&1&1&1&0&0\\
1&1&1&0&1&0\\
1&1&1&0&0&1\\
\ebbb
\parbox{7cm}{ 
\1\4\3\5\2\6\1\4\3\5\2\6\1\4\3\5\eeee
\1\4\3\5\2\6\1\4\3\5\2\6\1\4\3\5\eeee
\5\2\6\1\4\3\5\2\6\1\4\3\5\2\6\1\eeee
\5\2\6\1\4\3\5\2\6\1\4\3\5\2\6\1\eeee
\1\4\3\5\2\6\1\4\3\5\2\6\1\4\3\5\eeee
\1\4\3\5\2\6\1\4\3\5\2\6\1\4\3\5\eeee
\5\2\6\1\4\3\5\2\6\1\4\3\5\2\6\1\eeee
\5\2\6\1\4\3\5\2\6\1\4\3\5\2\6\1\eeee
\1\4\3\5\2\6\1\4\3\5\2\6\1\4\3\5\eeee
\1\4\3\5\2\6\1\4\3\5\2\6\1\4\3\5\eeee
\5\2\6\1\4\3\5\2\6\1\4\3\5\2\6\1\eeee
\5\2\6\1\4\3\5\2\6\1\4\3\5\2\6\1\eeee
} 

\baaa
6-89
\eaaa
\bbbb
1&0&0&1&1&1\\
0&1&1&0&1&1\\
0&1&1&1&0&1\\
1&0&1&1&1&0\\
1&1&0&1&1&0\\
1&1&1&0&0&1\\
\ebbb
\parbox{7cm}{ 
\1\4\3\6\1\4\3\6\1\4\3\6\1\4\3\6\eeee
\1\4\3\6\1\4\3\6\1\4\3\6\1\4\3\6\eeee
\5\5\2\2\5\5\2\2\5\5\2\2\5\5\2\2\eeee
\4\1\6\3\4\1\6\3\4\1\6\3\4\1\6\3\eeee
\4\1\6\3\4\1\6\3\4\1\6\3\4\1\6\3\eeee
\5\5\2\2\5\5\2\2\5\5\2\2\5\5\2\2\eeee
\1\4\3\6\1\4\3\6\1\4\3\6\1\4\3\6\eeee
\1\4\3\6\1\4\3\6\1\4\3\6\1\4\3\6\eeee
\5\5\2\2\5\5\2\2\5\5\2\2\5\5\2\2\eeee
\4\1\6\3\4\1\6\3\4\1\6\3\4\1\6\3\eeee
\4\1\6\3\4\1\6\3\4\1\6\3\4\1\6\3\eeee
\5\5\2\2\5\5\2\2\5\5\2\2\5\5\2\2\eeee
} 
\baaa
\phantom{0-0}\#
\eaaa
\mbox{}\phantom{\bbbb
0&0&0&0&0&0\\
\ebbb}
\parbox{7cm}{ 
\1\4\5\1\4\5\1\4\5\1\4\5\1\4\5\1\eeee
\1\4\5\1\4\5\1\4\5\1\4\5\1\4\5\1\eeee
\6\3\2\6\3\2\6\3\2\6\3\2\6\3\2\6\eeee
\6\3\2\6\3\2\6\3\2\6\3\2\6\3\2\6\eeee
\1\4\5\1\4\5\1\4\5\1\4\5\1\4\5\1\eeee
\1\4\5\1\4\5\1\4\5\1\4\5\1\4\5\1\eeee
\6\3\2\6\3\2\6\3\2\6\3\2\6\3\2\6\eeee
\6\3\2\6\3\2\6\3\2\6\3\2\6\3\2\6\eeee
\1\4\5\1\4\5\1\4\5\1\4\5\1\4\5\1\eeee
\1\4\5\1\4\5\1\4\5\1\4\5\1\4\5\1\eeee
\6\3\2\6\3\2\6\3\2\6\3\2\6\3\2\6\eeee
\6\3\2\6\3\2\6\3\2\6\3\2\6\3\2\6\eeee
} 
\baaa
\phantom{0-0}\#
\eaaa
\mbox{}\phantom{\bbbb
0&0&0&0&0&0\\
\ebbb}
\parbox{7cm}{ 
\1\4\3\6\1\4\3\6\1\4\3\6\1\4\3\6\eeee
\1\5\2\6\1\5\2\6\1\5\2\6\1\5\2\6\eeee
\4\5\2\3\4\5\2\3\4\5\2\3\4\5\2\3\eeee
\4\1\6\3\4\1\6\3\4\1\6\3\4\1\6\3\eeee
\5\1\6\2\5\1\6\2\5\1\6\2\5\1\6\2\eeee
\5\4\3\2\5\4\3\2\5\4\3\2\5\4\3\2\eeee
\1\4\3\6\1\4\3\6\1\4\3\6\1\4\3\6\eeee
\1\5\2\6\1\5\2\6\1\5\2\6\1\5\2\6\eeee
\4\5\2\3\4\5\2\3\4\5\2\3\4\5\2\3\eeee
\4\1\6\3\4\1\6\3\4\1\6\3\4\1\6\3\eeee
\5\1\6\2\5\1\6\2\5\1\6\2\5\1\6\2\eeee
\5\4\3\2\5\4\3\2\5\4\3\2\5\4\3\2\eeee
} 

\baaa
6-90
\eaaa
\bbbb
1&0&0&1&1&1\\
0&1&1&0&1&1\\
0&1&1&1&0&1\\
1&0&1&2&0&0\\
1&1&0&0&2&0\\
1&1&1&0&0&1\\
\ebbb
\parbox{7cm}{ 
\1\4\4\1\5\5\1\4\4\1\5\5\1\4\4\1\eeee
\1\4\4\1\5\5\1\4\4\1\5\5\1\4\4\1\eeee
\6\3\3\6\2\2\6\3\3\6\2\2\6\3\3\6\eeee
\6\2\2\6\3\3\6\2\2\6\3\3\6\2\2\6\eeee
\1\5\5\1\4\4\1\5\5\1\4\4\1\5\5\1\eeee
\1\5\5\1\4\4\1\5\5\1\4\4\1\5\5\1\eeee
\6\2\2\6\3\3\6\2\2\6\3\3\6\2\2\6\eeee
\6\3\3\6\2\2\6\3\3\6\2\2\6\3\3\6\eeee
\1\4\4\1\5\5\1\4\4\1\5\5\1\4\4\1\eeee
\1\4\4\1\5\5\1\4\4\1\5\5\1\4\4\1\eeee
\6\3\3\6\2\2\6\3\3\6\2\2\6\3\3\6\eeee
\6\2\2\6\3\3\6\2\2\6\3\3\6\2\2\6\eeee
} 

\baaa
6-91
\eaaa
\bbbb
1&0&0&1&1&1\\
0&1&1&0&1&1\\
0&1&2&0&0&1\\
1&0&0&2&1&0\\
1&1&0&1&1&0\\
1&1&1&0&0&1\\
\ebbb
\parbox{7cm}{ 
\1\4\5\2\3\6\1\4\5\2\3\6\1\4\5\2\eeee
\1\4\5\2\3\6\1\4\5\2\3\6\1\4\5\2\eeee
\5\4\1\6\3\2\5\4\1\6\3\2\5\4\1\6\eeee
\5\4\1\6\3\2\5\4\1\6\3\2\5\4\1\6\eeee
\1\4\5\2\3\6\1\4\5\2\3\6\1\4\5\2\eeee
\1\4\5\2\3\6\1\4\5\2\3\6\1\4\5\2\eeee
\5\4\1\6\3\2\5\4\1\6\3\2\5\4\1\6\eeee
\5\4\1\6\3\2\5\4\1\6\3\2\5\4\1\6\eeee
\1\4\5\2\3\6\1\4\5\2\3\6\1\4\5\2\eeee
\1\4\5\2\3\6\1\4\5\2\3\6\1\4\5\2\eeee
\5\4\1\6\3\2\5\4\1\6\3\2\5\4\1\6\eeee
\5\4\1\6\3\2\5\4\1\6\3\2\5\4\1\6\eeee
} 

\baaa
6-92
\eaaa
\bbbb
1&0&0&1&1&1\\
0&2&0&0&1&1\\
0&0&2&1&0&1\\
1&0&1&2&0&0\\
1&1&0&0&2&0\\
1&1&1&0&0&1\\
\ebbb
\parbox{7cm}{ 
\1\4\4\1\5\5\1\4\4\1\5\5\1\4\4\1\eeee
\1\4\4\1\5\5\1\4\4\1\5\5\1\4\4\1\eeee
\6\3\3\6\2\2\6\3\3\6\2\2\6\3\3\6\eeee
\6\3\3\6\2\2\6\3\3\6\2\2\6\3\3\6\eeee
\1\4\4\1\5\5\1\4\4\1\5\5\1\4\4\1\eeee
\1\4\4\1\5\5\1\4\4\1\5\5\1\4\4\1\eeee
\6\3\3\6\2\2\6\3\3\6\2\2\6\3\3\6\eeee
\6\3\3\6\2\2\6\3\3\6\2\2\6\3\3\6\eeee
\1\4\4\1\5\5\1\4\4\1\5\5\1\4\4\1\eeee
\1\4\4\1\5\5\1\4\4\1\5\5\1\4\4\1\eeee
\6\3\3\6\2\2\6\3\3\6\2\2\6\3\3\6\eeee
\6\3\3\6\2\2\6\3\3\6\2\2\6\3\3\6\eeee
} 

\baaa
6-93
\eaaa
\bbbb
2&0&0&0&0&2\\
0&2&0&0&1&1\\
0&0&2&1&1&0\\
0&0&1&3&0&0\\
0&1&1&0&2&0\\
1&1&0&0&0&2\\
\ebbb
\parbox{7cm}{ 
\1\6\2\5\3\4\4\3\5\2\6\1\6\2\5\3\eeee
\1\6\2\5\3\4\4\3\5\2\6\1\6\2\5\3\eeee
\1\6\2\5\3\4\4\3\5\2\6\1\6\2\5\3\eeee
\1\6\2\5\3\4\4\3\5\2\6\1\6\2\5\3\eeee
\1\6\2\5\3\4\4\3\5\2\6\1\6\2\5\3\eeee
\1\6\2\5\3\4\4\3\5\2\6\1\6\2\5\3\eeee
\1\6\2\5\3\4\4\3\5\2\6\1\6\2\5\3\eeee
\1\6\2\5\3\4\4\3\5\2\6\1\6\2\5\3\eeee
\1\6\2\5\3\4\4\3\5\2\6\1\6\2\5\3\eeee
\1\6\2\5\3\4\4\3\5\2\6\1\6\2\5\3\eeee
\1\6\2\5\3\4\4\3\5\2\6\1\6\2\5\3\eeee
\1\6\2\5\3\4\4\3\5\2\6\1\6\2\5\3\eeee
} 

\baaa
6-94
\eaaa
\bbbb
2&0&0&0&0&2\\
0&2&0&0&1&1\\
0&0&2&1&1&0\\
0&0&2&2&0&0\\
0&1&1&0&2&0\\
1&1&0&0&0&2\\
\ebbb
\parbox{7cm}{ 
\1\6\2\5\3\4\3\5\2\6\1\6\2\5\3\4\eeee
\1\6\2\5\3\4\3\5\2\6\1\6\2\5\3\4\eeee
\1\6\2\5\3\4\3\5\2\6\1\6\2\5\3\4\eeee
\1\6\2\5\3\4\3\5\2\6\1\6\2\5\3\4\eeee
\1\6\2\5\3\4\3\5\2\6\1\6\2\5\3\4\eeee
\1\6\2\5\3\4\3\5\2\6\1\6\2\5\3\4\eeee
\1\6\2\5\3\4\3\5\2\6\1\6\2\5\3\4\eeee
\1\6\2\5\3\4\3\5\2\6\1\6\2\5\3\4\eeee
\1\6\2\5\3\4\3\5\2\6\1\6\2\5\3\4\eeee
\1\6\2\5\3\4\3\5\2\6\1\6\2\5\3\4\eeee
\1\6\2\5\3\4\3\5\2\6\1\6\2\5\3\4\eeee
\1\6\2\5\3\4\3\5\2\6\1\6\2\5\3\4\eeee
} 

\baaa
6-95
\eaaa
\bbbb
2&0&0&0&1&1\\
0&2&0&1&0&1\\
0&0&2&1&1&0\\
0&1&1&2&0&0\\
1&0&1&0&2&0\\
1&1&0&0&0&2\\
\ebbb
\parbox{7cm}{ 
\1\5\3\4\2\6\1\5\3\4\2\6\1\5\3\4\eeee
\1\5\3\4\2\6\1\5\3\4\2\6\1\5\3\4\eeee
\1\5\3\4\2\6\1\5\3\4\2\6\1\5\3\4\eeee
\1\5\3\4\2\6\1\5\3\4\2\6\1\5\3\4\eeee
\1\5\3\4\2\6\1\5\3\4\2\6\1\5\3\4\eeee
\1\5\3\4\2\6\1\5\3\4\2\6\1\5\3\4\eeee
\1\5\3\4\2\6\1\5\3\4\2\6\1\5\3\4\eeee
\1\5\3\4\2\6\1\5\3\4\2\6\1\5\3\4\eeee
\1\5\3\4\2\6\1\5\3\4\2\6\1\5\3\4\eeee
\1\5\3\4\2\6\1\5\3\4\2\6\1\5\3\4\eeee
\1\5\3\4\2\6\1\5\3\4\2\6\1\5\3\4\eeee
\1\5\3\4\2\6\1\5\3\4\2\6\1\5\3\4\eeee
} 

\baaa
6-96
\eaaa
\bbbb
2&0&0&0&1&1\\
0&2&0&1&0&1\\
0&0&3&0&1&0\\
0&1&0&3&0&0\\
1&0&1&0&2&0\\
1&1&0&0&0&2\\
\ebbb
\parbox{7cm}{ 
\1\5\3\3\5\1\6\2\4\4\2\6\1\5\3\3\eeee
\1\5\3\3\5\1\6\2\4\4\2\6\1\5\3\3\eeee
\1\5\3\3\5\1\6\2\4\4\2\6\1\5\3\3\eeee
\1\5\3\3\5\1\6\2\4\4\2\6\1\5\3\3\eeee
\1\5\3\3\5\1\6\2\4\4\2\6\1\5\3\3\eeee
\1\5\3\3\5\1\6\2\4\4\2\6\1\5\3\3\eeee
\1\5\3\3\5\1\6\2\4\4\2\6\1\5\3\3\eeee
\1\5\3\3\5\1\6\2\4\4\2\6\1\5\3\3\eeee
\1\5\3\3\5\1\6\2\4\4\2\6\1\5\3\3\eeee
\1\5\3\3\5\1\6\2\4\4\2\6\1\5\3\3\eeee
\1\5\3\3\5\1\6\2\4\4\2\6\1\5\3\3\eeee
\1\5\3\3\5\1\6\2\4\4\2\6\1\5\3\3\eeee
} 



\baaa
7-1
\eaaa
\bbbb
0&0&0&0&0&0&4\\
0&0&0&0&0&0&4\\
0&0&0&0&0&2&2\\
0&0&0&0&0&2&2\\
0&0&0&0&0&4&0\\
0&0&1&1&2&0&0\\
1&1&1&1&0&0&0\\
\ebbb
\parbox{7cm}{ 
\1\7\3\6\5\6\3\7\1\7\3\6\5\6\3\7\eeee
\7\2\7\4\6\5\6\4\7\2\7\4\6\5\6\4\eeee
\3\7\1\7\3\6\5\6\3\7\1\7\3\6\5\6\eeee
\6\4\7\2\7\4\6\5\6\4\7\2\7\4\6\5\eeee
\5\6\3\7\1\7\3\6\5\6\3\7\1\7\3\6\eeee
\6\5\6\4\7\2\7\4\6\5\6\4\7\2\7\4\eeee
\3\6\5\6\3\7\1\7\3\6\5\6\3\7\1\7\eeee
\7\4\6\5\6\4\7\2\7\4\6\5\6\4\7\2\eeee
\1\7\3\6\5\6\3\7\1\7\3\6\5\6\3\7\eeee
\7\2\7\4\6\5\6\4\7\2\7\4\6\5\6\4\eeee
\3\7\1\7\3\6\5\6\3\7\1\7\3\6\5\6\eeee
\6\4\7\2\7\4\6\5\6\4\7\2\7\4\6\5\eeee
} 

\baaa
7-2
\eaaa
\bbbb
0&0&0&0&0&0&4\\
0&0&0&0&0&0&4\\
0&0&0&0&0&2&2\\
0&0&0&0&0&2&2\\
0&0&0&0&2&2&0\\
0&0&1&1&2&0&0\\
1&1&1&1&0&0&0\\
\ebbb
\parbox{7cm}{ 
\1\7\3\6\5\5\6\3\7\1\7\3\6\5\5\6\eeee
\7\2\7\4\6\5\5\6\4\7\2\7\4\6\5\5\eeee
\3\7\1\7\3\6\5\5\6\3\7\1\7\3\6\5\eeee
\6\4\7\2\7\4\6\5\5\6\4\7\2\7\4\6\eeee
\5\6\3\7\1\7\3\6\5\5\6\3\7\1\7\3\eeee
\5\5\6\4\7\2\7\4\6\5\5\6\4\7\2\7\eeee
\6\5\5\6\3\7\1\7\3\6\5\5\6\3\7\1\eeee
\3\6\5\5\6\4\7\2\7\4\6\5\5\6\4\7\eeee
\7\4\6\5\5\6\3\7\1\7\3\6\5\5\6\3\eeee
\1\7\3\6\5\5\6\4\7\2\7\4\6\5\5\6\eeee
\7\2\7\4\6\5\5\6\3\7\1\7\3\6\5\5\eeee
\3\7\1\7\3\6\5\5\6\4\7\2\7\4\6\5\eeee
} 

\baaa
7-3
\eaaa
\bbbb
0&0&0&0&0&0&4\\
0&0&0&0&0&0&4\\
0&0&0&0&0&2&2\\
0&0&0&0&0&4&0\\
0&0&0&0&0&4&0\\
0&0&2&1&1&0&0\\
1&1&2&0&0&0&0\\
\ebbb
\parbox{7cm}{ 
\1\7\3\6\4\6\3\7\1\7\3\6\4\6\3\7\eeee
\7\2\7\3\6\5\6\3\7\2\7\3\6\5\6\3\eeee
\3\7\1\7\3\6\4\6\3\7\1\7\3\6\4\6\eeee
\6\3\7\2\7\3\6\5\6\3\7\2\7\3\6\5\eeee
\4\6\3\7\1\7\3\6\4\6\3\7\1\7\3\6\eeee
\6\5\6\3\7\2\7\3\6\5\6\3\7\2\7\3\eeee
\3\6\4\6\3\7\1\7\3\6\4\6\3\7\1\7\eeee
\7\3\6\5\6\3\7\2\7\3\6\5\6\3\7\2\eeee
\1\7\3\6\4\6\3\7\1\7\3\6\4\6\3\7\eeee
\7\2\7\3\6\5\6\3\7\2\7\3\6\5\6\3\eeee
\3\7\1\7\3\6\4\6\3\7\1\7\3\6\4\6\eeee
\6\3\7\2\7\3\6\5\6\3\7\2\7\3\6\5\eeee
} 

\baaa
7-4
\eaaa
\bbbb
0&0&0&0&0&0&4\\
0&0&0&0&0&0&4\\
0&0&0&0&0&2&2\\
0&0&0&0&2&2&0\\
0&0&0&2&2&0&0\\
0&0&2&2&0&0&0\\
1&1&2&0&0&0&0\\
\ebbb
\parbox{7cm}{ 
\1\7\3\6\4\5\5\4\6\3\7\1\7\3\6\4\eeee
\7\2\7\3\6\4\5\5\4\6\3\7\2\7\3\6\eeee
\3\7\1\7\3\6\4\5\5\4\6\3\7\1\7\3\eeee
\6\3\7\2\7\3\6\4\5\5\4\6\3\7\2\7\eeee
\4\6\3\7\1\7\3\6\4\5\5\4\6\3\7\1\eeee
\5\4\6\3\7\2\7\3\6\4\5\5\4\6\3\7\eeee
\5\5\4\6\3\7\1\7\3\6\4\5\5\4\6\3\eeee
\4\5\5\4\6\3\7\2\7\3\6\4\5\5\4\6\eeee
\6\4\5\5\4\6\3\7\1\7\3\6\4\5\5\4\eeee
\3\6\4\5\5\4\6\3\7\2\7\3\6\4\5\5\eeee
\7\3\6\4\5\5\4\6\3\7\1\7\3\6\4\5\eeee
\1\7\3\6\4\5\5\4\6\3\7\2\7\3\6\4\eeee
} 

\baaa
7-5
\eaaa
\bbbb
0&0&0&0&0&0&4\\
0&0&0&0&0&0&4\\
0&0&0&0&0&2&2\\
0&0&0&0&2&2&0\\
0&0&0&4&0&0&0\\
0&0&2&2&0&0&0\\
1&1&2&0&0&0&0\\
\ebbb
\parbox{7cm}{ 
\1\7\3\6\4\5\4\6\3\7\1\7\3\6\4\5\eeee
\7\2\7\3\6\4\5\4\6\3\7\2\7\3\6\4\eeee
\3\7\1\7\3\6\4\5\4\6\3\7\1\7\3\6\eeee
\6\3\7\2\7\3\6\4\5\4\6\3\7\2\7\3\eeee
\4\6\3\7\1\7\3\6\4\5\4\6\3\7\1\7\eeee
\5\4\6\3\7\2\7\3\6\4\5\4\6\3\7\2\eeee
\4\5\4\6\3\7\1\7\3\6\4\5\4\6\3\7\eeee
\6\4\5\4\6\3\7\2\7\3\6\4\5\4\6\3\eeee
\3\6\4\5\4\6\3\7\1\7\3\6\4\5\4\6\eeee
\7\3\6\4\5\4\6\3\7\2\7\3\6\4\5\4\eeee
\1\7\3\6\4\5\4\6\3\7\1\7\3\6\4\5\eeee
\7\2\7\3\6\4\5\4\6\3\7\2\7\3\6\4\eeee
} 

\baaa
7-6
\eaaa
\bbbb
0&0&0&0&0&0&4\\
0&0&0&0&0&0&4\\
0&0&0&0&0&2&2\\
0&0&0&1&1&2&0\\
0&0&0&1&1&2&0\\
0&0&2&1&1&0&0\\
1&1&2&0&0&0&0\\
\ebbb
\parbox{7cm}{ 
\1\7\3\6\4\4\6\3\7\1\7\3\6\4\4\6\eeee
\7\2\7\3\6\5\5\6\3\7\2\7\3\6\5\5\eeee
\3\7\1\7\3\6\4\4\6\3\7\1\7\3\6\4\eeee
\6\3\7\2\7\3\6\5\5\6\3\7\2\7\3\6\eeee
\4\6\3\7\1\7\3\6\4\4\6\3\7\1\7\3\eeee
\4\5\6\3\7\2\7\3\6\5\5\6\3\7\2\7\eeee
\6\5\4\6\3\7\1\7\3\6\4\4\6\3\7\1\eeee
\3\6\4\5\6\3\7\2\7\3\6\5\5\6\3\7\eeee
\7\3\6\5\4\6\3\7\1\7\3\6\4\4\6\3\eeee
\1\7\3\6\4\5\6\3\7\2\7\3\6\5\5\6\eeee
\7\2\7\3\6\5\4\6\3\7\1\7\3\6\4\4\eeee
\3\7\1\7\3\6\4\5\6\3\7\2\7\3\6\5\eeee
} 

\baaa
7-7
\eaaa
\bbbb
0&0&0&0&0&0&4\\
0&0&0&0&0&0&4\\
0&0&0&0&1&1&2\\
0&0&0&0&1&1&2\\
0&0&1&1&1&1&0\\
0&0&1&1&1&1&0\\
1&1&1&1&0&0&0\\
\ebbb
\parbox{7cm}{ 
\1\7\3\5\5\3\7\1\7\3\5\5\3\7\1\7\eeee
\7\2\7\4\6\6\4\7\2\7\4\6\6\4\7\2\eeee
\3\7\1\7\3\5\5\3\7\1\7\3\5\5\3\7\eeee
\5\4\7\2\7\4\6\6\4\7\2\7\4\6\6\4\eeee
\5\6\3\7\1\7\3\5\5\3\7\1\7\3\5\5\eeee
\3\6\5\4\7\2\7\4\6\6\4\7\2\7\4\6\eeee
\7\4\5\6\3\7\1\7\3\5\5\3\7\1\7\3\eeee
\1\7\3\6\5\4\7\2\7\4\6\6\4\7\2\7\eeee
\7\2\7\4\5\6\3\7\1\7\3\5\5\3\7\1\eeee
\3\7\1\7\3\6\5\4\7\2\7\4\6\6\4\7\eeee
\5\4\7\2\7\4\5\6\3\7\1\7\3\5\5\3\eeee
\5\6\3\7\1\7\3\6\5\4\7\2\7\4\6\6\eeee
} 

\baaa
7-8
\eaaa
\bbbb
0&0&0&0&0&0&4\\
0&0&0&0&0&0&4\\
0&0&0&0&1&1&2\\
0&0&0&0&1&1&2\\
0&0&2&2&0&0&0\\
0&0&2&2&0&0&0\\
1&1&1&1&0&0&0\\
\ebbb
\parbox{7cm}{ 
\1\7\3\5\3\7\1\7\3\5\3\7\1\7\3\5\eeee
\7\2\7\4\6\4\7\2\7\4\6\4\7\2\7\4\eeee
\3\7\1\7\3\5\3\7\1\7\3\5\3\7\1\7\eeee
\5\4\7\2\7\4\6\4\7\2\7\4\6\4\7\2\eeee
\3\6\3\7\1\7\3\5\3\7\1\7\3\5\3\7\eeee
\7\4\5\4\7\2\7\4\6\4\7\2\7\4\6\4\eeee
\1\7\3\6\3\7\1\7\3\5\3\7\1\7\3\5\eeee
\7\2\7\4\5\4\7\2\7\4\6\4\7\2\7\4\eeee
\3\7\1\7\3\6\3\7\1\7\3\5\3\7\1\7\eeee
\5\4\7\2\7\4\5\4\7\2\7\4\6\4\7\2\eeee
\3\6\3\7\1\7\3\6\3\7\1\7\3\5\3\7\eeee
\7\4\5\4\7\2\7\4\5\4\7\2\7\4\6\4\eeee
} 

\baaa
7-9
\eaaa
\bbbb
0&0&0&0&0&0&4\\
0&0&0&0&0&0&4\\
0&0&0&0&1&1&2\\
0&0&0&0&2&2&0\\
0&0&2&2&0&0&0\\
0&0&2&2&0&0&0\\
1&1&2&0&0&0&0\\
\ebbb
\parbox{7cm}{ 
\1\7\3\5\4\5\3\7\1\7\3\5\4\5\3\7\eeee
\7\2\7\3\6\4\6\3\7\2\7\3\6\4\6\3\eeee
\3\7\1\7\3\5\4\5\3\7\1\7\3\5\4\5\eeee
\5\3\7\2\7\3\6\4\6\3\7\2\7\3\6\4\eeee
\4\6\3\7\1\7\3\5\4\5\3\7\1\7\3\5\eeee
\5\4\5\3\7\2\7\3\6\4\6\3\7\2\7\3\eeee
\3\6\4\6\3\7\1\7\3\5\4\5\3\7\1\7\eeee
\7\3\5\4\5\3\7\2\7\3\6\4\6\3\7\2\eeee
\1\7\3\6\4\6\3\7\1\7\3\5\4\5\3\7\eeee
\7\2\7\3\5\4\5\3\7\2\7\3\6\4\6\3\eeee
\3\7\1\7\3\6\4\6\3\7\1\7\3\5\4\5\eeee
\5\3\7\2\7\3\5\4\5\3\7\2\7\3\6\4\eeee
} 

\baaa
7-10
\eaaa
\bbbb
0&0&0&0&0&0&4\\
0&0&0&0&0&0&4\\
0&0&0&0&1&1&2\\
0&0&0&2&1&1&0\\
0&0&2&2&0&0&0\\
0&0&2&2&0&0&0\\
1&1&2&0&0&0&0\\
\ebbb
\parbox{7cm}{ 
\1\7\3\5\4\4\5\3\7\1\7\3\5\4\4\5\eeee
\7\2\7\3\6\4\4\6\3\7\2\7\3\6\4\4\eeee
\3\7\1\7\3\5\4\4\5\3\7\1\7\3\5\4\eeee
\5\3\7\2\7\3\6\4\4\6\3\7\2\7\3\6\eeee
\4\6\3\7\1\7\3\5\4\4\5\3\7\1\7\3\eeee
\4\4\5\3\7\2\7\3\6\4\4\6\3\7\2\7\eeee
\5\4\4\6\3\7\1\7\3\5\4\4\5\3\7\1\eeee
\3\6\4\4\5\3\7\2\7\3\6\4\4\6\3\7\eeee
\7\3\5\4\4\6\3\7\1\7\3\5\4\4\5\3\eeee
\1\7\3\6\4\4\5\3\7\2\7\3\6\4\4\6\eeee
\7\2\7\3\5\4\4\6\3\7\1\7\3\5\4\4\eeee
\3\7\1\7\3\6\4\4\5\3\7\2\7\3\6\4\eeee
} 

\baaa
7-11
\eaaa
\bbbb
0&0&0&0&0&0&4\\
0&0&0&0&0&2&2\\
0&0&0&0&0&2&2\\
0&0&0&0&0&2&2\\
0&0&0&0&0&4&0\\
0&1&1&1&1&0&0\\
1&1&1&1&0&0&0\\
\ebbb
\parbox{7cm}{ 
\1\7\3\7\1\7\3\7\1\7\3\7\1\7\3\7\eeee
\7\2\6\4\7\2\6\4\7\2\6\4\7\2\6\4\eeee
\3\6\5\6\3\6\5\6\3\6\5\6\3\6\5\6\eeee
\7\4\6\2\7\4\6\2\7\4\6\2\7\4\6\2\eeee
\1\7\3\7\1\7\3\7\1\7\3\7\1\7\3\7\eeee
\7\2\6\4\7\2\6\4\7\2\6\4\7\2\6\4\eeee
\3\6\5\6\3\6\5\6\3\6\5\6\3\6\5\6\eeee
\7\4\6\2\7\4\6\2\7\4\6\2\7\4\6\2\eeee
\1\7\3\7\1\7\3\7\1\7\3\7\1\7\3\7\eeee
\7\2\6\4\7\2\6\4\7\2\6\4\7\2\6\4\eeee
\3\6\5\6\3\6\5\6\3\6\5\6\3\6\5\6\eeee
\7\4\6\2\7\4\6\2\7\4\6\2\7\4\6\2\eeee
} 

\baaa
7-12
\eaaa
\bbbb
0&0&0&0&0&0&4\\
0&0&0&0&0&2&2\\
0&0&0&0&0&2&2\\
0&0&0&0&1&2&1\\
0&0&0&4&0&0&0\\
0&1&1&2&0&0&0\\
1&1&1&1&0&0&0\\
\ebbb
\parbox{7cm}{ 
\1\7\4\5\4\7\1\7\4\5\4\7\1\7\4\5\eeee
\7\2\6\4\6\3\7\2\6\4\6\3\7\2\6\4\eeee
\4\6\3\7\2\6\4\6\3\7\2\6\4\6\3\7\eeee
\5\4\7\1\7\4\5\4\7\1\7\4\5\4\7\1\eeee
\4\6\2\7\3\6\4\6\2\7\3\6\4\6\2\7\eeee
\7\3\6\4\6\2\7\3\6\4\6\2\7\3\6\4\eeee
\1\7\4\5\4\7\1\7\4\5\4\7\1\7\4\5\eeee
\7\2\6\4\6\3\7\2\6\4\6\3\7\2\6\4\eeee
\4\6\3\7\2\6\4\6\3\7\2\6\4\6\3\7\eeee
\5\4\7\1\7\4\5\4\7\1\7\4\5\4\7\1\eeee
\4\6\2\7\3\6\4\6\2\7\3\6\4\6\2\7\eeee
\7\3\6\4\6\2\7\3\6\4\6\2\7\3\6\4\eeee
} 

\baaa
7-13
\eaaa
\bbbb
0&0&0&0&0&0&4\\
0&0&0&0&0&2&2\\
0&0&0&0&0&2&2\\
0&0&0&0&2&2&0\\
0&0&0&2&2&0&0\\
0&1&1&2&0&0&0\\
2&1&1&0&0&0&0\\
\ebbb
\parbox{7cm}{ 
\1\7\2\6\4\5\5\4\6\2\7\1\7\2\6\4\eeee
\7\1\7\3\6\4\5\5\4\6\3\7\1\7\3\6\eeee
\2\7\1\7\2\6\4\5\5\4\6\2\7\1\7\2\eeee
\6\3\7\1\7\3\6\4\5\5\4\6\3\7\1\7\eeee
\4\6\2\7\1\7\2\6\4\5\5\4\6\2\7\1\eeee
\5\4\6\3\7\1\7\3\6\4\5\5\4\6\3\7\eeee
\5\5\4\6\2\7\1\7\2\6\4\5\5\4\6\2\eeee
\4\5\5\4\6\3\7\1\7\3\6\4\5\5\4\6\eeee
\6\4\5\5\4\6\2\7\1\7\2\6\4\5\5\4\eeee
\2\6\4\5\5\4\6\3\7\1\7\3\6\4\5\5\eeee
\7\3\6\4\5\5\4\6\2\7\1\7\2\6\4\5\eeee
\1\7\2\6\4\5\5\4\6\3\7\1\7\3\6\4\eeee
} 

\baaa
7-14
\eaaa
\bbbb
0&0&0&0&0&0&4\\
0&0&0&0&0&2&2\\
0&0&0&0&0&2&2\\
0&0&0&0&2&2&0\\
0&0&0&4&0&0&0\\
0&1&1&2&0&0&0\\
2&1&1&0&0&0&0\\
\ebbb
\parbox{7cm}{ 
\1\7\2\6\4\5\4\6\2\7\1\7\2\6\4\5\eeee
\7\1\7\3\6\4\5\4\6\3\7\1\7\3\6\4\eeee
\2\7\1\7\2\6\4\5\4\6\2\7\1\7\2\6\eeee
\6\3\7\1\7\3\6\4\5\4\6\3\7\1\7\3\eeee
\4\6\2\7\1\7\2\6\4\5\4\6\2\7\1\7\eeee
\5\4\6\3\7\1\7\3\6\4\5\4\6\3\7\1\eeee
\4\5\4\6\2\7\1\7\2\6\4\5\4\6\2\7\eeee
\6\4\5\4\6\3\7\1\7\3\6\4\5\4\6\3\eeee
\2\6\4\5\4\6\2\7\1\7\2\6\4\5\4\6\eeee
\7\3\6\4\5\4\6\3\7\1\7\3\6\4\5\4\eeee
\1\7\2\6\4\5\4\6\2\7\1\7\2\6\4\5\eeee
\7\1\7\3\6\4\5\4\6\3\7\1\7\3\6\4\eeee
} 

\baaa
7-15
\eaaa
\bbbb
0&0&0&0&0&0&4\\
0&0&0&0&0&2&2\\
0&0&0&0&0&2&2\\
0&0&0&1&1&2&0\\
0&0&0&1&1&2&0\\
0&1&1&1&1&0&0\\
2&1&1&0&0&0&0\\
\ebbb
\parbox{7cm}{ 
\1\7\2\6\4\4\6\2\7\1\7\2\6\4\4\6\eeee
\7\1\7\3\6\5\5\6\3\7\1\7\3\6\5\5\eeee
\2\7\1\7\2\6\4\4\6\2\7\1\7\2\6\4\eeee
\6\3\7\1\7\3\6\5\5\6\3\7\1\7\3\6\eeee
\4\6\2\7\1\7\2\6\4\4\6\2\7\1\7\2\eeee
\4\5\6\3\7\1\7\3\6\5\5\6\3\7\1\7\eeee
\6\5\4\6\2\7\1\7\2\6\4\4\6\2\7\1\eeee
\2\6\4\5\6\3\7\1\7\3\6\5\5\6\3\7\eeee
\7\3\6\5\4\6\2\7\1\7\2\6\4\4\6\2\eeee
\1\7\2\6\4\5\6\3\7\1\7\3\6\5\5\6\eeee
\7\1\7\3\6\5\4\6\2\7\1\7\2\6\4\4\eeee
\2\7\1\7\2\6\4\5\6\3\7\1\7\3\6\5\eeee
} 

\baaa
7-16
\eaaa
\bbbb
0&0&0&0&0&0&4\\
0&0&0&0&0&2&2\\
0&0&0&0&2&2&0\\
0&0&0&0&2&2&0\\
0&0&1&1&2&0&0\\
0&2&1&1&0&0&0\\
2&2&0&0&0&0&0\\
\ebbb
\parbox{7cm}{ 
\1\7\2\6\3\5\5\3\6\2\7\1\7\2\6\3\eeee
\7\1\7\2\6\4\5\5\4\6\2\7\1\7\2\6\eeee
\2\7\1\7\2\6\3\5\5\3\6\2\7\1\7\2\eeee
\6\2\7\1\7\2\6\4\5\5\4\6\2\7\1\7\eeee
\3\6\2\7\1\7\2\6\3\5\5\3\6\2\7\1\eeee
\5\4\6\2\7\1\7\2\6\4\5\5\4\6\2\7\eeee
\5\5\3\6\2\7\1\7\2\6\3\5\5\3\6\2\eeee
\3\5\5\4\6\2\7\1\7\2\6\4\5\5\4\6\eeee
\6\4\5\5\3\6\2\7\1\7\2\6\3\5\5\3\eeee
\2\6\3\5\5\4\6\2\7\1\7\2\6\4\5\5\eeee
\7\2\6\4\5\5\3\6\2\7\1\7\2\6\3\5\eeee
\1\7\2\6\3\5\5\4\6\2\7\1\7\2\6\4\eeee
} 

\baaa
7-17
\eaaa
\bbbb
0&0&0&0&0&0&4\\
0&0&0&0&0&2&2\\
0&0&0&0&2&2&0\\
0&0&0&0&2&2&0\\
0&0&2&2&0&0&0\\
0&2&1&1&0&0&0\\
2&2&0&0&0&0&0\\
\ebbb
\parbox{7cm}{ 
\1\7\2\6\3\5\3\6\2\7\1\7\2\6\3\5\eeee
\7\1\7\2\6\4\5\4\6\2\7\1\7\2\6\4\eeee
\2\7\1\7\2\6\3\5\3\6\2\7\1\7\2\6\eeee
\6\2\7\1\7\2\6\4\5\4\6\2\7\1\7\2\eeee
\3\6\2\7\1\7\2\6\3\5\3\6\2\7\1\7\eeee
\5\4\6\2\7\1\7\2\6\4\5\4\6\2\7\1\eeee
\3\5\3\6\2\7\1\7\2\6\3\5\3\6\2\7\eeee
\6\4\5\4\6\2\7\1\7\2\6\4\5\4\6\2\eeee
\2\6\3\5\3\6\2\7\1\7\2\6\3\5\3\6\eeee
\7\2\6\4\5\4\6\2\7\1\7\2\6\4\5\4\eeee
\1\7\2\6\3\5\3\6\2\7\1\7\2\6\3\5\eeee
\7\1\7\2\6\4\5\4\6\2\7\1\7\2\6\4\eeee
} 

\baaa
7-18
\eaaa
\bbbb
0&0&0&0&0&0&4\\
0&0&0&0&0&2&2\\
0&0&0&0&2&2&0\\
0&0&0&0&4&0&0\\
0&0&2&2&0&0&0\\
0&2&2&0&0&0&0\\
2&2&0&0&0&0&0\\
\ebbb
\parbox{7cm}{ 
\1\7\2\6\3\5\4\5\3\6\2\7\1\7\2\6\eeee
\7\1\7\2\6\3\5\4\5\3\6\2\7\1\7\2\eeee
\2\7\1\7\2\6\3\5\4\5\3\6\2\7\1\7\eeee
\6\2\7\1\7\2\6\3\5\4\5\3\6\2\7\1\eeee
\3\6\2\7\1\7\2\6\3\5\4\5\3\6\2\7\eeee
\5\3\6\2\7\1\7\2\6\3\5\4\5\3\6\2\eeee
\4\5\3\6\2\7\1\7\2\6\3\5\4\5\3\6\eeee
\5\4\5\3\6\2\7\1\7\2\6\3\5\4\5\3\eeee
\3\5\4\5\3\6\2\7\1\7\2\6\3\5\4\5\eeee
\6\3\5\4\5\3\6\2\7\1\7\2\6\3\5\4\eeee
\2\6\3\5\4\5\3\6\2\7\1\7\2\6\3\5\eeee
\7\2\6\3\5\4\5\3\6\2\7\1\7\2\6\3\eeee
} 

\baaa
7-19
\eaaa
\bbbb
0&0&0&0&0&0&4\\
0&0&0&0&0&2&2\\
0&0&0&0&2&2&0\\
0&0&0&2&2&0&0\\
0&0&2&2&0&0&0\\
0&2&2&0&0&0&0\\
2&2&0&0&0&0&0\\
\ebbb
\parbox{7cm}{ 
\1\7\2\6\3\5\4\4\5\3\6\2\7\1\7\2\eeee
\7\1\7\2\6\3\5\4\4\5\3\6\2\7\1\7\eeee
\2\7\1\7\2\6\3\5\4\4\5\3\6\2\7\1\eeee
\6\2\7\1\7\2\6\3\5\4\4\5\3\6\2\7\eeee
\3\6\2\7\1\7\2\6\3\5\4\4\5\3\6\2\eeee
\5\3\6\2\7\1\7\2\6\3\5\4\4\5\3\6\eeee
\4\5\3\6\2\7\1\7\2\6\3\5\4\4\5\3\eeee
\4\4\5\3\6\2\7\1\7\2\6\3\5\4\4\5\eeee
\5\4\4\5\3\6\2\7\1\7\2\6\3\5\4\4\eeee
\3\5\4\4\5\3\6\2\7\1\7\2\6\3\5\4\eeee
\6\3\5\4\4\5\3\6\2\7\1\7\2\6\3\5\eeee
\2\6\3\5\4\4\5\3\6\2\7\1\7\2\6\3\eeee
} 

\baaa
7-20
\eaaa
\bbbb
0&0&0&0&0&0&4\\
0&0&0&0&0&2&2\\
0&0&0&1&1&2&0\\
0&0&2&1&1&0&0\\
0&0&2&1&1&0&0\\
0&2&2&0&0&0&0\\
2&2&0&0&0&0&0\\
\ebbb
\parbox{7cm}{ 
\1\7\2\6\3\4\4\3\6\2\7\1\7\2\6\3\eeee
\7\1\7\2\6\3\5\5\3\6\2\7\1\7\2\6\eeee
\2\7\1\7\2\6\3\4\4\3\6\2\7\1\7\2\eeee
\6\2\7\1\7\2\6\3\5\5\3\6\2\7\1\7\eeee
\3\6\2\7\1\7\2\6\3\4\4\3\6\2\7\1\eeee
\4\3\6\2\7\1\7\2\6\3\5\5\3\6\2\7\eeee
\4\5\3\6\2\7\1\7\2\6\3\4\4\3\6\2\eeee
\3\5\4\3\6\2\7\1\7\2\6\3\5\5\3\6\eeee
\6\3\4\5\3\6\2\7\1\7\2\6\3\4\4\3\eeee
\2\6\3\5\4\3\6\2\7\1\7\2\6\3\5\5\eeee
\7\2\6\3\4\5\3\6\2\7\1\7\2\6\3\4\eeee
\1\7\2\6\3\5\4\3\6\2\7\1\7\2\6\3\eeee
} 

\baaa
7-21
\eaaa
\bbbb
0&0&0&0&0&0&4\\
0&0&0&0&0&4&0\\
0&0&0&0&1&1&2\\
0&0&0&0&1&2&1\\
0&0&2&2&0&0&0\\
0&1&1&2&0&0&0\\
1&0&2&1&0&0&0\\
\ebbb
\parbox{7cm}{ 
\1\7\4\6\4\5\4\6\4\7\1\7\4\6\4\5\eeee
\7\3\6\2\6\3\7\3\5\3\7\3\6\2\6\3\eeee
\4\5\4\6\4\7\1\7\4\6\4\5\4\6\4\7\eeee
\6\3\7\3\5\3\7\3\6\2\6\3\7\3\5\3\eeee
\4\7\1\7\4\6\4\5\4\6\4\7\1\7\4\6\eeee
\5\3\7\3\6\2\6\3\7\3\5\3\7\3\6\2\eeee
\4\6\4\5\4\6\4\7\1\7\4\6\4\5\4\6\eeee
\6\2\6\3\7\3\5\3\7\3\6\2\6\3\7\3\eeee
\4\6\4\7\1\7\4\6\4\5\4\6\4\7\1\7\eeee
\7\3\5\3\7\3\6\2\6\3\7\3\5\3\7\3\eeee
\1\7\4\6\4\5\4\6\4\7\1\7\4\6\4\5\eeee
\7\3\6\2\6\3\7\3\5\3\7\3\6\2\6\3\eeee
} 

\baaa
7-22
\eaaa
\bbbb
0&0&0&0&0&0&4\\
0&0&0&0&1&1&2\\
0&0&0&0&1&1&2\\
0&0&0&0&2&2&0\\
0&1&1&2&0&0&0\\
0&1&1&2&0&0&0\\
2&1&1&0&0&0&0\\
\ebbb
\parbox{7cm}{ 
\1\7\2\5\4\5\2\7\1\7\2\5\4\5\2\7\eeee
\7\1\7\3\6\4\6\3\7\1\7\3\6\4\6\3\eeee
\2\7\1\7\2\5\4\5\2\7\1\7\2\5\4\5\eeee
\5\3\7\1\7\3\6\4\6\3\7\1\7\3\6\4\eeee
\4\6\2\7\1\7\2\5\4\5\2\7\1\7\2\5\eeee
\5\4\5\3\7\1\7\3\6\4\6\3\7\1\7\3\eeee
\2\6\4\6\2\7\1\7\2\5\4\5\2\7\1\7\eeee
\7\3\5\4\5\3\7\1\7\3\6\4\6\3\7\1\eeee
\1\7\2\6\4\6\2\7\1\7\2\5\4\5\2\7\eeee
\7\1\7\3\5\4\5\3\7\1\7\3\6\4\6\3\eeee
\2\7\1\7\2\6\4\6\2\7\1\7\2\5\4\5\eeee
\5\3\7\1\7\3\5\4\5\3\7\1\7\3\6\4\eeee
} 

\baaa
7-23
\eaaa
\bbbb
0&0&0&0&0&0&4\\
0&0&0&0&1&1&2\\
0&0&0&0&1&1&2\\
0&0&0&0&2&2&0\\
0&1&2&1&0&0&0\\
0&1&2&1&0&0&0\\
1&1&2&0&0&0&0\\
\ebbb
\parbox{7cm}{ 
\1\7\3\7\1\7\3\7\1\7\3\7\1\7\3\7\eeee
\7\2\6\3\7\2\5\3\7\2\6\3\7\2\5\3\eeee
\3\5\4\5\3\6\4\6\3\5\4\5\3\6\4\6\eeee
\7\3\6\2\7\3\5\2\7\3\6\2\7\3\5\2\eeee
\1\7\3\7\1\7\3\7\1\7\3\7\1\7\3\7\eeee
\7\2\5\3\7\2\6\3\7\2\5\3\7\2\6\3\eeee
\3\6\4\6\3\5\4\5\3\6\4\6\3\5\4\5\eeee
\7\3\5\2\7\3\6\2\7\3\5\2\7\3\6\2\eeee
\1\7\3\7\1\7\3\7\1\7\3\7\1\7\3\7\eeee
\7\2\6\3\7\2\5\3\7\2\6\3\7\2\5\3\eeee
\3\5\4\5\3\6\4\6\3\5\4\5\3\6\4\6\eeee
\7\3\6\2\7\3\5\2\7\3\6\2\7\3\5\2\eeee
} 

\baaa
7-24
\eaaa
\bbbb
0&0&0&0&0&0&4\\
0&0&0&0&1&1&2\\
0&0&0&0&1&1&2\\
0&0&0&2&1&1&0\\
0&1&1&2&0&0&0\\
0&1&1&2&0&0&0\\
2&1&1&0&0&0&0\\
\ebbb
\parbox{7cm}{ 
\1\7\2\5\4\4\5\2\7\1\7\2\5\4\4\5\eeee
\7\1\7\3\6\4\4\6\3\7\1\7\3\6\4\4\eeee
\2\7\1\7\2\5\4\4\5\2\7\1\7\2\5\4\eeee
\5\3\7\1\7\3\6\4\4\6\3\7\1\7\3\6\eeee
\4\6\2\7\1\7\2\5\4\4\5\2\7\1\7\2\eeee
\4\4\5\3\7\1\7\3\6\4\4\6\3\7\1\7\eeee
\5\4\4\6\2\7\1\7\2\5\4\4\5\2\7\1\eeee
\2\6\4\4\5\3\7\1\7\3\6\4\4\6\3\7\eeee
\7\3\5\4\4\6\2\7\1\7\2\5\4\4\5\2\eeee
\1\7\2\6\4\4\5\3\7\1\7\3\6\4\4\6\eeee
\7\1\7\3\5\4\4\6\2\7\1\7\2\5\4\4\eeee
\2\7\1\7\2\6\4\4\5\3\7\1\7\3\6\4\eeee
} 

\baaa
7-25
\eaaa
\bbbb
0&0&0&0&0&0&4\\
0&0&0&0&1&1&2\\
0&0&0&0&2&2&0\\
0&0&0&0&2&2&0\\
0&2&1&1&0&0&0\\
0&2&1&1&0&0&0\\
2&2&0&0&0&0&0\\
\ebbb
\parbox{7cm}{ 
\1\7\2\5\3\5\2\7\1\7\2\5\3\5\2\7\eeee
\7\1\7\2\6\4\6\2\7\1\7\2\6\4\6\2\eeee
\2\7\1\7\2\5\3\5\2\7\1\7\2\5\3\5\eeee
\5\2\7\1\7\2\6\4\6\2\7\1\7\2\6\4\eeee
\3\6\2\7\1\7\2\5\3\5\2\7\1\7\2\5\eeee
\5\4\5\2\7\1\7\2\6\4\6\2\7\1\7\2\eeee
\2\6\3\6\2\7\1\7\2\5\3\5\2\7\1\7\eeee
\7\2\5\4\5\2\7\1\7\2\6\4\6\2\7\1\eeee
\1\7\2\6\3\6\2\7\1\7\2\5\3\5\2\7\eeee
\7\1\7\2\5\4\5\2\7\1\7\2\6\4\6\2\eeee
\2\7\1\7\2\6\3\6\2\7\1\7\2\5\3\5\eeee
\5\2\7\1\7\2\5\4\5\2\7\1\7\2\6\4\eeee
} 

\baaa
7-26
\eaaa
\bbbb
0&0&0&0&0&0&4\\
0&0&0&0&1&1&2\\
0&0&0&1&1&1&1\\
0&0&1&2&0&1&0\\
0&1&1&0&1&1&0\\
0&1&1&1&1&0&0\\
1&2&1&0&0&0&0\\
\ebbb
\parbox{7cm}{ 
\1\7\3\4\4\6\2\7\2\5\6\3\5\5\3\6\eeee
\7\2\6\4\4\3\7\1\7\3\4\4\6\2\7\2\eeee
\3\5\5\3\6\5\2\7\2\6\4\4\3\7\1\7\eeee
\4\6\2\7\2\5\6\3\5\5\3\6\5\2\7\2\eeee
\4\3\7\1\7\3\4\4\6\2\7\2\5\6\3\5\eeee
\6\5\2\7\2\6\4\4\3\7\1\7\3\4\4\6\eeee
\2\5\6\3\5\5\3\6\5\2\7\2\6\4\4\3\eeee
\7\3\4\4\6\2\7\2\5\6\3\5\5\3\6\5\eeee
\2\6\4\4\3\7\1\7\3\4\4\6\2\7\2\5\eeee
\5\5\3\6\5\2\7\2\6\4\4\3\7\1\7\3\eeee
\6\2\7\2\5\6\3\5\5\3\6\5\2\7\2\6\eeee
\3\7\1\7\3\4\4\6\2\7\2\5\6\3\5\5\eeee
} 

\baaa
7-27
\eaaa
\bbbb
0&0&0&0&0&0&4\\
0&0&0&0&1&1&2\\
0&0&0&2&1&1&0\\
0&0&2&2&0&0&0\\
0&2&2&0&0&0&0\\
0&2&2&0&0&0&0\\
2&2&0&0&0&0&0\\
\ebbb
\parbox{7cm}{ 
\1\7\2\5\3\4\4\3\5\2\7\1\7\2\5\3\eeee
\7\1\7\2\6\3\4\4\3\6\2\7\1\7\2\6\eeee
\2\7\1\7\2\5\3\4\4\3\5\2\7\1\7\2\eeee
\5\2\7\1\7\2\6\3\4\4\3\6\2\7\1\7\eeee
\3\6\2\7\1\7\2\5\3\4\4\3\5\2\7\1\eeee
\4\3\5\2\7\1\7\2\6\3\4\4\3\6\2\7\eeee
\4\4\3\6\2\7\1\7\2\5\3\4\4\3\5\2\eeee
\3\4\4\3\5\2\7\1\7\2\6\3\4\4\3\6\eeee
\5\3\4\4\3\6\2\7\1\7\2\5\3\4\4\3\eeee
\2\6\3\4\4\3\5\2\7\1\7\2\6\3\4\4\eeee
\7\2\5\3\4\4\3\6\2\7\1\7\2\5\3\4\eeee
\1\7\2\6\3\4\4\3\5\2\7\1\7\2\6\3\eeee
} 

\baaa
7-28
\eaaa
\bbbb
0&0&0&0&0&0&4\\
0&0&0&0&1&1&2\\
0&0&1&0&1&1&1\\
0&0&0&2&1&1&0\\
0&1&1&1&0&1&0\\
0&1&1&1&1&0&0\\
1&2&1&0&0&0&0\\
\ebbb
\parbox{7cm}{ 
\1\7\3\3\7\1\7\3\3\7\1\7\3\3\7\1\eeee
\7\2\6\5\2\7\2\6\5\2\7\2\6\5\2\7\eeee
\3\5\4\4\6\3\5\4\4\6\3\5\4\4\6\3\eeee
\3\6\4\4\5\3\6\4\4\5\3\6\4\4\5\3\eeee
\7\2\5\6\2\7\2\5\6\2\7\2\5\6\2\7\eeee
\1\7\3\3\7\1\7\3\3\7\1\7\3\3\7\1\eeee
\7\2\6\5\2\7\2\6\5\2\7\2\6\5\2\7\eeee
\3\5\4\4\6\3\5\4\4\6\3\5\4\4\6\3\eeee
\3\6\4\4\5\3\6\4\4\5\3\6\4\4\5\3\eeee
\7\2\5\6\2\7\2\5\6\2\7\2\5\6\2\7\eeee
\1\7\3\3\7\1\7\3\3\7\1\7\3\3\7\1\eeee
\7\2\6\5\2\7\2\6\5\2\7\2\6\5\2\7\eeee
} 

\baaa
7-29
\eaaa
\bbbb
0&0&0&0&0&0&4\\
0&0&0&0&1&1&2\\
0&0&1&0&1&1&1\\
0&0&0&2&1&1&0\\
0&1&1&1&1&0&0\\
0&1&1&1&0&1&0\\
1&2&1&0&0&0&0\\
\ebbb
\parbox{7cm}{ 
\1\7\3\3\7\1\7\3\3\7\1\7\3\3\7\1\eeee
\7\2\6\6\2\7\2\5\5\2\7\2\6\6\2\7\eeee
\3\5\4\4\5\3\6\4\4\6\3\5\4\4\5\3\eeee
\3\5\4\4\5\3\6\4\4\6\3\5\4\4\5\3\eeee
\7\2\6\6\2\7\2\5\5\2\7\2\6\6\2\7\eeee
\1\7\3\3\7\1\7\3\3\7\1\7\3\3\7\1\eeee
\7\2\5\5\2\7\2\6\6\2\7\2\5\5\2\7\eeee
\3\6\4\4\6\3\5\4\4\5\3\6\4\4\6\3\eeee
\3\6\4\4\6\3\5\4\4\5\3\6\4\4\6\3\eeee
\7\2\5\5\2\7\2\6\6\2\7\2\5\5\2\7\eeee
\1\7\3\3\7\1\7\3\3\7\1\7\3\3\7\1\eeee
\7\2\6\6\2\7\2\5\5\2\7\2\6\6\2\7\eeee
} 

\baaa
7-30
\eaaa
\bbbb
0&0&0&0&0&0&4\\
0&0&0&0&1&1&2\\
0&0&1&1&1&1&0\\
0&0&1&1&1&1&0\\
0&2&1&1&0&0&0\\
0&2&1&1&0&0&0\\
2&2&0&0&0&0&0\\
\ebbb
\parbox{7cm}{ 
\1\7\2\5\3\3\5\2\7\1\7\2\5\3\3\5\eeee
\7\1\7\2\6\4\4\6\2\7\1\7\2\6\4\4\eeee
\2\7\1\7\2\5\3\3\5\2\7\1\7\2\5\3\eeee
\5\2\7\1\7\2\6\4\4\6\2\7\1\7\2\6\eeee
\3\6\2\7\1\7\2\5\3\3\5\2\7\1\7\2\eeee
\3\4\5\2\7\1\7\2\6\4\4\6\2\7\1\7\eeee
\5\4\3\6\2\7\1\7\2\5\3\3\5\2\7\1\eeee
\2\6\3\4\5\2\7\1\7\2\6\4\4\6\2\7\eeee
\7\2\5\4\3\6\2\7\1\7\2\5\3\3\5\2\eeee
\1\7\2\6\3\4\5\2\7\1\7\2\6\4\4\6\eeee
\7\1\7\2\5\4\3\6\2\7\1\7\2\5\3\3\eeee
\2\7\1\7\2\6\3\4\5\2\7\1\7\2\6\4\eeee
} 

\baaa
7-31
\eaaa
\bbbb
0&0&0&0&0&1&3\\
0&0&0&0&0&1&3\\
0&0&0&0&1&0&3\\
0&0&0&0&1&0&3\\
0&0&2&2&0&0&0\\
2&2&0&0&0&0&0\\
1&1&1&1&0&0&0\\
\ebbb
\parbox{7cm}{ 
\1\7\3\7\2\7\4\7\1\7\3\7\2\7\4\7\eeee
\6\2\7\1\6\1\7\2\6\2\7\1\6\1\7\2\eeee
\1\7\4\7\2\7\3\7\1\7\4\7\2\7\3\7\eeee
\7\3\5\3\7\4\5\4\7\3\5\3\7\4\5\4\eeee
\2\7\4\7\1\7\3\7\2\7\4\7\1\7\3\7\eeee
\6\1\7\2\6\2\7\1\6\1\7\2\6\2\7\1\eeee
\2\7\3\7\1\7\4\7\2\7\3\7\1\7\4\7\eeee
\7\4\5\4\7\3\5\3\7\4\5\4\7\3\5\3\eeee
\1\7\3\7\2\7\4\7\1\7\3\7\2\7\4\7\eeee
\6\2\7\1\6\1\7\2\6\2\7\1\6\1\7\2\eeee
\1\7\4\7\2\7\3\7\1\7\4\7\2\7\3\7\eeee
\7\3\5\3\7\4\5\4\7\3\5\3\7\4\5\4\eeee
} 

\baaa
7-32
\eaaa
\bbbb
0&0&0&0&0&1&3\\
0&0&0&0&1&2&1\\
0&0&0&1&2&1&0\\
0&0&1&2&1&0&0\\
0&1&2&1&0&0&0\\
1&2&1&0&0&0&0\\
3&1&0&0&0&0&0\\
\ebbb
\parbox{7cm}{ 
\1\7\1\7\1\7\1\7\1\7\1\7\1\7\1\7\eeee
\6\2\6\2\6\2\6\2\6\2\6\2\6\2\6\2\eeee
\3\5\3\5\3\5\3\5\3\5\3\5\3\5\3\5\eeee
\4\4\4\4\4\4\4\4\4\4\4\4\4\4\4\4\eeee
\5\3\5\3\5\3\5\3\5\3\5\3\5\3\5\3\eeee
\2\6\2\6\2\6\2\6\2\6\2\6\2\6\2\6\eeee
\7\1\7\1\7\1\7\1\7\1\7\1\7\1\7\1\eeee
\1\7\1\7\1\7\1\7\1\7\1\7\1\7\1\7\eeee
\6\2\6\2\6\2\6\2\6\2\6\2\6\2\6\2\eeee
\3\5\3\5\3\5\3\5\3\5\3\5\3\5\3\5\eeee
\4\4\4\4\4\4\4\4\4\4\4\4\4\4\4\4\eeee
\5\3\5\3\5\3\5\3\5\3\5\3\5\3\5\3\eeee
} 

\baaa
7-33
\eaaa
\bbbb
0&0&0&0&0&2&2\\
0&0&0&0&0&2&2\\
0&0&0&0&2&0&2\\
0&0&0&0&2&0&2\\
0&0&1&1&0&2&0\\
1&1&0&0&2&0&0\\
1&1&1&1&0&0&0\\
\ebbb
\parbox{7cm}{ 
\1\6\5\3\7\2\6\5\3\7\2\6\5\4\7\1\eeee
\6\5\4\7\1\6\5\4\7\1\6\5\3\7\2\6\eeee
\5\3\7\2\6\5\3\7\2\6\5\4\7\1\6\5\eeee
\4\7\1\6\5\4\7\1\6\5\3\7\2\6\5\3\eeee
\7\2\6\5\3\7\2\6\5\4\7\1\6\5\4\7\eeee
\1\6\5\4\7\1\6\5\3\7\2\6\5\3\7\2\eeee
\6\5\3\7\2\6\5\4\7\1\6\5\4\7\1\6\eeee
\5\4\7\1\6\5\3\7\2\6\5\3\7\2\6\5\eeee
\3\7\2\6\5\4\7\1\6\5\4\7\1\6\5\3\eeee
\7\1\6\5\3\7\2\6\5\3\7\2\6\5\4\7\eeee
\2\6\5\4\7\1\6\5\4\7\1\6\5\3\7\2\eeee
\6\5\3\7\2\6\5\3\7\2\6\5\4\7\1\6\eeee
} 

\baaa
7-34
\eaaa
\bbbb
0&0&0&0&0&2&2\\
0&0&0&0&0&2&2\\
0&0&0&0&2&0&2\\
0&0&0&0&2&0&2\\
0&0&1&1&2&0&0\\
1&1&0&0&0&2&0\\
1&1&1&1&0&0&0\\
\ebbb
\parbox{7cm}{ 
\1\6\6\1\7\4\5\5\3\7\2\6\6\2\7\3\eeee
\6\6\2\7\3\5\5\4\7\1\6\6\1\7\4\5\eeee
\6\1\7\4\5\5\3\7\2\6\6\2\7\3\5\5\eeee
\2\7\3\5\5\4\7\1\6\6\1\7\4\5\5\3\eeee
\7\4\5\5\3\7\2\6\6\2\7\3\5\5\4\7\eeee
\3\5\5\4\7\1\6\6\1\7\4\5\5\3\7\2\eeee
\5\5\3\7\2\6\6\2\7\3\5\5\4\7\1\6\eeee
\5\4\7\1\6\6\1\7\4\5\5\3\7\2\6\6\eeee
\3\7\2\6\6\2\7\3\5\5\4\7\1\6\6\1\eeee
\7\1\6\6\1\7\4\5\5\3\7\2\6\6\2\7\eeee
\2\6\6\2\7\3\5\5\4\7\1\6\6\1\7\4\eeee
\6\6\1\7\4\5\5\3\7\2\6\6\2\7\3\5\eeee
} 

\baaa
7-35
\eaaa
\bbbb
0&0&0&0&0&2&2\\
0&0&0&0&0&2&2\\
0&0&0&0&2&0&2\\
0&0&0&0&2&0&2\\
0&0&1&1&2&0&0\\
2&2&0&0&0&0&0\\
1&1&1&1&0&0&0\\
\ebbb
\parbox{7cm}{ 
\1\6\2\7\3\5\5\3\7\2\6\1\7\4\5\5\eeee
\6\1\7\4\5\5\4\7\1\6\2\7\3\5\5\4\eeee
\2\7\3\5\5\3\7\2\6\1\7\4\5\5\3\7\eeee
\7\4\5\5\4\7\1\6\2\7\3\5\5\4\7\1\eeee
\3\5\5\3\7\2\6\1\7\4\5\5\3\7\2\6\eeee
\5\5\4\7\1\6\2\7\3\5\5\4\7\1\6\2\eeee
\5\3\7\2\6\1\7\4\5\5\3\7\2\6\1\7\eeee
\4\7\1\6\2\7\3\5\5\4\7\1\6\2\7\3\eeee
\7\2\6\1\7\4\5\5\3\7\2\6\1\7\4\5\eeee
\1\6\2\7\3\5\5\4\7\1\6\2\7\3\5\5\eeee
\6\1\7\4\5\5\3\7\2\6\1\7\4\5\5\4\eeee
\2\7\3\5\5\4\7\1\6\2\7\3\5\5\3\7\eeee
} 

\baaa
7-36
\eaaa
\bbbb
0&0&0&0&0&2&2\\
0&0&0&0&0&2&2\\
0&0&0&0&2&0&2\\
0&0&0&0&2&0&2\\
0&0&2&2&0&0&0\\
2&2&0&0&0&0&0\\
1&1&1&1&0&0&0\\
\ebbb
\parbox{7cm}{ 
\1\6\2\7\3\5\4\7\1\6\2\7\3\5\4\7\eeee
\6\1\7\4\5\3\7\2\6\1\7\4\5\3\7\2\eeee
\2\7\3\5\4\7\1\6\2\7\3\5\4\7\1\6\eeee
\7\4\5\3\7\2\6\1\7\4\5\3\7\2\6\1\eeee
\3\5\4\7\1\6\2\7\3\5\4\7\1\6\2\7\eeee
\5\3\7\2\6\1\7\4\5\3\7\2\6\1\7\4\eeee
\4\7\1\6\2\7\3\5\4\7\1\6\2\7\3\5\eeee
\7\2\6\1\7\4\5\3\7\2\6\1\7\4\5\3\eeee
\1\6\2\7\3\5\4\7\1\6\2\7\3\5\4\7\eeee
\6\1\7\4\5\3\7\2\6\1\7\4\5\3\7\2\eeee
\2\7\3\5\4\7\1\6\2\7\3\5\4\7\1\6\eeee
\7\4\5\3\7\2\6\1\7\4\5\3\7\2\6\1\eeee
} 

\baaa
7-37
\eaaa
\bbbb
0&0&0&0&0&2&2\\
0&0&0&0&0&2&2\\
0&0&0&0&2&0&2\\
0&0&0&0&2&2&0\\
0&0&2&2&0&0&0\\
1&1&0&2&0&0&0\\
1&1&2&0&0&0&0\\
\ebbb
\parbox{7cm}{ 
\1\6\4\5\3\7\2\6\4\5\3\7\1\6\4\5\eeee
\6\4\5\3\7\1\6\4\5\3\7\2\6\4\5\3\eeee
\4\5\3\7\2\6\4\5\3\7\1\6\4\5\3\7\eeee
\5\3\7\1\6\4\5\3\7\2\6\4\5\3\7\1\eeee
\3\7\2\6\4\5\3\7\1\6\4\5\3\7\2\6\eeee
\7\1\6\4\5\3\7\2\6\4\5\3\7\1\6\4\eeee
\2\6\4\5\3\7\1\6\4\5\3\7\2\6\4\5\eeee
\6\4\5\3\7\2\6\4\5\3\7\1\6\4\5\3\eeee
\4\5\3\7\1\6\4\5\3\7\2\6\4\5\3\7\eeee
\5\3\7\2\6\4\5\3\7\1\6\4\5\3\7\2\eeee
\3\7\1\6\4\5\3\7\2\6\4\5\3\7\1\6\eeee
\7\2\6\4\5\3\7\1\6\4\5\3\7\2\6\4\eeee
} 

\baaa
7-38
\eaaa
\bbbb
0&0&0&0&0&2&2\\
0&0&0&0&0&2&2\\
0&0&0&0&2&0&2\\
0&0&0&2&0&2&0\\
0&0&2&0&2&0&0\\
1&1&0&2&0&0&0\\
1&1&2&0&0&0&0\\
\ebbb
\parbox{7cm}{ 
\1\6\4\4\6\2\7\3\5\5\3\7\1\6\4\4\eeee
\6\4\4\6\1\7\3\5\5\3\7\2\6\4\4\6\eeee
\4\4\6\2\7\3\5\5\3\7\1\6\4\4\6\2\eeee
\4\6\1\7\3\5\5\3\7\2\6\4\4\6\1\7\eeee
\6\2\7\3\5\5\3\7\1\6\4\4\6\2\7\3\eeee
\1\7\3\5\5\3\7\2\6\4\4\6\1\7\3\5\eeee
\7\3\5\5\3\7\1\6\4\4\6\2\7\3\5\5\eeee
\3\5\5\3\7\2\6\4\4\6\1\7\3\5\5\3\eeee
\5\5\3\7\1\6\4\4\6\2\7\3\5\5\3\7\eeee
\5\3\7\2\6\4\4\6\1\7\3\5\5\3\7\2\eeee
\3\7\1\6\4\4\6\2\7\3\5\5\3\7\1\6\eeee
\7\2\6\4\4\6\1\7\3\5\5\3\7\2\6\4\eeee
} 

\baaa
7-39
\eaaa
\bbbb
0&0&0&0&0&2&2\\
0&0&0&0&0&2&2\\
0&0&0&0&2&0&2\\
0&0&0&2&2&0&0\\
0&0&2&2&0&0&0\\
1&1&0&0&0&2&0\\
1&1&2&0&0&0&0\\
\ebbb
\parbox{7cm}{ 
\1\6\6\1\7\3\5\4\4\5\3\7\1\6\6\1\eeee
\6\6\2\7\3\5\4\4\5\3\7\2\6\6\2\7\eeee
\6\1\7\3\5\4\4\5\3\7\1\6\6\1\7\3\eeee
\2\7\3\5\4\4\5\3\7\2\6\6\2\7\3\5\eeee
\7\3\5\4\4\5\3\7\1\6\6\1\7\3\5\4\eeee
\3\5\4\4\5\3\7\2\6\6\2\7\3\5\4\4\eeee
\5\4\4\5\3\7\1\6\6\1\7\3\5\4\4\5\eeee
\4\4\5\3\7\2\6\6\2\7\3\5\4\4\5\3\eeee
\4\5\3\7\1\6\6\1\7\3\5\4\4\5\3\7\eeee
\5\3\7\2\6\6\2\7\3\5\4\4\5\3\7\2\eeee
\3\7\1\6\6\1\7\3\5\4\4\5\3\7\1\6\eeee
\7\2\6\6\2\7\3\5\4\4\5\3\7\2\6\6\eeee
} 

\baaa
7-40
\eaaa
\bbbb
0&0&0&0&0&2&2\\
0&0&0&0&0&2&2\\
0&0&0&0&2&0&2\\
0&0&0&2&2&0&0\\
0&0&2&2&0&0&0\\
2&2&0&0&0&0&0\\
1&1&2&0&0&0&0\\
\ebbb
\parbox{7cm}{ 
\1\6\2\7\3\5\4\4\5\3\7\1\6\2\7\3\eeee
\6\1\7\3\5\4\4\5\3\7\2\6\1\7\3\5\eeee
\2\7\3\5\4\4\5\3\7\1\6\2\7\3\5\4\eeee
\7\3\5\4\4\5\3\7\2\6\1\7\3\5\4\4\eeee
\3\5\4\4\5\3\7\1\6\2\7\3\5\4\4\5\eeee
\5\4\4\5\3\7\2\6\1\7\3\5\4\4\5\3\eeee
\4\4\5\3\7\1\6\2\7\3\5\4\4\5\3\7\eeee
\4\5\3\7\2\6\1\7\3\5\4\4\5\3\7\1\eeee
\5\3\7\1\6\2\7\3\5\4\4\5\3\7\2\6\eeee
\3\7\2\6\1\7\3\5\4\4\5\3\7\1\6\2\eeee
\7\1\6\2\7\3\5\4\4\5\3\7\2\6\1\7\eeee
\2\6\1\7\3\5\4\4\5\3\7\1\6\2\7\3\eeee
} 

\baaa
7-41
\eaaa
\bbbb
0&0&0&0&0&2&2\\
0&0&0&0&0&2&2\\
0&0&0&0&2&1&1\\
0&0&0&0&2&1&1\\
0&0&1&1&2&0&0\\
1&1&1&1&0&0&0\\
1&1&1&1&0&0&0\\
\ebbb
\parbox{7cm}{ 
\1\6\3\5\5\3\6\1\6\3\5\5\3\6\1\6\eeee
\6\2\7\4\5\5\4\7\2\7\4\5\5\4\7\2\eeee
\3\7\1\6\3\5\5\3\6\1\6\3\5\5\3\6\eeee
\5\4\6\2\7\4\5\5\4\7\2\7\4\5\5\4\eeee
\5\5\3\7\1\6\3\5\5\3\6\1\6\3\5\5\eeee
\3\5\5\4\6\2\7\4\5\5\4\7\2\7\4\5\eeee
\6\4\5\5\3\7\1\6\3\5\5\3\6\1\6\3\eeee
\1\7\3\5\5\4\6\2\7\4\5\5\4\7\2\7\eeee
\6\2\6\4\5\5\3\7\1\6\3\5\5\3\6\1\eeee
\3\7\1\7\3\5\5\4\6\2\7\4\5\5\4\7\eeee
\5\4\6\2\6\4\5\5\3\7\1\6\3\5\5\3\eeee
\5\5\3\7\1\7\3\5\5\4\6\2\7\4\5\5\eeee
} 

\baaa
7-42
\eaaa
\bbbb
0&0&0&0&0&2&2\\
0&0&0&0&0&2&2\\
0&0&0&0&2&1&1\\
0&0&0&0&2&1&1\\
0&0&2&2&0&0&0\\
1&1&1&1&0&0&0\\
1&1&1&1&0&0&0\\
\ebbb
\parbox{7cm}{ 
\1\6\3\5\3\6\1\6\3\5\3\6\1\6\3\5\eeee
\6\2\7\4\5\4\7\2\7\4\5\4\7\2\7\4\eeee
\3\7\1\6\3\5\3\6\1\6\3\5\3\6\1\6\eeee
\5\4\6\2\7\4\5\4\7\2\7\4\5\4\7\2\eeee
\3\5\3\7\1\6\3\5\3\6\1\6\3\5\3\6\eeee
\6\4\5\4\6\2\7\4\5\4\7\2\7\4\5\4\eeee
\1\7\3\5\3\7\1\6\3\5\3\6\1\6\3\5\eeee
\6\2\6\4\5\4\6\2\7\4\5\4\7\2\7\4\eeee
\3\7\1\7\3\5\3\7\1\6\3\5\3\6\1\6\eeee
\5\4\6\2\6\4\5\4\6\2\7\4\5\4\7\2\eeee
\3\5\3\7\1\7\3\5\3\7\1\6\3\5\3\6\eeee
\6\4\5\4\6\2\6\4\5\4\6\2\7\4\5\4\eeee
} 

\baaa
7-43
\eaaa
\bbbb
0&0&0&0&0&2&2\\
0&0&0&0&0&2&2\\
0&0&0&0&2&1&1\\
0&0&0&2&2&0&0\\
0&0&2&2&0&0&0\\
1&1&2&0&0&0&0\\
1&1&2&0&0&0&0\\
\ebbb
\parbox{7cm}{ 
\1\6\3\5\4\4\5\3\6\1\6\3\5\4\4\5\eeee
\6\2\7\3\5\4\4\5\3\7\2\7\3\5\4\4\eeee
\3\7\1\6\3\5\4\4\5\3\6\1\6\3\5\4\eeee
\5\3\6\2\7\3\5\4\4\5\3\7\2\7\3\5\eeee
\4\5\3\7\1\6\3\5\4\4\5\3\6\1\6\3\eeee
\4\4\5\3\6\2\7\3\5\4\4\5\3\7\2\7\eeee
\5\4\4\5\3\7\1\6\3\5\4\4\5\3\6\1\eeee
\3\5\4\4\5\3\6\2\7\3\5\4\4\5\3\7\eeee
\6\3\5\4\4\5\3\7\1\6\3\5\4\4\5\3\eeee
\1\7\3\5\4\4\5\3\6\2\7\3\5\4\4\5\eeee
\6\2\6\3\5\4\4\5\3\7\1\6\3\5\4\4\eeee
\3\7\1\7\3\5\4\4\5\3\6\2\7\3\5\4\eeee
} 

\baaa
7-44
\eaaa
\bbbb
0&0&0&0&0&2&2\\
0&0&0&0&0&2&2\\
0&0&0&1&1&0&2\\
0&0&2&1&1&0&0\\
0&0&2&1&1&0&0\\
1&1&0&0&0&2&0\\
1&1&2&0&0&0&0\\
\ebbb
\parbox{7cm}{ 
\1\6\6\1\7\3\5\4\3\7\2\6\6\2\7\3\eeee
\6\6\2\7\3\4\5\3\7\1\6\6\1\7\3\5\eeee
\6\1\7\3\5\4\3\7\2\6\6\2\7\3\4\5\eeee
\2\7\3\4\5\3\7\1\6\6\1\7\3\5\4\3\eeee
\7\3\5\4\3\7\2\6\6\2\7\3\4\5\3\7\eeee
\3\4\5\3\7\1\6\6\1\7\3\5\4\3\7\2\eeee
\5\4\3\7\2\6\6\2\7\3\4\5\3\7\1\6\eeee
\5\3\7\1\6\6\1\7\3\5\4\3\7\2\6\6\eeee
\3\7\2\6\6\2\7\3\4\5\3\7\1\6\6\1\eeee
\7\1\6\6\1\7\3\5\4\3\7\2\6\6\2\7\eeee
\2\6\6\2\7\3\4\5\3\7\1\6\6\1\7\3\eeee
\6\6\1\7\3\5\4\3\7\2\6\6\2\7\3\4\eeee
} 

\baaa
7-45
\eaaa
\bbbb
0&0&0&0&0&2&2\\
0&0&0&0&0&2&2\\
0&0&0&1&1&0&2\\
0&0&2&1&1&0&0\\
0&0&2&1&1&0&0\\
2&2&0&0&0&0&0\\
1&1&2&0&0&0&0\\
\ebbb
\parbox{7cm}{ 
\1\6\2\7\3\5\5\3\7\2\6\1\7\3\5\4\eeee
\6\1\7\3\4\4\3\7\1\6\2\7\3\4\5\3\eeee
\2\7\3\5\5\3\7\2\6\1\7\3\5\4\3\7\eeee
\7\3\4\4\3\7\1\6\2\7\3\4\5\3\7\1\eeee
\3\5\5\3\7\2\6\1\7\3\5\4\3\7\2\6\eeee
\4\4\3\7\1\6\2\7\3\4\5\3\7\1\6\2\eeee
\5\3\7\2\6\1\7\3\5\4\3\7\2\6\1\7\eeee
\3\7\1\6\2\7\3\4\5\3\7\1\6\2\7\3\eeee
\7\2\6\1\7\3\5\4\3\7\2\6\1\7\3\4\eeee
\1\6\2\7\3\4\5\3\7\1\6\2\7\3\5\5\eeee
\6\1\7\3\5\4\3\7\2\6\1\7\3\4\4\3\eeee
\2\7\3\4\5\3\7\1\6\2\7\3\5\5\3\7\eeee
} 

\baaa
7-46
\eaaa
\bbbb
0&0&0&0&0&2&2\\
0&0&0&0&0&2&2\\
0&0&0&1&1&1&1\\
0&0&2&1&1&0&0\\
0&0&2&1&1&0&0\\
1&1&2&0&0&0&0\\
1&1&2&0&0&0&0\\
\ebbb
\parbox{7cm}{ 
\1\6\3\4\4\3\6\1\6\3\4\4\3\6\1\6\eeee
\6\2\7\3\5\5\3\7\2\7\3\5\5\3\7\2\eeee
\3\7\1\6\3\4\4\3\6\1\6\3\4\4\3\6\eeee
\4\3\6\2\7\3\5\5\3\7\2\7\3\5\5\3\eeee
\4\5\3\7\1\6\3\4\4\3\6\1\6\3\4\4\eeee
\3\5\4\3\6\2\7\3\5\5\3\7\2\7\3\5\eeee
\6\3\4\5\3\7\1\6\3\4\4\3\6\1\6\3\eeee
\1\7\3\5\4\3\6\2\7\3\5\5\3\7\2\7\eeee
\6\2\6\3\4\5\3\7\1\6\3\4\4\3\6\1\eeee
\3\7\1\7\3\5\4\3\6\2\7\3\5\5\3\7\eeee
\4\3\6\2\6\3\4\5\3\7\1\6\3\4\4\3\eeee
\4\5\3\7\1\7\3\5\4\3\6\2\7\3\5\5\eeee
} 

\baaa
7-47
\eaaa
\bbbb
0&0&0&0&0&2&2\\
0&0&0&0&0&2&2\\
0&0&1&0&1&0&2\\
0&0&0&2&0&2&0\\
0&0&1&0&1&0&2\\
1&1&0&2&0&0&0\\
1&1&1&0&1&0&0\\
\ebbb
\parbox{7cm}{ 
\1\6\4\4\6\2\7\3\3\7\2\6\4\4\6\1\eeee
\6\4\4\6\1\7\5\5\7\1\6\4\4\6\2\7\eeee
\4\4\6\2\7\3\3\7\2\6\4\4\6\1\7\5\eeee
\4\6\1\7\5\5\7\1\6\4\4\6\2\7\3\3\eeee
\6\2\7\3\3\7\2\6\4\4\6\1\7\5\5\7\eeee
\1\7\5\5\7\1\6\4\4\6\2\7\3\3\7\2\eeee
\7\3\3\7\2\6\4\4\6\1\7\5\5\7\1\6\eeee
\5\5\7\1\6\4\4\6\2\7\3\3\7\2\6\4\eeee
\3\7\2\6\4\4\6\1\7\5\5\7\1\6\4\4\eeee
\7\1\6\4\4\6\2\7\3\3\7\2\6\4\4\6\eeee
\2\6\4\4\6\1\7\5\5\7\1\6\4\4\6\2\eeee
\6\4\4\6\2\7\3\3\7\2\6\4\4\6\1\7\eeee
} 

\baaa
7-48
\eaaa
\bbbb
0&0&0&0&0&2&2\\
0&0&0&0&1&1&2\\
0&0&0&0&1&1&2\\
0&0&0&0&2&0&2\\
0&1&1&2&0&0&0\\
2&1&1&0&0&0&0\\
1&1&1&1&0&0&0\\
\ebbb
\parbox{7cm}{ 
\1\6\2\7\1\6\2\7\1\6\2\7\1\6\2\7\eeee
\6\1\7\3\6\1\7\3\6\1\7\3\6\1\7\3\eeee
\2\7\4\5\2\7\4\5\2\7\4\5\2\7\4\5\eeee
\7\3\5\4\7\3\5\4\7\3\5\4\7\3\5\4\eeee
\1\6\2\7\1\6\2\7\1\6\2\7\1\6\2\7\eeee
\6\1\7\3\6\1\7\3\6\1\7\3\6\1\7\3\eeee
\2\7\4\5\2\7\4\5\2\7\4\5\2\7\4\5\eeee
\7\3\5\4\7\3\5\4\7\3\5\4\7\3\5\4\eeee
\1\6\2\7\1\6\2\7\1\6\2\7\1\6\2\7\eeee
\6\1\7\3\6\1\7\3\6\1\7\3\6\1\7\3\eeee
\2\7\4\5\2\7\4\5\2\7\4\5\2\7\4\5\eeee
\7\3\5\4\7\3\5\4\7\3\5\4\7\3\5\4\eeee
} 

\baaa
7-49
\eaaa
\bbbb
0&0&0&0&0&2&2\\
0&0&0&0&1&1&2\\
0&0&0&0&2&0&2\\
0&0&0&0&2&2&0\\
0&1&1&2&0&0&0\\
1&1&0&2&0&0&0\\
1&2&1&0&0&0&0\\
\ebbb
\parbox{7cm}{ 
\1\6\4\6\1\7\2\6\4\5\2\7\3\5\4\5\eeee
\6\4\5\2\7\3\5\4\5\3\7\2\5\4\6\2\eeee
\4\5\3\7\2\5\4\6\2\7\1\6\4\6\1\7\eeee
\6\2\7\1\6\4\6\1\7\2\6\4\5\2\7\3\eeee
\1\7\2\6\4\5\2\7\3\5\4\5\3\7\2\5\eeee
\7\3\5\4\5\3\7\2\5\4\6\2\7\1\6\4\eeee
\2\5\4\6\2\7\1\6\4\6\1\7\2\6\4\5\eeee
\6\4\6\1\7\2\6\4\5\2\7\3\5\4\5\3\eeee
\4\5\2\7\3\5\4\5\3\7\2\5\4\6\2\7\eeee
\5\3\7\2\5\4\6\2\7\1\6\4\6\1\7\2\eeee
\2\7\1\6\4\6\1\7\2\6\4\5\2\7\3\5\eeee
\7\2\6\4\5\2\7\3\5\4\5\3\7\2\5\4\eeee
} 

\baaa
7-50
\eaaa
\bbbb
0&0&0&0&0&2&2\\
0&0&0&0&1&1&2\\
0&0&0&0&2&1&1\\
0&0&0&0&2&2&0\\
0&1&2&1&0&0&0\\
1&1&1&1&0&0&0\\
1&2&1&0&0&0&0\\
\ebbb
\parbox{7cm}{ 
\1\7\2\6\3\5\4\5\3\6\2\7\1\7\2\6\eeee
\6\2\7\1\7\2\6\3\5\4\5\3\6\2\7\1\eeee
\4\5\3\6\2\7\1\7\2\6\3\5\4\5\3\6\eeee
\6\3\5\4\5\3\6\2\7\1\7\2\6\3\5\4\eeee
\1\7\2\6\3\5\4\5\3\6\2\7\1\7\2\6\eeee
\6\2\7\1\7\2\6\3\5\4\5\3\6\2\7\1\eeee
\4\5\3\6\2\7\1\7\2\6\3\5\4\5\3\6\eeee
\6\3\5\4\5\3\6\2\7\1\7\2\6\3\5\4\eeee
\1\7\2\6\3\5\4\5\3\6\2\7\1\7\2\6\eeee
\6\2\7\1\7\2\6\3\5\4\5\3\6\2\7\1\eeee
\4\5\3\6\2\7\1\7\2\6\3\5\4\5\3\6\eeee
\6\3\5\4\5\3\6\2\7\1\7\2\6\3\5\4\eeee
} 

\baaa
7-51
\eaaa
\bbbb
0&0&0&0&0&2&2\\
0&0&0&0&1&1&2\\
0&0&0&1&1&1&1\\
0&0&1&1&1&1&0\\
0&1&1&1&1&0&0\\
1&1&1&1&0&0&0\\
1&2&1&0&0&0&0\\
\ebbb
\parbox{7cm}{ 
\1\7\2\6\3\5\4\4\5\3\6\2\7\1\7\2\eeee
\6\2\7\1\7\2\6\3\5\4\4\5\3\6\2\7\eeee
\4\5\3\6\2\7\1\7\2\6\3\5\4\4\5\3\eeee
\3\5\4\4\5\3\6\2\7\1\7\2\6\3\5\4\eeee
\7\2\6\3\5\4\4\5\3\6\2\7\1\7\2\6\eeee
\2\7\1\7\2\6\3\5\4\4\5\3\6\2\7\1\eeee
\5\3\6\2\7\1\7\2\6\3\5\4\4\5\3\6\eeee
\5\4\4\5\3\6\2\7\1\7\2\6\3\5\4\4\eeee
\2\6\3\5\4\4\5\3\6\2\7\1\7\2\6\3\eeee
\7\1\7\2\6\3\5\4\4\5\3\6\2\7\1\7\eeee
\3\6\2\7\1\7\2\6\3\5\4\4\5\3\6\2\eeee
\4\4\5\3\6\2\7\1\7\2\6\3\5\4\4\5\eeee
} 

\baaa
7-52
\eaaa
\bbbb
0&0&0&0&0&2&2\\
0&0&0&0&1&1&2\\
0&0&0&1&2&0&1\\
0&0&1&2&1&0&0\\
0&1&2&1&0&0&0\\
2&2&0&0&0&0&0\\
1&2&1&0&0&0&0\\
\ebbb
\parbox{7cm}{ 
\1\7\3\4\5\2\6\2\5\4\3\7\1\7\3\4\eeee
\6\2\5\4\3\7\1\7\3\4\5\2\6\2\5\4\eeee
\1\7\3\4\5\2\6\2\5\4\3\7\1\7\3\4\eeee
\6\2\5\4\3\7\1\7\3\4\5\2\6\2\5\4\eeee
\1\7\3\4\5\2\6\2\5\4\3\7\1\7\3\4\eeee
\6\2\5\4\3\7\1\7\3\4\5\2\6\2\5\4\eeee
\1\7\3\4\5\2\6\2\5\4\3\7\1\7\3\4\eeee
\6\2\5\4\3\7\1\7\3\4\5\2\6\2\5\4\eeee
\1\7\3\4\5\2\6\2\5\4\3\7\1\7\3\4\eeee
\6\2\5\4\3\7\1\7\3\4\5\2\6\2\5\4\eeee
\1\7\3\4\5\2\6\2\5\4\3\7\1\7\3\4\eeee
\6\2\5\4\3\7\1\7\3\4\5\2\6\2\5\4\eeee
} 

\baaa
7-53
\eaaa
\bbbb
0&0&0&0&0&2&2\\
0&0&0&0&2&0&2\\
0&0&0&1&1&1&1\\
0&0&1&0&1&1&1\\
0&1&1&1&0&1&0\\
1&0&1&1&1&0&0\\
1&1&1&1&0&0&0\\
\ebbb
\parbox{7cm}{ 
\1\7\2\7\1\7\2\7\1\7\2\7\1\7\2\7\eeee
\6\3\5\4\6\3\5\4\6\3\5\4\6\3\5\4\eeee
\5\4\6\3\5\4\6\3\5\4\6\3\5\4\6\3\eeee
\2\7\1\7\2\7\1\7\2\7\1\7\2\7\1\7\eeee
\5\3\6\4\5\3\6\4\5\3\6\4\5\3\6\4\eeee
\6\4\5\3\6\4\5\3\6\4\5\3\6\4\5\3\eeee
\1\7\2\7\1\7\2\7\1\7\2\7\1\7\2\7\eeee
\6\3\5\4\6\3\5\4\6\3\5\4\6\3\5\4\eeee
\5\4\6\3\5\4\6\3\5\4\6\3\5\4\6\3\eeee
\2\7\1\7\2\7\1\7\2\7\1\7\2\7\1\7\eeee
\5\3\6\4\5\3\6\4\5\3\6\4\5\3\6\4\eeee
\6\4\5\3\6\4\5\3\6\4\5\3\6\4\5\3\eeee
} 

\baaa
7-54
\eaaa
\bbbb
0&0&0&0&0&2&2\\
0&0&0&0&2&0&2\\
0&0&0&1&1&1&1\\
0&0&1&0&1&1&1\\
0&1&1&1&1&0&0\\
1&0&1&1&0&1&0\\
1&1&1&1&0&0&0\\
\ebbb
\parbox{7cm}{ 
\1\7\2\7\1\7\2\7\1\7\2\7\1\7\2\7\eeee
\6\3\5\4\6\3\5\4\6\3\5\4\6\3\5\4\eeee
\6\4\5\3\6\4\5\3\6\4\5\3\6\4\5\3\eeee
\1\7\2\7\1\7\2\7\1\7\2\7\1\7\2\7\eeee
\6\3\5\4\6\3\5\4\6\3\5\4\6\3\5\4\eeee
\6\4\5\3\6\4\5\3\6\4\5\3\6\4\5\3\eeee
\1\7\2\7\1\7\2\7\1\7\2\7\1\7\2\7\eeee
\6\3\5\4\6\3\5\4\6\3\5\4\6\3\5\4\eeee
\6\4\5\3\6\4\5\3\6\4\5\3\6\4\5\3\eeee
\1\7\2\7\1\7\2\7\1\7\2\7\1\7\2\7\eeee
\6\3\5\4\6\3\5\4\6\3\5\4\6\3\5\4\eeee
\6\4\5\3\6\4\5\3\6\4\5\3\6\4\5\3\eeee
} 

\baaa
7-55
\eaaa
\bbbb
0&0&0&0&0&2&2\\
0&0&0&0&2&0&2\\
0&0&0&2&0&2&0\\
0&0&2&0&2&0&0\\
0&2&0&2&0&0&0\\
2&0&2&0&0&0&0\\
2&2&0&0&0&0&0\\
\ebbb
\parbox{7cm}{ 
\1\6\3\4\5\2\7\1\6\3\4\5\2\7\1\6\eeee
\6\3\4\5\2\7\1\6\3\4\5\2\7\1\6\3\eeee
\3\4\5\2\7\1\6\3\4\5\2\7\1\6\3\4\eeee
\4\5\2\7\1\6\3\4\5\2\7\1\6\3\4\5\eeee
\5\2\7\1\6\3\4\5\2\7\1\6\3\4\5\2\eeee
\2\7\1\6\3\4\5\2\7\1\6\3\4\5\2\7\eeee
\7\1\6\3\4\5\2\7\1\6\3\4\5\2\7\1\eeee
\1\6\3\4\5\2\7\1\6\3\4\5\2\7\1\6\eeee
\6\3\4\5\2\7\1\6\3\4\5\2\7\1\6\3\eeee
\3\4\5\2\7\1\6\3\4\5\2\7\1\6\3\4\eeee
\4\5\2\7\1\6\3\4\5\2\7\1\6\3\4\5\eeee
\5\2\7\1\6\3\4\5\2\7\1\6\3\4\5\2\eeee
} 

\baaa
7-56
\eaaa
\bbbb
0&0&0&0&0&2&2\\
0&0&0&0&2&0&2\\
0&0&0&2&0&2&0\\
0&0&2&2&0&0&0\\
0&2&0&0&2&0&0\\
2&0&2&0&0&0&0\\
2&2&0&0&0&0&0\\
\ebbb
\parbox{7cm}{ 
\1\6\3\4\4\3\6\1\7\2\5\5\2\7\1\6\eeee
\6\3\4\4\3\6\1\7\2\5\5\2\7\1\6\3\eeee
\3\4\4\3\6\1\7\2\5\5\2\7\1\6\3\4\eeee
\4\4\3\6\1\7\2\5\5\2\7\1\6\3\4\4\eeee
\4\3\6\1\7\2\5\5\2\7\1\6\3\4\4\3\eeee
\3\6\1\7\2\5\5\2\7\1\6\3\4\4\3\6\eeee
\6\1\7\2\5\5\2\7\1\6\3\4\4\3\6\1\eeee
\1\7\2\5\5\2\7\1\6\3\4\4\3\6\1\7\eeee
\7\2\5\5\2\7\1\6\3\4\4\3\6\1\7\2\eeee
\2\5\5\2\7\1\6\3\4\4\3\6\1\7\2\5\eeee
\5\5\2\7\1\6\3\4\4\3\6\1\7\2\5\5\eeee
\5\2\7\1\6\3\4\4\3\6\1\7\2\5\5\2\eeee
} 

\baaa
7-57
\eaaa
\bbbb
0&0&0&0&0&2&2\\
0&0&0&0&2&0&2\\
0&0&1&0&1&1&1\\
0&0&0&1&1&1&1\\
0&1&1&1&0&1&0\\
1&0&1&1&1&0&0\\
1&1&1&1&0&0&0\\
\ebbb
\parbox{7cm}{ 
\1\7\2\7\1\7\2\7\1\7\2\7\1\7\2\7\eeee
\6\3\5\4\6\3\5\4\6\3\5\4\6\3\5\4\eeee
\5\3\6\4\5\3\6\4\5\3\6\4\5\3\6\4\eeee
\2\7\1\7\2\7\1\7\2\7\1\7\2\7\1\7\eeee
\5\4\6\3\5\4\6\3\5\4\6\3\5\4\6\3\eeee
\6\4\5\3\6\4\5\3\6\4\5\3\6\4\5\3\eeee
\1\7\2\7\1\7\2\7\1\7\2\7\1\7\2\7\eeee
\6\3\5\4\6\3\5\4\6\3\5\4\6\3\5\4\eeee
\5\3\6\4\5\3\6\4\5\3\6\4\5\3\6\4\eeee
\2\7\1\7\2\7\1\7\2\7\1\7\2\7\1\7\eeee
\5\4\6\3\5\4\6\3\5\4\6\3\5\4\6\3\eeee
\6\4\5\3\6\4\5\3\6\4\5\3\6\4\5\3\eeee
} 

\baaa
7-58
\eaaa
\bbbb
0&0&0&0&0&2&2\\
0&0&0&0&2&0&2\\
0&0&1&0&1&1&1\\
0&0&0&1&1&1&1\\
0&1&1&1&1&0&0\\
1&0&1&1&0&1&0\\
1&1&1&1&0&0&0\\
\ebbb
\parbox{7cm}{ 
\1\7\2\7\1\7\2\7\1\7\2\7\1\7\2\7\eeee
\6\3\5\4\6\3\5\4\6\3\5\4\6\3\5\4\eeee
\6\3\5\4\6\3\5\4\6\3\5\4\6\3\5\4\eeee
\1\7\2\7\1\7\2\7\1\7\2\7\1\7\2\7\eeee
\6\4\5\3\6\4\5\3\6\4\5\3\6\4\5\3\eeee
\6\4\5\3\6\4\5\3\6\4\5\3\6\4\5\3\eeee
\1\7\2\7\1\7\2\7\1\7\2\7\1\7\2\7\eeee
\6\3\5\4\6\3\5\4\6\3\5\4\6\3\5\4\eeee
\6\3\5\4\6\3\5\4\6\3\5\4\6\3\5\4\eeee
\1\7\2\7\1\7\2\7\1\7\2\7\1\7\2\7\eeee
\6\4\5\3\6\4\5\3\6\4\5\3\6\4\5\3\eeee
\6\4\5\3\6\4\5\3\6\4\5\3\6\4\5\3\eeee
} 

\baaa
7-59
\eaaa
\bbbb
0&0&0&0&0&2&2\\
0&0&0&0&2&0&2\\
0&0&1&1&0&2&0\\
0&0&1&1&0&2&0\\
0&2&0&0&2&0&0\\
2&0&1&1&0&0&0\\
2&2&0&0&0&0&0\\
\ebbb
\parbox{7cm}{ 
\1\6\4\3\6\1\7\2\5\5\2\7\1\6\4\3\eeee
\6\3\4\6\1\7\2\5\5\2\7\1\6\3\4\6\eeee
\4\3\6\1\7\2\5\5\2\7\1\6\4\3\6\1\eeee
\4\6\1\7\2\5\5\2\7\1\6\3\4\6\1\7\eeee
\6\1\7\2\5\5\2\7\1\6\4\3\6\1\7\2\eeee
\1\7\2\5\5\2\7\1\6\3\4\6\1\7\2\5\eeee
\7\2\5\5\2\7\1\6\4\3\6\1\7\2\5\5\eeee
\2\5\5\2\7\1\6\3\4\6\1\7\2\5\5\2\eeee
\5\5\2\7\1\6\4\3\6\1\7\2\5\5\2\7\eeee
\5\2\7\1\6\3\4\6\1\7\2\5\5\2\7\1\eeee
\2\7\1\6\4\3\6\1\7\2\5\5\2\7\1\6\eeee
\7\1\6\3\4\6\1\7\2\5\5\2\7\1\6\3\eeee
} 

\baaa
7-60
\eaaa
\bbbb
0&0&0&0&0&2&2\\
0&0&0&0&2&1&1\\
0&0&0&0&2&1&1\\
0&0&0&2&2&0&0\\
0&1&1&2&0&0&0\\
2&1&1&0&0&0&0\\
2&1&1&0&0&0&0\\
\ebbb
\parbox{7cm}{ 
\1\6\2\5\4\4\5\2\6\1\6\2\5\4\4\5\eeee
\6\1\7\3\5\4\4\5\3\7\1\7\3\5\4\4\eeee
\2\7\1\6\2\5\4\4\5\2\6\1\6\2\5\4\eeee
\5\3\6\1\7\3\5\4\4\5\3\7\1\7\3\5\eeee
\4\5\2\7\1\6\2\5\4\4\5\2\6\1\6\2\eeee
\4\4\5\3\6\1\7\3\5\4\4\5\3\7\1\7\eeee
\5\4\4\5\2\7\1\6\2\5\4\4\5\2\6\1\eeee
\2\5\4\4\5\3\6\1\7\3\5\4\4\5\3\7\eeee
\6\3\5\4\4\5\2\7\1\6\2\5\4\4\5\2\eeee
\1\7\2\5\4\4\5\3\6\1\7\3\5\4\4\5\eeee
\6\1\6\3\5\4\4\5\2\7\1\6\2\5\4\4\eeee
\2\7\1\7\2\5\4\4\5\3\6\1\7\3\5\4\eeee
} 

\baaa
7-61
\eaaa
\bbbb
0&0&0&0&0&2&2\\
0&0&0&1&1&0&2\\
0&0&0&1&1&0&2\\
0&1&1&0&0&2&0\\
0&1&1&0&0&2&0\\
2&0&0&1&1&0&0\\
2&1&1&0&0&0&0\\
\ebbb
\parbox{7cm}{ 
\1\6\5\2\7\1\6\4\2\7\1\6\4\3\7\1\eeee
\6\4\3\7\1\6\5\3\7\1\6\5\2\7\1\6\eeee
\5\2\7\1\6\4\2\7\1\6\4\3\7\1\6\5\eeee
\3\7\1\6\5\3\7\1\6\5\2\7\1\6\4\2\eeee
\7\1\6\4\2\7\1\6\4\3\7\1\6\5\3\7\eeee
\1\6\5\3\7\1\6\5\2\7\1\6\4\2\7\1\eeee
\6\4\2\7\1\6\4\3\7\1\6\5\3\7\1\6\eeee
\5\3\7\1\6\5\2\7\1\6\4\2\7\1\6\4\eeee
\2\7\1\6\4\3\7\1\6\5\3\7\1\6\5\2\eeee
\7\1\6\5\2\7\1\6\4\2\7\1\6\4\3\7\eeee
\1\6\4\3\7\1\6\5\3\7\1\6\5\2\7\1\eeee
\6\5\2\7\1\6\4\2\7\1\6\4\3\7\1\6\eeee
} 

\baaa
7-62
\eaaa
\bbbb
0&0&0&0&0&2&2\\
0&0&0&1&1&0&2\\
0&0&0&1&1&0&2\\
0&1&1&1&1&0&0\\
0&1&1&1&1&0&0\\
2&0&0&0&0&2&0\\
2&1&1&0&0&0&0\\
\ebbb
\parbox{7cm}{ 
\1\6\6\1\7\3\5\4\2\7\1\6\6\1\7\2\eeee
\6\6\1\7\2\4\5\3\7\1\6\6\1\7\3\5\eeee
\6\1\7\3\5\4\2\7\1\6\6\1\7\2\4\5\eeee
\1\7\2\4\5\3\7\1\6\6\1\7\3\5\4\2\eeee
\7\3\5\4\2\7\1\6\6\1\7\2\4\5\3\7\eeee
\2\4\5\3\7\1\6\6\1\7\3\5\4\2\7\1\eeee
\5\4\2\7\1\6\6\1\7\2\4\5\3\7\1\6\eeee
\5\3\7\1\6\6\1\7\3\5\4\2\7\1\6\6\eeee
\2\7\1\6\6\1\7\2\4\5\3\7\1\6\6\1\eeee
\7\1\6\6\1\7\3\5\4\2\7\1\6\6\1\7\eeee
\1\6\6\1\7\2\4\5\3\7\1\6\6\1\7\3\eeee
\6\6\1\7\3\5\4\2\7\1\6\6\1\7\2\4\eeee
} 

\baaa
7-63
\eaaa
\bbbb
0&0&0&0&0&2&2\\
0&0&0&1&1&1&1\\
0&0&0&1&1&1&1\\
0&1&1&1&1&0&0\\
0&1&1&1&1&0&0\\
2&1&1&0&0&0&0\\
2&1&1&0&0&0&0\\
\ebbb
\parbox{7cm}{ 
\1\6\2\4\4\2\6\1\6\2\4\4\2\6\1\6\eeee
\6\1\7\3\5\5\3\7\1\7\3\5\5\3\7\1\eeee
\2\7\1\6\2\4\4\2\6\1\6\2\4\4\2\6\eeee
\4\3\6\1\7\3\5\5\3\7\1\7\3\5\5\3\eeee
\4\5\2\7\1\6\2\4\4\2\6\1\6\2\4\4\eeee
\2\5\4\3\6\1\7\3\5\5\3\7\1\7\3\5\eeee
\6\3\4\5\2\7\1\6\2\4\4\2\6\1\6\2\eeee
\1\7\2\5\4\3\6\1\7\3\5\5\3\7\1\7\eeee
\6\1\6\3\4\5\2\7\1\6\2\4\4\2\6\1\eeee
\2\7\1\7\2\5\4\3\6\1\7\3\5\5\3\7\eeee
\4\3\6\1\6\3\4\5\2\7\1\6\2\4\4\2\eeee
\4\5\2\7\1\7\2\5\4\3\6\1\7\3\5\5\eeee
} 

\baaa
7-64
\eaaa
\bbbb
0&0&0&0&0&2&2\\
0&0&0&1&1&1&1\\
0&0&0&2&2&0&0\\
0&1&1&0&0&1&1\\
0&1&1&0&1&0&1\\
1&1&0&1&0&1&0\\
1&1&0&1&1&0&0\\
\ebbb
\parbox{7cm}{ 
\1\7\5\2\6\4\3\4\6\2\5\7\1\7\5\2\eeee
\6\2\5\7\1\7\5\2\6\4\3\4\6\2\5\7\eeee
\6\4\3\4\6\2\5\7\1\7\5\2\6\4\3\4\eeee
\1\7\5\2\6\4\3\4\6\2\5\7\1\7\5\2\eeee
\6\2\5\7\1\7\5\2\6\4\3\4\6\2\5\7\eeee
\6\4\3\4\6\2\5\7\1\7\5\2\6\4\3\4\eeee
\1\7\5\2\6\4\3\4\6\2\5\7\1\7\5\2\eeee
\6\2\5\7\1\7\5\2\6\4\3\4\6\2\5\7\eeee
\6\4\3\4\6\2\5\7\1\7\5\2\6\4\3\4\eeee
\1\7\5\2\6\4\3\4\6\2\5\7\1\7\5\2\eeee
\6\2\5\7\1\7\5\2\6\4\3\4\6\2\5\7\eeee
\6\4\3\4\6\2\5\7\1\7\5\2\6\4\3\4\eeee
} 

\baaa
7-65
\eaaa
\bbbb
0&0&0&0&0&2&2\\
0&0&0&1&1&1&1\\
0&0&0&2&2&0&0\\
0&1&1&0&1&0&1\\
0&1&1&1&1&0&0\\
1&1&0&0&0&1&1\\
1&1&0&1&0&1&0\\
\ebbb
\parbox{7cm}{ 
\1\7\4\3\4\7\1\7\4\3\4\7\1\7\4\2\eeee
\6\2\5\5\2\6\6\2\5\5\2\6\6\2\5\5\eeee
\7\4\3\4\7\1\7\4\3\4\7\1\7\4\3\4\eeee
\2\5\5\2\6\6\2\5\5\2\6\6\2\5\5\2\eeee
\4\3\4\7\1\7\4\3\4\7\1\7\4\3\4\7\eeee
\5\5\2\6\6\2\5\5\2\6\6\2\5\5\2\6\eeee
\3\4\7\1\7\4\3\4\7\1\7\4\3\4\7\1\eeee
\5\2\6\6\2\5\5\2\6\6\2\5\5\2\6\6\eeee
\4\7\1\7\4\3\4\7\1\7\4\3\4\7\1\7\eeee
\2\6\6\2\5\5\2\6\6\2\5\5\2\6\6\2\eeee
\7\1\7\4\3\4\7\1\7\4\3\4\7\1\7\4\eeee
\6\6\2\5\5\2\6\6\2\5\5\2\6\6\2\5\eeee
} 

\baaa
7-66
\eaaa
\bbbb
0&0&0&0&0&2&2\\
0&0&0&1&1&1&1\\
0&0&1&1&2&0&0\\
0&1&1&0&1&0&1\\
0&1&2&1&0&0&0\\
1&1&0&0&0&1&1\\
1&1&0&1&0&1&0\\
\ebbb
\parbox{7cm}{ 
\1\7\4\3\5\2\6\6\2\5\3\4\7\1\7\4\eeee
\6\2\5\3\4\7\1\7\4\3\5\2\6\6\2\5\eeee
\7\4\3\5\2\6\6\2\5\3\4\7\1\7\4\3\eeee
\2\5\3\4\7\1\7\4\3\5\2\6\6\2\5\3\eeee
\4\3\5\2\6\6\2\5\3\4\7\1\7\4\3\5\eeee
\5\3\4\7\1\7\4\3\5\2\6\6\2\5\3\4\eeee
\3\5\2\6\6\2\5\3\4\7\1\7\4\3\5\2\eeee
\3\4\7\1\7\4\3\5\2\6\6\2\5\3\4\7\eeee
\5\2\6\6\2\5\3\4\7\1\7\4\3\5\2\6\eeee
\4\7\1\7\4\3\5\2\6\6\2\5\3\4\7\1\eeee
\2\6\6\2\5\3\4\7\1\7\4\3\5\2\6\6\eeee
\7\1\7\4\3\5\2\6\6\2\5\3\4\7\1\7\eeee
} 

\baaa
7-67
\eaaa
\bbbb
0&0&0&0&0&2&2\\
0&0&0&1&1&1&1\\
0&0&2&0&1&0&1\\
0&1&0&2&0&1&0\\
0&1&1&0&0&1&1\\
1&1&0&1&1&0&0\\
1&1&1&0&1&0&0\\
\ebbb
\parbox{7cm}{ 
\1\7\3\3\5\6\2\5\7\2\4\4\6\1\6\4\eeee
\6\2\5\7\2\4\4\6\1\6\4\4\2\7\5\2\eeee
\4\4\6\1\6\4\4\2\7\5\2\6\5\3\3\7\eeee
\4\4\2\7\5\2\6\5\3\3\7\1\7\3\3\5\eeee
\2\6\5\3\3\7\1\7\3\3\5\6\2\5\7\2\eeee
\7\1\7\3\3\5\6\2\5\7\2\4\4\6\1\6\eeee
\5\6\2\5\7\2\4\4\6\1\6\4\4\2\7\5\eeee
\2\4\4\6\1\6\4\4\2\7\5\2\6\5\3\3\eeee
\6\4\4\2\7\5\2\6\5\3\3\7\1\7\3\3\eeee
\5\2\6\5\3\3\7\1\7\3\3\5\6\2\5\7\eeee
\3\7\1\7\3\3\5\6\2\5\7\2\4\4\6\1\eeee
\3\5\6\2\5\7\2\4\4\6\1\6\4\4\2\7\eeee
} 

\baaa
7-68
\eaaa
\bbbb
0&0&0&0&0&2&2\\
0&0&0&1&2&0&1\\
0&0&1&0&1&1&1\\
0&1&0&1&0&1&1\\
0&2&1&0&0&1&0\\
1&0&1&1&1&0&0\\
1&1&1&1&0&0&0\\
\ebbb
\parbox{7cm}{ 
\1\7\2\5\3\6\4\4\6\3\5\2\7\1\7\2\eeee
\6\3\5\2\7\1\7\2\5\3\6\4\4\6\3\5\eeee
\5\3\6\4\4\6\3\5\2\7\1\7\2\5\3\6\eeee
\2\7\1\7\2\5\3\6\4\4\6\3\5\2\7\1\eeee
\4\4\6\3\5\2\7\1\7\2\5\3\6\4\4\6\eeee
\7\2\5\3\6\4\4\6\3\5\2\7\1\7\2\5\eeee
\3\5\2\7\1\7\2\5\3\6\4\4\6\3\5\2\eeee
\3\6\4\4\6\3\5\2\7\1\7\2\5\3\6\4\eeee
\7\1\7\2\5\3\6\4\4\6\3\5\2\7\1\7\eeee
\4\6\3\5\2\7\1\7\2\5\3\6\4\4\6\3\eeee
\2\5\3\6\4\4\6\3\5\2\7\1\7\2\5\3\eeee
\5\2\7\1\7\2\5\3\6\4\4\6\3\5\2\7\eeee
} 

\baaa
7-69
\eaaa
\bbbb
0&0&0&0&0&2&2\\
0&1&0&0&1&0&2\\
0&0&1&1&0&2&0\\
0&0&1&1&0&2&0\\
0&1&0&0&1&0&2\\
2&0&1&1&0&0&0\\
2&1&0&0&1&0&0\\
\ebbb
\parbox{7cm}{ 
\1\6\4\3\6\1\7\2\2\7\1\6\3\4\6\1\eeee
\6\3\4\6\1\7\5\5\7\1\6\4\3\6\1\7\eeee
\4\3\6\1\7\2\2\7\1\6\3\4\6\1\7\5\eeee
\4\6\1\7\5\5\7\1\6\4\3\6\1\7\2\2\eeee
\6\1\7\2\2\7\1\6\3\4\6\1\7\5\5\7\eeee
\1\7\5\5\7\1\6\4\3\6\1\7\2\2\7\1\eeee
\7\2\2\7\1\6\3\4\6\1\7\5\5\7\1\6\eeee
\5\5\7\1\6\4\3\6\1\7\2\2\7\1\6\3\eeee
\2\7\1\6\3\4\6\1\7\5\5\7\1\6\4\3\eeee
\7\1\6\4\3\6\1\7\2\2\7\1\6\3\4\6\eeee
\1\6\3\4\6\1\7\5\5\7\1\6\4\3\6\1\eeee
\6\4\3\6\1\7\2\2\7\1\6\3\4\6\1\7\eeee
} 

\baaa
7-70
\eaaa
\bbbb
0&0&0&0&1&1&2\\
0&0&0&0&1&1&2\\
0&0&0&0&1&1&2\\
0&0&0&0&1&1&2\\
1&1&1&1&0&0&0\\
1&1&1&1&0&0&0\\
1&1&1&1&0&0&0\\
\ebbb
\parbox{7cm}{ 
\1\6\3\5\1\6\3\5\1\6\3\5\1\6\3\5\eeee
\5\2\7\4\7\2\7\4\7\2\7\4\7\2\7\4\eeee
\3\7\1\6\3\5\1\6\3\5\1\6\3\5\1\6\eeee
\6\4\5\2\7\4\7\2\7\4\7\2\7\4\7\2\eeee
\1\7\3\7\1\6\3\5\1\6\3\5\1\6\3\5\eeee
\5\2\6\4\5\2\7\4\7\2\7\4\7\2\7\4\eeee
\3\7\1\7\3\7\1\6\3\5\1\6\3\5\1\6\eeee
\6\4\5\2\6\4\5\2\7\4\7\2\7\4\7\2\eeee
\1\7\3\7\1\7\3\7\1\6\3\5\1\6\3\5\eeee
\5\2\6\4\5\2\6\4\5\2\7\4\7\2\7\4\eeee
\3\7\1\7\3\7\1\7\3\7\1\6\3\5\1\6\eeee
\6\4\5\2\6\4\5\2\6\4\5\2\7\4\7\2\eeee
} 

\baaa
7-71
\eaaa
\bbbb
0&0&0&0&1&1&2\\
0&0&0&0&1&1&2\\
0&0&0&2&1&1&0\\
0&0&2&2&0&0&0\\
1&1&2&0&0&0&0\\
1&1&2&0&0&0&0\\
1&1&0&0&0&0&2\\
\ebbb
\parbox{7cm}{ 
\1\6\3\4\4\3\6\1\7\7\2\5\3\4\4\3\eeee
\5\3\4\4\3\5\2\7\7\1\6\3\4\4\3\6\eeee
\3\4\4\3\6\1\7\7\2\5\3\4\4\3\5\2\eeee
\4\4\3\5\2\7\7\1\6\3\4\4\3\6\1\7\eeee
\4\3\6\1\7\7\2\5\3\4\4\3\5\2\7\7\eeee
\3\5\2\7\7\1\6\3\4\4\3\6\1\7\7\2\eeee
\6\1\7\7\2\5\3\4\4\3\5\2\7\7\1\6\eeee
\2\7\7\1\6\3\4\4\3\6\1\7\7\2\5\3\eeee
\7\7\2\5\3\4\4\3\5\2\7\7\1\6\3\4\eeee
\7\1\6\3\4\4\3\6\1\7\7\2\5\3\4\4\eeee
\2\5\3\4\4\3\5\2\7\7\1\6\3\4\4\3\eeee
\6\3\4\4\3\6\1\7\7\2\5\3\4\4\3\5\eeee
} 

\baaa
7-72
\eaaa
\bbbb
0&0&0&0&1&1&2\\
0&0&0&0&1&1&2\\
0&0&1&1&0&0&2\\
0&0&1&1&0&0&2\\
1&1&0&0&1&1&0\\
1&1&0&0&1&1&0\\
1&1&1&1&0&0&0\\
\ebbb
\parbox{7cm}{ 
\1\6\6\1\7\4\4\7\1\6\6\1\7\4\4\7\eeee
\5\5\2\7\3\3\7\2\5\5\2\7\3\3\7\2\eeee
\6\1\7\4\4\7\1\6\6\1\7\4\4\7\1\6\eeee
\2\7\3\3\7\2\5\5\2\7\3\3\7\2\5\5\eeee
\7\4\4\7\1\6\6\1\7\4\4\7\1\6\6\1\eeee
\3\3\7\2\5\5\2\7\3\3\7\2\5\5\2\7\eeee
\4\7\1\6\6\1\7\4\4\7\1\6\6\1\7\4\eeee
\7\2\5\5\2\7\3\3\7\2\5\5\2\7\3\3\eeee
\1\6\6\1\7\4\4\7\1\6\6\1\7\4\4\7\eeee
\5\5\2\7\3\3\7\2\5\5\2\7\3\3\7\2\eeee
\6\1\7\4\4\7\1\6\6\1\7\4\4\7\1\6\eeee
\2\7\3\3\7\2\5\5\2\7\3\3\7\2\5\5\eeee
} 

\baaa
7-73
\eaaa
\bbbb
0&0&0&0&1&1&2\\
0&0&0&0&1&1&2\\
0&0&1&1&1&1&0\\
0&0&1&1&1&1&0\\
1&1&1&1&0&0&0\\
1&1&1&1&0&0&0\\
1&1&0&0&0&0&2\\
\ebbb
\parbox{7cm}{ 
\1\6\4\3\5\2\7\7\1\6\4\3\5\2\7\7\eeee
\5\3\4\6\1\7\7\2\5\3\4\6\1\7\7\2\eeee
\4\3\5\2\7\7\1\6\4\3\5\2\7\7\1\6\eeee
\4\6\1\7\7\2\5\3\4\6\1\7\7\2\5\3\eeee
\5\2\7\7\1\6\4\3\5\2\7\7\1\6\4\3\eeee
\1\7\7\2\5\3\4\6\1\7\7\2\5\3\4\6\eeee
\7\7\1\6\4\3\5\2\7\7\1\6\4\3\5\2\eeee
\7\2\5\3\4\6\1\7\7\2\5\3\4\6\1\7\eeee
\1\6\4\3\5\2\7\7\1\6\4\3\5\2\7\7\eeee
\5\3\4\6\1\7\7\2\5\3\4\6\1\7\7\2\eeee
\4\3\5\2\7\7\1\6\4\3\5\2\7\7\1\6\eeee
\4\6\1\7\7\2\5\3\4\6\1\7\7\2\5\3\eeee
} 

\baaa
7-74
\eaaa
\bbbb
0&0&0&0&1&1&2\\
0&0&0&1&0&2&1\\
0&0&1&0&2&0&1\\
0&1&0&2&0&1&0\\
1&0&2&0&1&0&0\\
1&2&0&1&0&0&0\\
2&1&1&0&0&0&0\\
\ebbb
\parbox{7cm}{ 
\1\7\1\7\1\7\1\7\1\7\1\7\1\7\1\7\eeee
\5\3\5\3\5\3\5\3\5\3\5\3\5\3\5\3\eeee
\5\3\5\3\5\3\5\3\5\3\5\3\5\3\5\3\eeee
\1\7\1\7\1\7\1\7\1\7\1\7\1\7\1\7\eeee
\6\2\6\2\6\2\6\2\6\2\6\2\6\2\6\2\eeee
\4\4\4\4\4\4\4\4\4\4\4\4\4\4\4\4\eeee
\2\6\2\6\2\6\2\6\2\6\2\6\2\6\2\6\eeee
\7\1\7\1\7\1\7\1\7\1\7\1\7\1\7\1\eeee
\3\5\3\5\3\5\3\5\3\5\3\5\3\5\3\5\eeee
\3\5\3\5\3\5\3\5\3\5\3\5\3\5\3\5\eeee
\7\1\7\1\7\1\7\1\7\1\7\1\7\1\7\1\eeee
\2\6\2\6\2\6\2\6\2\6\2\6\2\6\2\6\eeee
} 

\baaa
7-75
\eaaa
\bbbb
0&0&0&0&1&1&2\\
0&0&0&1&1&2&0\\
0&0&1&1&0&1&1\\
0&1&1&0&1&0&1\\
1&1&0&1&1&0&0\\
1&2&1&0&0&0&0\\
2&0&1&1&0&0&0\\
\ebbb
\parbox{7cm}{ 
\1\6\2\5\4\3\7\1\6\2\5\4\3\7\1\6\eeee
\5\2\6\1\7\3\4\5\2\6\1\7\3\4\5\2\eeee
\5\4\3\7\1\6\2\5\4\3\7\1\6\2\5\4\eeee
\1\7\3\4\5\2\6\1\7\3\4\5\2\6\1\7\eeee
\7\1\6\2\5\4\3\7\1\6\2\5\4\3\7\1\eeee
\4\5\2\6\1\7\3\4\5\2\6\1\7\3\4\5\eeee
\2\5\4\3\7\1\6\2\5\4\3\7\1\6\2\5\eeee
\6\1\7\3\4\5\2\6\1\7\3\4\5\2\6\1\eeee
\3\7\1\6\2\5\4\3\7\1\6\2\5\4\3\7\eeee
\3\4\5\2\6\1\7\3\4\5\2\6\1\7\3\4\eeee
\6\2\5\4\3\7\1\6\2\5\4\3\7\1\6\2\eeee
\2\6\1\7\3\4\5\2\6\1\7\3\4\5\2\6\eeee
} 

\baaa
7-76
\eaaa
\bbbb
0&0&0&0&1&1&2\\
0&0&1&1&1&1&0\\
0&1&0&2&0&1&0\\
0&1&2&1&0&0&0\\
1&1&0&0&0&1&1\\
1&1&1&0&1&0&0\\
2&0&0&0&1&0&1\\
\ebbb
\parbox{7cm}{ 
\1\6\3\4\2\5\7\1\6\3\4\2\5\7\1\6\eeee
\5\2\4\3\6\1\7\5\2\4\3\6\1\7\5\2\eeee
\6\3\4\2\5\7\1\6\3\4\2\5\7\1\6\3\eeee
\2\4\3\6\1\7\5\2\4\3\6\1\7\5\2\4\eeee
\3\4\2\5\7\1\6\3\4\2\5\7\1\6\3\4\eeee
\4\3\6\1\7\5\2\4\3\6\1\7\5\2\4\3\eeee
\4\2\5\7\1\6\3\4\2\5\7\1\6\3\4\2\eeee
\3\6\1\7\5\2\4\3\6\1\7\5\2\4\3\6\eeee
\2\5\7\1\6\3\4\2\5\7\1\6\3\4\2\5\eeee
\6\1\7\5\2\4\3\6\1\7\5\2\4\3\6\1\eeee
\5\7\1\6\3\4\2\5\7\1\6\3\4\2\5\7\eeee
\1\7\5\2\4\3\6\1\7\5\2\4\3\6\1\7\eeee
} 

\baaa
7-77
\eaaa
\bbbb
0&0&0&1&1&1&1\\
0&0&1&0&1&1&1\\
0&1&0&1&0&1&1\\
1&0&1&0&1&0&1\\
1&1&0&1&0&1&0\\
1&1&1&0&1&0&0\\
1&1&1&1&0&0&0\\
\ebbb
\parbox{7cm}{ 
\1\6\3\4\5\2\7\1\6\3\4\5\2\7\1\6\eeee
\4\5\2\7\1\6\3\4\5\2\7\1\6\3\4\5\eeee
\7\1\6\3\4\5\2\7\1\6\3\4\5\2\7\1\eeee
\3\4\5\2\7\1\6\3\4\5\2\7\1\6\3\4\eeee
\2\7\1\6\3\4\5\2\7\1\6\3\4\5\2\7\eeee
\6\3\4\5\2\7\1\6\3\4\5\2\7\1\6\3\eeee
\5\2\7\1\6\3\4\5\2\7\1\6\3\4\5\2\eeee
\1\6\3\4\5\2\7\1\6\3\4\5\2\7\1\6\eeee
\4\5\2\7\1\6\3\4\5\2\7\1\6\3\4\5\eeee
\7\1\6\3\4\5\2\7\1\6\3\4\5\2\7\1\eeee
\3\4\5\2\7\1\6\3\4\5\2\7\1\6\3\4\eeee
\2\7\1\6\3\4\5\2\7\1\6\3\4\5\2\7\eeee
} 

\baaa
7-78
\eaaa
\bbbb
0&0&0&1&1&1&1\\
0&0&1&0&1&1&1\\
0&1&1&0&0&1&1\\
1&0&0&1&1&0&1\\
1&1&0&1&1&0&0\\
1&1&1&0&0&1&0\\
1&1&1&1&0&0&0\\
\ebbb
\parbox{7cm}{ 
\1\5\4\4\5\1\7\2\6\3\3\6\2\7\1\5\eeee
\4\5\1\7\2\6\3\3\6\2\7\1\5\4\4\5\eeee
\7\2\6\3\3\6\2\7\1\5\4\4\5\1\7\2\eeee
\3\3\6\2\7\1\5\4\4\5\1\7\2\6\3\3\eeee
\2\7\1\5\4\4\5\1\7\2\6\3\3\6\2\7\eeee
\5\4\4\5\1\7\2\6\3\3\6\2\7\1\5\4\eeee
\5\1\7\2\6\3\3\6\2\7\1\5\4\4\5\1\eeee
\2\6\3\3\6\2\7\1\5\4\4\5\1\7\2\6\eeee
\3\6\2\7\1\5\4\4\5\1\7\2\6\3\3\6\eeee
\7\1\5\4\4\5\1\7\2\6\3\3\6\2\7\1\eeee
\4\4\5\1\7\2\6\3\3\6\2\7\1\5\4\4\eeee
\1\7\2\6\3\3\6\2\7\1\5\4\4\5\1\7\eeee
} 

\baaa
7-79
\eaaa
\bbbb
1&0&0&0&0&1&2\\
0&1&0&0&1&1&1\\
0&0&1&1&1&0&1\\
0&0&1&2&1&0&0\\
0&1&1&1&1&0&0\\
1&2&0&0&0&1&0\\
1&1&1&0&0&0&1\\
\ebbb
\parbox{7cm}{ 
\1\7\3\4\5\2\6\2\5\4\3\7\1\7\3\4\eeee
\1\7\3\4\5\2\6\2\5\4\3\7\1\7\3\4\eeee
\6\2\5\4\3\7\1\7\3\4\5\2\6\2\5\4\eeee
\6\2\5\4\3\7\1\7\3\4\5\2\6\2\5\4\eeee
\1\7\3\4\5\2\6\2\5\4\3\7\1\7\3\4\eeee
\1\7\3\4\5\2\6\2\5\4\3\7\1\7\3\4\eeee
\6\2\5\4\3\7\1\7\3\4\5\2\6\2\5\4\eeee
\6\2\5\4\3\7\1\7\3\4\5\2\6\2\5\4\eeee
\1\7\3\4\5\2\6\2\5\4\3\7\1\7\3\4\eeee
\1\7\3\4\5\2\6\2\5\4\3\7\1\7\3\4\eeee
\6\2\5\4\3\7\1\7\3\4\5\2\6\2\5\4\eeee
\6\2\5\4\3\7\1\7\3\4\5\2\6\2\5\4\eeee
} 

\baaa
7-80
\eaaa
\bbbb
1&0&0&0&0&1&2\\
0&1&0&0&1&1&1\\
0&0&1&1&1&1&0\\
0&0&1&2&1&0&0\\
0&1&1&1&1&0&0\\
1&1&1&0&0&1&0\\
2&1&0&0&0&0&1\\
\ebbb
\parbox{7cm}{ 
\1\6\3\4\5\2\7\1\6\3\4\5\2\7\1\6\eeee
\1\6\3\4\5\2\7\1\6\3\4\5\2\7\1\6\eeee
\7\2\5\4\3\6\1\7\2\5\4\3\6\1\7\2\eeee
\7\2\5\4\3\6\1\7\2\5\4\3\6\1\7\2\eeee
\1\6\3\4\5\2\7\1\6\3\4\5\2\7\1\6\eeee
\1\6\3\4\5\2\7\1\6\3\4\5\2\7\1\6\eeee
\7\2\5\4\3\6\1\7\2\5\4\3\6\1\7\2\eeee
\7\2\5\4\3\6\1\7\2\5\4\3\6\1\7\2\eeee
\1\6\3\4\5\2\7\1\6\3\4\5\2\7\1\6\eeee
\1\6\3\4\5\2\7\1\6\3\4\5\2\7\1\6\eeee
\7\2\5\4\3\6\1\7\2\5\4\3\6\1\7\2\eeee
\7\2\5\4\3\6\1\7\2\5\4\3\6\1\7\2\eeee
} 

\baaa
7-81
\eaaa
\bbbb
1&0&0&0&1&1&1\\
0&1&0&1&0&1&1\\
0&0&2&0&1&0&1\\
0&1&0&2&0&1&0\\
1&0&1&0&2&0&0\\
1&1&0&1&0&1&0\\
1&1&1&0&0&0&1\\
\ebbb
\parbox{7cm}{ 
\1\5\5\1\6\4\2\7\3\3\7\2\4\6\1\5\eeee
\1\5\5\1\6\4\2\7\3\3\7\2\4\6\1\5\eeee
\7\3\3\7\2\4\6\1\5\5\1\6\4\2\7\3\eeee
\7\3\3\7\2\4\6\1\5\5\1\6\4\2\7\3\eeee
\1\5\5\1\6\4\2\7\3\3\7\2\4\6\1\5\eeee
\1\5\5\1\6\4\2\7\3\3\7\2\4\6\1\5\eeee
\7\3\3\7\2\4\6\1\5\5\1\6\4\2\7\3\eeee
\7\3\3\7\2\4\6\1\5\5\1\6\4\2\7\3\eeee
\1\5\5\1\6\4\2\7\3\3\7\2\4\6\1\5\eeee
\1\5\5\1\6\4\2\7\3\3\7\2\4\6\1\5\eeee
\7\3\3\7\2\4\6\1\5\5\1\6\4\2\7\3\eeee
\7\3\3\7\2\4\6\1\5\5\1\6\4\2\7\3\eeee
} 

\baaa
7-82
\eaaa
\bbbb
2&0&0&0&0&0&2\\
0&2&0&0&0&1&1\\
0&0&2&0&1&1&0\\
0&0&0&2&2&0&0\\
0&0&1&1&2&0&0\\
0&1&1&0&0&2&0\\
1&1&0&0&0&0&2\\
\ebbb
\parbox{7cm}{ 
\1\7\2\6\3\5\4\5\3\6\2\7\1\7\2\6\eeee
\1\7\2\6\3\5\4\5\3\6\2\7\1\7\2\6\eeee
\1\7\2\6\3\5\4\5\3\6\2\7\1\7\2\6\eeee
\1\7\2\6\3\5\4\5\3\6\2\7\1\7\2\6\eeee
\1\7\2\6\3\5\4\5\3\6\2\7\1\7\2\6\eeee
\1\7\2\6\3\5\4\5\3\6\2\7\1\7\2\6\eeee
\1\7\2\6\3\5\4\5\3\6\2\7\1\7\2\6\eeee
\1\7\2\6\3\5\4\5\3\6\2\7\1\7\2\6\eeee
\1\7\2\6\3\5\4\5\3\6\2\7\1\7\2\6\eeee
\1\7\2\6\3\5\4\5\3\6\2\7\1\7\2\6\eeee
\1\7\2\6\3\5\4\5\3\6\2\7\1\7\2\6\eeee
\1\7\2\6\3\5\4\5\3\6\2\7\1\7\2\6\eeee
} 

\baaa
7-83
\eaaa
\bbbb
2&0&0&0&0&0&2\\
0&2&0&0&0&1&1\\
0&0&2&0&1&1&0\\
0&0&0&3&1&0&0\\
0&0&1&1&2&0&0\\
0&1&1&0&0&2&0\\
1&1&0&0&0&0&2\\
\ebbb
\parbox{7cm}{ 
\1\7\2\6\3\5\4\4\5\3\6\2\7\1\7\2\eeee
\1\7\2\6\3\5\4\4\5\3\6\2\7\1\7\2\eeee
\1\7\2\6\3\5\4\4\5\3\6\2\7\1\7\2\eeee
\1\7\2\6\3\5\4\4\5\3\6\2\7\1\7\2\eeee
\1\7\2\6\3\5\4\4\5\3\6\2\7\1\7\2\eeee
\1\7\2\6\3\5\4\4\5\3\6\2\7\1\7\2\eeee
\1\7\2\6\3\5\4\4\5\3\6\2\7\1\7\2\eeee
\1\7\2\6\3\5\4\4\5\3\6\2\7\1\7\2\eeee
\1\7\2\6\3\5\4\4\5\3\6\2\7\1\7\2\eeee
\1\7\2\6\3\5\4\4\5\3\6\2\7\1\7\2\eeee
\1\7\2\6\3\5\4\4\5\3\6\2\7\1\7\2\eeee
\1\7\2\6\3\5\4\4\5\3\6\2\7\1\7\2\eeee
} 

\baaa
7-84
\eaaa
\bbbb
2&0&0&0&0&1&1\\
0&2&0&0&1&0&1\\
0&0&2&1&0&1&0\\
0&0&1&2&1&0&0\\
0&1&0&1&2&0&0\\
1&0&1&0&0&2&0\\
1&1&0&0&0&0&2\\
\ebbb
\parbox{7cm}{ 
\1\6\3\4\5\2\7\1\6\3\4\5\2\7\1\6\eeee
\1\6\3\4\5\2\7\1\6\3\4\5\2\7\1\6\eeee
\1\6\3\4\5\2\7\1\6\3\4\5\2\7\1\6\eeee
\1\6\3\4\5\2\7\1\6\3\4\5\2\7\1\6\eeee
\1\6\3\4\5\2\7\1\6\3\4\5\2\7\1\6\eeee
\1\6\3\4\5\2\7\1\6\3\4\5\2\7\1\6\eeee
\1\6\3\4\5\2\7\1\6\3\4\5\2\7\1\6\eeee
\1\6\3\4\5\2\7\1\6\3\4\5\2\7\1\6\eeee
\1\6\3\4\5\2\7\1\6\3\4\5\2\7\1\6\eeee
\1\6\3\4\5\2\7\1\6\3\4\5\2\7\1\6\eeee
\1\6\3\4\5\2\7\1\6\3\4\5\2\7\1\6\eeee
\1\6\3\4\5\2\7\1\6\3\4\5\2\7\1\6\eeee
} 

\baaa
7-85
\eaaa
\bbbb
2&0&0&0&0&1&1\\
0&2&0&0&1&0&1\\
0&0&2&1&0&1&0\\
0&0&1&3&0&0&0\\
0&1&0&0&3&0&0\\
1&0&1&0&0&2&0\\
1&1&0&0&0&0&2\\
\ebbb
\parbox{7cm}{ 
\1\6\3\4\4\3\6\1\7\2\5\5\2\7\1\6\eeee
\1\6\3\4\4\3\6\1\7\2\5\5\2\7\1\6\eeee
\1\6\3\4\4\3\6\1\7\2\5\5\2\7\1\6\eeee
\1\6\3\4\4\3\6\1\7\2\5\5\2\7\1\6\eeee
\1\6\3\4\4\3\6\1\7\2\5\5\2\7\1\6\eeee
\1\6\3\4\4\3\6\1\7\2\5\5\2\7\1\6\eeee
\1\6\3\4\4\3\6\1\7\2\5\5\2\7\1\6\eeee
\1\6\3\4\4\3\6\1\7\2\5\5\2\7\1\6\eeee
\1\6\3\4\4\3\6\1\7\2\5\5\2\7\1\6\eeee
\1\6\3\4\4\3\6\1\7\2\5\5\2\7\1\6\eeee
\1\6\3\4\4\3\6\1\7\2\5\5\2\7\1\6\eeee
\1\6\3\4\4\3\6\1\7\2\5\5\2\7\1\6\eeee
} 



\baaa
8-1
\eaaa
\bbbb
0&0&0&0&0&0&0&4\\
0&0&0&0&0&0&0&4\\
0&0&0&0&0&0&2&2\\
0&0&0&0&0&0&2&2\\
0&0&0&0&0&0&4&0\\
0&0&0&0&0&0&4&0\\
0&0&1&1&1&1&0&0\\
1&1&1&1&0&0&0&0\\
\ebbb
\parbox{7cm}{ 
\1\8\3\7\5\7\3\8\1\8\3\7\5\7\3\8\eeee
\8\2\8\4\7\6\7\4\8\2\8\4\7\6\7\4\eeee
\3\8\1\8\3\7\5\7\3\8\1\8\3\7\5\7\eeee
\7\4\8\2\8\4\7\6\7\4\8\2\8\4\7\6\eeee
\5\7\3\8\1\8\3\7\5\7\3\8\1\8\3\7\eeee
\7\6\7\4\8\2\8\4\7\6\7\4\8\2\8\4\eeee
\3\7\5\7\3\8\1\8\3\7\5\7\3\8\1\8\eeee
\8\4\7\6\7\4\8\2\8\4\7\6\7\4\8\2\eeee
\1\8\3\7\5\7\3\8\1\8\3\7\5\7\3\8\eeee
\8\2\8\4\7\6\7\4\8\2\8\4\7\6\7\4\eeee
\3\8\1\8\3\7\5\7\3\8\1\8\3\7\5\7\eeee
\7\4\8\2\8\4\7\6\7\4\8\2\8\4\7\6\eeee
} 

\baaa
8-2
\eaaa
\bbbb
0&0&0&0&0&0&0&4\\
0&0&0&0&0&0&0&4\\
0&0&0&0&0&0&2&2\\
0&0&0&0&0&0&2&2\\
0&0&0&0&0&2&2&0\\
0&0&0&0&2&2&0&0\\
0&0&1&1&2&0&0&0\\
1&1&1&1&0&0&0&0\\
\ebbb
\parbox{7cm}{ 
\1\8\3\7\5\6\6\5\7\3\8\1\8\3\7\5\eeee
\8\2\8\4\7\5\6\6\5\7\4\8\2\8\4\7\eeee
\3\8\1\8\3\7\5\6\6\5\7\3\8\1\8\3\eeee
\7\4\8\2\8\4\7\5\6\6\5\7\4\8\2\8\eeee
\5\7\3\8\1\8\3\7\5\6\6\5\7\3\8\1\eeee
\6\5\7\4\8\2\8\4\7\5\6\6\5\7\4\8\eeee
\6\6\5\7\3\8\1\8\3\7\5\6\6\5\7\3\eeee
\5\6\6\5\7\4\8\2\8\4\7\5\6\6\5\7\eeee
\7\5\6\6\5\7\3\8\1\8\3\7\5\6\6\5\eeee
\3\7\5\6\6\5\7\4\8\2\8\4\7\5\6\6\eeee
\8\4\7\5\6\6\5\7\3\8\1\8\3\7\5\6\eeee
\1\8\3\7\5\6\6\5\7\4\8\2\8\4\7\5\eeee
} 

\baaa
8-3
\eaaa
\bbbb
0&0&0&0&0&0&0&4\\
0&0&0&0&0&0&0&4\\
0&0&0&0&0&0&2&2\\
0&0&0&0&0&0&2&2\\
0&0&0&0&0&2&2&0\\
0&0&0&0&4&0&0&0\\
0&0&1&1&2&0&0&0\\
1&1&1&1&0&0&0&0\\
\ebbb
\parbox{7cm}{ 
\1\8\3\7\5\6\5\7\3\8\1\8\3\7\5\6\eeee
\8\2\8\4\7\5\6\5\7\4\8\2\8\4\7\5\eeee
\3\8\1\8\3\7\5\6\5\7\3\8\1\8\3\7\eeee
\7\4\8\2\8\4\7\5\6\5\7\4\8\2\8\4\eeee
\5\7\3\8\1\8\3\7\5\6\5\7\3\8\1\8\eeee
\6\5\7\4\8\2\8\4\7\5\6\5\7\4\8\2\eeee
\5\6\5\7\3\8\1\8\3\7\5\6\5\7\3\8\eeee
\7\5\6\5\7\4\8\2\8\4\7\5\6\5\7\4\eeee
\3\7\5\6\5\7\3\8\1\8\3\7\5\6\5\7\eeee
\8\4\7\5\6\5\7\4\8\2\8\4\7\5\6\5\eeee
\1\8\3\7\5\6\5\7\3\8\1\8\3\7\5\6\eeee
\8\2\8\4\7\5\6\5\7\4\8\2\8\4\7\5\eeee
} 

\baaa
8-4
\eaaa
\bbbb
0&0&0&0&0&0&0&4\\
0&0&0&0&0&0&0&4\\
0&0&0&0&0&0&2&2\\
0&0&0&0&0&0&2&2\\
0&0&0&0&1&1&2&0\\
0&0&0&0&1&1&2&0\\
0&0&1&1&1&1&0&0\\
1&1&1&1&0&0&0&0\\
\ebbb
\parbox{7cm}{ 
\1\8\3\7\5\5\7\3\8\1\8\3\7\5\5\7\eeee
\8\2\8\4\7\6\6\7\4\8\2\8\4\7\6\6\eeee
\3\8\1\8\3\7\5\5\7\3\8\1\8\3\7\5\eeee
\7\4\8\2\8\4\7\6\6\7\4\8\2\8\4\7\eeee
\5\7\3\8\1\8\3\7\5\5\7\3\8\1\8\3\eeee
\5\6\7\4\8\2\8\4\7\6\6\7\4\8\2\8\eeee
\7\6\5\7\3\8\1\8\3\7\5\5\7\3\8\1\eeee
\3\7\5\6\7\4\8\2\8\4\7\6\6\7\4\8\eeee
\8\4\7\6\5\7\3\8\1\8\3\7\5\5\7\3\eeee
\1\8\3\7\5\6\7\4\8\2\8\4\7\6\6\7\eeee
\8\2\8\4\7\6\5\7\3\8\1\8\3\7\5\5\eeee
\3\8\1\8\3\7\5\6\7\4\8\2\8\4\7\6\eeee
} 

\baaa
8-5
\eaaa
\bbbb
0&0&0&0&0&0&0&4\\
0&0&0&0&0&0&0&4\\
0&0&0&0&0&0&2&2\\
0&0&0&0&0&2&2&0\\
0&0&0&0&0&2&2&0\\
0&0&0&1&1&2&0&0\\
0&0&2&1&1&0&0&0\\
1&1&2&0&0&0&0&0\\
\ebbb
\parbox{7cm}{ 
\1\8\3\7\4\6\6\4\7\3\8\1\8\3\7\4\eeee
\8\2\8\3\7\5\6\6\5\7\3\8\2\8\3\7\eeee
\3\8\1\8\3\7\4\6\6\4\7\3\8\1\8\3\eeee
\7\3\8\2\8\3\7\5\6\6\5\7\3\8\2\8\eeee
\4\7\3\8\1\8\3\7\4\6\6\4\7\3\8\1\eeee
\6\5\7\3\8\2\8\3\7\5\6\6\5\7\3\8\eeee
\6\6\4\7\3\8\1\8\3\7\4\6\6\4\7\3\eeee
\4\6\6\5\7\3\8\2\8\3\7\5\6\6\5\7\eeee
\7\5\6\6\4\7\3\8\1\8\3\7\4\6\6\4\eeee
\3\7\4\6\6\5\7\3\8\2\8\3\7\5\6\6\eeee
\8\3\7\5\6\6\4\7\3\8\1\8\3\7\4\6\eeee
\1\8\3\7\4\6\6\5\7\3\8\2\8\3\7\5\eeee
} 

\baaa
8-6
\eaaa
\bbbb
0&0&0&0&0&0&0&4\\
0&0&0&0&0&0&0&4\\
0&0&0&0&0&0&2&2\\
0&0&0&0&0&2&2&0\\
0&0&0&0&0&2&2&0\\
0&0&0&2&2&0&0&0\\
0&0&2&1&1&0&0&0\\
1&1&2&0&0&0&0&0\\
\ebbb
\parbox{7cm}{ 
\1\8\3\7\4\6\4\7\3\8\1\8\3\7\4\6\eeee
\8\2\8\3\7\5\6\5\7\3\8\2\8\3\7\5\eeee
\3\8\1\8\3\7\4\6\4\7\3\8\1\8\3\7\eeee
\7\3\8\2\8\3\7\5\6\5\7\3\8\2\8\3\eeee
\4\7\3\8\1\8\3\7\4\6\4\7\3\8\1\8\eeee
\6\5\7\3\8\2\8\3\7\5\6\5\7\3\8\2\eeee
\4\6\4\7\3\8\1\8\3\7\4\6\4\7\3\8\eeee
\7\5\6\5\7\3\8\2\8\3\7\5\6\5\7\3\eeee
\3\7\4\6\4\7\3\8\1\8\3\7\4\6\4\7\eeee
\8\3\7\5\6\5\7\3\8\2\8\3\7\5\6\5\eeee
\1\8\3\7\4\6\4\7\3\8\1\8\3\7\4\6\eeee
\8\2\8\3\7\5\6\5\7\3\8\2\8\3\7\5\eeee
} 

\baaa
8-7
\eaaa
\bbbb
0&0&0&0&0&0&0&4\\
0&0&0&0&0&0&0&4\\
0&0&0&0&0&0&2&2\\
0&0&0&0&0&2&2&0\\
0&0&0&0&0&4&0&0\\
0&0&0&2&2&0&0&0\\
0&0&2&2&0&0&0&0\\
1&1&2&0&0&0&0&0\\
\ebbb
\parbox{7cm}{ 
\1\8\3\7\4\6\5\6\4\7\3\8\1\8\3\7\eeee
\8\2\8\3\7\4\6\5\6\4\7\3\8\2\8\3\eeee
\3\8\1\8\3\7\4\6\5\6\4\7\3\8\1\8\eeee
\7\3\8\2\8\3\7\4\6\5\6\4\7\3\8\2\eeee
\4\7\3\8\1\8\3\7\4\6\5\6\4\7\3\8\eeee
\6\4\7\3\8\2\8\3\7\4\6\5\6\4\7\3\eeee
\5\6\4\7\3\8\1\8\3\7\4\6\5\6\4\7\eeee
\6\5\6\4\7\3\8\2\8\3\7\4\6\5\6\4\eeee
\4\6\5\6\4\7\3\8\1\8\3\7\4\6\5\6\eeee
\7\4\6\5\6\4\7\3\8\2\8\3\7\4\6\5\eeee
\3\7\4\6\5\6\4\7\3\8\1\8\3\7\4\6\eeee
\8\3\7\4\6\5\6\4\7\3\8\2\8\3\7\4\eeee
} 

\baaa
8-8
\eaaa
\bbbb
0&0&0&0&0&0&0&4\\
0&0&0&0&0&0&0&4\\
0&0&0&0&0&0&2&2\\
0&0&0&0&0&2&2&0\\
0&0&0&0&2&2&0&0\\
0&0&0&2&2&0&0&0\\
0&0&2&2&0&0&0&0\\
1&1&2&0&0&0&0&0\\
\ebbb
\parbox{7cm}{ 
\1\8\3\7\4\6\5\5\6\4\7\3\8\1\8\3\eeee
\8\2\8\3\7\4\6\5\5\6\4\7\3\8\2\8\eeee
\3\8\1\8\3\7\4\6\5\5\6\4\7\3\8\1\eeee
\7\3\8\2\8\3\7\4\6\5\5\6\4\7\3\8\eeee
\4\7\3\8\1\8\3\7\4\6\5\5\6\4\7\3\eeee
\6\4\7\3\8\2\8\3\7\4\6\5\5\6\4\7\eeee
\5\6\4\7\3\8\1\8\3\7\4\6\5\5\6\4\eeee
\5\5\6\4\7\3\8\2\8\3\7\4\6\5\5\6\eeee
\6\5\5\6\4\7\3\8\1\8\3\7\4\6\5\5\eeee
\4\6\5\5\6\4\7\3\8\2\8\3\7\4\6\5\eeee
\7\4\6\5\5\6\4\7\3\8\1\8\3\7\4\6\eeee
\3\7\4\6\5\5\6\4\7\3\8\2\8\3\7\4\eeee
} 

\baaa
8-9
\eaaa
\bbbb
0&0&0&0&0&0&0&4\\
0&0&0&0&0&0&0&4\\
0&0&0&0&0&0&2&2\\
0&0&0&0&0&4&0&0\\
0&0&0&0&0&4&0&0\\
0&0&0&1&1&0&2&0\\
0&0&2&0&0&2&0&0\\
1&1&2&0&0&0&0&0\\
\ebbb
\parbox{7cm}{ 
\1\8\3\7\6\4\6\7\3\8\1\8\3\7\6\4\eeee
\8\2\8\3\7\6\5\6\7\3\8\2\8\3\7\6\eeee
\3\8\1\8\3\7\6\4\6\7\3\8\1\8\3\7\eeee
\7\3\8\2\8\3\7\6\5\6\7\3\8\2\8\3\eeee
\6\7\3\8\1\8\3\7\6\4\6\7\3\8\1\8\eeee
\4\6\7\3\8\2\8\3\7\6\5\6\7\3\8\2\eeee
\6\5\6\7\3\8\1\8\3\7\6\4\6\7\3\8\eeee
\7\6\4\6\7\3\8\2\8\3\7\6\5\6\7\3\eeee
\3\7\6\5\6\7\3\8\1\8\3\7\6\4\6\7\eeee
\8\3\7\6\4\6\7\3\8\2\8\3\7\6\5\6\eeee
\1\8\3\7\6\5\6\7\3\8\1\8\3\7\6\4\eeee
\8\2\8\3\7\6\4\6\7\3\8\2\8\3\7\6\eeee
} 

\baaa
8-10
\eaaa
\bbbb
0&0&0&0&0&0&0&4\\
0&0&0&0&0&0&0&4\\
0&0&0&0&0&0&2&2\\
0&0&0&0&1&1&2&0\\
0&0&0&2&1&1&0&0\\
0&0&0&2&1&1&0&0\\
0&0&2&2&0&0&0&0\\
1&1&2&0&0&0&0&0\\
\ebbb
\parbox{7cm}{ 
\1\8\3\7\4\5\5\4\7\3\8\1\8\3\7\4\eeee
\8\2\8\3\7\4\6\6\4\7\3\8\2\8\3\7\eeee
\3\8\1\8\3\7\4\5\5\4\7\3\8\1\8\3\eeee
\7\3\8\2\8\3\7\4\6\6\4\7\3\8\2\8\eeee
\4\7\3\8\1\8\3\7\4\5\5\4\7\3\8\1\eeee
\5\4\7\3\8\2\8\3\7\4\6\6\4\7\3\8\eeee
\5\6\4\7\3\8\1\8\3\7\4\5\5\4\7\3\eeee
\4\6\5\4\7\3\8\2\8\3\7\4\6\6\4\7\eeee
\7\4\5\6\4\7\3\8\1\8\3\7\4\5\5\4\eeee
\3\7\4\6\5\4\7\3\8\2\8\3\7\4\6\6\eeee
\8\3\7\4\5\6\4\7\3\8\1\8\3\7\4\5\eeee
\1\8\3\7\4\6\5\4\7\3\8\2\8\3\7\4\eeee
} 

\baaa
8-11
\eaaa
\bbbb
0&0&0&0&0&0&0&4\\
0&0&0&0&0&0&0&4\\
0&0&0&0&0&0&4&0\\
0&0&0&0&0&2&2&0\\
0&0&0&0&0&2&2&0\\
0&0&0&1&1&0&0&2\\
0&0&2&1&1&0&0&0\\
1&1&0&0&0&2&0&0\\
\ebbb
\parbox{7cm}{ 
\1\8\6\4\7\3\7\4\6\8\1\8\6\4\7\3\eeee
\8\2\8\6\5\7\3\7\5\6\8\2\8\6\5\7\eeee
\6\8\1\8\6\4\7\3\7\4\6\8\1\8\6\4\eeee
\4\6\8\2\8\6\5\7\3\7\5\6\8\2\8\6\eeee
\7\5\6\8\1\8\6\4\7\3\7\4\6\8\1\8\eeee
\3\7\4\6\8\2\8\6\5\7\3\7\5\6\8\2\eeee
\7\3\7\5\6\8\1\8\6\4\7\3\7\4\6\8\eeee
\4\7\3\7\4\6\8\2\8\6\5\7\3\7\5\6\eeee
\6\5\7\3\7\5\6\8\1\8\6\4\7\3\7\4\eeee
\8\6\4\7\3\7\4\6\8\2\8\6\5\7\3\7\eeee
\1\8\6\5\7\3\7\5\6\8\1\8\6\4\7\3\eeee
\8\2\8\6\4\7\3\7\4\6\8\2\8\6\5\7\eeee
} 

\baaa
8-12
\eaaa
\bbbb
0&0&0&0&0&0&0&4\\
0&0&0&0&0&0&0&4\\
0&0&0&0&0&1&1&2\\
0&0&0&0&0&1&1&2\\
0&0&0&0&0&2&2&0\\
0&0&1&1&2&0&0&0\\
0&0&1&1&2&0&0&0\\
1&1&1&1&0&0&0&0\\
\ebbb
\parbox{7cm}{ 
\1\8\3\6\5\6\3\8\1\8\3\6\5\6\3\8\eeee
\8\2\8\4\7\5\7\4\8\2\8\4\7\5\7\4\eeee
\3\8\1\8\3\6\5\6\3\8\1\8\3\6\5\6\eeee
\6\4\8\2\8\4\7\5\7\4\8\2\8\4\7\5\eeee
\5\7\3\8\1\8\3\6\5\6\3\8\1\8\3\6\eeee
\6\5\6\4\8\2\8\4\7\5\7\4\8\2\8\4\eeee
\3\7\5\7\3\8\1\8\3\6\5\6\3\8\1\8\eeee
\8\4\6\5\6\4\8\2\8\4\7\5\7\4\8\2\eeee
\1\8\3\7\5\7\3\8\1\8\3\6\5\6\3\8\eeee
\8\2\8\4\6\5\6\4\8\2\8\4\7\5\7\4\eeee
\3\8\1\8\3\7\5\7\3\8\1\8\3\6\5\6\eeee
\6\4\8\2\8\4\6\5\6\4\8\2\8\4\7\5\eeee
} 

\baaa
8-13
\eaaa
\bbbb
0&0&0&0&0&0&0&4\\
0&0&0&0&0&0&0&4\\
0&0&0&0&0&1&1&2\\
0&0&0&0&0&1&1&2\\
0&0&0&0&2&1&1&0\\
0&0&1&1&2&0&0&0\\
0&0&1&1&2&0&0&0\\
1&1&1&1&0&0&0&0\\
\ebbb
\parbox{7cm}{ 
\1\8\3\6\5\5\6\3\8\1\8\3\6\5\5\6\eeee
\8\2\8\4\7\5\5\7\4\8\2\8\4\7\5\5\eeee
\3\8\1\8\3\6\5\5\6\3\8\1\8\3\6\5\eeee
\6\4\8\2\8\4\7\5\5\7\4\8\2\8\4\7\eeee
\5\7\3\8\1\8\3\6\5\5\6\3\8\1\8\3\eeee
\5\5\6\4\8\2\8\4\7\5\5\7\4\8\2\8\eeee
\6\5\5\7\3\8\1\8\3\6\5\5\6\3\8\1\eeee
\3\7\5\5\6\4\8\2\8\4\7\5\5\7\4\8\eeee
\8\4\6\5\5\7\3\8\1\8\3\6\5\5\6\3\eeee
\1\8\3\7\5\5\6\4\8\2\8\4\7\5\5\7\eeee
\8\2\8\4\6\5\5\7\3\8\1\8\3\6\5\5\eeee
\3\8\1\8\3\7\5\5\6\4\8\2\8\4\7\5\eeee
} 

\baaa
8-14
\eaaa
\bbbb
0&0&0&0&0&0&0&4\\
0&0&0&0&0&0&0&4\\
0&0&0&0&0&1&1&2\\
0&0&0&0&0&2&2&0\\
0&0&0&0&0&2&2&0\\
0&0&2&1&1&0&0&0\\
0&0&2&1&1&0&0&0\\
1&1&2&0&0&0&0&0\\
\ebbb
\parbox{7cm}{ 
\1\8\3\6\4\6\3\8\1\8\3\6\4\6\3\8\eeee
\8\2\8\3\7\5\7\3\8\2\8\3\7\5\7\3\eeee
\3\8\1\8\3\6\4\6\3\8\1\8\3\6\4\6\eeee
\6\3\8\2\8\3\7\5\7\3\8\2\8\3\7\5\eeee
\4\7\3\8\1\8\3\6\4\6\3\8\1\8\3\6\eeee
\6\5\6\3\8\2\8\3\7\5\7\3\8\2\8\3\eeee
\3\7\4\7\3\8\1\8\3\6\4\6\3\8\1\8\eeee
\8\3\6\5\6\3\8\2\8\3\7\5\7\3\8\2\eeee
\1\8\3\7\4\7\3\8\1\8\3\6\4\6\3\8\eeee
\8\2\8\3\6\5\6\3\8\2\8\3\7\5\7\3\eeee
\3\8\1\8\3\7\4\7\3\8\1\8\3\6\4\6\eeee
\6\3\8\2\8\3\6\5\6\3\8\2\8\3\7\5\eeee
} 

\baaa
8-15
\eaaa
\bbbb
0&0&0&0&0&0&0&4\\
0&0&0&0&0&0&0&4\\
0&0&0&0&0&1&1&2\\
0&0&0&0&2&1&1&0\\
0&0&0&2&2&0&0&0\\
0&0&2&2&0&0&0&0\\
0&0&2&2&0&0&0&0\\
1&1&2&0&0&0&0&0\\
\ebbb
\parbox{7cm}{ 
\1\8\3\6\4\5\5\4\6\3\8\1\8\3\6\4\eeee
\8\2\8\3\7\4\5\5\4\7\3\8\2\8\3\7\eeee
\3\8\1\8\3\6\4\5\5\4\6\3\8\1\8\3\eeee
\6\3\8\2\8\3\7\4\5\5\4\7\3\8\2\8\eeee
\4\7\3\8\1\8\3\6\4\5\5\4\6\3\8\1\eeee
\5\4\6\3\8\2\8\3\7\4\5\5\4\7\3\8\eeee
\5\5\4\7\3\8\1\8\3\6\4\5\5\4\6\3\eeee
\4\5\5\4\6\3\8\2\8\3\7\4\5\5\4\7\eeee
\6\4\5\5\4\7\3\8\1\8\3\6\4\5\5\4\eeee
\3\7\4\5\5\4\6\3\8\2\8\3\7\4\5\5\eeee
\8\3\6\4\5\5\4\7\3\8\1\8\3\6\4\5\eeee
\1\8\3\7\4\5\5\4\6\3\8\2\8\3\7\4\eeee
} 

\baaa
8-16
\eaaa
\bbbb
0&0&0&0&0&0&0&4\\
0&0&0&0&0&0&0&4\\
0&0&0&0&0&1&1&2\\
0&0&0&1&1&1&1&0\\
0&0&0&1&1&1&1&0\\
0&0&2&1&1&0&0&0\\
0&0&2&1&1&0&0&0\\
1&1&2&0&0&0&0&0\\
\ebbb
\parbox{7cm}{ 
\1\8\3\6\4\4\6\3\8\1\8\3\6\4\4\6\eeee
\8\2\8\3\7\5\5\7\3\8\2\8\3\7\5\5\eeee
\3\8\1\8\3\6\4\4\6\3\8\1\8\3\6\4\eeee
\6\3\8\2\8\3\7\5\5\7\3\8\2\8\3\7\eeee
\4\7\3\8\1\8\3\6\4\4\6\3\8\1\8\3\eeee
\4\5\6\3\8\2\8\3\7\5\5\7\3\8\2\8\eeee
\6\5\4\7\3\8\1\8\3\6\4\4\6\3\8\1\eeee
\3\7\4\5\6\3\8\2\8\3\7\5\5\7\3\8\eeee
\8\3\6\5\4\7\3\8\1\8\3\6\4\4\6\3\eeee
\1\8\3\7\4\5\6\3\8\2\8\3\7\5\5\7\eeee
\8\2\8\3\6\5\4\7\3\8\1\8\3\6\4\4\eeee
\3\8\1\8\3\7\4\5\6\3\8\2\8\3\7\5\eeee
} 

\baaa
8-17
\eaaa
\bbbb
0&0&0&0&0&0&0&4\\
0&0&0&0&0&0&2&2\\
0&0&0&0&0&0&2&2\\
0&0&0&0&0&2&2&0\\
0&0&0&0&0&2&2&0\\
0&0&0&1&1&2&0&0\\
0&1&1&1&1&0&0&0\\
2&1&1&0&0&0&0&0\\
\ebbb
\parbox{7cm}{ 
\1\8\2\7\4\6\6\4\7\2\8\1\8\2\7\4\eeee
\8\1\8\3\7\5\6\6\5\7\3\8\1\8\3\7\eeee
\2\8\1\8\2\7\4\6\6\4\7\2\8\1\8\2\eeee
\7\3\8\1\8\3\7\5\6\6\5\7\3\8\1\8\eeee
\4\7\2\8\1\8\2\7\4\6\6\4\7\2\8\1\eeee
\6\5\7\3\8\1\8\3\7\5\6\6\5\7\3\8\eeee
\6\6\4\7\2\8\1\8\2\7\4\6\6\4\7\2\eeee
\4\6\6\5\7\3\8\1\8\3\7\5\6\6\5\7\eeee
\7\5\6\6\4\7\2\8\1\8\2\7\4\6\6\4\eeee
\2\7\4\6\6\5\7\3\8\1\8\3\7\5\6\6\eeee
\8\3\7\5\6\6\4\7\2\8\1\8\2\7\4\6\eeee
\1\8\2\7\4\6\6\5\7\3\8\1\8\3\7\5\eeee
} 

\baaa
8-18
\eaaa
\bbbb
0&0&0&0&0&0&0&4\\
0&0&0&0&0&0&2&2\\
0&0&0&0&0&0&2&2\\
0&0&0&0&0&2&2&0\\
0&0&0&0&0&2&2&0\\
0&0&0&2&2&0&0&0\\
0&1&1&1&1&0&0&0\\
2&1&1&0&0&0&0&0\\
\ebbb
\parbox{7cm}{ 
\1\8\2\7\4\6\4\7\2\8\1\8\2\7\4\6\eeee
\8\1\8\3\7\5\6\5\7\3\8\1\8\3\7\5\eeee
\2\8\1\8\2\7\4\6\4\7\2\8\1\8\2\7\eeee
\7\3\8\1\8\3\7\5\6\5\7\3\8\1\8\3\eeee
\4\7\2\8\1\8\2\7\4\6\4\7\2\8\1\8\eeee
\6\5\7\3\8\1\8\3\7\5\6\5\7\3\8\1\eeee
\4\6\4\7\2\8\1\8\2\7\4\6\4\7\2\8\eeee
\7\5\6\5\7\3\8\1\8\3\7\5\6\5\7\3\eeee
\2\7\4\6\4\7\2\8\1\8\2\7\4\6\4\7\eeee
\8\3\7\5\6\5\7\3\8\1\8\3\7\5\6\5\eeee
\1\8\2\7\4\6\4\7\2\8\1\8\2\7\4\6\eeee
\8\1\8\3\7\5\6\5\7\3\8\1\8\3\7\5\eeee
} 

\baaa
8-19
\eaaa
\bbbb
0&0&0&0&0&0&0&4\\
0&0&0&0&0&0&2&2\\
0&0&0&0&0&0&2&2\\
0&0&0&0&0&2&2&0\\
0&0&0&0&0&4&0&0\\
0&0&0&2&2&0&0&0\\
0&1&1&2&0&0&0&0\\
2&1&1&0&0&0&0&0\\
\ebbb
\parbox{7cm}{ 
\1\8\2\7\4\6\5\6\4\7\2\8\1\8\2\7\eeee
\8\1\8\3\7\4\6\5\6\4\7\3\8\1\8\3\eeee
\2\8\1\8\2\7\4\6\5\6\4\7\2\8\1\8\eeee
\7\3\8\1\8\3\7\4\6\5\6\4\7\3\8\1\eeee
\4\7\2\8\1\8\2\7\4\6\5\6\4\7\2\8\eeee
\6\4\7\3\8\1\8\3\7\4\6\5\6\4\7\3\eeee
\5\6\4\7\2\8\1\8\2\7\4\6\5\6\4\7\eeee
\6\5\6\4\7\3\8\1\8\3\7\4\6\5\6\4\eeee
\4\6\5\6\4\7\2\8\1\8\2\7\4\6\5\6\eeee
\7\4\6\5\6\4\7\3\8\1\8\3\7\4\6\5\eeee
\2\7\4\6\5\6\4\7\2\8\1\8\2\7\4\6\eeee
\8\3\7\4\6\5\6\4\7\3\8\1\8\3\7\4\eeee
} 

\baaa
8-20
\eaaa
\bbbb
0&0&0&0&0&0&0&4\\
0&0&0&0&0&0&2&2\\
0&0&0&0&0&0&2&2\\
0&0&0&0&0&2&2&0\\
0&0&0&0&2&2&0&0\\
0&0&0&2&2&0&0&0\\
0&1&1&2&0&0&0&0\\
2&1&1&0&0&0&0&0\\
\ebbb
\parbox{7cm}{ 
\1\8\2\7\4\6\5\5\6\4\7\2\8\1\8\2\eeee
\8\1\8\3\7\4\6\5\5\6\4\7\3\8\1\8\eeee
\2\8\1\8\2\7\4\6\5\5\6\4\7\2\8\1\eeee
\7\3\8\1\8\3\7\4\6\5\5\6\4\7\3\8\eeee
\4\7\2\8\1\8\2\7\4\6\5\5\6\4\7\2\eeee
\6\4\7\3\8\1\8\3\7\4\6\5\5\6\4\7\eeee
\5\6\4\7\2\8\1\8\2\7\4\6\5\5\6\4\eeee
\5\5\6\4\7\3\8\1\8\3\7\4\6\5\5\6\eeee
\6\5\5\6\4\7\2\8\1\8\2\7\4\6\5\5\eeee
\4\6\5\5\6\4\7\3\8\1\8\3\7\4\6\5\eeee
\7\4\6\5\5\6\4\7\2\8\1\8\2\7\4\6\eeee
\2\7\4\6\5\5\6\4\7\3\8\1\8\3\7\4\eeee
} 

\baaa
8-21
\eaaa
\bbbb
0&0&0&0&0&0&0&4\\
0&0&0&0&0&0&2&2\\
0&0&0&0&0&0&2&2\\
0&0&0&0&1&1&2&0\\
0&0&0&2&1&1&0&0\\
0&0&0&2&1&1&0&0\\
0&1&1&2&0&0&0&0\\
2&1&1&0&0&0&0&0\\
\ebbb
\parbox{7cm}{ 
\1\8\2\7\4\5\5\4\7\2\8\1\8\2\7\4\eeee
\8\1\8\3\7\4\6\6\4\7\3\8\1\8\3\7\eeee
\2\8\1\8\2\7\4\5\5\4\7\2\8\1\8\2\eeee
\7\3\8\1\8\3\7\4\6\6\4\7\3\8\1\8\eeee
\4\7\2\8\1\8\2\7\4\5\5\4\7\2\8\1\eeee
\5\4\7\3\8\1\8\3\7\4\6\6\4\7\3\8\eeee
\5\6\4\7\2\8\1\8\2\7\4\5\5\4\7\2\eeee
\4\6\5\4\7\3\8\1\8\3\7\4\6\6\4\7\eeee
\7\4\5\6\4\7\2\8\1\8\2\7\4\5\5\4\eeee
\2\7\4\6\5\4\7\3\8\1\8\3\7\4\6\6\eeee
\8\3\7\4\5\6\4\7\2\8\1\8\2\7\4\5\eeee
\1\8\2\7\4\6\5\4\7\3\8\1\8\3\7\4\eeee
} 

\baaa
8-22
\eaaa
\bbbb
0&0&0&0&0&0&0&4\\
0&0&0&0&0&0&2&2\\
0&0&0&0&0&2&0&2\\
0&0&0&0&1&1&1&1\\
0&0&0&2&0&0&2&0\\
0&0&2&2&0&0&0&0\\
0&1&0&1&1&0&1&0\\
1&1&1&1&0&0&0&0\\
\ebbb
\parbox{7cm}{ 
\1\8\4\5\7\4\6\3\8\2\7\7\2\8\3\6\eeee
\8\2\7\7\2\8\3\6\4\7\5\4\8\1\8\4\eeee
\4\7\5\4\8\1\8\4\5\7\4\6\3\8\2\7\eeee
\5\7\4\6\3\8\2\7\7\2\8\3\6\4\7\5\eeee
\7\2\8\3\6\4\7\5\4\8\1\8\4\5\7\4\eeee
\4\8\1\8\4\5\7\4\6\3\8\2\7\7\2\8\eeee
\6\3\8\2\7\7\2\8\3\6\4\7\5\4\8\1\eeee
\3\6\4\7\5\4\8\1\8\4\5\7\4\6\3\8\eeee
\8\4\5\7\4\6\3\8\2\7\7\2\8\3\6\4\eeee
\2\7\7\2\8\3\6\4\7\5\4\8\1\8\4\5\eeee
\7\5\4\8\1\8\4\5\7\4\6\3\8\2\7\7\eeee
\7\4\6\3\8\2\7\7\2\8\3\6\4\7\5\4\eeee
} 

\baaa
8-23
\eaaa
\bbbb
0&0&0&0&0&0&0&4\\
0&0&0&0&0&0&2&2\\
0&0&0&0&0&2&0&2\\
0&0&0&0&1&1&1&1\\
0&0&0&4&0&0&0&0\\
0&0&2&2&0&0&0&0\\
0&2&0&2&0&0&0&0\\
1&1&1&1&0&0&0&0\\
\ebbb
\parbox{7cm}{ 
\1\8\4\5\4\8\1\8\4\5\4\8\1\8\4\5\eeee
\8\2\7\4\6\3\8\2\7\4\6\3\8\2\7\4\eeee
\4\7\2\8\3\6\4\7\2\8\3\6\4\7\2\8\eeee
\5\4\8\1\8\4\5\4\8\1\8\4\5\4\8\1\eeee
\4\6\3\8\2\7\4\6\3\8\2\7\4\6\3\8\eeee
\8\3\6\4\7\2\8\3\6\4\7\2\8\3\6\4\eeee
\1\8\4\5\4\8\1\8\4\5\4\8\1\8\4\5\eeee
\8\2\7\4\6\3\8\2\7\4\6\3\8\2\7\4\eeee
\4\7\2\8\3\6\4\7\2\8\3\6\4\7\2\8\eeee
\5\4\8\1\8\4\5\4\8\1\8\4\5\4\8\1\eeee
\4\6\3\8\2\7\4\6\3\8\2\7\4\6\3\8\eeee
\8\3\6\4\7\2\8\3\6\4\7\2\8\3\6\4\eeee
} 

\baaa
8-24
\eaaa
\bbbb
0&0&0&0&0&0&0&4\\
0&0&0&0&0&0&2&2\\
0&0&0&0&0&2&2&0\\
0&0&0&0&0&2&2&0\\
0&0&0&0&2&2&0&0\\
0&0&1&1&2&0&0&0\\
0&2&1&1&0&0&0&0\\
2&2&0&0&0&0&0&0\\
\ebbb
\parbox{7cm}{ 
\1\8\2\7\3\6\5\5\6\3\7\2\8\1\8\2\eeee
\8\1\8\2\7\4\6\5\5\6\4\7\2\8\1\8\eeee
\2\8\1\8\2\7\3\6\5\5\6\3\7\2\8\1\eeee
\7\2\8\1\8\2\7\4\6\5\5\6\4\7\2\8\eeee
\3\7\2\8\1\8\2\7\3\6\5\5\6\3\7\2\eeee
\6\4\7\2\8\1\8\2\7\4\6\5\5\6\4\7\eeee
\5\6\3\7\2\8\1\8\2\7\3\6\5\5\6\3\eeee
\5\5\6\4\7\2\8\1\8\2\7\4\6\5\5\6\eeee
\6\5\5\6\3\7\2\8\1\8\2\7\3\6\5\5\eeee
\3\6\5\5\6\4\7\2\8\1\8\2\7\4\6\5\eeee
\7\4\6\5\5\6\3\7\2\8\1\8\2\7\3\6\eeee
\2\7\3\6\5\5\6\4\7\2\8\1\8\2\7\4\eeee
} 

\baaa
8-25
\eaaa
\bbbb
0&0&0&0&0&0&0&4\\
0&0&0&0&0&0&2&2\\
0&0&0&0&0&2&2&0\\
0&0&0&0&2&2&0&0\\
0&0&0&2&2&0&0&0\\
0&0&2&2&0&0&0&0\\
0&2&2&0&0&0&0&0\\
2&2&0&0&0&0&0&0\\
\ebbb
\parbox{7cm}{ 
\1\8\2\7\3\6\4\5\5\4\6\3\7\2\8\1\eeee
\8\1\8\2\7\3\6\4\5\5\4\6\3\7\2\8\eeee
\2\8\1\8\2\7\3\6\4\5\5\4\6\3\7\2\eeee
\7\2\8\1\8\2\7\3\6\4\5\5\4\6\3\7\eeee
\3\7\2\8\1\8\2\7\3\6\4\5\5\4\6\3\eeee
\6\3\7\2\8\1\8\2\7\3\6\4\5\5\4\6\eeee
\4\6\3\7\2\8\1\8\2\7\3\6\4\5\5\4\eeee
\5\4\6\3\7\2\8\1\8\2\7\3\6\4\5\5\eeee
\5\5\4\6\3\7\2\8\1\8\2\7\3\6\4\5\eeee
\4\5\5\4\6\3\7\2\8\1\8\2\7\3\6\4\eeee
\6\4\5\5\4\6\3\7\2\8\1\8\2\7\3\6\eeee
\3\6\4\5\5\4\6\3\7\2\8\1\8\2\7\3\eeee
} 

\baaa
8-26
\eaaa
\bbbb
0&0&0&0&0&0&0&4\\
0&0&0&0&0&0&2&2\\
0&0&0&0&0&2&2&0\\
0&0&0&0&2&2&0&0\\
0&0&0&4&0&0&0&0\\
0&0&2&2&0&0&0&0\\
0&2&2&0&0&0&0&0\\
2&2&0&0&0&0&0&0\\
\ebbb
\parbox{7cm}{ 
\1\8\2\7\3\6\4\5\4\6\3\7\2\8\1\8\eeee
\8\1\8\2\7\3\6\4\5\4\6\3\7\2\8\1\eeee
\2\8\1\8\2\7\3\6\4\5\4\6\3\7\2\8\eeee
\7\2\8\1\8\2\7\3\6\4\5\4\6\3\7\2\eeee
\3\7\2\8\1\8\2\7\3\6\4\5\4\6\3\7\eeee
\6\3\7\2\8\1\8\2\7\3\6\4\5\4\6\3\eeee
\4\6\3\7\2\8\1\8\2\7\3\6\4\5\4\6\eeee
\5\4\6\3\7\2\8\1\8\2\7\3\6\4\5\4\eeee
\4\5\4\6\3\7\2\8\1\8\2\7\3\6\4\5\eeee
\6\4\5\4\6\3\7\2\8\1\8\2\7\3\6\4\eeee
\3\6\4\5\4\6\3\7\2\8\1\8\2\7\3\6\eeee
\7\3\6\4\5\4\6\3\7\2\8\1\8\2\7\3\eeee
} 

\baaa
8-27
\eaaa
\bbbb
0&0&0&0&0&0&0&4\\
0&0&0&0&0&0&2&2\\
0&0&0&0&0&2&2&0\\
0&0&0&1&1&2&0&0\\
0&0&0&1&1&2&0&0\\
0&0&2&1&1&0&0&0\\
0&2&2&0&0&0&0&0\\
2&2&0&0&0&0&0&0\\
\ebbb
\parbox{7cm}{ 
\1\8\2\7\3\6\4\4\6\3\7\2\8\1\8\2\eeee
\8\1\8\2\7\3\6\5\5\6\3\7\2\8\1\8\eeee
\2\8\1\8\2\7\3\6\4\4\6\3\7\2\8\1\eeee
\7\2\8\1\8\2\7\3\6\5\5\6\3\7\2\8\eeee
\3\7\2\8\1\8\2\7\3\6\4\4\6\3\7\2\eeee
\6\3\7\2\8\1\8\2\7\3\6\5\5\6\3\7\eeee
\4\6\3\7\2\8\1\8\2\7\3\6\4\4\6\3\eeee
\4\5\6\3\7\2\8\1\8\2\7\3\6\5\5\6\eeee
\6\5\4\6\3\7\2\8\1\8\2\7\3\6\4\4\eeee
\3\6\4\5\6\3\7\2\8\1\8\2\7\3\6\5\eeee
\7\3\6\5\4\6\3\7\2\8\1\8\2\7\3\6\eeee
\2\7\3\6\4\5\6\3\7\2\8\1\8\2\7\3\eeee
} 

\baaa
8-28
\eaaa
\bbbb
0&0&0&0&0&0&0&4\\
0&0&0&0&0&0&2&2\\
0&0&0&0&1&1&2&0\\
0&0&0&0&1&1&2&0\\
0&0&1&1&1&1&0&0\\
0&0&1&1&1&1&0&0\\
0&2&1&1&0&0&0&0\\
2&2&0&0&0&0&0&0\\
\ebbb
\parbox{7cm}{ 
\1\8\2\7\3\5\5\3\7\2\8\1\8\2\7\3\eeee
\8\1\8\2\7\4\6\6\4\7\2\8\1\8\2\7\eeee
\2\8\1\8\2\7\3\5\5\3\7\2\8\1\8\2\eeee
\7\2\8\1\8\2\7\4\6\6\4\7\2\8\1\8\eeee
\3\7\2\8\1\8\2\7\3\5\5\3\7\2\8\1\eeee
\5\4\7\2\8\1\8\2\7\4\6\6\4\7\2\8\eeee
\5\6\3\7\2\8\1\8\2\7\3\5\5\3\7\2\eeee
\3\6\5\4\7\2\8\1\8\2\7\4\6\6\4\7\eeee
\7\4\5\6\3\7\2\8\1\8\2\7\3\5\5\3\eeee
\2\7\3\6\5\4\7\2\8\1\8\2\7\4\6\6\eeee
\8\2\7\4\5\6\3\7\2\8\1\8\2\7\3\5\eeee
\1\8\2\7\3\6\5\4\7\2\8\1\8\2\7\4\eeee
} 

\baaa
8-29
\eaaa
\bbbb
0&0&0&0&0&0&0&4\\
0&0&0&0&0&0&2&2\\
0&0&0&0&1&1&2&0\\
0&0&0&0&1&1&2&0\\
0&0&2&2&0&0&0&0\\
0&0&2&2&0&0&0&0\\
0&2&1&1&0&0&0&0\\
2&2&0&0&0&0&0&0\\
\ebbb
\parbox{7cm}{ 
\1\8\2\7\3\5\3\7\2\8\1\8\2\7\3\5\eeee
\8\1\8\2\7\4\6\4\7\2\8\1\8\2\7\4\eeee
\2\8\1\8\2\7\3\5\3\7\2\8\1\8\2\7\eeee
\7\2\8\1\8\2\7\4\6\4\7\2\8\1\8\2\eeee
\3\7\2\8\1\8\2\7\3\5\3\7\2\8\1\8\eeee
\5\4\7\2\8\1\8\2\7\4\6\4\7\2\8\1\eeee
\3\6\3\7\2\8\1\8\2\7\3\5\3\7\2\8\eeee
\7\4\5\4\7\2\8\1\8\2\7\4\6\4\7\2\eeee
\2\7\3\6\3\7\2\8\1\8\2\7\3\5\3\7\eeee
\8\2\7\4\5\4\7\2\8\1\8\2\7\4\6\4\eeee
\1\8\2\7\3\6\3\7\2\8\1\8\2\7\3\5\eeee
\8\1\8\2\7\4\5\4\7\2\8\1\8\2\7\4\eeee
} 

\baaa
8-30
\eaaa
\bbbb
0&0&0&0&0&0&0&4\\
0&0&0&0&0&0&2&2\\
0&0&0&0&1&1&2&0\\
0&0&0&0&2&2&0&0\\
0&0&2&2&0&0&0&0\\
0&0&2&2&0&0&0&0\\
0&2&2&0&0&0&0&0\\
2&2&0&0&0&0&0&0\\
\ebbb
\parbox{7cm}{ 
\1\8\2\7\3\5\4\5\3\7\2\8\1\8\2\7\eeee
\8\1\8\2\7\3\6\4\6\3\7\2\8\1\8\2\eeee
\2\8\1\8\2\7\3\5\4\5\3\7\2\8\1\8\eeee
\7\2\8\1\8\2\7\3\6\4\6\3\7\2\8\1\eeee
\3\7\2\8\1\8\2\7\3\5\4\5\3\7\2\8\eeee
\5\3\7\2\8\1\8\2\7\3\6\4\6\3\7\2\eeee
\4\6\3\7\2\8\1\8\2\7\3\5\4\5\3\7\eeee
\5\4\5\3\7\2\8\1\8\2\7\3\6\4\6\3\eeee
\3\6\4\6\3\7\2\8\1\8\2\7\3\5\4\5\eeee
\7\3\5\4\5\3\7\2\8\1\8\2\7\3\6\4\eeee
\2\7\3\6\4\6\3\7\2\8\1\8\2\7\3\5\eeee
\8\2\7\3\5\4\5\3\7\2\8\1\8\2\7\3\eeee
} 

\baaa
8-31
\eaaa
\bbbb
0&0&0&0&0&0&0&4\\
0&0&0&0&0&0&2&2\\
0&0&0&0&1&1&2&0\\
0&0&0&2&1&1&0&0\\
0&0&2&2&0&0&0&0\\
0&0&2&2&0&0&0&0\\
0&2&2&0&0&0&0&0\\
2&2&0&0&0&0&0&0\\
\ebbb
\parbox{7cm}{ 
\1\8\2\7\3\5\4\4\5\3\7\2\8\1\8\2\eeee
\8\1\8\2\7\3\6\4\4\6\3\7\2\8\1\8\eeee
\2\8\1\8\2\7\3\5\4\4\5\3\7\2\8\1\eeee
\7\2\8\1\8\2\7\3\6\4\4\6\3\7\2\8\eeee
\3\7\2\8\1\8\2\7\3\5\4\4\5\3\7\2\eeee
\5\3\7\2\8\1\8\2\7\3\6\4\4\6\3\7\eeee
\4\6\3\7\2\8\1\8\2\7\3\5\4\4\5\3\eeee
\4\4\5\3\7\2\8\1\8\2\7\3\6\4\4\6\eeee
\5\4\4\6\3\7\2\8\1\8\2\7\3\5\4\4\eeee
\3\6\4\4\5\3\7\2\8\1\8\2\7\3\6\4\eeee
\7\3\5\4\4\6\3\7\2\8\1\8\2\7\3\5\eeee
\2\7\3\6\4\4\5\3\7\2\8\1\8\2\7\3\eeee
} 

\baaa
8-32
\eaaa
\bbbb
0&0&0&0&0&0&0&4\\
0&0&0&0&0&0&4&0\\
0&0&0&0&1&1&0&2\\
0&0&0&0&1&1&0&2\\
0&0&1&1&0&0&2&0\\
0&0&1&1&0&0&2&0\\
0&1&0&0&1&1&0&1\\
1&0&1&1&0&0&1&0\\
\ebbb
\parbox{7cm}{ 
\1\8\7\2\7\8\1\8\7\2\7\8\1\8\7\2\eeee
\8\3\6\7\5\4\8\3\6\7\5\4\8\3\6\7\eeee
\7\5\4\8\3\6\7\5\4\8\3\6\7\5\4\8\eeee
\2\7\8\1\8\7\2\7\8\1\8\7\2\7\8\1\eeee
\7\6\3\8\4\5\7\6\3\8\4\5\7\6\3\8\eeee
\8\4\5\7\6\3\8\4\5\7\6\3\8\4\5\7\eeee
\1\8\7\2\7\8\1\8\7\2\7\8\1\8\7\2\eeee
\8\3\6\7\5\4\8\3\6\7\5\4\8\3\6\7\eeee
\7\5\4\8\3\6\7\5\4\8\3\6\7\5\4\8\eeee
\2\7\8\1\8\7\2\7\8\1\8\7\2\7\8\1\eeee
\7\6\3\8\4\5\7\6\3\8\4\5\7\6\3\8\eeee
\8\4\5\7\6\3\8\4\5\7\6\3\8\4\5\7\eeee
} 

\baaa
8-33
\eaaa
\bbbb
0&0&0&0&0&0&0&4\\
0&0&0&0&0&0&4&0\\
0&0&0&0&1&1&0&2\\
0&0&0&0&1&1&0&2\\
0&0&1&1&0&0&2&0\\
0&0&1&1&0&0&2&0\\
0&2&0&0&1&1&0&0\\
2&0&1&1&0&0&0&0\\
\ebbb
\parbox{7cm}{ 
\1\8\3\5\7\2\7\5\3\8\1\8\3\5\7\2\eeee
\8\1\8\4\6\7\2\7\6\4\8\1\8\4\6\7\eeee
\3\8\1\8\3\5\7\2\7\5\3\8\1\8\3\5\eeee
\5\4\8\1\8\4\6\7\2\7\6\4\8\1\8\4\eeee
\7\6\3\8\1\8\3\5\7\2\7\5\3\8\1\8\eeee
\2\7\5\4\8\1\8\4\6\7\2\7\6\4\8\1\eeee
\7\2\7\6\3\8\1\8\3\5\7\2\7\5\3\8\eeee
\5\7\2\7\5\4\8\1\8\4\6\7\2\7\6\4\eeee
\3\6\7\2\7\6\3\8\1\8\3\5\7\2\7\5\eeee
\8\4\5\7\2\7\5\4\8\1\8\4\6\7\2\7\eeee
\1\8\3\6\7\2\7\6\3\8\1\8\3\5\7\2\eeee
\8\1\8\4\5\7\2\7\5\4\8\1\8\4\6\7\eeee
} 

\baaa
8-34
\eaaa
\bbbb
0&0&0&0&0&0&0&4\\
0&0&0&0&0&0&4&0\\
0&0&0&0&1&1&0&2\\
0&0&0&0&1&1&2&0\\
0&0&2&2&0&0&0&0\\
0&0&2&2&0&0&0&0\\
0&2&0&2&0&0&0&0\\
2&0&2&0&0&0&0&0\\
\ebbb
\parbox{7cm}{ 
\1\8\3\5\4\7\2\7\4\5\3\8\1\8\3\5\eeee
\8\1\8\3\6\4\7\2\7\4\6\3\8\1\8\3\eeee
\3\8\1\8\3\5\4\7\2\7\4\5\3\8\1\8\eeee
\5\3\8\1\8\3\6\4\7\2\7\4\6\3\8\1\eeee
\4\6\3\8\1\8\3\5\4\7\2\7\4\5\3\8\eeee
\7\4\5\3\8\1\8\3\6\4\7\2\7\4\6\3\eeee
\2\7\4\6\3\8\1\8\3\5\4\7\2\7\4\5\eeee
\7\2\7\4\5\3\8\1\8\3\6\4\7\2\7\4\eeee
\4\7\2\7\4\6\3\8\1\8\3\5\4\7\2\7\eeee
\5\4\7\2\7\4\5\3\8\1\8\3\6\4\7\2\eeee
\3\6\4\7\2\7\4\6\3\8\1\8\3\5\4\7\eeee
\8\3\5\4\7\2\7\4\5\3\8\1\8\3\6\4\eeee
} 

\baaa
8-35
\eaaa
\bbbb
0&0&0&0&0&0&0&4\\
0&0&0&0&0&0&4&0\\
0&0&0&0&1&1&0&2\\
0&0&0&0&1&2&0&1\\
0&0&1&1&0&0&2&0\\
0&0&1&2&0&0&1&0\\
0&1&0&0&2&1&0&0\\
1&0&2&1&0&0&0&0\\
\ebbb
\parbox{7cm}{ 
\1\8\4\6\3\8\3\5\7\5\4\6\7\2\7\6\eeee
\8\3\6\4\8\1\8\4\6\3\8\3\5\7\5\4\eeee
\4\5\7\5\3\8\3\6\4\8\1\8\4\6\3\8\eeee
\6\7\2\7\6\4\5\7\5\3\8\3\6\4\8\1\eeee
\3\5\7\5\4\6\7\2\7\6\4\5\7\5\3\8\eeee
\8\4\6\3\8\3\5\7\5\4\6\7\2\7\6\4\eeee
\3\6\4\8\1\8\4\6\3\8\3\5\7\5\4\6\eeee
\5\7\5\3\8\3\6\4\8\1\8\4\6\3\8\3\eeee
\7\2\7\6\4\5\7\5\3\8\3\6\4\8\1\8\eeee
\5\7\5\4\6\7\2\7\6\4\5\7\5\3\8\3\eeee
\4\6\3\8\3\5\7\5\4\6\7\2\7\6\4\5\eeee
\6\4\8\1\8\4\6\3\8\3\5\7\5\4\6\7\eeee
} 

\baaa
8-36
\eaaa
\bbbb
0&0&0&0&0&0&0&4\\
0&0&0&0&0&1&1&2\\
0&0&0&0&0&1&1&2\\
0&0&0&0&0&1&1&2\\
0&0&0&0&0&2&2&0\\
0&1&1&1&1&0&0&0\\
0&1&1&1&1&0&0&0\\
1&1&1&1&0&0&0&0\\
\ebbb
\parbox{7cm}{ 
\1\8\3\8\1\8\3\8\1\8\3\8\1\8\3\8\eeee
\8\2\7\4\8\2\6\4\8\2\7\4\8\2\6\4\eeee
\3\6\5\6\3\7\5\7\3\6\5\6\3\7\5\7\eeee
\8\4\7\2\8\4\6\2\8\4\7\2\8\4\6\2\eeee
\1\8\3\8\1\8\3\8\1\8\3\8\1\8\3\8\eeee
\8\2\6\4\8\2\7\4\8\2\6\4\8\2\7\4\eeee
\3\7\5\7\3\6\5\6\3\7\5\7\3\6\5\6\eeee
\8\4\6\2\8\4\7\2\8\4\6\2\8\4\7\2\eeee
\1\8\3\8\1\8\3\8\1\8\3\8\1\8\3\8\eeee
\8\2\7\4\8\2\6\4\8\2\7\4\8\2\6\4\eeee
\3\6\5\6\3\7\5\7\3\6\5\6\3\7\5\7\eeee
\8\4\7\2\8\4\6\2\8\4\7\2\8\4\6\2\eeee
} 

\baaa
8-37
\eaaa
\bbbb
0&0&0&0&0&0&0&4\\
0&0&0&0&0&1&1&2\\
0&0&0&0&0&1&1&2\\
0&0&0&0&0&2&2&0\\
0&0&0&0&0&2&2&0\\
0&1&1&1&1&0&0&0\\
0&1&1&1&1&0&0&0\\
2&1&1&0&0&0&0&0\\
\ebbb
\parbox{7cm}{ 
\1\8\2\6\4\6\2\8\1\8\2\6\4\6\2\8\eeee
\8\1\8\3\7\5\7\3\8\1\8\3\7\5\7\3\eeee
\2\8\1\8\2\6\4\6\2\8\1\8\2\6\4\6\eeee
\6\3\8\1\8\3\7\5\7\3\8\1\8\3\7\5\eeee
\4\7\2\8\1\8\2\6\4\6\2\8\1\8\2\6\eeee
\6\5\6\3\8\1\8\3\7\5\7\3\8\1\8\3\eeee
\2\7\4\7\2\8\1\8\2\6\4\6\2\8\1\8\eeee
\8\3\6\5\6\3\8\1\8\3\7\5\7\3\8\1\eeee
\1\8\2\7\4\7\2\8\1\8\2\6\4\6\2\8\eeee
\8\1\8\3\6\5\6\3\8\1\8\3\7\5\7\3\eeee
\2\8\1\8\2\7\4\7\2\8\1\8\2\6\4\6\eeee
\6\3\8\1\8\3\6\5\6\3\8\1\8\3\7\5\eeee
} 

\baaa
8-38
\eaaa
\bbbb
0&0&0&0&0&0&0&4\\
0&0&0&0&0&1&1&2\\
0&0&0&0&0&1&1&2\\
0&0&0&0&2&1&1&0\\
0&0&0&2&2&0&0&0\\
0&1&1&2&0&0&0&0\\
0&1&1&2&0&0&0&0\\
2&1&1&0&0&0&0&0\\
\ebbb
\parbox{7cm}{ 
\1\8\2\6\4\5\5\4\6\2\8\1\8\2\6\4\eeee
\8\1\8\3\7\4\5\5\4\7\3\8\1\8\3\7\eeee
\2\8\1\8\2\6\4\5\5\4\6\2\8\1\8\2\eeee
\6\3\8\1\8\3\7\4\5\5\4\7\3\8\1\8\eeee
\4\7\2\8\1\8\2\6\4\5\5\4\6\2\8\1\eeee
\5\4\6\3\8\1\8\3\7\4\5\5\4\7\3\8\eeee
\5\5\4\7\2\8\1\8\2\6\4\5\5\4\6\2\eeee
\4\5\5\4\6\3\8\1\8\3\7\4\5\5\4\7\eeee
\6\4\5\5\4\7\2\8\1\8\2\6\4\5\5\4\eeee
\2\7\4\5\5\4\6\3\8\1\8\3\7\4\5\5\eeee
\8\3\6\4\5\5\4\7\2\8\1\8\2\6\4\5\eeee
\1\8\2\7\4\5\5\4\6\3\8\1\8\3\7\4\eeee
} 

\baaa
8-39
\eaaa
\bbbb
0&0&0&0&0&0&0&4\\
0&0&0&0&0&1&1&2\\
0&0&0&0&0&1&1&2\\
0&0&0&1&1&1&1&0\\
0&0&0&1&1&1&1&0\\
0&1&1&1&1&0&0&0\\
0&1&1&1&1&0&0&0\\
2&1&1&0&0&0&0&0\\
\ebbb
\parbox{7cm}{ 
\1\8\2\6\4\4\6\2\8\1\8\2\6\4\4\6\eeee
\8\1\8\3\7\5\5\7\3\8\1\8\3\7\5\5\eeee
\2\8\1\8\2\6\4\4\6\2\8\1\8\2\6\4\eeee
\6\3\8\1\8\3\7\5\5\7\3\8\1\8\3\7\eeee
\4\7\2\8\1\8\2\6\4\4\6\2\8\1\8\2\eeee
\4\5\6\3\8\1\8\3\7\5\5\7\3\8\1\8\eeee
\6\5\4\7\2\8\1\8\2\6\4\4\6\2\8\1\eeee
\2\7\4\5\6\3\8\1\8\3\7\5\5\7\3\8\eeee
\8\3\6\5\4\7\2\8\1\8\2\6\4\4\6\2\eeee
\1\8\2\7\4\5\6\3\8\1\8\3\7\5\5\7\eeee
\8\1\8\3\6\5\4\7\2\8\1\8\2\6\4\4\eeee
\2\8\1\8\2\7\4\5\6\3\8\1\8\3\7\5\eeee
} 

\baaa
8-40
\eaaa
\bbbb
0&0&0&0&0&0&0&4\\
0&0&0&0&0&1&1&2\\
0&0&0&0&1&1&1&1\\
0&0&0&2&1&0&1&0\\
0&0&1&1&0&1&1&0\\
0&1&1&0&1&1&0&0\\
0&1&1&1&1&0&0&0\\
1&2&1&0&0&0&0&0\\
\ebbb
\parbox{7cm}{ 
\1\8\3\5\7\3\6\5\4\4\7\2\8\2\6\6\eeee
\8\2\7\4\4\5\6\3\7\5\3\8\1\8\3\5\eeee
\3\6\5\4\4\7\2\8\2\6\6\2\8\2\7\4\eeee
\5\6\3\7\5\3\8\1\8\3\5\7\3\6\5\4\eeee
\7\2\8\2\6\6\2\8\2\7\4\4\5\6\3\7\eeee
\3\8\1\8\3\5\7\3\6\5\4\4\7\2\8\2\eeee
\6\2\8\2\7\4\4\5\6\3\7\5\3\8\1\8\eeee
\5\7\3\6\5\4\4\7\2\8\2\6\6\2\8\2\eeee
\4\4\5\6\3\7\5\3\8\1\8\3\5\7\3\6\eeee
\4\4\7\2\8\2\6\6\2\8\2\7\4\4\5\6\eeee
\7\5\3\8\1\8\3\5\7\3\6\5\4\4\7\2\eeee
\2\6\6\2\8\2\7\4\4\5\6\3\7\5\3\8\eeee
} 

\baaa
8-41
\eaaa
\bbbb
0&0&0&0&0&0&0&4\\
0&0&0&0&0&1&1&2\\
0&0&0&0&2&1&1&0\\
0&0&0&0&2&1&1&0\\
0&0&1&1&2&0&0&0\\
0&2&1&1&0&0&0&0\\
0&2&1&1&0&0&0&0\\
2&2&0&0&0&0&0&0\\
\ebbb
\parbox{7cm}{ 
\1\8\2\6\3\5\5\3\6\2\8\1\8\2\6\3\eeee
\8\1\8\2\7\4\5\5\4\7\2\8\1\8\2\7\eeee
\2\8\1\8\2\6\3\5\5\3\6\2\8\1\8\2\eeee
\6\2\8\1\8\2\7\4\5\5\4\7\2\8\1\8\eeee
\3\7\2\8\1\8\2\6\3\5\5\3\6\2\8\1\eeee
\5\4\6\2\8\1\8\2\7\4\5\5\4\7\2\8\eeee
\5\5\3\7\2\8\1\8\2\6\3\5\5\3\6\2\eeee
\3\5\5\4\6\2\8\1\8\2\7\4\5\5\4\7\eeee
\6\4\5\5\3\7\2\8\1\8\2\6\3\5\5\3\eeee
\2\7\3\5\5\4\6\2\8\1\8\2\7\4\5\5\eeee
\8\2\6\4\5\5\3\7\2\8\1\8\2\6\3\5\eeee
\1\8\2\7\3\5\5\4\6\2\8\1\8\2\7\4\eeee
} 

\baaa
8-42
\eaaa
\bbbb
0&0&0&0&0&0&0&4\\
0&0&0&0&0&1&1&2\\
0&0&0&0&2&1&1&0\\
0&0&0&0&2&1&1&0\\
0&0&2&2&0&0&0&0\\
0&2&1&1&0&0&0&0\\
0&2&1&1&0&0&0&0\\
2&2&0&0&0&0&0&0\\
\ebbb
\parbox{7cm}{ 
\1\8\2\6\3\5\3\6\2\8\1\8\2\6\3\5\eeee
\8\1\8\2\7\4\5\4\7\2\8\1\8\2\7\4\eeee
\2\8\1\8\2\6\3\5\3\6\2\8\1\8\2\6\eeee
\6\2\8\1\8\2\7\4\5\4\7\2\8\1\8\2\eeee
\3\7\2\8\1\8\2\6\3\5\3\6\2\8\1\8\eeee
\5\4\6\2\8\1\8\2\7\4\5\4\7\2\8\1\eeee
\3\5\3\7\2\8\1\8\2\6\3\5\3\6\2\8\eeee
\6\4\5\4\6\2\8\1\8\2\7\4\5\4\7\2\eeee
\2\7\3\5\3\7\2\8\1\8\2\6\3\5\3\6\eeee
\8\2\6\4\5\4\6\2\8\1\8\2\7\4\5\4\eeee
\1\8\2\7\3\5\3\7\2\8\1\8\2\6\3\5\eeee
\8\1\8\2\6\4\5\4\6\2\8\1\8\2\7\4\eeee
} 

\baaa
8-43
\eaaa
\bbbb
0&0&0&0&0&0&0&4\\
0&0&0&0&0&1&1&2\\
0&0&0&0&2&1&1&0\\
0&0&0&2&2&0&0&0\\
0&0&2&2&0&0&0&0\\
0&2&2&0&0&0&0&0\\
0&2&2&0&0&0&0&0\\
2&2&0&0&0&0&0&0\\
\ebbb
\parbox{7cm}{ 
\1\8\2\6\3\5\4\4\5\3\6\2\8\1\8\2\eeee
\8\1\8\2\7\3\5\4\4\5\3\7\2\8\1\8\eeee
\2\8\1\8\2\6\3\5\4\4\5\3\6\2\8\1\eeee
\6\2\8\1\8\2\7\3\5\4\4\5\3\7\2\8\eeee
\3\7\2\8\1\8\2\6\3\5\4\4\5\3\6\2\eeee
\5\3\6\2\8\1\8\2\7\3\5\4\4\5\3\7\eeee
\4\5\3\7\2\8\1\8\2\6\3\5\4\4\5\3\eeee
\4\4\5\3\6\2\8\1\8\2\7\3\5\4\4\5\eeee
\5\4\4\5\3\7\2\8\1\8\2\6\3\5\4\4\eeee
\3\5\4\4\5\3\6\2\8\1\8\2\7\3\5\4\eeee
\6\3\5\4\4\5\3\7\2\8\1\8\2\6\3\5\eeee
\2\7\3\5\4\4\5\3\6\2\8\1\8\2\7\3\eeee
} 

\baaa
8-44
\eaaa
\bbbb
0&0&0&0&0&0&0&4\\
0&0&0&0&0&1&1&2\\
0&0&0&1&1&1&1&0\\
0&0&2&1&1&0&0&0\\
0&0&2&1&1&0&0&0\\
0&2&2&0&0&0&0&0\\
0&2&2&0&0&0&0&0\\
2&2&0&0&0&0&0&0\\
\ebbb
\parbox{7cm}{ 
\1\8\2\6\3\4\4\3\6\2\8\1\8\2\6\3\eeee
\8\1\8\2\7\3\5\5\3\7\2\8\1\8\2\7\eeee
\2\8\1\8\2\6\3\4\4\3\6\2\8\1\8\2\eeee
\6\2\8\1\8\2\7\3\5\5\3\7\2\8\1\8\eeee
\3\7\2\8\1\8\2\6\3\4\4\3\6\2\8\1\eeee
\4\3\6\2\8\1\8\2\7\3\5\5\3\7\2\8\eeee
\4\5\3\7\2\8\1\8\2\6\3\4\4\3\6\2\eeee
\3\5\4\3\6\2\8\1\8\2\7\3\5\5\3\7\eeee
\6\3\4\5\3\7\2\8\1\8\2\6\3\4\4\3\eeee
\2\7\3\5\4\3\6\2\8\1\8\2\7\3\5\5\eeee
\8\2\6\3\4\5\3\7\2\8\1\8\2\6\3\4\eeee
\1\8\2\7\3\5\4\3\6\2\8\1\8\2\7\3\eeee
} 

\baaa
8-45
\eaaa
\bbbb
0&0&0&0&0&0&1&3\\
0&0&0&0&0&1&0&3\\
0&0&0&0&1&0&0&3\\
0&0&0&0&1&1&1&1\\
0&0&1&3&0&0&0&0\\
0&1&0&3&0&0&0&0\\
1&0&0&3&0&0&0&0\\
1&1&1&1&0&0&0&0\\
\ebbb
\parbox{7cm}{ 
\1\8\2\8\3\8\1\8\2\8\3\8\1\8\2\8\eeee
\7\4\6\4\5\4\7\4\6\4\5\4\7\4\6\4\eeee
\4\5\4\7\4\6\4\5\4\7\4\6\4\5\4\7\eeee
\8\3\8\1\8\2\8\3\8\1\8\2\8\3\8\1\eeee
\1\8\2\8\3\8\1\8\2\8\3\8\1\8\2\8\eeee
\7\4\6\4\5\4\7\4\6\4\5\4\7\4\6\4\eeee
\4\5\4\7\4\6\4\5\4\7\4\6\4\5\4\7\eeee
\8\3\8\1\8\2\8\3\8\1\8\2\8\3\8\1\eeee
\1\8\2\8\3\8\1\8\2\8\3\8\1\8\2\8\eeee
\7\4\6\4\5\4\7\4\6\4\5\4\7\4\6\4\eeee
\4\5\4\7\4\6\4\5\4\7\4\6\4\5\4\7\eeee
\8\3\8\1\8\2\8\3\8\1\8\2\8\3\8\1\eeee
} 

\baaa
8-46
\eaaa
\bbbb
0&0&0&0&0&0&1&3\\
0&0&0&0&0&1&2&1\\
0&0&0&0&1&2&1&0\\
0&0&0&0&2&2&0&0\\
0&0&2&2&0&0&0&0\\
0&1&2&1&0&0&0&0\\
1&2&1&0&0&0&0&0\\
3&1&0&0&0&0&0&0\\
\ebbb
\parbox{7cm}{ 
\1\8\1\8\1\8\1\8\1\8\1\8\1\8\1\8\eeee
\7\2\7\2\7\2\7\2\7\2\7\2\7\2\7\2\eeee
\3\6\3\6\3\6\3\6\3\6\3\6\3\6\3\6\eeee
\5\4\5\4\5\4\5\4\5\4\5\4\5\4\5\4\eeee
\3\6\3\6\3\6\3\6\3\6\3\6\3\6\3\6\eeee
\7\2\7\2\7\2\7\2\7\2\7\2\7\2\7\2\eeee
\1\8\1\8\1\8\1\8\1\8\1\8\1\8\1\8\eeee
\8\1\8\1\8\1\8\1\8\1\8\1\8\1\8\1\eeee
\2\7\2\7\2\7\2\7\2\7\2\7\2\7\2\7\eeee
\6\3\6\3\6\3\6\3\6\3\6\3\6\3\6\3\eeee
\4\5\4\5\4\5\4\5\4\5\4\5\4\5\4\5\eeee
\6\3\6\3\6\3\6\3\6\3\6\3\6\3\6\3\eeee
} 

\baaa
8-47
\eaaa
\bbbb
0&0&0&0&0&0&1&3\\
0&0&0&0&0&1&2&1\\
0&0&0&0&1&2&1&0\\
0&0&0&0&3&1&0&0\\
0&0&1&3&0&0&0&0\\
0&1&2&1&0&0&0&0\\
1&2&1&0&0&0&0&0\\
3&1&0&0&0&0&0&0\\
\ebbb
\parbox{7cm}{ 
\1\8\1\8\1\8\1\8\1\8\1\8\1\8\1\8\eeee
\7\2\7\2\7\2\7\2\7\2\7\2\7\2\7\2\eeee
\3\6\3\6\3\6\3\6\3\6\3\6\3\6\3\6\eeee
\5\4\5\4\5\4\5\4\5\4\5\4\5\4\5\4\eeee
\4\5\4\5\4\5\4\5\4\5\4\5\4\5\4\5\eeee
\6\3\6\3\6\3\6\3\6\3\6\3\6\3\6\3\eeee
\2\7\2\7\2\7\2\7\2\7\2\7\2\7\2\7\eeee
\8\1\8\1\8\1\8\1\8\1\8\1\8\1\8\1\eeee
\1\8\1\8\1\8\1\8\1\8\1\8\1\8\1\8\eeee
\7\2\7\2\7\2\7\2\7\2\7\2\7\2\7\2\eeee
\3\6\3\6\3\6\3\6\3\6\3\6\3\6\3\6\eeee
\5\4\5\4\5\4\5\4\5\4\5\4\5\4\5\4\eeee
} 

\baaa
8-48
\eaaa
\bbbb
0&0&0&0&0&0&1&3\\
0&0&0&0&0&1&2&1\\
0&0&0&0&1&2&1&0\\
0&0&0&1&2&1&0&0\\
0&0&1&2&1&0&0&0\\
0&1&2&1&0&0&0&0\\
1&2&1&0&0&0&0&0\\
3&1&0&0&0&0&0&0\\
\ebbb
\parbox{7cm}{ 
\1\8\1\8\1\8\1\8\1\8\1\8\1\8\1\8\eeee
\7\2\7\2\7\2\7\2\7\2\7\2\7\2\7\2\eeee
\3\6\3\6\3\6\3\6\3\6\3\6\3\6\3\6\eeee
\5\4\5\4\5\4\5\4\5\4\5\4\5\4\5\4\eeee
\5\4\5\4\5\4\5\4\5\4\5\4\5\4\5\4\eeee
\3\6\3\6\3\6\3\6\3\6\3\6\3\6\3\6\eeee
\7\2\7\2\7\2\7\2\7\2\7\2\7\2\7\2\eeee
\1\8\1\8\1\8\1\8\1\8\1\8\1\8\1\8\eeee
\8\1\8\1\8\1\8\1\8\1\8\1\8\1\8\1\eeee
\2\7\2\7\2\7\2\7\2\7\2\7\2\7\2\7\eeee
\6\3\6\3\6\3\6\3\6\3\6\3\6\3\6\3\eeee
\4\5\4\5\4\5\4\5\4\5\4\5\4\5\4\5\eeee
} 

\baaa
8-49
\eaaa
\bbbb
0&0&0&0&0&0&1&3\\
0&0&0&0&0&2&1&1\\
0&0&0&0&1&1&0&2\\
0&0&0&0&1&3&0&0\\
0&0&3&1&0&0&0&0\\
0&2&1&1&0&0&0&0\\
1&3&0&0&0&0&0&0\\
1&1&2&0&0&0&0&0\\
\ebbb
\parbox{7cm}{ 
\1\8\3\8\1\8\3\8\1\8\3\8\1\8\3\8\eeee
\7\2\6\2\7\2\6\2\7\2\6\2\7\2\6\2\eeee
\2\6\4\6\2\6\4\6\2\6\4\6\2\6\4\6\eeee
\8\3\5\3\8\3\5\3\8\3\5\3\8\3\5\3\eeee
\1\8\3\8\1\8\3\8\1\8\3\8\1\8\3\8\eeee
\7\2\6\2\7\2\6\2\7\2\6\2\7\2\6\2\eeee
\2\6\4\6\2\6\4\6\2\6\4\6\2\6\4\6\eeee
\8\3\5\3\8\3\5\3\8\3\5\3\8\3\5\3\eeee
\1\8\3\8\1\8\3\8\1\8\3\8\1\8\3\8\eeee
\7\2\6\2\7\2\6\2\7\2\6\2\7\2\6\2\eeee
\2\6\4\6\2\6\4\6\2\6\4\6\2\6\4\6\eeee
\8\3\5\3\8\3\5\3\8\3\5\3\8\3\5\3\eeee
} 

\baaa
8-50
\eaaa
\bbbb
0&0&0&0&0&0&2&2\\
0&0&0&0&0&0&2&2\\
0&0&0&0&0&2&0&2\\
0&0&0&0&0&2&0&2\\
0&0&0&0&0&2&2&0\\
0&0&1&1&2&0&0&0\\
1&1&0&0&2&0&0&0\\
1&1&1&1&0&0&0&0\\
\ebbb
\parbox{7cm}{ 
\1\7\5\6\3\8\2\7\5\6\4\8\1\7\5\6\eeee
\7\5\6\4\8\1\7\5\6\3\8\2\7\5\6\4\eeee
\5\6\3\8\2\7\5\6\4\8\1\7\5\6\3\8\eeee
\6\4\8\1\7\5\6\3\8\2\7\5\6\4\8\1\eeee
\3\8\2\7\5\6\4\8\1\7\5\6\3\8\2\7\eeee
\8\1\7\5\6\3\8\2\7\5\6\4\8\1\7\5\eeee
\2\7\5\6\4\8\1\7\5\6\3\8\2\7\5\6\eeee
\7\5\6\3\8\2\7\5\6\4\8\1\7\5\6\3\eeee
\5\6\4\8\1\7\5\6\3\8\2\7\5\6\4\8\eeee
\6\3\8\2\7\5\6\4\8\1\7\5\6\3\8\2\eeee
\4\8\1\7\5\6\3\8\2\7\5\6\4\8\1\7\eeee
\8\2\7\5\6\4\8\1\7\5\6\3\8\2\7\5\eeee
} 

\baaa
8-51
\eaaa
\bbbb
0&0&0&0&0&0&2&2\\
0&0&0&0&0&0&2&2\\
0&0&0&0&0&2&0&2\\
0&0&0&0&0&2&0&2\\
0&0&0&0&2&0&2&0\\
0&0&1&1&0&2&0&0\\
1&1&0&0&2&0&0&0\\
1&1&1&1&0&0&0&0\\
\ebbb
\parbox{7cm}{ 
\1\7\5\5\7\2\8\3\6\6\4\8\1\7\5\5\eeee
\7\5\5\7\1\8\4\6\6\3\8\2\7\5\5\7\eeee
\5\5\7\2\8\3\6\6\4\8\1\7\5\5\7\2\eeee
\5\7\1\8\4\6\6\3\8\2\7\5\5\7\1\8\eeee
\7\2\8\3\6\6\4\8\1\7\5\5\7\2\8\3\eeee
\1\8\4\6\6\3\8\2\7\5\5\7\1\8\4\6\eeee
\8\3\6\6\4\8\1\7\5\5\7\2\8\3\6\6\eeee
\4\6\6\3\8\2\7\5\5\7\1\8\4\6\6\3\eeee
\6\6\4\8\1\7\5\5\7\2\8\3\6\6\4\8\eeee
\6\3\8\2\7\5\5\7\1\8\4\6\6\3\8\2\eeee
\4\8\1\7\5\5\7\2\8\3\6\6\4\8\1\7\eeee
\8\2\7\5\5\7\1\8\4\6\6\3\8\2\7\5\eeee
} 

\baaa
8-52
\eaaa
\bbbb
0&0&0&0&0&0&2&2\\
0&0&0&0&0&0&2&2\\
0&0&0&0&0&2&0&2\\
0&0&0&0&0&2&0&2\\
0&0&0&0&2&0&2&0\\
0&0&2&2&0&0&0&0\\
1&1&0&0&2&0&0&0\\
1&1&1&1&0&0&0&0\\
\ebbb
\parbox{7cm}{ 
\1\7\5\5\7\2\8\3\6\4\8\1\7\5\5\7\eeee
\7\5\5\7\1\8\4\6\3\8\2\7\5\5\7\2\eeee
\5\5\7\2\8\3\6\4\8\1\7\5\5\7\1\8\eeee
\5\7\1\8\4\6\3\8\2\7\5\5\7\2\8\3\eeee
\7\2\8\3\6\4\8\1\7\5\5\7\1\8\4\6\eeee
\1\8\4\6\3\8\2\7\5\5\7\2\8\3\6\4\eeee
\8\3\6\4\8\1\7\5\5\7\1\8\4\6\3\8\eeee
\4\6\3\8\2\7\5\5\7\2\8\3\6\4\8\1\eeee
\6\4\8\1\7\5\5\7\1\8\4\6\3\8\2\7\eeee
\3\8\2\7\5\5\7\2\8\3\6\4\8\1\7\5\eeee
\8\1\7\5\5\7\1\8\4\6\3\8\2\7\5\5\eeee
\2\7\5\5\7\2\8\3\6\4\8\1\7\5\5\7\eeee
} 

\baaa
8-53
\eaaa
\bbbb
0&0&0&0&0&0&2&2\\
0&0&0&0&0&0&2&2\\
0&0&0&0&0&2&0&2\\
0&0&0&0&0&2&1&1\\
0&0&0&0&0&2&2&0\\
0&0&1&2&1&0&0&0\\
1&1&0&1&1&0&0&0\\
1&1&1&1&0&0&0&0\\
\ebbb
\parbox{7cm}{ 
\1\7\4\6\5\7\2\7\5\6\4\7\1\8\4\6\eeee
\7\5\6\4\7\1\8\4\6\3\8\2\8\3\6\4\eeee
\4\6\3\8\2\8\3\6\4\8\1\7\4\6\5\7\eeee
\6\4\8\1\7\4\6\5\7\2\7\5\6\4\7\1\eeee
\5\7\2\7\5\6\4\7\1\8\4\6\3\8\2\8\eeee
\7\1\8\4\6\3\8\2\8\3\6\4\8\1\7\4\eeee
\2\8\3\6\4\8\1\7\4\6\5\7\2\7\5\6\eeee
\7\4\6\5\7\2\7\5\6\4\7\1\8\4\6\3\eeee
\5\6\4\7\1\8\4\6\3\8\2\8\3\6\4\8\eeee
\6\3\8\2\8\3\6\4\8\1\7\4\6\5\7\2\eeee
\4\8\1\7\4\6\5\7\2\7\5\6\4\7\1\8\eeee
\7\2\7\5\6\4\7\1\8\4\6\3\8\2\8\3\eeee
} 

\baaa
8-54
\eaaa
\bbbb
0&0&0&0&0&0&2&2\\
0&0&0&0&0&0&2&2\\
0&0&0&0&0&2&0&2\\
0&0&0&0&2&0&2&0\\
0&0&0&2&0&2&0&0\\
0&0&2&0&2&0&0&0\\
1&1&0&2&0&0&0&0\\
1&1&2&0&0&0&0&0\\
\ebbb
\parbox{7cm}{ 
\1\7\4\5\6\3\8\1\7\4\5\6\3\8\2\7\eeee
\7\4\5\6\3\8\2\7\4\5\6\3\8\1\7\4\eeee
\4\5\6\3\8\1\7\4\5\6\3\8\2\7\4\5\eeee
\5\6\3\8\2\7\4\5\6\3\8\1\7\4\5\6\eeee
\6\3\8\1\7\4\5\6\3\8\2\7\4\5\6\3\eeee
\3\8\2\7\4\5\6\3\8\1\7\4\5\6\3\8\eeee
\8\1\7\4\5\6\3\8\2\7\4\5\6\3\8\2\eeee
\2\7\4\5\6\3\8\1\7\4\5\6\3\8\1\7\eeee
\7\4\5\6\3\8\2\7\4\5\6\3\8\2\7\4\eeee
\4\5\6\3\8\1\7\4\5\6\3\8\1\7\4\5\eeee
\5\6\3\8\2\7\4\5\6\3\8\2\7\4\5\6\eeee
\6\3\8\1\7\4\5\6\3\8\1\7\4\5\6\3\eeee
} 

\baaa
8-55
\eaaa
\bbbb
0&0&0&0&0&0&2&2\\
0&0&0&0&0&0&2&2\\
0&0&0&0&0&2&0&2\\
0&0&0&0&2&0&2&0\\
0&0&0&2&2&0&0&0\\
0&0&2&0&0&2&0&0\\
1&1&0&2&0&0&0&0\\
1&1&2&0&0&0&0&0\\
\ebbb
\parbox{7cm}{ 
\1\7\4\5\5\4\7\1\8\3\6\6\3\8\2\7\eeee
\7\4\5\5\4\7\2\8\3\6\6\3\8\1\7\4\eeee
\4\5\5\4\7\1\8\3\6\6\3\8\2\7\4\5\eeee
\5\5\4\7\2\8\3\6\6\3\8\1\7\4\5\5\eeee
\5\4\7\1\8\3\6\6\3\8\2\7\4\5\5\4\eeee
\4\7\2\8\3\6\6\3\8\1\7\4\5\5\4\7\eeee
\7\1\8\3\6\6\3\8\2\7\4\5\5\4\7\2\eeee
\2\8\3\6\6\3\8\1\7\4\5\5\4\7\1\8\eeee
\8\3\6\6\3\8\2\7\4\5\5\4\7\2\8\3\eeee
\3\6\6\3\8\1\7\4\5\5\4\7\1\8\3\6\eeee
\6\6\3\8\2\7\4\5\5\4\7\2\8\3\6\6\eeee
\6\3\8\1\7\4\5\5\4\7\1\8\3\6\6\3\eeee
} 

\baaa
8-56
\eaaa
\bbbb
0&0&0&0&0&0&2&2\\
0&0&0&0&0&0&2&2\\
0&0&0&0&0&2&0&2\\
0&0&0&0&2&2&0&0\\
0&0&0&2&2&0&0&0\\
0&0&2&2&0&0&0&0\\
1&1&0&0&0&0&2&0\\
1&1&2&0&0&0&0&0\\
\ebbb
\parbox{7cm}{ 
\1\7\7\1\8\3\6\4\5\5\4\6\3\8\2\7\eeee
\7\7\2\8\3\6\4\5\5\4\6\3\8\1\7\7\eeee
\7\1\8\3\6\4\5\5\4\6\3\8\2\7\7\2\eeee
\2\8\3\6\4\5\5\4\6\3\8\1\7\7\1\8\eeee
\8\3\6\4\5\5\4\6\3\8\2\7\7\2\8\3\eeee
\3\6\4\5\5\4\6\3\8\1\7\7\1\8\3\6\eeee
\6\4\5\5\4\6\3\8\2\7\7\2\8\3\6\4\eeee
\4\5\5\4\6\3\8\1\7\7\1\8\3\6\4\5\eeee
\5\5\4\6\3\8\2\7\7\2\8\3\6\4\5\5\eeee
\5\4\6\3\8\1\7\7\1\8\3\6\4\5\5\4\eeee
\4\6\3\8\2\7\7\2\8\3\6\4\5\5\4\6\eeee
\6\3\8\1\7\7\1\8\3\6\4\5\5\4\6\3\eeee
} 

\baaa
8-57
\eaaa
\bbbb
0&0&0&0&0&0&2&2\\
0&0&0&0&0&0&2&2\\
0&0&0&0&0&2&0&2\\
0&0&0&0&2&2&0&0\\
0&0&0&2&2&0&0&0\\
0&0&2&2&0&0&0&0\\
2&2&0&0&0&0&0&0\\
1&1&2&0&0&0&0&0\\
\ebbb
\parbox{7cm}{ 
\1\7\2\8\3\6\4\5\5\4\6\3\8\2\7\1\eeee
\7\1\8\3\6\4\5\5\4\6\3\8\1\7\2\8\eeee
\2\8\3\6\4\5\5\4\6\3\8\2\7\1\8\3\eeee
\8\3\6\4\5\5\4\6\3\8\1\7\2\8\3\6\eeee
\3\6\4\5\5\4\6\3\8\2\7\1\8\3\6\4\eeee
\6\4\5\5\4\6\3\8\1\7\2\8\3\6\4\5\eeee
\4\5\5\4\6\3\8\2\7\1\8\3\6\4\5\5\eeee
\5\5\4\6\3\8\1\7\2\8\3\6\4\5\5\4\eeee
\5\4\6\3\8\2\7\1\8\3\6\4\5\5\4\6\eeee
\4\6\3\8\1\7\2\8\3\6\4\5\5\4\6\3\eeee
\6\3\8\2\7\1\8\3\6\4\5\5\4\6\3\8\eeee
\3\8\1\7\2\8\3\6\4\5\5\4\6\3\8\1\eeee
} 

\baaa
8-58
\eaaa
\bbbb
0&0&0&0&0&0&2&2\\
0&0&0&0&0&0&2&2\\
0&0&0&0&0&2&0&2\\
0&0&0&1&1&0&2&0\\
0&0&0&1&1&0&2&0\\
0&0&2&0&0&2&0&0\\
1&1&0&1&1&0&0&0\\
1&1&2&0&0&0&0&0\\
\ebbb
\parbox{7cm}{ 
\1\7\5\4\7\2\8\3\6\6\3\8\1\7\5\4\eeee
\7\4\5\7\1\8\3\6\6\3\8\2\7\4\5\7\eeee
\5\4\7\2\8\3\6\6\3\8\1\7\5\4\7\2\eeee
\5\7\1\8\3\6\6\3\8\2\7\4\5\7\1\8\eeee
\7\2\8\3\6\6\3\8\1\7\5\4\7\2\8\3\eeee
\1\8\3\6\6\3\8\2\7\4\5\7\1\8\3\6\eeee
\8\3\6\6\3\8\1\7\5\4\7\2\8\3\6\6\eeee
\3\6\6\3\8\2\7\4\5\7\1\8\3\6\6\3\eeee
\6\6\3\8\1\7\5\4\7\2\8\3\6\6\3\8\eeee
\6\3\8\2\7\4\5\7\1\8\3\6\6\3\8\2\eeee
\3\8\1\7\5\4\7\2\8\3\6\6\3\8\1\7\eeee
\8\2\7\4\5\7\1\8\3\6\6\3\8\2\7\4\eeee
} 

\baaa
8-59
\eaaa
\bbbb
0&0&0&0&0&0&2&2\\
0&0&0&0&0&0&2&2\\
0&0&0&0&0&2&0&2\\
0&0&0&1&1&2&0&0\\
0&0&0&1&1&2&0&0\\
0&0&2&1&1&0&0&0\\
1&1&0&0&0&0&2&0\\
1&1&2&0&0&0&0&0\\
\ebbb
\parbox{7cm}{ 
\1\7\7\1\8\3\6\4\4\6\3\8\1\7\7\1\eeee
\7\7\2\8\3\6\5\5\6\3\8\2\7\7\2\8\eeee
\7\1\8\3\6\4\4\6\3\8\1\7\7\1\8\3\eeee
\2\8\3\6\5\5\6\3\8\2\7\7\2\8\3\6\eeee
\8\3\6\4\4\6\3\8\1\7\7\1\8\3\6\4\eeee
\3\6\5\5\6\3\8\2\7\7\2\8\3\6\5\5\eeee
\6\4\4\6\3\8\1\7\7\1\8\3\6\4\4\6\eeee
\5\5\6\3\8\2\7\7\2\8\3\6\5\5\6\3\eeee
\4\6\3\8\1\7\7\1\8\3\6\4\4\6\3\8\eeee
\6\3\8\2\7\7\2\8\3\6\5\5\6\3\8\2\eeee
\3\8\1\7\7\1\8\3\6\4\4\6\3\8\1\7\eeee
\8\2\7\7\2\8\3\6\5\5\6\3\8\2\7\7\eeee
} 

\baaa
8-60
\eaaa
\bbbb
0&0&0&0&0&0&2&2\\
0&0&0&0&0&0&2&2\\
0&0&0&0&0&2&0&2\\
0&0&0&1&1&2&0&0\\
0&0&0&1&1&2&0&0\\
0&0&2&1&1&0&0&0\\
2&2&0&0&0&0&0&0\\
1&1&2&0&0&0&0&0\\
\ebbb
\parbox{7cm}{ 
\1\7\2\8\3\6\5\4\6\3\8\1\7\2\8\3\eeee
\7\1\8\3\6\4\5\6\3\8\2\7\1\8\3\6\eeee
\2\8\3\6\5\4\6\3\8\1\7\2\8\3\6\5\eeee
\8\3\6\4\5\6\3\8\2\7\1\8\3\6\4\4\eeee
\3\6\5\4\6\3\8\1\7\2\8\3\6\5\5\6\eeee
\6\4\5\6\3\8\2\7\1\8\3\6\4\4\6\3\eeee
\5\4\6\3\8\1\7\2\8\3\6\5\5\6\3\8\eeee
\5\6\3\8\2\7\1\8\3\6\4\4\6\3\8\1\eeee
\6\3\8\1\7\2\8\3\6\5\5\6\3\8\2\7\eeee
\3\8\2\7\1\8\3\6\4\4\6\3\8\1\7\2\eeee
\8\1\7\2\8\3\6\5\5\6\3\8\2\7\1\8\eeee
\2\7\1\8\3\6\4\4\6\3\8\1\7\2\8\3\eeee
} 

\baaa
8-61
\eaaa
\bbbb
0&0&0&0&0&0&2&2\\
0&0&0&0&0&0&2&2\\
0&0&0&0&0&2&0&2\\
0&0&0&2&0&0&2&0\\
0&0&0&0&2&2&0&0\\
0&0&2&0&2&0&0&0\\
1&1&0&2&0&0&0&0\\
1&1&2&0&0&0&0&0\\
\ebbb
\parbox{7cm}{ 
\1\7\4\4\7\2\8\3\6\5\5\6\3\8\2\7\eeee
\7\4\4\7\1\8\3\6\5\5\6\3\8\1\7\4\eeee
\4\4\7\2\8\3\6\5\5\6\3\8\2\7\4\4\eeee
\4\7\1\8\3\6\5\5\6\3\8\1\7\4\4\7\eeee
\7\2\8\3\6\5\5\6\3\8\2\7\4\4\7\1\eeee
\1\8\3\6\5\5\6\3\8\1\7\4\4\7\2\8\eeee
\8\3\6\5\5\6\3\8\2\7\4\4\7\1\8\3\eeee
\3\6\5\5\6\3\8\1\7\4\4\7\2\8\3\6\eeee
\6\5\5\6\3\8\2\7\4\4\7\1\8\3\6\5\eeee
\5\5\6\3\8\1\7\4\4\7\2\8\3\6\5\5\eeee
\5\6\3\8\2\7\4\4\7\1\8\3\6\5\5\6\eeee
\6\3\8\1\7\4\4\7\2\8\3\6\5\5\6\3\eeee
} 

\baaa
8-62
\eaaa
\bbbb
0&0&0&0&0&0&2&2\\
0&0&0&0&0&0&2&2\\
0&0&0&0&0&2&1&1\\
0&0&0&0&0&2&1&1\\
0&0&0&0&2&2&0&0\\
0&0&1&1&2&0&0&0\\
1&1&1&1&0&0&0&0\\
1&1&1&1&0&0&0&0\\
\ebbb
\parbox{7cm}{ 
\1\7\3\6\5\5\6\3\7\1\7\3\6\5\5\6\eeee
\7\2\8\4\6\5\5\6\4\8\2\8\4\6\5\5\eeee
\3\8\1\7\3\6\5\5\6\3\7\1\7\3\6\5\eeee
\6\4\7\2\8\4\6\5\5\6\4\8\2\8\4\6\eeee
\5\6\3\8\1\7\3\6\5\5\6\3\7\1\7\3\eeee
\5\5\6\4\7\2\8\4\6\5\5\6\4\8\2\8\eeee
\6\5\5\6\3\8\1\7\3\6\5\5\6\3\7\1\eeee
\3\6\5\5\6\4\7\2\8\4\6\5\5\6\4\8\eeee
\7\4\6\5\5\6\3\8\1\7\3\6\5\5\6\3\eeee
\1\8\3\6\5\5\6\4\7\2\8\4\6\5\5\6\eeee
\7\2\7\4\6\5\5\6\3\8\1\7\3\6\5\5\eeee
\3\8\1\8\3\6\5\5\6\4\7\2\8\4\6\5\eeee
} 

\baaa
8-63
\eaaa
\bbbb
0&0&0&0&0&0&2&2\\
0&0&0&0&0&0&2&2\\
0&0&0&0&0&2&1&1\\
0&0&0&0&2&2&0&0\\
0&0&0&2&2&0&0&0\\
0&0&2&2&0&0&0&0\\
1&1&2&0&0&0&0&0\\
1&1&2&0&0&0&0&0\\
\ebbb
\parbox{7cm}{ 
\1\7\3\6\4\5\5\4\6\3\7\1\7\3\6\4\eeee
\7\2\8\3\6\4\5\5\4\6\3\8\2\8\3\6\eeee
\3\8\1\7\3\6\4\5\5\4\6\3\7\1\7\3\eeee
\6\3\7\2\8\3\6\4\5\5\4\6\3\8\2\8\eeee
\4\6\3\8\1\7\3\6\4\5\5\4\6\3\7\1\eeee
\5\4\6\3\7\2\8\3\6\4\5\5\4\6\3\8\eeee
\5\5\4\6\3\8\1\7\3\6\4\5\5\4\6\3\eeee
\4\5\5\4\6\3\7\2\8\3\6\4\5\5\4\6\eeee
\6\4\5\5\4\6\3\8\1\7\3\6\4\5\5\4\eeee
\3\6\4\5\5\4\6\3\7\2\8\3\6\4\5\5\eeee
\7\3\6\4\5\5\4\6\3\8\1\7\3\6\4\5\eeee
\1\8\3\6\4\5\5\4\6\3\7\2\8\3\6\4\eeee
} 

\baaa
8-64
\eaaa
\bbbb
0&0&0&0&0&0&2&2\\
0&0&0&0&0&0&2&2\\
0&0&0&0&0&2&1&1\\
0&0&0&1&1&2&0&0\\
0&0&0&1&1&2&0&0\\
0&0&2&1&1&0&0&0\\
1&1&2&0&0&0&0&0\\
1&1&2&0&0&0&0&0\\
\ebbb
\parbox{7cm}{ 
\1\7\3\6\4\4\6\3\7\1\7\3\6\4\4\6\eeee
\7\2\8\3\6\5\5\6\3\8\2\8\3\6\5\5\eeee
\3\8\1\7\3\6\4\4\6\3\7\1\7\3\6\4\eeee
\6\3\7\2\8\3\6\5\5\6\3\8\2\8\3\6\eeee
\4\6\3\8\1\7\3\6\4\4\6\3\7\1\7\3\eeee
\4\5\6\3\7\2\8\3\6\5\5\6\3\8\2\8\eeee
\6\5\4\6\3\8\1\7\3\6\4\4\6\3\7\1\eeee
\3\6\4\5\6\3\7\2\8\3\6\5\5\6\3\8\eeee
\7\3\6\5\4\6\3\8\1\7\3\6\4\4\6\3\eeee
\1\8\3\6\4\5\6\3\7\2\8\3\6\5\5\6\eeee
\7\2\7\3\6\5\4\6\3\8\1\7\3\6\4\4\eeee
\3\8\1\8\3\6\4\5\6\3\7\2\8\3\6\5\eeee
} 

\baaa
8-65
\eaaa
\bbbb
0&0&0&0&0&0&2&2\\
0&0&0&0&0&0&2&2\\
0&0&0&0&1&1&0&2\\
0&0&0&0&1&1&0&2\\
0&0&1&1&0&0&2&0\\
0&0&1&1&0&0&2&0\\
1&1&0&0&1&1&0&0\\
1&1&1&1&0&0&0&0\\
\ebbb
\parbox{7cm}{ 
\1\7\6\3\8\2\7\5\3\8\2\7\5\4\8\1\eeee
\7\5\4\8\1\7\6\4\8\1\7\6\3\8\2\7\eeee
\6\3\8\2\7\5\3\8\2\7\5\4\8\1\7\6\eeee
\4\8\1\7\6\4\8\1\7\6\3\8\2\7\5\3\eeee
\8\2\7\5\3\8\2\7\5\4\8\1\7\6\4\8\eeee
\1\7\6\4\8\1\7\6\3\8\2\7\5\3\8\2\eeee
\7\5\3\8\2\7\5\4\8\1\7\6\4\8\1\7\eeee
\6\4\8\1\7\6\3\8\2\7\5\3\8\2\7\5\eeee
\3\8\2\7\5\4\8\1\7\6\4\8\1\7\6\3\eeee
\8\1\7\6\3\8\2\7\5\3\8\2\7\5\4\8\eeee
\2\7\5\4\8\1\7\6\4\8\1\7\6\3\8\2\eeee
\7\6\3\8\2\7\5\3\8\2\7\5\4\8\1\7\eeee
} 

\baaa
8-66
\eaaa
\bbbb
0&0&0&0&0&0&2&2\\
0&0&0&0&0&0&2&2\\
0&0&0&0&1&1&0&2\\
0&0&0&0&1&1&0&2\\
0&0&1&1&1&1&0&0\\
0&0&1&1&1&1&0&0\\
1&1&0&0&0&0&2&0\\
1&1&1&1&0&0&0&0\\
\ebbb
\parbox{7cm}{ 
\1\7\7\1\8\4\6\5\3\8\2\7\7\2\8\3\eeee
\7\7\2\8\3\5\6\4\8\1\7\7\1\8\4\6\eeee
\7\1\8\4\6\5\3\8\2\7\7\2\8\3\5\6\eeee
\2\8\3\5\6\4\8\1\7\7\1\8\4\6\5\3\eeee
\8\4\6\5\3\8\2\7\7\2\8\3\5\6\4\8\eeee
\3\5\6\4\8\1\7\7\1\8\4\6\5\3\8\2\eeee
\6\5\3\8\2\7\7\2\8\3\5\6\4\8\1\7\eeee
\6\4\8\1\7\7\1\8\4\6\5\3\8\2\7\7\eeee
\3\8\2\7\7\2\8\3\5\6\4\8\1\7\7\1\eeee
\8\1\7\7\1\8\4\6\5\3\8\2\7\7\2\8\eeee
\2\7\7\2\8\3\5\6\4\8\1\7\7\1\8\4\eeee
\7\7\1\8\4\6\5\3\8\2\7\7\2\8\3\5\eeee
} 

\baaa
8-67
\eaaa
\bbbb
0&0&0&0&0&0&2&2\\
0&0&0&0&0&0&2&2\\
0&0&0&0&1&1&0&2\\
0&0&0&0&1&1&0&2\\
0&0&1&1&1&1&0&0\\
0&0&1&1&1&1&0&0\\
2&2&0&0&0&0&0&0\\
1&1&1&1&0&0&0&0\\
\ebbb
\parbox{7cm}{ 
\1\7\2\8\3\6\6\3\8\2\7\1\8\4\6\5\eeee
\7\1\8\4\5\5\4\8\1\7\2\8\3\5\6\4\eeee
\2\8\3\6\6\3\8\2\7\1\8\4\6\5\3\8\eeee
\8\4\5\5\4\8\1\7\2\8\3\5\6\4\8\1\eeee
\3\6\6\3\8\2\7\1\8\4\6\5\3\8\2\7\eeee
\5\5\4\8\1\7\2\8\3\5\6\4\8\1\7\2\eeee
\6\3\8\2\7\1\8\4\6\5\3\8\2\7\1\8\eeee
\4\8\1\7\2\8\3\5\6\4\8\1\7\2\8\3\eeee
\8\2\7\1\8\4\6\5\3\8\2\7\1\8\4\5\eeee
\1\7\2\8\3\5\6\4\8\1\7\2\8\3\6\6\eeee
\7\1\8\4\6\5\3\8\2\7\1\8\4\5\5\4\eeee
\2\8\3\5\6\4\8\1\7\2\8\3\6\6\3\8\eeee
} 

\baaa
8-68
\eaaa
\bbbb
0&0&0&0&0&0&2&2\\
0&0&0&0&0&0&2&2\\
0&0&0&0&1&1&0&2\\
0&0&0&0&1&1&0&2\\
0&0&2&2&0&0&0&0\\
0&0&2&2&0&0&0&0\\
1&1&0&0&0&0&2&0\\
1&1&1&1&0&0&0&0\\
\ebbb
\parbox{7cm}{ 
\1\7\7\1\8\4\6\3\8\2\7\7\1\8\4\5\eeee
\7\7\2\8\3\5\4\8\1\7\7\2\8\3\6\4\eeee
\7\1\8\4\6\3\8\2\7\7\1\8\4\5\3\8\eeee
\2\8\3\5\4\8\1\7\7\2\8\3\6\4\8\1\eeee
\8\4\6\3\8\2\7\7\1\8\4\5\3\8\2\7\eeee
\3\5\4\8\1\7\7\2\8\3\6\4\8\1\7\7\eeee
\6\3\8\2\7\7\1\8\4\5\3\8\2\7\7\2\eeee
\4\8\1\7\7\2\8\3\6\4\8\1\7\7\1\8\eeee
\8\2\7\7\1\8\4\5\3\8\2\7\7\2\8\3\eeee
\1\7\7\2\8\3\6\4\8\1\7\7\1\8\4\6\eeee
\7\7\1\8\4\5\3\8\2\7\7\2\8\3\5\4\eeee
\7\2\8\3\6\4\8\1\7\7\1\8\4\6\3\8\eeee
} 

\baaa
8-69
\eaaa
\bbbb
0&0&0&0&0&0&2&2\\
0&0&0&0&0&0&2&2\\
0&0&0&0&1&1&0&2\\
0&0&0&0&1&1&0&2\\
0&0&2&2&0&0&0&0\\
0&0&2&2&0&0&0&0\\
2&2&0&0&0&0&0&0\\
1&1&1&1&0&0&0&0\\
\ebbb
\parbox{7cm}{ 
\1\7\2\8\3\6\4\8\1\7\2\8\3\6\4\8\eeee
\7\1\8\4\5\3\8\2\7\1\8\4\5\3\8\2\eeee
\2\8\3\6\4\8\1\7\2\8\3\6\4\8\1\7\eeee
\8\4\5\3\8\2\7\1\8\4\5\3\8\2\7\1\eeee
\3\6\4\8\1\7\2\8\3\6\4\8\1\7\2\8\eeee
\5\3\8\2\7\1\8\4\5\3\8\2\7\1\8\4\eeee
\4\8\1\7\2\8\3\6\4\8\1\7\2\8\3\6\eeee
\8\2\7\1\8\4\5\3\8\2\7\1\8\4\5\3\eeee
\1\7\2\8\3\6\4\8\1\7\2\8\3\6\4\8\eeee
\7\1\8\4\5\3\8\2\7\1\8\4\5\3\8\2\eeee
\2\8\3\6\4\8\1\7\2\8\3\6\4\8\1\7\eeee
\8\4\5\3\8\2\7\1\8\4\5\3\8\2\7\1\eeee
} 

\baaa
8-70
\eaaa
\bbbb
0&0&0&0&0&0&2&2\\
0&0&0&0&0&0&2&2\\
0&0&0&0&1&1&0&2\\
0&0&0&0&1&1&2&0\\
0&0&2&2&0&0&0&0\\
0&0&2&2&0&0&0&0\\
1&1&0&2&0&0&0&0\\
1&1&2&0&0&0&0&0\\
\ebbb
\parbox{7cm}{ 
\1\7\4\5\3\8\2\7\4\6\3\8\1\7\4\5\eeee
\7\4\6\3\8\1\7\4\5\3\8\2\7\4\6\3\eeee
\4\5\3\8\2\7\4\6\3\8\1\7\4\5\3\8\eeee
\6\3\8\1\7\4\5\3\8\2\7\4\6\3\8\1\eeee
\3\8\2\7\4\6\3\8\1\7\4\5\3\8\2\7\eeee
\8\1\7\4\5\3\8\2\7\4\6\3\8\1\7\4\eeee
\2\7\4\6\3\8\1\7\4\5\3\8\2\7\4\6\eeee
\7\4\5\3\8\2\7\4\6\3\8\1\7\4\5\3\eeee
\4\6\3\8\1\7\4\5\3\8\2\7\4\6\3\8\eeee
\5\3\8\2\7\4\6\3\8\1\7\4\5\3\8\2\eeee
\3\8\1\7\4\5\3\8\2\7\4\6\3\8\1\7\eeee
\8\2\7\4\6\3\8\1\7\4\5\3\8\2\7\4\eeee
} 

\baaa
8-71
\eaaa
\bbbb
0&0&0&0&0&0&2&2\\
0&0&0&0&0&0&2&2\\
0&0&0&0&1&1&0&2\\
0&0&0&0&2&2&0&0\\
0&0&2&2&0&0&0&0\\
0&0&2&2&0&0&0&0\\
1&1&0&0&0&0&2&0\\
1&1&2&0&0&0&0&0\\
\ebbb
\parbox{7cm}{ 
\1\7\7\1\8\3\6\4\5\3\8\1\7\7\2\8\eeee
\7\7\2\8\3\5\4\6\3\8\2\7\7\1\8\3\eeee
\7\1\8\3\6\4\5\3\8\1\7\7\2\8\3\6\eeee
\2\8\3\5\4\6\3\8\2\7\7\1\8\3\5\4\eeee
\8\3\6\4\5\3\8\1\7\7\2\8\3\6\4\5\eeee
\3\5\4\6\3\8\2\7\7\1\8\3\5\4\6\3\eeee
\6\4\5\3\8\1\7\7\2\8\3\6\4\5\3\8\eeee
\4\6\3\8\2\7\7\1\8\3\5\4\6\3\8\1\eeee
\5\3\8\1\7\7\2\8\3\6\4\5\3\8\2\7\eeee
\3\8\2\7\7\1\8\3\5\4\6\3\8\1\7\7\eeee
\8\1\7\7\2\8\3\6\4\5\3\8\2\7\7\2\eeee
\2\7\7\1\8\3\5\4\6\3\8\1\7\7\1\8\eeee
} 

\baaa
8-72
\eaaa
\bbbb
0&0&0&0&0&0&2&2\\
0&0&0&0&0&0&2&2\\
0&0&0&0&1&1&0&2\\
0&0&0&0&2&2&0&0\\
0&0&2&2&0&0&0&0\\
0&0&2&2&0&0&0&0\\
2&2&0&0&0&0&0&0\\
1&1&2&0&0&0&0&0\\
\ebbb
\parbox{7cm}{ 
\1\7\2\8\3\6\4\5\3\8\2\7\1\8\3\5\eeee
\7\1\8\3\5\4\6\3\8\1\7\2\8\3\6\4\eeee
\2\8\3\6\4\5\3\8\2\7\1\8\3\5\4\6\eeee
\8\3\5\4\6\3\8\1\7\2\8\3\6\4\5\3\eeee
\3\6\4\5\3\8\2\7\1\8\3\5\4\6\3\8\eeee
\5\4\6\3\8\1\7\2\8\3\6\4\5\3\8\2\eeee
\4\5\3\8\2\7\1\8\3\5\4\6\3\8\1\7\eeee
\6\3\8\1\7\2\8\3\6\4\5\3\8\2\7\1\eeee
\3\8\2\7\1\8\3\5\4\6\3\8\1\7\2\8\eeee
\8\1\7\2\8\3\6\4\5\3\8\2\7\1\8\3\eeee
\2\7\1\8\3\5\4\6\3\8\1\7\2\8\3\6\eeee
\7\2\8\3\6\4\5\3\8\2\7\1\8\3\5\4\eeee
} 

\baaa
8-73
\eaaa
\bbbb
0&0&0&0&0&0&2&2\\
0&0&0&0&0&0&2&2\\
0&0&0&0&1&1&0&2\\
0&0&0&2&0&0&2&0\\
0&0&2&0&1&1&0&0\\
0&0&2&0&1&1&0&0\\
1&1&0&2&0&0&0&0\\
1&1&2&0&0&0&0&0\\
\ebbb
\parbox{7cm}{ 
\1\7\4\4\7\2\8\3\5\6\3\8\1\7\4\4\eeee
\7\4\4\7\1\8\3\6\5\3\8\2\7\4\4\7\eeee
\4\4\7\2\8\3\5\6\3\8\1\7\4\4\7\2\eeee
\4\7\1\8\3\6\5\3\8\2\7\4\4\7\1\8\eeee
\7\2\8\3\5\6\3\8\1\7\4\4\7\2\8\3\eeee
\1\8\3\6\5\3\8\2\7\4\4\7\1\8\3\6\eeee
\8\3\5\6\3\8\1\7\4\4\7\2\8\3\5\6\eeee
\3\6\5\3\8\2\7\4\4\7\1\8\3\6\5\3\eeee
\5\6\3\8\1\7\4\4\7\2\8\3\5\6\3\8\eeee
\5\3\8\2\7\4\4\7\1\8\3\6\5\3\8\2\eeee
\3\8\1\7\4\4\7\2\8\3\5\6\3\8\1\7\eeee
\8\2\7\4\4\7\1\8\3\6\5\3\8\2\7\4\eeee
} 

\baaa
8-74
\eaaa
\bbbb
0&0&0&0&0&0&2&2\\
0&0&0&0&0&0&2&2\\
0&0&0&0&1&1&0&2\\
0&0&0&2&1&1&0&0\\
0&0&2&2&0&0&0&0\\
0&0&2&2&0&0&0&0\\
1&1&0&0&0&0&2&0\\
1&1&2&0&0&0&0&0\\
\ebbb
\parbox{7cm}{ 
\1\7\7\1\8\3\6\4\4\6\3\8\1\7\7\1\eeee
\7\7\2\8\3\5\4\4\5\3\8\2\7\7\2\8\eeee
\7\1\8\3\6\4\4\6\3\8\1\7\7\1\8\3\eeee
\2\8\3\5\4\4\5\3\8\2\7\7\2\8\3\5\eeee
\8\3\6\4\4\6\3\8\1\7\7\1\8\3\6\4\eeee
\3\5\4\4\5\3\8\2\7\7\2\8\3\5\4\4\eeee
\6\4\4\6\3\8\1\7\7\1\8\3\6\4\4\6\eeee
\4\4\5\3\8\2\7\7\2\8\3\5\4\4\5\3\eeee
\4\6\3\8\1\7\7\1\8\3\6\4\4\6\3\8\eeee
\5\3\8\2\7\7\2\8\3\5\4\4\5\3\8\2\eeee
\3\8\1\7\7\1\8\3\6\4\4\6\3\8\1\7\eeee
\8\2\7\7\2\8\3\5\4\4\5\3\8\2\7\7\eeee
} 

\baaa
8-75
\eaaa
\bbbb
0&0&0&0&0&0&2&2\\
0&0&0&0&0&0&2&2\\
0&0&0&0&1&1&1&1\\
0&0&0&0&1&1&1&1\\
0&0&1&1&1&1&0&0\\
0&0&1&1&1&1&0&0\\
1&1&1&1&0&0&0&0\\
1&1&1&1&0&0&0&0\\
\ebbb
\parbox{7cm}{ 
\1\7\3\5\5\3\7\1\7\3\5\5\3\7\1\7\eeee
\7\2\8\4\6\6\4\8\2\8\4\6\6\4\8\2\eeee
\3\8\1\7\3\5\5\3\7\1\7\3\5\5\3\7\eeee
\5\4\7\2\8\4\6\6\4\8\2\8\4\6\6\4\eeee
\5\6\3\8\1\7\3\5\5\3\7\1\7\3\5\5\eeee
\3\6\5\4\7\2\8\4\6\6\4\8\2\8\4\6\eeee
\7\4\5\6\3\8\1\7\3\5\5\3\7\1\7\3\eeee
\1\8\3\6\5\4\7\2\8\4\6\6\4\8\2\8\eeee
\7\2\7\4\5\6\3\8\1\7\3\5\5\3\7\1\eeee
\3\8\1\8\3\6\5\4\7\2\8\4\6\6\4\8\eeee
\5\4\7\2\7\4\5\6\3\8\1\7\3\5\5\3\eeee
\5\6\3\8\1\8\3\6\5\4\7\2\8\4\6\6\eeee
} 

\baaa
8-76
\eaaa
\bbbb
0&0&0&0&0&0&2&2\\
0&0&0&0&0&0&2&2\\
0&0&0&0&1&1&1&1\\
0&0&0&0&1&1&1&1\\
0&0&2&2&0&0&0&0\\
0&0&2&2&0&0&0&0\\
1&1&1&1&0&0&0&0\\
1&1&1&1&0&0&0&0\\
\ebbb
\parbox{7cm}{ 
\1\7\3\5\3\7\1\7\3\5\3\7\1\7\3\5\eeee
\7\2\8\4\6\4\8\2\8\4\6\4\8\2\8\4\eeee
\3\8\1\7\3\5\3\7\1\7\3\5\3\7\1\7\eeee
\5\4\7\2\8\4\6\4\8\2\8\4\6\4\8\2\eeee
\3\6\3\8\1\7\3\5\3\7\1\7\3\5\3\7\eeee
\7\4\5\4\7\2\8\4\6\4\8\2\8\4\6\4\eeee
\1\8\3\6\3\8\1\7\3\5\3\7\1\7\3\5\eeee
\7\2\7\4\5\4\7\2\8\4\6\4\8\2\8\4\eeee
\3\8\1\8\3\6\3\8\1\7\3\5\3\7\1\7\eeee
\5\4\7\2\7\4\5\4\7\2\8\4\6\4\8\2\eeee
\3\6\3\8\1\8\3\6\3\8\1\7\3\5\3\7\eeee
\7\4\5\4\7\2\7\4\5\4\7\2\8\4\6\4\eeee
} 

\baaa
8-77
\eaaa
\bbbb
0&0&0&0&0&0&2&2\\
0&0&0&0&0&0&2&2\\
0&0&1&0&0&1&0&2\\
0&0&0&1&1&0&2&0\\
0&0&0&1&1&0&2&0\\
0&0&1&0&0&1&0&2\\
1&1&0&1&1&0&0&0\\
1&1&1&0&0&1&0&0\\
\ebbb
\parbox{7cm}{ 
\1\7\5\4\7\2\8\3\3\8\2\7\4\5\7\1\eeee
\7\4\5\7\1\8\6\6\8\1\7\5\4\7\2\8\eeee
\5\4\7\2\8\3\3\8\2\7\4\5\7\1\8\6\eeee
\5\7\1\8\6\6\8\1\7\5\4\7\2\8\3\3\eeee
\7\2\8\3\3\8\2\7\4\5\7\1\8\6\6\8\eeee
\1\8\6\6\8\1\7\5\4\7\2\8\3\3\8\2\eeee
\8\3\3\8\2\7\4\5\7\1\8\6\6\8\1\7\eeee
\6\6\8\1\7\5\4\7\2\8\3\3\8\2\7\4\eeee
\3\8\2\7\4\5\7\1\8\6\6\8\1\7\5\4\eeee
\8\1\7\5\4\7\2\8\3\3\8\2\7\4\5\7\eeee
\2\7\4\5\7\1\8\6\6\8\1\7\5\4\7\2\eeee
\7\5\4\7\2\8\3\3\8\2\7\4\5\7\1\8\eeee
} 

\baaa
8-78
\eaaa
\bbbb
0&0&0&0&0&0&2&2\\
0&0&0&0&0&1&1&2\\
0&0&0&0&1&0&1&2\\
0&0&0&0&1&1&0&2\\
0&0&2&2&0&0&0&0\\
0&2&0&2&0&0&0&0\\
2&1&1&0&0&0&0&0\\
1&1&1&1&0&0&0&0\\
\ebbb
\parbox{7cm}{ 
\1\8\4\8\1\8\4\8\1\8\4\8\1\8\4\8\eeee
\7\2\6\2\7\3\5\3\7\2\6\2\7\3\5\3\eeee
\1\8\4\8\1\8\4\8\1\8\4\8\1\8\4\8\eeee
\7\3\5\3\7\2\6\2\7\3\5\3\7\2\6\2\eeee
\1\8\4\8\1\8\4\8\1\8\4\8\1\8\4\8\eeee
\7\2\6\2\7\3\5\3\7\2\6\2\7\3\5\3\eeee
\1\8\4\8\1\8\4\8\1\8\4\8\1\8\4\8\eeee
\7\3\5\3\7\2\6\2\7\3\5\3\7\2\6\2\eeee
\1\8\4\8\1\8\4\8\1\8\4\8\1\8\4\8\eeee
\7\2\6\2\7\3\5\3\7\2\6\2\7\3\5\3\eeee
\1\8\4\8\1\8\4\8\1\8\4\8\1\8\4\8\eeee
\7\3\5\3\7\2\6\2\7\3\5\3\7\2\6\2\eeee
} 

\baaa
8-79
\eaaa
\bbbb
0&0&0&0&0&0&2&2\\
0&0&0&0&0&1&1&2\\
0&0&0&0&1&0&1&2\\
0&0&0&1&1&1&0&1\\
0&0&1&2&0&1&0&0\\
0&1&0&2&1&0&0&0\\
2&1&1&0&0&0&0&0\\
1&1&1&1&0&0&0&0\\
\ebbb
\parbox{7cm}{ 
\1\8\4\4\8\1\8\4\4\8\1\8\4\4\8\1\eeee
\7\2\6\5\3\7\2\6\5\3\7\2\6\5\3\7\eeee
\1\8\4\4\8\1\8\4\4\8\1\8\4\4\8\1\eeee
\7\3\5\6\2\7\3\5\6\2\7\3\5\6\2\7\eeee
\1\8\4\4\8\1\8\4\4\8\1\8\4\4\8\1\eeee
\7\2\6\5\3\7\2\6\5\3\7\2\6\5\3\7\eeee
\1\8\4\4\8\1\8\4\4\8\1\8\4\4\8\1\eeee
\7\3\5\6\2\7\3\5\6\2\7\3\5\6\2\7\eeee
\1\8\4\4\8\1\8\4\4\8\1\8\4\4\8\1\eeee
\7\2\6\5\3\7\2\6\5\3\7\2\6\5\3\7\eeee
\1\8\4\4\8\1\8\4\4\8\1\8\4\4\8\1\eeee
\7\3\5\6\2\7\3\5\6\2\7\3\5\6\2\7\eeee
} 

\baaa
8-80
\eaaa
\bbbb
0&0&0&0&0&0&2&2\\
0&0&0&0&0&1&1&2\\
0&0&0&0&1&0&1&2\\
0&0&0&1&1&1&0&1\\
0&0&1&2&1&0&0&0\\
0&1&0&2&0&1&0&0\\
2&1&1&0&0&0&0&0\\
1&1&1&1&0&0&0&0\\
\ebbb
\parbox{7cm}{ 
\1\8\4\4\8\1\8\4\4\8\1\8\4\4\8\1\eeee
\7\2\6\6\2\7\3\5\5\3\7\2\6\6\2\7\eeee
\1\8\4\4\8\1\8\4\4\8\1\8\4\4\8\1\eeee
\7\3\5\5\3\7\2\6\6\2\7\3\5\5\3\7\eeee
\1\8\4\4\8\1\8\4\4\8\1\8\4\4\8\1\eeee
\7\2\6\6\2\7\3\5\5\3\7\2\6\6\2\7\eeee
\1\8\4\4\8\1\8\4\4\8\1\8\4\4\8\1\eeee
\7\3\5\5\3\7\2\6\6\2\7\3\5\5\3\7\eeee
\1\8\4\4\8\1\8\4\4\8\1\8\4\4\8\1\eeee
\7\2\6\6\2\7\3\5\5\3\7\2\6\6\2\7\eeee
\1\8\4\4\8\1\8\4\4\8\1\8\4\4\8\1\eeee
\7\3\5\5\3\7\2\6\6\2\7\3\5\5\3\7\eeee
} 

\baaa
8-81
\eaaa
\bbbb
0&0&0&0&0&0&2&2\\
0&0&0&0&0&1&1&2\\
0&0&0&0&1&1&1&1\\
0&0&0&0&1&2&1&0\\
0&0&2&2&0&0&0&0\\
0&1&1&2&0&0&0&0\\
1&1&1&1&0&0&0&0\\
1&2&1&0&0&0&0&0\\
\ebbb
\parbox{7cm}{ 
\1\8\2\7\3\6\4\5\4\6\3\7\2\8\1\8\eeee
\7\2\8\1\8\2\7\3\6\4\5\4\6\3\7\2\eeee
\4\6\3\7\2\8\1\8\2\7\3\6\4\5\4\6\eeee
\6\4\5\4\6\3\7\2\8\1\8\2\7\3\6\4\eeee
\2\7\3\6\4\5\4\6\3\7\2\8\1\8\2\7\eeee
\8\1\8\2\7\3\6\4\5\4\6\3\7\2\8\1\eeee
\3\7\2\8\1\8\2\7\3\6\4\5\4\6\3\7\eeee
\5\4\6\3\7\2\8\1\8\2\7\3\6\4\5\4\eeee
\3\6\4\5\4\6\3\7\2\8\1\8\2\7\3\6\eeee
\8\2\7\3\6\4\5\4\6\3\7\2\8\1\8\2\eeee
\2\8\1\8\2\7\3\6\4\5\4\6\3\7\2\8\eeee
\6\3\7\2\8\1\8\2\7\3\6\4\5\4\6\3\eeee
} 

\baaa
8-82
\eaaa
\bbbb
0&0&0&0&0&0&2&2\\
0&0&0&0&0&1&1&2\\
0&0&0&0&1&1&1&1\\
0&0&0&1&1&1&1&0\\
0&0&1&1&1&1&0&0\\
0&1&1&1&1&0&0&0\\
1&1&1&1&0&0&0&0\\
1&2&1&0&0&0&0&0\\
\ebbb
\parbox{7cm}{ 
\1\8\2\7\3\6\4\5\5\4\6\3\7\2\8\1\eeee
\7\2\8\1\8\2\7\3\6\4\5\5\4\6\3\7\eeee
\4\6\3\7\2\8\1\8\2\7\3\6\4\5\5\4\eeee
\4\5\5\4\6\3\7\2\8\1\8\2\7\3\6\4\eeee
\7\3\6\4\5\5\4\6\3\7\2\8\1\8\2\7\eeee
\1\8\2\7\3\6\4\5\5\4\6\3\7\2\8\1\eeee
\7\2\8\1\8\2\7\3\6\4\5\5\4\6\3\7\eeee
\4\6\3\7\2\8\1\8\2\7\3\6\4\5\5\4\eeee
\4\5\5\4\6\3\7\2\8\1\8\2\7\3\6\4\eeee
\7\3\6\4\5\5\4\6\3\7\2\8\1\8\2\7\eeee
\1\8\2\7\3\6\4\5\5\4\6\3\7\2\8\1\eeee
\7\2\8\1\8\2\7\3\6\4\5\5\4\6\3\7\eeee
} 

\baaa
8-83
\eaaa
\bbbb
0&0&0&0&0&0&2&2\\
0&0&0&0&0&1&1&2\\
0&0&0&0&1&2&0&1\\
0&0&0&0&2&2&0&0\\
0&0&2&2&0&0&0&0\\
0&1&2&1&0&0&0&0\\
2&2&0&0&0&0&0&0\\
1&2&1&0&0&0&0&0\\
\ebbb
\parbox{7cm}{ 
\1\8\3\5\3\8\1\8\3\5\3\8\1\8\3\5\eeee
\7\2\6\4\6\2\7\2\6\4\6\2\7\2\6\4\eeee
\1\8\3\5\3\8\1\8\3\5\3\8\1\8\3\5\eeee
\7\2\6\4\6\2\7\2\6\4\6\2\7\2\6\4\eeee
\1\8\3\5\3\8\1\8\3\5\3\8\1\8\3\5\eeee
\7\2\6\4\6\2\7\2\6\4\6\2\7\2\6\4\eeee
\1\8\3\5\3\8\1\8\3\5\3\8\1\8\3\5\eeee
\7\2\6\4\6\2\7\2\6\4\6\2\7\2\6\4\eeee
\1\8\3\5\3\8\1\8\3\5\3\8\1\8\3\5\eeee
\7\2\6\4\6\2\7\2\6\4\6\2\7\2\6\4\eeee
\1\8\3\5\3\8\1\8\3\5\3\8\1\8\3\5\eeee
\7\2\6\4\6\2\7\2\6\4\6\2\7\2\6\4\eeee
} 

\baaa
8-84
\eaaa
\bbbb
0&0&0&0&0&0&2&2\\
0&0&0&0&0&1&1&2\\
0&0&0&0&1&2&0&1\\
0&0&0&1&2&1&0&0\\
0&0&1&2&1&0&0&0\\
0&1&2&1&0&0&0&0\\
2&2&0&0&0&0&0&0\\
1&2&1&0&0&0&0&0\\
\ebbb
\parbox{7cm}{ 
\1\8\3\5\5\3\8\1\8\3\5\5\3\8\1\8\eeee
\7\2\6\4\4\6\2\7\2\6\4\4\6\2\7\2\eeee
\1\8\3\5\5\3\8\1\8\3\5\5\3\8\1\8\eeee
\7\2\6\4\4\6\2\7\2\6\4\4\6\2\7\2\eeee
\1\8\3\5\5\3\8\1\8\3\5\5\3\8\1\8\eeee
\7\2\6\4\4\6\2\7\2\6\4\4\6\2\7\2\eeee
\1\8\3\5\5\3\8\1\8\3\5\5\3\8\1\8\eeee
\7\2\6\4\4\6\2\7\2\6\4\4\6\2\7\2\eeee
\1\8\3\5\5\3\8\1\8\3\5\5\3\8\1\8\eeee
\7\2\6\4\4\6\2\7\2\6\4\4\6\2\7\2\eeee
\1\8\3\5\5\3\8\1\8\3\5\5\3\8\1\8\eeee
\7\2\6\4\4\6\2\7\2\6\4\4\6\2\7\2\eeee
} 

\baaa
8-85
\eaaa
\bbbb
0&0&0&0&0&0&2&2\\
0&0&0&0&0&2&0&2\\
0&0&0&0&1&1&1&1\\
0&0&0&0&2&0&2&0\\
0&0&2&2&0&0&0&0\\
0&2&2&0&0&0&0&0\\
2&0&1&1&0&0&0&0\\
2&1&1&0&0&0&0&0\\
\ebbb
\parbox{7cm}{ 
\1\7\3\6\2\8\1\8\2\6\3\7\1\8\3\5\eeee
\7\4\5\3\8\1\7\3\6\2\8\1\8\2\6\3\eeee
\3\5\4\7\1\7\4\5\3\8\1\7\3\6\2\8\eeee
\6\3\7\1\8\3\5\4\7\1\7\4\5\3\8\1\eeee
\2\8\1\8\2\6\3\7\1\8\3\5\4\7\1\7\eeee
\8\1\7\3\6\2\8\1\8\2\6\3\7\1\8\3\eeee
\1\7\4\5\3\8\1\7\3\6\2\8\1\8\2\6\eeee
\8\3\5\4\7\1\7\4\5\3\8\1\7\3\6\2\eeee
\2\6\3\7\1\8\3\5\4\7\1\7\4\5\3\8\eeee
\6\2\8\1\8\2\6\3\7\1\8\3\5\4\7\1\eeee
\3\8\1\7\3\6\2\8\1\8\2\6\3\7\1\8\eeee
\7\1\7\4\5\3\8\1\7\3\6\2\8\1\8\2\eeee
} 

\baaa
8-86
\eaaa
\bbbb
0&0&0&0&0&0&2&2\\
0&0&0&0&0&2&0&2\\
0&0&0&0&2&0&2&0\\
0&0&0&0&2&2&0&0\\
0&0&1&1&0&0&1&1\\
0&1&0&1&0&1&0&1\\
1&0&1&0&1&0&1&0\\
1&1&0&0&1&1&0&0\\
\ebbb
\parbox{7cm}{ 
\1\7\7\1\8\6\2\8\5\3\7\5\4\6\6\4\eeee
\7\3\5\8\2\6\8\1\7\7\1\8\6\2\8\5\eeee
\7\5\4\6\6\4\5\7\3\5\8\2\6\8\1\7\eeee
\1\8\6\2\8\5\3\7\5\4\6\6\4\5\7\3\eeee
\8\2\6\8\1\7\7\1\8\6\2\8\5\3\7\5\eeee
\6\6\4\5\7\3\5\8\2\6\8\1\7\7\1\8\eeee
\2\8\5\3\7\5\4\6\6\4\5\7\3\5\8\2\eeee
\8\1\7\7\1\8\6\2\8\5\3\7\5\4\6\6\eeee
\5\7\3\5\8\2\6\8\1\7\7\1\8\6\2\8\eeee
\3\7\5\4\6\6\4\5\7\3\5\8\2\6\8\1\eeee
\7\1\8\6\2\8\5\3\7\5\4\6\6\4\5\7\eeee
\5\8\2\6\8\1\7\7\1\8\6\2\8\5\3\7\eeee
} 

\baaa
8-87
\eaaa
\bbbb
0&0&0&0&0&0&2&2\\
0&0&0&0&0&2&0&2\\
0&0&0&0&2&0&2&0\\
0&0&0&0&2&2&0&0\\
0&0&2&2&0&0&0&0\\
0&2&0&2&0&0&0&0\\
2&0&2&0&0&0&0&0\\
2&2&0&0&0&0&0&0\\
\ebbb
\parbox{7cm}{ 
\1\7\3\5\4\6\2\8\1\7\3\5\4\6\2\8\eeee
\7\3\5\4\6\2\8\1\7\3\5\4\6\2\8\1\eeee
\3\5\4\6\2\8\1\7\3\5\4\6\2\8\1\7\eeee
\5\4\6\2\8\1\7\3\5\4\6\2\8\1\7\3\eeee
\4\6\2\8\1\7\3\5\4\6\2\8\1\7\3\5\eeee
\6\2\8\1\7\3\5\4\6\2\8\1\7\3\5\4\eeee
\2\8\1\7\3\5\4\6\2\8\1\7\3\5\4\6\eeee
\8\1\7\3\5\4\6\2\8\1\7\3\5\4\6\2\eeee
\1\7\3\5\4\6\2\8\1\7\3\5\4\6\2\8\eeee
\7\3\5\4\6\2\8\1\7\3\5\4\6\2\8\1\eeee
\3\5\4\6\2\8\1\7\3\5\4\6\2\8\1\7\eeee
\5\4\6\2\8\1\7\3\5\4\6\2\8\1\7\3\eeee
} 

\baaa
8-88
\eaaa
\bbbb
0&0&0&0&0&0&2&2\\
0&0&0&0&0&2&0&2\\
0&0&0&0&2&0&2&0\\
0&0&0&2&0&2&0&0\\
0&0&2&0&2&0&0&0\\
0&2&0&2&0&0&0&0\\
2&0&2&0&0&0&0&0\\
2&2&0&0&0&0&0&0\\
\ebbb
\parbox{7cm}{ 
\1\7\3\5\5\3\7\1\8\2\6\4\4\6\2\8\eeee
\7\3\5\5\3\7\1\8\2\6\4\4\6\2\8\1\eeee
\3\5\5\3\7\1\8\2\6\4\4\6\2\8\1\7\eeee
\5\5\3\7\1\8\2\6\4\4\6\2\8\1\7\3\eeee
\5\3\7\1\8\2\6\4\4\6\2\8\1\7\3\5\eeee
\3\7\1\8\2\6\4\4\6\2\8\1\7\3\5\5\eeee
\7\1\8\2\6\4\4\6\2\8\1\7\3\5\5\3\eeee
\1\8\2\6\4\4\6\2\8\1\7\3\5\5\3\7\eeee
\8\2\6\4\4\6\2\8\1\7\3\5\5\3\7\1\eeee
\2\6\4\4\6\2\8\1\7\3\5\5\3\7\1\8\eeee
\6\4\4\6\2\8\1\7\3\5\5\3\7\1\8\2\eeee
\4\4\6\2\8\1\7\3\5\5\3\7\1\8\2\6\eeee
} 

\baaa
8-89
\eaaa
\bbbb
0&0&0&0&0&0&2&2\\
0&0&0&0&0&2&0&2\\
0&0&0&1&1&0&2&0\\
0&0&2&1&1&0&0&0\\
0&0&2&1&1&0&0&0\\
0&2&0&0&0&2&0&0\\
2&0&2&0&0&0&0&0\\
2&2&0&0&0&0&0&0\\
\ebbb
\parbox{7cm}{ 
\1\7\3\4\4\3\7\1\8\2\6\6\2\8\1\7\eeee
\7\3\5\5\3\7\1\8\2\6\6\2\8\1\7\3\eeee
\3\4\4\3\7\1\8\2\6\6\2\8\1\7\3\5\eeee
\5\5\3\7\1\8\2\6\6\2\8\1\7\3\4\4\eeee
\4\3\7\1\8\2\6\6\2\8\1\7\3\5\5\3\eeee
\3\7\1\8\2\6\6\2\8\1\7\3\4\4\3\7\eeee
\7\1\8\2\6\6\2\8\1\7\3\5\5\3\7\1\eeee
\1\8\2\6\6\2\8\1\7\3\4\4\3\7\1\8\eeee
\8\2\6\6\2\8\1\7\3\5\5\3\7\1\8\2\eeee
\2\6\6\2\8\1\7\3\4\4\3\7\1\8\2\6\eeee
\6\6\2\8\1\7\3\5\5\3\7\1\8\2\6\6\eeee
\6\2\8\1\7\3\4\4\3\7\1\8\2\6\6\2\eeee
} 

\baaa
8-90
\eaaa
\bbbb
0&0&0&0&0&0&2&2\\
0&0&0&0&0&2&1&1\\
0&0&0&0&0&2&1&1\\
0&0&0&0&2&2&0&0\\
0&0&0&2&2&0&0&0\\
0&1&1&2&0&0&0&0\\
2&1&1&0&0&0&0&0\\
2&1&1&0&0&0&0&0\\
\ebbb
\parbox{7cm}{ 
\1\7\2\6\4\5\5\4\6\2\7\1\7\2\6\4\eeee
\7\1\8\3\6\4\5\5\4\6\3\8\1\8\3\6\eeee
\2\8\1\7\2\6\4\5\5\4\6\2\7\1\7\2\eeee
\6\3\7\1\8\3\6\4\5\5\4\6\3\8\1\8\eeee
\4\6\2\8\1\7\2\6\4\5\5\4\6\2\7\1\eeee
\5\4\6\3\7\1\8\3\6\4\5\5\4\6\3\8\eeee
\5\5\4\6\2\8\1\7\2\6\4\5\5\4\6\2\eeee
\4\5\5\4\6\3\7\1\8\3\6\4\5\5\4\6\eeee
\6\4\5\5\4\6\2\8\1\7\2\6\4\5\5\4\eeee
\2\6\4\5\5\4\6\3\7\1\8\3\6\4\5\5\eeee
\7\3\6\4\5\5\4\6\2\8\1\7\2\6\4\5\eeee
\1\8\2\6\4\5\5\4\6\3\7\1\8\3\6\4\eeee
} 

\baaa
8-91
\eaaa
\bbbb
0&0&0&0&0&0&2&2\\
0&0&0&0&0&2&1&1\\
0&0&0&0&0&2&1&1\\
0&0&0&1&1&2&0&0\\
0&0&0&1&1&2&0&0\\
0&1&1&1&1&0&0&0\\
2&1&1&0&0&0&0&0\\
2&1&1&0&0&0&0&0\\
\ebbb
\parbox{7cm}{ 
\1\7\2\6\4\4\6\2\7\1\7\2\6\4\4\6\eeee
\7\1\8\3\6\5\5\6\3\8\1\8\3\6\5\5\eeee
\2\8\1\7\2\6\4\4\6\2\7\1\7\2\6\4\eeee
\6\3\7\1\8\3\6\5\5\6\3\8\1\8\3\6\eeee
\4\6\2\8\1\7\2\6\4\4\6\2\7\1\7\2\eeee
\4\5\6\3\7\1\8\3\6\5\5\6\3\8\1\8\eeee
\6\5\4\6\2\8\1\7\2\6\4\4\6\2\7\1\eeee
\2\6\4\5\6\3\7\1\8\3\6\5\5\6\3\8\eeee
\7\3\6\5\4\6\2\8\1\7\2\6\4\4\6\2\eeee
\1\8\2\6\4\5\6\3\7\1\8\3\6\5\5\6\eeee
\7\1\7\3\6\5\4\6\2\8\1\7\2\6\4\4\eeee
\2\8\1\8\2\6\4\5\6\3\7\1\8\3\6\5\eeee
} 

\baaa
8-92
\eaaa
\bbbb
0&0&0&0&0&0&2&2\\
0&0&0&0&0&2&1&1\\
0&0&0&0&1&2&0&1\\
0&0&0&0&1&2&1&0\\
0&0&2&2&0&0&0&0\\
0&2&1&1&0&0&0&0\\
1&2&0&1&0&0&0&0\\
1&2&1&0&0&0&0&0\\
\ebbb
\parbox{7cm}{ 
\1\8\3\5\3\8\1\8\3\5\3\8\1\8\3\5\eeee
\7\2\6\4\6\2\7\2\6\4\6\2\7\2\6\4\eeee
\4\6\2\7\2\6\4\6\2\7\2\6\4\6\2\7\eeee
\5\3\8\1\8\3\5\3\8\1\8\3\5\3\8\1\eeee
\4\6\2\7\2\6\4\6\2\7\2\6\4\6\2\7\eeee
\7\2\6\4\6\2\7\2\6\4\6\2\7\2\6\4\eeee
\1\8\3\5\3\8\1\8\3\5\3\8\1\8\3\5\eeee
\7\2\6\4\6\2\7\2\6\4\6\2\7\2\6\4\eeee
\4\6\2\7\2\6\4\6\2\7\2\6\4\6\2\7\eeee
\5\3\8\1\8\3\5\3\8\1\8\3\5\3\8\1\eeee
\4\6\2\7\2\6\4\6\2\7\2\6\4\6\2\7\eeee
\7\2\6\4\6\2\7\2\6\4\6\2\7\2\6\4\eeee
} 

\baaa
8-93
\eaaa
\bbbb
0&0&0&0&0&0&2&2\\
0&0&0&0&0&2&1&1\\
0&0&0&1&1&0&1&1\\
0&0&1&0&1&1&0&1\\
0&0&1&1&0&1&1&0\\
0&2&0&1&1&0&0&0\\
1&1&1&0&1&0&0&0\\
1&1&1&1&0&0&0&0\\
\ebbb
\parbox{7cm}{ 
\1\8\4\3\8\2\6\4\5\6\2\7\3\5\7\1\eeee
\7\2\6\5\4\6\2\8\3\4\8\1\8\4\3\8\eeee
\5\6\2\7\3\5\7\1\7\5\3\7\2\6\5\4\eeee
\3\4\8\1\8\4\3\8\2\6\4\5\6\2\7\3\eeee
\7\5\3\7\2\6\5\4\6\2\8\3\4\8\1\8\eeee
\2\6\4\5\6\2\7\3\5\7\1\7\5\3\7\2\eeee
\6\2\8\3\4\8\1\8\4\3\8\2\6\4\5\6\eeee
\5\7\1\7\5\3\7\2\6\5\4\6\2\8\3\4\eeee
\4\3\8\2\6\4\5\6\2\7\3\5\7\1\7\5\eeee
\6\5\4\6\2\8\3\4\8\1\8\4\3\8\2\6\eeee
\2\7\3\5\7\1\7\5\3\7\2\6\5\4\6\2\eeee
\8\1\8\4\3\8\2\6\4\5\6\2\7\3\5\7\eeee
} 

\baaa
8-94
\eaaa
\bbbb
0&0&0&0&0&0&2&2\\
0&0&0&0&0&2&1&1\\
0&0&0&1&1&0&1&1\\
0&0&1&1&0&1&0&1\\
0&0&1&0&1&1&1&0\\
0&2&0&1&1&0&0&0\\
1&1&1&0&1&0&0&0\\
1&1&1&1&0&0&0&0\\
\ebbb
\parbox{7cm}{ 
\1\8\4\3\7\2\6\5\5\6\2\7\3\4\8\1\eeee
\7\2\6\5\5\6\2\7\3\4\8\1\8\4\3\7\eeee
\5\6\2\7\3\4\8\1\8\4\3\7\2\6\5\5\eeee
\3\4\8\1\8\4\3\7\2\6\5\5\6\2\7\3\eeee
\8\4\3\7\2\6\5\5\6\2\7\3\4\8\1\8\eeee
\2\6\5\5\6\2\7\3\4\8\1\8\4\3\7\2\eeee
\6\2\7\3\4\8\1\8\4\3\7\2\6\5\5\6\eeee
\4\8\1\8\4\3\7\2\6\5\5\6\2\7\3\4\eeee
\4\3\7\2\6\5\5\6\2\7\3\4\8\1\8\4\eeee
\6\5\5\6\2\7\3\4\8\1\8\4\3\7\2\6\eeee
\2\7\3\4\8\1\8\4\3\7\2\6\5\5\6\2\eeee
\8\1\8\4\3\7\2\6\5\5\6\2\7\3\4\8\eeee
} 

\baaa
8-95
\eaaa
\bbbb
0&0&0&0&0&0&2&2\\
0&0&0&0&0&2&1&1\\
0&0&0&1&1&1&0&1\\
0&0&1&1&0&0&1&1\\
0&0&1&0&1&1&1&0\\
0&2&1&0&1&0&0&0\\
1&1&0&1&1&0&0&0\\
1&1&1&1&0&0&0&0\\
\ebbb
\parbox{7cm}{ 
\1\8\3\5\6\2\7\4\4\7\2\6\5\3\8\1\eeee
\7\2\6\5\3\8\1\8\3\5\6\2\7\4\4\7\eeee
\5\6\2\7\4\4\7\2\6\5\3\8\1\8\3\5\eeee
\5\3\8\1\8\3\5\6\2\7\4\4\7\2\6\5\eeee
\7\4\4\7\2\6\5\3\8\1\8\3\5\6\2\7\eeee
\1\8\3\5\6\2\7\4\4\7\2\6\5\3\8\1\eeee
\7\2\6\5\3\8\1\8\3\5\6\2\7\4\4\7\eeee
\5\6\2\7\4\4\7\2\6\5\3\8\1\8\3\5\eeee
\5\3\8\1\8\3\5\6\2\7\4\4\7\2\6\5\eeee
\7\4\4\7\2\6\5\3\8\1\8\3\5\6\2\7\eeee
\1\8\3\5\6\2\7\4\4\7\2\6\5\3\8\1\eeee
\7\2\6\5\3\8\1\8\3\5\6\2\7\4\4\7\eeee
} 

\baaa
8-96
\eaaa
\bbbb
0&0&0&0&0&0&2&2\\
0&0&0&0&1&1&0&2\\
0&0&0&0&1&1&0&2\\
0&0&0&0&1&1&2&0\\
0&1&1&2&0&0&0&0\\
0&1&1&2&0&0&0&0\\
2&0&0&2&0&0&0&0\\
2&1&1&0&0&0&0&0\\
\ebbb
\parbox{7cm}{ 
\1\7\4\5\2\8\1\7\4\6\3\8\1\7\4\5\eeee
\7\4\6\3\8\1\7\4\5\2\8\1\7\4\6\3\eeee
\4\5\2\8\1\7\4\6\3\8\1\7\4\5\2\8\eeee
\6\3\8\1\7\4\5\2\8\1\7\4\6\3\8\1\eeee
\2\8\1\7\4\6\3\8\1\7\4\5\2\8\1\7\eeee
\8\1\7\4\5\2\8\1\7\4\6\3\8\1\7\4\eeee
\1\7\4\6\3\8\1\7\4\5\2\8\1\7\4\6\eeee
\7\4\5\2\8\1\7\4\6\3\8\1\7\4\5\2\eeee
\4\6\3\8\1\7\4\5\2\8\1\7\4\6\3\8\eeee
\5\2\8\1\7\4\6\3\8\1\7\4\5\2\8\1\eeee
\3\8\1\7\4\5\2\8\1\7\4\6\3\8\1\7\eeee
\8\1\7\4\6\3\8\1\7\4\5\2\8\1\7\4\eeee
} 

\baaa
8-97
\eaaa
\bbbb
0&0&0&0&0&0&2&2\\
0&0&0&0&1&1&0&2\\
0&0&0&0&1&1&0&2\\
0&0&0&2&0&0&2&0\\
0&1&1&0&1&1&0&0\\
0&1&1&0&1&1&0&0\\
2&0&0&2&0&0&0&0\\
2&1&1&0&0&0&0&0\\
\ebbb
\parbox{7cm}{ 
\1\7\4\4\7\1\8\2\5\6\3\8\1\7\4\4\eeee
\7\4\4\7\1\8\3\6\5\2\8\1\7\4\4\7\eeee
\4\4\7\1\8\2\5\6\3\8\1\7\4\4\7\1\eeee
\4\7\1\8\3\6\5\2\8\1\7\4\4\7\1\8\eeee
\7\1\8\2\5\6\3\8\1\7\4\4\7\1\8\2\eeee
\1\8\3\6\5\2\8\1\7\4\4\7\1\8\3\6\eeee
\8\2\5\6\3\8\1\7\4\4\7\1\8\2\5\6\eeee
\3\6\5\2\8\1\7\4\4\7\1\8\3\6\5\2\eeee
\5\6\3\8\1\7\4\4\7\1\8\2\5\6\3\8\eeee
\5\2\8\1\7\4\4\7\1\8\3\6\5\2\8\1\eeee
\3\8\1\7\4\4\7\1\8\2\5\6\3\8\1\7\eeee
\8\1\7\4\4\7\1\8\3\6\5\2\8\1\7\4\eeee
} 

\baaa
8-98
\eaaa
\bbbb
0&0&0&0&0&0&2&2\\
0&0&0&0&1&1&0&2\\
0&0&0&0&1&1&0&2\\
0&0&0&2&1&1&0&0\\
0&1&1&2&0&0&0&0\\
0&1&1&2&0&0&0&0\\
2&0&0&0&0&0&2&0\\
2&1&1&0&0&0&0&0\\
\ebbb
\parbox{7cm}{ 
\1\7\7\1\8\3\6\4\4\6\3\8\1\7\7\1\eeee
\7\7\1\8\2\5\4\4\5\2\8\1\7\7\1\8\eeee
\7\1\8\3\6\4\4\6\3\8\1\7\7\1\8\3\eeee
\1\8\2\5\4\4\5\2\8\1\7\7\1\8\2\5\eeee
\8\3\6\4\4\6\3\8\1\7\7\1\8\3\6\4\eeee
\2\5\4\4\5\2\8\1\7\7\1\8\2\5\4\4\eeee
\6\4\4\6\3\8\1\7\7\1\8\3\6\4\4\6\eeee
\4\4\5\2\8\1\7\7\1\8\2\5\4\4\5\2\eeee
\4\6\3\8\1\7\7\1\8\3\6\4\4\6\3\8\eeee
\5\2\8\1\7\7\1\8\2\5\4\4\5\2\8\1\eeee
\3\8\1\7\7\1\8\3\6\4\4\6\3\8\1\7\eeee
\8\1\7\7\1\8\2\5\4\4\5\2\8\1\7\7\eeee
} 

\baaa
8-99
\eaaa
\bbbb
0&0&0&0&0&0&2&2\\
0&0&0&0&1&1&0&2\\
0&0&0&0&1&1&1&1\\
0&0&0&0&1&1&2&0\\
0&1&2&1&0&0&0&0\\
0&1&2&1&0&0&0&0\\
2&0&1&1&0&0&0&0\\
2&1&1&0&0&0&0&0\\
\ebbb
\parbox{7cm}{ 
\1\7\3\5\4\7\1\7\4\6\3\7\1\8\3\5\eeee
\7\4\6\3\7\1\8\3\5\2\8\1\8\2\6\3\eeee
\3\5\2\8\1\8\2\6\3\8\1\7\3\5\4\7\eeee
\6\3\8\1\7\3\5\4\7\1\7\4\6\3\7\1\eeee
\4\7\1\7\4\6\3\7\1\8\3\5\2\8\1\8\eeee
\7\1\8\3\5\2\8\1\8\2\6\3\8\1\7\3\eeee
\1\8\2\6\3\8\1\7\3\5\4\7\1\7\4\6\eeee
\7\3\5\4\7\1\7\4\6\3\7\1\8\3\5\2\eeee
\4\6\3\7\1\8\3\5\2\8\1\8\2\6\3\8\eeee
\5\2\8\1\8\2\6\3\8\1\7\3\5\4\7\1\eeee
\3\8\1\7\3\5\4\7\1\7\4\6\3\7\1\8\eeee
\7\1\7\4\6\3\7\1\8\3\5\2\8\1\8\2\eeee
} 

\baaa
8-100
\eaaa
\bbbb
0&0&0&0&0&0&2&2\\
0&0&0&0&1&1&0&2\\
0&0&0&1&0&1&1&1\\
0&0&1&0&1&0&1&1\\
0&1&0&2&0&1&0&0\\
0&1&2&0&1&0&0&0\\
1&0&1&1&0&0&1&0\\
1&1&1&1&0&0&0&0\\
\ebbb
\parbox{7cm}{ 
\1\8\2\8\1\8\2\8\1\8\2\8\1\8\2\8\eeee
\7\3\6\3\7\4\5\4\7\3\6\3\7\4\5\4\eeee
\7\4\5\4\7\3\6\3\7\4\5\4\7\3\6\3\eeee
\1\8\2\8\1\8\2\8\1\8\2\8\1\8\2\8\eeee
\7\3\6\3\7\4\5\4\7\3\6\3\7\4\5\4\eeee
\7\4\5\4\7\3\6\3\7\4\5\4\7\3\6\3\eeee
\1\8\2\8\1\8\2\8\1\8\2\8\1\8\2\8\eeee
\7\3\6\3\7\4\5\4\7\3\6\3\7\4\5\4\eeee
\7\4\5\4\7\3\6\3\7\4\5\4\7\3\6\3\eeee
\1\8\2\8\1\8\2\8\1\8\2\8\1\8\2\8\eeee
\7\3\6\3\7\4\5\4\7\3\6\3\7\4\5\4\eeee
\7\4\5\4\7\3\6\3\7\4\5\4\7\3\6\3\eeee
} 

\baaa
8-101
\eaaa
\bbbb
0&0&0&0&0&0&2&2\\
0&0&0&0&1&1&0&2\\
0&0&1&0&0&1&1&1\\
0&0&0&1&1&0&1&1\\
0&1&0&2&1&0&0&0\\
0&1&2&0&0&1&0&0\\
1&0&1&1&0&0&1&0\\
1&1&1&1&0&0&0&0\\
\ebbb
\parbox{7cm}{ 
\1\8\2\8\1\8\2\8\1\8\2\8\1\8\2\8\eeee
\7\3\6\3\7\4\5\4\7\3\6\3\7\4\5\4\eeee
\7\3\6\3\7\4\5\4\7\3\6\3\7\4\5\4\eeee
\1\8\2\8\1\8\2\8\1\8\2\8\1\8\2\8\eeee
\7\4\5\4\7\3\6\3\7\4\5\4\7\3\6\3\eeee
\7\4\5\4\7\3\6\3\7\4\5\4\7\3\6\3\eeee
\1\8\2\8\1\8\2\8\1\8\2\8\1\8\2\8\eeee
\7\3\6\3\7\4\5\4\7\3\6\3\7\4\5\4\eeee
\7\3\6\3\7\4\5\4\7\3\6\3\7\4\5\4\eeee
\1\8\2\8\1\8\2\8\1\8\2\8\1\8\2\8\eeee
\7\4\5\4\7\3\6\3\7\4\5\4\7\3\6\3\eeee
\7\4\5\4\7\3\6\3\7\4\5\4\7\3\6\3\eeee
} 

\baaa
8-102
\eaaa
\bbbb
0&0&0&0&0&0&2&2\\
0&0&0&0&1&1&0&2\\
0&0&1&1&0&0&2&0\\
0&0&1&1&0&0&2&0\\
0&2&0&0&1&1&0&0\\
0&2&0&0&1&1&0&0\\
2&0&1&1&0&0&0&0\\
2&2&0&0&0&0&0&0\\
\ebbb
\parbox{7cm}{ 
\1\7\4\3\7\1\8\2\5\6\2\8\1\7\4\3\eeee
\7\3\4\7\1\8\2\6\5\2\8\1\7\3\4\7\eeee
\4\3\7\1\8\2\5\6\2\8\1\7\4\3\7\1\eeee
\4\7\1\8\2\6\5\2\8\1\7\3\4\7\1\8\eeee
\7\1\8\2\5\6\2\8\1\7\4\3\7\1\8\2\eeee
\1\8\2\6\5\2\8\1\7\3\4\7\1\8\2\6\eeee
\8\2\5\6\2\8\1\7\4\3\7\1\8\2\5\6\eeee
\2\6\5\2\8\1\7\3\4\7\1\8\2\6\5\2\eeee
\5\6\2\8\1\7\4\3\7\1\8\2\5\6\2\8\eeee
\5\2\8\1\7\3\4\7\1\8\2\6\5\2\8\1\eeee
\2\8\1\7\4\3\7\1\8\2\5\6\2\8\1\7\eeee
\8\1\7\3\4\7\1\8\2\6\5\2\8\1\7\3\eeee
} 

\baaa
8-103
\eaaa
\bbbb
0&0&0&0&0&0&2&2\\
0&0&0&0&1&1&1&1\\
0&0&0&0&1&1&1&1\\
0&0&0&0&2&2&0&0\\
0&1&1&2&0&0&0&0\\
0&1&1&2&0&0&0&0\\
2&1&1&0&0&0&0&0\\
2&1&1&0&0&0&0&0\\
\ebbb
\parbox{7cm}{ 
\1\7\2\5\4\5\2\7\1\7\2\5\4\5\2\7\eeee
\7\1\8\3\6\4\6\3\8\1\8\3\6\4\6\3\eeee
\2\8\1\7\2\5\4\5\2\7\1\7\2\5\4\5\eeee
\5\3\7\1\8\3\6\4\6\3\8\1\8\3\6\4\eeee
\4\6\2\8\1\7\2\5\4\5\2\7\1\7\2\5\eeee
\5\4\5\3\7\1\8\3\6\4\6\3\8\1\8\3\eeee
\2\6\4\6\2\8\1\7\2\5\4\5\2\7\1\7\eeee
\7\3\5\4\5\3\7\1\8\3\6\4\6\3\8\1\eeee
\1\8\2\6\4\6\2\8\1\7\2\5\4\5\2\7\eeee
\7\1\7\3\5\4\5\3\7\1\8\3\6\4\6\3\eeee
\2\8\1\8\2\6\4\6\2\8\1\7\2\5\4\5\eeee
\5\3\7\1\7\3\5\4\5\3\7\1\8\3\6\4\eeee
} 

\baaa
8-104
\eaaa
\bbbb
0&0&0&0&0&0&2&2\\
0&0&0&0&1&1&1&1\\
0&0&0&0&1&1&1&1\\
0&0&0&0&2&2&0&0\\
0&1&2&1&0&0&0&0\\
0&1&2&1&0&0&0&0\\
1&1&2&0&0&0&0&0\\
1&1&2&0&0&0&0&0\\
\ebbb
\parbox{7cm}{ 
\1\8\3\7\1\7\3\8\1\8\3\7\1\7\3\8\eeee
\7\2\6\3\8\2\5\3\7\2\6\3\8\2\5\3\eeee
\3\5\4\5\3\6\4\6\3\5\4\5\3\6\4\6\eeee
\8\3\6\2\7\3\5\2\8\3\6\2\7\3\5\2\eeee
\1\7\3\8\1\8\3\7\1\7\3\8\1\8\3\7\eeee
\8\2\5\3\7\2\6\3\8\2\5\3\7\2\6\3\eeee
\3\6\4\6\3\5\4\5\3\6\4\6\3\5\4\5\eeee
\7\3\5\2\8\3\6\2\7\3\5\2\8\3\6\2\eeee
\1\8\3\7\1\7\3\8\1\8\3\7\1\7\3\8\eeee
\7\2\6\3\8\2\5\3\7\2\6\3\8\2\5\3\eeee
\3\5\4\5\3\6\4\6\3\5\4\5\3\6\4\6\eeee
\8\3\6\2\7\3\5\2\8\3\6\2\7\3\5\2\eeee
} 

\baaa
8-105
\eaaa
\bbbb
0&0&0&0&0&0&2&2\\
0&0&0&0&1&1&1&1\\
0&0&0&0&1&1&1&1\\
0&0&0&2&1&1&0&0\\
0&1&1&2&0&0&0&0\\
0&1&1&2&0&0&0&0\\
2&1&1&0&0&0&0&0\\
2&1&1&0&0&0&0&0\\
\ebbb
\parbox{7cm}{ 
\1\7\2\5\4\4\5\2\7\1\7\2\5\4\4\5\eeee
\7\1\8\3\6\4\4\6\3\8\1\8\3\6\4\4\eeee
\2\8\1\7\2\5\4\4\5\2\7\1\7\2\5\4\eeee
\5\3\7\1\8\3\6\4\4\6\3\8\1\8\3\6\eeee
\4\6\2\8\1\7\2\5\4\4\5\2\7\1\7\2\eeee
\4\4\5\3\7\1\8\3\6\4\4\6\3\8\1\8\eeee
\5\4\4\6\2\8\1\7\2\5\4\4\5\2\7\1\eeee
\2\6\4\4\5\3\7\1\8\3\6\4\4\6\3\8\eeee
\7\3\5\4\4\6\2\8\1\7\2\5\4\4\5\2\eeee
\1\8\2\6\4\4\5\3\7\1\8\3\6\4\4\6\eeee
\7\1\7\3\5\4\4\6\2\8\1\7\2\5\4\4\eeee
\2\8\1\8\2\6\4\4\5\3\7\1\8\3\6\4\eeee
} 

\baaa
8-106
\eaaa
\bbbb
0&0&0&0&0&0&2&2\\
0&0&0&0&1&1&1&1\\
0&0&0&1&0&2&0&1\\
0&0&1&1&1&1&0&0\\
0&1&0&1&0&0&1&1\\
0&1&2&1&0&0&0&0\\
1&1&0&0&1&0&1&0\\
1&1&1&0&1&0&0&0\\
\ebbb
\parbox{7cm}{ 
\1\8\3\6\2\7\5\4\4\5\7\2\6\3\8\1\eeee
\7\2\6\3\8\1\8\3\6\2\7\5\4\4\5\7\eeee
\7\5\4\4\5\7\2\6\3\8\1\8\3\6\2\7\eeee
\1\8\3\6\2\7\5\4\4\5\7\2\6\3\8\1\eeee
\7\2\6\3\8\1\8\3\6\2\7\5\4\4\5\7\eeee
\7\5\4\4\5\7\2\6\3\8\1\8\3\6\2\7\eeee
\1\8\3\6\2\7\5\4\4\5\7\2\6\3\8\1\eeee
\7\2\6\3\8\1\8\3\6\2\7\5\4\4\5\7\eeee
\7\5\4\4\5\7\2\6\3\8\1\8\3\6\2\7\eeee
\1\8\3\6\2\7\5\4\4\5\7\2\6\3\8\1\eeee
\7\2\6\3\8\1\8\3\6\2\7\5\4\4\5\7\eeee
\7\5\4\4\5\7\2\6\3\8\1\8\3\6\2\7\eeee
} 

\baaa
8-107
\eaaa
\bbbb
0&0&0&0&0&0&2&2\\
0&0&0&0&1&1&1&1\\
0&0&0&1&1&1&0&1\\
0&0&1&1&0&0&1&1\\
0&1&1&0&1&0&1&0\\
0&2&2&0&0&0&0&0\\
1&1&0&1&1&0&0&0\\
1&1&1&1&0&0&0&0\\
\ebbb
\parbox{7cm}{ 
\1\8\3\5\5\3\8\1\8\3\5\5\3\8\1\8\eeee
\7\2\6\2\7\4\4\7\2\6\2\7\4\4\7\2\eeee
\5\5\3\8\1\8\3\5\5\3\8\1\8\3\5\5\eeee
\2\7\4\4\7\2\6\2\7\4\4\7\2\6\2\7\eeee
\8\1\8\3\5\5\3\8\1\8\3\5\5\3\8\1\eeee
\4\7\2\6\2\7\4\4\7\2\6\2\7\4\4\7\eeee
\3\5\5\3\8\1\8\3\5\5\3\8\1\8\3\5\eeee
\6\2\7\4\4\7\2\6\2\7\4\4\7\2\6\2\eeee
\3\8\1\8\3\5\5\3\8\1\8\3\5\5\3\8\eeee
\4\4\7\2\6\2\7\4\4\7\2\6\2\7\4\4\eeee
\8\3\5\5\3\8\1\8\3\5\5\3\8\1\8\3\eeee
\2\6\2\7\4\4\7\2\6\2\7\4\4\7\2\6\eeee
} 

\baaa
8-108
\eaaa
\bbbb
0&0&0&0&0&0&2&2\\
0&0&0&0&1&1&1&1\\
0&0&0&1&1&2&0&0\\
0&0&1&1&0&1&0&1\\
0&1&1&0&0&0&1&1\\
0&1&2&1&0&0&0&0\\
1&1&0&0&1&0&1&0\\
1&1&0&1&1&0&0&0\\
\ebbb
\parbox{7cm}{ 
\1\8\4\4\8\1\8\4\4\8\1\8\4\4\8\1\eeee
\7\2\6\3\5\7\2\6\3\5\7\2\6\3\5\7\eeee
\7\5\3\6\2\7\5\3\6\2\7\5\3\6\2\7\eeee
\1\8\4\4\8\1\8\4\4\8\1\8\4\4\8\1\eeee
\7\2\6\3\5\7\2\6\3\5\7\2\6\3\5\7\eeee
\7\5\3\6\2\7\5\3\6\2\7\5\3\6\2\7\eeee
\1\8\4\4\8\1\8\4\4\8\1\8\4\4\8\1\eeee
\7\2\6\3\5\7\2\6\3\5\7\2\6\3\5\7\eeee
\7\5\3\6\2\7\5\3\6\2\7\5\3\6\2\7\eeee
\1\8\4\4\8\1\8\4\4\8\1\8\4\4\8\1\eeee
\7\2\6\3\5\7\2\6\3\5\7\2\6\3\5\7\eeee
\7\5\3\6\2\7\5\3\6\2\7\5\3\6\2\7\eeee
} 

\baaa
8-109
\eaaa
\bbbb
0&0&0&0&0&0&2&2\\
0&0&0&0&1&1&1&1\\
0&0&0&2&0&2&0&0\\
0&0&1&1&1&1&0&0\\
0&1&0&1&0&1&0&1\\
0&1&1&1&1&0&0&0\\
1&1&0&0&0&0&1&1\\
1&1&0&0&1&0&1&0\\
\ebbb
\parbox{7cm}{ 
\1\8\5\4\4\5\8\1\8\5\4\4\5\8\1\8\eeee
\7\2\6\3\6\2\7\7\2\6\3\6\2\7\7\2\eeee
\8\5\4\4\5\8\1\8\5\4\4\5\8\1\8\5\eeee
\2\6\3\6\2\7\7\2\6\3\6\2\7\7\2\6\eeee
\5\4\4\5\8\1\8\5\4\4\5\8\1\8\5\4\eeee
\6\3\6\2\7\7\2\6\3\6\2\7\7\2\6\3\eeee
\4\4\5\8\1\8\5\4\4\5\8\1\8\5\4\4\eeee
\3\6\2\7\7\2\6\3\6\2\7\7\2\6\3\6\eeee
\4\5\8\1\8\5\4\4\5\8\1\8\5\4\4\5\eeee
\6\2\7\7\2\6\3\6\2\7\7\2\6\3\6\2\eeee
\5\8\1\8\5\4\4\5\8\1\8\5\4\4\5\8\eeee
\2\7\7\2\6\3\6\2\7\7\2\6\3\6\2\7\eeee
} 

\baaa
8-110
\eaaa
\bbbb
0&0&0&0&0&0&2&2\\
0&0&0&0&1&1&1&1\\
0&0&0&2&1&1&0&0\\
0&0&2&1&0&1&0&0\\
0&1&1&0&0&1&0&1\\
0&1&1&1&1&0&0&0\\
1&1&0&0&0&0&1&1\\
1&1&0&0&1&0&1&0\\
\ebbb
\parbox{7cm}{ 
\1\8\5\3\4\6\2\7\7\2\6\4\3\5\8\1\eeee
\7\2\6\4\3\5\8\1\8\5\3\4\6\2\7\7\eeee
\8\5\3\4\6\2\7\7\2\6\4\3\5\8\1\8\eeee
\2\6\4\3\5\8\1\8\5\3\4\6\2\7\7\2\eeee
\5\3\4\6\2\7\7\2\6\4\3\5\8\1\8\5\eeee
\6\4\3\5\8\1\8\5\3\4\6\2\7\7\2\6\eeee
\3\4\6\2\7\7\2\6\4\3\5\8\1\8\5\3\eeee
\4\3\5\8\1\8\5\3\4\6\2\7\7\2\6\4\eeee
\4\6\2\7\7\2\6\4\3\5\8\1\8\5\3\4\eeee
\3\5\8\1\8\5\3\4\6\2\7\7\2\6\4\3\eeee
\6\2\7\7\2\6\4\3\5\8\1\8\5\3\4\6\eeee
\5\8\1\8\5\3\4\6\2\7\7\2\6\4\3\5\eeee
} 

\baaa
8-111
\eaaa
\bbbb
0&0&0&0&0&0&2&2\\
0&0&0&0&1&1&1&1\\
0&0&1&1&0&1&0&1\\
0&0&1&1&1&1&0&0\\
0&1&0&1&0&0&1&1\\
0&1&1&1&0&1&0&0\\
1&1&0&0&1&0&1&0\\
1&1&1&0&1&0&0&0\\
\ebbb
\parbox{7cm}{ 
\1\8\3\3\8\1\8\3\3\8\1\8\3\3\8\1\eeee
\7\2\6\6\2\7\5\4\4\5\7\2\6\6\2\7\eeee
\7\5\4\4\5\7\2\6\6\2\7\5\4\4\5\7\eeee
\1\8\3\3\8\1\8\3\3\8\1\8\3\3\8\1\eeee
\7\2\6\6\2\7\5\4\4\5\7\2\6\6\2\7\eeee
\7\5\4\4\5\7\2\6\6\2\7\5\4\4\5\7\eeee
\1\8\3\3\8\1\8\3\3\8\1\8\3\3\8\1\eeee
\7\2\6\6\2\7\5\4\4\5\7\2\6\6\2\7\eeee
\7\5\4\4\5\7\2\6\6\2\7\5\4\4\5\7\eeee
\1\8\3\3\8\1\8\3\3\8\1\8\3\3\8\1\eeee
\7\2\6\6\2\7\5\4\4\5\7\2\6\6\2\7\eeee
\7\5\4\4\5\7\2\6\6\2\7\5\4\4\5\7\eeee
} 

\baaa
8-112
\eaaa
\bbbb
0&0&0&0&0&0&2&2\\
0&0&0&0&1&2&0&1\\
0&0&1&0&0&1&1&1\\
0&0&0&1&1&2&0&0\\
0&2&0&2&0&0&0&0\\
0&1&1&1&0&1&0&0\\
1&0&2&0&0&0&1&0\\
1&1&2&0&0&0&0&0\\
\ebbb
\parbox{7cm}{ 
\1\8\2\5\2\8\1\8\2\5\2\8\1\8\2\5\eeee
\7\3\6\4\6\3\7\3\6\4\6\3\7\3\6\4\eeee
\7\3\6\4\6\3\7\3\6\4\6\3\7\3\6\4\eeee
\1\8\2\5\2\8\1\8\2\5\2\8\1\8\2\5\eeee
\7\3\6\4\6\3\7\3\6\4\6\3\7\3\6\4\eeee
\7\3\6\4\6\3\7\3\6\4\6\3\7\3\6\4\eeee
\1\8\2\5\2\8\1\8\2\5\2\8\1\8\2\5\eeee
\7\3\6\4\6\3\7\3\6\4\6\3\7\3\6\4\eeee
\7\3\6\4\6\3\7\3\6\4\6\3\7\3\6\4\eeee
\1\8\2\5\2\8\1\8\2\5\2\8\1\8\2\5\eeee
\7\3\6\4\6\3\7\3\6\4\6\3\7\3\6\4\eeee
\7\3\6\4\6\3\7\3\6\4\6\3\7\3\6\4\eeee
} 

\baaa
8-113
\eaaa
\bbbb
0&0&0&0&0&0&2&2\\
0&0&0&0&1&2&0&1\\
0&0&1&0&0&1&1&1\\
0&0&0&2&1&1&0&0\\
0&1&0&2&1&0&0&0\\
0&1&1&1&0&1&0&0\\
1&0&2&0&0&0&1&0\\
1&1&2&0&0&0&0&0\\
\ebbb
\parbox{7cm}{ 
\1\8\2\5\5\2\8\1\8\2\5\5\2\8\1\8\eeee
\7\3\6\4\4\6\3\7\3\6\4\4\6\3\7\3\eeee
\7\3\6\4\4\6\3\7\3\6\4\4\6\3\7\3\eeee
\1\8\2\5\5\2\8\1\8\2\5\5\2\8\1\8\eeee
\7\3\6\4\4\6\3\7\3\6\4\4\6\3\7\3\eeee
\7\3\6\4\4\6\3\7\3\6\4\4\6\3\7\3\eeee
\1\8\2\5\5\2\8\1\8\2\5\5\2\8\1\8\eeee
\7\3\6\4\4\6\3\7\3\6\4\4\6\3\7\3\eeee
\7\3\6\4\4\6\3\7\3\6\4\4\6\3\7\3\eeee
\1\8\2\5\5\2\8\1\8\2\5\5\2\8\1\8\eeee
\7\3\6\4\4\6\3\7\3\6\4\4\6\3\7\3\eeee
\7\3\6\4\4\6\3\7\3\6\4\4\6\3\7\3\eeee
} 

\baaa
8-114
\eaaa
\bbbb
0&0&0&0&0&0&2&2\\
0&1&0&0&0&1&1&1\\
0&0&1&0&1&0&1&1\\
0&0&0&1&1&1&0&1\\
0&0&1&1&1&1&0&0\\
0&1&0&1&1&1&0&0\\
1&1&1&0&0&0&1&0\\
1&1&1&1&0&0&0&0\\
\ebbb
\parbox{7cm}{ 
\1\8\4\4\8\1\8\4\4\8\1\8\4\4\8\1\eeee
\7\2\6\5\3\7\2\6\5\3\7\2\6\5\3\7\eeee
\7\2\6\5\3\7\2\6\5\3\7\2\6\5\3\7\eeee
\1\8\4\4\8\1\8\4\4\8\1\8\4\4\8\1\eeee
\7\3\5\6\2\7\3\5\6\2\7\3\5\6\2\7\eeee
\7\3\5\6\2\7\3\5\6\2\7\3\5\6\2\7\eeee
\1\8\4\4\8\1\8\4\4\8\1\8\4\4\8\1\eeee
\7\2\6\5\3\7\2\6\5\3\7\2\6\5\3\7\eeee
\7\2\6\5\3\7\2\6\5\3\7\2\6\5\3\7\eeee
\1\8\4\4\8\1\8\4\4\8\1\8\4\4\8\1\eeee
\7\3\5\6\2\7\3\5\6\2\7\3\5\6\2\7\eeee
\7\3\5\6\2\7\3\5\6\2\7\3\5\6\2\7\eeee
} 

\baaa
8-115
\eaaa
\bbbb
0&0&0&0&0&0&2&2\\
0&1&0&0&0&1&1&1\\
0&0&1&0&1&0&1&1\\
0&0&0&1&1&1&0&1\\
0&0&1&1&2&0&0&0\\
0&1&0&1&0&2&0&0\\
1&1&1&0&0&0&1&0\\
1&1&1&1&0&0&0&0\\
\ebbb
\parbox{7cm}{ 
\1\8\4\4\8\1\8\4\4\8\1\8\4\4\8\1\eeee
\7\2\6\6\2\7\3\5\5\3\7\2\6\6\2\7\eeee
\7\2\6\6\2\7\3\5\5\3\7\2\6\6\2\7\eeee
\1\8\4\4\8\1\8\4\4\8\1\8\4\4\8\1\eeee
\7\3\5\5\3\7\2\6\6\2\7\3\5\5\3\7\eeee
\7\3\5\5\3\7\2\6\6\2\7\3\5\5\3\7\eeee
\1\8\4\4\8\1\8\4\4\8\1\8\4\4\8\1\eeee
\7\2\6\6\2\7\3\5\5\3\7\2\6\6\2\7\eeee
\7\2\6\6\2\7\3\5\5\3\7\2\6\6\2\7\eeee
\1\8\4\4\8\1\8\4\4\8\1\8\4\4\8\1\eeee
\7\3\5\5\3\7\2\6\6\2\7\3\5\5\3\7\eeee
\7\3\5\5\3\7\2\6\6\2\7\3\5\5\3\7\eeee
} 

\baaa
8-116
\eaaa
\bbbb
0&0&0&0&0&1&1&2\\
0&0&0&0&0&1&1&2\\
0&0&0&0&2&1&1&0\\
0&0&0&0&2&1&1&0\\
0&0&1&1&0&0&0&2\\
1&1&1&1&0&0&0&0\\
1&1&1&1&0&0&0&0\\
1&1&0&0&2&0&0&0\\
\ebbb
\parbox{7cm}{ 
\1\7\4\5\8\2\7\3\5\8\2\6\3\5\8\1\eeee
\6\3\5\8\1\6\4\5\8\1\7\4\5\8\2\7\eeee
\4\5\8\2\7\3\5\8\2\6\3\5\8\1\6\4\eeee
\5\8\1\6\4\5\8\1\7\4\5\8\2\7\3\5\eeee
\8\2\7\3\5\8\2\6\3\5\8\1\6\4\5\8\eeee
\1\6\4\5\8\1\7\4\5\8\2\7\3\5\8\2\eeee
\7\3\5\8\2\6\3\5\8\1\6\4\5\8\1\7\eeee
\4\5\8\1\7\4\5\8\2\7\3\5\8\2\6\3\eeee
\5\8\2\6\3\5\8\1\6\4\5\8\1\7\4\5\eeee
\8\1\7\4\5\8\2\7\3\5\8\2\6\3\5\8\eeee
\2\6\3\5\8\1\6\4\5\8\1\7\4\5\8\2\eeee
\7\4\5\8\2\7\3\5\8\2\6\3\5\8\1\6\eeee
} 

\baaa
8-117
\eaaa
\bbbb
0&0&0&0&0&1&1&2\\
0&0&0&0&0&1&1&2\\
0&0&0&0&2&1&1&0\\
0&0&0&0&2&1&1&0\\
0&0&1&1&2&0&0&0\\
1&1&1&1&0&0&0&0\\
1&1&1&1&0&0&0&0\\
1&1&0&0&0&0&0&2\\
\ebbb
\parbox{7cm}{ 
\1\7\4\5\5\4\7\1\8\8\2\6\3\5\5\3\eeee
\6\3\5\5\3\6\2\8\8\1\7\4\5\5\4\7\eeee
\4\5\5\4\7\1\8\8\2\6\3\5\5\3\6\2\eeee
\5\5\3\6\2\8\8\1\7\4\5\5\4\7\1\8\eeee
\5\4\7\1\8\8\2\6\3\5\5\3\6\2\8\8\eeee
\3\6\2\8\8\1\7\4\5\5\4\7\1\8\8\2\eeee
\7\1\8\8\2\6\3\5\5\3\6\2\8\8\1\7\eeee
\2\8\8\1\7\4\5\5\4\7\1\8\8\2\6\3\eeee
\8\8\2\6\3\5\5\3\6\2\8\8\1\7\4\5\eeee
\8\1\7\4\5\5\4\7\1\8\8\2\6\3\5\5\eeee
\2\6\3\5\5\3\6\2\8\8\1\7\4\5\5\4\eeee
\7\4\5\5\4\7\1\8\8\2\6\3\5\5\3\6\eeee
} 

\baaa
8-118
\eaaa
\bbbb
0&0&0&0&0&1&1&2\\
0&0&0&0&0&1&1&2\\
0&0&0&0&2&1&1&0\\
0&0&0&2&0&0&0&2\\
0&0&2&0&2&0&0&0\\
1&1&2&0&0&0&0&0\\
1&1&2&0&0&0&0&0\\
1&1&0&2&0&0&0&0\\
\ebbb
\parbox{7cm}{ 
\1\7\3\5\5\3\7\1\8\4\4\8\1\7\3\5\eeee
\6\3\5\5\3\6\2\8\4\4\8\2\6\3\5\5\eeee
\3\5\5\3\7\1\8\4\4\8\1\7\3\5\5\3\eeee
\5\5\3\6\2\8\4\4\8\2\6\3\5\5\3\6\eeee
\5\3\7\1\8\4\4\8\1\7\3\5\5\3\7\1\eeee
\3\6\2\8\4\4\8\2\6\3\5\5\3\6\2\8\eeee
\7\1\8\4\4\8\1\7\3\5\5\3\7\1\8\4\eeee
\2\8\4\4\8\2\6\3\5\5\3\6\2\8\4\4\eeee
\8\4\4\8\1\7\3\5\5\3\7\1\8\4\4\8\eeee
\4\4\8\2\6\3\5\5\3\6\2\8\4\4\8\2\eeee
\4\8\1\7\3\5\5\3\7\1\8\4\4\8\1\7\eeee
\8\2\6\3\5\5\3\6\2\8\4\4\8\2\6\3\eeee
} 

\baaa
8-119
\eaaa
\bbbb
0&0&0&0&0&1&1&2\\
0&0&0&0&0&1&1&2\\
0&0&0&1&1&0&0&2\\
0&0&2&1&1&0&0&0\\
0&0&2&1&1&0&0&0\\
1&1&0&0&0&1&1&0\\
1&1&0&0&0&1&1&0\\
1&1&2&0&0&0&0&0\\
\ebbb
\parbox{7cm}{ 
\1\7\7\1\8\3\5\4\3\8\2\6\6\2\8\3\eeee
\6\6\2\8\3\4\5\3\8\1\7\7\1\8\3\5\eeee
\7\1\8\3\5\4\3\8\2\6\6\2\8\3\4\5\eeee
\2\8\3\4\5\3\8\1\7\7\1\8\3\5\4\3\eeee
\8\3\5\4\3\8\2\6\6\2\8\3\4\5\3\8\eeee
\3\4\5\3\8\1\7\7\1\8\3\5\4\3\8\2\eeee
\5\4\3\8\2\6\6\2\8\3\4\5\3\8\1\7\eeee
\5\3\8\1\7\7\1\8\3\5\4\3\8\2\6\6\eeee
\3\8\2\6\6\2\8\3\4\5\3\8\1\7\7\1\eeee
\8\1\7\7\1\8\3\5\4\3\8\2\6\6\2\8\eeee
\2\6\6\2\8\3\4\5\3\8\1\7\7\1\8\3\eeee
\7\7\1\8\3\5\4\3\8\2\6\6\2\8\3\4\eeee
} 

\baaa
8-120
\eaaa
\bbbb
0&0&0&0&0&1&1&2\\
0&0&0&0&0&1&1&2\\
0&0&0&1&1&1&1&0\\
0&0&2&1&1&0&0&0\\
0&0&2&1&1&0&0&0\\
1&1&2&0&0&0&0&0\\
1&1&2&0&0&0&0&0\\
1&1&0&0&0&0&0&2\\
\ebbb
\parbox{7cm}{ 
\1\7\3\4\4\3\7\1\8\8\2\6\3\5\5\3\eeee
\6\3\5\5\3\6\2\8\8\1\7\3\4\4\3\7\eeee
\3\4\4\3\7\1\8\8\2\6\3\5\5\3\6\2\eeee
\5\5\3\6\2\8\8\1\7\3\4\4\3\7\1\8\eeee
\4\3\7\1\8\8\2\6\3\5\5\3\6\2\8\8\eeee
\3\6\2\8\8\1\7\3\4\4\3\7\1\8\8\2\eeee
\7\1\8\8\2\6\3\5\5\3\6\2\8\8\1\7\eeee
\2\8\8\1\7\3\4\4\3\7\1\8\8\2\6\3\eeee
\8\8\2\6\3\5\5\3\6\2\8\8\1\7\3\4\eeee
\8\1\7\3\4\4\3\7\1\8\8\2\6\3\5\5\eeee
\2\6\3\5\5\3\6\2\8\8\1\7\3\4\4\3\eeee
\7\3\4\4\3\7\1\8\8\2\6\3\5\5\3\6\eeee
} 

\baaa
8-121
\eaaa
\bbbb
0&0&0&0&0&1&1&2\\
0&0&0&0&0&1&1&2\\
0&0&1&0&1&1&1&0\\
0&0&0&2&0&0&0&2\\
0&0&1&0&1&1&1&0\\
1&1&1&0&1&0&0&0\\
1&1&1&0&1&0&0&0\\
1&1&0&2&0&0&0&0\\
\ebbb
\parbox{7cm}{ 
\1\7\5\3\6\2\8\4\4\8\2\6\3\5\7\1\eeee
\6\3\5\7\1\8\4\4\8\1\7\5\3\6\2\8\eeee
\5\3\6\2\8\4\4\8\2\6\3\5\7\1\8\4\eeee
\5\7\1\8\4\4\8\1\7\5\3\6\2\8\4\4\eeee
\6\2\8\4\4\8\2\6\3\5\7\1\8\4\4\8\eeee
\1\8\4\4\8\1\7\5\3\6\2\8\4\4\8\2\eeee
\8\4\4\8\2\6\3\5\7\1\8\4\4\8\1\7\eeee
\4\4\8\1\7\5\3\6\2\8\4\4\8\2\6\3\eeee
\4\8\2\6\3\5\7\1\8\4\4\8\1\7\5\3\eeee
\8\1\7\5\3\6\2\8\4\4\8\2\6\3\5\7\eeee
\2\6\3\5\7\1\8\4\4\8\1\7\5\3\6\2\eeee
\7\5\3\6\2\8\4\4\8\2\6\3\5\7\1\8\eeee
} 

\baaa
8-122
\eaaa
\bbbb
0&0&0&0&0&1&1&2\\
0&0&0&0&1&0&1&2\\
0&0&0&0&1&1&0&2\\
0&0&0&0&1&1&1&1\\
0&1&1&2&0&0&0&0\\
1&0&1&2&0&0&0&0\\
1&1&0&2&0&0&0&0\\
1&1&1&1&0&0&0&0\\
\ebbb
\parbox{7cm}{ 
\1\7\4\6\4\5\3\8\2\8\1\7\4\6\4\5\eeee
\6\4\5\3\8\2\8\1\7\4\6\4\5\3\8\2\eeee
\3\8\2\8\1\7\4\6\4\5\3\8\2\8\1\7\eeee
\8\1\7\4\6\4\5\3\8\2\8\1\7\4\6\4\eeee
\4\6\4\5\3\8\2\8\1\7\4\6\4\5\3\8\eeee
\5\3\8\2\8\1\7\4\6\4\5\3\8\2\8\1\eeee
\2\8\1\7\4\6\4\5\3\8\2\8\1\7\4\6\eeee
\7\4\6\4\5\3\8\2\8\1\7\4\6\4\5\3\eeee
\4\5\3\8\2\8\1\7\4\6\4\5\3\8\2\8\eeee
\8\2\8\1\7\4\6\4\5\3\8\2\8\1\7\4\eeee
\1\7\4\6\4\5\3\8\2\8\1\7\4\6\4\5\eeee
\6\4\5\3\8\2\8\1\7\4\6\4\5\3\8\2\eeee
} 

\baaa
8-123
\eaaa
\bbbb
0&0&0&0&0&1&1&2\\
0&0&0&0&1&0&1&2\\
0&0&0&0&1&2&0&1\\
0&0&0&0&1&2&1&0\\
0&1&2&1&0&0&0&0\\
1&0&2&1&0&0&0&0\\
2&1&0&1&0&0&0&0\\
2&1&1&0&0&0&0&0\\
\ebbb
\parbox{7cm}{ 
\1\7\1\8\2\8\1\7\1\8\2\8\1\7\1\8\eeee
\6\4\6\3\5\3\6\4\6\3\5\3\6\4\6\3\eeee
\3\5\3\6\4\6\3\5\3\6\4\6\3\5\3\6\eeee
\8\2\8\1\7\1\8\2\8\1\7\1\8\2\8\1\eeee
\1\7\1\8\2\8\1\7\1\8\2\8\1\7\1\8\eeee
\6\4\6\3\5\3\6\4\6\3\5\3\6\4\6\3\eeee
\3\5\3\6\4\6\3\5\3\6\4\6\3\5\3\6\eeee
\8\2\8\1\7\1\8\2\8\1\7\1\8\2\8\1\eeee
\1\7\1\8\2\8\1\7\1\8\2\8\1\7\1\8\eeee
\6\4\6\3\5\3\6\4\6\3\5\3\6\4\6\3\eeee
\3\5\3\6\4\6\3\5\3\6\4\6\3\5\3\6\eeee
\8\2\8\1\7\1\8\2\8\1\7\1\8\2\8\1\eeee
} 
\baaa
\phantom{0-00}\#
\eaaa
\mbox{}\phantom{\bbbb
0&0&0&0&0&0&0&0\\
\ebbb}
\parbox{7cm}{ 
\1\7\2\8\1\8\1\7\2\8\1\8\1\7\2\8\eeee
\6\4\5\3\6\3\6\4\5\3\6\3\6\4\5\3\eeee
\3\6\3\6\4\5\3\6\3\6\4\5\3\6\3\6\eeee
\8\1\8\1\7\2\8\1\8\1\7\2\8\1\8\1\eeee
\1\7\2\8\1\8\1\7\2\8\1\8\1\7\2\8\eeee
\6\4\5\3\6\3\6\4\5\3\6\3\6\4\5\3\eeee
\3\6\3\6\4\5\3\6\3\6\4\5\3\6\3\6\eeee
\8\1\8\1\7\2\8\1\8\1\7\2\8\1\8\1\eeee
\1\7\2\8\1\8\1\7\2\8\1\8\1\7\2\8\eeee
\6\4\5\3\6\3\6\4\5\3\6\3\6\4\5\3\eeee
\3\6\3\6\4\5\3\6\3\6\4\5\3\6\3\6\eeee
\8\1\8\1\7\2\8\1\8\1\7\2\8\1\8\1\eeee
} 

\baaa
8-124
\eaaa
\bbbb
0&0&0&0&0&1&1&2\\
0&0&0&0&1&0&1&2\\
0&0&1&0&0&2&1&0\\
0&0&0&1&1&1&0&1\\
0&1&0&2&1&0&0&0\\
1&0&1&1&0&1&0&0\\
2&1&1&0&0&0&0&0\\
2&1&0&1&0&0&0&0\\
\ebbb
\parbox{7cm}{ 
\1\7\1\8\2\8\1\7\1\8\2\8\1\7\1\8\eeee
\6\3\6\4\5\4\6\3\6\4\5\4\6\3\6\4\eeee
\6\3\6\4\5\4\6\3\6\4\5\4\6\3\6\4\eeee
\1\7\1\8\2\8\1\7\1\8\2\8\1\7\1\8\eeee
\8\2\8\1\7\1\8\2\8\1\7\1\8\2\8\1\eeee
\4\5\4\6\3\6\4\5\4\6\3\6\4\5\4\6\eeee
\4\5\4\6\3\6\4\5\4\6\3\6\4\5\4\6\eeee
\8\2\8\1\7\1\8\2\8\1\7\1\8\2\8\1\eeee
\1\7\1\8\2\8\1\7\1\8\2\8\1\7\1\8\eeee
\6\3\6\4\5\4\6\3\6\4\5\4\6\3\6\4\eeee
\6\3\6\4\5\4\6\3\6\4\5\4\6\3\6\4\eeee
\1\7\1\8\2\8\1\7\1\8\2\8\1\7\1\8\eeee
} 

\baaa
8-125
\eaaa
\bbbb
0&0&0&0&0&1&1&2\\
0&0&0&0&1&0&2&1\\
0&0&0&0&1&2&0&1\\
0&0&0&0&2&1&1&0\\
0&1&1&2&0&0&0&0\\
1&0&2&1&0&0&0&0\\
1&2&0&1&0&0&0&0\\
2&1&1&0&0&0&0&0\\
\ebbb
\parbox{7cm}{ 
\1\7\2\8\1\7\2\8\1\7\2\8\1\7\2\8\eeee
\6\4\5\3\6\4\5\3\6\4\5\3\6\4\5\3\eeee
\3\5\4\6\3\5\4\6\3\5\4\6\3\5\4\6\eeee
\8\2\7\1\8\2\7\1\8\2\7\1\8\2\7\1\eeee
\1\7\2\8\1\7\2\8\1\7\2\8\1\7\2\8\eeee
\6\4\5\3\6\4\5\3\6\4\5\3\6\4\5\3\eeee
\3\5\4\6\3\5\4\6\3\5\4\6\3\5\4\6\eeee
\8\2\7\1\8\2\7\1\8\2\7\1\8\2\7\1\eeee
\1\7\2\8\1\7\2\8\1\7\2\8\1\7\2\8\eeee
\6\4\5\3\6\4\5\3\6\4\5\3\6\4\5\3\eeee
\3\5\4\6\3\5\4\6\3\5\4\6\3\5\4\6\eeee
\8\2\7\1\8\2\7\1\8\2\7\1\8\2\7\1\eeee
} 

\baaa
\phantom{0-00}\#
\eaaa
\mbox{}\phantom{\bbbb
0&0&0&0&0&0&0&0\\
\ebbb}
\parbox{7cm}{ 
\1\8\1\8\1\8\1\8\1\8\1\8\1\8\1\8\eeee
\6\3\6\3\6\3\6\3\6\3\6\3\6\3\6\3\eeee
\4\5\4\5\4\5\4\5\4\5\4\5\4\5\4\5\eeee
\7\2\7\2\7\2\7\2\7\2\7\2\7\2\7\2\eeee
\1\8\1\8\1\8\1\8\1\8\1\8\1\8\1\8\eeee
\6\3\6\3\6\3\6\3\6\3\6\3\6\3\6\3\eeee
\4\5\4\5\4\5\4\5\4\5\4\5\4\5\4\5\eeee
\7\2\7\2\7\2\7\2\7\2\7\2\7\2\7\2\eeee
\1\8\1\8\1\8\1\8\1\8\1\8\1\8\1\8\eeee
\6\3\6\3\6\3\6\3\6\3\6\3\6\3\6\3\eeee
\4\5\4\5\4\5\4\5\4\5\4\5\4\5\4\5\eeee
\7\2\7\2\7\2\7\2\7\2\7\2\7\2\7\2\eeee
} 

\baaa
8-126
\eaaa
\bbbb
0&0&0&0&0&1&1&2\\
0&0&0&0&1&0&2&1\\
0&0&0&0&1&2&1&0\\
0&0&0&0&2&1&0&1\\
0&1&1&2&0&0&0&0\\
1&0&2&1&0&0&0&0\\
1&2&1&0&0&0&0&0\\
2&1&0&1&0&0&0&0\\
\ebbb
\parbox{7cm}{ 
\1\7\2\8\1\7\2\8\1\7\2\8\1\7\2\8\eeee
\6\3\5\4\6\3\5\4\6\3\5\4\6\3\5\4\eeee
\3\6\4\5\3\6\4\5\3\6\4\5\3\6\4\5\eeee
\7\1\8\2\7\1\8\2\7\1\8\2\7\1\8\2\eeee
\2\8\1\7\2\8\1\7\2\8\1\7\2\8\1\7\eeee
\5\4\6\3\5\4\6\3\5\4\6\3\5\4\6\3\eeee
\4\5\3\6\4\5\3\6\4\5\3\6\4\5\3\6\eeee
\8\2\7\1\8\2\7\1\8\2\7\1\8\2\7\1\eeee
\1\7\2\8\1\7\2\8\1\7\2\8\1\7\2\8\eeee
\6\3\5\4\6\3\5\4\6\3\5\4\6\3\5\4\eeee
\3\6\4\5\3\6\4\5\3\6\4\5\3\6\4\5\eeee
\7\1\8\2\7\1\8\2\7\1\8\2\7\1\8\2\eeee
} 

\baaa
8-127
\eaaa
\bbbb
0&0&0&0&0&1&1&2\\
0&0&0&0&1&0&2&1\\
0&0&0&1&0&2&0&1\\
0&0&1&0&2&0&1&0\\
0&1&0&2&0&1&0&0\\
1&0&2&0&1&0&0&0\\
1&2&0&1&0&0&0&0\\
2&1&1&0&0&0&0&0\\
\ebbb
\parbox{7cm}{ 
\1\8\1\8\1\8\1\8\1\8\1\8\1\8\1\8\eeee
\6\3\6\3\6\3\6\3\6\3\6\3\6\3\6\3\eeee
\5\4\5\4\5\4\5\4\5\4\5\4\5\4\5\4\eeee
\2\7\2\7\2\7\2\7\2\7\2\7\2\7\2\7\eeee
\8\1\8\1\8\1\8\1\8\1\8\1\8\1\8\1\eeee
\3\6\3\6\3\6\3\6\3\6\3\6\3\6\3\6\eeee
\4\5\4\5\4\5\4\5\4\5\4\5\4\5\4\5\eeee
\7\2\7\2\7\2\7\2\7\2\7\2\7\2\7\2\eeee
\1\8\1\8\1\8\1\8\1\8\1\8\1\8\1\8\eeee
\6\3\6\3\6\3\6\3\6\3\6\3\6\3\6\3\eeee
\5\4\5\4\5\4\5\4\5\4\5\4\5\4\5\4\eeee
\2\7\2\7\2\7\2\7\2\7\2\7\2\7\2\7\eeee
} 

\baaa
8-128
\eaaa
\bbbb
0&0&0&0&0&1&1&2\\
0&0&0&0&1&0&2&1\\
0&0&0&1&0&2&0&1\\
0&0&1&2&0&1&0&0\\
0&1&0&0&2&0&1&0\\
1&0&2&1&0&0&0&0\\
1&2&0&0&1&0&0&0\\
2&1&1&0&0&0&0&0\\
\ebbb
\parbox{7cm}{ 
\1\8\1\8\1\8\1\8\1\8\1\8\1\8\1\8\eeee
\6\3\6\3\6\3\6\3\6\3\6\3\6\3\6\3\eeee
\4\4\4\4\4\4\4\4\4\4\4\4\4\4\4\4\eeee
\3\6\3\6\3\6\3\6\3\6\3\6\3\6\3\6\eeee
\8\1\8\1\8\1\8\1\8\1\8\1\8\1\8\1\eeee
\2\7\2\7\2\7\2\7\2\7\2\7\2\7\2\7\eeee
\5\5\5\5\5\5\5\5\5\5\5\5\5\5\5\5\eeee
\7\2\7\2\7\2\7\2\7\2\7\2\7\2\7\2\eeee
\1\8\1\8\1\8\1\8\1\8\1\8\1\8\1\8\eeee
\6\3\6\3\6\3\6\3\6\3\6\3\6\3\6\3\eeee
\4\4\4\4\4\4\4\4\4\4\4\4\4\4\4\4\eeee
\3\6\3\6\3\6\3\6\3\6\3\6\3\6\3\6\eeee
} 

\baaa
8-129
\eaaa
\bbbb
0&0&0&0&0&1&1&2\\
0&0&0&0&1&0&2&1\\
0&0&0&1&1&1&0&1\\
0&0&1&0&1&1&1&0\\
0&1&1&1&0&1&0&0\\
1&0&1&1&1&0&0&0\\
1&2&0&1&0&0&0&0\\
2&1&1&0&0&0&0&0\\
\ebbb
\parbox{7cm}{ 
\1\7\2\8\1\7\2\8\1\7\2\8\1\7\2\8\eeee
\6\4\5\3\6\4\5\3\6\4\5\3\6\4\5\3\eeee
\5\3\6\4\5\3\6\4\5\3\6\4\5\3\6\4\eeee
\2\8\1\7\2\8\1\7\2\8\1\7\2\8\1\7\eeee
\7\1\8\2\7\1\8\2\7\1\8\2\7\1\8\2\eeee
\4\6\3\5\4\6\3\5\4\6\3\5\4\6\3\5\eeee
\3\5\4\6\3\5\4\6\3\5\4\6\3\5\4\6\eeee
\8\2\7\1\8\2\7\1\8\2\7\1\8\2\7\1\eeee
\1\7\2\8\1\7\2\8\1\7\2\8\1\7\2\8\eeee
\6\4\5\3\6\4\5\3\6\4\5\3\6\4\5\3\eeee
\5\3\6\4\5\3\6\4\5\3\6\4\5\3\6\4\eeee
\2\8\1\7\2\8\1\7\2\8\1\7\2\8\1\7\eeee
} 

\baaa
8-130
\eaaa
\bbbb
0&0&0&0&0&1&1&2\\
0&0&0&0&1&0&2&1\\
0&0&0&1&1&1&0&1\\
0&0&1&0&1&1&1&0\\
0&1&1&1&1&0&0&0\\
1&0&1&1&0&1&0&0\\
1&2&0&1&0&0&0&0\\
2&1&1&0&0&0&0&0\\
\ebbb
\parbox{7cm}{ 
\1\7\2\8\1\7\2\8\1\7\2\8\1\7\2\8\eeee
\6\4\5\3\6\4\5\3\6\4\5\3\6\4\5\3\eeee
\6\3\5\4\6\3\5\4\6\3\5\4\6\3\5\4\eeee
\1\8\2\7\1\8\2\7\1\8\2\7\1\8\2\7\eeee
\8\1\7\2\8\1\7\2\8\1\7\2\8\1\7\2\eeee
\3\6\4\5\3\6\4\5\3\6\4\5\3\6\4\5\eeee
\4\6\3\5\4\6\3\5\4\6\3\5\4\6\3\5\eeee
\7\1\8\2\7\1\8\2\7\1\8\2\7\1\8\2\eeee
\2\8\1\7\2\8\1\7\2\8\1\7\2\8\1\7\eeee
\5\3\6\4\5\3\6\4\5\3\6\4\5\3\6\4\eeee
\5\4\6\3\5\4\6\3\5\4\6\3\5\4\6\3\eeee
\2\7\1\8\2\7\1\8\2\7\1\8\2\7\1\8\eeee
} 

\baaa
8-131
\eaaa
\bbbb
0&0&0&0&0&1&1&2\\
0&0&0&0&1&0&2&1\\
0&0&1&0&0&2&0&1\\
0&0&0&1&2&0&1&0\\
0&1&0&2&1&0&0&0\\
1&0&2&0&0&1&0&0\\
1&2&0&1&0&0&0&0\\
2&1&1&0&0&0&0&0\\
\ebbb
\parbox{7cm}{ 
\1\8\1\8\1\8\1\8\1\8\1\8\1\8\1\8\eeee
\6\3\6\3\6\3\6\3\6\3\6\3\6\3\6\3\eeee
\6\3\6\3\6\3\6\3\6\3\6\3\6\3\6\3\eeee
\1\8\1\8\1\8\1\8\1\8\1\8\1\8\1\8\eeee
\7\2\7\2\7\2\7\2\7\2\7\2\7\2\7\2\eeee
\4\5\4\5\4\5\4\5\4\5\4\5\4\5\4\5\eeee
\4\5\4\5\4\5\4\5\4\5\4\5\4\5\4\5\eeee
\7\2\7\2\7\2\7\2\7\2\7\2\7\2\7\2\eeee
\1\8\1\8\1\8\1\8\1\8\1\8\1\8\1\8\eeee
\6\3\6\3\6\3\6\3\6\3\6\3\6\3\6\3\eeee
\6\3\6\3\6\3\6\3\6\3\6\3\6\3\6\3\eeee
\1\8\1\8\1\8\1\8\1\8\1\8\1\8\1\8\eeee
} 

\baaa
8-132
\eaaa
\bbbb
0&0&0&0&0&1&1&2\\
0&0&0&0&1&0&2&1\\
0&0&1&0&1&1&0&1\\
0&0&0&1&1&1&1&0\\
0&1&1&1&1&0&0&0\\
1&0&1&1&0&1&0&0\\
1&2&0&1&0&0&0&0\\
2&1&1&0&0&0&0&0\\
\ebbb
\parbox{7cm}{ 
\1\7\2\8\1\7\2\8\1\7\2\8\1\7\2\8\eeee
\6\4\5\3\6\4\5\3\6\4\5\3\6\4\5\3\eeee
\6\4\5\3\6\4\5\3\6\4\5\3\6\4\5\3\eeee
\1\7\2\8\1\7\2\8\1\7\2\8\1\7\2\8\eeee
\8\2\7\1\8\2\7\1\8\2\7\1\8\2\7\1\eeee
\3\5\4\6\3\5\4\6\3\5\4\6\3\5\4\6\eeee
\3\5\4\6\3\5\4\6\3\5\4\6\3\5\4\6\eeee
\8\2\7\1\8\2\7\1\8\2\7\1\8\2\7\1\eeee
\1\7\2\8\1\7\2\8\1\7\2\8\1\7\2\8\eeee
\6\4\5\3\6\4\5\3\6\4\5\3\6\4\5\3\eeee
\6\4\5\3\6\4\5\3\6\4\5\3\6\4\5\3\eeee
\1\7\2\8\1\7\2\8\1\7\2\8\1\7\2\8\eeee
} 

\baaa
8-133
\eaaa
\bbbb
0&0&0&0&0&1&1&2\\
0&0&0&0&1&1&1&1\\
0&0&0&0&1&1&1&1\\
0&0&0&0&2&1&1&0\\
0&1&1&2&0&0&0&0\\
1&1&1&1&0&0&0&0\\
1&1&1&1&0&0&0&0\\
2&1&1&0&0&0&0&0\\
\ebbb
\parbox{7cm}{ 
\1\7\2\8\1\7\2\8\1\7\2\8\1\7\2\6\eeee
\6\4\5\3\6\4\5\3\6\4\5\3\6\4\5\3\eeee
\2\5\4\7\2\5\4\7\2\5\4\7\2\5\4\7\eeee
\8\3\6\1\8\3\6\1\8\3\6\1\8\3\6\1\eeee
\1\7\2\8\1\7\2\8\1\7\2\8\1\7\2\8\eeee
\6\4\5\3\6\4\5\3\6\4\5\3\6\4\5\3\eeee
\2\5\4\7\2\5\4\7\2\5\4\7\2\5\4\7\eeee
\8\3\6\1\8\3\6\1\8\3\6\1\8\3\6\1\eeee
\1\7\2\8\1\7\2\8\1\7\2\8\1\7\2\8\eeee
\6\4\5\3\6\4\5\3\6\4\5\3\6\4\5\3\eeee
\2\5\4\7\2\5\4\7\2\5\4\7\2\5\4\7\eeee
\8\3\6\1\8\3\6\1\8\3\6\1\8\3\6\1\eeee
} 

\baaa
8-134
\eaaa
\bbbb
0&0&0&0&0&1&1&2\\
0&0&0&0&1&1&1&1\\
0&0&0&1&0&1&1&1\\
0&0&1&1&1&0&1&0\\
0&1&0&1&1&1&0&0\\
1&1&1&0&1&0&0&0\\
1&1&1&1&0&0&0&0\\
2&1&1&0&0&0&0&0\\
\ebbb
\parbox{7cm}{ 
\1\8\1\8\1\8\1\8\1\8\1\8\1\8\1\6\eeee
\6\2\7\3\6\2\7\3\6\2\7\3\6\2\7\3\eeee
\5\5\4\4\5\5\4\4\5\5\4\4\5\5\4\4\eeee
\2\6\3\7\2\6\3\7\2\6\3\7\2\6\3\7\eeee
\8\1\8\1\8\1\8\1\8\1\8\1\8\1\8\1\eeee
\3\7\2\6\3\7\2\6\3\7\2\6\3\7\2\6\eeee
\4\4\5\5\4\4\5\5\4\4\5\5\4\4\5\5\eeee
\7\3\6\2\7\3\6\2\7\3\6\2\7\3\6\2\eeee
\1\8\1\8\1\8\1\8\1\8\1\8\1\8\1\8\eeee
\6\2\7\3\6\2\7\3\6\2\7\3\6\2\7\3\eeee
\5\5\4\4\5\5\4\4\5\5\4\4\5\5\4\4\eeee
\2\6\3\7\2\6\3\7\2\6\3\7\2\6\3\7\eeee
} 

\baaa
8-135
\eaaa
\bbbb
0&0&0&0&0&1&1&2\\
0&0&0&0&1&1&2&0\\
0&0&0&1&1&0&1&1\\
0&0&1&1&0&1&0&1\\
0&1&1&0&1&1&0&0\\
1&1&0&1&1&0&0&0\\
1&2&1&0&0&0&0&0\\
2&0&1&1&0&0&0&0\\
\ebbb
\parbox{7cm}{ 
\1\7\2\6\4\4\6\2\7\1\8\3\5\5\3\8\eeee
\6\2\7\1\8\3\5\5\3\8\1\7\2\6\4\4\eeee
\5\5\3\8\1\7\2\6\4\4\6\2\7\1\8\3\eeee
\2\6\4\4\6\2\7\1\8\3\5\5\3\8\1\7\eeee
\7\1\8\3\5\5\3\8\1\7\2\6\4\4\6\2\eeee
\3\8\1\7\2\6\4\4\6\2\7\1\8\3\5\5\eeee
\4\4\6\2\7\1\8\3\5\5\3\8\1\7\2\6\eeee
\8\3\5\5\3\8\1\7\2\6\4\4\6\2\7\1\eeee
\1\7\2\6\4\4\6\2\7\1\8\3\5\5\3\8\eeee
\6\2\7\1\8\3\5\5\3\8\1\7\2\6\4\4\eeee
\5\5\3\8\1\7\2\6\4\4\6\2\7\1\8\3\eeee
\2\6\4\4\6\2\7\1\8\3\5\5\3\8\1\7\eeee
} 

\baaa
8-136
\eaaa
\bbbb
0&0&0&0&0&1&1&2\\
0&0&0&0&2&1&1&0\\
0&0&0&1&1&0&1&1\\
0&0&1&0&1&1&0&1\\
0&1&1&1&0&0&0&1\\
1&1&0&2&0&0&0&0\\
1&1&2&0&0&0&0&0\\
1&0&1&1&1&0&0&0\\
\ebbb
\parbox{7cm}{ 
\1\8\5\2\5\8\1\8\5\2\5\8\1\8\5\2\eeee
\6\4\3\7\3\4\6\4\3\7\3\4\6\4\3\7\eeee
\2\5\8\1\8\5\2\5\8\1\8\5\2\5\8\1\eeee
\7\3\4\6\4\3\7\3\4\6\4\3\7\3\4\6\eeee
\1\8\5\2\5\8\1\8\5\2\5\8\1\8\5\2\eeee
\6\4\3\7\3\4\6\4\3\7\3\4\6\4\3\7\eeee
\2\5\8\1\8\5\2\5\8\1\8\5\2\5\8\1\eeee
\7\3\4\6\4\3\7\3\4\6\4\3\7\3\4\6\eeee
\1\8\5\2\5\8\1\8\5\2\5\8\1\8\5\2\eeee
\6\4\3\7\3\4\6\4\3\7\3\4\6\4\3\7\eeee
\2\5\8\1\8\5\2\5\8\1\8\5\2\5\8\1\eeee
\7\3\4\6\4\3\7\3\4\6\4\3\7\3\4\6\eeee
} 

\baaa
8-137
\eaaa
\bbbb
0&0&0&0&0&1&1&2\\
0&0&0&0&2&1&1&0\\
0&0&0&1&1&0&1&1\\
0&0&1&0&1&1&0&1\\
0&1&1&1&1&0&0&0\\
1&1&0&2&0&0&0&0\\
1&1&2&0&0&0&0&0\\
1&0&1&1&0&0&0&1\\
\ebbb
\parbox{7cm}{ 
\1\8\8\1\8\8\1\8\8\1\8\8\1\8\8\1\eeee
\6\4\3\7\3\4\6\4\3\7\3\4\6\4\3\7\eeee
\2\5\5\2\5\5\2\5\5\2\5\5\2\5\5\2\eeee
\7\3\4\6\4\3\7\3\4\6\4\3\7\3\4\6\eeee
\1\8\8\1\8\8\1\8\8\1\8\8\1\8\8\1\eeee
\6\4\3\7\3\4\6\4\3\7\3\4\6\4\3\7\eeee
\2\5\5\2\5\5\2\5\5\2\5\5\2\5\5\2\eeee
\7\3\4\6\4\3\7\3\4\6\4\3\7\3\4\6\eeee
\1\8\8\1\8\8\1\8\8\1\8\8\1\8\8\1\eeee
\6\4\3\7\3\4\6\4\3\7\3\4\6\4\3\7\eeee
\2\5\5\2\5\5\2\5\5\2\5\5\2\5\5\2\eeee
\7\3\4\6\4\3\7\3\4\6\4\3\7\3\4\6\eeee
} 

\baaa
8-138
\eaaa
\bbbb
0&0&0&0&0&1&1&2\\
0&0&0&0&2&1&1&0\\
0&0&1&0&1&0&1&1\\
0&0&0&1&1&1&0&1\\
0&1&1&1&1&0&0&0\\
1&1&0&2&0&0&0&0\\
1&1&2&0&0&0&0&0\\
1&0&1&1&0&0&0&1\\
\ebbb
\parbox{7cm}{ 
\1\8\8\1\8\8\1\8\8\1\8\8\1\8\8\1\eeee
\6\4\4\6\4\4\6\4\4\6\4\4\6\4\4\6\eeee
\2\5\5\2\5\5\2\5\5\2\5\5\2\5\5\2\eeee
\7\3\3\7\3\3\7\3\3\7\3\3\7\3\3\7\eeee
\1\8\8\1\8\8\1\8\8\1\8\8\1\8\8\1\eeee
\6\4\4\6\4\4\6\4\4\6\4\4\6\4\4\6\eeee
\2\5\5\2\5\5\2\5\5\2\5\5\2\5\5\2\eeee
\7\3\3\7\3\3\7\3\3\7\3\3\7\3\3\7\eeee
\1\8\8\1\8\8\1\8\8\1\8\8\1\8\8\1\eeee
\6\4\4\6\4\4\6\4\4\6\4\4\6\4\4\6\eeee
\2\5\5\2\5\5\2\5\5\2\5\5\2\5\5\2\eeee
\7\3\3\7\3\3\7\3\3\7\3\3\7\3\3\7\eeee
} 

\baaa
8-139
\eaaa
\bbbb
0&0&0&0&0&1&1&2\\
0&0&0&1&1&0&1&1\\
0&0&2&0&0&1&1&0\\
0&1&0&0&1&1&0&1\\
0&1&0&1&2&0&0&0\\
1&0&1&1&0&0&1&0\\
1&1&1&0&0&1&0&0\\
2&1&0&1&0&0&0&0\\
\ebbb
\parbox{7cm}{ 
\1\7\6\1\8\4\2\8\1\7\6\1\8\4\2\8\eeee
\6\3\3\7\2\5\5\4\6\3\3\7\2\5\5\4\eeee
\7\3\3\6\4\5\5\2\7\3\3\6\4\5\5\2\eeee
\1\6\7\1\8\2\4\8\1\6\7\1\8\2\4\8\eeee
\8\4\2\8\1\7\6\1\8\4\2\8\1\7\6\1\eeee
\2\5\5\4\6\3\3\7\2\5\5\4\6\3\3\7\eeee
\4\5\5\2\7\3\3\6\4\5\5\2\7\3\3\6\eeee
\8\2\4\8\1\6\7\1\8\2\4\8\1\6\7\1\eeee
\1\7\6\1\8\4\2\8\1\7\6\1\8\4\2\8\eeee
\6\3\3\7\2\5\5\4\6\3\3\7\2\5\5\4\eeee
\7\3\3\6\4\5\5\2\7\3\3\6\4\5\5\2\eeee
\1\6\7\1\8\2\4\8\1\6\7\1\8\2\4\8\eeee
} 

\baaa
8-140
\eaaa
\bbbb
0&0&0&0&0&1&1&2\\
0&0&0&1&1&1&1&0\\
0&0&1&1&2&0&0&0\\
0&1&1&0&1&0&1&0\\
0&1&2&1&0&0&0&0\\
1&1&0&0&0&0&1&1\\
1&1&0&1&0&1&0&0\\
2&0&0&0&0&1&0&1\\
\ebbb
\parbox{7cm}{ 
\1\7\4\3\5\2\6\8\1\7\4\3\5\2\6\8\eeee
\6\2\5\3\4\7\1\8\6\2\5\3\4\7\1\8\eeee
\7\4\3\5\2\6\8\1\7\4\3\5\2\6\8\1\eeee
\2\5\3\4\7\1\8\6\2\5\3\4\7\1\8\6\eeee
\4\3\5\2\6\8\1\7\4\3\5\2\6\8\1\7\eeee
\5\3\4\7\1\8\6\2\5\3\4\7\1\8\6\2\eeee
\3\5\2\6\8\1\7\4\3\5\2\6\8\1\7\4\eeee
\3\4\7\1\8\6\2\5\3\4\7\1\8\6\2\5\eeee
\5\2\6\8\1\7\4\3\5\2\6\8\1\7\4\3\eeee
\4\7\1\8\6\2\5\3\4\7\1\8\6\2\5\3\eeee
\2\6\8\1\7\4\3\5\2\6\8\1\7\4\3\5\eeee
\7\1\8\6\2\5\3\4\7\1\8\6\2\5\3\4\eeee
} 

\baaa
8-141
\eaaa
\bbbb
0&0&0&0&0&1&1&2\\
0&1&0&0&1&0&1&1\\
0&0&1&0&1&1&0&1\\
0&0&0&2&0&1&1&0\\
0&1&1&0&2&0&0&0\\
1&0&1&1&0&1&0&0\\
1&1&0&1&0&0&1&0\\
2&1&1&0&0&0&0&0\\
\ebbb
\parbox{7cm}{ 
\1\7\7\1\8\2\2\8\1\7\7\1\8\2\2\8\eeee
\6\4\4\6\3\5\5\3\6\4\4\6\3\5\5\3\eeee
\6\4\4\6\3\5\5\3\6\4\4\6\3\5\5\3\eeee
\1\7\7\1\8\2\2\8\1\7\7\1\8\2\2\8\eeee
\8\2\2\8\1\7\7\1\8\2\2\8\1\7\7\1\eeee
\3\5\5\3\6\4\4\6\3\5\5\3\6\4\4\6\eeee
\3\5\5\3\6\4\4\6\3\5\5\3\6\4\4\6\eeee
\8\2\2\8\1\7\7\1\8\2\2\8\1\7\7\1\eeee
\1\7\7\1\8\2\2\8\1\7\7\1\8\2\2\8\eeee
\6\4\4\6\3\5\5\3\6\4\4\6\3\5\5\3\eeee
\6\4\4\6\3\5\5\3\6\4\4\6\3\5\5\3\eeee
\1\7\7\1\8\2\2\8\1\7\7\1\8\2\2\8\eeee
} 

\baaa
8-142
\eaaa
\bbbb
0&0&0&0&0&1&1&2\\
0&1&0&0&1&0&1&1\\
0&0&1&0&2&0&1&0\\
0&0&0&2&0&1&0&1\\
0&1&1&0&2&0&0&0\\
1&0&0&2&0&1&0&0\\
1&2&1&0&0&0&0&0\\
1&1&0&1&0&0&0&1\\
\ebbb
\parbox{7cm}{ 
\1\8\8\1\8\8\1\8\8\1\8\8\1\8\8\1\eeee
\6\4\4\6\4\4\6\4\4\6\4\4\6\4\4\6\eeee
\6\4\4\6\4\4\6\4\4\6\4\4\6\4\4\6\eeee
\1\8\8\1\8\8\1\8\8\1\8\8\1\8\8\1\eeee
\7\2\2\7\2\2\7\2\2\7\2\2\7\2\2\7\eeee
\3\5\5\3\5\5\3\5\5\3\5\5\3\5\5\3\eeee
\3\5\5\3\5\5\3\5\5\3\5\5\3\5\5\3\eeee
\7\2\2\7\2\2\7\2\2\7\2\2\7\2\2\7\eeee
\1\8\8\1\8\8\1\8\8\1\8\8\1\8\8\1\eeee
\6\4\4\6\4\4\6\4\4\6\4\4\6\4\4\6\eeee
\6\4\4\6\4\4\6\4\4\6\4\4\6\4\4\6\eeee
\1\8\8\1\8\8\1\8\8\1\8\8\1\8\8\1\eeee
} 

\baaa
8-143
\eaaa
\bbbb
0&0&0&0&1&1&1&1\\
0&0&0&0&1&1&1&1\\
0&0&0&0&1&1&1&1\\
0&0&0&0&1&1&1&1\\
1&1&1&1&0&0&0&0\\
1&1&1&1&0&0&0&0\\
1&1&1&1&0&0&0&0\\
1&1&1&1&0&0&0&0\\
\ebbb
\parbox{7cm}{ 
\1\6\3\5\1\6\3\5\1\6\3\5\1\6\3\5\eeee
\5\2\8\4\7\2\8\4\7\2\8\4\7\2\8\4\eeee
\3\7\1\6\3\5\1\6\3\5\1\6\3\5\1\6\eeee
\6\4\5\2\8\4\7\2\8\4\7\2\8\4\7\2\eeee
\1\8\3\7\1\6\3\5\1\6\3\5\1\6\3\5\eeee
\5\2\6\4\5\2\8\4\7\2\8\4\7\2\8\4\eeee
\3\7\1\8\3\7\1\6\3\5\1\6\3\5\1\6\eeee
\6\4\5\2\6\4\5\2\8\4\7\2\8\4\7\2\eeee
\1\8\3\7\1\8\3\7\1\6\3\5\1\6\3\5\eeee
\5\2\6\4\5\2\6\4\5\2\8\4\7\2\8\4\eeee
\3\7\1\8\3\7\1\8\3\7\1\6\3\5\1\6\eeee
\6\4\5\2\6\4\5\2\6\4\5\2\8\4\7\2\eeee
} 

\baaa
8-144
\eaaa
\bbbb
0&0&0&0&1&1&1&1\\
0&0&0&0&1&1&1&1\\
0&0&1&1&0&0&1&1\\
0&0&1&1&0&0&1&1\\
1&1&0&0&1&1&0&0\\
1&1&0&0&1&1&0&0\\
1&1&1&1&0&0&0&0\\
1&1&1&1&0&0&0&0\\
\ebbb
\parbox{7cm}{ 
\1\6\6\1\7\4\4\7\1\6\6\1\7\4\4\7\eeee
\5\5\2\8\3\3\8\2\5\5\2\8\3\3\8\2\eeee
\6\1\7\4\4\7\1\6\6\1\7\4\4\7\1\6\eeee
\2\8\3\3\8\2\5\5\2\8\3\3\8\2\5\5\eeee
\7\4\4\7\1\6\6\1\7\4\4\7\1\6\6\1\eeee
\3\3\8\2\5\5\2\8\3\3\8\2\5\5\2\8\eeee
\4\7\1\6\6\1\7\4\4\7\1\6\6\1\7\4\eeee
\8\2\5\5\2\8\3\3\8\2\5\5\2\8\3\3\eeee
\1\6\6\1\7\4\4\7\1\6\6\1\7\4\4\7\eeee
\5\5\2\8\3\3\8\2\5\5\2\8\3\3\8\2\eeee
\6\1\7\4\4\7\1\6\6\1\7\4\4\7\1\6\eeee
\2\8\3\3\8\2\5\5\2\8\3\3\8\2\5\5\eeee
} 

\baaa
8-145
\eaaa
\bbbb
0&0&0&0&1&1&1&1\\
0&0&0&1&0&1&1&1\\
0&0&1&0&1&0&1&1\\
0&1&0&1&0&1&0&1\\
1&0&1&0&1&0&1&0\\
1&1&0&1&0&1&0&0\\
1&1&1&0&1&0&0&0\\
1&1&1&1&0&0&0&0\\
\ebbb
\parbox{7cm}{ 
\1\7\3\5\5\3\7\1\8\2\6\4\4\6\2\8\eeee
\5\5\3\7\1\8\2\6\4\4\6\2\8\1\7\3\eeee
\7\1\8\2\6\4\4\6\2\8\1\7\3\5\5\3\eeee
\2\6\4\4\6\2\8\1\7\3\5\5\3\7\1\8\eeee
\4\6\2\8\1\7\3\5\5\3\7\1\8\2\6\4\eeee
\8\1\7\3\5\5\3\7\1\8\2\6\4\4\6\2\eeee
\3\5\5\3\7\1\8\2\6\4\4\6\2\8\1\7\eeee
\3\7\1\8\2\6\4\4\6\2\8\1\7\3\5\5\eeee
\8\2\6\4\4\6\2\8\1\7\3\5\5\3\7\1\eeee
\4\4\6\2\8\1\7\3\5\5\3\7\1\8\2\6\eeee
\2\8\1\7\3\5\5\3\7\1\8\2\6\4\4\6\eeee
\7\3\5\5\3\7\1\8\2\6\4\4\6\2\8\1\eeee
} 

\baaa
8-146
\eaaa
\bbbb
0&0&0&0&1&1&1&1\\
0&0&1&1&0&0&1&1\\
0&1&0&1&0&1&0&1\\
0&1&1&0&1&0&0&1\\
1&0&0&1&0&1&1&0\\
1&0&1&0&1&0&1&0\\
1&1&0&0&1&1&0&0\\
1&1&1&1&0&0&0&0\\
\ebbb
\parbox{7cm}{ 
\1\6\5\7\1\6\5\7\1\6\5\7\1\6\5\7\eeee
\5\7\1\6\5\7\1\6\5\7\1\6\5\7\1\6\eeee
\4\2\8\3\4\2\8\3\4\2\8\3\4\2\8\3\eeee
\8\3\4\2\8\3\4\2\8\3\4\2\8\3\4\2\eeee
\1\6\5\7\1\6\5\7\1\6\5\7\1\6\5\7\eeee
\5\7\1\6\5\7\1\6\5\7\1\6\5\7\1\6\eeee
\4\2\8\3\4\2\8\3\4\2\8\3\4\2\8\3\eeee
\8\3\4\2\8\3\4\2\8\3\4\2\8\3\4\2\eeee
\1\6\5\7\1\6\5\7\1\6\5\7\1\6\5\7\eeee
\5\7\1\6\5\7\1\6\5\7\1\6\5\7\1\6\eeee
\4\2\8\3\4\2\8\3\4\2\8\3\4\2\8\3\eeee
\8\3\4\2\8\3\4\2\8\3\4\2\8\3\4\2\eeee
} 

\baaa
8-147
\eaaa
\bbbb
0&0&0&0&1&1&1&1\\
0&0&1&1&0&0&1&1\\
0&1&0&1&0&1&0&1\\
0&1&1&0&1&0&0&1\\
1&0&0&1&0&1&1&0\\
1&0&1&0&1&1&0&0\\
1&1&0&0&1&0&1&0\\
1&1&1&1&0&0&0&0\\
\ebbb
\parbox{7cm}{ 
\1\6\5\7\1\6\5\7\1\6\5\7\1\6\5\7\eeee
\5\6\1\7\5\6\1\7\5\6\1\7\5\6\1\7\eeee
\4\3\8\2\4\3\8\2\4\3\8\2\4\3\8\2\eeee
\8\2\4\3\8\2\4\3\8\2\4\3\8\2\4\3\eeee
\1\7\5\6\1\7\5\6\1\7\5\6\1\7\5\6\eeee
\5\7\1\6\5\7\1\6\5\7\1\6\5\7\1\6\eeee
\4\2\8\3\4\2\8\3\4\2\8\3\4\2\8\3\eeee
\8\3\4\2\8\3\4\2\8\3\4\2\8\3\4\2\eeee
\1\6\5\7\1\6\5\7\1\6\5\7\1\6\5\7\eeee
\5\6\1\7\5\6\1\7\5\6\1\7\5\6\1\7\eeee
\4\3\8\2\4\3\8\2\4\3\8\2\4\3\8\2\eeee
\8\2\4\3\8\2\4\3\8\2\4\3\8\2\4\3\eeee
} 

\baaa
8-148
\eaaa
\bbbb
0&0&0&0&1&1&1&1\\
0&0&1&1&0&0&1&1\\
0&1&0&1&0&1&0&1\\
0&1&1&0&1&0&1&0\\
1&0&0&1&0&1&1&0\\
1&0&1&0&1&0&0&1\\
1&1&0&1&1&0&0&0\\
1&1&1&0&0&1&0&0\\
\ebbb
\parbox{7cm}{ 
\1\7\2\8\1\7\2\8\1\7\2\8\1\7\2\8\eeee
\5\4\3\6\5\4\3\6\5\4\3\6\5\4\3\6\eeee
\7\2\8\1\7\2\8\1\7\2\8\1\7\2\8\1\eeee
\4\3\6\5\4\3\6\5\4\3\6\5\4\3\6\5\eeee
\2\8\1\7\2\8\1\7\2\8\1\7\2\8\1\7\eeee
\3\6\5\4\3\6\5\4\3\6\5\4\3\6\5\4\eeee
\8\1\7\2\8\1\7\2\8\1\7\2\8\1\7\2\eeee
\6\5\4\3\6\5\4\3\6\5\4\3\6\5\4\3\eeee
\1\7\2\8\1\7\2\8\1\7\2\8\1\7\2\8\eeee
\5\4\3\6\5\4\3\6\5\4\3\6\5\4\3\6\eeee
\7\2\8\1\7\2\8\1\7\2\8\1\7\2\8\1\eeee
\4\3\6\5\4\3\6\5\4\3\6\5\4\3\6\5\eeee
} 

\baaa
8-149
\eaaa
\bbbb
0&0&0&0&1&1&1&1\\
0&0&1&1&0&0&1&1\\
0&1&0&1&0&1&0&1\\
0&1&1&1&1&0&0&0\\
1&0&0&1&0&1&1&0\\
1&0&1&0&1&1&0&0\\
1&1&0&0&1&0&1&0\\
1&1&1&0&0&0&0&1\\
\ebbb
\parbox{7cm}{ 
\1\6\5\7\1\6\5\7\1\6\5\7\1\6\5\7\eeee
\5\6\1\7\5\6\1\7\5\6\1\7\5\6\1\7\eeee
\4\3\8\2\4\3\8\2\4\3\8\2\4\3\8\2\eeee
\4\2\8\3\4\2\8\3\4\2\8\3\4\2\8\3\eeee
\5\7\1\6\5\7\1\6\5\7\1\6\5\7\1\6\eeee
\1\7\5\6\1\7\5\6\1\7\5\6\1\7\5\6\eeee
\8\2\4\3\8\2\4\3\8\2\4\3\8\2\4\3\eeee
\8\3\4\2\8\3\4\2\8\3\4\2\8\3\4\2\eeee
\1\6\5\7\1\6\5\7\1\6\5\7\1\6\5\7\eeee
\5\6\1\7\5\6\1\7\5\6\1\7\5\6\1\7\eeee
\4\3\8\2\4\3\8\2\4\3\8\2\4\3\8\2\eeee
\4\2\8\3\4\2\8\3\4\2\8\3\4\2\8\3\eeee
} 

\baaa
8-150
\eaaa
\bbbb
0&0&0&0&1&1&1&1\\
0&0&1&1&0&0&1&1\\
0&1&1&0&0&1&0&1\\
0&1&0&1&1&0&0&1\\
1&0&0&1&1&0&1&0\\
1&0&1&0&0&1&1&0\\
1&1&0&0&1&1&0&0\\
1&1&1&1&0&0&0&0\\
\ebbb
\parbox{7cm}{ 
\1\6\6\7\1\6\6\7\1\6\6\7\1\6\6\7\eeee
\5\7\1\5\5\7\1\5\5\7\1\5\5\7\1\5\eeee
\4\2\8\4\4\2\8\4\4\2\8\4\4\2\8\4\eeee
\8\3\3\2\8\3\3\2\8\3\3\2\8\3\3\2\eeee
\1\6\6\7\1\6\6\7\1\6\6\7\1\6\6\7\eeee
\5\7\1\5\5\7\1\5\5\7\1\5\5\7\1\5\eeee
\4\2\8\4\4\2\8\4\4\2\8\4\4\2\8\4\eeee
\8\3\3\2\8\3\3\2\8\3\3\2\8\3\3\2\eeee
\1\6\6\7\1\6\6\7\1\6\6\7\1\6\6\7\eeee
\5\7\1\5\5\7\1\5\5\7\1\5\5\7\1\5\eeee
\4\2\8\4\4\2\8\4\4\2\8\4\4\2\8\4\eeee
\8\3\3\2\8\3\3\2\8\3\3\2\8\3\3\2\eeee
} 

\baaa
\phantom{0-00}\#
\eaaa
\mbox{}\phantom{\bbbb
0&0&0&0&0&0&0&0\\
\ebbb}\parbox{7cm}{ 
\1\7\2\8\1\7\2\8\1\7\2\8\1\7\2\8\eeee
\5\5\4\4\5\5\4\4\5\5\4\4\5\5\4\4\eeee
\7\1\8\2\7\1\8\2\7\1\8\2\7\1\8\2\eeee
\6\6\3\3\6\6\3\3\6\6\3\3\6\6\3\3\eeee
\1\7\2\8\1\7\2\8\1\7\2\8\1\7\2\8\eeee
\5\5\4\4\5\5\4\4\5\5\4\4\5\5\4\4\eeee
\7\1\8\2\7\1\8\2\7\1\8\2\7\1\8\2\eeee
\6\6\3\3\6\6\3\3\6\6\3\3\6\6\3\3\eeee
\1\7\2\8\1\7\2\8\1\7\2\8\1\7\2\8\eeee
\5\5\4\4\5\5\4\4\5\5\4\4\5\5\4\4\eeee
\7\1\8\2\7\1\8\2\7\1\8\2\7\1\8\2\eeee
\6\6\3\3\6\6\3\3\6\6\3\3\6\6\3\3\eeee
} 

\baaa
8-151
\eaaa
\bbbb
0&0&0&0&1&1&1&1\\
0&0&1&1&0&0&1&1\\
0&1&1&0&0&1&0&1\\
0&1&0&1&1&0&1&0\\
1&0&0&1&1&0&0&1\\
1&0&1&0&0&1&1&0\\
1&1&0&1&0&1&0&0\\
1&1&1&0&1&0&0&0\\
\ebbb
\parbox{7cm}{ 
\1\7\2\8\1\7\2\8\1\7\2\8\1\7\2\8\eeee
\5\4\4\5\5\4\4\5\5\4\4\5\5\4\4\5\eeee
\8\2\7\1\8\2\7\1\8\2\7\1\8\2\7\1\eeee
\3\3\6\6\3\3\6\6\3\3\6\6\3\3\6\6\eeee
\2\8\1\7\2\8\1\7\2\8\1\7\2\8\1\7\eeee
\4\5\5\4\4\5\5\4\4\5\5\4\4\5\5\4\eeee
\7\1\8\2\7\1\8\2\7\1\8\2\7\1\8\2\eeee
\6\6\3\3\6\6\3\3\6\6\3\3\6\6\3\3\eeee
\1\7\2\8\1\7\2\8\1\7\2\8\1\7\2\8\eeee
\5\4\4\5\5\4\4\5\5\4\4\5\5\4\4\5\eeee
\8\2\7\1\8\2\7\1\8\2\7\1\8\2\7\1\eeee
\3\3\6\6\3\3\6\6\3\3\6\6\3\3\6\6\eeee
} 

\baaa
8-152
\eaaa
\bbbb
0&0&0&0&1&1&1&1\\
0&1&0&1&0&0&1&1\\
0&0&1&1&0&1&0&1\\
0&1&1&1&1&0&0&0\\
1&0&0&1&0&1&1&0\\
1&0&1&0&1&0&1&0\\
1&1&0&0&1&1&0&0\\
1&1&1&0&0&0&0&1\\
\ebbb
\parbox{7cm}{ 
\1\6\5\7\1\6\5\7\1\6\5\7\1\6\5\7\eeee
\5\7\1\6\5\7\1\6\5\7\1\6\5\7\1\6\eeee
\4\2\8\3\4\2\8\3\4\2\8\3\4\2\8\3\eeee
\4\2\8\3\4\2\8\3\4\2\8\3\4\2\8\3\eeee
\5\7\1\6\5\7\1\6\5\7\1\6\5\7\1\6\eeee
\1\6\5\7\1\6\5\7\1\6\5\7\1\6\5\7\eeee
\8\3\4\2\8\3\4\2\8\3\4\2\8\3\4\2\eeee
\8\3\4\2\8\3\4\2\8\3\4\2\8\3\4\2\eeee
\1\6\5\7\1\6\5\7\1\6\5\7\1\6\5\7\eeee
\5\7\1\6\5\7\1\6\5\7\1\6\5\7\1\6\eeee
\4\2\8\3\4\2\8\3\4\2\8\3\4\2\8\3\eeee
\4\2\8\3\4\2\8\3\4\2\8\3\4\2\8\3\eeee
} 

\baaa
8-153
\eaaa
\bbbb
0&0&0&0&1&1&1&1\\
0&1&0&1&0&0&1&1\\
0&0&1&1&0&1&0&1\\
0&1&1&1&1&0&0&0\\
1&0&0&1&0&1&1&0\\
1&0&1&0&1&1&0&0\\
1&1&0&0&1&0&1&0\\
1&1&1&0&0&0&0&1\\
\ebbb
\parbox{7cm}{ 
\1\6\5\7\1\6\5\7\1\6\5\7\1\6\5\7\eeee
\5\6\1\7\5\6\1\7\5\6\1\7\5\6\1\7\eeee
\4\3\8\2\4\3\8\2\4\3\8\2\4\3\8\2\eeee
\4\3\8\2\4\3\8\2\4\3\8\2\4\3\8\2\eeee
\5\6\1\7\5\6\1\7\5\6\1\7\5\6\1\7\eeee
\1\6\5\7\1\6\5\7\1\6\5\7\1\6\5\7\eeee
\8\3\4\2\8\3\4\2\8\3\4\2\8\3\4\2\eeee
\8\3\4\2\8\3\4\2\8\3\4\2\8\3\4\2\eeee
\1\6\5\7\1\6\5\7\1\6\5\7\1\6\5\7\eeee
\5\6\1\7\5\6\1\7\5\6\1\7\5\6\1\7\eeee
\4\3\8\2\4\3\8\2\4\3\8\2\4\3\8\2\eeee
\4\3\8\2\4\3\8\2\4\3\8\2\4\3\8\2\eeee
} 

\baaa
8-154
\eaaa
\bbbb
1&0&0&0&0&0&1&2\\
0&1&0&0&0&1&1&1\\
0&0&1&0&1&1&0&1\\
0&0&0&1&1&2&0&0\\
0&0&2&1&1&0&0&0\\
0&1&1&1&0&1&0&0\\
1&2&0&0&0&0&1&0\\
1&1&1&0&0&0&0&1\\
\ebbb
\parbox{7cm}{ 
\1\8\3\5\3\8\1\8\3\5\3\8\1\8\3\5\eeee
\1\8\3\5\3\8\1\8\3\5\3\8\1\8\3\5\eeee
\7\2\6\4\6\2\7\2\6\4\6\2\7\2\6\4\eeee
\7\2\6\4\6\2\7\2\6\4\6\2\7\2\6\4\eeee
\1\8\3\5\3\8\1\8\3\5\3\8\1\8\3\5\eeee
\1\8\3\5\3\8\1\8\3\5\3\8\1\8\3\5\eeee
\7\2\6\4\6\2\7\2\6\4\6\2\7\2\6\4\eeee
\7\2\6\4\6\2\7\2\6\4\6\2\7\2\6\4\eeee
\1\8\3\5\3\8\1\8\3\5\3\8\1\8\3\5\eeee
\1\8\3\5\3\8\1\8\3\5\3\8\1\8\3\5\eeee
\7\2\6\4\6\2\7\2\6\4\6\2\7\2\6\4\eeee
\7\2\6\4\6\2\7\2\6\4\6\2\7\2\6\4\eeee
} 

\baaa
8-155
\eaaa
\bbbb
1&0&0&0&0&0&1&2\\
0&1&0&0&0&1&1&1\\
0&0&1&0&1&1&0&1\\
0&0&0&1&2&1&0&0\\
0&0&1&2&1&0&0&0\\
0&1&1&1&0&1&0&0\\
1&2&0&0&0&0&1&0\\
1&1&1&0&0&0&0&1\\
\ebbb
\parbox{7cm}{ 
\1\8\3\5\4\6\2\7\2\6\4\5\3\8\1\8\eeee
\1\8\3\5\4\6\2\7\2\6\4\5\3\8\1\8\eeee
\7\2\6\4\5\3\8\1\8\3\5\4\6\2\7\2\eeee
\7\2\6\4\5\3\8\1\8\3\5\4\6\2\7\2\eeee
\1\8\3\5\4\6\2\7\2\6\4\5\3\8\1\8\eeee
\1\8\3\5\4\6\2\7\2\6\4\5\3\8\1\8\eeee
\7\2\6\4\5\3\8\1\8\3\5\4\6\2\7\2\eeee
\7\2\6\4\5\3\8\1\8\3\5\4\6\2\7\2\eeee
\1\8\3\5\4\6\2\7\2\6\4\5\3\8\1\8\eeee
\1\8\3\5\4\6\2\7\2\6\4\5\3\8\1\8\eeee
\7\2\6\4\5\3\8\1\8\3\5\4\6\2\7\2\eeee
\7\2\6\4\5\3\8\1\8\3\5\4\6\2\7\2\eeee
} 

\baaa
8-156
\eaaa
\bbbb
1&0&0&0&0&0&1&2\\
0&1&0&0&0&1&1&1\\
0&0&1&0&1&1&0&1\\
0&0&0&2&1&1&0&0\\
0&0&1&1&2&0&0&0\\
0&1&1&1&0&1&0&0\\
1&2&0&0&0&0&1&0\\
1&1&1&0&0&0&0&1\\
\ebbb
\parbox{7cm}{ 
\1\8\3\5\5\3\8\1\8\3\5\5\3\8\1\8\eeee
\1\8\3\5\5\3\8\1\8\3\5\5\3\8\1\8\eeee
\7\2\6\4\4\6\2\7\2\6\4\4\6\2\7\2\eeee
\7\2\6\4\4\6\2\7\2\6\4\4\6\2\7\2\eeee
\1\8\3\5\5\3\8\1\8\3\5\5\3\8\1\8\eeee
\1\8\3\5\5\3\8\1\8\3\5\5\3\8\1\8\eeee
\7\2\6\4\4\6\2\7\2\6\4\4\6\2\7\2\eeee
\7\2\6\4\4\6\2\7\2\6\4\4\6\2\7\2\eeee
\1\8\3\5\5\3\8\1\8\3\5\5\3\8\1\8\eeee
\1\8\3\5\5\3\8\1\8\3\5\5\3\8\1\8\eeee
\7\2\6\4\4\6\2\7\2\6\4\4\6\2\7\2\eeee
\7\2\6\4\4\6\2\7\2\6\4\4\6\2\7\2\eeee
} 

\baaa
8-157
\eaaa
\bbbb
1&0&0&0&0&0&1&2\\
0&1&0&0&0&1&1&1\\
0&0&1&0&1&1&1&0\\
0&0&0&1&2&1&0&0\\
0&0&1&2&1&0&0&0\\
0&1&1&1&0&1&0&0\\
1&1&1&0&0&0&1&0\\
2&1&0&0&0&0&0&1\\
\ebbb
\parbox{7cm}{ 
\1\7\3\5\4\6\2\8\1\7\3\5\4\6\2\8\eeee
\1\7\3\5\4\6\2\8\1\7\3\5\4\6\2\8\eeee
\8\2\6\4\5\3\7\1\8\2\6\4\5\3\7\1\eeee
\8\2\6\4\5\3\7\1\8\2\6\4\5\3\7\1\eeee
\1\7\3\5\4\6\2\8\1\7\3\5\4\6\2\8\eeee
\1\7\3\5\4\6\2\8\1\7\3\5\4\6\2\8\eeee
\8\2\6\4\5\3\7\1\8\2\6\4\5\3\7\1\eeee
\8\2\6\4\5\3\7\1\8\2\6\4\5\3\7\1\eeee
\1\7\3\5\4\6\2\8\1\7\3\5\4\6\2\8\eeee
\1\7\3\5\4\6\2\8\1\7\3\5\4\6\2\8\eeee
\8\2\6\4\5\3\7\1\8\2\6\4\5\3\7\1\eeee
\8\2\6\4\5\3\7\1\8\2\6\4\5\3\7\1\eeee
} 

\baaa
8-158
\eaaa
\bbbb
1&0&0&0&0&0&1&2\\
0&1&0&0&0&1&1&1\\
0&0&1&0&1&1&1&0\\
0&0&0&2&1&1&0&0\\
0&0&1&1&2&0&0&0\\
0&1&1&1&0&1&0&0\\
1&1&1&0&0&0&1&0\\
2&1&0&0&0&0&0&1\\
\ebbb
\parbox{7cm}{ 
\1\7\3\5\5\3\7\1\8\2\6\4\4\6\2\8\eeee
\1\7\3\5\5\3\7\1\8\2\6\4\4\6\2\8\eeee
\8\2\6\4\4\6\2\8\1\7\3\5\5\3\7\1\eeee
\8\2\6\4\4\6\2\8\1\7\3\5\5\3\7\1\eeee
\1\7\3\5\5\3\7\1\8\2\6\4\4\6\2\8\eeee
\1\7\3\5\5\3\7\1\8\2\6\4\4\6\2\8\eeee
\8\2\6\4\4\6\2\8\1\7\3\5\5\3\7\1\eeee
\8\2\6\4\4\6\2\8\1\7\3\5\5\3\7\1\eeee
\1\7\3\5\5\3\7\1\8\2\6\4\4\6\2\8\eeee
\1\7\3\5\5\3\7\1\8\2\6\4\4\6\2\8\eeee
\8\2\6\4\4\6\2\8\1\7\3\5\5\3\7\1\eeee
\8\2\6\4\4\6\2\8\1\7\3\5\5\3\7\1\eeee
} 

\baaa
8-159
\eaaa
\bbbb
1&0&0&0&0&1&1&1\\
0&1&0&0&1&0&1&1\\
0&0&1&0&1&1&0&1\\
0&0&0&1&1&1&1&0\\
0&1&1&1&1&0&0&0\\
1&0&1&1&0&1&0&0\\
1&1&0&1&0&0&1&0\\
1&1&1&0&0&0&0&1\\
\ebbb
\parbox{7cm}{ 
\1\6\3\8\1\6\3\8\1\6\3\8\1\6\3\8\eeee
\1\6\3\8\1\6\3\8\1\6\3\8\1\6\3\8\eeee
\7\4\5\2\7\4\5\2\7\4\5\2\7\4\5\2\eeee
\7\4\5\2\7\4\5\2\7\4\5\2\7\4\5\2\eeee
\1\6\3\8\1\6\3\8\1\6\3\8\1\6\3\8\eeee
\1\6\3\8\1\6\3\8\1\6\3\8\1\6\3\8\eeee
\7\4\5\2\7\4\5\2\7\4\5\2\7\4\5\2\eeee
\7\4\5\2\7\4\5\2\7\4\5\2\7\4\5\2\eeee
\1\6\3\8\1\6\3\8\1\6\3\8\1\6\3\8\eeee
\1\6\3\8\1\6\3\8\1\6\3\8\1\6\3\8\eeee
\7\4\5\2\7\4\5\2\7\4\5\2\7\4\5\2\eeee
\7\4\5\2\7\4\5\2\7\4\5\2\7\4\5\2\eeee
} 

\baaa
8-160
\eaaa
\bbbb
1&0&0&0&0&1&1&1\\
0&1&0&0&1&0&1&1\\
0&0&1&1&0&1&0&1\\
0&0&1&1&1&0&1&0\\
0&1&0&1&1&1&0&0\\
1&0&1&0&1&1&0&0\\
1&1&0&1&0&0&1&0\\
1&1&1&0&0&0&0&1\\
\ebbb
\parbox{7cm}{ 
\1\6\5\2\8\3\4\7\1\6\5\2\8\3\4\7\eeee
\1\6\5\2\8\3\4\7\1\6\5\2\8\3\4\7\eeee
\8\3\4\7\1\6\5\2\8\3\4\7\1\6\5\2\eeee
\8\3\4\7\1\6\5\2\8\3\4\7\1\6\5\2\eeee
\1\6\5\2\8\3\4\7\1\6\5\2\8\3\4\7\eeee
\1\6\5\2\8\3\4\7\1\6\5\2\8\3\4\7\eeee
\8\3\4\7\1\6\5\2\8\3\4\7\1\6\5\2\eeee
\8\3\4\7\1\6\5\2\8\3\4\7\1\6\5\2\eeee
\1\6\5\2\8\3\4\7\1\6\5\2\8\3\4\7\eeee
\1\6\5\2\8\3\4\7\1\6\5\2\8\3\4\7\eeee
\8\3\4\7\1\6\5\2\8\3\4\7\1\6\5\2\eeee
\8\3\4\7\1\6\5\2\8\3\4\7\1\6\5\2\eeee
} 

\baaa
8-161
\eaaa
\bbbb
1&0&0&0&0&1&1&1\\
0&1&0&0&1&0&1&1\\
0&0&1&1&0&1&0&1\\
0&0&1&2&0&1&0&0\\
0&1&0&0&2&0&1&0\\
1&0&1&1&0&1&0&0\\
1&1&0&0&1&0&1&0\\
1&1&1&0&0&0&0&1\\
\ebbb
\parbox{7cm}{ 
\1\6\4\3\8\2\5\7\1\6\4\3\8\2\5\7\eeee
\1\6\4\3\8\2\5\7\1\6\4\3\8\2\5\7\eeee
\8\3\4\6\1\7\5\2\8\3\4\6\1\7\5\2\eeee
\8\3\4\6\1\7\5\2\8\3\4\6\1\7\5\2\eeee
\1\6\4\3\8\2\5\7\1\6\4\3\8\2\5\7\eeee
\1\6\4\3\8\2\5\7\1\6\4\3\8\2\5\7\eeee
\8\3\4\6\1\7\5\2\8\3\4\6\1\7\5\2\eeee
\8\3\4\6\1\7\5\2\8\3\4\6\1\7\5\2\eeee
\1\6\4\3\8\2\5\7\1\6\4\3\8\2\5\7\eeee
\1\6\4\3\8\2\5\7\1\6\4\3\8\2\5\7\eeee
\8\3\4\6\1\7\5\2\8\3\4\6\1\7\5\2\eeee
\8\3\4\6\1\7\5\2\8\3\4\6\1\7\5\2\eeee
} 

\baaa
8-162
\eaaa
\bbbb
1&0&0&0&0&1&1&1\\
0&1&0&0&1&0&1&1\\
0&0&2&0&0&1&0&1\\
0&0&0&2&1&0&1&0\\
0&1&0&1&2&0&0&0\\
1&0&1&0&0&2&0&0\\
1&1&0&1&0&0&1&0\\
1&1&1&0&0&0&0&1\\
\ebbb
\parbox{7cm}{ 
\1\6\6\1\7\4\4\7\1\6\6\1\7\4\4\7\eeee
\1\6\6\1\7\4\4\7\1\6\6\1\7\4\4\7\eeee
\8\3\3\8\2\5\5\2\8\3\3\8\2\5\5\2\eeee
\8\3\3\8\2\5\5\2\8\3\3\8\2\5\5\2\eeee
\1\6\6\1\7\4\4\7\1\6\6\1\7\4\4\7\eeee
\1\6\6\1\7\4\4\7\1\6\6\1\7\4\4\7\eeee
\8\3\3\8\2\5\5\2\8\3\3\8\2\5\5\2\eeee
\8\3\3\8\2\5\5\2\8\3\3\8\2\5\5\2\eeee
\1\6\6\1\7\4\4\7\1\6\6\1\7\4\4\7\eeee
\1\6\6\1\7\4\4\7\1\6\6\1\7\4\4\7\eeee
\8\3\3\8\2\5\5\2\8\3\3\8\2\5\5\2\eeee
\8\3\3\8\2\5\5\2\8\3\3\8\2\5\5\2\eeee
} 

\baaa
8-163
\eaaa
\bbbb
2&0&0&0&0&0&0&2\\
0&2&0&0&0&0&1&1\\
0&0&2&0&0&1&1&0\\
0&0&0&2&1&1&0&0\\
0&0&0&1&3&0&0&0\\
0&0&1&1&0&2&0&0\\
0&1&1&0&0&0&2&0\\
1&1&0&0&0&0&0&2\\
\ebbb
\parbox{7cm}{ 
\1\8\2\7\3\6\4\5\5\4\6\3\7\2\8\1\eeee
\1\8\2\7\3\6\4\5\5\4\6\3\7\2\8\1\eeee
\1\8\2\7\3\6\4\5\5\4\6\3\7\2\8\1\eeee
\1\8\2\7\3\6\4\5\5\4\6\3\7\2\8\1\eeee
\1\8\2\7\3\6\4\5\5\4\6\3\7\2\8\1\eeee
\1\8\2\7\3\6\4\5\5\4\6\3\7\2\8\1\eeee
\1\8\2\7\3\6\4\5\5\4\6\3\7\2\8\1\eeee
\1\8\2\7\3\6\4\5\5\4\6\3\7\2\8\1\eeee
\1\8\2\7\3\6\4\5\5\4\6\3\7\2\8\1\eeee
\1\8\2\7\3\6\4\5\5\4\6\3\7\2\8\1\eeee
\1\8\2\7\3\6\4\5\5\4\6\3\7\2\8\1\eeee
\1\8\2\7\3\6\4\5\5\4\6\3\7\2\8\1\eeee
} 

\baaa
8-164
\eaaa
\bbbb
2&0&0&0&0&0&0&2\\
0&2&0&0&0&0&1&1\\
0&0&2&0&0&1&1&0\\
0&0&0&2&1&1&0&0\\
0&0&0&2&2&0&0&0\\
0&0&1&1&0&2&0&0\\
0&1&1&0&0&0&2&0\\
1&1&0&0&0&0&0&2\\
\ebbb
\parbox{7cm}{ 
\1\8\2\7\3\6\4\5\4\6\3\7\2\8\1\8\eeee
\1\8\2\7\3\6\4\5\4\6\3\7\2\8\1\8\eeee
\1\8\2\7\3\6\4\5\4\6\3\7\2\8\1\8\eeee
\1\8\2\7\3\6\4\5\4\6\3\7\2\8\1\8\eeee
\1\8\2\7\3\6\4\5\4\6\3\7\2\8\1\8\eeee
\1\8\2\7\3\6\4\5\4\6\3\7\2\8\1\8\eeee
\1\8\2\7\3\6\4\5\4\6\3\7\2\8\1\8\eeee
\1\8\2\7\3\6\4\5\4\6\3\7\2\8\1\8\eeee
\1\8\2\7\3\6\4\5\4\6\3\7\2\8\1\8\eeee
\1\8\2\7\3\6\4\5\4\6\3\7\2\8\1\8\eeee
\1\8\2\7\3\6\4\5\4\6\3\7\2\8\1\8\eeee
\1\8\2\7\3\6\4\5\4\6\3\7\2\8\1\8\eeee
} 

\baaa
8-165
\eaaa
\bbbb
2&0&0&0&0&0&1&1\\
0&2&0&0&0&1&0&1\\
0&0&2&0&1&0&1&0\\
0&0&0&2&1&1&0&0\\
0&0&1&1&2&0&0&0\\
0&1&0&1&0&2&0&0\\
1&0&1&0&0&0&2&0\\
1&1&0&0&0&0&0&2\\
\ebbb
\parbox{7cm}{ 
\1\7\3\5\4\6\2\8\1\7\3\5\4\6\2\8\eeee
\1\7\3\5\4\6\2\8\1\7\3\5\4\6\2\8\eeee
\1\7\3\5\4\6\2\8\1\7\3\5\4\6\2\8\eeee
\1\7\3\5\4\6\2\8\1\7\3\5\4\6\2\8\eeee
\1\7\3\5\4\6\2\8\1\7\3\5\4\6\2\8\eeee
\1\7\3\5\4\6\2\8\1\7\3\5\4\6\2\8\eeee
\1\7\3\5\4\6\2\8\1\7\3\5\4\6\2\8\eeee
\1\7\3\5\4\6\2\8\1\7\3\5\4\6\2\8\eeee
\1\7\3\5\4\6\2\8\1\7\3\5\4\6\2\8\eeee
\1\7\3\5\4\6\2\8\1\7\3\5\4\6\2\8\eeee
\1\7\3\5\4\6\2\8\1\7\3\5\4\6\2\8\eeee
\1\7\3\5\4\6\2\8\1\7\3\5\4\6\2\8\eeee
} 

\baaa
8-166
\eaaa
\bbbb
2&0&0&0&0&0&1&1\\
0&2&0&0&0&1&0&1\\
0&0&2&0&1&0&1&0\\
0&0&0&3&0&1&0&0\\
0&0&1&0&3&0&0&0\\
0&1&0&1&0&2&0&0\\
1&0&1&0&0&0&2&0\\
1&1&0&0&0&0&0&2\\
\ebbb
\parbox{7cm}{ 
\1\7\3\5\5\3\7\1\8\2\6\4\4\6\2\8\eeee
\1\7\3\5\5\3\7\1\8\2\6\4\4\6\2\8\eeee
\1\7\3\5\5\3\7\1\8\2\6\4\4\6\2\8\eeee
\1\7\3\5\5\3\7\1\8\2\6\4\4\6\2\8\eeee
\1\7\3\5\5\3\7\1\8\2\6\4\4\6\2\8\eeee
\1\7\3\5\5\3\7\1\8\2\6\4\4\6\2\8\eeee
\1\7\3\5\5\3\7\1\8\2\6\4\4\6\2\8\eeee
\1\7\3\5\5\3\7\1\8\2\6\4\4\6\2\8\eeee
\1\7\3\5\5\3\7\1\8\2\6\4\4\6\2\8\eeee
\1\7\3\5\5\3\7\1\8\2\6\4\4\6\2\8\eeee
\1\7\3\5\5\3\7\1\8\2\6\4\4\6\2\8\eeee
\1\7\3\5\5\3\7\1\8\2\6\4\4\6\2\8\eeee
} 



\baaa
9-1
\eaaa
\bbbb
0&0&0&0&0&0&0&0&4\\
0&0&0&0&0&0&0&0&4\\
0&0&0&0&0&0&0&2&2\\
0&0&0&0&0&0&0&2&2\\
0&0&0&0&0&0&2&2&0\\
0&0&0&0&0&0&2&2&0\\
0&0&0&0&1&1&2&0&0\\
0&0&1&1&1&1&0&0&0\\
1&1&1&1&0&0&0&0&0\\
\ebbb
\parbox{7cm}{ 
\1\9\3\8\5\7\7\5\8\3\9\1\9\3\8\5\eeee
\9\2\9\4\8\6\7\7\6\8\4\9\2\9\4\8\eeee
\3\9\1\9\3\8\5\7\7\5\8\3\9\1\9\3\eeee
\8\4\9\2\9\4\8\6\7\7\6\8\4\9\2\9\eeee
\5\8\3\9\1\9\3\8\5\7\7\5\8\3\9\1\eeee
\7\6\8\4\9\2\9\4\8\6\7\7\6\8\4\9\eeee
\7\7\5\8\3\9\1\9\3\8\5\7\7\5\8\3\eeee
\5\7\7\6\8\4\9\2\9\4\8\6\7\7\6\8\eeee
\8\6\7\7\5\8\3\9\1\9\3\8\5\7\7\5\eeee
\3\8\5\7\7\6\8\4\9\2\9\4\8\6\7\7\eeee
\9\4\8\6\7\7\5\8\3\9\1\9\3\8\5\7\eeee
\1\9\3\8\5\7\7\6\8\4\9\2\9\4\8\6\eeee
} 

\baaa
9-2
\eaaa
\bbbb
0&0&0&0&0&0&0&0&4\\
0&0&0&0&0&0&0&0&4\\
0&0&0&0&0&0&0&2&2\\
0&0&0&0&0&0&0&2&2\\
0&0&0&0&0&0&2&2&0\\
0&0&0&0&0&0&2&2&0\\
0&0&0&0&2&2&0&0&0\\
0&0&1&1&1&1&0&0&0\\
1&1&1&1&0&0&0&0&0\\
\ebbb
\parbox{7cm}{ 
\1\9\3\8\5\7\5\8\3\9\1\9\3\8\5\7\eeee
\9\2\9\4\8\6\7\6\8\4\9\2\9\4\8\6\eeee
\3\9\1\9\3\8\5\7\5\8\3\9\1\9\3\8\eeee
\8\4\9\2\9\4\8\6\7\6\8\4\9\2\9\4\eeee
\5\8\3\9\1\9\3\8\5\7\5\8\3\9\1\9\eeee
\7\6\8\4\9\2\9\4\8\6\7\6\8\4\9\2\eeee
\5\7\5\8\3\9\1\9\3\8\5\7\5\8\3\9\eeee
\8\6\7\6\8\4\9\2\9\4\8\6\7\6\8\4\eeee
\3\8\5\7\5\8\3\9\1\9\3\8\5\7\5\8\eeee
\9\4\8\6\7\6\8\4\9\2\9\4\8\6\7\6\eeee
\1\9\3\8\5\7\5\8\3\9\1\9\3\8\5\7\eeee
\9\2\9\4\8\6\7\6\8\4\9\2\9\4\8\6\eeee
} 

\baaa
9-3
\eaaa
\bbbb
0&0&0&0&0&0&0&0&4\\
0&0&0&0&0&0&0&0&4\\
0&0&0&0&0&0&0&2&2\\
0&0&0&0&0&0&0&2&2\\
0&0&0&0&0&0&2&2&0\\
0&0&0&0&0&0&4&0&0\\
0&0&0&0&2&2&0&0&0\\
0&0&1&1&2&0&0&0&0\\
1&1&1&1&0&0&0&0&0\\
\ebbb
\parbox{7cm}{ 
\1\9\3\8\5\7\6\7\5\8\3\9\1\9\3\8\eeee
\9\2\9\4\8\5\7\6\7\5\8\4\9\2\9\4\eeee
\3\9\1\9\3\8\5\7\6\7\5\8\3\9\1\9\eeee
\8\4\9\2\9\4\8\5\7\6\7\5\8\4\9\2\eeee
\5\8\3\9\1\9\3\8\5\7\6\7\5\8\3\9\eeee
\7\5\8\4\9\2\9\4\8\5\7\6\7\5\8\4\eeee
\6\7\5\8\3\9\1\9\3\8\5\7\6\7\5\8\eeee
\7\6\7\5\8\4\9\2\9\4\8\5\7\6\7\5\eeee
\5\7\6\7\5\8\3\9\1\9\3\8\5\7\6\7\eeee
\8\5\7\6\7\5\8\4\9\2\9\4\8\5\7\6\eeee
\3\8\5\7\6\7\5\8\3\9\1\9\3\8\5\7\eeee
\9\4\8\5\7\6\7\5\8\4\9\2\9\4\8\5\eeee
} 

\baaa
9-4
\eaaa
\bbbb
0&0&0&0&0&0&0&0&4\\
0&0&0&0&0&0&0&0&4\\
0&0&0&0&0&0&0&2&2\\
0&0&0&0&0&0&0&2&2\\
0&0&0&0&0&0&2&2&0\\
0&0&0&0&0&2&2&0&0\\
0&0&0&0&2&2&0&0&0\\
0&0&1&1&2&0&0&0&0\\
1&1&1&1&0&0&0&0&0\\
\ebbb
\parbox{7cm}{ 
\1\9\3\8\5\7\6\6\7\5\8\3\9\1\9\3\eeee
\9\2\9\4\8\5\7\6\6\7\5\8\4\9\2\9\eeee
\3\9\1\9\3\8\5\7\6\6\7\5\8\3\9\1\eeee
\8\4\9\2\9\4\8\5\7\6\6\7\5\8\4\9\eeee
\5\8\3\9\1\9\3\8\5\7\6\6\7\5\8\3\eeee
\7\5\8\4\9\2\9\4\8\5\7\6\6\7\5\8\eeee
\6\7\5\8\3\9\1\9\3\8\5\7\6\6\7\5\eeee
\6\6\7\5\8\4\9\2\9\4\8\5\7\6\6\7\eeee
\7\6\6\7\5\8\3\9\1\9\3\8\5\7\6\6\eeee
\5\7\6\6\7\5\8\4\9\2\9\4\8\5\7\6\eeee
\8\5\7\6\6\7\5\8\3\9\1\9\3\8\5\7\eeee
\3\8\5\7\6\6\7\5\8\4\9\2\9\4\8\5\eeee
} 

\baaa
9-5
\eaaa
\bbbb
0&0&0&0&0&0&0&0&4\\
0&0&0&0&0&0&0&0&4\\
0&0&0&0&0&0&0&2&2\\
0&0&0&0&0&0&0&2&2\\
0&0&0&0&0&0&4&0&0\\
0&0&0&0&0&0&4&0&0\\
0&0&0&0&1&1&0&2&0\\
0&0&1&1&0&0&2&0&0\\
1&1&1&1&0&0&0&0&0\\
\ebbb
\parbox{7cm}{ 
\1\9\3\8\7\5\7\8\3\9\1\9\3\8\7\5\eeee
\9\2\9\4\8\7\6\7\8\4\9\2\9\4\8\7\eeee
\3\9\1\9\3\8\7\5\7\8\3\9\1\9\3\8\eeee
\8\4\9\2\9\4\8\7\6\7\8\4\9\2\9\4\eeee
\7\8\3\9\1\9\3\8\7\5\7\8\3\9\1\9\eeee
\5\7\8\4\9\2\9\4\8\7\6\7\8\4\9\2\eeee
\7\6\7\8\3\9\1\9\3\8\7\5\7\8\3\9\eeee
\8\7\5\7\8\4\9\2\9\4\8\7\6\7\8\4\eeee
\3\8\7\6\7\8\3\9\1\9\3\8\7\5\7\8\eeee
\9\4\8\7\5\7\8\4\9\2\9\4\8\7\6\7\eeee
\1\9\3\8\7\6\7\8\3\9\1\9\3\8\7\5\eeee
\9\2\9\4\8\7\5\7\8\4\9\2\9\4\8\7\eeee
} 

\baaa
9-6
\eaaa
\bbbb
0&0&0&0&0&0&0&0&4\\
0&0&0&0&0&0&0&0&4\\
0&0&0&0&0&0&0&2&2\\
0&0&0&0&0&0&0&2&2\\
0&0&0&0&0&1&1&2&0\\
0&0&0&0&2&1&1&0&0\\
0&0&0&0&2&1&1&0&0\\
0&0&1&1&2&0&0&0&0\\
1&1&1&1&0&0&0&0&0\\
\ebbb
\parbox{7cm}{ 
\1\9\3\8\5\6\6\5\8\3\9\1\9\3\8\5\eeee
\9\2\9\4\8\5\7\7\5\8\4\9\2\9\4\8\eeee
\3\9\1\9\3\8\5\6\6\5\8\3\9\1\9\3\eeee
\8\4\9\2\9\4\8\5\7\7\5\8\4\9\2\9\eeee
\5\8\3\9\1\9\3\8\5\6\6\5\8\3\9\1\eeee
\6\5\8\4\9\2\9\4\8\5\7\7\5\8\4\9\eeee
\6\7\5\8\3\9\1\9\3\8\5\6\6\5\8\3\eeee
\5\7\6\5\8\4\9\2\9\4\8\5\7\7\5\8\eeee
\8\5\6\7\5\8\3\9\1\9\3\8\5\6\6\5\eeee
\3\8\5\7\6\5\8\4\9\2\9\4\8\5\7\7\eeee
\9\4\8\5\6\7\5\8\3\9\1\9\3\8\5\6\eeee
\1\9\3\8\5\7\6\5\8\4\9\2\9\4\8\5\eeee
} 

\baaa
9-7
\eaaa
\bbbb
0&0&0&0&0&0&0&0&4\\
0&0&0&0&0&0&0&0&4\\
0&0&0&0&0&0&0&2&2\\
0&0&0&0&0&0&2&2&0\\
0&0&0&0&0&0&2&2&0\\
0&0&0&0&0&0&4&0&0\\
0&0&0&1&1&2&0&0&0\\
0&0&2&1&1&0&0&0&0\\
1&1&2&0&0&0&0&0&0\\
\ebbb
\parbox{7cm}{ 
\1\9\3\8\4\7\6\7\4\8\3\9\1\9\3\8\eeee
\9\2\9\3\8\5\7\6\7\5\8\3\9\2\9\3\eeee
\3\9\1\9\3\8\4\7\6\7\4\8\3\9\1\9\eeee
\8\3\9\2\9\3\8\5\7\6\7\5\8\3\9\2\eeee
\4\8\3\9\1\9\3\8\4\7\6\7\4\8\3\9\eeee
\7\5\8\3\9\2\9\3\8\5\7\6\7\5\8\3\eeee
\6\7\4\8\3\9\1\9\3\8\4\7\6\7\4\8\eeee
\7\6\7\5\8\3\9\2\9\3\8\5\7\6\7\5\eeee
\4\7\6\7\4\8\3\9\1\9\3\8\4\7\6\7\eeee
\8\5\7\6\7\5\8\3\9\2\9\3\8\5\7\6\eeee
\3\8\4\7\6\7\4\8\3\9\1\9\3\8\4\7\eeee
\9\3\8\5\7\6\7\5\8\3\9\2\9\3\8\5\eeee
} 

\baaa
9-8
\eaaa
\bbbb
0&0&0&0&0&0&0&0&4\\
0&0&0&0&0&0&0&0&4\\
0&0&0&0&0&0&0&2&2\\
0&0&0&0&0&0&2&2&0\\
0&0&0&0&0&0&2&2&0\\
0&0&0&0&0&2&2&0&0\\
0&0&0&1&1&2&0&0&0\\
0&0&2&1&1&0&0&0&0\\
1&1&2&0&0&0&0&0&0\\
\ebbb
\parbox{7cm}{ 
\1\9\3\8\4\7\6\6\7\4\8\3\9\1\9\3\eeee
\9\2\9\3\8\5\7\6\6\7\5\8\3\9\2\9\eeee
\3\9\1\9\3\8\4\7\6\6\7\4\8\3\9\1\eeee
\8\3\9\2\9\3\8\5\7\6\6\7\5\8\3\9\eeee
\4\8\3\9\1\9\3\8\4\7\6\6\7\4\8\3\eeee
\7\5\8\3\9\2\9\3\8\5\7\6\6\7\5\8\eeee
\6\7\4\8\3\9\1\9\3\8\4\7\6\6\7\4\eeee
\6\6\7\5\8\3\9\2\9\3\8\5\7\6\6\7\eeee
\7\6\6\7\4\8\3\9\1\9\3\8\4\7\6\6\eeee
\4\7\6\6\7\5\8\3\9\2\9\3\8\5\7\6\eeee
\8\5\7\6\6\7\4\8\3\9\1\9\3\8\4\7\eeee
\3\8\4\7\6\6\7\5\8\3\9\2\9\3\8\5\eeee
} 

\baaa
9-9
\eaaa
\bbbb
0&0&0&0&0&0&0&0&4\\
0&0&0&0&0&0&0&0&4\\
0&0&0&0&0&0&0&2&2\\
0&0&0&0&0&0&2&2&0\\
0&0&0&0&0&0&4&0&0\\
0&0&0&0&0&0&4&0&0\\
0&0&0&2&1&1&0&0&0\\
0&0&2&2&0&0&0&0&0\\
1&1&2&0&0&0&0&0&0\\
\ebbb
\parbox{7cm}{ 
\1\9\3\8\4\7\5\7\4\8\3\9\1\9\3\8\eeee
\9\2\9\3\8\4\7\6\7\4\8\3\9\2\9\3\eeee
\3\9\1\9\3\8\4\7\5\7\4\8\3\9\1\9\eeee
\8\3\9\2\9\3\8\4\7\6\7\4\8\3\9\2\eeee
\4\8\3\9\1\9\3\8\4\7\5\7\4\8\3\9\eeee
\7\4\8\3\9\2\9\3\8\4\7\6\7\4\8\3\eeee
\5\7\4\8\3\9\1\9\3\8\4\7\5\7\4\8\eeee
\7\6\7\4\8\3\9\2\9\3\8\4\7\6\7\4\eeee
\4\7\5\7\4\8\3\9\1\9\3\8\4\7\5\7\eeee
\8\4\7\6\7\4\8\3\9\2\9\3\8\4\7\6\eeee
\3\8\4\7\5\7\4\8\3\9\1\9\3\8\4\7\eeee
\9\3\8\4\7\6\7\4\8\3\9\2\9\3\8\4\eeee
} 

\baaa
9-10
\eaaa
\bbbb
0&0&0&0&0&0&0&0&4\\
0&0&0&0&0&0&0&0&4\\
0&0&0&0&0&0&0&2&2\\
0&0&0&0&0&0&2&2&0\\
0&0&0&0&0&2&2&0&0\\
0&0&0&0&2&2&0&0&0\\
0&0&0&2&2&0&0&0&0\\
0&0&2&2&0&0&0&0&0\\
1&1&2&0&0&0&0&0&0\\
\ebbb
\parbox{7cm}{ 
\1\9\3\8\4\7\5\6\6\5\7\4\8\3\9\1\eeee
\9\2\9\3\8\4\7\5\6\6\5\7\4\8\3\9\eeee
\3\9\1\9\3\8\4\7\5\6\6\5\7\4\8\3\eeee
\8\3\9\2\9\3\8\4\7\5\6\6\5\7\4\8\eeee
\4\8\3\9\1\9\3\8\4\7\5\6\6\5\7\4\eeee
\7\4\8\3\9\2\9\3\8\4\7\5\6\6\5\7\eeee
\5\7\4\8\3\9\1\9\3\8\4\7\5\6\6\5\eeee
\6\5\7\4\8\3\9\2\9\3\8\4\7\5\6\6\eeee
\6\6\5\7\4\8\3\9\1\9\3\8\4\7\5\6\eeee
\5\6\6\5\7\4\8\3\9\2\9\3\8\4\7\5\eeee
\7\5\6\6\5\7\4\8\3\9\1\9\3\8\4\7\eeee
\4\7\5\6\6\5\7\4\8\3\9\2\9\3\8\4\eeee
} 

\baaa
9-11
\eaaa
\bbbb
0&0&0&0&0&0&0&0&4\\
0&0&0&0&0&0&0&0&4\\
0&0&0&0&0&0&0&2&2\\
0&0&0&0&0&0&2&2&0\\
0&0&0&0&0&2&2&0&0\\
0&0&0&0&4&0&0&0&0\\
0&0&0&2&2&0&0&0&0\\
0&0&2&2&0&0&0&0&0\\
1&1&2&0&0&0&0&0&0\\
\ebbb
\parbox{7cm}{ 
\1\9\3\8\4\7\5\6\5\7\4\8\3\9\1\9\eeee
\9\2\9\3\8\4\7\5\6\5\7\4\8\3\9\2\eeee
\3\9\1\9\3\8\4\7\5\6\5\7\4\8\3\9\eeee
\8\3\9\2\9\3\8\4\7\5\6\5\7\4\8\3\eeee
\4\8\3\9\1\9\3\8\4\7\5\6\5\7\4\8\eeee
\7\4\8\3\9\2\9\3\8\4\7\5\6\5\7\4\eeee
\5\7\4\8\3\9\1\9\3\8\4\7\5\6\5\7\eeee
\6\5\7\4\8\3\9\2\9\3\8\4\7\5\6\5\eeee
\5\6\5\7\4\8\3\9\1\9\3\8\4\7\5\6\eeee
\7\5\6\5\7\4\8\3\9\2\9\3\8\4\7\5\eeee
\4\7\5\6\5\7\4\8\3\9\1\9\3\8\4\7\eeee
\8\4\7\5\6\5\7\4\8\3\9\2\9\3\8\4\eeee
} 

\baaa
9-12
\eaaa
\bbbb
0&0&0&0&0&0&0&0&4\\
0&0&0&0&0&0&0&0&4\\
0&0&0&0&0&0&0&2&2\\
0&0&0&0&0&0&2&2&0\\
0&0&0&0&1&1&2&0&0\\
0&0&0&0&1&1&2&0&0\\
0&0&0&2&1&1&0&0&0\\
0&0&2&2&0&0&0&0&0\\
1&1&2&0&0&0&0&0&0\\
\ebbb
\parbox{7cm}{ 
\1\9\3\8\4\7\5\5\7\4\8\3\9\1\9\3\eeee
\9\2\9\3\8\4\7\6\6\7\4\8\3\9\2\9\eeee
\3\9\1\9\3\8\4\7\5\5\7\4\8\3\9\1\eeee
\8\3\9\2\9\3\8\4\7\6\6\7\4\8\3\9\eeee
\4\8\3\9\1\9\3\8\4\7\5\5\7\4\8\3\eeee
\7\4\8\3\9\2\9\3\8\4\7\6\6\7\4\8\eeee
\5\7\4\8\3\9\1\9\3\8\4\7\5\5\7\4\eeee
\5\6\7\4\8\3\9\2\9\3\8\4\7\6\6\7\eeee
\7\6\5\7\4\8\3\9\1\9\3\8\4\7\5\5\eeee
\4\7\5\6\7\4\8\3\9\2\9\3\8\4\7\6\eeee
\8\4\7\6\5\7\4\8\3\9\1\9\3\8\4\7\eeee
\3\8\4\7\5\6\7\4\8\3\9\2\9\3\8\4\eeee
} 

\baaa
9-13
\eaaa
\bbbb
0&0&0&0&0&0&0&0&4\\
0&0&0&0&0&0&0&0&4\\
0&0&0&0&0&0&0&2&2\\
0&0&0&0&0&1&1&2&0\\
0&0&0&0&0&1&1&2&0\\
0&0&0&1&1&1&1&0&0\\
0&0&0&1&1&1&1&0&0\\
0&0&2&1&1&0&0&0&0\\
1&1&2&0&0&0&0&0&0\\
\ebbb
\parbox{7cm}{ 
\1\9\3\8\4\6\6\4\8\3\9\1\9\3\8\4\eeee
\9\2\9\3\8\5\7\7\5\8\3\9\2\9\3\8\eeee
\3\9\1\9\3\8\4\6\6\4\8\3\9\1\9\3\eeee
\8\3\9\2\9\3\8\5\7\7\5\8\3\9\2\9\eeee
\4\8\3\9\1\9\3\8\4\6\6\4\8\3\9\1\eeee
\6\5\8\3\9\2\9\3\8\5\7\7\5\8\3\9\eeee
\6\7\4\8\3\9\1\9\3\8\4\6\6\4\8\3\eeee
\4\7\6\5\8\3\9\2\9\3\8\5\7\7\5\8\eeee
\8\5\6\7\4\8\3\9\1\9\3\8\4\6\6\4\eeee
\3\8\4\7\6\5\8\3\9\2\9\3\8\5\7\7\eeee
\9\3\8\5\6\7\4\8\3\9\1\9\3\8\4\6\eeee
\1\9\3\8\4\7\6\5\8\3\9\2\9\3\8\5\eeee
} 

\baaa
9-14
\eaaa
\bbbb
0&0&0&0&0&0&0&0&4\\
0&0&0&0&0&0&0&0&4\\
0&0&0&0&0&0&0&2&2\\
0&0&0&0&0&1&1&2&0\\
0&0&0&0&0&1&1&2&0\\
0&0&0&2&2&0&0&0&0\\
0&0&0&2&2&0&0&0&0\\
0&0&2&1&1&0&0&0&0\\
1&1&2&0&0&0&0&0&0\\
\ebbb
\parbox{7cm}{ 
\1\9\3\8\4\6\4\8\3\9\1\9\3\8\4\6\eeee
\9\2\9\3\8\5\7\5\8\3\9\2\9\3\8\5\eeee
\3\9\1\9\3\8\4\6\4\8\3\9\1\9\3\8\eeee
\8\3\9\2\9\3\8\5\7\5\8\3\9\2\9\3\eeee
\4\8\3\9\1\9\3\8\4\6\4\8\3\9\1\9\eeee
\6\5\8\3\9\2\9\3\8\5\7\5\8\3\9\2\eeee
\4\7\4\8\3\9\1\9\3\8\4\6\4\8\3\9\eeee
\8\5\6\5\8\3\9\2\9\3\8\5\7\5\8\3\eeee
\3\8\4\7\4\8\3\9\1\9\3\8\4\6\4\8\eeee
\9\3\8\5\6\5\8\3\9\2\9\3\8\5\7\5\eeee
\1\9\3\8\4\7\4\8\3\9\1\9\3\8\4\6\eeee
\9\2\9\3\8\5\6\5\8\3\9\2\9\3\8\5\eeee
} 

\baaa
9-15
\eaaa
\bbbb
0&0&0&0&0&0&0&0&4\\
0&0&0&0&0&0&0&0&4\\
0&0&0&0&0&0&0&2&2\\
0&0&0&0&0&1&1&2&0\\
0&0&0&0&0&2&2&0&0\\
0&0&0&2&2&0&0&0&0\\
0&0&0&2&2&0&0&0&0\\
0&0&2&2&0&0&0&0&0\\
1&1&2&0&0&0&0&0&0\\
\ebbb
\parbox{7cm}{ 
\1\9\3\8\4\6\5\6\4\8\3\9\1\9\3\8\eeee
\9\2\9\3\8\4\7\5\7\4\8\3\9\2\9\3\eeee
\3\9\1\9\3\8\4\6\5\6\4\8\3\9\1\9\eeee
\8\3\9\2\9\3\8\4\7\5\7\4\8\3\9\2\eeee
\4\8\3\9\1\9\3\8\4\6\5\6\4\8\3\9\eeee
\6\4\8\3\9\2\9\3\8\4\7\5\7\4\8\3\eeee
\5\7\4\8\3\9\1\9\3\8\4\6\5\6\4\8\eeee
\6\5\6\4\8\3\9\2\9\3\8\4\7\5\7\4\eeee
\4\7\5\7\4\8\3\9\1\9\3\8\4\6\5\6\eeee
\8\4\6\5\6\4\8\3\9\2\9\3\8\4\7\5\eeee
\3\8\4\7\5\7\4\8\3\9\1\9\3\8\4\6\eeee
\9\3\8\4\6\5\6\4\8\3\9\2\9\3\8\4\eeee
} 

\baaa
9-16
\eaaa
\bbbb
0&0&0&0&0&0&0&0&4\\
0&0&0&0&0&0&0&0&4\\
0&0&0&0&0&0&0&2&2\\
0&0&0&0&0&1&1&2&0\\
0&0&0&0&2&1&1&0&0\\
0&0&0&2&2&0&0&0&0\\
0&0&0&2&2&0&0&0&0\\
0&0&2&2&0&0&0&0&0\\
1&1&2&0&0&0&0&0&0\\
\ebbb
\parbox{7cm}{ 
\1\9\3\8\4\6\5\5\6\4\8\3\9\1\9\3\eeee
\9\2\9\3\8\4\7\5\5\7\4\8\3\9\2\9\eeee
\3\9\1\9\3\8\4\6\5\5\6\4\8\3\9\1\eeee
\8\3\9\2\9\3\8\4\7\5\5\7\4\8\3\9\eeee
\4\8\3\9\1\9\3\8\4\6\5\5\6\4\8\3\eeee
\6\4\8\3\9\2\9\3\8\4\7\5\5\7\4\8\eeee
\5\7\4\8\3\9\1\9\3\8\4\6\5\5\6\4\eeee
\5\5\6\4\8\3\9\2\9\3\8\4\7\5\5\7\eeee
\6\5\5\7\4\8\3\9\1\9\3\8\4\6\5\5\eeee
\4\7\5\5\6\4\8\3\9\2\9\3\8\4\7\5\eeee
\8\4\6\5\5\7\4\8\3\9\1\9\3\8\4\6\eeee
\3\8\4\7\5\5\6\4\8\3\9\2\9\3\8\4\eeee
} 

\baaa
9-17
\eaaa
\bbbb
0&0&0&0&0&0&0&0&4\\
0&0&0&0&0&0&0&0&4\\
0&0&0&0&0&0&0&4&0\\
0&0&0&0&0&1&1&0&2\\
0&0&0&0&0&1&1&0&2\\
0&0&0&1&1&0&0&2&0\\
0&0&0&1&1&0&0&2&0\\
0&0&2&0&0&1&1&0&0\\
1&1&0&1&1&0&0&0&0\\
\ebbb
\parbox{7cm}{ 
\1\9\4\6\8\3\8\6\4\9\1\9\4\6\8\3\eeee
\9\2\9\5\7\8\3\8\7\5\9\2\9\5\7\8\eeee
\4\9\1\9\4\6\8\3\8\6\4\9\1\9\4\6\eeee
\6\5\9\2\9\5\7\8\3\8\7\5\9\2\9\5\eeee
\8\7\4\9\1\9\4\6\8\3\8\6\4\9\1\9\eeee
\3\8\6\5\9\2\9\5\7\8\3\8\7\5\9\2\eeee
\8\3\8\7\4\9\1\9\4\6\8\3\8\6\4\9\eeee
\6\8\3\8\6\5\9\2\9\5\7\8\3\8\7\5\eeee
\4\7\8\3\8\7\4\9\1\9\4\6\8\3\8\6\eeee
\9\5\6\8\3\8\6\5\9\2\9\5\7\8\3\8\eeee
\1\9\4\7\8\3\8\7\4\9\1\9\4\6\8\3\eeee
\9\2\9\5\6\8\3\8\6\5\9\2\9\5\7\8\eeee
} 

\baaa
9-18
\eaaa
\bbbb
0&0&0&0&0&0&0&0&4\\
0&0&0&0&0&0&0&0&4\\
0&0&0&0&0&0&0&4&0\\
0&0&0&0&0&1&1&0&2\\
0&0&0&0&0&1&1&2&0\\
0&0&0&2&2&0&0&0&0\\
0&0&0&2&2&0&0&0&0\\
0&0&2&0&2&0&0&0&0\\
1&1&0&2&0&0&0&0&0\\
\ebbb
\parbox{7cm}{ 
\1\9\4\6\5\8\3\8\5\6\4\9\1\9\4\6\eeee
\9\2\9\4\7\5\8\3\8\5\7\4\9\2\9\4\eeee
\4\9\1\9\4\6\5\8\3\8\5\6\4\9\1\9\eeee
\6\4\9\2\9\4\7\5\8\3\8\5\7\4\9\2\eeee
\5\7\4\9\1\9\4\6\5\8\3\8\5\6\4\9\eeee
\8\5\6\4\9\2\9\4\7\5\8\3\8\5\7\4\eeee
\3\8\5\7\4\9\1\9\4\6\5\8\3\8\5\6\eeee
\8\3\8\5\6\4\9\2\9\4\7\5\8\3\8\5\eeee
\5\8\3\8\5\7\4\9\1\9\4\6\5\8\3\8\eeee
\6\5\8\3\8\5\6\4\9\2\9\4\7\5\8\3\eeee
\4\7\5\8\3\8\5\7\4\9\1\9\4\6\5\8\eeee
\9\4\6\5\8\3\8\5\6\4\9\2\9\4\7\5\eeee
} 

\baaa
9-19
\eaaa
\bbbb
0&0&0&0&0&0&0&0&4\\
0&0&0&0&0&0&0&0&4\\
0&0&0&0&0&0&1&1&2\\
0&0&0&0&0&0&1&1&2\\
0&0&0&0&0&0&2&2&0\\
0&0&0&0&0&0&2&2&0\\
0&0&1&1&1&1&0&0&0\\
0&0&1&1&1&1&0&0&0\\
1&1&1&1&0&0&0&0&0\\
\ebbb
\parbox{7cm}{ 
\1\9\3\7\5\7\3\9\1\9\3\7\5\7\3\9\eeee
\9\2\9\4\8\6\8\4\9\2\9\4\8\6\8\4\eeee
\3\9\1\9\3\7\5\7\3\9\1\9\3\7\5\7\eeee
\7\4\9\2\9\4\8\6\8\4\9\2\9\4\8\6\eeee
\5\8\3\9\1\9\3\7\5\7\3\9\1\9\3\7\eeee
\7\6\7\4\9\2\9\4\8\6\8\4\9\2\9\4\eeee
\3\8\5\8\3\9\1\9\3\7\5\7\3\9\1\9\eeee
\9\4\7\6\7\4\9\2\9\4\8\6\8\4\9\2\eeee
\1\9\3\8\5\8\3\9\1\9\3\7\5\7\3\9\eeee
\9\2\9\4\7\6\7\4\9\2\9\4\8\6\8\4\eeee
\3\9\1\9\3\8\5\8\3\9\1\9\3\7\5\7\eeee
\7\4\9\2\9\4\7\6\7\4\9\2\9\4\8\6\eeee
} 

\baaa
9-20
\eaaa
\bbbb
0&0&0&0&0&0&0&0&4\\
0&0&0&0&0&0&0&0&4\\
0&0&0&0&0&0&1&1&2\\
0&0&0&0&0&0&1&1&2\\
0&0&0&0&0&2&1&1&0\\
0&0&0&0&2&2&0&0&0\\
0&0&1&1&2&0&0&0&0\\
0&0&1&1&2&0&0&0&0\\
1&1&1&1&0&0&0&0&0\\
\ebbb
\parbox{7cm}{ 
\1\9\3\7\5\6\6\5\7\3\9\1\9\3\7\5\eeee
\9\2\9\4\8\5\6\6\5\8\4\9\2\9\4\8\eeee
\3\9\1\9\3\7\5\6\6\5\7\3\9\1\9\3\eeee
\7\4\9\2\9\4\8\5\6\6\5\8\4\9\2\9\eeee
\5\8\3\9\1\9\3\7\5\6\6\5\7\3\9\1\eeee
\6\5\7\4\9\2\9\4\8\5\6\6\5\8\4\9\eeee
\6\6\5\8\3\9\1\9\3\7\5\6\6\5\7\3\eeee
\5\6\6\5\7\4\9\2\9\4\8\5\6\6\5\8\eeee
\7\5\6\6\5\8\3\9\1\9\3\7\5\6\6\5\eeee
\3\8\5\6\6\5\7\4\9\2\9\4\8\5\6\6\eeee
\9\4\7\5\6\6\5\8\3\9\1\9\3\7\5\6\eeee
\1\9\3\8\5\6\6\5\7\4\9\2\9\4\8\5\eeee
} 

\baaa
9-21
\eaaa
\bbbb
0&0&0&0&0&0&0&0&4\\
0&0&0&0&0&0&0&0&4\\
0&0&0&0&0&0&1&1&2\\
0&0&0&0&0&0&1&1&2\\
0&0&0&0&1&1&1&1&0\\
0&0&0&0&1&1&1&1&0\\
0&0&1&1&1&1&0&0&0\\
0&0&1&1&1&1&0&0&0\\
1&1&1&1&0&0&0&0&0\\
\ebbb
\parbox{7cm}{ 
\1\9\3\7\5\5\7\3\9\1\9\3\7\5\5\7\eeee
\9\2\9\4\8\6\6\8\4\9\2\9\4\8\6\6\eeee
\3\9\1\9\3\7\5\5\7\3\9\1\9\3\7\5\eeee
\7\4\9\2\9\4\8\6\6\8\4\9\2\9\4\8\eeee
\5\8\3\9\1\9\3\7\5\5\7\3\9\1\9\3\eeee
\5\6\7\4\9\2\9\4\8\6\6\8\4\9\2\9\eeee
\7\6\5\8\3\9\1\9\3\7\5\5\7\3\9\1\eeee
\3\8\5\6\7\4\9\2\9\4\8\6\6\8\4\9\eeee
\9\4\7\6\5\8\3\9\1\9\3\7\5\5\7\3\eeee
\1\9\3\8\5\6\7\4\9\2\9\4\8\6\6\8\eeee
\9\2\9\4\7\6\5\8\3\9\1\9\3\7\5\5\eeee
\3\9\1\9\3\8\5\6\7\4\9\2\9\4\8\6\eeee
} 

\baaa
9-22
\eaaa
\bbbb
0&0&0&0&0&0&0&0&4\\
0&0&0&0&0&0&0&0&4\\
0&0&0&0&0&0&1&1&2\\
0&0&0&0&0&2&1&1&0\\
0&0&0&0&0&2&1&1&0\\
0&0&0&1&1&2&0&0&0\\
0&0&2&1&1&0&0&0&0\\
0&0&2&1&1&0&0&0&0\\
1&1&2&0&0&0&0&0&0\\
\ebbb
\parbox{7cm}{ 
\1\9\3\7\4\6\6\4\7\3\9\1\9\3\7\4\eeee
\9\2\9\3\8\5\6\6\5\8\3\9\2\9\3\8\eeee
\3\9\1\9\3\7\4\6\6\4\7\3\9\1\9\3\eeee
\7\3\9\2\9\3\8\5\6\6\5\8\3\9\2\9\eeee
\4\8\3\9\1\9\3\7\4\6\6\4\7\3\9\1\eeee
\6\5\7\3\9\2\9\3\8\5\6\6\5\8\3\9\eeee
\6\6\4\8\3\9\1\9\3\7\4\6\6\4\7\3\eeee
\4\6\6\5\7\3\9\2\9\3\8\5\6\6\5\8\eeee
\7\5\6\6\4\8\3\9\1\9\3\7\4\6\6\4\eeee
\3\8\4\6\6\5\7\3\9\2\9\3\8\5\6\6\eeee
\9\3\7\5\6\6\4\8\3\9\1\9\3\7\4\6\eeee
\1\9\3\8\4\6\6\5\7\3\9\2\9\3\8\5\eeee
} 

\baaa
9-23
\eaaa
\bbbb
0&0&0&0&0&0&0&0&4\\
0&0&0&0&0&0&0&0&4\\
0&0&0&0&0&0&1&1&2\\
0&0&0&0&0&2&1&1&0\\
0&0&0&0&0&2&1&1&0\\
0&0&0&2&2&0&0&0&0\\
0&0&2&1&1&0&0&0&0\\
0&0&2&1&1&0&0&0&0\\
1&1&2&0&0&0&0&0&0\\
\ebbb
\parbox{7cm}{ 
\1\9\3\7\4\6\4\7\3\9\1\9\3\7\4\6\eeee
\9\2\9\3\8\5\6\5\8\3\9\2\9\3\8\5\eeee
\3\9\1\9\3\7\4\6\4\7\3\9\1\9\3\7\eeee
\7\3\9\2\9\3\8\5\6\5\8\3\9\2\9\3\eeee
\4\8\3\9\1\9\3\7\4\6\4\7\3\9\1\9\eeee
\6\5\7\3\9\2\9\3\8\5\6\5\8\3\9\2\eeee
\4\6\4\8\3\9\1\9\3\7\4\6\4\7\3\9\eeee
\7\5\6\5\7\3\9\2\9\3\8\5\6\5\8\3\eeee
\3\8\4\6\4\8\3\9\1\9\3\7\4\6\4\7\eeee
\9\3\7\5\6\5\7\3\9\2\9\3\8\5\6\5\eeee
\1\9\3\8\4\6\4\8\3\9\1\9\3\7\4\6\eeee
\9\2\9\3\7\5\6\5\7\3\9\2\9\3\8\5\eeee
} 

\baaa
9-24
\eaaa
\bbbb
0&0&0&0&0&0&0&0&4\\
0&0&0&0&0&0&0&0&4\\
0&0&0&0&0&0&1&1&2\\
0&0&0&0&0&2&1&1&0\\
0&0&0&0&2&2&0&0&0\\
0&0&0&2&2&0&0&0&0\\
0&0&2&2&0&0&0&0&0\\
0&0&2&2&0&0&0&0&0\\
1&1&2&0&0&0&0&0&0\\
\ebbb
\parbox{7cm}{ 
\1\9\3\7\4\6\5\5\6\4\7\3\9\1\9\3\eeee
\9\2\9\3\8\4\6\5\5\6\4\8\3\9\2\9\eeee
\3\9\1\9\3\7\4\6\5\5\6\4\7\3\9\1\eeee
\7\3\9\2\9\3\8\4\6\5\5\6\4\8\3\9\eeee
\4\8\3\9\1\9\3\7\4\6\5\5\6\4\7\3\eeee
\6\4\7\3\9\2\9\3\8\4\6\5\5\6\4\8\eeee
\5\6\4\8\3\9\1\9\3\7\4\6\5\5\6\4\eeee
\5\5\6\4\7\3\9\2\9\3\8\4\6\5\5\6\eeee
\6\5\5\6\4\8\3\9\1\9\3\7\4\6\5\5\eeee
\4\6\5\5\6\4\7\3\9\2\9\3\8\4\6\5\eeee
\7\4\6\5\5\6\4\8\3\9\1\9\3\7\4\6\eeee
\3\8\4\6\5\5\6\4\7\3\9\2\9\3\8\4\eeee
} 

\baaa
9-25
\eaaa
\bbbb
0&0&0&0&0&0&0&0&4\\
0&0&0&0&0&0&0&0&4\\
0&0&0&0&0&0&1&1&2\\
0&0&0&0&1&1&1&1&0\\
0&0&0&2&1&1&0&0&0\\
0&0&0&2&1&1&0&0&0\\
0&0&2&2&0&0&0&0&0\\
0&0&2&2&0&0&0&0&0\\
1&1&2&0&0&0&0&0&0\\
\ebbb
\parbox{7cm}{ 
\1\9\3\7\4\5\5\4\7\3\9\1\9\3\7\4\eeee
\9\2\9\3\8\4\6\6\4\8\3\9\2\9\3\8\eeee
\3\9\1\9\3\7\4\5\5\4\7\3\9\1\9\3\eeee
\7\3\9\2\9\3\8\4\6\6\4\8\3\9\2\9\eeee
\4\8\3\9\1\9\3\7\4\5\5\4\7\3\9\1\eeee
\5\4\7\3\9\2\9\3\8\4\6\6\4\8\3\9\eeee
\5\6\4\8\3\9\1\9\3\7\4\5\5\4\7\3\eeee
\4\6\5\4\7\3\9\2\9\3\8\4\6\6\4\8\eeee
\7\4\5\6\4\8\3\9\1\9\3\7\4\5\5\4\eeee
\3\8\4\6\5\4\7\3\9\2\9\3\8\4\6\6\eeee
\9\3\7\4\5\6\4\8\3\9\1\9\3\7\4\5\eeee
\1\9\3\8\4\6\5\4\7\3\9\2\9\3\8\4\eeee
} 

\baaa
9-26
\eaaa
\bbbb
0&0&0&0&0&0&0&0&4\\
0&0&0&0&0&0&0&2&2\\
0&0&0&0&0&0&0&2&2\\
0&0&0&0&0&0&2&2&0\\
0&0&0&0&0&0&2&2&0\\
0&0&0&0&0&0&4&0&0\\
0&0&0&1&1&2&0&0&0\\
0&1&1&1&1&0&0&0&0\\
2&1&1&0&0&0&0&0&0\\
\ebbb
\parbox{7cm}{ 
\1\9\2\8\4\7\6\7\4\8\2\9\1\9\2\8\eeee
\9\1\9\3\8\5\7\6\7\5\8\3\9\1\9\3\eeee
\2\9\1\9\2\8\4\7\6\7\4\8\2\9\1\9\eeee
\8\3\9\1\9\3\8\5\7\6\7\5\8\3\9\1\eeee
\4\8\2\9\1\9\2\8\4\7\6\7\4\8\2\9\eeee
\7\5\8\3\9\1\9\3\8\5\7\6\7\5\8\3\eeee
\6\7\4\8\2\9\1\9\2\8\4\7\6\7\4\8\eeee
\7\6\7\5\8\3\9\1\9\3\8\5\7\6\7\5\eeee
\4\7\6\7\4\8\2\9\1\9\2\8\4\7\6\7\eeee
\8\5\7\6\7\5\8\3\9\1\9\3\8\5\7\6\eeee
\2\8\4\7\6\7\4\8\2\9\1\9\2\8\4\7\eeee
\9\3\8\5\7\6\7\5\8\3\9\1\9\3\8\5\eeee
} 

\baaa
9-27
\eaaa
\bbbb
0&0&0&0&0&0&0&0&4\\
0&0&0&0&0&0&0&2&2\\
0&0&0&0&0&0&0&2&2\\
0&0&0&0&0&0&2&2&0\\
0&0&0&0&0&0&2&2&0\\
0&0&0&0&0&2&2&0&0\\
0&0&0&1&1&2&0&0&0\\
0&1&1&1&1&0&0&0&0\\
2&1&1&0&0&0&0&0&0\\
\ebbb
\parbox{7cm}{ 
\1\9\2\8\4\7\6\6\7\4\8\2\9\1\9\2\eeee
\9\1\9\3\8\5\7\6\6\7\5\8\3\9\1\9\eeee
\2\9\1\9\2\8\4\7\6\6\7\4\8\2\9\1\eeee
\8\3\9\1\9\3\8\5\7\6\6\7\5\8\3\9\eeee
\4\8\2\9\1\9\2\8\4\7\6\6\7\4\8\2\eeee
\7\5\8\3\9\1\9\3\8\5\7\6\6\7\5\8\eeee
\6\7\4\8\2\9\1\9\2\8\4\7\6\6\7\4\eeee
\6\6\7\5\8\3\9\1\9\3\8\5\7\6\6\7\eeee
\7\6\6\7\4\8\2\9\1\9\2\8\4\7\6\6\eeee
\4\7\6\6\7\5\8\3\9\1\9\3\8\5\7\6\eeee
\8\5\7\6\6\7\4\8\2\9\1\9\2\8\4\7\eeee
\2\8\4\7\6\6\7\5\8\3\9\1\9\3\8\5\eeee
} 

\baaa
9-28
\eaaa
\bbbb
0&0&0&0&0&0&0&0&4\\
0&0&0&0&0&0&0&2&2\\
0&0&0&0&0&0&0&2&2\\
0&0&0&0&0&0&2&2&0\\
0&0&0&0&0&2&2&0&0\\
0&0&0&0&2&2&0&0&0\\
0&0&0&2&2&0&0&0&0\\
0&1&1&2&0&0&0&0&0\\
2&1&1&0&0&0&0&0&0\\
\ebbb
\parbox{7cm}{ 
\1\9\2\8\4\7\5\6\6\5\7\4\8\2\9\1\eeee
\9\1\9\3\8\4\7\5\6\6\5\7\4\8\3\9\eeee
\2\9\1\9\2\8\4\7\5\6\6\5\7\4\8\2\eeee
\8\3\9\1\9\3\8\4\7\5\6\6\5\7\4\8\eeee
\4\8\2\9\1\9\2\8\4\7\5\6\6\5\7\4\eeee
\7\4\8\3\9\1\9\3\8\4\7\5\6\6\5\7\eeee
\5\7\4\8\2\9\1\9\2\8\4\7\5\6\6\5\eeee
\6\5\7\4\8\3\9\1\9\3\8\4\7\5\6\6\eeee
\6\6\5\7\4\8\2\9\1\9\2\8\4\7\5\6\eeee
\5\6\6\5\7\4\8\3\9\1\9\3\8\4\7\5\eeee
\7\5\6\6\5\7\4\8\2\9\1\9\2\8\4\7\eeee
\4\7\5\6\6\5\7\4\8\3\9\1\9\3\8\4\eeee
} 

\baaa
9-29
\eaaa
\bbbb
0&0&0&0&0&0&0&0&4\\
0&0&0&0&0&0&0&2&2\\
0&0&0&0&0&0&0&2&2\\
0&0&0&0&0&0&2&2&0\\
0&0&0&0&0&2&2&0&0\\
0&0&0&0&4&0&0&0&0\\
0&0&0&2&2&0&0&0&0\\
0&1&1&2&0&0&0&0&0\\
2&1&1&0&0&0&0&0&0\\
\ebbb
\parbox{7cm}{ 
\1\9\2\8\4\7\5\6\5\7\4\8\2\9\1\9\eeee
\9\1\9\3\8\4\7\5\6\5\7\4\8\3\9\1\eeee
\2\9\1\9\2\8\4\7\5\6\5\7\4\8\2\9\eeee
\8\3\9\1\9\3\8\4\7\5\6\5\7\4\8\3\eeee
\4\8\2\9\1\9\2\8\4\7\5\6\5\7\4\8\eeee
\7\4\8\3\9\1\9\3\8\4\7\5\6\5\7\4\eeee
\5\7\4\8\2\9\1\9\2\8\4\7\5\6\5\7\eeee
\6\5\7\4\8\3\9\1\9\3\8\4\7\5\6\5\eeee
\5\6\5\7\4\8\2\9\1\9\2\8\4\7\5\6\eeee
\7\5\6\5\7\4\8\3\9\1\9\3\8\4\7\5\eeee
\4\7\5\6\5\7\4\8\2\9\1\9\2\8\4\7\eeee
\8\4\7\5\6\5\7\4\8\3\9\1\9\3\8\4\eeee
} 

\baaa
9-30
\eaaa
\bbbb
0&0&0&0&0&0&0&0&4\\
0&0&0&0&0&0&0&2&2\\
0&0&0&0&0&0&0&2&2\\
0&0&0&0&0&0&2&2&0\\
0&0&0&0&1&1&2&0&0\\
0&0&0&0&1&1&2&0&0\\
0&0&0&2&1&1&0&0&0\\
0&1&1&2&0&0&0&0&0\\
2&1&1&0&0&0&0&0&0\\
\ebbb
\parbox{7cm}{ 
\1\9\2\8\4\7\5\5\7\4\8\2\9\1\9\2\eeee
\9\1\9\3\8\4\7\6\6\7\4\8\3\9\1\9\eeee
\2\9\1\9\2\8\4\7\5\5\7\4\8\2\9\1\eeee
\8\3\9\1\9\3\8\4\7\6\6\7\4\8\3\9\eeee
\4\8\2\9\1\9\2\8\4\7\5\5\7\4\8\2\eeee
\7\4\8\3\9\1\9\3\8\4\7\6\6\7\4\8\eeee
\5\7\4\8\2\9\1\9\2\8\4\7\5\5\7\4\eeee
\5\6\7\4\8\3\9\1\9\3\8\4\7\6\6\7\eeee
\7\6\5\7\4\8\2\9\1\9\2\8\4\7\5\5\eeee
\4\7\5\6\7\4\8\3\9\1\9\3\8\4\7\6\eeee
\8\4\7\6\5\7\4\8\2\9\1\9\2\8\4\7\eeee
\2\8\4\7\5\6\7\4\8\3\9\1\9\3\8\4\eeee
} 

\baaa
9-31
\eaaa
\bbbb
0&0&0&0&0&0&0&0&4\\
0&0&0&0&0&0&0&2&2\\
0&0&0&0&0&0&0&2&2\\
0&0&0&0&0&1&1&2&0\\
0&0&0&0&0&1&1&2&0\\
0&0&0&1&1&1&1&0&0\\
0&0&0&1&1&1&1&0&0\\
0&1&1&1&1&0&0&0&0\\
2&1&1&0&0&0&0&0&0\\
\ebbb
\parbox{7cm}{ 
\1\9\2\8\4\6\6\4\8\2\9\1\9\2\8\4\eeee
\9\1\9\3\8\5\7\7\5\8\3\9\1\9\3\8\eeee
\2\9\1\9\2\8\4\6\6\4\8\2\9\1\9\2\eeee
\8\3\9\1\9\3\8\5\7\7\5\8\3\9\1\9\eeee
\4\8\2\9\1\9\2\8\4\6\6\4\8\2\9\1\eeee
\6\5\8\3\9\1\9\3\8\5\7\7\5\8\3\9\eeee
\6\7\4\8\2\9\1\9\2\8\4\6\6\4\8\2\eeee
\4\7\6\5\8\3\9\1\9\3\8\5\7\7\5\8\eeee
\8\5\6\7\4\8\2\9\1\9\2\8\4\6\6\4\eeee
\2\8\4\7\6\5\8\3\9\1\9\3\8\5\7\7\eeee
\9\3\8\5\6\7\4\8\2\9\1\9\2\8\4\6\eeee
\1\9\2\8\4\7\6\5\8\3\9\1\9\3\8\5\eeee
} 

\baaa
9-32
\eaaa
\bbbb
0&0&0&0&0&0&0&0&4\\
0&0&0&0&0&0&0&2&2\\
0&0&0&0&0&0&0&2&2\\
0&0&0&0&0&1&1&2&0\\
0&0&0&0&0&1&1&2&0\\
0&0&0&2&2&0&0&0&0\\
0&0&0&2&2&0&0&0&0\\
0&1&1&1&1&0&0&0&0\\
2&1&1&0&0&0&0&0&0\\
\ebbb
\parbox{7cm}{ 
\1\9\2\8\4\6\4\8\2\9\1\9\2\8\4\6\eeee
\9\1\9\3\8\5\7\5\8\3\9\1\9\3\8\5\eeee
\2\9\1\9\2\8\4\6\4\8\2\9\1\9\2\8\eeee
\8\3\9\1\9\3\8\5\7\5\8\3\9\1\9\3\eeee
\4\8\2\9\1\9\2\8\4\6\4\8\2\9\1\9\eeee
\6\5\8\3\9\1\9\3\8\5\7\5\8\3\9\1\eeee
\4\7\4\8\2\9\1\9\2\8\4\6\4\8\2\9\eeee
\8\5\6\5\8\3\9\1\9\3\8\5\7\5\8\3\eeee
\2\8\4\7\4\8\2\9\1\9\2\8\4\6\4\8\eeee
\9\3\8\5\6\5\8\3\9\1\9\3\8\5\7\5\eeee
\1\9\2\8\4\7\4\8\2\9\1\9\2\8\4\6\eeee
\9\1\9\3\8\5\6\5\8\3\9\1\9\3\8\5\eeee
} 

\baaa
9-33
\eaaa
\bbbb
0&0&0&0&0&0&0&0&4\\
0&0&0&0&0&0&0&2&2\\
0&0&0&0&0&0&0&2&2\\
0&0&0&0&0&1&1&2&0\\
0&0&0&0&0&2&2&0&0\\
0&0&0&2&2&0&0&0&0\\
0&0&0&2&2&0&0&0&0\\
0&1&1&2&0&0&0&0&0\\
2&1&1&0&0&0&0&0&0\\
\ebbb
\parbox{7cm}{ 
\1\9\2\8\4\6\5\6\4\8\2\9\1\9\2\8\eeee
\9\1\9\3\8\4\7\5\7\4\8\3\9\1\9\3\eeee
\2\9\1\9\2\8\4\6\5\6\4\8\2\9\1\9\eeee
\8\3\9\1\9\3\8\4\7\5\7\4\8\3\9\1\eeee
\4\8\2\9\1\9\2\8\4\6\5\6\4\8\2\9\eeee
\6\4\8\3\9\1\9\3\8\4\7\5\7\4\8\3\eeee
\5\7\4\8\2\9\1\9\2\8\4\6\5\6\4\8\eeee
\6\5\6\4\8\3\9\1\9\3\8\4\7\5\7\4\eeee
\4\7\5\7\4\8\2\9\1\9\2\8\4\6\5\6\eeee
\8\4\6\5\6\4\8\3\9\1\9\3\8\4\7\5\eeee
\2\8\4\7\5\7\4\8\2\9\1\9\2\8\4\6\eeee
\9\3\8\4\6\5\6\4\8\3\9\1\9\3\8\4\eeee
} 

\baaa
9-34
\eaaa
\bbbb
0&0&0&0&0&0&0&0&4\\
0&0&0&0&0&0&0&2&2\\
0&0&0&0&0&0&0&2&2\\
0&0&0&0&0&1&1&2&0\\
0&0&0&0&2&1&1&0&0\\
0&0&0&2&2&0&0&0&0\\
0&0&0&2&2&0&0&0&0\\
0&1&1&2&0&0&0&0&0\\
2&1&1&0&0&0&0&0&0\\
\ebbb
\parbox{7cm}{ 
\1\9\2\8\4\6\5\5\6\4\8\2\9\1\9\2\eeee
\9\1\9\3\8\4\7\5\5\7\4\8\3\9\1\9\eeee
\2\9\1\9\2\8\4\6\5\5\6\4\8\2\9\1\eeee
\8\3\9\1\9\3\8\4\7\5\5\7\4\8\3\9\eeee
\4\8\2\9\1\9\2\8\4\6\5\5\6\4\8\2\eeee
\6\4\8\3\9\1\9\3\8\4\7\5\5\7\4\8\eeee
\5\7\4\8\2\9\1\9\2\8\4\6\5\5\6\4\eeee
\5\5\6\4\8\3\9\1\9\3\8\4\7\5\5\7\eeee
\6\5\5\7\4\8\2\9\1\9\2\8\4\6\5\5\eeee
\4\7\5\5\6\4\8\3\9\1\9\3\8\4\7\5\eeee
\8\4\6\5\5\7\4\8\2\9\1\9\2\8\4\6\eeee
\2\8\4\7\5\5\6\4\8\3\9\1\9\3\8\4\eeee
} 

\baaa
9-35
\eaaa
\bbbb
0&0&0&0&0&0&0&0&4\\
0&0&0&0&0&0&0&2&2\\
0&0&0&0&0&0&0&4&0\\
0&0&0&0&0&0&1&1&2\\
0&0&0&0&0&0&2&0&2\\
0&0&0&0&0&0&4&0&0\\
0&0&0&1&2&1&0&0&0\\
0&2&1&1&0&0&0&0&0\\
1&1&0&1&1&0&0&0&0\\
\ebbb
\parbox{7cm}{ 
\1\9\4\9\1\9\4\9\1\9\4\9\1\9\4\9\eeee
\9\2\8\2\9\5\7\5\9\2\8\2\9\5\7\5\eeee
\4\8\3\8\4\7\6\7\4\8\3\8\4\7\6\7\eeee
\9\2\8\2\9\5\7\5\9\2\8\2\9\5\7\5\eeee
\1\9\4\9\1\9\4\9\1\9\4\9\1\9\4\9\eeee
\9\5\7\5\9\2\8\2\9\5\7\5\9\2\8\2\eeee
\4\7\6\7\4\8\3\8\4\7\6\7\4\8\3\8\eeee
\9\5\7\5\9\2\8\2\9\5\7\5\9\2\8\2\eeee
\1\9\4\9\1\9\4\9\1\9\4\9\1\9\4\9\eeee
\9\2\8\2\9\5\7\5\9\2\8\2\9\5\7\5\eeee
\4\8\3\8\4\7\6\7\4\8\3\8\4\7\6\7\eeee
\9\2\8\2\9\5\7\5\9\2\8\2\9\5\7\5\eeee
} 

\baaa
9-36
\eaaa
\bbbb
0&0&0&0&0&0&0&0&4\\
0&0&0&0&0&0&0&2&2\\
0&0&0&0&0&0&0&4&0\\
0&0&0&0&0&0&2&0&2\\
0&0&0&0&0&0&2&1&1\\
0&0&0&0&0&0&2&2&0\\
0&0&0&1&2&1&0&0&0\\
0&1&1&0&1&1&0&0&0\\
1&1&0&1&1&0&0&0&0\\
\ebbb
\parbox{7cm}{ 
\1\9\5\7\4\9\2\8\6\7\5\8\3\8\5\7\eeee
\9\2\8\6\7\5\8\3\8\5\7\6\8\2\9\4\eeee
\5\8\3\8\5\7\6\8\2\9\4\7\5\9\1\9\eeee
\7\6\8\2\9\4\7\5\9\1\9\5\7\4\9\2\eeee
\4\7\5\9\1\9\5\7\4\9\2\8\6\7\5\8\eeee
\9\5\7\4\9\2\8\6\7\5\8\3\8\5\7\6\eeee
\2\8\6\7\5\8\3\8\5\7\6\8\2\9\4\7\eeee
\8\3\8\5\7\6\8\2\9\4\7\5\9\1\9\5\eeee
\6\8\2\9\4\7\5\9\1\9\5\7\4\9\2\8\eeee
\7\5\9\1\9\5\7\4\9\2\8\6\7\5\8\3\eeee
\5\7\4\9\2\8\6\7\5\8\3\8\5\7\6\8\eeee
\8\6\7\5\8\3\8\5\7\6\8\2\9\4\7\5\eeee
} 

\baaa
9-37
\eaaa
\bbbb
0&0&0&0&0&0&0&0&4\\
0&0&0&0&0&0&0&2&2\\
0&0&0&0&0&0&2&0&2\\
0&0&0&0&0&2&0&2&0\\
0&0&0&0&1&0&1&1&1\\
0&0&0&2&0&0&2&0&0\\
0&0&1&0&1&1&0&1&0\\
0&1&0&1&1&0&1&0&0\\
1&1&1&0&1&0&0&0&0\\
\ebbb
\parbox{7cm}{ 
\1\9\5\5\9\1\9\5\5\9\1\9\5\5\9\1\eeee
\9\2\8\7\3\9\2\8\7\3\9\2\8\7\3\9\eeee
\5\8\4\6\7\5\8\4\6\7\5\8\4\6\7\5\eeee
\5\7\6\4\8\5\7\6\4\8\5\7\6\4\8\5\eeee
\9\3\7\8\2\9\3\7\8\2\9\3\7\8\2\9\eeee
\1\9\5\5\9\1\9\5\5\9\1\9\5\5\9\1\eeee
\9\2\8\7\3\9\2\8\7\3\9\2\8\7\3\9\eeee
\5\8\4\6\7\5\8\4\6\7\5\8\4\6\7\5\eeee
\5\7\6\4\8\5\7\6\4\8\5\7\6\4\8\5\eeee
\9\3\7\8\2\9\3\7\8\2\9\3\7\8\2\9\eeee
\1\9\5\5\9\1\9\5\5\9\1\9\5\5\9\1\eeee
\9\2\8\7\3\9\2\8\7\3\9\2\8\7\3\9\eeee
} 

\baaa
9-38
\eaaa
\bbbb
0&0&0&0&0&0&0&0&4\\
0&0&0&0&0&0&0&2&2\\
0&0&0&0&0&0&2&0&2\\
0&0&0&1&0&0&1&1&1\\
0&0&0&0&2&0&0&2&0\\
0&0&0&0&0&2&2&0&0\\
0&0&1&1&0&1&1&0&0\\
0&1&0&1&1&0&0&1&0\\
1&1&1&1&0&0&0&0&0\\
\ebbb
\parbox{7cm}{ 
\1\9\4\4\9\1\9\4\4\9\1\9\4\4\9\1\eeee
\9\2\8\8\2\9\3\7\7\3\9\2\8\8\2\9\eeee
\4\8\5\5\8\4\7\6\6\7\4\8\5\5\8\4\eeee
\4\8\5\5\8\4\7\6\6\7\4\8\5\5\8\4\eeee
\9\2\8\8\2\9\3\7\7\3\9\2\8\8\2\9\eeee
\1\9\4\4\9\1\9\4\4\9\1\9\4\4\9\1\eeee
\9\3\7\7\3\9\2\8\8\2\9\3\7\7\3\9\eeee
\4\7\6\6\7\4\8\5\5\8\4\7\6\6\7\4\eeee
\4\7\6\6\7\4\8\5\5\8\4\7\6\6\7\4\eeee
\9\3\7\7\3\9\2\8\8\2\9\3\7\7\3\9\eeee
\1\9\4\4\9\1\9\4\4\9\1\9\4\4\9\1\eeee
\9\2\8\8\2\9\3\7\7\3\9\2\8\8\2\9\eeee
} 

\baaa
9-39
\eaaa
\bbbb
0&0&0&0&0&0&0&0&4\\
0&0&0&0&0&0&0&2&2\\
0&0&0&0&0&0&2&2&0\\
0&0&0&0&0&0&2&2&0\\
0&0&0&0&0&2&2&0&0\\
0&0&0&0&2&2&0&0&0\\
0&0&1&1&2&0&0&0&0\\
0&2&1&1&0&0&0&0&0\\
2&2&0&0&0&0&0&0&0\\
\ebbb
\parbox{7cm}{ 
\1\9\2\8\3\7\5\6\6\5\7\3\8\2\9\1\eeee
\9\1\9\2\8\4\7\5\6\6\5\7\4\8\2\9\eeee
\2\9\1\9\2\8\3\7\5\6\6\5\7\3\8\2\eeee
\8\2\9\1\9\2\8\4\7\5\6\6\5\7\4\8\eeee
\3\8\2\9\1\9\2\8\3\7\5\6\6\5\7\3\eeee
\7\4\8\2\9\1\9\2\8\4\7\5\6\6\5\7\eeee
\5\7\3\8\2\9\1\9\2\8\3\7\5\6\6\5\eeee
\6\5\7\4\8\2\9\1\9\2\8\4\7\5\6\6\eeee
\6\6\5\7\3\8\2\9\1\9\2\8\3\7\5\6\eeee
\5\6\6\5\7\4\8\2\9\1\9\2\8\4\7\5\eeee
\7\5\6\6\5\7\3\8\2\9\1\9\2\8\3\7\eeee
\3\7\5\6\6\5\7\4\8\2\9\1\9\2\8\4\eeee
} 

\baaa
9-40
\eaaa
\bbbb
0&0&0&0&0&0&0&0&4\\
0&0&0&0&0&0&0&2&2\\
0&0&0&0&0&0&2&2&0\\
0&0&0&0&0&0&2&2&0\\
0&0&0&0&0&2&2&0&0\\
0&0&0&0&4&0&0&0&0\\
0&0&1&1&2&0&0&0&0\\
0&2&1&1&0&0&0&0&0\\
2&2&0&0&0&0&0&0&0\\
\ebbb
\parbox{7cm}{ 
\1\9\2\8\3\7\5\6\5\7\3\8\2\9\1\9\eeee
\9\1\9\2\8\4\7\5\6\5\7\4\8\2\9\1\eeee
\2\9\1\9\2\8\3\7\5\6\5\7\3\8\2\9\eeee
\8\2\9\1\9\2\8\4\7\5\6\5\7\4\8\2\eeee
\3\8\2\9\1\9\2\8\3\7\5\6\5\7\3\8\eeee
\7\4\8\2\9\1\9\2\8\4\7\5\6\5\7\4\eeee
\5\7\3\8\2\9\1\9\2\8\3\7\5\6\5\7\eeee
\6\5\7\4\8\2\9\1\9\2\8\4\7\5\6\5\eeee
\5\6\5\7\3\8\2\9\1\9\2\8\3\7\5\6\eeee
\7\5\6\5\7\4\8\2\9\1\9\2\8\4\7\5\eeee
\3\7\5\6\5\7\3\8\2\9\1\9\2\8\3\7\eeee
\8\4\7\5\6\5\7\4\8\2\9\1\9\2\8\4\eeee
} 

\baaa
9-41
\eaaa
\bbbb
0&0&0&0&0&0&0&0&4\\
0&0&0&0&0&0&0&2&2\\
0&0&0&0&0&0&2&2&0\\
0&0&0&0&0&0&2&2&0\\
0&0&0&0&1&1&2&0&0\\
0&0&0&0&1&1&2&0&0\\
0&0&1&1&1&1&0&0&0\\
0&2&1&1&0&0&0&0&0\\
2&2&0&0&0&0&0&0&0\\
\ebbb
\parbox{7cm}{ 
\1\9\2\8\3\7\5\5\7\3\8\2\9\1\9\2\eeee
\9\1\9\2\8\4\7\6\6\7\4\8\2\9\1\9\eeee
\2\9\1\9\2\8\3\7\5\5\7\3\8\2\9\1\eeee
\8\2\9\1\9\2\8\4\7\6\6\7\4\8\2\9\eeee
\3\8\2\9\1\9\2\8\3\7\5\5\7\3\8\2\eeee
\7\4\8\2\9\1\9\2\8\4\7\6\6\7\4\8\eeee
\5\7\3\8\2\9\1\9\2\8\3\7\5\5\7\3\eeee
\5\6\7\4\8\2\9\1\9\2\8\4\7\6\6\7\eeee
\7\6\5\7\3\8\2\9\1\9\2\8\3\7\5\5\eeee
\3\7\5\6\7\4\8\2\9\1\9\2\8\4\7\6\eeee
\8\4\7\6\5\7\3\8\2\9\1\9\2\8\3\7\eeee
\2\8\3\7\5\6\7\4\8\2\9\1\9\2\8\4\eeee
} 

\baaa
9-42
\eaaa
\bbbb
0&0&0&0&0&0&0&0&4\\
0&0&0&0&0&0&0&2&2\\
0&0&0&0&0&0&2&2&0\\
0&0&0&0&0&2&2&0&0\\
0&0&0&0&0&2&2&0&0\\
0&0&0&1&1&2&0&0&0\\
0&0&2&1&1&0&0&0&0\\
0&2&2&0&0&0&0&0&0\\
2&2&0&0&0&0&0&0&0\\
\ebbb
\parbox{7cm}{ 
\1\9\2\8\3\7\4\6\6\4\7\3\8\2\9\1\eeee
\9\1\9\2\8\3\7\5\6\6\5\7\3\8\2\9\eeee
\2\9\1\9\2\8\3\7\4\6\6\4\7\3\8\2\eeee
\8\2\9\1\9\2\8\3\7\5\6\6\5\7\3\8\eeee
\3\8\2\9\1\9\2\8\3\7\4\6\6\4\7\3\eeee
\7\3\8\2\9\1\9\2\8\3\7\5\6\6\5\7\eeee
\4\7\3\8\2\9\1\9\2\8\3\7\4\6\6\4\eeee
\6\5\7\3\8\2\9\1\9\2\8\3\7\5\6\6\eeee
\6\6\4\7\3\8\2\9\1\9\2\8\3\7\4\6\eeee
\4\6\6\5\7\3\8\2\9\1\9\2\8\3\7\5\eeee
\7\5\6\6\4\7\3\8\2\9\1\9\2\8\3\7\eeee
\3\7\4\6\6\5\7\3\8\2\9\1\9\2\8\3\eeee
} 

\baaa
9-43
\eaaa
\bbbb
0&0&0&0&0&0&0&0&4\\
0&0&0&0&0&0&0&2&2\\
0&0&0&0&0&0&2&2&0\\
0&0&0&0&0&2&2&0&0\\
0&0&0&0&0&2&2&0&0\\
0&0&0&2&2&0&0&0&0\\
0&0&2&1&1&0&0&0&0\\
0&2&2&0&0&0&0&0&0\\
2&2&0&0&0&0&0&0&0\\
\ebbb
\parbox{7cm}{ 
\1\9\2\8\3\7\4\6\4\7\3\8\2\9\1\9\eeee
\9\1\9\2\8\3\7\5\6\5\7\3\8\2\9\1\eeee
\2\9\1\9\2\8\3\7\4\6\4\7\3\8\2\9\eeee
\8\2\9\1\9\2\8\3\7\5\6\5\7\3\8\2\eeee
\3\8\2\9\1\9\2\8\3\7\4\6\4\7\3\8\eeee
\7\3\8\2\9\1\9\2\8\3\7\5\6\5\7\3\eeee
\4\7\3\8\2\9\1\9\2\8\3\7\4\6\4\7\eeee
\6\5\7\3\8\2\9\1\9\2\8\3\7\5\6\5\eeee
\4\6\4\7\3\8\2\9\1\9\2\8\3\7\4\6\eeee
\7\5\6\5\7\3\8\2\9\1\9\2\8\3\7\5\eeee
\3\7\4\6\4\7\3\8\2\9\1\9\2\8\3\7\eeee
\8\3\7\5\6\5\7\3\8\2\9\1\9\2\8\3\eeee
} 

\baaa
9-44
\eaaa
\bbbb
0&0&0&0&0&0&0&0&4\\
0&0&0&0&0&0&0&2&2\\
0&0&0&0&0&0&2&2&0\\
0&0&0&0&0&2&2&0&0\\
0&0&0&0&0&4&0&0&0\\
0&0&0&2&2&0&0&0&0\\
0&0&2&2&0&0&0&0&0\\
0&2&2&0&0&0&0&0&0\\
2&2&0&0&0&0&0&0&0\\
\ebbb
\parbox{7cm}{ 
\1\9\2\8\3\7\4\6\5\6\4\7\3\8\2\9\eeee
\9\1\9\2\8\3\7\4\6\5\6\4\7\3\8\2\eeee
\2\9\1\9\2\8\3\7\4\6\5\6\4\7\3\8\eeee
\8\2\9\1\9\2\8\3\7\4\6\5\6\4\7\3\eeee
\3\8\2\9\1\9\2\8\3\7\4\6\5\6\4\7\eeee
\7\3\8\2\9\1\9\2\8\3\7\4\6\5\6\4\eeee
\4\7\3\8\2\9\1\9\2\8\3\7\4\6\5\6\eeee
\6\4\7\3\8\2\9\1\9\2\8\3\7\4\6\5\eeee
\5\6\4\7\3\8\2\9\1\9\2\8\3\7\4\6\eeee
\6\5\6\4\7\3\8\2\9\1\9\2\8\3\7\4\eeee
\4\6\5\6\4\7\3\8\2\9\1\9\2\8\3\7\eeee
\7\4\6\5\6\4\7\3\8\2\9\1\9\2\8\3\eeee
} 

\baaa
9-45
\eaaa
\bbbb
0&0&0&0&0&0&0&0&4\\
0&0&0&0&0&0&0&2&2\\
0&0&0&0&0&0&2&2&0\\
0&0&0&0&0&2&2&0&0\\
0&0&0&0&2&2&0&0&0\\
0&0&0&2&2&0&0&0&0\\
0&0&2&2&0&0&0&0&0\\
0&2&2&0&0&0&0&0&0\\
2&2&0&0&0&0&0&0&0\\
\ebbb
\parbox{7cm}{ 
\1\9\2\8\3\7\4\6\5\5\6\4\7\3\8\2\eeee
\9\1\9\2\8\3\7\4\6\5\5\6\4\7\3\8\eeee
\2\9\1\9\2\8\3\7\4\6\5\5\6\4\7\3\eeee
\8\2\9\1\9\2\8\3\7\4\6\5\5\6\4\7\eeee
\3\8\2\9\1\9\2\8\3\7\4\6\5\5\6\4\eeee
\7\3\8\2\9\1\9\2\8\3\7\4\6\5\5\6\eeee
\4\7\3\8\2\9\1\9\2\8\3\7\4\6\5\5\eeee
\6\4\7\3\8\2\9\1\9\2\8\3\7\4\6\5\eeee
\5\6\4\7\3\8\2\9\1\9\2\8\3\7\4\6\eeee
\5\5\6\4\7\3\8\2\9\1\9\2\8\3\7\4\eeee
\6\5\5\6\4\7\3\8\2\9\1\9\2\8\3\7\eeee
\4\6\5\5\6\4\7\3\8\2\9\1\9\2\8\3\eeee
} 

\baaa
9-46
\eaaa
\bbbb
0&0&0&0&0&0&0&0&4\\
0&0&0&0&0&0&0&2&2\\
0&0&0&0&0&0&2&2&0\\
0&0&0&0&1&1&2&0&0\\
0&0&0&2&1&1&0&0&0\\
0&0&0&2&1&1&0&0&0\\
0&0&2&2&0&0&0&0&0\\
0&2&2&0&0&0&0&0&0\\
2&2&0&0&0&0&0&0&0\\
\ebbb
\parbox{7cm}{ 
\1\9\2\8\3\7\4\5\5\4\7\3\8\2\9\1\eeee
\9\1\9\2\8\3\7\4\6\6\4\7\3\8\2\9\eeee
\2\9\1\9\2\8\3\7\4\5\5\4\7\3\8\2\eeee
\8\2\9\1\9\2\8\3\7\4\6\6\4\7\3\8\eeee
\3\8\2\9\1\9\2\8\3\7\4\5\5\4\7\3\eeee
\7\3\8\2\9\1\9\2\8\3\7\4\6\6\4\7\eeee
\4\7\3\8\2\9\1\9\2\8\3\7\4\5\5\4\eeee
\5\4\7\3\8\2\9\1\9\2\8\3\7\4\6\6\eeee
\5\6\4\7\3\8\2\9\1\9\2\8\3\7\4\5\eeee
\4\6\5\4\7\3\8\2\9\1\9\2\8\3\7\4\eeee
\7\4\5\6\4\7\3\8\2\9\1\9\2\8\3\7\eeee
\3\7\4\6\5\4\7\3\8\2\9\1\9\2\8\3\eeee
} 

\baaa
9-47
\eaaa
\bbbb
0&0&0&0&0&0&0&0&4\\
0&0&0&0&0&0&0&2&2\\
0&0&0&0&0&1&1&2&0\\
0&0&0&0&0&1&1&2&0\\
0&0&0&0&0&2&2&0&0\\
0&0&1&1&2&0&0&0&0\\
0&0&1&1&2&0&0&0&0\\
0&2&1&1&0&0&0&0&0\\
2&2&0&0&0&0&0&0&0\\
\ebbb
\parbox{7cm}{ 
\1\9\2\8\3\6\5\6\3\8\2\9\1\9\2\8\eeee
\9\1\9\2\8\4\7\5\7\4\8\2\9\1\9\2\eeee
\2\9\1\9\2\8\3\6\5\6\3\8\2\9\1\9\eeee
\8\2\9\1\9\2\8\4\7\5\7\4\8\2\9\1\eeee
\3\8\2\9\1\9\2\8\3\6\5\6\3\8\2\9\eeee
\6\4\8\2\9\1\9\2\8\4\7\5\7\4\8\2\eeee
\5\7\3\8\2\9\1\9\2\8\3\6\5\6\3\8\eeee
\6\5\6\4\8\2\9\1\9\2\8\4\7\5\7\4\eeee
\3\7\5\7\3\8\2\9\1\9\2\8\3\6\5\6\eeee
\8\4\6\5\6\4\8\2\9\1\9\2\8\4\7\5\eeee
\2\8\3\7\5\7\3\8\2\9\1\9\2\8\3\6\eeee
\9\2\8\4\6\5\6\4\8\2\9\1\9\2\8\4\eeee
} 

\baaa
9-48
\eaaa
\bbbb
0&0&0&0&0&0&0&0&4\\
0&0&0&0&0&0&0&2&2\\
0&0&0&0&0&1&1&2&0\\
0&0&0&0&0&1&1&2&0\\
0&0&0&0&2&1&1&0&0\\
0&0&1&1&2&0&0&0&0\\
0&0&1&1&2&0&0&0&0\\
0&2&1&1&0&0&0&0&0\\
2&2&0&0&0&0&0&0&0\\
\ebbb
\parbox{7cm}{ 
\1\9\2\8\3\6\5\5\6\3\8\2\9\1\9\2\eeee
\9\1\9\2\8\4\7\5\5\7\4\8\2\9\1\9\eeee
\2\9\1\9\2\8\3\6\5\5\6\3\8\2\9\1\eeee
\8\2\9\1\9\2\8\4\7\5\5\7\4\8\2\9\eeee
\3\8\2\9\1\9\2\8\3\6\5\5\6\3\8\2\eeee
\6\4\8\2\9\1\9\2\8\4\7\5\5\7\4\8\eeee
\5\7\3\8\2\9\1\9\2\8\3\6\5\5\6\3\eeee
\5\5\6\4\8\2\9\1\9\2\8\4\7\5\5\7\eeee
\6\5\5\7\3\8\2\9\1\9\2\8\3\6\5\5\eeee
\3\7\5\5\6\4\8\2\9\1\9\2\8\4\7\5\eeee
\8\4\6\5\5\7\3\8\2\9\1\9\2\8\3\6\eeee
\2\8\3\7\5\5\6\4\8\2\9\1\9\2\8\4\eeee
} 

\baaa
9-49
\eaaa
\bbbb
0&0&0&0&0&0&0&0&4\\
0&0&0&0&0&0&0&2&2\\
0&0&0&0&0&1&1&2&0\\
0&0&0&0&0&2&2&0&0\\
0&0&0&0&0&2&2&0&0\\
0&0&2&1&1&0&0&0&0\\
0&0&2&1&1&0&0&0&0\\
0&2&2&0&0&0&0&0&0\\
2&2&0&0&0&0&0&0&0\\
\ebbb
\parbox{7cm}{ 
\1\9\2\8\3\6\4\6\3\8\2\9\1\9\2\8\eeee
\9\1\9\2\8\3\7\5\7\3\8\2\9\1\9\2\eeee
\2\9\1\9\2\8\3\6\4\6\3\8\2\9\1\9\eeee
\8\2\9\1\9\2\8\3\7\5\7\3\8\2\9\1\eeee
\3\8\2\9\1\9\2\8\3\6\4\6\3\8\2\9\eeee
\6\3\8\2\9\1\9\2\8\3\7\5\7\3\8\2\eeee
\4\7\3\8\2\9\1\9\2\8\3\6\4\6\3\8\eeee
\6\5\6\3\8\2\9\1\9\2\8\3\7\5\7\3\eeee
\3\7\4\7\3\8\2\9\1\9\2\8\3\6\4\6\eeee
\8\3\6\5\6\3\8\2\9\1\9\2\8\3\7\5\eeee
\2\8\3\7\4\7\3\8\2\9\1\9\2\8\3\6\eeee
\9\2\8\3\6\5\6\3\8\2\9\1\9\2\8\3\eeee
} 

\baaa
9-50
\eaaa
\bbbb
0&0&0&0&0&0&0&0&4\\
0&0&0&0&0&0&0&2&2\\
0&0&0&0&0&1&1&2&0\\
0&0&0&0&2&1&1&0&0\\
0&0&0&2&2&0&0&0&0\\
0&0&2&2&0&0&0&0&0\\
0&0&2&2&0&0&0&0&0\\
0&2&2&0&0&0&0&0&0\\
2&2&0&0&0&0&0&0&0\\
\ebbb
\parbox{7cm}{ 
\1\9\2\8\3\6\4\5\5\4\6\3\8\2\9\1\eeee
\9\1\9\2\8\3\7\4\5\5\4\7\3\8\2\9\eeee
\2\9\1\9\2\8\3\6\4\5\5\4\6\3\8\2\eeee
\8\2\9\1\9\2\8\3\7\4\5\5\4\7\3\8\eeee
\3\8\2\9\1\9\2\8\3\6\4\5\5\4\6\3\eeee
\6\3\8\2\9\1\9\2\8\3\7\4\5\5\4\7\eeee
\4\7\3\8\2\9\1\9\2\8\3\6\4\5\5\4\eeee
\5\4\6\3\8\2\9\1\9\2\8\3\7\4\5\5\eeee
\5\5\4\7\3\8\2\9\1\9\2\8\3\6\4\5\eeee
\4\5\5\4\6\3\8\2\9\1\9\2\8\3\7\4\eeee
\6\4\5\5\4\7\3\8\2\9\1\9\2\8\3\6\eeee
\3\7\4\5\5\4\6\3\8\2\9\1\9\2\8\3\eeee
} 

\baaa
9-51
\eaaa
\bbbb
0&0&0&0&0&0&0&0&4\\
0&0&0&0&0&0&0&2&2\\
0&0&0&0&0&1&1&2&0\\
0&0&0&1&1&1&1&0&0\\
0&0&0&1&1&1&1&0&0\\
0&0&2&1&1&0&0&0&0\\
0&0&2&1&1&0&0&0&0\\
0&2&2&0&0&0&0&0&0\\
2&2&0&0&0&0&0&0&0\\
\ebbb
\parbox{7cm}{ 
\1\9\2\8\3\6\4\4\6\3\8\2\9\1\9\2\eeee
\9\1\9\2\8\3\7\5\5\7\3\8\2\9\1\9\eeee
\2\9\1\9\2\8\3\6\4\4\6\3\8\2\9\1\eeee
\8\2\9\1\9\2\8\3\7\5\5\7\3\8\2\9\eeee
\3\8\2\9\1\9\2\8\3\6\4\4\6\3\8\2\eeee
\6\3\8\2\9\1\9\2\8\3\7\5\5\7\3\8\eeee
\4\7\3\8\2\9\1\9\2\8\3\6\4\4\6\3\eeee
\4\5\6\3\8\2\9\1\9\2\8\3\7\5\5\7\eeee
\6\5\4\7\3\8\2\9\1\9\2\8\3\6\4\4\eeee
\3\7\4\5\6\3\8\2\9\1\9\2\8\3\7\5\eeee
\8\3\6\5\4\7\3\8\2\9\1\9\2\8\3\6\eeee
\2\8\3\7\4\5\6\3\8\2\9\1\9\2\8\3\eeee
} 

\baaa
9-52
\eaaa
\bbbb
0&0&0&0&0&0&0&0&4\\
0&0&0&0&0&0&0&4&0\\
0&0&0&0&0&1&1&0&2\\
0&0&0&0&0&1&1&0&2\\
0&0&0&0&0&1&1&2&0\\
0&0&1&1&2&0&0&0&0\\
0&0&1&1&2&0&0&0&0\\
0&2&0&0&2&0&0&0&0\\
2&0&1&1&0&0&0&0&0\\
\ebbb
\parbox{7cm}{ 
\1\9\3\6\5\8\2\8\5\6\3\9\1\9\3\6\eeee
\9\1\9\4\7\5\8\2\8\5\7\4\9\1\9\4\eeee
\3\9\1\9\3\6\5\8\2\8\5\6\3\9\1\9\eeee
\6\4\9\1\9\4\7\5\8\2\8\5\7\4\9\1\eeee
\5\7\3\9\1\9\3\6\5\8\2\8\5\6\3\9\eeee
\8\5\6\4\9\1\9\4\7\5\8\2\8\5\7\4\eeee
\2\8\5\7\3\9\1\9\3\6\5\8\2\8\5\6\eeee
\8\2\8\5\6\4\9\1\9\4\7\5\8\2\8\5\eeee
\5\8\2\8\5\7\3\9\1\9\3\6\5\8\2\8\eeee
\6\5\8\2\8\5\6\4\9\1\9\4\7\5\8\2\eeee
\3\7\5\8\2\8\5\7\3\9\1\9\3\6\5\8\eeee
\9\4\6\5\8\2\8\5\6\4\9\1\9\4\7\5\eeee
} 

\baaa
9-53
\eaaa
\bbbb
0&0&0&0&0&0&0&0&4\\
0&0&0&0&0&0&1&1&2\\
0&0&0&0&0&0&1&1&2\\
0&0&0&0&0&2&1&1&0\\
0&0&0&0&0&2&1&1&0\\
0&0&0&1&1&2&0&0&0\\
0&1&1&1&1&0&0&0&0\\
0&1&1&1&1&0&0&0&0\\
2&1&1&0&0&0&0&0&0\\
\ebbb
\parbox{7cm}{ 
\1\9\2\7\4\6\6\4\7\2\9\1\9\2\7\4\eeee
\9\1\9\3\8\5\6\6\5\8\3\9\1\9\3\8\eeee
\2\9\1\9\2\7\4\6\6\4\7\2\9\1\9\2\eeee
\7\3\9\1\9\3\8\5\6\6\5\8\3\9\1\9\eeee
\4\8\2\9\1\9\2\7\4\6\6\4\7\2\9\1\eeee
\6\5\7\3\9\1\9\3\8\5\6\6\5\8\3\9\eeee
\6\6\4\8\2\9\1\9\2\7\4\6\6\4\7\2\eeee
\4\6\6\5\7\3\9\1\9\3\8\5\6\6\5\8\eeee
\7\5\6\6\4\8\2\9\1\9\2\7\4\6\6\4\eeee
\2\8\4\6\6\5\7\3\9\1\9\3\8\5\6\6\eeee
\9\3\7\5\6\6\4\8\2\9\1\9\2\7\4\6\eeee
\1\9\2\8\4\6\6\5\7\3\9\1\9\3\8\5\eeee
} 

\baaa
9-54
\eaaa
\bbbb
0&0&0&0&0&0&0&0&4\\
0&0&0&0&0&0&1&1&2\\
0&0&0&0&0&0&1&1&2\\
0&0&0&0&0&2&1&1&0\\
0&0&0&0&0&2&1&1&0\\
0&0&0&2&2&0&0&0&0\\
0&1&1&1&1&0&0&0&0\\
0&1&1&1&1&0&0&0&0\\
2&1&1&0&0&0&0&0&0\\
\ebbb
\parbox{7cm}{ 
\1\9\2\7\4\6\4\7\2\9\1\9\2\7\4\6\eeee
\9\1\9\3\8\5\6\5\8\3\9\1\9\3\8\5\eeee
\2\9\1\9\2\7\4\6\4\7\2\9\1\9\2\7\eeee
\7\3\9\1\9\3\8\5\6\5\8\3\9\1\9\3\eeee
\4\8\2\9\1\9\2\7\4\6\4\7\2\9\1\9\eeee
\6\5\7\3\9\1\9\3\8\5\6\5\8\3\9\1\eeee
\4\6\4\8\2\9\1\9\2\7\4\6\4\7\2\9\eeee
\7\5\6\5\7\3\9\1\9\3\8\5\6\5\8\3\eeee
\2\8\4\6\4\8\2\9\1\9\2\7\4\6\4\7\eeee
\9\3\7\5\6\5\7\3\9\1\9\3\8\5\6\5\eeee
\1\9\2\8\4\6\4\8\2\9\1\9\2\7\4\6\eeee
\9\1\9\3\7\5\6\5\7\3\9\1\9\3\8\5\eeee
} 

\baaa
9-55
\eaaa
\bbbb
0&0&0&0&0&0&0&0&4\\
0&0&0&0&0&0&1&1&2\\
0&0&0&0&0&0&1&1&2\\
0&0&0&0&0&2&1&1&0\\
0&0&0&0&2&2&0&0&0\\
0&0&0&2&2&0&0&0&0\\
0&1&1&2&0&0&0&0&0\\
0&1&1&2&0&0&0&0&0\\
2&1&1&0&0&0&0&0&0\\
\ebbb
\parbox{7cm}{ 
\1\9\2\7\4\6\5\5\6\4\7\2\9\1\9\2\eeee
\9\1\9\3\8\4\6\5\5\6\4\8\3\9\1\9\eeee
\2\9\1\9\2\7\4\6\5\5\6\4\7\2\9\1\eeee
\7\3\9\1\9\3\8\4\6\5\5\6\4\8\3\9\eeee
\4\8\2\9\1\9\2\7\4\6\5\5\6\4\7\2\eeee
\6\4\7\3\9\1\9\3\8\4\6\5\5\6\4\8\eeee
\5\6\4\8\2\9\1\9\2\7\4\6\5\5\6\4\eeee
\5\5\6\4\7\3\9\1\9\3\8\4\6\5\5\6\eeee
\6\5\5\6\4\8\2\9\1\9\2\7\4\6\5\5\eeee
\4\6\5\5\6\4\7\3\9\1\9\3\8\4\6\5\eeee
\7\4\6\5\5\6\4\8\2\9\1\9\2\7\4\6\eeee
\2\8\4\6\5\5\6\4\7\3\9\1\9\3\8\4\eeee
} 

\baaa
9-56
\eaaa
\bbbb
0&0&0&0&0&0&0&0&4\\
0&0&0&0&0&0&1&1&2\\
0&0&0&0&0&0&1&1&2\\
0&0&0&0&1&1&1&1&0\\
0&0&0&2&1&1&0&0&0\\
0&0&0&2&1&1&0&0&0\\
0&1&1&2&0&0&0&0&0\\
0&1&1&2&0&0&0&0&0\\
2&1&1&0&0&0&0&0&0\\
\ebbb
\parbox{7cm}{ 
\1\9\2\7\4\5\5\4\7\2\9\1\9\2\7\4\eeee
\9\1\9\3\8\4\6\6\4\8\3\9\1\9\3\8\eeee
\2\9\1\9\2\7\4\5\5\4\7\2\9\1\9\2\eeee
\7\3\9\1\9\3\8\4\6\6\4\8\3\9\1\9\eeee
\4\8\2\9\1\9\2\7\4\5\5\4\7\2\9\1\eeee
\5\4\7\3\9\1\9\3\8\4\6\6\4\8\3\9\eeee
\5\6\4\8\2\9\1\9\2\7\4\5\5\4\7\2\eeee
\4\6\5\4\7\3\9\1\9\3\8\4\6\6\4\8\eeee
\7\4\5\6\4\8\2\9\1\9\2\7\4\5\5\4\eeee
\2\8\4\6\5\4\7\3\9\1\9\3\8\4\6\6\eeee
\9\3\7\4\5\6\4\8\2\9\1\9\2\7\4\5\eeee
\1\9\2\8\4\6\5\4\7\3\9\1\9\3\8\4\eeee
} 

\baaa
9-57
\eaaa
\bbbb
0&0&0&0&0&0&0&0&4\\
0&0&0&0&0&0&1&1&2\\
0&0&0&0&0&2&1&1&0\\
0&0&0&0&0&2&1&1&0\\
0&0&0&0&2&2&0&0&0\\
0&0&1&1&2&0&0&0&0\\
0&2&1&1&0&0&0&0&0\\
0&2&1&1&0&0&0&0&0\\
2&2&0&0&0&0&0&0&0\\
\ebbb
\parbox{7cm}{ 
\1\9\2\7\3\6\5\5\6\3\7\2\9\1\9\2\eeee
\9\1\9\2\8\4\6\5\5\6\4\8\2\9\1\9\eeee
\2\9\1\9\2\7\3\6\5\5\6\3\7\2\9\1\eeee
\7\2\9\1\9\2\8\4\6\5\5\6\4\8\2\9\eeee
\3\8\2\9\1\9\2\7\3\6\5\5\6\3\7\2\eeee
\6\4\7\2\9\1\9\2\8\4\6\5\5\6\4\8\eeee
\5\6\3\8\2\9\1\9\2\7\3\6\5\5\6\3\eeee
\5\5\6\4\7\2\9\1\9\2\8\4\6\5\5\6\eeee
\6\5\5\6\3\8\2\9\1\9\2\7\3\6\5\5\eeee
\3\6\5\5\6\4\7\2\9\1\9\2\8\4\6\5\eeee
\7\4\6\5\5\6\3\8\2\9\1\9\2\7\3\6\eeee
\2\8\3\6\5\5\6\4\7\2\9\1\9\2\8\4\eeee
} 

\baaa
9-58
\eaaa
\bbbb
0&0&0&0&0&0&0&0&4\\
0&0&0&0&0&0&1&1&2\\
0&0&0&0&0&2&1&1&0\\
0&0&0&0&2&2&0&0&0\\
0&0&0&2&2&0&0&0&0\\
0&0&2&2&0&0&0&0&0\\
0&2&2&0&0&0&0&0&0\\
0&2&2&0&0&0&0&0&0\\
2&2&0&0&0&0&0&0&0\\
\ebbb
\parbox{7cm}{ 
\1\9\2\7\3\6\4\5\5\4\6\3\7\2\9\1\eeee
\9\1\9\2\8\3\6\4\5\5\4\6\3\8\2\9\eeee
\2\9\1\9\2\7\3\6\4\5\5\4\6\3\7\2\eeee
\7\2\9\1\9\2\8\3\6\4\5\5\4\6\3\8\eeee
\3\8\2\9\1\9\2\7\3\6\4\5\5\4\6\3\eeee
\6\3\7\2\9\1\9\2\8\3\6\4\5\5\4\6\eeee
\4\6\3\8\2\9\1\9\2\7\3\6\4\5\5\4\eeee
\5\4\6\3\7\2\9\1\9\2\8\3\6\4\5\5\eeee
\5\5\4\6\3\8\2\9\1\9\2\7\3\6\4\5\eeee
\4\5\5\4\6\3\7\2\9\1\9\2\8\3\6\4\eeee
\6\4\5\5\4\6\3\8\2\9\1\9\2\7\3\6\eeee
\3\6\4\5\5\4\6\3\7\2\9\1\9\2\8\3\eeee
} 

\baaa
9-59
\eaaa
\bbbb
0&0&0&0&0&0&0&0&4\\
0&0&0&0&0&0&1&1&2\\
0&0&0&0&0&2&1&1&0\\
0&0&0&1&1&2&0&0&0\\
0&0&0&1&1&2&0&0&0\\
0&0&2&1&1&0&0&0&0\\
0&2&2&0&0&0&0&0&0\\
0&2&2&0&0&0&0&0&0\\
2&2&0&0&0&0&0&0&0\\
\ebbb
\parbox{7cm}{ 
\1\9\2\7\3\6\4\4\6\3\7\2\9\1\9\2\eeee
\9\1\9\2\8\3\6\5\5\6\3\8\2\9\1\9\eeee
\2\9\1\9\2\7\3\6\4\4\6\3\7\2\9\1\eeee
\7\2\9\1\9\2\8\3\6\5\5\6\3\8\2\9\eeee
\3\8\2\9\1\9\2\7\3\6\4\4\6\3\7\2\eeee
\6\3\7\2\9\1\9\2\8\3\6\5\5\6\3\8\eeee
\4\6\3\8\2\9\1\9\2\7\3\6\4\4\6\3\eeee
\4\5\6\3\7\2\9\1\9\2\8\3\6\5\5\6\eeee
\6\5\4\6\3\8\2\9\1\9\2\7\3\6\4\4\eeee
\3\6\4\5\6\3\7\2\9\1\9\2\8\3\6\5\eeee
\7\3\6\5\4\6\3\8\2\9\1\9\2\7\3\6\eeee
\2\8\3\6\4\5\6\3\7\2\9\1\9\2\8\3\eeee
} 

\baaa
9-60
\eaaa
\bbbb
0&0&0&0&0&0&0&0&4\\
0&0&0&0&0&0&1&1&2\\
0&0&0&0&1&1&1&1&0\\
0&0&0&0&1&1&1&1&0\\
0&0&1&1&1&1&0&0&0\\
0&0&1&1&1&1&0&0&0\\
0&2&1&1&0&0&0&0&0\\
0&2&1&1&0&0&0&0&0\\
2&2&0&0&0&0&0&0&0\\
\ebbb
\parbox{7cm}{ 
\1\9\2\7\3\5\5\3\7\2\9\1\9\2\7\3\eeee
\9\1\9\2\8\4\6\6\4\8\2\9\1\9\2\8\eeee
\2\9\1\9\2\7\3\5\5\3\7\2\9\1\9\2\eeee
\7\2\9\1\9\2\8\4\6\6\4\8\2\9\1\9\eeee
\3\8\2\9\1\9\2\7\3\5\5\3\7\2\9\1\eeee
\5\4\7\2\9\1\9\2\8\4\6\6\4\8\2\9\eeee
\5\6\3\8\2\9\1\9\2\7\3\5\5\3\7\2\eeee
\3\6\5\4\7\2\9\1\9\2\8\4\6\6\4\8\eeee
\7\4\5\6\3\8\2\9\1\9\2\7\3\5\5\3\eeee
\2\8\3\6\5\4\7\2\9\1\9\2\8\4\6\6\eeee
\9\2\7\4\5\6\3\8\2\9\1\9\2\7\3\5\eeee
\1\9\2\8\3\6\5\4\7\2\9\1\9\2\8\4\eeee
} 

\baaa
9-61
\eaaa
\bbbb
0&0&0&0&0&0&0&0&4\\
0&0&0&0&0&0&1&1&2\\
0&0&0&0&1&1&1&1&0\\
0&0&0&0&1&1&1&1&0\\
0&0&2&2&0&0&0&0&0\\
0&0&2&2&0&0&0&0&0\\
0&2&1&1&0&0&0&0&0\\
0&2&1&1&0&0&0&0&0\\
2&2&0&0&0&0&0&0&0\\
\ebbb
\parbox{7cm}{ 
\1\9\2\7\3\5\3\7\2\9\1\9\2\7\3\5\eeee
\9\1\9\2\8\4\6\4\8\2\9\1\9\2\8\4\eeee
\2\9\1\9\2\7\3\5\3\7\2\9\1\9\2\7\eeee
\7\2\9\1\9\2\8\4\6\4\8\2\9\1\9\2\eeee
\3\8\2\9\1\9\2\7\3\5\3\7\2\9\1\9\eeee
\5\4\7\2\9\1\9\2\8\4\6\4\8\2\9\1\eeee
\3\6\3\8\2\9\1\9\2\7\3\5\3\7\2\9\eeee
\7\4\5\4\7\2\9\1\9\2\8\4\6\4\8\2\eeee
\2\8\3\6\3\8\2\9\1\9\2\7\3\5\3\7\eeee
\9\2\7\4\5\4\7\2\9\1\9\2\8\4\6\4\eeee
\1\9\2\8\3\6\3\8\2\9\1\9\2\7\3\5\eeee
\9\1\9\2\7\4\5\4\7\2\9\1\9\2\8\4\eeee
} 

\baaa
9-62
\eaaa
\bbbb
0&0&0&0&0&0&0&0&4\\
0&0&0&0&0&0&1&1&2\\
0&0&0&0&1&1&1&1&0\\
0&0&0&0&2&2&0&0&0\\
0&0&2&2&0&0&0&0&0\\
0&0&2&2&0&0&0&0&0\\
0&2&2&0&0&0&0&0&0\\
0&2&2&0&0&0&0&0&0\\
2&2&0&0&0&0&0&0&0\\
\ebbb
\parbox{7cm}{ 
\1\9\2\7\3\5\4\5\3\7\2\9\1\9\2\7\eeee
\9\1\9\2\8\3\6\4\6\3\8\2\9\1\9\2\eeee
\2\9\1\9\2\7\3\5\4\5\3\7\2\9\1\9\eeee
\7\2\9\1\9\2\8\3\6\4\6\3\8\2\9\1\eeee
\3\8\2\9\1\9\2\7\3\5\4\5\3\7\2\9\eeee
\5\3\7\2\9\1\9\2\8\3\6\4\6\3\8\2\eeee
\4\6\3\8\2\9\1\9\2\7\3\5\4\5\3\7\eeee
\5\4\5\3\7\2\9\1\9\2\8\3\6\4\6\3\eeee
\3\6\4\6\3\8\2\9\1\9\2\7\3\5\4\5\eeee
\7\3\5\4\5\3\7\2\9\1\9\2\8\3\6\4\eeee
\2\8\3\6\4\6\3\8\2\9\1\9\2\7\3\5\eeee
\9\2\7\3\5\4\5\3\7\2\9\1\9\2\8\3\eeee
} 

\baaa
9-63
\eaaa
\bbbb
0&0&0&0&0&0&0&0&4\\
0&0&0&0&0&0&1&1&2\\
0&0&0&0&1&1&1&1&0\\
0&0&0&2&1&1&0&0&0\\
0&0&2&2&0&0&0&0&0\\
0&0&2&2&0&0&0&0&0\\
0&2&2&0&0&0&0&0&0\\
0&2&2&0&0&0&0&0&0\\
2&2&0&0&0&0&0&0&0\\
\ebbb
\parbox{7cm}{ 
\1\9\2\7\3\5\4\4\5\3\7\2\9\1\9\2\eeee
\9\1\9\2\8\3\6\4\4\6\3\8\2\9\1\9\eeee
\2\9\1\9\2\7\3\5\4\4\5\3\7\2\9\1\eeee
\7\2\9\1\9\2\8\3\6\4\4\6\3\8\2\9\eeee
\3\8\2\9\1\9\2\7\3\5\4\4\5\3\7\2\eeee
\5\3\7\2\9\1\9\2\8\3\6\4\4\6\3\8\eeee
\4\6\3\8\2\9\1\9\2\7\3\5\4\4\5\3\eeee
\4\4\5\3\7\2\9\1\9\2\8\3\6\4\4\6\eeee
\5\4\4\6\3\8\2\9\1\9\2\7\3\5\4\4\eeee
\3\6\4\4\5\3\7\2\9\1\9\2\8\3\6\4\eeee
\7\3\5\4\4\6\3\8\2\9\1\9\2\7\3\5\eeee
\2\8\3\6\4\4\5\3\7\2\9\1\9\2\8\3\eeee
} 

\baaa
9-64
\eaaa
\bbbb
0&0&0&0&0&0&0&1&3\\
0&0&0&0&0&0&1&2&1\\
0&0&0&0&0&1&2&1&0\\
0&0&0&0&1&2&1&0&0\\
0&0&0&1&2&1&0&0&0\\
0&0&1&2&1&0&0&0&0\\
0&1&2&1&0&0&0&0&0\\
1&2&1&0&0&0&0&0&0\\
3&1&0&0&0&0&0&0&0\\
\ebbb
\parbox{7cm}{ 
\1\9\1\9\1\9\1\9\1\9\1\9\1\9\1\9\eeee
\8\2\8\2\8\2\8\2\8\2\8\2\8\2\8\2\eeee
\3\7\3\7\3\7\3\7\3\7\3\7\3\7\3\7\eeee
\6\4\6\4\6\4\6\4\6\4\6\4\6\4\6\4\eeee
\5\5\5\5\5\5\5\5\5\5\5\5\5\5\5\5\eeee
\4\6\4\6\4\6\4\6\4\6\4\6\4\6\4\6\eeee
\7\3\7\3\7\3\7\3\7\3\7\3\7\3\7\3\eeee
\2\8\2\8\2\8\2\8\2\8\2\8\2\8\2\8\eeee
\9\1\9\1\9\1\9\1\9\1\9\1\9\1\9\1\eeee
\1\9\1\9\1\9\1\9\1\9\1\9\1\9\1\9\eeee
\8\2\8\2\8\2\8\2\8\2\8\2\8\2\8\2\eeee
\3\7\3\7\3\7\3\7\3\7\3\7\3\7\3\7\eeee
} 

\baaa
9-65
\eaaa
\bbbb
0&0&0&0&0&0&0&2&2\\
0&0&0&0&0&0&0&2&2\\
0&0&0&0&0&0&2&0&2\\
0&0&0&0&0&0&2&0&2\\
0&0&0&0&0&0&2&2&0\\
0&0&0&0&0&0&2&2&0\\
0&0&1&1&1&1&0&0&0\\
1&1&0&0&1&1&0&0&0\\
1&1&1&1&0&0&0&0&0\\
\ebbb
\parbox{7cm}{ 
\1\8\6\7\3\9\2\8\5\7\4\9\1\8\6\7\eeee
\8\5\7\4\9\1\8\6\7\3\9\2\8\5\7\4\eeee
\6\7\3\9\2\8\5\7\4\9\1\8\6\7\3\9\eeee
\7\4\9\1\8\6\7\3\9\2\8\5\7\4\9\1\eeee
\3\9\2\8\5\7\4\9\1\8\6\7\3\9\2\8\eeee
\9\1\8\6\7\3\9\2\8\5\7\4\9\1\8\6\eeee
\2\8\5\7\4\9\1\8\6\7\3\9\2\8\5\7\eeee
\8\6\7\3\9\2\8\5\7\4\9\1\8\6\7\3\eeee
\5\7\4\9\1\8\6\7\3\9\2\8\5\7\4\9\eeee
\7\3\9\2\8\5\7\4\9\1\8\6\7\3\9\2\eeee
\4\9\1\8\6\7\3\9\2\8\5\7\4\9\1\8\eeee
\9\2\8\5\7\4\9\1\8\6\7\3\9\2\8\5\eeee
} 

\baaa
9-66
\eaaa
\bbbb
0&0&0&0&0&0&0&2&2\\
0&0&0&0&0&0&0&2&2\\
0&0&0&0&0&0&2&0&2\\
0&0&0&0&0&0&2&0&2\\
0&0&0&0&0&2&0&2&0\\
0&0&0&0&2&0&2&0&0\\
0&0&1&1&0&2&0&0&0\\
1&1&0&0&2&0&0&0&0\\
1&1&1&1&0&0&0&0&0\\
\ebbb
\parbox{7cm}{ 
\1\8\5\6\7\4\9\1\8\5\6\7\3\9\2\8\eeee
\8\5\6\7\3\9\2\8\5\6\7\4\9\1\8\5\eeee
\5\6\7\4\9\1\8\5\6\7\3\9\2\8\5\6\eeee
\6\7\3\9\2\8\5\6\7\4\9\1\8\5\6\7\eeee
\7\4\9\1\8\5\6\7\3\9\2\8\5\6\7\3\eeee
\3\9\2\8\5\6\7\4\9\1\8\5\6\7\4\9\eeee
\9\1\8\5\6\7\3\9\2\8\5\6\7\3\9\2\eeee
\2\8\5\6\7\4\9\1\8\5\6\7\4\9\1\8\eeee
\8\5\6\7\3\9\2\8\5\6\7\3\9\2\8\5\eeee
\5\6\7\4\9\1\8\5\6\7\4\9\1\8\5\6\eeee
\6\7\3\9\2\8\5\6\7\3\9\2\8\5\6\7\eeee
\7\4\9\1\8\5\6\7\4\9\1\8\5\6\7\3\eeee
} 

\baaa
9-67
\eaaa
\bbbb
0&0&0&0&0&0&0&2&2\\
0&0&0&0&0&0&0&2&2\\
0&0&0&0&0&0&2&0&2\\
0&0&0&0&0&0&2&0&2\\
0&0&0&0&0&2&0&2&0\\
0&0&0&0&2&2&0&0&0\\
0&0&1&1&0&0&2&0&0\\
1&1&0&0&2&0&0&0&0\\
1&1&1&1&0&0&0&0&0\\
\ebbb
\parbox{7cm}{ 
\1\8\5\6\6\5\8\1\9\4\7\7\3\9\2\8\eeee
\8\5\6\6\5\8\2\9\3\7\7\4\9\1\8\5\eeee
\5\6\6\5\8\1\9\4\7\7\3\9\2\8\5\6\eeee
\6\6\5\8\2\9\3\7\7\4\9\1\8\5\6\6\eeee
\6\5\8\1\9\4\7\7\3\9\2\8\5\6\6\5\eeee
\5\8\2\9\3\7\7\4\9\1\8\5\6\6\5\8\eeee
\8\1\9\4\7\7\3\9\2\8\5\6\6\5\8\2\eeee
\2\9\3\7\7\4\9\1\8\5\6\6\5\8\1\9\eeee
\9\4\7\7\3\9\2\8\5\6\6\5\8\2\9\3\eeee
\3\7\7\4\9\1\8\5\6\6\5\8\1\9\4\7\eeee
\7\7\3\9\2\8\5\6\6\5\8\2\9\3\7\7\eeee
\7\4\9\1\8\5\6\6\5\8\1\9\4\7\7\3\eeee
} 

\baaa
9-68
\eaaa
\bbbb
0&0&0&0&0&0&0&2&2\\
0&0&0&0&0&0&0&2&2\\
0&0&0&0&0&0&2&0&2\\
0&0&0&0&0&0&2&0&2\\
0&0&0&0&0&2&0&2&0\\
0&0&0&0&2&2&0&0&0\\
0&0&2&2&0&0&0&0&0\\
1&1&0&0&2&0&0&0&0\\
1&1&1&1&0&0&0&0&0\\
\ebbb
\parbox{7cm}{ 
\1\8\5\6\6\5\8\1\9\4\7\3\9\2\8\5\eeee
\8\5\6\6\5\8\2\9\3\7\4\9\1\8\5\6\eeee
\5\6\6\5\8\1\9\4\7\3\9\2\8\5\6\6\eeee
\6\6\5\8\2\9\3\7\4\9\1\8\5\6\6\5\eeee
\6\5\8\1\9\4\7\3\9\2\8\5\6\6\5\8\eeee
\5\8\2\9\3\7\4\9\1\8\5\6\6\5\8\2\eeee
\8\1\9\4\7\3\9\2\8\5\6\6\5\8\1\9\eeee
\2\9\3\7\4\9\1\8\5\6\6\5\8\2\9\3\eeee
\9\4\7\3\9\2\8\5\6\6\5\8\1\9\4\7\eeee
\3\7\4\9\1\8\5\6\6\5\8\2\9\3\7\4\eeee
\7\3\9\2\8\5\6\6\5\8\1\9\4\7\3\9\eeee
\4\9\1\8\5\6\6\5\8\2\9\3\7\4\9\1\eeee
} 

\baaa
9-69
\eaaa
\bbbb
0&0&0&0&0&0&0&2&2\\
0&0&0&0&0&0&0&2&2\\
0&0&0&0&0&0&2&0&2\\
0&0&0&0&0&0&2&0&2\\
0&0&0&0&1&1&0&2&0\\
0&0&0&0&1&1&0&2&0\\
0&0&1&1&0&0&2&0&0\\
1&1&0&0&1&1&0&0&0\\
1&1&1&1&0&0&0&0&0\\
\ebbb
\parbox{7cm}{ 
\1\8\6\5\8\2\9\3\7\7\4\9\1\8\6\5\eeee
\8\5\6\8\1\9\4\7\7\3\9\2\8\5\6\8\eeee
\6\5\8\2\9\3\7\7\4\9\1\8\6\5\8\2\eeee
\6\8\1\9\4\7\7\3\9\2\8\5\6\8\1\9\eeee
\8\2\9\3\7\7\4\9\1\8\6\5\8\2\9\3\eeee
\1\9\4\7\7\3\9\2\8\5\6\8\1\9\4\7\eeee
\9\3\7\7\4\9\1\8\6\5\8\2\9\3\7\7\eeee
\4\7\7\3\9\2\8\5\6\8\1\9\4\7\7\3\eeee
\7\7\4\9\1\8\6\5\8\2\9\3\7\7\4\9\eeee
\7\3\9\2\8\5\6\8\1\9\4\7\7\3\9\2\eeee
\4\9\1\8\6\5\8\2\9\3\7\7\4\9\1\8\eeee
\9\2\8\5\6\8\1\9\4\7\7\3\9\2\8\5\eeee
} 

\baaa
9-70
\eaaa
\bbbb
0&0&0&0&0&0&0&2&2\\
0&0&0&0&0&0&0&2&2\\
0&0&0&0&0&0&2&0&2\\
0&0&0&0&0&0&2&0&2\\
0&0&0&0&1&1&0&2&0\\
0&0&0&0&1&1&0&2&0\\
0&0&2&2&0&0&0&0&0\\
1&1&0&0&1&1&0&0&0\\
1&1&1&1&0&0&0&0&0\\
\ebbb
\parbox{7cm}{ 
\1\8\6\5\8\2\9\3\7\4\9\1\8\6\6\8\eeee
\8\5\6\8\1\9\4\7\3\9\2\8\5\5\8\2\eeee
\6\5\8\2\9\3\7\4\9\1\8\6\6\8\1\9\eeee
\6\8\1\9\4\7\3\9\2\8\5\5\8\2\9\3\eeee
\8\2\9\3\7\4\9\1\8\6\6\8\1\9\4\7\eeee
\1\9\4\7\3\9\2\8\5\5\8\2\9\3\7\4\eeee
\9\3\7\4\9\1\8\6\6\8\1\9\4\7\3\9\eeee
\4\7\3\9\2\8\5\5\8\2\9\3\7\4\9\1\eeee
\7\4\9\1\8\6\6\8\1\9\4\7\3\9\2\8\eeee
\3\9\2\8\5\5\8\2\9\3\7\4\9\1\8\6\eeee
\9\1\8\6\6\8\1\9\4\7\3\9\2\8\5\6\eeee
\2\8\5\5\8\2\9\3\7\4\9\1\8\6\5\8\eeee
} 

\baaa
9-71
\eaaa
\bbbb
0&0&0&0&0&0&0&2&2\\
0&0&0&0&0&0&0&2&2\\
0&0&0&0&0&0&2&0&2\\
0&0&0&0&0&0&2&0&2\\
0&0&0&0&2&0&0&2&0\\
0&0&0&0&0&2&2&0&0\\
0&0&1&1&0&2&0&0&0\\
1&1&0&0&2&0&0&0&0\\
1&1&1&1&0&0&0&0&0\\
\ebbb
\parbox{7cm}{ 
\1\8\5\5\8\2\9\3\7\6\6\7\3\9\2\8\eeee
\8\5\5\8\1\9\4\7\6\6\7\4\9\1\8\5\eeee
\5\5\8\2\9\3\7\6\6\7\3\9\2\8\5\5\eeee
\5\8\1\9\4\7\6\6\7\4\9\1\8\5\5\8\eeee
\8\2\9\3\7\6\6\7\3\9\2\8\5\5\8\1\eeee
\1\9\4\7\6\6\7\4\9\1\8\5\5\8\2\9\eeee
\9\3\7\6\6\7\3\9\2\8\5\5\8\1\9\4\eeee
\4\7\6\6\7\4\9\1\8\5\5\8\2\9\3\7\eeee
\7\6\6\7\3\9\2\8\5\5\8\1\9\4\7\6\eeee
\6\6\7\4\9\1\8\5\5\8\2\9\3\7\6\6\eeee
\6\7\3\9\2\8\5\5\8\1\9\4\7\6\6\7\eeee
\7\4\9\1\8\5\5\8\2\9\3\7\6\6\7\3\eeee
} 

\baaa
9-72
\eaaa
\bbbb
0&0&0&0&0&0&0&2&2\\
0&0&0&0&0&0&0&2&2\\
0&0&0&0&0&0&2&0&2\\
0&0&0&0&0&1&1&1&1\\
0&0&0&0&0&2&0&2&0\\
0&0&0&2&2&0&0&0&0\\
0&0&2&2&0&0&0&0&0\\
1&1&0&1&1&0&0&0&0\\
1&1&1&1&0&0&0&0&0\\
\ebbb
\parbox{7cm}{ 
\1\8\4\7\3\9\2\9\3\7\4\8\1\9\4\6\eeee
\8\5\6\4\9\1\8\4\7\3\9\2\9\3\7\4\eeee
\4\6\5\8\2\8\5\6\4\9\1\8\4\7\3\9\eeee
\7\4\8\1\9\4\6\5\8\2\8\5\6\4\9\1\eeee
\3\9\2\9\3\7\4\8\1\9\4\6\5\8\2\8\eeee
\9\1\8\4\7\3\9\2\9\3\7\4\8\1\9\4\eeee
\2\8\5\6\4\9\1\8\4\7\3\9\2\9\3\7\eeee
\9\4\6\5\8\2\8\5\6\4\9\1\8\4\7\3\eeee
\3\7\4\8\1\9\4\6\5\8\2\8\5\6\4\9\eeee
\7\3\9\2\9\3\7\4\8\1\9\4\6\5\8\2\eeee
\4\9\1\8\4\7\3\9\2\9\3\7\4\8\1\9\eeee
\8\2\8\5\6\4\9\1\8\4\7\3\9\2\9\3\eeee
} 

\baaa
9-73
\eaaa
\bbbb
0&0&0&0&0&0&0&2&2\\
0&0&0&0&0&0&0&2&2\\
0&0&0&0&0&0&2&0&2\\
0&0&0&0&0&2&0&2&0\\
0&0&0&0&0&2&2&0&0\\
0&0&0&2&2&0&0&0&0\\
0&0&2&0&2&0&0&0&0\\
1&1&0&2&0&0&0&0&0\\
1&1&2&0&0&0&0&0&0\\
\ebbb
\parbox{7cm}{ 
\1\8\4\6\5\7\3\9\1\8\4\6\5\7\3\9\eeee
\8\4\6\5\7\3\9\2\8\4\6\5\7\3\9\2\eeee
\4\6\5\7\3\9\1\8\4\6\5\7\3\9\1\8\eeee
\6\5\7\3\9\2\8\4\6\5\7\3\9\2\8\4\eeee
\5\7\3\9\1\8\4\6\5\7\3\9\1\8\4\6\eeee
\7\3\9\2\8\4\6\5\7\3\9\2\8\4\6\5\eeee
\3\9\1\8\4\6\5\7\3\9\1\8\4\6\5\7\eeee
\9\2\8\4\6\5\7\3\9\2\8\4\6\5\7\3\eeee
\1\8\4\6\5\7\3\9\1\8\4\6\5\7\3\9\eeee
\8\4\6\5\7\3\9\2\8\4\6\5\7\3\9\2\eeee
\4\6\5\7\3\9\1\8\4\6\5\7\3\9\1\8\eeee
\6\5\7\3\9\2\8\4\6\5\7\3\9\2\8\4\eeee
} 

\baaa
9-74
\eaaa
\bbbb
0&0&0&0&0&0&0&2&2\\
0&0&0&0&0&0&0&2&2\\
0&0&0&0&0&0&2&0&2\\
0&0&0&0&0&2&0&2&0\\
0&0&0&0&2&0&2&0&0\\
0&0&0&2&0&2&0&0&0\\
0&0&2&0&2&0&0&0&0\\
1&1&0&2&0&0&0&0&0\\
1&1&2&0&0&0&0&0&0\\
\ebbb
\parbox{7cm}{ 
\1\8\4\6\6\4\8\1\9\3\7\5\5\7\3\9\eeee
\8\4\6\6\4\8\2\9\3\7\5\5\7\3\9\2\eeee
\4\6\6\4\8\1\9\3\7\5\5\7\3\9\1\8\eeee
\6\6\4\8\2\9\3\7\5\5\7\3\9\2\8\4\eeee
\6\4\8\1\9\3\7\5\5\7\3\9\1\8\4\6\eeee
\4\8\2\9\3\7\5\5\7\3\9\2\8\4\6\6\eeee
\8\1\9\3\7\5\5\7\3\9\1\8\4\6\6\4\eeee
\2\9\3\7\5\5\7\3\9\2\8\4\6\6\4\8\eeee
\9\3\7\5\5\7\3\9\1\8\4\6\6\4\8\1\eeee
\3\7\5\5\7\3\9\2\8\4\6\6\4\8\2\9\eeee
\7\5\5\7\3\9\1\8\4\6\6\4\8\1\9\3\eeee
\5\5\7\3\9\2\8\4\6\6\4\8\2\9\3\7\eeee
} 

\baaa
9-75
\eaaa
\bbbb
0&0&0&0&0&0&0&2&2\\
0&0&0&0&0&0&0&2&2\\
0&0&0&0&0&0&2&0&2\\
0&0&0&0&0&2&2&0&0\\
0&0&0&0&0&2&2&0&0\\
0&0&0&1&1&0&0&2&0\\
0&0&2&1&1&0&0&0&0\\
1&1&0&0&0&2&0&0&0\\
1&1&2&0&0&0&0&0&0\\
\ebbb
\parbox{7cm}{ 
\1\8\6\4\7\3\9\1\8\6\5\7\3\9\2\8\eeee
\8\6\5\7\3\9\2\8\6\4\7\3\9\1\8\6\eeee
\6\4\7\3\9\1\8\6\5\7\3\9\2\8\6\5\eeee
\5\7\3\9\2\8\6\4\7\3\9\1\8\6\4\7\eeee
\7\3\9\1\8\6\5\7\3\9\2\8\6\5\7\3\eeee
\3\9\2\8\6\4\7\3\9\1\8\6\4\7\3\9\eeee
\9\1\8\6\5\7\3\9\2\8\6\5\7\3\9\2\eeee
\2\8\6\4\7\3\9\1\8\6\4\7\3\9\1\8\eeee
\8\6\5\7\3\9\2\8\6\5\7\3\9\2\8\6\eeee
\6\4\7\3\9\1\8\6\4\7\3\9\1\8\6\5\eeee
\5\7\3\9\2\8\6\5\7\3\9\2\8\6\4\7\eeee
\7\3\9\1\8\6\4\7\3\9\1\8\6\5\7\3\eeee
} 

\baaa
9-76
\eaaa
\bbbb
0&0&0&0&0&0&0&2&2\\
0&0&0&0&0&0&0&2&2\\
0&0&0&0&0&0&2&0&2\\
0&0&0&0&0&2&2&0&0\\
0&0&0&0&0&2&2&0&0\\
0&0&0&1&1&2&0&0&0\\
0&0&2&1&1&0&0&0&0\\
1&1&0&0&0&0&0&2&0\\
1&1&2&0&0&0&0&0&0\\
\ebbb
\parbox{7cm}{ 
\1\8\8\1\9\3\7\4\6\6\5\7\3\9\2\8\eeee
\8\8\2\9\3\7\5\6\6\4\7\3\9\1\8\8\eeee
\8\1\9\3\7\4\6\6\5\7\3\9\2\8\8\2\eeee
\2\9\3\7\5\6\6\4\7\3\9\1\8\8\1\9\eeee
\9\3\7\4\6\6\5\7\3\9\2\8\8\2\9\3\eeee
\3\7\5\6\6\4\7\3\9\1\8\8\1\9\3\7\eeee
\7\4\6\6\5\7\3\9\2\8\8\2\9\3\7\5\eeee
\5\6\6\4\7\3\9\1\8\8\1\9\3\7\4\6\eeee
\6\6\5\7\3\9\2\8\8\2\9\3\7\5\6\6\eeee
\6\4\7\3\9\1\8\8\1\9\3\7\4\6\6\5\eeee
\5\7\3\9\2\8\8\2\9\3\7\5\6\6\4\7\eeee
\7\3\9\1\8\8\1\9\3\7\4\6\6\5\7\3\eeee
} 

\baaa
9-77
\eaaa
\bbbb
0&0&0&0&0&0&0&2&2\\
0&0&0&0&0&0&0&2&2\\
0&0&0&0&0&0&2&0&2\\
0&0&0&0&0&2&2&0&0\\
0&0&0&0&0&2&2&0&0\\
0&0&0&1&1&2&0&0&0\\
0&0&2&1&1&0&0&0&0\\
2&2&0&0&0&0&0&0&0\\
1&1&2&0&0&0&0&0&0\\
\ebbb
\parbox{7cm}{ 
\1\8\2\9\3\7\5\6\6\5\7\3\9\2\8\1\eeee
\8\1\9\3\7\4\6\6\4\7\3\9\1\8\2\9\eeee
\2\9\3\7\5\6\6\5\7\3\9\2\8\1\9\3\eeee
\9\3\7\4\6\6\4\7\3\9\1\8\2\9\3\7\eeee
\3\7\5\6\6\5\7\3\9\2\8\1\9\3\7\4\eeee
\7\4\6\6\4\7\3\9\1\8\2\9\3\7\5\6\eeee
\5\6\6\5\7\3\9\2\8\1\9\3\7\4\6\6\eeee
\6\6\4\7\3\9\1\8\2\9\3\7\5\6\6\4\eeee
\6\5\7\3\9\2\8\1\9\3\7\4\6\6\5\7\eeee
\4\7\3\9\1\8\2\9\3\7\5\6\6\4\7\3\eeee
\7\3\9\2\8\1\9\3\7\4\6\6\5\7\3\9\eeee
\3\9\1\8\2\9\3\7\5\6\6\4\7\3\9\1\eeee
} 

\baaa
9-78
\eaaa
\bbbb
0&0&0&0&0&0&0&2&2\\
0&0&0&0&0&0&0&2&2\\
0&0&0&0&0&0&2&0&2\\
0&0&0&0&0&2&2&0&0\\
0&0&0&0&0&2&2&0&0\\
0&0&0&2&2&0&0&0&0\\
0&0&2&1&1&0&0&0&0\\
2&2&0&0&0&0&0&0&0\\
1&1&2&0&0&0&0&0&0\\
\ebbb
\parbox{7cm}{ 
\1\8\2\9\3\7\5\6\4\7\3\9\1\8\2\9\eeee
\8\1\9\3\7\4\6\5\7\3\9\2\8\1\9\3\eeee
\2\9\3\7\5\6\4\7\3\9\1\8\2\9\3\7\eeee
\9\3\7\4\6\5\7\3\9\2\8\1\9\3\7\4\eeee
\3\7\5\6\4\7\3\9\1\8\2\9\3\7\5\6\eeee
\7\4\6\5\7\3\9\2\8\1\9\3\7\4\6\5\eeee
\5\6\4\7\3\9\1\8\2\9\3\7\5\6\4\7\eeee
\6\5\7\3\9\2\8\1\9\3\7\4\6\5\7\3\eeee
\4\7\3\9\1\8\2\9\3\7\5\6\4\7\3\9\eeee
\7\3\9\2\8\1\9\3\7\4\6\5\7\3\9\2\eeee
\3\9\1\8\2\9\3\7\5\6\4\7\3\9\1\8\eeee
\9\2\8\1\9\3\7\4\6\5\7\3\9\2\8\1\eeee
} 

\baaa
9-79
\eaaa
\bbbb
0&0&0&0&0&0&0&2&2\\
0&0&0&0&0&0&0&2&2\\
0&0&0&0&0&0&2&0&2\\
0&0&0&0&0&2&2&0&0\\
0&0&0&0&2&0&0&2&0\\
0&0&0&2&0&2&0&0&0\\
0&0&2&2&0&0&0&0&0\\
1&1&0&0&2&0&0&0&0\\
1&1&2&0&0&0&0&0&0\\
\ebbb
\parbox{7cm}{ 
\1\8\5\5\8\2\9\3\7\4\6\6\4\7\3\9\eeee
\8\5\5\8\1\9\3\7\4\6\6\4\7\3\9\2\eeee
\5\5\8\2\9\3\7\4\6\6\4\7\3\9\1\8\eeee
\5\8\1\9\3\7\4\6\6\4\7\3\9\2\8\5\eeee
\8\2\9\3\7\4\6\6\4\7\3\9\1\8\5\5\eeee
\1\9\3\7\4\6\6\4\7\3\9\2\8\5\5\8\eeee
\9\3\7\4\6\6\4\7\3\9\1\8\5\5\8\2\eeee
\3\7\4\6\6\4\7\3\9\2\8\5\5\8\1\9\eeee
\7\4\6\6\4\7\3\9\1\8\5\5\8\2\9\3\eeee
\4\6\6\4\7\3\9\2\8\5\5\8\1\9\3\7\eeee
\6\6\4\7\3\9\1\8\5\5\8\2\9\3\7\4\eeee
\6\4\7\3\9\2\8\5\5\8\1\9\3\7\4\6\eeee
} 

\baaa
9-80
\eaaa
\bbbb
0&0&0&0&0&0&0&2&2\\
0&0&0&0&0&0&0&2&2\\
0&0&0&0&0&0&2&0&2\\
0&0&0&0&0&2&2&0&0\\
0&0&0&0&2&2&0&0&0\\
0&0&0&2&2&0&0&0&0\\
0&0&2&2&0&0&0&0&0\\
1&1&0&0&0&0&0&2&0\\
1&1&2&0&0&0&0&0&0\\
\ebbb
\parbox{7cm}{ 
\1\8\8\1\9\3\7\4\6\5\5\6\4\7\3\9\eeee
\8\8\2\9\3\7\4\6\5\5\6\4\7\3\9\2\eeee
\8\1\9\3\7\4\6\5\5\6\4\7\3\9\1\8\eeee
\2\9\3\7\4\6\5\5\6\4\7\3\9\2\8\8\eeee
\9\3\7\4\6\5\5\6\4\7\3\9\1\8\8\1\eeee
\3\7\4\6\5\5\6\4\7\3\9\2\8\8\2\9\eeee
\7\4\6\5\5\6\4\7\3\9\1\8\8\1\9\3\eeee
\4\6\5\5\6\4\7\3\9\2\8\8\2\9\3\7\eeee
\6\5\5\6\4\7\3\9\1\8\8\1\9\3\7\4\eeee
\5\5\6\4\7\3\9\2\8\8\2\9\3\7\4\6\eeee
\5\6\4\7\3\9\1\8\8\1\9\3\7\4\6\5\eeee
\6\4\7\3\9\2\8\8\2\9\3\7\4\6\5\5\eeee
} 

\baaa
9-81
\eaaa
\bbbb
0&0&0&0&0&0&0&2&2\\
0&0&0&0&0&0&0&2&2\\
0&0&0&0&0&0&2&0&2\\
0&0&0&0&0&2&2&0&0\\
0&0&0&0&2&2&0&0&0\\
0&0&0&2&2&0&0&0&0\\
0&0&2&2&0&0&0&0&0\\
2&2&0&0&0&0&0&0&0\\
1&1&2&0&0&0&0&0&0\\
\ebbb
\parbox{7cm}{ 
\1\8\2\9\3\7\4\6\5\5\6\4\7\3\9\1\eeee
\8\1\9\3\7\4\6\5\5\6\4\7\3\9\2\8\eeee
\2\9\3\7\4\6\5\5\6\4\7\3\9\1\8\2\eeee
\9\3\7\4\6\5\5\6\4\7\3\9\2\8\1\9\eeee
\3\7\4\6\5\5\6\4\7\3\9\1\8\2\9\3\eeee
\7\4\6\5\5\6\4\7\3\9\2\8\1\9\3\7\eeee
\4\6\5\5\6\4\7\3\9\1\8\2\9\3\7\4\eeee
\6\5\5\6\4\7\3\9\2\8\1\9\3\7\4\6\eeee
\5\5\6\4\7\3\9\1\8\2\9\3\7\4\6\5\eeee
\5\6\4\7\3\9\2\8\1\9\3\7\4\6\5\5\eeee
\6\4\7\3\9\1\8\2\9\3\7\4\6\5\5\6\eeee
\4\7\3\9\2\8\1\9\3\7\4\6\5\5\6\4\eeee
} 

\baaa
9-82
\eaaa
\bbbb
0&0&0&0&0&0&0&2&2\\
0&0&0&0&0&0&0&2&2\\
0&0&0&0&0&0&2&0&2\\
0&0&0&0&1&1&0&2&0\\
0&0&0&2&1&1&0&0&0\\
0&0&0&2&1&1&0&0&0\\
0&0&2&0&0&0&2&0&0\\
1&1&0&2&0&0&0&0&0\\
1&1&2&0&0&0&0&0&0\\
\ebbb
\parbox{7cm}{ 
\1\8\4\5\5\4\8\1\9\3\7\7\3\9\2\8\eeee
\8\4\6\6\4\8\2\9\3\7\7\3\9\1\8\4\eeee
\4\5\5\4\8\1\9\3\7\7\3\9\2\8\4\6\eeee
\6\6\4\8\2\9\3\7\7\3\9\1\8\4\5\5\eeee
\5\4\8\1\9\3\7\7\3\9\2\8\4\6\6\4\eeee
\4\8\2\9\3\7\7\3\9\1\8\4\5\5\4\8\eeee
\8\1\9\3\7\7\3\9\2\8\4\6\6\4\8\2\eeee
\2\9\3\7\7\3\9\1\8\4\5\5\4\8\1\9\eeee
\9\3\7\7\3\9\2\8\4\6\6\4\8\2\9\3\eeee
\3\7\7\3\9\1\8\4\5\5\4\8\1\9\3\7\eeee
\7\7\3\9\2\8\4\6\6\4\8\2\9\3\7\7\eeee
\7\3\9\1\8\4\5\5\4\8\1\9\3\7\7\3\eeee
} 

\baaa
9-83
\eaaa
\bbbb
0&0&0&0&0&0&0&2&2\\
0&0&0&0&0&0&0&2&2\\
0&0&0&0&0&0&2&0&2\\
0&0&0&0&1&1&2&0&0\\
0&0&0&2&1&1&0&0&0\\
0&0&0&2&1&1&0&0&0\\
0&0&2&2&0&0&0&0&0\\
1&1&0&0&0&0&0&2&0\\
1&1&2&0&0&0&0&0&0\\
\ebbb
\parbox{7cm}{ 
\1\8\8\1\9\3\7\4\5\6\4\7\3\9\2\8\eeee
\8\8\2\9\3\7\4\6\5\4\7\3\9\1\8\8\eeee
\8\1\9\3\7\4\5\6\4\7\3\9\2\8\8\2\eeee
\2\9\3\7\4\6\5\4\7\3\9\1\8\8\1\9\eeee
\9\3\7\4\5\6\4\7\3\9\2\8\8\2\9\3\eeee
\3\7\4\6\5\4\7\3\9\1\8\8\1\9\3\7\eeee
\7\4\5\6\4\7\3\9\2\8\8\2\9\3\7\4\eeee
\4\6\5\4\7\3\9\1\8\8\1\9\3\7\4\5\eeee
\5\6\4\7\3\9\2\8\8\2\9\3\7\4\6\5\eeee
\5\4\7\3\9\1\8\8\1\9\3\7\4\5\6\4\eeee
\4\7\3\9\2\8\8\2\9\3\7\4\6\5\4\7\eeee
\7\3\9\1\8\8\1\9\3\7\4\5\6\4\7\3\eeee
} 

\baaa
9-84
\eaaa
\bbbb
0&0&0&0&0&0&0&2&2\\
0&0&0&0&0&0&0&2&2\\
0&0&0&0&0&0&2&0&2\\
0&0&0&0&1&1&2&0&0\\
0&0&0&2&1&1&0&0&0\\
0&0&0&2&1&1&0&0&0\\
0&0&2&2&0&0&0&0&0\\
2&2&0&0&0&0&0&0&0\\
1&1&2&0&0&0&0&0&0\\
\ebbb
\parbox{7cm}{ 
\1\8\2\9\3\7\4\5\5\4\7\3\9\2\8\1\eeee
\8\1\9\3\7\4\6\6\4\7\3\9\1\8\2\9\eeee
\2\9\3\7\4\5\5\4\7\3\9\2\8\1\9\3\eeee
\9\3\7\4\6\6\4\7\3\9\1\8\2\9\3\7\eeee
\3\7\4\5\5\4\7\3\9\2\8\1\9\3\7\4\eeee
\7\4\6\6\4\7\3\9\1\8\2\9\3\7\4\6\eeee
\4\5\5\4\7\3\9\2\8\1\9\3\7\4\5\6\eeee
\6\6\4\7\3\9\1\8\2\9\3\7\4\6\5\4\eeee
\5\4\7\3\9\2\8\1\9\3\7\4\5\6\4\7\eeee
\4\7\3\9\1\8\2\9\3\7\4\6\5\4\7\3\eeee
\7\3\9\2\8\1\9\3\7\4\5\6\4\7\3\9\eeee
\3\9\1\8\2\9\3\7\4\6\5\4\7\3\9\1\eeee
} 

\baaa
9-85
\eaaa
\bbbb
0&0&0&0&0&0&0&2&2\\
0&0&0&0&0&0&0&2&2\\
0&0&0&0&0&0&2&0&2\\
0&0&0&1&0&1&0&2&0\\
0&0&0&0&2&0&2&0&0\\
0&0&0&1&0&1&0&2&0\\
0&0&2&0&2&0&0&0&0\\
1&1&0&1&0&1&0&0&0\\
1&1&2&0&0&0&0&0&0\\
\ebbb
\parbox{7cm}{ 
\1\8\6\4\8\2\9\3\7\5\5\7\3\9\2\8\eeee
\8\4\6\8\1\9\3\7\5\5\7\3\9\1\8\6\eeee
\6\4\8\2\9\3\7\5\5\7\3\9\2\8\4\6\eeee
\6\8\1\9\3\7\5\5\7\3\9\1\8\6\4\8\eeee
\8\2\9\3\7\5\5\7\3\9\2\8\4\6\8\1\eeee
\1\9\3\7\5\5\7\3\9\1\8\6\4\8\2\9\eeee
\9\3\7\5\5\7\3\9\2\8\4\6\8\1\9\3\eeee
\3\7\5\5\7\3\9\1\8\6\4\8\2\9\3\7\eeee
\7\5\5\7\3\9\2\8\4\6\8\1\9\3\7\5\eeee
\5\5\7\3\9\1\8\6\4\8\2\9\3\7\5\5\eeee
\5\7\3\9\2\8\4\6\8\1\9\3\7\5\5\7\eeee
\7\3\9\1\8\6\4\8\2\9\3\7\5\5\7\3\eeee
} 

\baaa
9-86
\eaaa
\bbbb
0&0&0&0&0&0&0&2&2\\
0&0&0&0&0&0&0&2&2\\
0&0&0&0&0&0&2&0&2\\
0&0&0&1&0&1&2&0&0\\
0&0&0&0&2&0&0&2&0\\
0&0&0&1&0&1&2&0&0\\
0&0&2&1&0&1&0&0&0\\
1&1&0&0&2&0&0&0&0\\
1&1&2&0&0&0&0&0&0\\
\ebbb
\parbox{7cm}{ 
\1\8\5\5\8\2\9\3\7\6\6\7\3\9\2\8\eeee
\8\5\5\8\1\9\3\7\4\4\7\3\9\1\8\5\eeee
\5\5\8\2\9\3\7\6\6\7\3\9\2\8\5\5\eeee
\5\8\1\9\3\7\4\4\7\3\9\1\8\5\5\8\eeee
\8\2\9\3\7\6\6\7\3\9\2\8\5\5\8\1\eeee
\1\9\3\7\4\4\7\3\9\1\8\5\5\8\2\9\eeee
\9\3\7\6\6\7\3\9\2\8\5\5\8\1\9\3\eeee
\3\7\4\4\7\3\9\1\8\5\5\8\2\9\3\7\eeee
\7\6\6\7\3\9\2\8\5\5\8\1\9\3\7\4\eeee
\4\4\7\3\9\1\8\5\5\8\2\9\3\7\6\6\eeee
\6\7\3\9\2\8\5\5\8\1\9\3\7\4\4\7\eeee
\7\3\9\1\8\5\5\8\2\9\3\7\6\6\7\3\eeee
} 

\baaa
9-87
\eaaa
\bbbb
0&0&0&0&0&0&0&2&2\\
0&0&0&0&0&0&0&2&2\\
0&0&0&0&0&0&2&1&1\\
0&0&0&0&0&0&2&1&1\\
0&0&0&0&0&2&2&0&0\\
0&0&0&0&2&2&0&0&0\\
0&0&1&1&2&0&0&0&0\\
1&1&1&1&0&0&0&0&0\\
1&1&1&1&0&0&0&0&0\\
\ebbb
\parbox{7cm}{ 
\1\8\3\7\5\6\6\5\7\3\8\1\8\3\7\5\eeee
\8\2\9\4\7\5\6\6\5\7\4\9\2\9\4\7\eeee
\3\9\1\8\3\7\5\6\6\5\7\3\8\1\8\3\eeee
\7\4\8\2\9\4\7\5\6\6\5\7\4\9\2\9\eeee
\5\7\3\9\1\8\3\7\5\6\6\5\7\3\8\1\eeee
\6\5\7\4\8\2\9\4\7\5\6\6\5\7\4\9\eeee
\6\6\5\7\3\9\1\8\3\7\5\6\6\5\7\3\eeee
\5\6\6\5\7\4\8\2\9\4\7\5\6\6\5\7\eeee
\7\5\6\6\5\7\3\9\1\8\3\7\5\6\6\5\eeee
\3\7\5\6\6\5\7\4\8\2\9\4\7\5\6\6\eeee
\8\4\7\5\6\6\5\7\3\9\1\8\3\7\5\6\eeee
\1\9\3\7\5\6\6\5\7\4\8\2\9\4\7\5\eeee
} 

\baaa
9-88
\eaaa
\bbbb
0&0&0&0&0&0&0&2&2\\
0&0&0&0&0&0&0&2&2\\
0&0&0&0&0&0&2&1&1\\
0&0&0&0&0&0&2&1&1\\
0&0&0&0&1&1&2&0&0\\
0&0&0&0&1&1&2&0&0\\
0&0&1&1&1&1&0&0&0\\
1&1&1&1&0&0&0&0&0\\
1&1&1&1&0&0&0&0&0\\
\ebbb
\parbox{7cm}{ 
\1\8\3\7\5\5\7\3\8\1\8\3\7\5\5\7\eeee
\8\2\9\4\7\6\6\7\4\9\2\9\4\7\6\6\eeee
\3\9\1\8\3\7\5\5\7\3\8\1\8\3\7\5\eeee
\7\4\8\2\9\4\7\6\6\7\4\9\2\9\4\7\eeee
\5\7\3\9\1\8\3\7\5\5\7\3\8\1\8\3\eeee
\5\6\7\4\8\2\9\4\7\6\6\7\4\9\2\9\eeee
\7\6\5\7\3\9\1\8\3\7\5\5\7\3\8\1\eeee
\3\7\5\6\7\4\8\2\9\4\7\6\6\7\4\9\eeee
\8\4\7\6\5\7\3\9\1\8\3\7\5\5\7\3\eeee
\1\9\3\7\5\6\7\4\8\2\9\4\7\6\6\7\eeee
\8\2\8\4\7\6\5\7\3\9\1\8\3\7\5\5\eeee
\3\9\1\9\3\7\5\6\7\4\8\2\9\4\7\6\eeee
} 

\baaa
9-89
\eaaa
\bbbb
0&0&0&0&0&0&0&2&2\\
0&0&0&0&0&0&0&2&2\\
0&0&0&0&0&0&2&1&1\\
0&0&0&0&0&2&2&0&0\\
0&0&0&0&0&2&2&0&0\\
0&0&0&1&1&2&0&0&0\\
0&0&2&1&1&0&0&0&0\\
1&1&2&0&0&0&0&0&0\\
1&1&2&0&0&0&0&0&0\\
\ebbb
\parbox{7cm}{ 
\1\8\3\7\4\6\6\4\7\3\8\1\8\3\7\4\eeee
\8\2\9\3\7\5\6\6\5\7\3\9\2\9\3\7\eeee
\3\9\1\8\3\7\4\6\6\4\7\3\8\1\8\3\eeee
\7\3\8\2\9\3\7\5\6\6\5\7\3\9\2\9\eeee
\4\7\3\9\1\8\3\7\4\6\6\4\7\3\8\1\eeee
\6\5\7\3\8\2\9\3\7\5\6\6\5\7\3\9\eeee
\6\6\4\7\3\9\1\8\3\7\4\6\6\4\7\3\eeee
\4\6\6\5\7\3\8\2\9\3\7\5\6\6\5\7\eeee
\7\5\6\6\4\7\3\9\1\8\3\7\4\6\6\4\eeee
\3\7\4\6\6\5\7\3\8\2\9\3\7\5\6\6\eeee
\8\3\7\5\6\6\4\7\3\9\1\8\3\7\4\6\eeee
\1\9\3\7\4\6\6\5\7\3\8\2\9\3\7\5\eeee
} 

\baaa
9-90
\eaaa
\bbbb
0&0&0&0&0&0&0&2&2\\
0&0&0&0&0&0&0&2&2\\
0&0&0&0&0&0&2&1&1\\
0&0&0&0&0&2&2&0&0\\
0&0&0&0&0&2&2&0&0\\
0&0&0&2&2&0&0&0&0\\
0&0&2&1&1&0&0&0&0\\
1&1&2&0&0&0&0&0&0\\
1&1&2&0&0&0&0&0&0\\
\ebbb
\parbox{7cm}{ 
\1\8\3\7\4\6\4\7\3\8\1\8\3\7\4\6\eeee
\8\2\9\3\7\5\6\5\7\3\9\2\9\3\7\5\eeee
\3\9\1\8\3\7\4\6\4\7\3\8\1\8\3\7\eeee
\7\3\8\2\9\3\7\5\6\5\7\3\9\2\9\3\eeee
\4\7\3\9\1\8\3\7\4\6\4\7\3\8\1\8\eeee
\6\5\7\3\8\2\9\3\7\5\6\5\7\3\9\2\eeee
\4\6\4\7\3\9\1\8\3\7\4\6\4\7\3\8\eeee
\7\5\6\5\7\3\8\2\9\3\7\5\6\5\7\3\eeee
\3\7\4\6\4\7\3\9\1\8\3\7\4\6\4\7\eeee
\8\3\7\5\6\5\7\3\8\2\9\3\7\5\6\5\eeee
\1\9\3\7\4\6\4\7\3\9\1\8\3\7\4\6\eeee
\8\2\8\3\7\5\6\5\7\3\8\2\9\3\7\5\eeee
} 

\baaa
9-91
\eaaa
\bbbb
0&0&0&0&0&0&0&2&2\\
0&0&0&0&0&0&0&2&2\\
0&0&0&0&0&0&2&1&1\\
0&0&0&0&0&2&2&0&0\\
0&0&0&0&2&2&0&0&0\\
0&0&0&2&2&0&0&0&0\\
0&0&2&2&0&0&0&0&0\\
1&1&2&0&0&0&0&0&0\\
1&1&2&0&0&0&0&0&0\\
\ebbb
\parbox{7cm}{ 
\1\8\3\7\4\6\5\5\6\4\7\3\8\1\8\3\eeee
\8\2\9\3\7\4\6\5\5\6\4\7\3\9\2\9\eeee
\3\9\1\8\3\7\4\6\5\5\6\4\7\3\8\1\eeee
\7\3\8\2\9\3\7\4\6\5\5\6\4\7\3\9\eeee
\4\7\3\9\1\8\3\7\4\6\5\5\6\4\7\3\eeee
\6\4\7\3\8\2\9\3\7\4\6\5\5\6\4\7\eeee
\5\6\4\7\3\9\1\8\3\7\4\6\5\5\6\4\eeee
\5\5\6\4\7\3\8\2\9\3\7\4\6\5\5\6\eeee
\6\5\5\6\4\7\3\9\1\8\3\7\4\6\5\5\eeee
\4\6\5\5\6\4\7\3\8\2\9\3\7\4\6\5\eeee
\7\4\6\5\5\6\4\7\3\9\1\8\3\7\4\6\eeee
\3\7\4\6\5\5\6\4\7\3\8\2\9\3\7\4\eeee
} 

\baaa
9-92
\eaaa
\bbbb
0&0&0&0&0&0&0&2&2\\
0&0&0&0&0&0&0&2&2\\
0&0&0&0&0&0&2&1&1\\
0&0&0&0&1&1&2&0&0\\
0&0&0&2&1&1&0&0&0\\
0&0&0&2&1&1&0&0&0\\
0&0&2&2&0&0&0&0&0\\
1&1&2&0&0&0&0&0&0\\
1&1&2&0&0&0&0&0&0\\
\ebbb
\parbox{7cm}{ 
\1\8\3\7\4\5\5\4\7\3\8\1\8\3\7\4\eeee
\8\2\9\3\7\4\6\6\4\7\3\9\2\9\3\7\eeee
\3\9\1\8\3\7\4\5\5\4\7\3\8\1\8\3\eeee
\7\3\8\2\9\3\7\4\6\6\4\7\3\9\2\9\eeee
\4\7\3\9\1\8\3\7\4\5\5\4\7\3\8\1\eeee
\5\4\7\3\8\2\9\3\7\4\6\6\4\7\3\9\eeee
\5\6\4\7\3\9\1\8\3\7\4\5\5\4\7\3\eeee
\4\6\5\4\7\3\8\2\9\3\7\4\6\6\4\7\eeee
\7\4\5\6\4\7\3\9\1\8\3\7\4\5\5\4\eeee
\3\7\4\6\5\4\7\3\8\2\9\3\7\4\6\6\eeee
\8\3\7\4\5\6\4\7\3\9\1\8\3\7\4\5\eeee
\1\9\3\7\4\6\5\4\7\3\8\2\9\3\7\4\eeee
} 

\baaa
9-93
\eaaa
\bbbb
0&0&0&0&0&0&0&2&2\\
0&0&0&0&0&0&0&2&2\\
0&0&0&0&0&1&1&0&2\\
0&0&0&0&0&1&1&0&2\\
0&0&0&0&0&1&1&2&0\\
0&0&1&1&2&0&0&0&0\\
0&0&1&1&2&0&0&0&0\\
1&1&0&0&2&0&0&0&0\\
1&1&1&1&0&0&0&0&0\\
\ebbb
\parbox{7cm}{ 
\1\8\5\6\3\9\2\8\5\7\4\9\1\8\5\6\eeee
\8\5\7\4\9\1\8\5\6\3\9\2\8\5\7\4\eeee
\5\6\3\9\2\8\5\7\4\9\1\8\5\6\3\9\eeee
\7\4\9\1\8\5\6\3\9\2\8\5\7\4\9\1\eeee
\3\9\2\8\5\7\4\9\1\8\5\6\3\9\2\8\eeee
\9\1\8\5\6\3\9\2\8\5\7\4\9\1\8\5\eeee
\2\8\5\7\4\9\1\8\5\6\3\9\2\8\5\7\eeee
\8\5\6\3\9\2\8\5\7\4\9\1\8\5\6\3\eeee
\5\7\4\9\1\8\5\6\3\9\2\8\5\7\4\9\eeee
\6\3\9\2\8\5\7\4\9\1\8\5\6\3\9\2\eeee
\4\9\1\8\5\6\3\9\2\8\5\7\4\9\1\8\eeee
\9\2\8\5\7\4\9\1\8\5\6\3\9\2\8\5\eeee
} 

\baaa
9-94
\eaaa
\bbbb
0&0&0&0&0&0&0&2&2\\
0&0&0&0&0&0&0&2&2\\
0&0&0&0&0&1&1&0&2\\
0&0&0&0&0&1&1&0&2\\
0&0&0&0&0&2&2&0&0\\
0&0&1&1&2&0&0&0&0\\
0&0&1&1&2&0&0&0&0\\
1&1&0&0&0&0&0&2&0\\
1&1&1&1&0&0&0&0&0\\
\ebbb
\parbox{7cm}{ 
\1\8\8\1\9\4\7\5\6\4\9\1\8\8\2\9\eeee
\8\8\2\9\3\6\5\7\3\9\2\8\8\1\9\4\eeee
\8\1\9\4\7\5\6\4\9\1\8\8\2\9\3\7\eeee
\2\9\3\6\5\7\3\9\2\8\8\1\9\4\6\5\eeee
\9\4\7\5\6\4\9\1\8\8\2\9\3\7\5\6\eeee
\3\6\5\7\3\9\2\8\8\1\9\4\6\5\7\4\eeee
\7\5\6\4\9\1\8\8\2\9\3\7\5\6\3\9\eeee
\5\7\3\9\2\8\8\1\9\4\6\5\7\4\9\1\eeee
\6\4\9\1\8\8\2\9\3\7\5\6\3\9\2\8\eeee
\3\9\2\8\8\1\9\4\6\5\7\4\9\1\8\8\eeee
\9\1\8\8\2\9\3\7\5\6\3\9\2\8\8\2\eeee
\2\8\8\1\9\4\6\5\7\4\9\1\8\8\1\9\eeee
} 

\baaa
9-95
\eaaa
\bbbb
0&0&0&0&0&0&0&2&2\\
0&0&0&0&0&0&0&2&2\\
0&0&0&0&0&1&1&0&2\\
0&0&0&0&0&1&1&0&2\\
0&0&0&0&0&2&2&0&0\\
0&0&1&1&2&0&0&0&0\\
0&0&1&1&2&0&0&0&0\\
2&2&0&0&0&0&0&0&0\\
1&1&1&1&0&0&0&0&0\\
\ebbb
\parbox{7cm}{ 
\1\8\2\9\3\7\5\6\3\9\2\8\1\9\4\6\eeee
\8\1\9\4\6\5\7\4\9\1\8\2\9\3\7\5\eeee
\2\9\3\7\5\6\3\9\2\8\1\9\4\6\5\7\eeee
\9\4\6\5\7\4\9\1\8\2\9\3\7\5\6\3\eeee
\3\7\5\6\3\9\2\8\1\9\4\6\5\7\4\9\eeee
\6\5\7\4\9\1\8\2\9\3\7\5\6\3\9\2\eeee
\5\6\3\9\2\8\1\9\4\6\5\7\4\9\1\8\eeee
\7\4\9\1\8\2\9\3\7\5\6\3\9\2\8\1\eeee
\3\9\2\8\1\9\4\6\5\7\4\9\1\8\2\9\eeee
\9\1\8\2\9\3\7\5\6\3\9\2\8\1\9\4\eeee
\2\8\1\9\4\6\5\7\4\9\1\8\2\9\3\7\eeee
\8\2\9\3\7\5\6\3\9\2\8\1\9\4\6\5\eeee
} 

\baaa
9-96
\eaaa
\bbbb
0&0&0&0&0&0&0&2&2\\
0&0&0&0&0&0&0&2&2\\
0&0&0&0&0&1&1&0&2\\
0&0&0&0&0&1&1&0&2\\
0&0&0&0&2&0&0&2&0\\
0&0&1&1&0&1&1&0&0\\
0&0&1&1&0&1&1&0&0\\
1&1&0&0&2&0&0&0&0\\
1&1&1&1&0&0&0&0&0\\
\ebbb
\parbox{7cm}{ 
\1\8\5\5\8\2\9\3\6\7\4\9\1\8\5\5\eeee
\8\5\5\8\1\9\4\7\6\3\9\2\8\5\5\8\eeee
\5\5\8\2\9\3\6\7\4\9\1\8\5\5\8\2\eeee
\5\8\1\9\4\7\6\3\9\2\8\5\5\8\1\9\eeee
\8\2\9\3\6\7\4\9\1\8\5\5\8\2\9\3\eeee
\1\9\4\7\6\3\9\2\8\5\5\8\1\9\4\7\eeee
\9\3\6\7\4\9\1\8\5\5\8\2\9\3\6\7\eeee
\4\7\6\3\9\2\8\5\5\8\1\9\4\7\6\3\eeee
\6\7\4\9\1\8\5\5\8\2\9\3\6\7\4\9\eeee
\6\3\9\2\8\5\5\8\1\9\4\7\6\3\9\2\eeee
\4\9\1\8\5\5\8\2\9\3\6\7\4\9\1\8\eeee
\9\2\8\5\5\8\1\9\4\7\6\3\9\2\8\5\eeee
} 

\baaa
9-97
\eaaa
\bbbb
0&0&0&0&0&0&0&2&2\\
0&0&0&0&0&0&0&2&2\\
0&0&0&0&0&1&1&0&2\\
0&0&0&0&0&1&1&0&2\\
0&0&0&0&2&0&0&2&0\\
0&0&2&2&0&0&0&0&0\\
0&0&2&2&0&0&0&0&0\\
1&1&0&0&2&0&0&0&0\\
1&1&1&1&0&0&0&0&0\\
\ebbb
\parbox{7cm}{ 
\1\8\5\5\8\2\9\3\6\4\9\1\8\5\5\8\eeee
\8\5\5\8\1\9\4\7\3\9\2\8\5\5\8\2\eeee
\5\5\8\2\9\3\6\4\9\1\8\5\5\8\1\9\eeee
\5\8\1\9\4\7\3\9\2\8\5\5\8\2\9\3\eeee
\8\2\9\3\6\4\9\1\8\5\5\8\1\9\4\6\eeee
\1\9\4\7\3\9\2\8\5\5\8\2\9\3\7\4\eeee
\9\3\6\4\9\1\8\5\5\8\1\9\4\6\3\9\eeee
\4\7\3\9\2\8\5\5\8\2\9\3\7\4\9\1\eeee
\6\4\9\1\8\5\5\8\1\9\4\6\3\9\2\8\eeee
\3\9\2\8\5\5\8\2\9\3\7\4\9\1\8\5\eeee
\9\1\8\5\5\8\1\9\4\6\3\9\2\8\5\5\eeee
\2\8\5\5\8\2\9\3\7\4\9\1\8\5\5\8\eeee
} 

\baaa
9-98
\eaaa
\bbbb
0&0&0&0&0&0&0&2&2\\
0&0&0&0&0&0&0&2&2\\
0&0&0&0&0&1&1&0&2\\
0&0&0&0&0&1&1&0&2\\
0&0&0&0&2&1&1&0&0\\
0&0&1&1&2&0&0&0&0\\
0&0&1&1&2&0&0&0&0\\
1&1&0&0&0&0&0&2&0\\
1&1&1&1&0&0&0&0&0\\
\ebbb
\parbox{7cm}{ 
\1\8\8\1\9\4\7\5\5\7\4\9\1\8\8\1\eeee
\8\8\2\9\3\6\5\5\6\3\9\2\8\8\2\9\eeee
\8\1\9\4\7\5\5\7\4\9\1\8\8\1\9\4\eeee
\2\9\3\6\5\5\6\3\9\2\8\8\2\9\3\6\eeee
\9\4\7\5\5\7\4\9\1\8\8\1\9\4\7\5\eeee
\3\6\5\5\6\3\9\2\8\8\2\9\3\6\5\5\eeee
\7\5\5\7\4\9\1\8\8\1\9\4\7\5\5\7\eeee
\5\5\6\3\9\2\8\8\2\9\3\6\5\5\6\3\eeee
\5\7\4\9\1\8\8\1\9\4\7\5\5\7\4\9\eeee
\6\3\9\2\8\8\2\9\3\6\5\5\6\3\9\2\eeee
\4\9\1\8\8\1\9\4\7\5\5\7\4\9\1\8\eeee
\9\2\8\8\2\9\3\6\5\5\6\3\9\2\8\8\eeee
} 

\baaa
9-99
\eaaa
\bbbb
0&0&0&0&0&0&0&2&2\\
0&0&0&0&0&0&0&2&2\\
0&0&0&0&0&1&1&0&2\\
0&0&0&0&0&1&1&0&2\\
0&0&0&0&2&1&1&0&0\\
0&0&1&1&2&0&0&0&0\\
0&0&1&1&2&0&0&0&0\\
2&2&0&0&0&0&0&0&0\\
1&1&1&1&0&0&0&0&0\\
\ebbb
\parbox{7cm}{ 
\1\8\2\9\3\7\5\5\6\4\9\1\8\2\9\3\eeee
\8\1\9\4\6\5\5\7\3\9\2\8\1\9\4\7\eeee
\2\9\3\7\5\5\6\4\9\1\8\2\9\3\6\5\eeee
\9\4\6\5\5\7\3\9\2\8\1\9\4\7\5\5\eeee
\3\7\5\5\6\4\9\1\8\2\9\3\6\5\5\6\eeee
\6\5\5\7\3\9\2\8\1\9\4\7\5\5\7\4\eeee
\5\5\6\4\9\1\8\2\9\3\6\5\5\6\3\9\eeee
\5\7\3\9\2\8\1\9\4\7\5\5\7\4\9\1\eeee
\6\4\9\1\8\2\9\3\6\5\5\6\3\9\2\8\eeee
\3\9\2\8\1\9\4\7\5\5\7\4\9\1\8\2\eeee
\9\1\8\2\9\3\6\5\5\6\3\9\2\8\1\9\eeee
\2\8\1\9\4\7\5\5\7\4\9\1\8\2\9\3\eeee
} 

\baaa
9-100
\eaaa
\bbbb
0&0&0&0&0&0&0&2&2\\
0&0&0&0&0&0&0&2&2\\
0&0&0&0&0&1&1&0&2\\
0&0&0&0&0&1&1&1&1\\
0&0&0&0&0&1&1&2&0\\
0&0&1&2&1&0&0&0&0\\
0&0&1&2&1&0&0&0&0\\
1&1&0&1&1&0&0&0&0\\
1&1&1&1&0&0&0&0&0\\
\ebbb
\parbox{7cm}{ 
\1\8\4\6\5\8\2\8\5\7\4\8\1\9\4\6\eeee
\8\5\7\4\8\1\9\4\6\3\9\2\9\3\7\4\eeee
\4\6\3\9\2\9\3\7\4\9\1\8\4\6\5\8\eeee
\7\4\9\1\8\4\6\5\8\2\8\5\7\4\8\1\eeee
\5\8\2\8\5\7\4\8\1\9\4\6\3\9\2\9\eeee
\8\1\9\4\6\3\9\2\9\3\7\4\9\1\8\4\eeee
\2\9\3\7\4\9\1\8\4\6\5\8\2\8\5\7\eeee
\8\4\6\5\8\2\8\5\7\4\8\1\9\4\6\3\eeee
\5\7\4\8\1\9\4\6\3\9\2\9\3\7\4\9\eeee
\6\3\9\2\9\3\7\4\9\1\8\4\6\5\8\2\eeee
\4\9\1\8\4\6\5\8\2\8\5\7\4\8\1\9\eeee
\8\2\8\5\7\4\8\1\9\4\6\3\9\2\9\3\eeee
} 

\baaa
9-101
\eaaa
\bbbb
0&0&0&0&0&0&0&2&2\\
0&0&0&0&0&0&0&2&2\\
0&0&0&0&0&1&1&0&2\\
0&0&0&0&0&2&2&0&0\\
0&0&0&0&2&0&0&2&0\\
0&0&2&2&0&0&0&0&0\\
0&0&2&2&0&0&0&0&0\\
1&1&0&0&2&0&0&0&0\\
1&1&2&0&0&0&0&0&0\\
\ebbb
\parbox{7cm}{ 
\1\8\5\5\8\2\9\3\6\4\7\3\9\2\8\5\eeee
\8\5\5\8\1\9\3\7\4\6\3\9\1\8\5\5\eeee
\5\5\8\2\9\3\6\4\7\3\9\2\8\5\5\8\eeee
\5\8\1\9\3\7\4\6\3\9\1\8\5\5\8\1\eeee
\8\2\9\3\6\4\7\3\9\2\8\5\5\8\2\9\eeee
\1\9\3\7\4\6\3\9\1\8\5\5\8\1\9\3\eeee
\9\3\6\4\7\3\9\2\8\5\5\8\2\9\3\7\eeee
\3\7\4\6\3\9\1\8\5\5\8\1\9\3\6\4\eeee
\6\4\7\3\9\2\8\5\5\8\2\9\3\7\4\6\eeee
\4\6\3\9\1\8\5\5\8\1\9\3\6\4\7\3\eeee
\7\3\9\2\8\5\5\8\2\9\3\7\4\6\3\9\eeee
\3\9\1\8\5\5\8\1\9\3\6\4\7\3\9\1\eeee
} 

\baaa
9-102
\eaaa
\bbbb
0&0&0&0&0&0&0&2&2\\
0&0&0&0&0&0&0&2&2\\
0&0&0&0&0&1&1&0&2\\
0&0&0&0&2&1&1&0&0\\
0&0&0&2&2&0&0&0&0\\
0&0&2&2&0&0&0&0&0\\
0&0&2&2&0&0&0&0&0\\
1&1&0&0&0&0&0&2&0\\
1&1&2&0&0&0&0&0&0\\
\ebbb
\parbox{7cm}{ 
\1\8\8\1\9\3\7\4\5\5\4\6\3\9\2\8\eeee
\8\8\2\9\3\6\4\5\5\4\7\3\9\1\8\8\eeee
\8\1\9\3\7\4\5\5\4\6\3\9\2\8\8\2\eeee
\2\9\3\6\4\5\5\4\7\3\9\1\8\8\1\9\eeee
\9\3\7\4\5\5\4\6\3\9\2\8\8\2\9\3\eeee
\3\6\4\5\5\4\7\3\9\1\8\8\1\9\3\7\eeee
\7\4\5\5\4\6\3\9\2\8\8\2\9\3\6\4\eeee
\4\5\5\4\7\3\9\1\8\8\1\9\3\7\4\5\eeee
\5\5\4\6\3\9\2\8\8\2\9\3\6\4\5\5\eeee
\5\4\7\3\9\1\8\8\1\9\3\7\4\5\5\4\eeee
\4\6\3\9\2\8\8\2\9\3\6\4\5\5\4\7\eeee
\7\3\9\1\8\8\1\9\3\7\4\5\5\4\6\3\eeee
} 

\baaa
9-103
\eaaa
\bbbb
0&0&0&0&0&0&0&2&2\\
0&0&0&0&0&0&0&2&2\\
0&0&0&0&0&1&1&0&2\\
0&0&0&1&1&0&0&2&0\\
0&0&0&1&1&0&0&2&0\\
0&0&2&0&0&1&1&0&0\\
0&0&2&0&0&1&1&0&0\\
1&1&0&1&1&0&0&0&0\\
1&1&2&0&0&0&0&0&0\\
\ebbb
\parbox{7cm}{ 
\1\8\5\4\8\2\9\3\6\7\3\9\1\8\5\4\eeee
\8\4\5\8\1\9\3\7\6\3\9\2\8\4\5\8\eeee
\5\4\8\2\9\3\6\7\3\9\1\8\5\4\8\2\eeee
\5\8\1\9\3\7\6\3\9\2\8\4\5\8\1\9\eeee
\8\2\9\3\6\7\3\9\1\8\5\4\8\2\9\3\eeee
\1\9\3\7\6\3\9\2\8\4\5\8\1\9\3\7\eeee
\9\3\6\7\3\9\1\8\5\4\8\2\9\3\6\7\eeee
\3\7\6\3\9\2\8\4\5\8\1\9\3\7\6\3\eeee
\6\7\3\9\1\8\5\4\8\2\9\3\6\7\3\9\eeee
\6\3\9\2\8\4\5\8\1\9\3\7\6\3\9\2\eeee
\3\9\1\8\5\4\8\2\9\3\6\7\3\9\1\8\eeee
\9\2\8\4\5\8\1\9\3\7\6\3\9\2\8\4\eeee
} 

\baaa
9-104
\eaaa
\bbbb
0&0&0&0&0&0&0&2&2\\
0&0&0&0&0&0&0&2&2\\
0&0&0&0&0&1&1&0&2\\
0&0&0&1&1&1&1&0&0\\
0&0&0&1&1&1&1&0&0\\
0&0&2&1&1&0&0&0&0\\
0&0&2&1&1&0&0&0&0\\
1&1&0&0&0&0&0&2&0\\
1&1&2&0&0&0&0&0&0\\
\ebbb
\parbox{7cm}{ 
\1\8\8\1\9\3\7\4\4\7\3\9\1\8\8\1\eeee
\8\8\2\9\3\6\5\5\6\3\9\2\8\8\2\9\eeee
\8\1\9\3\7\4\4\7\3\9\1\8\8\1\9\3\eeee
\2\9\3\6\5\5\6\3\9\2\8\8\2\9\3\6\eeee
\9\3\7\4\4\7\3\9\1\8\8\1\9\3\7\4\eeee
\3\6\5\5\6\3\9\2\8\8\2\9\3\6\5\5\eeee
\7\4\4\7\3\9\1\8\8\1\9\3\7\4\4\7\eeee
\5\5\6\3\9\2\8\8\2\9\3\6\5\5\6\3\eeee
\4\7\3\9\1\8\8\1\9\3\7\4\4\7\3\9\eeee
\6\3\9\2\8\8\2\9\3\6\5\5\6\3\9\2\eeee
\3\9\1\8\8\1\9\3\7\4\4\7\3\9\1\8\eeee
\9\2\8\8\2\9\3\6\5\5\6\3\9\2\8\8\eeee
} 

\baaa
9-105
\eaaa
\bbbb
0&0&0&0&0&0&0&2&2\\
0&0&0&0&0&0&0&2&2\\
0&0&0&0&0&1&1&0&2\\
0&0&0&1&1&1&1&0&0\\
0&0&0&1&1&1&1&0&0\\
0&0&2&1&1&0&0&0&0\\
0&0&2&1&1&0&0&0&0\\
2&2&0&0&0&0&0&0&0\\
1&1&2&0&0&0&0&0&0\\
\ebbb
\parbox{7cm}{ 
\1\8\2\9\3\7\5\4\6\3\9\1\8\2\9\3\eeee
\8\1\9\3\6\4\5\7\3\9\2\8\1\9\3\7\eeee
\2\9\3\7\5\4\6\3\9\1\8\2\9\3\6\5\eeee
\9\3\6\4\5\7\3\9\2\8\1\9\3\7\4\4\eeee
\3\7\5\4\6\3\9\1\8\2\9\3\6\5\5\6\eeee
\6\4\5\7\3\9\2\8\1\9\3\7\4\4\7\3\eeee
\5\4\6\3\9\1\8\2\9\3\6\5\5\6\3\9\eeee
\5\7\3\9\2\8\1\9\3\7\4\4\7\3\9\1\eeee
\6\3\9\1\8\2\9\3\6\5\5\6\3\9\2\8\eeee
\3\9\2\8\1\9\3\7\4\4\7\3\9\1\8\2\eeee
\9\1\8\2\9\3\6\5\5\6\3\9\2\8\1\9\eeee
\2\8\1\9\3\7\4\4\7\3\9\1\8\2\9\3\eeee
} 

\baaa
9-106
\eaaa
\bbbb
0&0&0&0&0&0&0&2&2\\
0&0&0&0&0&0&0&2&2\\
0&0&0&0&0&1&1&1&1\\
0&0&0&0&0&1&1&1&1\\
0&0&0&0&0&2&2&0&0\\
0&0&1&1&2&0&0&0&0\\
0&0&1&1&2&0&0&0&0\\
1&1&1&1&0&0&0&0&0\\
1&1&1&1&0&0&0&0&0\\
\ebbb
\parbox{7cm}{ 
\1\8\3\6\5\6\3\8\1\8\3\6\5\6\3\8\eeee
\8\2\9\4\7\5\7\4\9\2\9\4\7\5\7\4\eeee
\3\9\1\8\3\6\5\6\3\8\1\8\3\6\5\6\eeee
\6\4\8\2\9\4\7\5\7\4\9\2\9\4\7\5\eeee
\5\7\3\9\1\8\3\6\5\6\3\8\1\8\3\6\eeee
\6\5\6\4\8\2\9\4\7\5\7\4\9\2\9\4\eeee
\3\7\5\7\3\9\1\8\3\6\5\6\3\8\1\8\eeee
\8\4\6\5\6\4\8\2\9\4\7\5\7\4\9\2\eeee
\1\9\3\7\5\7\3\9\1\8\3\6\5\6\3\8\eeee
\8\2\8\4\6\5\6\4\8\2\9\4\7\5\7\4\eeee
\3\9\1\9\3\7\5\7\3\9\1\8\3\6\5\6\eeee
\6\4\8\2\8\4\6\5\6\4\8\2\9\4\7\5\eeee
} 

\baaa
9-107
\eaaa
\bbbb
0&0&0&0&0&0&0&2&2\\
0&0&0&0&0&0&0&2&2\\
0&0&0&0&0&1&1&1&1\\
0&0&0&0&0&1&1&1&1\\
0&0&0&0&2&1&1&0&0\\
0&0&1&1&2&0&0&0&0\\
0&0&1&1&2&0&0&0&0\\
1&1&1&1&0&0&0&0&0\\
1&1&1&1&0&0&0&0&0\\
\ebbb
\parbox{7cm}{ 
\1\8\3\6\5\5\6\3\8\1\8\3\6\5\5\6\eeee
\8\2\9\4\7\5\5\7\4\9\2\9\4\7\5\5\eeee
\3\9\1\8\3\6\5\5\6\3\8\1\8\3\6\5\eeee
\6\4\8\2\9\4\7\5\5\7\4\9\2\9\4\7\eeee
\5\7\3\9\1\8\3\6\5\5\6\3\8\1\8\3\eeee
\5\5\6\4\8\2\9\4\7\5\5\7\4\9\2\9\eeee
\6\5\5\7\3\9\1\8\3\6\5\5\6\3\8\1\eeee
\3\7\5\5\6\4\8\2\9\4\7\5\5\7\4\9\eeee
\8\4\6\5\5\7\3\9\1\8\3\6\5\5\6\3\eeee
\1\9\3\7\5\5\6\4\8\2\9\4\7\5\5\7\eeee
\8\2\8\4\6\5\5\7\3\9\1\8\3\6\5\5\eeee
\3\9\1\9\3\7\5\5\6\4\8\2\9\4\7\5\eeee
} 

\baaa
9-108
\eaaa
\bbbb
0&0&0&0&0&0&0&2&2\\
0&0&0&0&0&0&0&2&2\\
0&0&0&0&0&1&1&1&1\\
0&0&0&0&0&2&2&0&0\\
0&0&0&0&0&2&2&0&0\\
0&0&2&1&1&0&0&0&0\\
0&0&2&1&1&0&0&0&0\\
1&1&2&0&0&0&0&0&0\\
1&1&2&0&0&0&0&0&0\\
\ebbb
\parbox{7cm}{ 
\1\8\3\6\4\6\3\8\1\8\3\6\4\6\3\8\eeee
\8\2\9\3\7\5\7\3\9\2\9\3\7\5\7\3\eeee
\3\9\1\8\3\6\4\6\3\8\1\8\3\6\4\6\eeee
\6\3\8\2\9\3\7\5\7\3\9\2\9\3\7\5\eeee
\4\7\3\9\1\8\3\6\4\6\3\8\1\8\3\6\eeee
\6\5\6\3\8\2\9\3\7\5\7\3\9\2\9\3\eeee
\3\7\4\7\3\9\1\8\3\6\4\6\3\8\1\8\eeee
\8\3\6\5\6\3\8\2\9\3\7\5\7\3\9\2\eeee
\1\9\3\7\4\7\3\9\1\8\3\6\4\6\3\8\eeee
\8\2\8\3\6\5\6\3\8\2\9\3\7\5\7\3\eeee
\3\9\1\9\3\7\4\7\3\9\1\8\3\6\4\6\eeee
\6\3\8\2\8\3\6\5\6\3\8\2\9\3\7\5\eeee
} 

\baaa
9-109
\eaaa
\bbbb
0&0&0&0&0&0&0&2&2\\
0&0&0&0&0&0&0&2&2\\
0&0&0&0&0&1&1&1&1\\
0&0&0&1&1&1&1&0&0\\
0&0&0&1&1&1&1&0&0\\
0&0&2&1&1&0&0&0&0\\
0&0&2&1&1&0&0&0&0\\
1&1&2&0&0&0&0&0&0\\
1&1&2&0&0&0&0&0&0\\
\ebbb
\parbox{7cm}{ 
\1\8\3\6\4\4\6\3\8\1\8\3\6\4\4\6\eeee
\8\2\9\3\7\5\5\7\3\9\2\9\3\7\5\5\eeee
\3\9\1\8\3\6\4\4\6\3\8\1\8\3\6\4\eeee
\6\3\8\2\9\3\7\5\5\7\3\9\2\9\3\7\eeee
\4\7\3\9\1\8\3\6\4\4\6\3\8\1\8\3\eeee
\4\5\6\3\8\2\9\3\7\5\5\7\3\9\2\9\eeee
\6\5\4\7\3\9\1\8\3\6\4\4\6\3\8\1\eeee
\3\7\4\5\6\3\8\2\9\3\7\5\5\7\3\9\eeee
\8\3\6\5\4\7\3\9\1\8\3\6\4\4\6\3\eeee
\1\9\3\7\4\5\6\3\8\2\9\3\7\5\5\7\eeee
\8\2\8\3\6\5\4\7\3\9\1\8\3\6\4\4\eeee
\3\9\1\9\3\7\4\5\6\3\8\2\9\3\7\5\eeee
} 

\baaa
9-110
\eaaa
\bbbb
0&0&0&0&0&0&0&2&2\\
0&0&0&0&0&0&1&1&2\\
0&0&0&0&0&1&1&1&1\\
0&0&0&0&0&2&1&1&0\\
0&0&0&0&0&2&2&0&0\\
0&0&1&2&1&0&0&0&0\\
0&1&1&1&1&0&0&0&0\\
1&1&1&1&0&0&0&0&0\\
1&2&1&0&0&0&0&0&0\\
\ebbb
\parbox{7cm}{ 
\1\9\2\8\3\7\4\6\5\6\4\7\3\8\2\9\eeee
\8\2\9\1\9\2\8\3\7\4\6\5\6\4\7\3\eeee
\4\7\3\8\2\9\1\9\2\8\3\7\4\6\5\6\eeee
\6\5\6\4\7\3\8\2\9\1\9\2\8\3\7\4\eeee
\3\7\4\6\5\6\4\7\3\8\2\9\1\9\2\8\eeee
\9\2\8\3\7\4\6\5\6\4\7\3\8\2\9\1\eeee
\2\9\1\9\2\8\3\7\4\6\5\6\4\7\3\8\eeee
\7\3\8\2\9\1\9\2\8\3\7\4\6\5\6\4\eeee
\5\6\4\7\3\8\2\9\1\9\2\8\3\7\4\6\eeee
\7\4\6\5\6\4\7\3\8\2\9\1\9\2\8\3\eeee
\2\8\3\7\4\6\5\6\4\7\3\8\2\9\1\9\eeee
\9\1\9\2\8\3\7\4\6\5\6\4\7\3\8\2\eeee
} 

\baaa
9-111
\eaaa
\bbbb
0&0&0&0&0&0&0&2&2\\
0&0&0&0&0&0&1&1&2\\
0&0&0&0&0&1&1&1&1\\
0&0&0&0&1&1&1&1&0\\
0&0&0&1&1&1&1&0&0\\
0&0&1&1&1&1&0&0&0\\
0&1&1&1&1&0&0&0&0\\
1&1&1&1&0&0&0&0&0\\
1&2&1&0&0&0&0&0&0\\
\ebbb
\parbox{7cm}{ 
\1\9\2\8\3\7\4\6\5\5\6\4\7\3\8\2\eeee
\8\2\9\1\9\2\8\3\7\4\6\5\5\6\4\7\eeee
\4\7\3\8\2\9\1\9\2\8\3\7\4\6\5\5\eeee
\5\5\6\4\7\3\8\2\9\1\9\2\8\3\7\4\eeee
\7\4\6\5\5\6\4\7\3\8\2\9\1\9\2\8\eeee
\2\8\3\7\4\6\5\5\6\4\7\3\8\2\9\1\eeee
\9\1\9\2\8\3\7\4\6\5\5\6\4\7\3\8\eeee
\3\8\2\9\1\9\2\8\3\7\4\6\5\5\6\4\eeee
\6\4\7\3\8\2\9\1\9\2\8\3\7\4\6\5\eeee
\6\5\5\6\4\7\3\8\2\9\1\9\2\8\3\7\eeee
\3\7\4\6\5\5\6\4\7\3\8\2\9\1\9\2\eeee
\9\2\8\3\7\4\6\5\5\6\4\7\3\8\2\9\eeee
} 

\baaa
9-112
\eaaa
\bbbb
0&0&0&0&0&0&0&2&2\\
0&0&0&0&0&0&1&1&2\\
0&0&0&0&0&1&2&0&1\\
0&0&0&0&1&2&1&0&0\\
0&0&0&1&2&1&0&0&0\\
0&0&1&2&1&0&0&0&0\\
0&1&2&1&0&0&0&0&0\\
2&2&0&0&0&0&0&0&0\\
1&2&1&0&0&0&0&0&0\\
\ebbb
\parbox{7cm}{ 
\1\9\3\6\5\4\7\2\8\2\7\4\5\6\3\9\eeee
\8\2\7\4\5\6\3\9\1\9\3\6\5\4\7\2\eeee
\1\9\3\6\5\4\7\2\8\2\7\4\5\6\3\9\eeee
\8\2\7\4\5\6\3\9\1\9\3\6\5\4\7\2\eeee
\1\9\3\6\5\4\7\2\8\2\7\4\5\6\3\9\eeee
\8\2\7\4\5\6\3\9\1\9\3\6\5\4\7\2\eeee
\1\9\3\6\5\4\7\2\8\2\7\4\5\6\3\9\eeee
\8\2\7\4\5\6\3\9\1\9\3\6\5\4\7\2\eeee
\1\9\3\6\5\4\7\2\8\2\7\4\5\6\3\9\eeee
\8\2\7\4\5\6\3\9\1\9\3\6\5\4\7\2\eeee
\1\9\3\6\5\4\7\2\8\2\7\4\5\6\3\9\eeee
\8\2\7\4\5\6\3\9\1\9\3\6\5\4\7\2\eeee
} 

\baaa
9-113
\eaaa
\bbbb
0&0&0&0&0&0&0&2&2\\
0&0&0&0&0&0&2&0&2\\
0&0&0&0&0&0&2&2&0\\
0&0&0&0&1&1&0&1&1\\
0&0&0&1&0&1&1&0&1\\
0&0&0&1&1&0&1&1&0\\
0&1&1&0&1&1&0&0&0\\
1&0&1&1&0&1&0&0&0\\
1&1&0&1&1&0&0&0&0\\
\ebbb
\parbox{7cm}{ 
\1\8\6\4\8\3\7\5\6\7\2\9\4\5\9\1\eeee
\8\3\7\5\6\7\2\9\4\5\9\1\8\6\4\8\eeee
\6\7\2\9\4\5\9\1\8\6\4\8\3\7\5\6\eeee
\4\5\9\1\8\6\4\8\3\7\5\6\7\2\9\4\eeee
\8\6\4\8\3\7\5\6\7\2\9\4\5\9\1\8\eeee
\3\7\5\6\7\2\9\4\5\9\1\8\6\4\8\3\eeee
\7\2\9\4\5\9\1\8\6\4\8\3\7\5\6\7\eeee
\5\9\1\8\6\4\8\3\7\5\6\7\2\9\4\5\eeee
\6\4\8\3\7\5\6\7\2\9\4\5\9\1\8\6\eeee
\7\5\6\7\2\9\4\5\9\1\8\6\4\8\3\7\eeee
\2\9\4\5\9\1\8\6\4\8\3\7\5\6\7\2\eeee
\9\1\8\6\4\8\3\7\5\6\7\2\9\4\5\9\eeee
} 

\baaa
9-114
\eaaa
\bbbb
0&0&0&0&0&0&0&2&2\\
0&0&0&0&0&0&2&0&2\\
0&0&0&0&0&0&2&2&0\\
0&0&0&0&1&1&0&1&1\\
0&0&0&1&1&0&1&0&1\\
0&0&0&1&0&1&1&1&0\\
0&1&1&0&1&1&0&0&0\\
1&0&1&1&0&1&0&0&0\\
1&1&0&1&1&0&0&0&0\\
\ebbb
\parbox{7cm}{ 
\1\8\6\4\9\2\7\6\6\7\2\9\4\6\8\1\eeee
\8\3\7\5\5\7\3\8\4\5\9\1\8\6\4\9\eeee
\6\7\2\9\4\6\8\1\9\5\4\8\3\7\5\5\eeee
\4\5\9\1\8\6\4\9\2\7\6\6\7\2\9\4\eeee
\9\5\4\8\3\7\5\5\7\3\8\4\5\9\1\8\eeee
\2\7\6\6\7\2\9\4\6\8\1\9\5\4\8\3\eeee
\7\3\8\4\5\9\1\8\6\4\9\2\7\6\6\7\eeee
\6\8\1\9\5\4\8\3\7\5\5\7\3\8\4\5\eeee
\6\4\9\2\7\6\6\7\2\9\4\6\8\1\9\5\eeee
\7\5\5\7\3\8\4\5\9\1\8\6\4\9\2\7\eeee
\2\9\4\6\8\1\9\5\4\8\3\7\5\5\7\3\eeee
\9\1\8\6\4\9\2\7\6\6\7\2\9\4\6\8\eeee
} 

\baaa
9-115
\eaaa
\bbbb
0&0&0&0&0&0&0&2&2\\
0&0&0&0&0&0&2&0&2\\
0&0&0&0&0&1&1&1&1\\
0&0&0&0&0&1&1&1&1\\
0&0&0&0&0&2&1&1&0\\
0&0&1&1&2&0&0&0&0\\
0&1&1&1&1&0&0&0&0\\
1&0&1&1&1&0&0&0&0\\
1&1&1&1&0&0&0&0&0\\
\ebbb
\parbox{7cm}{ 
\1\9\2\9\1\9\2\9\1\9\2\9\1\9\2\9\eeee
\8\3\7\4\8\3\7\4\8\3\7\4\8\3\7\4\eeee
\5\6\5\6\5\6\5\6\5\6\5\6\5\6\5\6\eeee
\7\4\8\3\7\4\8\3\7\4\8\3\7\4\8\3\eeee
\2\9\1\9\2\9\1\9\2\9\1\9\2\9\1\9\eeee
\7\3\8\4\7\3\8\4\7\3\8\4\7\3\8\4\eeee
\5\6\5\6\5\6\5\6\5\6\5\6\5\6\5\6\eeee
\8\4\7\3\8\4\7\3\8\4\7\3\8\4\7\3\eeee
\1\9\2\9\1\9\2\9\1\9\2\9\1\9\2\9\eeee
\8\3\7\4\8\3\7\4\8\3\7\4\8\3\7\4\eeee
\5\6\5\6\5\6\5\6\5\6\5\6\5\6\5\6\eeee
\7\4\8\3\7\4\8\3\7\4\8\3\7\4\8\3\eeee
} 

\baaa
9-116
\eaaa
\bbbb
0&0&0&0&0&0&0&2&2\\
0&0&0&0&0&0&2&0&2\\
0&0&0&0&0&1&1&1&1\\
0&0&0&0&0&2&0&2&0\\
0&0&0&0&0&2&2&0&0\\
0&0&1&1&2&0&0&0&0\\
0&1&1&0&2&0&0&0&0\\
2&0&1&1&0&0&0&0&0\\
2&1&1&0&0&0&0&0&0\\
\ebbb
\parbox{7cm}{ 
\1\8\3\7\5\6\3\9\1\8\3\7\5\6\3\9\eeee
\8\4\6\5\6\4\8\1\8\4\6\5\6\4\8\1\eeee
\3\6\5\7\3\8\1\9\3\6\5\7\3\8\1\9\eeee
\7\5\7\2\9\1\9\2\7\5\7\2\9\1\9\2\eeee
\5\6\3\9\1\8\3\7\5\6\3\9\1\8\3\7\eeee
\6\4\8\1\8\4\6\5\6\4\8\1\8\4\6\5\eeee
\3\8\1\9\3\6\5\7\3\8\1\9\3\6\5\7\eeee
\9\1\9\2\7\5\7\2\9\1\9\2\7\5\7\2\eeee
\1\8\3\7\5\6\3\9\1\8\3\7\5\6\3\9\eeee
\8\4\6\5\6\4\8\1\8\4\6\5\6\4\8\1\eeee
\3\6\5\7\3\8\1\9\3\6\5\7\3\8\1\9\eeee
\7\5\7\2\9\1\9\2\7\5\7\2\9\1\9\2\eeee
} 

\baaa
9-117
\eaaa
\bbbb
0&0&0&0&0&0&0&2&2\\
0&0&0&0&0&0&2&0&2\\
0&0&0&0&0&1&1&1&1\\
0&0&0&0&0&2&0&2&0\\
0&0&0&0&0&2&2&0&0\\
0&0&2&1&1&0&0&0&0\\
0&1&2&0&1&0&0&0&0\\
1&0&2&1&0&0&0&0&0\\
1&1&2&0&0&0&0&0&0\\
\ebbb
\parbox{7cm}{ 
\1\8\3\9\1\8\3\9\1\8\3\9\1\8\3\9\eeee
\8\4\6\3\8\4\6\3\8\4\6\3\8\4\6\3\eeee
\3\6\5\7\3\6\5\7\3\6\5\7\3\6\5\7\eeee
\9\3\7\2\9\3\7\2\9\3\7\2\9\3\7\2\eeee
\1\8\3\9\1\8\3\9\1\8\3\9\1\8\3\9\eeee
\8\4\6\3\8\4\6\3\8\4\6\3\8\4\6\3\eeee
\3\6\5\7\3\6\5\7\3\6\5\7\3\6\5\7\eeee
\9\3\7\2\9\3\7\2\9\3\7\2\9\3\7\2\eeee
\1\8\3\9\1\8\3\9\1\8\3\9\1\8\3\9\eeee
\8\4\6\3\8\4\6\3\8\4\6\3\8\4\6\3\eeee
\3\6\5\7\3\6\5\7\3\6\5\7\3\6\5\7\eeee
\9\3\7\2\9\3\7\2\9\3\7\2\9\3\7\2\eeee
} 

\baaa
\phantom{9-11}\#
\eaaa
\mbox{}\phantom{\bbbb
0&0&0&0&0&0&0&0&0\\
\ebbb}
\parbox{7cm}{ 
\1\9\2\9\1\9\2\9\1\9\2\9\1\9\2\9\eeee
\8\3\7\3\8\3\7\3\8\3\7\3\8\3\7\3\eeee
\4\6\5\6\4\6\5\6\4\6\5\6\4\6\5\6\eeee
\8\3\7\3\8\3\7\3\8\3\7\3\8\3\7\3\eeee
\1\9\2\9\1\9\2\9\1\9\2\9\1\9\2\9\eeee
\8\3\7\3\8\3\7\3\8\3\7\3\8\3\7\3\eeee
\4\6\5\6\4\6\5\6\4\6\5\6\4\6\5\6\eeee
\8\3\7\3\8\3\7\3\8\3\7\3\8\3\7\3\eeee
\1\9\2\9\1\9\2\9\1\9\2\9\1\9\2\9\eeee
\8\3\7\3\8\3\7\3\8\3\7\3\8\3\7\3\eeee
\4\6\5\6\4\6\5\6\4\6\5\6\4\6\5\6\eeee
\8\3\7\3\8\3\7\3\8\3\7\3\8\3\7\3\eeee
} 

\baaa
9-118
\eaaa
\bbbb
0&0&0&0&0&0&0&2&2\\
0&0&0&0&0&0&2&0&2\\
0&0&0&0&0&1&1&1&1\\
0&0&0&0&1&0&1&1&1\\
0&0&0&1&1&1&0&1&0\\
0&0&1&0&1&1&1&0&0\\
0&1&1&1&0&1&0&0&0\\
1&0&1&1&1&0&0&0&0\\
1&1&1&1&0&0&0&0&0\\
\ebbb
\parbox{7cm}{ 
\1\9\2\9\1\9\2\9\1\9\2\9\1\9\2\9\eeee
\8\3\7\4\8\3\7\4\8\3\7\4\8\3\7\4\eeee
\5\6\6\5\5\6\6\5\5\6\6\5\5\6\6\5\eeee
\4\7\3\8\4\7\3\8\4\7\3\8\4\7\3\8\eeee
\9\2\9\1\9\2\9\1\9\2\9\1\9\2\9\1\eeee
\3\7\4\8\3\7\4\8\3\7\4\8\3\7\4\8\eeee
\6\6\5\5\6\6\5\5\6\6\5\5\6\6\5\5\eeee
\7\3\8\4\7\3\8\4\7\3\8\4\7\3\8\4\eeee
\2\9\1\9\2\9\1\9\2\9\1\9\2\9\1\9\eeee
\7\4\8\3\7\4\8\3\7\4\8\3\7\4\8\3\eeee
\6\5\5\6\6\5\5\6\6\5\5\6\6\5\5\6\eeee
\3\8\4\7\3\8\4\7\3\8\4\7\3\8\4\7\eeee
} 

\baaa
9-119
\eaaa
\bbbb
0&0&0&0&0&0&0&2&2\\
0&0&0&0&0&0&2&0&2\\
0&0&0&0&0&1&1&1&1\\
0&0&0&0&1&2&0&1&0\\
0&0&0&1&0&2&1&0&0\\
0&0&1&1&1&1&0&0&0\\
0&1&2&0&1&0&0&0&0\\
1&0&2&1&0&0&0&0&0\\
1&1&2&0&0&0&0&0&0\\
\ebbb
\parbox{7cm}{ 
\1\9\2\9\1\9\2\9\1\9\2\9\1\9\2\9\eeee
\8\3\7\3\8\3\7\3\8\3\7\3\8\3\7\3\eeee
\4\6\5\6\4\6\5\6\4\6\5\6\4\6\5\6\eeee
\5\6\4\6\5\6\4\6\5\6\4\6\5\6\4\6\eeee
\7\3\8\3\7\3\8\3\7\3\8\3\7\3\8\3\eeee
\2\9\1\9\2\9\1\9\2\9\1\9\2\9\1\9\eeee
\7\3\8\3\7\3\8\3\7\3\8\3\7\3\8\3\eeee
\5\6\4\6\5\6\4\6\5\6\4\6\5\6\4\6\eeee
\4\6\5\6\4\6\5\6\4\6\5\6\4\6\5\6\eeee
\8\3\7\3\8\3\7\3\8\3\7\3\8\3\7\3\eeee
\1\9\2\9\1\9\2\9\1\9\2\9\1\9\2\9\eeee
\8\3\7\3\8\3\7\3\8\3\7\3\8\3\7\3\eeee
} 

\baaa
9-120
\eaaa
\bbbb
0&0&0&0&0&0&0&2&2\\
0&0&0&0&0&0&2&0&2\\
0&0&0&0&0&1&1&1&1\\
0&0&0&1&0&2&0&1&0\\
0&0&0&0&1&2&1&0&0\\
0&0&1&1&1&1&0&0&0\\
0&1&2&0&1&0&0&0&0\\
1&0&2&1&0&0&0&0&0\\
1&1&2&0&0&0&0&0&0\\
\ebbb
\parbox{7cm}{ 
\1\9\2\9\1\9\2\9\1\9\2\9\1\9\2\9\eeee
\8\3\7\3\8\3\7\3\8\3\7\3\8\3\7\3\eeee
\4\6\5\6\4\6\5\6\4\6\5\6\4\6\5\6\eeee
\4\6\5\6\4\6\5\6\4\6\5\6\4\6\5\6\eeee
\8\3\7\3\8\3\7\3\8\3\7\3\8\3\7\3\eeee
\1\9\2\9\1\9\2\9\1\9\2\9\1\9\2\9\eeee
\8\3\7\3\8\3\7\3\8\3\7\3\8\3\7\3\eeee
\4\6\5\6\4\6\5\6\4\6\5\6\4\6\5\6\eeee
\4\6\5\6\4\6\5\6\4\6\5\6\4\6\5\6\eeee
\8\3\7\3\8\3\7\3\8\3\7\3\8\3\7\3\eeee
\1\9\2\9\1\9\2\9\1\9\2\9\1\9\2\9\eeee
\8\3\7\3\8\3\7\3\8\3\7\3\8\3\7\3\eeee
} 

\baaa
9-121
\eaaa
\bbbb
0&0&0&0&0&0&0&2&2\\
0&0&0&0&0&0&2&0&2\\
0&0&0&0&0&2&0&2&0\\
0&0&0&0&2&0&2&0&0\\
0&0&0&2&0&2&0&0&0\\
0&0&2&0&2&0&0&0&0\\
0&2&0&2&0&0&0&0&0\\
2&0&2&0&0&0&0&0&0\\
2&2&0&0&0&0&0&0&0\\
\ebbb
\parbox{7cm}{ 
\1\8\3\6\5\4\7\2\9\1\8\3\6\5\4\7\eeee
\8\3\6\5\4\7\2\9\1\8\3\6\5\4\7\2\eeee
\3\6\5\4\7\2\9\1\8\3\6\5\4\7\2\9\eeee
\6\5\4\7\2\9\1\8\3\6\5\4\7\2\9\1\eeee
\5\4\7\2\9\1\8\3\6\5\4\7\2\9\1\8\eeee
\4\7\2\9\1\8\3\6\5\4\7\2\9\1\8\3\eeee
\7\2\9\1\8\3\6\5\4\7\2\9\1\8\3\6\eeee
\2\9\1\8\3\6\5\4\7\2\9\1\8\3\6\5\eeee
\9\1\8\3\6\5\4\7\2\9\1\8\3\6\5\4\eeee
\1\8\3\6\5\4\7\2\9\1\8\3\6\5\4\7\eeee
\8\3\6\5\4\7\2\9\1\8\3\6\5\4\7\2\eeee
\3\6\5\4\7\2\9\1\8\3\6\5\4\7\2\9\eeee
} 

\baaa
9-122
\eaaa
\bbbb
0&0&0&0&0&0&0&2&2\\
0&0&0&0&0&0&2&0&2\\
0&0&0&0&0&2&0&2&0\\
0&0&0&0&2&0&2&0&0\\
0&0&0&2&2&0&0&0&0\\
0&0&2&0&0&2&0&0&0\\
0&2&0&2&0&0&0&0&0\\
2&0&2&0&0&0&0&0&0\\
2&2&0&0&0&0&0&0&0\\
\ebbb
\parbox{7cm}{ 
\1\8\3\6\6\3\8\1\9\2\7\4\5\5\4\7\eeee
\8\3\6\6\3\8\1\9\2\7\4\5\5\4\7\2\eeee
\3\6\6\3\8\1\9\2\7\4\5\5\4\7\2\9\eeee
\6\6\3\8\1\9\2\7\4\5\5\4\7\2\9\1\eeee
\6\3\8\1\9\2\7\4\5\5\4\7\2\9\1\8\eeee
\3\8\1\9\2\7\4\5\5\4\7\2\9\1\8\3\eeee
\8\1\9\2\7\4\5\5\4\7\2\9\1\8\3\6\eeee
\1\9\2\7\4\5\5\4\7\2\9\1\8\3\6\6\eeee
\9\2\7\4\5\5\4\7\2\9\1\8\3\6\6\3\eeee
\2\7\4\5\5\4\7\2\9\1\8\3\6\6\3\8\eeee
\7\4\5\5\4\7\2\9\1\8\3\6\6\3\8\1\eeee
\4\5\5\4\7\2\9\1\8\3\6\6\3\8\1\9\eeee
} 

\baaa
9-123
\eaaa
\bbbb
0&0&0&0&0&0&0&2&2\\
0&0&0&0&0&0&2&0&2\\
0&0&0&0&0&2&0&2&0\\
0&0&0&1&1&0&2&0&0\\
0&0&0&1&1&0&2&0&0\\
0&0&2&0&0&2&0&0&0\\
0&2&0&1&1&0&0&0&0\\
2&0&2&0&0&0&0&0&0\\
2&2&0&0&0&0&0&0&0\\
\ebbb
\parbox{7cm}{ 
\1\8\3\6\6\3\8\1\9\2\7\4\4\7\2\9\eeee
\8\3\6\6\3\8\1\9\2\7\5\5\7\2\9\1\eeee
\3\6\6\3\8\1\9\2\7\4\4\7\2\9\1\8\eeee
\6\6\3\8\1\9\2\7\5\5\7\2\9\1\8\3\eeee
\6\3\8\1\9\2\7\4\4\7\2\9\1\8\3\6\eeee
\3\8\1\9\2\7\5\5\7\2\9\1\8\3\6\6\eeee
\8\1\9\2\7\4\4\7\2\9\1\8\3\6\6\3\eeee
\1\9\2\7\5\5\7\2\9\1\8\3\6\6\3\8\eeee
\9\2\7\4\4\7\2\9\1\8\3\6\6\3\8\1\eeee
\2\7\5\5\7\2\9\1\8\3\6\6\3\8\1\9\eeee
\7\4\4\7\2\9\1\8\3\6\6\3\8\1\9\2\eeee
\5\5\7\2\9\1\8\3\6\6\3\8\1\9\2\7\eeee
} 

\baaa
9-124
\eaaa
\bbbb
0&0&0&0&0&0&0&2&2\\
0&0&0&0&0&0&2&0&2\\
0&0&0&0&1&1&0&2&0\\
0&0&0&0&1&1&0&2&0\\
0&0&1&1&0&0&2&0&0\\
0&0&1&1&0&0&2&0&0\\
0&2&0&0&1&1&0&0&0\\
2&0&1&1&0&0&0&0&0\\
2&2&0&0&0&0&0&0&0\\
\ebbb
\parbox{7cm}{ 
\1\8\4\5\7\2\9\1\8\4\6\7\2\9\1\8\eeee
\8\3\6\7\2\9\1\8\3\5\7\2\9\1\8\4\eeee
\4\5\7\2\9\1\8\4\6\7\2\9\1\8\3\6\eeee
\6\7\2\9\1\8\3\5\7\2\9\1\8\4\5\7\eeee
\7\2\9\1\8\4\6\7\2\9\1\8\3\6\7\2\eeee
\2\9\1\8\3\5\7\2\9\1\8\4\5\7\2\9\eeee
\9\1\8\4\6\7\2\9\1\8\3\6\7\2\9\1\eeee
\1\8\3\5\7\2\9\1\8\4\5\7\2\9\1\8\eeee
\8\4\6\7\2\9\1\8\3\6\7\2\9\1\8\3\eeee
\3\5\7\2\9\1\8\4\5\7\2\9\1\8\4\6\eeee
\6\7\2\9\1\8\3\6\7\2\9\1\8\3\5\7\eeee
\7\2\9\1\8\4\5\7\2\9\1\8\4\6\7\2\eeee
} 

\baaa
9-125
\eaaa
\bbbb
0&0&0&0&0&0&0&2&2\\
0&0&0&0&0&0&2&0&2\\
0&0&0&0&1&1&0&2&0\\
0&0&0&0&1&1&0&2&0\\
0&0&1&1&1&1&0&0&0\\
0&0&1&1&1&1&0&0&0\\
0&2&0&0&0&0&2&0&0\\
2&0&1&1&0&0&0&0&0\\
2&2&0&0&0&0&0&0&0\\
\ebbb
\parbox{7cm}{ 
\1\8\4\5\5\4\8\1\9\2\7\7\2\9\1\8\eeee
\8\3\6\6\3\8\1\9\2\7\7\2\9\1\8\4\eeee
\4\5\5\4\8\1\9\2\7\7\2\9\1\8\3\6\eeee
\6\6\3\8\1\9\2\7\7\2\9\1\8\4\5\5\eeee
\5\4\8\1\9\2\7\7\2\9\1\8\3\6\6\3\eeee
\3\8\1\9\2\7\7\2\9\1\8\4\5\5\4\8\eeee
\8\1\9\2\7\7\2\9\1\8\3\6\6\3\8\1\eeee
\1\9\2\7\7\2\9\1\8\4\5\5\4\8\1\9\eeee
\9\2\7\7\2\9\1\8\3\6\6\3\8\1\9\2\eeee
\2\7\7\2\9\1\8\4\5\5\4\8\1\9\2\7\eeee
\7\7\2\9\1\8\3\6\6\3\8\1\9\2\7\7\eeee
\7\2\9\1\8\4\5\5\4\8\1\9\2\7\7\2\eeee
} 

\baaa
9-126
\eaaa
\bbbb
0&0&0&0&0&0&0&2&2\\
0&0&0&0&0&0&2&0&2\\
0&0&1&0&0&1&0&2&0\\
0&0&0&1&1&0&2&0&0\\
0&0&0&1&1&0&2&0&0\\
0&0&1&0&0&1&0&2&0\\
0&2&0&1&1&0&0&0&0\\
2&0&1&0&0&1&0&0&0\\
2&2&0&0&0&0&0&0&0\\
\ebbb
\parbox{7cm}{ 
\1\8\6\3\8\1\9\2\7\5\5\7\2\9\1\8\eeee
\8\3\6\8\1\9\2\7\4\4\7\2\9\1\8\6\eeee
\6\3\8\1\9\2\7\5\5\7\2\9\1\8\3\6\eeee
\6\8\1\9\2\7\4\4\7\2\9\1\8\6\3\8\eeee
\8\1\9\2\7\5\5\7\2\9\1\8\3\6\8\1\eeee
\1\9\2\7\4\4\7\2\9\1\8\6\3\8\1\9\eeee
\9\2\7\5\5\7\2\9\1\8\3\6\8\1\9\2\eeee
\2\7\4\4\7\2\9\1\8\6\3\8\1\9\2\7\eeee
\7\5\5\7\2\9\1\8\3\6\8\1\9\2\7\4\eeee
\4\4\7\2\9\1\8\6\3\8\1\9\2\7\5\5\eeee
\5\7\2\9\1\8\3\6\8\1\9\2\7\4\4\7\eeee
\7\2\9\1\8\6\3\8\1\9\2\7\5\5\7\2\eeee
} 

\baaa
9-127
\eaaa
\bbbb
0&0&0&0&0&0&0&2&2\\
0&0&0&0&0&0&2&1&1\\
0&0&0&0&0&0&2&1&1\\
0&0&0&0&0&1&2&0&1\\
0&0&0&0&0&1&2&1&0\\
0&0&0&2&2&0&0&0&0\\
0&1&1&1&1&0&0&0&0\\
1&1&1&0&1&0&0&0&0\\
1&1&1&1&0&0&0&0&0\\
\ebbb
\parbox{7cm}{ 
\1\9\4\6\4\9\1\9\4\6\4\9\1\9\4\6\eeee
\8\2\7\5\7\3\8\2\7\5\7\3\8\2\7\5\eeee
\5\7\3\8\2\7\5\7\3\8\2\7\5\7\3\8\eeee
\6\4\9\1\9\4\6\4\9\1\9\4\6\4\9\1\eeee
\5\7\2\8\3\7\5\7\2\8\3\7\5\7\2\8\eeee
\8\3\7\5\7\2\8\3\7\5\7\2\8\3\7\5\eeee
\1\9\4\6\4\9\1\9\4\6\4\9\1\9\4\6\eeee
\8\2\7\5\7\3\8\2\7\5\7\3\8\2\7\5\eeee
\5\7\3\8\2\7\5\7\3\8\2\7\5\7\3\8\eeee
\6\4\9\1\9\4\6\4\9\1\9\4\6\4\9\1\eeee
\5\7\2\8\3\7\5\7\2\8\3\7\5\7\2\8\eeee
\8\3\7\5\7\2\8\3\7\5\7\2\8\3\7\5\eeee
} 

\baaa
9-128
\eaaa
\bbbb
0&0&0&0&0&0&0&2&2\\
0&0&0&0&0&0&2&1&1\\
0&0&0&0&0&0&2&1&1\\
0&0&0&0&0&2&2&0&0\\
0&0&0&0&2&2&0&0&0\\
0&0&0&2&2&0&0&0&0\\
0&1&1&2&0&0&0&0&0\\
2&1&1&0&0&0&0&0&0\\
2&1&1&0&0&0&0&0&0\\
\ebbb
\parbox{7cm}{ 
\1\8\2\7\4\6\5\5\6\4\7\2\8\1\8\2\eeee
\8\1\9\3\7\4\6\5\5\6\4\7\3\9\1\9\eeee
\2\9\1\8\2\7\4\6\5\5\6\4\7\2\8\1\eeee
\7\3\8\1\9\3\7\4\6\5\5\6\4\7\3\9\eeee
\4\7\2\9\1\8\2\7\4\6\5\5\6\4\7\2\eeee
\6\4\7\3\8\1\9\3\7\4\6\5\5\6\4\7\eeee
\5\6\4\7\2\9\1\8\2\7\4\6\5\5\6\4\eeee
\5\5\6\4\7\3\8\1\9\3\7\4\6\5\5\6\eeee
\6\5\5\6\4\7\2\9\1\8\2\7\4\6\5\5\eeee
\4\6\5\5\6\4\7\3\8\1\9\3\7\4\6\5\eeee
\7\4\6\5\5\6\4\7\2\9\1\8\2\7\4\6\eeee
\2\7\4\6\5\5\6\4\7\3\8\1\9\3\7\4\eeee
} 

\baaa
9-129
\eaaa
\bbbb
0&0&0&0&0&0&0&2&2\\
0&0&0&0&0&0&2&1&1\\
0&0&0&0&0&0&2&1&1\\
0&0&0&0&1&1&2&0&0\\
0&0&0&2&1&1&0&0&0\\
0&0&0&2&1&1&0&0&0\\
0&1&1&2&0&0&0&0&0\\
2&1&1&0&0&0&0&0&0\\
2&1&1&0&0&0&0&0&0\\
\ebbb
\parbox{7cm}{ 
\1\8\2\7\4\5\5\4\7\2\8\1\8\2\7\4\eeee
\8\1\9\3\7\4\6\6\4\7\3\9\1\9\3\7\eeee
\2\9\1\8\2\7\4\5\5\4\7\2\8\1\8\2\eeee
\7\3\8\1\9\3\7\4\6\6\4\7\3\9\1\9\eeee
\4\7\2\9\1\8\2\7\4\5\5\4\7\2\8\1\eeee
\5\4\7\3\8\1\9\3\7\4\6\6\4\7\3\9\eeee
\5\6\4\7\2\9\1\8\2\7\4\5\5\4\7\2\eeee
\4\6\5\4\7\3\8\1\9\3\7\4\6\6\4\7\eeee
\7\4\5\6\4\7\2\9\1\8\2\7\4\5\5\4\eeee
\2\7\4\6\5\4\7\3\8\1\9\3\7\4\6\6\eeee
\8\3\7\4\5\6\4\7\2\9\1\8\2\7\4\5\eeee
\1\9\2\7\4\6\5\4\7\3\8\1\9\3\7\4\eeee
} 

\baaa
9-130
\eaaa
\bbbb
0&0&0&0&0&0&0&2&2\\
0&0&0&0&0&1&1&0&2\\
0&0&0&0&0&1&1&0&2\\
0&0&0&0&0&1&1&2&0\\
0&0&0&0&0&1&1&2&0\\
0&1&1&1&1&0&0&0&0\\
0&1&1&1&1&0&0&0&0\\
2&0&0&1&1&0&0&0&0\\
2&1&1&0&0&0&0&0&0\\
\ebbb
\parbox{7cm}{ 
\1\8\5\6\2\9\1\8\4\7\3\9\1\8\5\6\eeee
\8\4\7\3\9\1\8\5\6\2\9\1\8\4\7\3\eeee
\5\6\2\9\1\8\4\7\3\9\1\8\5\6\2\9\eeee
\7\3\9\1\8\5\6\2\9\1\8\4\7\3\9\1\eeee
\2\9\1\8\4\7\3\9\1\8\5\6\2\9\1\8\eeee
\9\1\8\5\6\2\9\1\8\4\7\3\9\1\8\5\eeee
\1\8\4\7\3\9\1\8\5\6\2\9\1\8\4\7\eeee
\8\5\6\2\9\1\8\4\7\3\9\1\8\5\6\2\eeee
\4\7\3\9\1\8\5\6\2\9\1\8\4\7\3\9\eeee
\6\2\9\1\8\4\7\3\9\1\8\5\6\2\9\1\eeee
\3\9\1\8\5\6\2\9\1\8\4\7\3\9\1\8\eeee
\9\1\8\4\7\3\9\1\8\5\6\2\9\1\8\4\eeee
} 

\baaa
9-131
\eaaa
\bbbb
0&0&0&0&0&0&0&2&2\\
0&0&0&0&0&1&1&0&2\\
0&0&0&0&0&1&1&0&2\\
0&0&0&0&2&1&1&0&0\\
0&0&0&2&2&0&0&0&0\\
0&1&1&2&0&0&0&0&0\\
0&1&1&2&0&0&0&0&0\\
2&0&0&0&0&0&0&2&0\\
2&1&1&0&0&0&0&0&0\\
\ebbb
\parbox{7cm}{ 
\1\8\8\1\9\3\7\4\5\5\4\6\2\9\1\8\eeee
\8\8\1\9\2\6\4\5\5\4\7\3\9\1\8\8\eeee
\8\1\9\3\7\4\5\5\4\6\2\9\1\8\8\1\eeee
\1\9\2\6\4\5\5\4\7\3\9\1\8\8\1\9\eeee
\9\3\7\4\5\5\4\6\2\9\1\8\8\1\9\2\eeee
\2\6\4\5\5\4\7\3\9\1\8\8\1\9\3\7\eeee
\7\4\5\5\4\6\2\9\1\8\8\1\9\2\6\4\eeee
\4\5\5\4\7\3\9\1\8\8\1\9\3\7\4\5\eeee
\5\5\4\6\2\9\1\8\8\1\9\2\6\4\5\5\eeee
\5\4\7\3\9\1\8\8\1\9\3\7\4\5\5\4\eeee
\4\6\2\9\1\8\8\1\9\2\6\4\5\5\4\7\eeee
\7\3\9\1\8\8\1\9\3\7\4\5\5\4\6\2\eeee
} 

\baaa
9-132
\eaaa
\bbbb
0&0&0&0&0&0&0&2&2\\
0&0&0&0&0&1&1&0&2\\
0&0&0&0&0&1&1&0&2\\
0&0&0&1&1&0&0&2&0\\
0&0&0&1&1&0&0&2&0\\
0&1&1&0&0&1&1&0&0\\
0&1&1&0&0&1&1&0&0\\
2&0&0&1&1&0&0&0&0\\
2&1&1&0&0&0&0&0&0\\
\ebbb
\parbox{7cm}{ 
\1\8\5\4\8\1\9\2\6\7\3\9\1\8\5\4\eeee
\8\4\5\8\1\9\3\7\6\2\9\1\8\4\5\8\eeee
\5\4\8\1\9\2\6\7\3\9\1\8\5\4\8\1\eeee
\5\8\1\9\3\7\6\2\9\1\8\4\5\8\1\9\eeee
\8\1\9\2\6\7\3\9\1\8\5\4\8\1\9\2\eeee
\1\9\3\7\6\2\9\1\8\4\5\8\1\9\3\7\eeee
\9\2\6\7\3\9\1\8\5\4\8\1\9\2\6\7\eeee
\3\7\6\2\9\1\8\4\5\8\1\9\3\7\6\2\eeee
\6\7\3\9\1\8\5\4\8\1\9\2\6\7\3\9\eeee
\6\2\9\1\8\4\5\8\1\9\3\7\6\2\9\1\eeee
\3\9\1\8\5\4\8\1\9\2\6\7\3\9\1\8\eeee
\9\1\8\4\5\8\1\9\3\7\6\2\9\1\8\4\eeee
} 

\baaa
9-133
\eaaa
\bbbb
0&0&0&0&0&0&0&2&2\\
0&0&0&0&0&1&1&0&2\\
0&0&0&0&0&1&1&0&2\\
0&0&0&1&1&1&1&0&0\\
0&0&0&1&1&1&1&0&0\\
0&1&1&1&1&0&0&0&0\\
0&1&1&1&1&0&0&0&0\\
2&0&0&0&0&0&0&2&0\\
2&1&1&0&0&0&0&0&0\\
\ebbb
\parbox{7cm}{ 
\1\8\8\1\9\3\7\4\4\7\3\9\1\8\8\1\eeee
\8\8\1\9\2\6\5\5\6\2\9\1\8\8\1\9\eeee
\8\1\9\3\7\4\4\7\3\9\1\8\8\1\9\3\eeee
\1\9\2\6\5\5\6\2\9\1\8\8\1\9\2\6\eeee
\9\3\7\4\4\7\3\9\1\8\8\1\9\3\7\4\eeee
\2\6\5\5\6\2\9\1\8\8\1\9\2\6\5\5\eeee
\7\4\4\7\3\9\1\8\8\1\9\3\7\4\4\7\eeee
\5\5\6\2\9\1\8\8\1\9\2\6\5\5\6\2\eeee
\4\7\3\9\1\8\8\1\9\3\7\4\4\7\3\9\eeee
\6\2\9\1\8\8\1\9\2\6\5\5\6\2\9\1\eeee
\3\9\1\8\8\1\9\3\7\4\4\7\3\9\1\8\eeee
\9\1\8\8\1\9\2\6\5\5\6\2\9\1\8\8\eeee
} 

\baaa
9-134
\eaaa
\bbbb
0&0&0&0&0&0&0&2&2\\
0&0&0&0&0&1&1&0&2\\
0&0&0&0&0&1&1&0&2\\
0&0&0&2&0&0&0&2&0\\
0&0&0&0&2&1&1&0&0\\
0&1&1&0&2&0&0&0&0\\
0&1&1&0&2&0&0&0&0\\
2&0&0&2&0&0&0&0&0\\
2&1&1&0&0&0&0&0&0\\
\ebbb
\parbox{7cm}{ 
\1\8\4\4\8\1\9\2\6\5\5\6\2\9\1\8\eeee
\8\4\4\8\1\9\3\7\5\5\7\3\9\1\8\4\eeee
\4\4\8\1\9\2\6\5\5\6\2\9\1\8\4\4\eeee
\4\8\1\9\3\7\5\5\7\3\9\1\8\4\4\8\eeee
\8\1\9\2\6\5\5\6\2\9\1\8\4\4\8\1\eeee
\1\9\3\7\5\5\7\3\9\1\8\4\4\8\1\9\eeee
\9\2\6\5\5\6\2\9\1\8\4\4\8\1\9\3\eeee
\3\7\5\5\7\3\9\1\8\4\4\8\1\9\2\6\eeee
\6\5\5\6\2\9\1\8\4\4\8\1\9\3\7\5\eeee
\5\5\7\3\9\1\8\4\4\8\1\9\2\6\5\5\eeee
\5\6\2\9\1\8\4\4\8\1\9\3\7\5\5\7\eeee
\7\3\9\1\8\4\4\8\1\9\2\6\5\5\6\2\eeee
} 

\baaa
9-135
\eaaa
\bbbb
0&0&0&0&0&0&0&2&2\\
0&0&0&0&0&1&1&1&1\\
0&0&0&0&0&1&1&1&1\\
0&0&0&0&0&1&1&1&1\\
0&0&0&0&0&2&2&0&0\\
0&1&1&1&1&0&0&0&0\\
0&1&1&1&1&0&0&0&0\\
1&1&1&1&0&0&0&0&0\\
1&1&1&1&0&0&0&0&0\\
\ebbb
\parbox{7cm}{ 
\1\9\3\8\1\8\3\9\1\9\3\8\1\8\3\9\eeee
\8\2\7\4\9\2\6\4\8\2\7\4\9\2\6\4\eeee
\3\6\5\6\3\7\5\7\3\6\5\6\3\7\5\7\eeee
\9\4\7\2\8\4\6\2\9\4\7\2\8\4\6\2\eeee
\1\8\3\9\1\9\3\8\1\8\3\9\1\9\3\8\eeee
\9\2\6\4\8\2\7\4\9\2\6\4\8\2\7\4\eeee
\3\7\5\7\3\6\5\6\3\7\5\7\3\6\5\6\eeee
\8\4\6\2\9\4\7\2\8\4\6\2\9\4\7\2\eeee
\1\9\3\8\1\8\3\9\1\9\3\8\1\8\3\9\eeee
\8\2\7\4\9\2\6\4\8\2\7\4\9\2\6\4\eeee
\3\6\5\6\3\7\5\7\3\6\5\6\3\7\5\7\eeee
\9\4\7\2\8\4\6\2\9\4\7\2\8\4\6\2\eeee
} 

\baaa
9-136
\eaaa
\bbbb
0&0&0&0&0&0&0&2&2\\
0&0&0&0&0&1&1&1&1\\
0&0&0&0&0&1&1&1&1\\
0&0&0&0&2&1&1&0&0\\
0&0&0&2&2&0&0&0&0\\
0&1&1&2&0&0&0&0&0\\
0&1&1&2&0&0&0&0&0\\
2&1&1&0&0&0&0&0&0\\
2&1&1&0&0&0&0&0&0\\
\ebbb
\parbox{7cm}{ 
\1\8\2\6\4\5\5\4\6\2\8\1\8\2\6\4\eeee
\8\1\9\3\7\4\5\5\4\7\3\9\1\9\3\7\eeee
\2\9\1\8\2\6\4\5\5\4\6\2\8\1\8\2\eeee
\6\3\8\1\9\3\7\4\5\5\4\7\3\9\1\9\eeee
\4\7\2\9\1\8\2\6\4\5\5\4\6\2\8\1\eeee
\5\4\6\3\8\1\9\3\7\4\5\5\4\7\3\9\eeee
\5\5\4\7\2\9\1\8\2\6\4\5\5\4\6\2\eeee
\4\5\5\4\6\3\8\1\9\3\7\4\5\5\4\7\eeee
\6\4\5\5\4\7\2\9\1\8\2\6\4\5\5\4\eeee
\2\7\4\5\5\4\6\3\8\1\9\3\7\4\5\5\eeee
\8\3\6\4\5\5\4\7\2\9\1\8\2\6\4\5\eeee
\1\9\2\7\4\5\5\4\6\3\8\1\9\3\7\4\eeee
} 

\baaa
9-137
\eaaa
\bbbb
0&0&0&0&0&0&0&2&2\\
0&0&0&0&0&1&1&1&1\\
0&0&0&0&0&1&1&1&1\\
0&0&0&1&1&1&1&0&0\\
0&0&0&1&1&1&1&0&0\\
0&1&1&1&1&0&0&0&0\\
0&1&1&1&1&0&0&0&0\\
2&1&1&0&0&0&0&0&0\\
2&1&1&0&0&0&0&0&0\\
\ebbb
\parbox{7cm}{ 
\1\8\2\6\4\4\6\2\8\1\8\2\6\4\4\6\eeee
\8\1\9\3\7\5\5\7\3\9\1\9\3\7\5\5\eeee
\2\9\1\8\2\6\4\4\6\2\8\1\8\2\6\4\eeee
\6\3\8\1\9\3\7\5\5\7\3\9\1\9\3\7\eeee
\4\7\2\9\1\8\2\6\4\4\6\2\8\1\8\2\eeee
\4\5\6\3\8\1\9\3\7\5\5\7\3\9\1\9\eeee
\6\5\4\7\2\9\1\8\2\6\4\4\6\2\8\1\eeee
\2\7\4\5\6\3\8\1\9\3\7\5\5\7\3\9\eeee
\8\3\6\5\4\7\2\9\1\8\2\6\4\4\6\2\eeee
\1\9\2\7\4\5\6\3\8\1\9\3\7\5\5\7\eeee
\8\1\8\3\6\5\4\7\2\9\1\8\2\6\4\4\eeee
\2\9\1\9\2\7\4\5\6\3\8\1\9\3\7\5\eeee
} 

\baaa
9-138
\eaaa
\bbbb
0&0&0&0&0&0&0&2&2\\
0&0&0&0&0&1&1&1&1\\
0&0&0&0&1&0&2&0&1\\
0&0&0&1&1&0&0&1&1\\
0&0&1&1&0&1&0&1&0\\
0&1&0&0&1&1&1&0&0\\
0&1&2&0&0&1&0&0&0\\
1&1&0&1&1&0&0&0&0\\
1&1&1&1&0&0&0&0&0\\
\ebbb
\parbox{7cm}{ 
\1\9\3\7\2\8\4\5\6\6\5\4\8\2\7\3\eeee
\8\2\7\3\9\1\9\3\7\2\8\4\5\6\6\5\eeee
\5\6\6\5\4\8\2\7\3\9\1\9\3\7\2\8\eeee
\3\7\2\8\4\5\6\6\5\4\8\2\7\3\9\1\eeee
\7\3\9\1\9\3\7\2\8\4\5\6\6\5\4\8\eeee
\6\5\4\8\2\7\3\9\1\9\3\7\2\8\4\5\eeee
\2\8\4\5\6\6\5\4\8\2\7\3\9\1\9\3\eeee
\9\1\9\3\7\2\8\4\5\6\6\5\4\8\2\7\eeee
\4\8\2\7\3\9\1\9\3\7\2\8\4\5\6\6\eeee
\4\5\6\6\5\4\8\2\7\3\9\1\9\3\7\2\eeee
\9\3\7\2\8\4\5\6\6\5\4\8\2\7\3\9\eeee
\2\7\3\9\1\9\3\7\2\8\4\5\6\6\5\4\eeee
} 

\baaa
9-139
\eaaa
\bbbb
0&0&0&0&0&0&0&2&2\\
0&0&0&0&0&1&1&1&1\\
0&0&0&0&1&0&2&0&1\\
0&0&0&2&0&1&1&0&0\\
0&0&1&0&0&2&0&1&0\\
0&1&0&1&1&0&1&0&0\\
0&1&1&1&0&1&0&0&0\\
1&2&0&0&1&0&0&0&0\\
1&2&1&0&0&0&0&0&0\\
\ebbb
\parbox{7cm}{ 
\1\9\3\5\8\1\8\5\3\9\1\9\3\5\8\1\eeee
\8\2\7\6\2\9\2\6\7\2\8\2\7\6\2\9\eeee
\5\6\4\4\7\3\7\4\4\6\5\6\4\4\7\3\eeee
\3\7\4\4\6\5\6\4\4\7\3\7\4\4\6\5\eeee
\9\2\6\7\2\8\2\7\6\2\9\2\6\7\2\8\eeee
\1\8\5\3\9\1\9\3\5\8\1\8\5\3\9\1\eeee
\9\2\6\7\2\8\2\7\6\2\9\2\6\7\2\8\eeee
\3\7\4\4\6\5\6\4\4\7\3\7\4\4\6\5\eeee
\5\6\4\4\7\3\7\4\4\6\5\6\4\4\7\3\eeee
\8\2\7\6\2\9\2\6\7\2\8\2\7\6\2\9\eeee
\1\9\3\5\8\1\8\5\3\9\1\9\3\5\8\1\eeee
\8\2\7\6\2\9\2\6\7\2\8\2\7\6\2\9\eeee
} 

\baaa
9-140
\eaaa
\bbbb
0&0&0&0&0&0&0&2&2\\
0&0&0&0&0&1&1&1&1\\
0&0&0&0&1&1&1&0&1\\
0&0&0&1&1&0&0&1&1\\
0&0&1&1&0&0&1&1&0\\
0&2&2&0&0&0&0&0&0\\
0&1&1&0&1&0&1&0&0\\
1&1&0&1&1&0&0&0&0\\
1&1&1&1&0&0&0&0&0\\
\ebbb
\parbox{7cm}{ 
\1\9\3\7\5\4\8\2\6\2\8\4\5\7\3\9\eeee
\8\2\6\2\8\4\5\7\3\9\1\9\3\7\5\4\eeee
\5\7\3\9\1\9\3\7\5\4\8\2\6\2\8\4\eeee
\3\7\5\4\8\2\6\2\8\4\5\7\3\9\1\9\eeee
\6\2\8\4\5\7\3\9\1\9\3\7\5\4\8\2\eeee
\3\9\1\9\3\7\5\4\8\2\6\2\8\4\5\7\eeee
\5\4\8\2\6\2\8\4\5\7\3\9\1\9\3\7\eeee
\8\4\5\7\3\9\1\9\3\7\5\4\8\2\6\2\eeee
\1\9\3\7\5\4\8\2\6\2\8\4\5\7\3\9\eeee
\8\2\6\2\8\4\5\7\3\9\1\9\3\7\5\4\eeee
\5\7\3\9\1\9\3\7\5\4\8\2\6\2\8\4\eeee
\3\7\5\4\8\2\6\2\8\4\5\7\3\9\1\9\eeee
} 

\baaa
9-141
\eaaa
\bbbb
0&0&0&0&0&0&0&2&2\\
0&0&0&0&0&1&1&1&1\\
0&0&0&0&1&1&2&0&0\\
0&0&0&1&1&0&1&0&1\\
0&0&1&1&0&0&0&1&1\\
0&1&1&0&0&1&0&1&0\\
0&1&2&1&0&0&0&0&0\\
1&1&0&0&1&1&0&0&0\\
1&1&0&1&1&0&0&0&0\\
\ebbb
\parbox{7cm}{ 
\1\9\4\5\8\2\7\3\6\6\3\7\2\8\5\4\eeee
\8\2\7\3\6\6\3\7\2\8\5\4\9\1\9\4\eeee
\6\6\3\7\2\8\5\4\9\1\9\4\5\8\2\7\eeee
\2\8\5\4\9\1\9\4\5\8\2\7\3\6\6\3\eeee
\9\1\9\4\5\8\2\7\3\6\6\3\7\2\8\5\eeee
\5\8\2\7\3\6\6\3\7\2\8\5\4\9\1\9\eeee
\3\6\6\3\7\2\8\5\4\9\1\9\4\5\8\2\eeee
\7\2\8\5\4\9\1\9\4\5\8\2\7\3\6\6\eeee
\4\9\1\9\4\5\8\2\7\3\6\6\3\7\2\8\eeee
\4\5\8\2\7\3\6\6\3\7\2\8\5\4\9\1\eeee
\7\3\6\6\3\7\2\8\5\4\9\1\9\4\5\8\eeee
\3\7\2\8\5\4\9\1\9\4\5\8\2\7\3\6\eeee
} 

\baaa
9-142
\eaaa
\bbbb
0&0&0&0&0&0&0&2&2\\
0&0&0&0&0&1&1&1&1\\
0&0&0&1&1&0&0&1&1\\
0&0&1&0&1&0&1&0&1\\
0&0&1&1&2&0&0&0&0\\
0&1&0&0&0&2&1&0&0\\
0&1&0&1&0&1&0&1&0\\
1&1&1&0&0&0&1&0&0\\
1&1&1&1&0&0&0&0&0\\
\ebbb
\parbox{7cm}{ 
\1\9\4\3\9\2\6\6\7\4\5\5\3\8\2\7\eeee
\8\2\7\8\1\8\7\2\8\3\5\5\4\7\6\6\eeee
\7\6\6\2\9\3\4\9\1\9\4\3\9\2\6\6\eeee
\2\6\6\7\4\5\5\3\8\2\7\8\1\8\7\2\eeee
\8\7\2\8\3\5\5\4\7\6\6\2\9\3\4\9\eeee
\3\4\9\1\9\4\3\9\2\6\6\7\4\5\5\3\eeee
\5\5\3\8\2\7\8\1\8\7\2\8\3\5\5\4\eeee
\5\5\4\7\6\6\2\9\3\4\9\1\9\4\3\9\eeee
\4\3\9\2\6\6\7\4\5\5\3\8\2\7\8\1\eeee
\7\8\1\8\7\2\8\3\5\5\4\7\6\6\2\9\eeee
\6\2\9\3\4\9\1\9\4\3\9\2\6\6\7\4\eeee
\6\7\4\5\5\3\8\2\7\8\1\8\7\2\8\3\eeee
} 

\baaa
9-143
\eaaa
\bbbb
0&0&0&0&0&0&0&2&2\\
0&0&0&0&0&1&1&1&1\\
0&0&0&1&1&1&1&0&0\\
0&0&1&0&2&0&1&0&0\\
0&0&1&2&1&0&0&0&0\\
0&1&1&0&0&0&1&0&1\\
0&1&1&1&0&1&0&0&0\\
1&1&0&0&0&0&0&1&1\\
1&1&0&0&0&1&0&1&0\\
\ebbb
\parbox{7cm}{ 
\1\9\6\3\5\4\7\2\8\8\2\7\4\5\3\6\eeee
\8\2\7\4\5\3\6\9\1\9\6\3\5\4\7\2\eeee
\9\6\3\5\4\7\2\8\8\2\7\4\5\3\6\9\eeee
\2\7\4\5\3\6\9\1\9\6\3\5\4\7\2\8\eeee
\6\3\5\4\7\2\8\8\2\7\4\5\3\6\9\1\eeee
\7\4\5\3\6\9\1\9\6\3\5\4\7\2\8\8\eeee
\3\5\4\7\2\8\8\2\7\4\5\3\6\9\1\9\eeee
\4\5\3\6\9\1\9\6\3\5\4\7\2\8\8\2\eeee
\5\4\7\2\8\8\2\7\4\5\3\6\9\1\9\6\eeee
\5\3\6\9\1\9\6\3\5\4\7\2\8\8\2\7\eeee
\4\7\2\8\8\2\7\4\5\3\6\9\1\9\6\3\eeee
\3\6\9\1\9\6\3\5\4\7\2\8\8\2\7\4\eeee
} 

\baaa
9-144
\eaaa
\bbbb
0&0&0&0&0&0&0&2&2\\
0&0&0&0&0&1&1&1&1\\
0&0&0&1&1&1&1&0&0\\
0&0&1&1&1&0&1&0&0\\
0&0&2&2&0&0&0&0&0\\
0&1&1&0&0&0&1&0&1\\
0&1&1&1&0&1&0&0&0\\
1&1&0&0&0&0&0&1&1\\
1&1&0&0&0&1&0&1&0\\
\ebbb
\parbox{7cm}{ 
\1\9\6\3\5\3\6\9\1\9\6\3\5\3\6\9\eeee
\8\2\7\4\4\7\2\8\8\2\7\4\4\7\2\8\eeee
\9\6\3\5\3\6\9\1\9\6\3\5\3\6\9\1\eeee
\2\7\4\4\7\2\8\8\2\7\4\4\7\2\8\8\eeee
\6\3\5\3\6\9\1\9\6\3\5\3\6\9\1\9\eeee
\7\4\4\7\2\8\8\2\7\4\4\7\2\8\8\2\eeee
\3\5\3\6\9\1\9\6\3\5\3\6\9\1\9\6\eeee
\4\4\7\2\8\8\2\7\4\4\7\2\8\8\2\7\eeee
\5\3\6\9\1\9\6\3\5\3\6\9\1\9\6\3\eeee
\4\7\2\8\8\2\7\4\4\7\2\8\8\2\7\4\eeee
\3\6\9\1\9\6\3\5\3\6\9\1\9\6\3\5\eeee
\7\2\8\8\2\7\4\4\7\2\8\8\2\7\4\4\eeee
} 

\baaa
9-145
\eaaa
\bbbb
0&0&0&0&0&0&0&2&2\\
0&0&0&0&0&1&1&1&1\\
0&0&1&0&0&0&2&0&1\\
0&0&0&1&0&2&0&1&0\\
0&0&0&0&2&1&1&0&0\\
0&1&0&1&1&1&0&0&0\\
0&1&1&0&1&0&1&0&0\\
1&2&0&1&0&0&0&0&0\\
1&2&1&0&0&0&0&0&0\\
\ebbb
\parbox{7cm}{ 
\1\9\3\3\9\1\9\3\3\9\1\9\3\3\9\1\eeee
\8\2\7\7\2\8\2\7\7\2\8\2\7\7\2\8\eeee
\4\6\5\5\6\4\6\5\5\6\4\6\5\5\6\4\eeee
\4\6\5\5\6\4\6\5\5\6\4\6\5\5\6\4\eeee
\8\2\7\7\2\8\2\7\7\2\8\2\7\7\2\8\eeee
\1\9\3\3\9\1\9\3\3\9\1\9\3\3\9\1\eeee
\8\2\7\7\2\8\2\7\7\2\8\2\7\7\2\8\eeee
\4\6\5\5\6\4\6\5\5\6\4\6\5\5\6\4\eeee
\4\6\5\5\6\4\6\5\5\6\4\6\5\5\6\4\eeee
\8\2\7\7\2\8\2\7\7\2\8\2\7\7\2\8\eeee
\1\9\3\3\9\1\9\3\3\9\1\9\3\3\9\1\eeee
\8\2\7\7\2\8\2\7\7\2\8\2\7\7\2\8\eeee
} 

\baaa
9-146
\eaaa
\bbbb
0&0&0&0&0&0&1&1&2\\
0&0&0&0&0&0&1&1&2\\
0&0&0&0&0&2&1&1&0\\
0&0&0&0&0&2&1&1&0\\
0&0&0&0&2&0&0&0&2\\
0&0&1&1&0&2&0&0&0\\
1&1&1&1&0&0&0&0&0\\
1&1&1&1&0&0&0&0&0\\
1&1&0&0&2&0&0&0&0\\
\ebbb
\parbox{7cm}{ 
\1\8\4\6\6\4\8\1\9\5\5\9\1\8\4\6\eeee
\7\3\6\6\3\7\2\9\5\5\9\2\7\3\6\6\eeee
\4\6\6\4\8\1\9\5\5\9\1\8\4\6\6\4\eeee
\6\6\3\7\2\9\5\5\9\2\7\3\6\6\3\7\eeee
\6\4\8\1\9\5\5\9\1\8\4\6\6\4\8\1\eeee
\3\7\2\9\5\5\9\2\7\3\6\6\3\7\2\9\eeee
\8\1\9\5\5\9\1\8\4\6\6\4\8\1\9\5\eeee
\2\9\5\5\9\2\7\3\6\6\3\7\2\9\5\5\eeee
\9\5\5\9\1\8\4\6\6\4\8\1\9\5\5\9\eeee
\5\5\9\2\7\3\6\6\3\7\2\9\5\5\9\2\eeee
\5\9\1\8\4\6\6\4\8\1\9\5\5\9\1\8\eeee
\9\2\7\3\6\6\3\7\2\9\5\5\9\2\7\3\eeee
} 

\baaa
9-147
\eaaa
\bbbb
0&0&0&0&0&0&1&1&2\\
0&0&0&0&0&0&1&1&2\\
0&0&0&0&0&2&1&1&0\\
0&0&0&1&1&0&0&0&2\\
0&0&0&1&1&0&0&0&2\\
0&0&2&0&0&2&0&0&0\\
1&1&2&0&0&0&0&0&0\\
1&1&2&0&0&0&0&0&0\\
1&1&0&1&1&0&0&0&0\\
\ebbb
\parbox{7cm}{ 
\1\8\3\6\6\3\8\1\9\5\5\9\1\8\3\6\eeee
\7\3\6\6\3\7\2\9\4\4\9\2\7\3\6\6\eeee
\3\6\6\3\8\1\9\5\5\9\1\8\3\6\6\3\eeee
\6\6\3\7\2\9\4\4\9\2\7\3\6\6\3\7\eeee
\6\3\8\1\9\5\5\9\1\8\3\6\6\3\8\1\eeee
\3\7\2\9\4\4\9\2\7\3\6\6\3\7\2\9\eeee
\8\1\9\5\5\9\1\8\3\6\6\3\8\1\9\5\eeee
\2\9\4\4\9\2\7\3\6\6\3\7\2\9\4\4\eeee
\9\5\5\9\1\8\3\6\6\3\8\1\9\5\5\9\eeee
\4\4\9\2\7\3\6\6\3\7\2\9\4\4\9\2\eeee
\5\9\1\8\3\6\6\3\8\1\9\5\5\9\1\8\eeee
\9\2\7\3\6\6\3\7\2\9\4\4\9\2\7\3\eeee
} 

\baaa
9-148
\eaaa
\bbbb
0&0&0&0&0&0&1&1&2\\
0&0&0&0&0&0&1&1&2\\
0&0&0&0&0&2&1&1&0\\
0&0&0&1&1&2&0&0&0\\
0&0&0&1&1&2&0&0&0\\
0&0&2&1&1&0&0&0&0\\
1&1&2&0&0&0&0&0&0\\
1&1&2&0&0&0&0&0&0\\
1&1&0&0&0&0&0&0&2\\
\ebbb
\parbox{7cm}{ 
\1\8\3\6\4\5\6\3\7\2\9\9\1\8\3\6\eeee
\7\3\6\5\4\6\3\8\1\9\9\2\7\3\6\5\eeee
\3\6\4\5\6\3\7\2\9\9\1\8\3\6\4\5\eeee
\6\5\4\6\3\8\1\9\9\2\7\3\6\5\4\6\eeee
\4\5\6\3\7\2\9\9\1\8\3\6\4\5\6\3\eeee
\4\6\3\8\1\9\9\2\7\3\6\5\4\6\3\8\eeee
\6\3\7\2\9\9\1\8\3\6\4\5\6\3\7\2\eeee
\3\8\1\9\9\2\7\3\6\5\4\6\3\8\1\9\eeee
\7\2\9\9\1\8\3\6\4\5\6\3\7\2\9\9\eeee
\1\9\9\2\7\3\6\5\4\6\3\8\1\9\9\2\eeee
\9\9\1\8\3\6\4\5\6\3\7\2\9\9\1\8\eeee
\9\2\7\3\6\5\4\6\3\8\1\9\9\2\7\3\eeee
} 

\baaa
9-149
\eaaa
\bbbb
0&0&0&0&0&0&1&1&2\\
0&0&0&0&0&0&1&1&2\\
0&0&0&0&1&1&0&0&2\\
0&0&0&0&1&1&0&0&2\\
0&0&1&1&0&0&1&1&0\\
0&0&1&1&0&0&1&1&0\\
1&1&0&0&1&1&0&0&0\\
1&1&0&0&1&1&0&0&0\\
1&1&1&1&0&0&0&0&0\\
\ebbb
\parbox{7cm}{ 
\1\8\6\3\9\2\8\5\3\9\2\7\5\4\9\1\eeee
\7\5\4\9\1\7\6\4\9\1\8\6\3\9\2\8\eeee
\6\3\9\2\8\5\3\9\2\7\5\4\9\1\7\6\eeee
\4\9\1\7\6\4\9\1\8\6\3\9\2\8\5\3\eeee
\9\2\8\5\3\9\2\7\5\4\9\1\7\6\4\9\eeee
\1\7\6\4\9\1\8\6\3\9\2\8\5\3\9\2\eeee
\8\5\3\9\2\7\5\4\9\1\7\6\4\9\1\8\eeee
\6\4\9\1\8\6\3\9\2\8\5\3\9\2\7\5\eeee
\3\9\2\7\5\4\9\1\7\6\4\9\1\8\6\3\eeee
\9\1\8\6\3\9\2\8\5\3\9\2\7\5\4\9\eeee
\2\7\5\4\9\1\7\6\4\9\1\8\6\3\9\2\eeee
\8\6\3\9\2\8\5\3\9\2\7\5\4\9\1\7\eeee
} 

\baaa
9-150
\eaaa
\bbbb
0&0&0&0&0&0&1&1&2\\
0&0&0&0&0&0&1&1&2\\
0&0&0&0&1&1&0&0&2\\
0&0&0&0&1&1&0&0&2\\
0&0&1&1&1&1&0&0&0\\
0&0&1&1&1&1&0&0&0\\
1&1&0&0&0&0&1&1&0\\
1&1&0&0&0&0&1&1&0\\
1&1&1&1&0&0&0&0&0\\
\ebbb
\parbox{7cm}{ 
\1\8\8\1\9\4\6\5\3\9\2\7\7\2\9\3\eeee
\7\7\2\9\3\5\6\4\9\1\8\8\1\9\4\6\eeee
\8\1\9\4\6\5\3\9\2\7\7\2\9\3\5\6\eeee
\2\9\3\5\6\4\9\1\8\8\1\9\4\6\5\3\eeee
\9\4\6\5\3\9\2\7\7\2\9\3\5\6\4\9\eeee
\3\5\6\4\9\1\8\8\1\9\4\6\5\3\9\2\eeee
\6\5\3\9\2\7\7\2\9\3\5\6\4\9\1\8\eeee
\6\4\9\1\8\8\1\9\4\6\5\3\9\2\7\7\eeee
\3\9\2\7\7\2\9\3\5\6\4\9\1\8\8\1\eeee
\9\1\8\8\1\9\4\6\5\3\9\2\7\7\2\9\eeee
\2\7\7\2\9\3\5\6\4\9\1\8\8\1\9\4\eeee
\8\8\1\9\4\6\5\3\9\2\7\7\2\9\3\5\eeee
} 

\baaa
9-151
\eaaa
\bbbb
0&0&0&0&0&0&1&1&2\\
0&0&0&0&0&0&1&1&2\\
0&0&0&0&1&1&1&1&0\\
0&0&0&0&1&1&1&1&0\\
0&0&1&1&1&1&0&0&0\\
0&0&1&1&1&1&0&0&0\\
1&1&1&1&0&0&0&0&0\\
1&1&1&1&0&0&0&0&0\\
1&1&0&0&0&0&0&0&2\\
\ebbb
\parbox{7cm}{ 
\1\8\4\5\5\4\8\1\9\9\2\7\3\6\6\3\eeee
\7\3\6\6\3\7\2\9\9\1\8\4\5\5\4\8\eeee
\4\5\5\4\8\1\9\9\2\7\3\6\6\3\7\2\eeee
\6\6\3\7\2\9\9\1\8\4\5\5\4\8\1\9\eeee
\5\4\8\1\9\9\2\7\3\6\6\3\7\2\9\9\eeee
\3\7\2\9\9\1\8\4\5\5\4\8\1\9\9\2\eeee
\8\1\9\9\2\7\3\6\6\3\7\2\9\9\1\8\eeee
\2\9\9\1\8\4\5\5\4\8\1\9\9\2\7\3\eeee
\9\9\2\7\3\6\6\3\7\2\9\9\1\8\4\5\eeee
\9\1\8\4\5\5\4\8\1\9\9\2\7\3\6\6\eeee
\2\7\3\6\6\3\7\2\9\9\1\8\4\5\5\4\eeee
\8\4\5\5\4\8\1\9\9\2\7\3\6\6\3\7\eeee
} 

\baaa
9-152
\eaaa
\bbbb
0&0&0&0&0&0&1&1&2\\
0&0&0&0&0&0&1&1&2\\
0&0&1&0&0&1&0&0&2\\
0&0&0&1&1&0&1&1&0\\
0&0&0&1&1&0&1&1&0\\
0&0&1&0&0&1&0&0&2\\
1&1&0&1&1&0&0&0&0\\
1&1&0&1&1&0&0&0&0\\
1&1&1&0&0&1&0&0&0\\
\ebbb
\parbox{7cm}{ 
\1\8\5\4\7\2\9\3\3\9\2\7\4\5\8\1\eeee
\7\4\5\8\1\9\6\6\9\1\8\5\4\7\2\9\eeee
\5\4\7\2\9\3\3\9\2\7\4\5\8\1\9\6\eeee
\5\8\1\9\6\6\9\1\8\5\4\7\2\9\3\3\eeee
\7\2\9\3\3\9\2\7\4\5\8\1\9\6\6\9\eeee
\1\9\6\6\9\1\8\5\4\7\2\9\3\3\9\2\eeee
\9\3\3\9\2\7\4\5\8\1\9\6\6\9\1\8\eeee
\6\6\9\1\8\5\4\7\2\9\3\3\9\2\7\4\eeee
\3\9\2\7\4\5\8\1\9\6\6\9\1\8\5\4\eeee
\9\1\8\5\4\7\2\9\3\3\9\2\7\4\5\8\eeee
\2\7\4\5\8\1\9\6\6\9\1\8\5\4\7\2\eeee
\8\5\4\7\2\9\3\3\9\2\7\4\5\8\1\9\eeee
} 

\baaa
9-153
\eaaa
\bbbb
0&0&0&0&0&0&1&1&2\\
0&0&0&0&0&1&0&1&2\\
0&0&0&0&1&0&1&0&2\\
0&0&0&0&1&1&0&0&2\\
0&0&1&1&0&0&1&1&0\\
0&1&0&1&0&1&0&1&0\\
1&0&1&0&1&0&1&0&0\\
1&1&0&0&1&1&0&0&0\\
1&1&1&1&0&0&0&0&0\\
\ebbb
\parbox{7cm}{ 
\1\8\6\2\9\4\6\6\4\9\2\6\8\1\9\2\eeee
\7\5\4\9\3\5\8\2\9\1\8\6\2\9\4\6\eeee
\7\3\9\1\7\7\1\9\3\7\5\4\9\3\5\8\eeee
\1\9\2\8\5\3\9\4\5\7\3\9\1\7\7\1\eeee
\9\4\6\6\4\9\2\6\8\1\9\2\8\5\3\9\eeee
\3\5\8\2\9\1\8\6\2\9\4\6\6\4\9\2\eeee
\7\7\1\9\3\7\5\4\9\3\5\8\2\9\1\8\eeee
\5\3\9\4\5\7\3\9\1\7\7\1\9\3\7\5\eeee
\4\9\2\6\8\1\9\2\8\5\3\9\4\5\7\3\eeee
\9\1\8\6\2\9\4\6\6\4\9\2\6\8\1\9\eeee
\3\7\5\4\9\3\5\8\2\9\1\8\6\2\9\4\eeee
\5\7\3\9\1\7\7\1\9\3\7\5\4\9\3\5\eeee
} 

\baaa
9-154
\eaaa
\bbbb
0&0&0&0&0&0&1&1&2\\
0&0&0&0&0&1&0&2&1\\
0&0&0&0&1&0&2&0&1\\
0&0&0&1&0&2&0&1&0\\
0&0&1&0&2&0&1&0&0\\
0&1&0&2&0&1&0&0&0\\
1&0&2&0&1&0&0&0&0\\
1&2&0&1&0&0&0&0&0\\
2&1&1&0&0&0&0&0&0\\
\ebbb
\parbox{7cm}{ 
\1\9\1\9\1\9\1\9\1\9\1\9\1\9\1\9\eeee
\7\3\7\3\7\3\7\3\7\3\7\3\7\3\7\3\eeee
\5\5\5\5\5\5\5\5\5\5\5\5\5\5\5\5\eeee
\3\7\3\7\3\7\3\7\3\7\3\7\3\7\3\7\eeee
\9\1\9\1\9\1\9\1\9\1\9\1\9\1\9\1\eeee
\2\8\2\8\2\8\2\8\2\8\2\8\2\8\2\8\eeee
\6\4\6\4\6\4\6\4\6\4\6\4\6\4\6\4\eeee
\6\4\6\4\6\4\6\4\6\4\6\4\6\4\6\4\eeee
\2\8\2\8\2\8\2\8\2\8\2\8\2\8\2\8\eeee
\9\1\9\1\9\1\9\1\9\1\9\1\9\1\9\1\eeee
\3\7\3\7\3\7\3\7\3\7\3\7\3\7\3\7\eeee
\5\5\5\5\5\5\5\5\5\5\5\5\5\5\5\5\eeee
} 

\baaa
9-155
\eaaa
\bbbb
0&0&0&0&0&0&1&1&2\\
0&0&0&0&0&1&1&2&0\\
0&0&0&0&1&1&0&1&1\\
0&0&0&1&1&0&1&0&1\\
0&0&1&1&0&1&1&0&0\\
0&1&1&0&1&1&0&0&0\\
1&1&0&1&1&0&0&0&0\\
1&2&1&0&0&0&0&0&0\\
2&0&1&1&0&0&0&0&0\\
\ebbb
\parbox{7cm}{ 
\1\8\2\7\4\5\6\3\9\1\8\2\7\4\5\6\eeee
\7\2\8\1\9\3\6\5\4\7\2\8\1\9\3\6\eeee
\5\6\3\9\1\8\2\7\4\5\6\3\9\1\8\2\eeee
\3\6\5\4\7\2\8\1\9\3\6\5\4\7\2\8\eeee
\8\2\7\4\5\6\3\9\1\8\2\7\4\5\6\3\eeee
\2\8\1\9\3\6\5\4\7\2\8\1\9\3\6\5\eeee
\6\3\9\1\8\2\7\4\5\6\3\9\1\8\2\7\eeee
\6\5\4\7\2\8\1\9\3\6\5\4\7\2\8\1\eeee
\2\7\4\5\6\3\9\1\8\2\7\4\5\6\3\9\eeee
\8\1\9\3\6\5\4\7\2\8\1\9\3\6\5\4\eeee
\3\9\1\8\2\7\4\5\6\3\9\1\8\2\7\4\eeee
\5\4\7\2\8\1\9\3\6\5\4\7\2\8\1\9\eeee
} 

\baaa
9-156
\eaaa
\bbbb
0&0&0&0&0&0&1&1&2\\
0&0&0&0&1&1&0&0&2\\
0&0&0&1&0&1&0&1&1\\
0&0&1&0&1&0&1&0&1\\
0&1&0&2&0&0&0&1&0\\
0&1&2&0&0&0&1&0&0\\
1&0&0&2&0&1&0&0&0\\
1&0&2&0&1&0&0&0&0\\
1&1&1&1&0&0&0&0&0\\
\ebbb
\parbox{7cm}{ 
\1\9\2\9\1\9\2\9\1\9\2\9\1\9\2\9\eeee
\7\4\5\4\7\4\5\4\7\4\5\4\7\4\5\4\eeee
\6\3\8\3\6\3\8\3\6\3\8\3\6\3\8\3\eeee
\2\9\1\9\2\9\1\9\2\9\1\9\2\9\1\9\eeee
\5\4\7\4\5\4\7\4\5\4\7\4\5\4\7\4\eeee
\8\3\6\3\8\3\6\3\8\3\6\3\8\3\6\3\eeee
\1\9\2\9\1\9\2\9\1\9\2\9\1\9\2\9\eeee
\7\4\5\4\7\4\5\4\7\4\5\4\7\4\5\4\eeee
\6\3\8\3\6\3\8\3\6\3\8\3\6\3\8\3\eeee
\2\9\1\9\2\9\1\9\2\9\1\9\2\9\1\9\eeee
\5\4\7\4\5\4\7\4\5\4\7\4\5\4\7\4\eeee
\8\3\6\3\8\3\6\3\8\3\6\3\8\3\6\3\eeee
} 

\baaa
9-157
\eaaa
\bbbb
0&0&0&0&0&0&1&1&2\\
0&0&0&0&1&1&0&0&2\\
0&0&0&1&0&1&0&1&1\\
0&0&1&0&1&0&1&0&1\\
0&1&0&2&0&1&0&0&0\\
0&1&2&0&1&0&0&0&0\\
1&0&0&2&0&0&0&1&0\\
1&0&2&0&0&0&1&0&0\\
1&1&1&1&0&0&0&0&0\\
\ebbb
\parbox{7cm}{ 
\1\9\2\9\1\9\2\9\1\9\2\9\1\9\2\9\eeee
\7\4\5\4\7\4\5\4\7\4\5\4\7\4\5\4\eeee
\8\3\6\3\8\3\6\3\8\3\6\3\8\3\6\3\eeee
\1\9\2\9\1\9\2\9\1\9\2\9\1\9\2\9\eeee
\7\4\5\4\7\4\5\4\7\4\5\4\7\4\5\4\eeee
\8\3\6\3\8\3\6\3\8\3\6\3\8\3\6\3\eeee
\1\9\2\9\1\9\2\9\1\9\2\9\1\9\2\9\eeee
\7\4\5\4\7\4\5\4\7\4\5\4\7\4\5\4\eeee
\8\3\6\3\8\3\6\3\8\3\6\3\8\3\6\3\eeee
\1\9\2\9\1\9\2\9\1\9\2\9\1\9\2\9\eeee
\7\4\5\4\7\4\5\4\7\4\5\4\7\4\5\4\eeee
\8\3\6\3\8\3\6\3\8\3\6\3\8\3\6\3\eeee
} 

\baaa
9-158
\eaaa
\bbbb
0&0&0&0&0&0&1&1&2\\
0&0&0&0&1&1&0&0&2\\
0&0&1&0&0&1&0&1&1\\
0&0&0&1&1&0&1&0&1\\
0&1&0&2&0&0&1&0&0\\
0&1&2&0&0&0&0&1&0\\
1&0&0&2&1&0&0&0&0\\
1&0&2&0&0&1&0&0&0\\
1&1&1&1&0&0&0&0&0\\
\ebbb
\parbox{7cm}{ 
\1\9\2\9\1\9\2\9\1\9\2\9\1\9\2\9\eeee
\7\4\5\4\7\4\5\4\7\4\5\4\7\4\5\4\eeee
\5\4\7\4\5\4\7\4\5\4\7\4\5\4\7\4\eeee
\2\9\1\9\2\9\1\9\2\9\1\9\2\9\1\9\eeee
\6\3\8\3\6\3\8\3\6\3\8\3\6\3\8\3\eeee
\8\3\6\3\8\3\6\3\8\3\6\3\8\3\6\3\eeee
\1\9\2\9\1\9\2\9\1\9\2\9\1\9\2\9\eeee
\7\4\5\4\7\4\5\4\7\4\5\4\7\4\5\4\eeee
\5\4\7\4\5\4\7\4\5\4\7\4\5\4\7\4\eeee
\2\9\1\9\2\9\1\9\2\9\1\9\2\9\1\9\eeee
\6\3\8\3\6\3\8\3\6\3\8\3\6\3\8\3\eeee
\8\3\6\3\8\3\6\3\8\3\6\3\8\3\6\3\eeee
} 

\baaa
9-159
\eaaa
\bbbb
0&0&0&0&0&0&1&1&2\\
0&0&0&0&1&1&0&0&2\\
0&0&1&0&0&1&0&1&1\\
0&0&0&1&1&0&1&0&1\\
0&1&0&2&0&0&1&0&0\\
0&1&2&0&0&1&0&0&0\\
1&0&0&2&1&0&0&0&0\\
1&0&2&0&0&0&0&1&0\\
1&1&1&1&0&0&0&0&0\\
\ebbb
\parbox{7cm}{ 
\1\9\2\9\1\9\2\9\1\9\2\9\1\9\2\9\eeee
\7\4\5\4\7\4\5\4\7\4\5\4\7\4\5\4\eeee
\5\4\7\4\5\4\7\4\5\4\7\4\5\4\7\4\eeee
\2\9\1\9\2\9\1\9\2\9\1\9\2\9\1\9\eeee
\6\3\8\3\6\3\8\3\6\3\8\3\6\3\8\3\eeee
\6\3\8\3\6\3\8\3\6\3\8\3\6\3\8\3\eeee
\2\9\1\9\2\9\1\9\2\9\1\9\2\9\1\9\eeee
\5\4\7\4\5\4\7\4\5\4\7\4\5\4\7\4\eeee
\7\4\5\4\7\4\5\4\7\4\5\4\7\4\5\4\eeee
\1\9\2\9\1\9\2\9\1\9\2\9\1\9\2\9\eeee
\8\3\6\3\8\3\6\3\8\3\6\3\8\3\6\3\eeee
\8\3\6\3\8\3\6\3\8\3\6\3\8\3\6\3\eeee
} 

\baaa
9-160
\eaaa
\bbbb
0&0&0&0&0&0&1&1&2\\
0&0&0&0&1&1&0&0&2\\
0&0&1&0&0&1&0&1&1\\
0&0&0&1&1&0&1&0&1\\
0&1&0&2&1&0&0&0&0\\
0&1&2&0&0&1&0&0&0\\
1&0&0&2&0&0&1&0&0\\
1&0&2&0&0&0&0&1&0\\
1&1&1&1&0&0&0&0&0\\
\ebbb
\parbox{7cm}{ 
\1\9\2\9\1\9\2\9\1\9\2\9\1\9\2\9\eeee
\7\4\5\4\7\4\5\4\7\4\5\4\7\4\5\4\eeee
\7\4\5\4\7\4\5\4\7\4\5\4\7\4\5\4\eeee
\1\9\2\9\1\9\2\9\1\9\2\9\1\9\2\9\eeee
\8\3\6\3\8\3\6\3\8\3\6\3\8\3\6\3\eeee
\8\3\6\3\8\3\6\3\8\3\6\3\8\3\6\3\eeee
\1\9\2\9\1\9\2\9\1\9\2\9\1\9\2\9\eeee
\7\4\5\4\7\4\5\4\7\4\5\4\7\4\5\4\eeee
\7\4\5\4\7\4\5\4\7\4\5\4\7\4\5\4\eeee
\1\9\2\9\1\9\2\9\1\9\2\9\1\9\2\9\eeee
\8\3\6\3\8\3\6\3\8\3\6\3\8\3\6\3\eeee
\8\3\6\3\8\3\6\3\8\3\6\3\8\3\6\3\eeee
} 

\baaa
9-161
\eaaa
\bbbb
0&0&0&0&0&0&1&1&2\\
0&0&0&0&1&1&0&1&1\\
0&0&0&1&0&1&0&2&0\\
0&0&1&0&2&1&0&0&0\\
0&1&0&1&0&0&1&1&0\\
0&2&1&1&0&0&0&0&0\\
1&0&0&0&1&0&1&0&1\\
1&1&1&0&1&0&0&0&0\\
2&1&0&0&0&0&1&0&0\\
\ebbb
\parbox{7cm}{ 
\1\8\3\8\1\9\2\6\2\9\1\8\3\8\1\9\eeee
\7\5\4\5\7\7\5\4\5\7\7\5\4\5\7\7\eeee
\9\2\6\2\9\1\8\3\8\1\9\2\6\2\9\1\eeee
\1\8\3\8\1\9\2\6\2\9\1\8\3\8\1\9\eeee
\7\5\4\5\7\7\5\4\5\7\7\5\4\5\7\7\eeee
\9\2\6\2\9\1\8\3\8\1\9\2\6\2\9\1\eeee
\1\8\3\8\1\9\2\6\2\9\1\8\3\8\1\9\eeee
\7\5\4\5\7\7\5\4\5\7\7\5\4\5\7\7\eeee
\9\2\6\2\9\1\8\3\8\1\9\2\6\2\9\1\eeee
\1\8\3\8\1\9\2\6\2\9\1\8\3\8\1\9\eeee
\7\5\4\5\7\7\5\4\5\7\7\5\4\5\7\7\eeee
\9\2\6\2\9\1\8\3\8\1\9\2\6\2\9\1\eeee
} 

\baaa
9-162
\eaaa
\bbbb
0&0&0&0&0&0&1&1&2\\
0&0&0&0&1&1&0&1&1\\
0&0&0&1&0&1&1&0&1\\
0&0&1&0&2&0&0&1&0\\
0&1&0&2&0&0&1&0&0\\
0&1&1&0&0&0&1&1&0\\
1&0&1&0&1&1&0&0&0\\
1&1&0&1&0&1&0&0&0\\
2&1&1&0&0&0&0&0&0\\
\ebbb
\parbox{7cm}{ 
\1\8\4\5\2\9\1\8\4\5\2\9\1\8\4\5\eeee
\7\6\3\7\6\3\7\6\3\7\6\3\7\6\3\7\eeee
\5\2\9\1\8\4\5\2\9\1\8\4\5\2\9\1\eeee
\4\8\1\9\2\5\4\8\1\9\2\5\4\8\1\9\eeee
\3\6\7\3\6\7\3\6\7\3\6\7\3\6\7\3\eeee
\9\2\5\4\8\1\9\2\5\4\8\1\9\2\5\4\eeee
\1\8\4\5\2\9\1\8\4\5\2\9\1\8\4\5\eeee
\7\6\3\7\6\3\7\6\3\7\6\3\7\6\3\7\eeee
\5\2\9\1\8\4\5\2\9\1\8\4\5\2\9\1\eeee
\4\8\1\9\2\5\4\8\1\9\2\5\4\8\1\9\eeee
\3\6\7\3\6\7\3\6\7\3\6\7\3\6\7\3\eeee
\9\2\5\4\8\1\9\2\5\4\8\1\9\2\5\4\eeee
} 

\baaa
9-163
\eaaa
\bbbb
0&0&0&0&0&0&1&1&2\\
0&0&0&0&1&1&0&1&1\\
0&0&0&1&0&1&1&0&1\\
0&0&1&0&2&0&0&1&0\\
0&1&0&2&0&1&0&0&0\\
0&1&1&0&1&1&0&0&0\\
1&0&1&0&0&0&1&1&0\\
1&1&0&1&0&0&1&0&0\\
2&1&1&0&0&0&0&0&0\\
\ebbb
\parbox{7cm}{ 
\1\8\4\5\2\9\1\8\4\5\2\9\1\8\4\5\eeee
\7\7\3\6\6\3\7\7\3\6\6\3\7\7\3\6\eeee
\8\1\9\2\5\4\8\1\9\2\5\4\8\1\9\2\eeee
\2\9\1\8\4\5\2\9\1\8\4\5\2\9\1\8\eeee
\6\3\7\7\3\6\6\3\7\7\3\6\6\3\7\7\eeee
\5\4\8\1\9\2\5\4\8\1\9\2\5\4\8\1\eeee
\4\5\2\9\1\8\4\5\2\9\1\8\4\5\2\9\eeee
\3\6\6\3\7\7\3\6\6\3\7\7\3\6\6\3\eeee
\9\2\5\4\8\1\9\2\5\4\8\1\9\2\5\4\eeee
\1\8\4\5\2\9\1\8\4\5\2\9\1\8\4\5\eeee
\7\7\3\6\6\3\7\7\3\6\6\3\7\7\3\6\eeee
\8\1\9\2\5\4\8\1\9\2\5\4\8\1\9\2\eeee
} 

\baaa
9-164
\eaaa
\bbbb
0&0&0&0&0&0&1&1&2\\
0&0&0&0&1&1&0&1&1\\
0&0&0&1&0&1&1&0&1\\
0&0&1&2&0&0&1&0&0\\
0&1&0&0&2&0&0&1&0\\
0&1&1&0&0&0&1&1&0\\
1&0&1&1&0&1&0&0&0\\
1&1&0&0&1&1&0&0&0\\
2&1&1&0&0&0&0&0&0\\
\ebbb
\parbox{7cm}{ 
\1\8\5\5\2\9\1\8\5\5\2\9\1\8\5\5\eeee
\7\6\2\8\6\3\7\6\2\8\6\3\7\6\2\8\eeee
\4\3\9\1\7\4\4\3\9\1\7\4\4\3\9\1\eeee
\4\7\1\9\3\4\4\7\1\9\3\4\4\7\1\9\eeee
\3\6\8\2\6\7\3\6\8\2\6\7\3\6\8\2\eeee
\9\2\5\5\8\1\9\2\5\5\8\1\9\2\5\5\eeee
\1\8\5\5\2\9\1\8\5\5\2\9\1\8\5\5\eeee
\7\6\2\8\6\3\7\6\2\8\6\3\7\6\2\8\eeee
\4\3\9\1\7\4\4\3\9\1\7\4\4\3\9\1\eeee
\4\7\1\9\3\4\4\7\1\9\3\4\4\7\1\9\eeee
\3\6\8\2\6\7\3\6\8\2\6\7\3\6\8\2\eeee
\9\2\5\5\8\1\9\2\5\5\8\1\9\2\5\5\eeee
} 

\baaa
9-165
\eaaa
\bbbb
0&0&0&0&0&0&1&1&2\\
0&0&0&0&1&1&0&1&1\\
0&0&0&1&0&2&0&1&0\\
0&0&1&1&1&1&0&0&0\\
0&1&0&1&0&0&1&1&0\\
0&1&2&1&0&0&0&0&0\\
1&0&0&0&1&0&1&0&1\\
1&1&1&0&1&0&0&0&0\\
2&1&0&0&0&0&1&0&0\\
\ebbb
\parbox{7cm}{ 
\1\8\3\6\2\9\1\8\3\6\2\9\1\8\3\6\eeee
\7\5\4\4\5\7\7\5\4\4\5\7\7\5\4\4\eeee
\9\2\6\3\8\1\9\2\6\3\8\1\9\2\6\3\eeee
\1\8\3\6\2\9\1\8\3\6\2\9\1\8\3\6\eeee
\7\5\4\4\5\7\7\5\4\4\5\7\7\5\4\4\eeee
\9\2\6\3\8\1\9\2\6\3\8\1\9\2\6\3\eeee
\1\8\3\6\2\9\1\8\3\6\2\9\1\8\3\6\eeee
\7\5\4\4\5\7\7\5\4\4\5\7\7\5\4\4\eeee
\9\2\6\3\8\1\9\2\6\3\8\1\9\2\6\3\eeee
\1\8\3\6\2\9\1\8\3\6\2\9\1\8\3\6\eeee
\7\5\4\4\5\7\7\5\4\4\5\7\7\5\4\4\eeee
\9\2\6\3\8\1\9\2\6\3\8\1\9\2\6\3\eeee
} 

\baaa
9-166
\eaaa
\bbbb
0&0&0&0&0&0&1&1&2\\
0&0&0&0&1&1&0&1&1\\
0&0&0&1&1&2&0&0&0\\
0&0&1&1&0&1&0&1&0\\
0&1&1&0&0&0&1&1&0\\
0&1&2&1&0&0&0&0&0\\
1&0&0&0&1&0&1&0&1\\
1&1&0&1&1&0&0&0&0\\
2&1&0&0&0&0&1&0&0\\
\ebbb
\parbox{7cm}{ 
\1\8\4\4\8\1\9\2\6\3\5\7\7\5\3\6\eeee
\7\5\3\6\2\9\1\8\4\4\8\1\9\2\6\3\eeee
\9\2\6\3\5\7\7\5\3\6\2\9\1\8\4\4\eeee
\1\8\4\4\8\1\9\2\6\3\5\7\7\5\3\6\eeee
\7\5\3\6\2\9\1\8\4\4\8\1\9\2\6\3\eeee
\9\2\6\3\5\7\7\5\3\6\2\9\1\8\4\4\eeee
\1\8\4\4\8\1\9\2\6\3\5\7\7\5\3\6\eeee
\7\5\3\6\2\9\1\8\4\4\8\1\9\2\6\3\eeee
\9\2\6\3\5\7\7\5\3\6\2\9\1\8\4\4\eeee
\1\8\4\4\8\1\9\2\6\3\5\7\7\5\3\6\eeee
\7\5\3\6\2\9\1\8\4\4\8\1\9\2\6\3\eeee
\9\2\6\3\5\7\7\5\3\6\2\9\1\8\4\4\eeee
} 

\baaa
9-167
\eaaa
\bbbb
0&0&0&0&0&0&1&1&2\\
0&0&0&0&1&1&0&1&1\\
0&0&1&1&0&1&0&0&1\\
0&0&1&1&1&1&0&0&0\\
0&1&0&1&0&0&1&0&1\\
0&1&1&1&0&1&0&0&0\\
1&0&0&0&2&0&0&1&0\\
1&2&0&0&0&0&1&0&0\\
1&1&1&0&1&0&0&0&0\\
\ebbb
\parbox{7cm}{ 
\1\9\3\3\9\1\9\3\3\9\1\9\3\3\9\1\eeee
\7\5\4\4\5\7\5\4\4\5\7\5\4\4\5\7\eeee
\8\2\6\6\2\8\2\6\6\2\8\2\6\6\2\8\eeee
\1\9\3\3\9\1\9\3\3\9\1\9\3\3\9\1\eeee
\7\5\4\4\5\7\5\4\4\5\7\5\4\4\5\7\eeee
\8\2\6\6\2\8\2\6\6\2\8\2\6\6\2\8\eeee
\1\9\3\3\9\1\9\3\3\9\1\9\3\3\9\1\eeee
\7\5\4\4\5\7\5\4\4\5\7\5\4\4\5\7\eeee
\8\2\6\6\2\8\2\6\6\2\8\2\6\6\2\8\eeee
\1\9\3\3\9\1\9\3\3\9\1\9\3\3\9\1\eeee
\7\5\4\4\5\7\5\4\4\5\7\5\4\4\5\7\eeee
\8\2\6\6\2\8\2\6\6\2\8\2\6\6\2\8\eeee
} 

\baaa
9-168
\eaaa
\bbbb
0&0&0&0&0&0&1&1&2\\
0&0&0&0&1&1&0&1&1\\
0&0&1&1&0&1&0&1&0\\
0&0&1&1&1&1&0&0&0\\
0&1&0&1&0&0&1&1&0\\
0&1&1&1&0&1&0&0&0\\
1&0&0&0&1&0&1&0&1\\
1&1&1&0&1&0&0&0&0\\
2&1&0&0&0&0&1&0&0\\
\ebbb
\parbox{7cm}{ 
\1\8\3\3\8\1\9\2\6\6\2\9\1\8\3\3\eeee
\7\5\4\4\5\7\7\5\4\4\5\7\7\5\4\4\eeee
\9\2\6\6\2\9\1\8\3\3\8\1\9\2\6\6\eeee
\1\8\3\3\8\1\9\2\6\6\2\9\1\8\3\3\eeee
\7\5\4\4\5\7\7\5\4\4\5\7\7\5\4\4\eeee
\9\2\6\6\2\9\1\8\3\3\8\1\9\2\6\6\eeee
\1\8\3\3\8\1\9\2\6\6\2\9\1\8\3\3\eeee
\7\5\4\4\5\7\7\5\4\4\5\7\7\5\4\4\eeee
\9\2\6\6\2\9\1\8\3\3\8\1\9\2\6\6\eeee
\1\8\3\3\8\1\9\2\6\6\2\9\1\8\3\3\eeee
\7\5\4\4\5\7\7\5\4\4\5\7\7\5\4\4\eeee
\9\2\6\6\2\9\1\8\3\3\8\1\9\2\6\6\eeee
} 

\baaa
9-169
\eaaa
\bbbb
0&0&0&0&0&0&1&1&2\\
0&0&0&0&1&1&1&1&0\\
0&0&0&2&1&1&0&0&0\\
0&0&2&1&0&1&0&0&0\\
0&1&1&0&0&1&0&1&0\\
0&1&1&1&1&0&0&0&0\\
1&1&0&0&0&0&0&1&1\\
1&1&0&0&1&0&1&0&0\\
2&0&0&0&0&0&1&0&1\\
\ebbb
\parbox{7cm}{ 
\1\8\5\3\4\6\2\7\9\1\8\5\3\4\6\2\eeee
\7\2\6\4\3\5\8\1\9\7\2\6\4\3\5\8\eeee
\8\5\3\4\6\2\7\9\1\8\5\3\4\6\2\7\eeee
\2\6\4\3\5\8\1\9\7\2\6\4\3\5\8\1\eeee
\5\3\4\6\2\7\9\1\8\5\3\4\6\2\7\9\eeee
\6\4\3\5\8\1\9\7\2\6\4\3\5\8\1\9\eeee
\3\4\6\2\7\9\1\8\5\3\4\6\2\7\9\1\eeee
\4\3\5\8\1\9\7\2\6\4\3\5\8\1\9\7\eeee
\4\6\2\7\9\1\8\5\3\4\6\2\7\9\1\8\eeee
\3\5\8\1\9\7\2\6\4\3\5\8\1\9\7\2\eeee
\6\2\7\9\1\8\5\3\4\6\2\7\9\1\8\5\eeee
\5\8\1\9\7\2\6\4\3\5\8\1\9\7\2\6\eeee
} 

\baaa
9-170
\eaaa
\bbbb
0&0&0&0&0&0&1&1&2\\
0&1&0&0&0&1&0&1&1\\
0&0&1&0&1&0&1&0&1\\
0&0&0&1&1&1&0&0&1\\
0&0&1&1&1&1&0&0&0\\
0&1&0&1&1&1&0&0&0\\
1&0&2&0&0&0&1&0&0\\
1&2&0&0&0&0&0&1&0\\
1&1&1&1&0&0&0&0&0\\
\ebbb
\parbox{7cm}{ 
\1\9\4\4\9\1\9\4\4\9\1\9\4\4\9\1\eeee
\7\3\5\6\2\8\2\6\5\3\7\3\5\6\2\8\eeee
\7\3\5\6\2\8\2\6\5\3\7\3\5\6\2\8\eeee
\1\9\4\4\9\1\9\4\4\9\1\9\4\4\9\1\eeee
\8\2\6\5\3\7\3\5\6\2\8\2\6\5\3\7\eeee
\8\2\6\5\3\7\3\5\6\2\8\2\6\5\3\7\eeee
\1\9\4\4\9\1\9\4\4\9\1\9\4\4\9\1\eeee
\7\3\5\6\2\8\2\6\5\3\7\3\5\6\2\8\eeee
\7\3\5\6\2\8\2\6\5\3\7\3\5\6\2\8\eeee
\1\9\4\4\9\1\9\4\4\9\1\9\4\4\9\1\eeee
\8\2\6\5\3\7\3\5\6\2\8\2\6\5\3\7\eeee
\8\2\6\5\3\7\3\5\6\2\8\2\6\5\3\7\eeee
} 

\baaa
9-171
\eaaa
\bbbb
0&0&0&0&0&0&1&1&2\\
0&1&0&0&0&1&0&1&1\\
0&0&1&0&1&0&1&0&1\\
0&0&0&1&1&1&0&0&1\\
0&0&1&1&2&0&0&0&0\\
0&1&0&1&0&2&0&0&0\\
1&0&2&0&0&0&1&0&0\\
1&2&0&0&0&0&0&1&0\\
1&1&1&1&0&0&0&0&0\\
\ebbb
\parbox{7cm}{ 
\1\9\4\4\9\1\9\4\4\9\1\9\4\4\9\1\eeee
\7\3\5\5\3\7\3\5\5\3\7\3\5\5\3\7\eeee
\7\3\5\5\3\7\3\5\5\3\7\3\5\5\3\7\eeee
\1\9\4\4\9\1\9\4\4\9\1\9\4\4\9\1\eeee
\8\2\6\6\2\8\2\6\6\2\8\2\6\6\2\8\eeee
\8\2\6\6\2\8\2\6\6\2\8\2\6\6\2\8\eeee
\1\9\4\4\9\1\9\4\4\9\1\9\4\4\9\1\eeee
\7\3\5\5\3\7\3\5\5\3\7\3\5\5\3\7\eeee
\7\3\5\5\3\7\3\5\5\3\7\3\5\5\3\7\eeee
\1\9\4\4\9\1\9\4\4\9\1\9\4\4\9\1\eeee
\8\2\6\6\2\8\2\6\6\2\8\2\6\6\2\8\eeee
\8\2\6\6\2\8\2\6\6\2\8\2\6\6\2\8\eeee
} 

\baaa
9-172
\eaaa
\bbbb
0&0&0&0&0&1&1&1&1\\
0&0&0&0&1&0&1&1&1\\
0&0&0&1&0&0&1&1&1\\
0&0&1&0&1&1&0&0&1\\
0&1&0&1&0&1&0&1&0\\
1&0&0&1&1&0&1&0&0\\
1&1&1&0&0&1&0&0&0\\
1&1&1&0&1&0&0&0&0\\
1&1&1&1&0&0&0&0&0\\
\ebbb
\parbox{7cm}{ 
\1\8\3\7\2\9\1\8\3\7\2\9\1\8\3\7\eeee
\6\5\4\6\5\4\6\5\4\6\5\4\6\5\4\6\eeee
\7\2\9\1\8\3\7\2\9\1\8\3\7\2\9\1\eeee
\1\8\3\7\2\9\1\8\3\7\2\9\1\8\3\7\eeee
\6\5\4\6\5\4\6\5\4\6\5\4\6\5\4\6\eeee
\7\2\9\1\8\3\7\2\9\1\8\3\7\2\9\1\eeee
\1\8\3\7\2\9\1\8\3\7\2\9\1\8\3\7\eeee
\6\5\4\6\5\4\6\5\4\6\5\4\6\5\4\6\eeee
\7\2\9\1\8\3\7\2\9\1\8\3\7\2\9\1\eeee
\1\8\3\7\2\9\1\8\3\7\2\9\1\8\3\7\eeee
\6\5\4\6\5\4\6\5\4\6\5\4\6\5\4\6\eeee
\7\2\9\1\8\3\7\2\9\1\8\3\7\2\9\1\eeee
} 

\baaa
9-173
\eaaa
\bbbb
0&0&0&0&0&1&1&1&1\\
0&0&0&0&1&0&1&1&1\\
0&0&0&1&0&0&1&1&1\\
0&0&1&0&1&1&0&0&1\\
0&1&0&1&1&0&0&1&0\\
1&0&0&1&0&1&1&0&0\\
1&1&1&0&0&1&0&0&0\\
1&1&1&0&1&0&0&0&0\\
1&1&1&1&0&0&0&0&0\\
\ebbb
\parbox{7cm}{ 
\1\7\3\8\2\9\1\7\3\8\2\9\1\7\3\8\eeee
\6\6\4\5\5\4\6\6\4\5\5\4\6\6\4\5\eeee
\7\1\9\2\8\3\7\1\9\2\8\3\7\1\9\2\eeee
\2\8\3\7\1\9\2\8\3\7\1\9\2\8\3\7\eeee
\5\5\4\6\6\4\5\5\4\6\6\4\5\5\4\6\eeee
\8\2\9\1\7\3\8\2\9\1\7\3\8\2\9\1\eeee
\1\7\3\8\2\9\1\7\3\8\2\9\1\7\3\8\eeee
\6\6\4\5\5\4\6\6\4\5\5\4\6\6\4\5\eeee
\7\1\9\2\8\3\7\1\9\2\8\3\7\1\9\2\eeee
\2\8\3\7\1\9\2\8\3\7\1\9\2\8\3\7\eeee
\5\5\4\6\6\4\5\5\4\6\6\4\5\5\4\6\eeee
\8\2\9\1\7\3\8\2\9\1\7\3\8\2\9\1\eeee
} 

\baaa
9-174
\eaaa
\bbbb
0&0&0&0&0&1&1&1&1\\
0&0&0&0&1&0&1&1&1\\
0&0&0&1&0&1&0&1&1\\
0&0&1&0&1&0&1&0&1\\
0&1&0&1&0&1&0&1&0\\
1&0&1&0&1&0&1&0&0\\
1&1&0&1&0&1&0&0&0\\
1&1&1&0&1&0&0&0&0\\
1&1&1&1&0&0&0&0&0\\
\ebbb
\parbox{7cm}{ 
\1\8\3\6\5\4\7\2\9\1\8\3\6\5\4\7\eeee
\6\5\4\7\2\9\1\8\3\6\5\4\7\2\9\1\eeee
\7\2\9\1\8\3\6\5\4\7\2\9\1\8\3\6\eeee
\1\8\3\6\5\4\7\2\9\1\8\3\6\5\4\7\eeee
\6\5\4\7\2\9\1\8\3\6\5\4\7\2\9\1\eeee
\7\2\9\1\8\3\6\5\4\7\2\9\1\8\3\6\eeee
\1\8\3\6\5\4\7\2\9\1\8\3\6\5\4\7\eeee
\6\5\4\7\2\9\1\8\3\6\5\4\7\2\9\1\eeee
\7\2\9\1\8\3\6\5\4\7\2\9\1\8\3\6\eeee
\1\8\3\6\5\4\7\2\9\1\8\3\6\5\4\7\eeee
\6\5\4\7\2\9\1\8\3\6\5\4\7\2\9\1\eeee
\7\2\9\1\8\3\6\5\4\7\2\9\1\8\3\6\eeee
} 

\baaa
9-175
\eaaa
\bbbb
0&0&0&0&0&1&1&1&1\\
0&0&0&0&1&0&1&1&1\\
0&0&0&1&0&1&0&1&1\\
0&0&1&1&0&1&0&0&1\\
0&1&0&0&1&0&1&1&0\\
1&0&1&1&0&1&0&0&0\\
1&1&0&0&1&0&1&0&0\\
1&1&1&0&1&0&0&0&0\\
1&1&1&1&0&0&0&0&0\\
\ebbb
\parbox{7cm}{ 
\1\8\2\7\5\5\7\2\8\1\9\3\6\4\4\6\eeee
\6\3\9\1\8\2\7\5\5\7\2\8\1\9\3\6\eeee
\6\4\4\6\3\9\1\8\2\7\5\5\7\2\8\1\eeee
\1\9\3\6\4\4\6\3\9\1\8\2\7\5\5\7\eeee
\7\2\8\1\9\3\6\4\4\6\3\9\1\8\2\7\eeee
\7\5\5\7\2\8\1\9\3\6\4\4\6\3\9\1\eeee
\1\8\2\7\5\5\7\2\8\1\9\3\6\4\4\6\eeee
\6\3\9\1\8\2\7\5\5\7\2\8\1\9\3\6\eeee
\6\4\4\6\3\9\1\8\2\7\5\5\7\2\8\1\eeee
\1\9\3\6\4\4\6\3\9\1\8\2\7\5\5\7\eeee
\7\2\8\1\9\3\6\4\4\6\3\9\1\8\2\7\eeee
\7\5\5\7\2\8\1\9\3\6\4\4\6\3\9\1\eeee
} 

\baaa
9-176
\eaaa
\bbbb
0&0&0&0&0&1&1&1&1\\
0&0&0&1&1&0&0&1&1\\
0&0&0&1&1&1&1&0&0\\
0&1&1&0&0&0&1&0&1\\
0&1&1&0&0&1&0&1&0\\
1&0&1&0&1&0&0&0&1\\
1&0&1&1&0&0&0&1&0\\
1&1&0&0&1&0&1&0&0\\
1&1&0&1&0&1&0&0&0\\
\ebbb
\parbox{7cm}{ 
\1\7\8\1\7\8\1\7\8\1\7\8\1\7\8\1\eeee
\6\3\5\6\3\5\6\3\5\6\3\5\6\3\5\6\eeee
\9\4\2\9\4\2\9\4\2\9\4\2\9\4\2\9\eeee
\1\7\8\1\7\8\1\7\8\1\7\8\1\7\8\1\eeee
\6\3\5\6\3\5\6\3\5\6\3\5\6\3\5\6\eeee
\9\4\2\9\4\2\9\4\2\9\4\2\9\4\2\9\eeee
\1\7\8\1\7\8\1\7\8\1\7\8\1\7\8\1\eeee
\6\3\5\6\3\5\6\3\5\6\3\5\6\3\5\6\eeee
\9\4\2\9\4\2\9\4\2\9\4\2\9\4\2\9\eeee
\1\7\8\1\7\8\1\7\8\1\7\8\1\7\8\1\eeee
\6\3\5\6\3\5\6\3\5\6\3\5\6\3\5\6\eeee
\9\4\2\9\4\2\9\4\2\9\4\2\9\4\2\9\eeee
} 

\baaa
9-177
\eaaa
\bbbb
0&0&0&0&0&1&1&1&1\\
0&0&0&1&1&0&0&1&1\\
0&0&0&1&1&1&1&0&0\\
0&1&1&0&0&0&1&0&1\\
0&1&1&0&0&1&0&1&0\\
1&0&1&0&1&0&0&1&0\\
1&0&1&1&0&0&0&0&1\\
1&1&0&0&1&1&0&0&0\\
1&1&0&1&0&0&1&0&0\\
\ebbb
\parbox{7cm}{ 
\1\7\9\1\6\8\1\7\9\1\6\8\1\7\9\1\eeee
\6\3\4\7\3\5\6\3\4\7\3\5\6\3\4\7\eeee
\8\5\2\9\4\2\8\5\2\9\4\2\8\5\2\9\eeee
\1\6\8\1\7\9\1\6\8\1\7\9\1\6\8\1\eeee
\7\3\5\6\3\4\7\3\5\6\3\4\7\3\5\6\eeee
\9\4\2\8\5\2\9\4\2\8\5\2\9\4\2\8\eeee
\1\7\9\1\6\8\1\7\9\1\6\8\1\7\9\1\eeee
\6\3\4\7\3\5\6\3\4\7\3\5\6\3\4\7\eeee
\8\5\2\9\4\2\8\5\2\9\4\2\8\5\2\9\eeee
\1\6\8\1\7\9\1\6\8\1\7\9\1\6\8\1\eeee
\7\3\5\6\3\4\7\3\5\6\3\4\7\3\5\6\eeee
\9\4\2\8\5\2\9\4\2\8\5\2\9\4\2\8\eeee
} 

\baaa
9-178
\eaaa
\bbbb
0&0&0&0&0&1&1&1&1\\
0&0&0&1&1&0&0&1&1\\
0&0&0&1&1&1&1&0&0\\
0&1&1&0&0&0&1&0&1\\
0&1&1&0&0&1&0&1&0\\
1&0&1&0&1&0&0&1&0\\
1&0&1&1&0&0&1&0&0\\
1&1&0&0&1&1&0&0&0\\
1&1&0&1&0&0&0&0&1\\
\ebbb
\parbox{7cm}{ 
\1\7\7\1\8\6\1\9\9\1\6\8\1\7\7\1\eeee
\6\3\4\9\2\5\8\2\4\7\3\5\6\3\4\9\eeee
\8\5\2\9\4\3\6\5\3\7\4\2\8\5\2\9\eeee
\1\6\8\1\7\7\1\8\6\1\9\9\1\6\8\1\eeee
\7\3\5\6\3\4\9\2\5\8\2\4\7\3\5\6\eeee
\7\4\2\8\5\2\9\4\3\6\5\3\7\4\2\8\eeee
\1\9\9\1\6\8\1\7\7\1\8\6\1\9\9\1\eeee
\8\2\4\7\3\5\6\3\4\9\2\5\8\2\4\7\eeee
\6\5\3\7\4\2\8\5\2\9\4\3\6\5\3\7\eeee
\1\8\6\1\9\9\1\6\8\1\7\7\1\8\6\1\eeee
\9\2\5\8\2\4\7\3\5\6\3\4\9\2\5\8\eeee
\9\4\3\6\5\3\7\4\2\8\5\2\9\4\3\6\eeee
} 

\baaa
9-179
\eaaa
\bbbb
0&0&0&0&0&1&1&1&1\\
0&0&0&1&1&0&0&1&1\\
0&0&0&1&1&1&1&0&0\\
0&1&1&0&0&0&1&0&1\\
0&1&1&0&0&1&0&1&0\\
1&0&1&0&1&1&0&0&0\\
1&0&1&1&0&0&1&0&0\\
1&1&0&0&1&0&0&1&0\\
1&1&0&1&0&0&0&0&1\\
\ebbb
\parbox{7cm}{ 
\1\7\7\1\8\8\1\7\7\1\8\8\1\7\7\1\eeee
\6\3\4\9\2\5\6\3\4\9\2\5\6\3\4\9\eeee
\6\5\2\9\4\3\6\5\2\9\4\3\6\5\2\9\eeee
\1\8\8\1\7\7\1\8\8\1\7\7\1\8\8\1\eeee
\9\2\5\6\3\4\9\2\5\6\3\4\9\2\5\6\eeee
\9\4\3\6\5\2\9\4\3\6\5\2\9\4\3\6\eeee
\1\7\7\1\8\8\1\7\7\1\8\8\1\7\7\1\eeee
\6\3\4\9\2\5\6\3\4\9\2\5\6\3\4\9\eeee
\6\5\2\9\4\3\6\5\2\9\4\3\6\5\2\9\eeee
\1\8\8\1\7\7\1\8\8\1\7\7\1\8\8\1\eeee
\9\2\5\6\3\4\9\2\5\6\3\4\9\2\5\6\eeee
\9\4\3\6\5\2\9\4\3\6\5\2\9\4\3\6\eeee
} 

\baaa
9-180
\eaaa
\bbbb
0&0&0&0&0&1&1&1&1\\
0&0&0&1&1&0&0&1&1\\
0&0&0&1&1&1&1&0&0\\
0&1&1&0&0&0&1&0&1\\
0&1&1&0&1&0&0&1&0\\
1&0&1&0&0&1&0&1&0\\
1&0&1&1&0&0&1&0&0\\
1&1&0&0&1&1&0&0&0\\
1&1&0&1&0&0&0&0&1\\
\ebbb
\parbox{7cm}{ 
\1\7\7\1\8\5\2\9\4\3\6\6\3\4\9\2\eeee
\6\3\4\9\2\5\8\1\7\7\1\8\5\2\9\4\eeee
\8\5\2\9\4\3\6\6\3\4\9\2\5\8\1\7\eeee
\2\5\8\1\7\7\1\8\5\2\9\4\3\6\6\3\eeee
\4\3\6\6\3\4\9\2\5\8\1\7\7\1\8\5\eeee
\7\7\1\8\5\2\9\4\3\6\6\3\4\9\2\5\eeee
\3\4\9\2\5\8\1\7\7\1\8\5\2\9\4\3\eeee
\5\2\9\4\3\6\6\3\4\9\2\5\8\1\7\7\eeee
\5\8\1\7\7\1\8\5\2\9\4\3\6\6\3\4\eeee
\3\6\6\3\4\9\2\5\8\1\7\7\1\8\5\2\eeee
\7\1\8\5\2\9\4\3\6\6\3\4\9\2\5\8\eeee
\4\9\2\5\8\1\7\7\1\8\5\2\9\4\3\6\eeee
} 

\baaa
9-181
\eaaa
\bbbb
0&0&0&0&0&1&1&1&1\\
0&0&0&1&1&0&0&1&1\\
0&0&0&1&1&1&1&0&0\\
0&1&1&0&1&0&0&0&1\\
0&1&1&1&0&0&1&0&0\\
1&0&1&0&0&0&1&1&0\\
1&0&1&0&1&1&0&0&0\\
1&1&0&0&0&1&0&0&1\\
1&1&0&1&0&0&0&1&0\\
\ebbb
\parbox{7cm}{ 
\1\7\5\2\8\6\3\4\9\1\7\5\2\8\6\3\eeee
\6\3\4\9\1\7\5\2\8\6\3\4\9\1\7\5\eeee
\7\5\2\8\6\3\4\9\1\7\5\2\8\6\3\4\eeee
\3\4\9\1\7\5\2\8\6\3\4\9\1\7\5\2\eeee
\5\2\8\6\3\4\9\1\7\5\2\8\6\3\4\9\eeee
\4\9\1\7\5\2\8\6\3\4\9\1\7\5\2\8\eeee
\2\8\6\3\4\9\1\7\5\2\8\6\3\4\9\1\eeee
\9\1\7\5\2\8\6\3\4\9\1\7\5\2\8\6\eeee
\8\6\3\4\9\1\7\5\2\8\6\3\4\9\1\7\eeee
\1\7\5\2\8\6\3\4\9\1\7\5\2\8\6\3\eeee
\6\3\4\9\1\7\5\2\8\6\3\4\9\1\7\5\eeee
\7\5\2\8\6\3\4\9\1\7\5\2\8\6\3\4\eeee
} 

\baaa
9-182
\eaaa
\bbbb
0&0&0&0&0&1&1&1&1\\
0&1&0&0&1&0&0&1&1\\
0&0&1&1&0&0&1&0&1\\
0&0&1&1&0&1&0&1&0\\
0&1&0&0&1&1&1&0&0\\
1&0&0&1&1&0&0&0&1\\
1&0&1&0&1&0&1&0&0\\
1&1&0&1&0&0&0&1&0\\
1&1&1&0&0&1&0&0&0\\
\ebbb
\parbox{7cm}{ 
\1\7\7\1\8\8\1\7\7\1\8\8\1\7\7\1\eeee
\6\5\5\6\4\4\6\5\5\6\4\4\6\5\5\6\eeee
\9\2\2\9\3\3\9\2\2\9\3\3\9\2\2\9\eeee
\1\8\8\1\7\7\1\8\8\1\7\7\1\8\8\1\eeee
\6\4\4\6\5\5\6\4\4\6\5\5\6\4\4\6\eeee
\9\3\3\9\2\2\9\3\3\9\2\2\9\3\3\9\eeee
\1\7\7\1\8\8\1\7\7\1\8\8\1\7\7\1\eeee
\6\5\5\6\4\4\6\5\5\6\4\4\6\5\5\6\eeee
\9\2\2\9\3\3\9\2\2\9\3\3\9\2\2\9\eeee
\1\8\8\1\7\7\1\8\8\1\7\7\1\8\8\1\eeee
\6\4\4\6\5\5\6\4\4\6\5\5\6\4\4\6\eeee
\9\3\3\9\2\2\9\3\3\9\2\2\9\3\3\9\eeee
} 

\baaa
9-183
\eaaa
\bbbb
0&0&0&0&0&1&1&1&1\\
0&1&0&0&1&0&0&1&1\\
0&0&1&1&0&0&1&0&1\\
0&0&1&1&0&1&1&0&0\\
0&1&0&0&1&1&0&1&0\\
1&0&0&1&1&0&0&0&1\\
1&0&1&1&0&0&1&0&0\\
1&1&0&0&1&0&0&1&0\\
1&1&1&0&0&1&0&0&0\\
\ebbb
\parbox{7cm}{ 
\1\7\7\1\8\8\1\7\7\1\8\8\1\7\7\1\eeee
\6\4\4\6\5\5\6\4\4\6\5\5\6\4\4\6\eeee
\9\3\3\9\2\2\9\3\3\9\2\2\9\3\3\9\eeee
\1\7\7\1\8\8\1\7\7\1\8\8\1\7\7\1\eeee
\6\4\4\6\5\5\6\4\4\6\5\5\6\4\4\6\eeee
\9\3\3\9\2\2\9\3\3\9\2\2\9\3\3\9\eeee
\1\7\7\1\8\8\1\7\7\1\8\8\1\7\7\1\eeee
\6\4\4\6\5\5\6\4\4\6\5\5\6\4\4\6\eeee
\9\3\3\9\2\2\9\3\3\9\2\2\9\3\3\9\eeee
\1\7\7\1\8\8\1\7\7\1\8\8\1\7\7\1\eeee
\6\4\4\6\5\5\6\4\4\6\5\5\6\4\4\6\eeee
\9\3\3\9\2\2\9\3\3\9\2\2\9\3\3\9\eeee
} 

\baaa
9-184
\eaaa
\bbbb
0&0&0&0&0&1&1&1&1\\
0&1&0&0&1&0&0&1&1\\
0&0&1&1&0&0&1&0&1\\
0&0&1&1&0&1&1&0&0\\
0&1&0&0&1&1&0&1&0\\
1&0&0&1&1&1&0&0&0\\
1&0&1&1&0&0&1&0&0\\
1&1&0&0&1&0&0&1&0\\
1&1&1&0&0&0&0&0&1\\
\ebbb
\parbox{7cm}{ 
\1\7\7\1\8\8\1\7\7\1\8\8\1\7\7\1\eeee
\6\4\3\9\2\5\6\4\3\9\2\5\6\4\3\9\eeee
\6\4\3\9\2\5\6\4\3\9\2\5\6\4\3\9\eeee
\1\7\7\1\8\8\1\7\7\1\8\8\1\7\7\1\eeee
\9\3\4\6\5\2\9\3\4\6\5\2\9\3\4\6\eeee
\9\3\4\6\5\2\9\3\4\6\5\2\9\3\4\6\eeee
\1\7\7\1\8\8\1\7\7\1\8\8\1\7\7\1\eeee
\6\4\3\9\2\5\6\4\3\9\2\5\6\4\3\9\eeee
\6\4\3\9\2\5\6\4\3\9\2\5\6\4\3\9\eeee
\1\7\7\1\8\8\1\7\7\1\8\8\1\7\7\1\eeee
\9\3\4\6\5\2\9\3\4\6\5\2\9\3\4\6\eeee
\9\3\4\6\5\2\9\3\4\6\5\2\9\3\4\6\eeee
} 

\baaa
9-185
\eaaa
\bbbb
0&0&0&0&0&1&1&1&1\\
0&1&0&0&1&0&0&1&1\\
0&0&1&1&0&0&1&0&1\\
0&0&1&1&1&0&0&1&0\\
0&1&0&1&2&0&0&0&0\\
1&0&0&0&0&2&1&0&0\\
1&0&1&0&0&1&0&1&0\\
1&1&0&1&0&0&1&0&0\\
1&1&1&0&0&0&0&0&1\\
\ebbb
\parbox{7cm}{ 
\1\7\8\2\2\8\7\1\8\4\4\8\1\7\8\2\eeee
\6\6\1\9\9\1\6\6\7\3\3\7\6\6\1\9\eeee
\6\6\7\3\3\7\6\6\1\9\9\1\6\6\7\3\eeee
\7\1\8\4\4\8\1\7\8\2\2\8\7\1\8\4\eeee
\3\9\2\5\5\2\9\3\4\5\5\4\3\9\2\5\eeee
\3\9\2\5\5\2\9\3\4\5\5\4\3\9\2\5\eeee
\7\1\8\4\4\8\1\7\8\2\2\8\7\1\8\4\eeee
\6\6\7\3\3\7\6\6\1\9\9\1\6\6\7\3\eeee
\6\6\1\9\9\1\6\6\7\3\3\7\6\6\1\9\eeee
\1\7\8\2\2\8\7\1\8\4\4\8\1\7\8\2\eeee
\9\3\4\5\5\4\3\9\2\5\5\2\9\3\4\5\eeee
\9\3\4\5\5\4\3\9\2\5\5\2\9\3\4\5\eeee
} 

\baaa
9-186
\eaaa
\bbbb
0&0&0&0&0&1&1&1&1\\
0&1&0&0&1&0&0&1&1\\
0&0&2&0&0&0&1&0&1\\
0&0&0&2&0&1&1&0&0\\
0&1&0&0&1&1&0&1&0\\
1&0&0&1&1&1&0&0&0\\
1&0&1&1&0&0&1&0&0\\
1&1&0&0&1&0&0&1&0\\
1&1&1&0&0&0&0&0&1\\
\ebbb
\parbox{7cm}{ 
\1\7\7\1\8\8\1\7\7\1\8\8\1\7\7\1\eeee
\6\4\4\6\5\2\9\3\3\9\2\5\6\4\4\6\eeee
\6\4\4\6\5\2\9\3\3\9\2\5\6\4\4\6\eeee
\1\7\7\1\8\8\1\7\7\1\8\8\1\7\7\1\eeee
\9\3\3\9\2\5\6\4\4\6\5\2\9\3\3\9\eeee
\9\3\3\9\2\5\6\4\4\6\5\2\9\3\3\9\eeee
\1\7\7\1\8\8\1\7\7\1\8\8\1\7\7\1\eeee
\6\4\4\6\5\2\9\3\3\9\2\5\6\4\4\6\eeee
\6\4\4\6\5\2\9\3\3\9\2\5\6\4\4\6\eeee
\1\7\7\1\8\8\1\7\7\1\8\8\1\7\7\1\eeee
\9\3\3\9\2\5\6\4\4\6\5\2\9\3\3\9\eeee
\9\3\3\9\2\5\6\4\4\6\5\2\9\3\3\9\eeee
} 

\baaa
9-187
\eaaa
\bbbb
0&0&0&0&0&1&1&1&1\\
0&2&0&0&0&0&0&1&1\\
0&0&2&0&0&0&1&0&1\\
0&0&0&2&0&1&0&1&0\\
0&0&0&0&2&1&1&0&0\\
1&0&0&1&1&1&0&0&0\\
1&0&1&0&1&0&1&0&0\\
1&1&0&1&0&0&0&1&0\\
1&1&1&0&0&0&0&0&1\\
\ebbb
\parbox{7cm}{ 
\1\7\7\1\8\8\1\7\7\1\8\8\1\7\7\1\eeee
\6\5\5\6\4\4\6\5\5\6\4\4\6\5\5\6\eeee
\6\5\5\6\4\4\6\5\5\6\4\4\6\5\5\6\eeee
\1\7\7\1\8\8\1\7\7\1\8\8\1\7\7\1\eeee
\9\3\3\9\2\2\9\3\3\9\2\2\9\3\3\9\eeee
\9\3\3\9\2\2\9\3\3\9\2\2\9\3\3\9\eeee
\1\7\7\1\8\8\1\7\7\1\8\8\1\7\7\1\eeee
\6\5\5\6\4\4\6\5\5\6\4\4\6\5\5\6\eeee
\6\5\5\6\4\4\6\5\5\6\4\4\6\5\5\6\eeee
\1\7\7\1\8\8\1\7\7\1\8\8\1\7\7\1\eeee
\9\3\3\9\2\2\9\3\3\9\2\2\9\3\3\9\eeee
\9\3\3\9\2\2\9\3\3\9\2\2\9\3\3\9\eeee
} 

\baaa
9-188
\eaaa
\bbbb
1&0&0&0&0&0&0&1&2\\
0&1&0&0&0&0&1&1&1\\
0&0&1&0&0&1&1&0&1\\
0&0&0&1&1&1&1&0&0\\
0&0&0&1&2&1&0&0&0\\
0&0&1&1&1&1&0&0&0\\
0&1&1&1&0&0&1&0&0\\
1&2&0&0&0&0&0&1&0\\
1&1&1&0&0&0&0&0&1\\
\ebbb
\parbox{7cm}{ 
\1\9\3\6\5\4\7\2\8\2\7\4\5\6\3\9\eeee
\1\9\3\6\5\4\7\2\8\2\7\4\5\6\3\9\eeee
\8\2\7\4\5\6\3\9\1\9\3\6\5\4\7\2\eeee
\8\2\7\4\5\6\3\9\1\9\3\6\5\4\7\2\eeee
\1\9\3\6\5\4\7\2\8\2\7\4\5\6\3\9\eeee
\1\9\3\6\5\4\7\2\8\2\7\4\5\6\3\9\eeee
\8\2\7\4\5\6\3\9\1\9\3\6\5\4\7\2\eeee
\8\2\7\4\5\6\3\9\1\9\3\6\5\4\7\2\eeee
\1\9\3\6\5\4\7\2\8\2\7\4\5\6\3\9\eeee
\1\9\3\6\5\4\7\2\8\2\7\4\5\6\3\9\eeee
\8\2\7\4\5\6\3\9\1\9\3\6\5\4\7\2\eeee
\8\2\7\4\5\6\3\9\1\9\3\6\5\4\7\2\eeee
} 

\baaa
9-189
\eaaa
\bbbb
1&0&0&0&0&0&0&1&2\\
0&1&0&0&0&0&1&1&1\\
0&0&1&0&0&1&1&1&0\\
0&0&0&1&1&1&1&0&0\\
0&0&0&1&2&1&0&0&0\\
0&0&1&1&1&1&0&0&0\\
0&1&1&1&0&0&1&0&0\\
1&1&1&0&0&0&0&1&0\\
2&1&0&0&0&0&0&0&1\\
\ebbb
\parbox{7cm}{ 
\1\8\3\6\5\4\7\2\9\1\8\3\6\5\4\7\eeee
\1\8\3\6\5\4\7\2\9\1\8\3\6\5\4\7\eeee
\9\2\7\4\5\6\3\8\1\9\2\7\4\5\6\3\eeee
\9\2\7\4\5\6\3\8\1\9\2\7\4\5\6\3\eeee
\1\8\3\6\5\4\7\2\9\1\8\3\6\5\4\7\eeee
\1\8\3\6\5\4\7\2\9\1\8\3\6\5\4\7\eeee
\9\2\7\4\5\6\3\8\1\9\2\7\4\5\6\3\eeee
\9\2\7\4\5\6\3\8\1\9\2\7\4\5\6\3\eeee
\1\8\3\6\5\4\7\2\9\1\8\3\6\5\4\7\eeee
\1\8\3\6\5\4\7\2\9\1\8\3\6\5\4\7\eeee
\9\2\7\4\5\6\3\8\1\9\2\7\4\5\6\3\eeee
\9\2\7\4\5\6\3\8\1\9\2\7\4\5\6\3\eeee
} 

\baaa
9-190
\eaaa
\bbbb
1&0&0&0&0&0&1&1&1\\
0&1&0&0&0&1&0&1&1\\
0&0&1&0&1&0&1&0&1\\
0&0&0&2&0&1&0&1&0\\
0&0&1&0&2&0&1&0&0\\
0&1&0&1&0&2&0&0&0\\
1&0&1&0&1&0&1&0&0\\
1&1&0&1&0&0&0&1&0\\
1&1&1&0&0&0&0&0&1\\
\ebbb
\parbox{7cm}{ 
\1\7\5\3\9\2\6\6\2\9\3\5\7\1\8\4\eeee
\1\7\5\3\9\2\6\6\2\9\3\5\7\1\8\4\eeee
\9\3\5\7\1\8\4\4\8\1\7\5\3\9\2\6\eeee
\9\3\5\7\1\8\4\4\8\1\7\5\3\9\2\6\eeee
\1\7\5\3\9\2\6\6\2\9\3\5\7\1\8\4\eeee
\1\7\5\3\9\2\6\6\2\9\3\5\7\1\8\4\eeee
\9\3\5\7\1\8\4\4\8\1\7\5\3\9\2\6\eeee
\9\3\5\7\1\8\4\4\8\1\7\5\3\9\2\6\eeee
\1\7\5\3\9\2\6\6\2\9\3\5\7\1\8\4\eeee
\1\7\5\3\9\2\6\6\2\9\3\5\7\1\8\4\eeee
\9\3\5\7\1\8\4\4\8\1\7\5\3\9\2\6\eeee
\9\3\5\7\1\8\4\4\8\1\7\5\3\9\2\6\eeee
} 

\baaa
9-191
\eaaa
\bbbb
2&0&0&0&0&0&0&0&2\\
0&2&0&0&0&0&0&1&1\\
0&0&2&0&0&0&1&1&0\\
0&0&0&2&0&1&1&0&0\\
0&0&0&0&2&2&0&0&0\\
0&0&0&1&1&2&0&0&0\\
0&0&1&1&0&0&2&0&0\\
0&1&1&0&0&0&0&2&0\\
1&1&0&0&0&0&0&0&2\\
\ebbb
\parbox{7cm}{ 
\1\9\2\8\3\7\4\6\5\6\4\7\3\8\2\9\eeee
\1\9\2\8\3\7\4\6\5\6\4\7\3\8\2\9\eeee
\1\9\2\8\3\7\4\6\5\6\4\7\3\8\2\9\eeee
\1\9\2\8\3\7\4\6\5\6\4\7\3\8\2\9\eeee
\1\9\2\8\3\7\4\6\5\6\4\7\3\8\2\9\eeee
\1\9\2\8\3\7\4\6\5\6\4\7\3\8\2\9\eeee
\1\9\2\8\3\7\4\6\5\6\4\7\3\8\2\9\eeee
\1\9\2\8\3\7\4\6\5\6\4\7\3\8\2\9\eeee
\1\9\2\8\3\7\4\6\5\6\4\7\3\8\2\9\eeee
\1\9\2\8\3\7\4\6\5\6\4\7\3\8\2\9\eeee
\1\9\2\8\3\7\4\6\5\6\4\7\3\8\2\9\eeee
\1\9\2\8\3\7\4\6\5\6\4\7\3\8\2\9\eeee
} 

\baaa
9-192
\eaaa
\bbbb
2&0&0&0&0&0&0&0&2\\
0&2&0&0&0&0&0&1&1\\
0&0&2&0&0&0&1&1&0\\
0&0&0&2&0&1&1&0&0\\
0&0&0&0&3&1&0&0&0\\
0&0&0&1&1&2&0&0&0\\
0&0&1&1&0&0&2&0&0\\
0&1&1&0&0&0&0&2&0\\
1&1&0&0&0&0&0&0&2\\
\ebbb
\parbox{7cm}{ 
\1\9\2\8\3\7\4\6\5\5\6\4\7\3\8\2\eeee
\1\9\2\8\3\7\4\6\5\5\6\4\7\3\8\2\eeee
\1\9\2\8\3\7\4\6\5\5\6\4\7\3\8\2\eeee
\1\9\2\8\3\7\4\6\5\5\6\4\7\3\8\2\eeee
\1\9\2\8\3\7\4\6\5\5\6\4\7\3\8\2\eeee
\1\9\2\8\3\7\4\6\5\5\6\4\7\3\8\2\eeee
\1\9\2\8\3\7\4\6\5\5\6\4\7\3\8\2\eeee
\1\9\2\8\3\7\4\6\5\5\6\4\7\3\8\2\eeee
\1\9\2\8\3\7\4\6\5\5\6\4\7\3\8\2\eeee
\1\9\2\8\3\7\4\6\5\5\6\4\7\3\8\2\eeee
\1\9\2\8\3\7\4\6\5\5\6\4\7\3\8\2\eeee
\1\9\2\8\3\7\4\6\5\5\6\4\7\3\8\2\eeee
} 

\baaa
9-193
\eaaa
\bbbb
2&0&0&0&0&0&0&1&1\\
0&2&0&0&0&0&1&0&1\\
0&0&2&0&0&1&0&1&0\\
0&0&0&2&1&0&1&0&0\\
0&0&0&1&2&1&0&0&0\\
0&0&1&0&1&2&0&0&0\\
0&1&0&1&0&0&2&0&0\\
1&0&1&0&0&0&0&2&0\\
1&1&0&0&0&0&0&0&2\\
\ebbb
\parbox{7cm}{ 
\1\8\3\6\5\4\7\2\9\1\8\3\6\5\4\7\eeee
\1\8\3\6\5\4\7\2\9\1\8\3\6\5\4\7\eeee
\1\8\3\6\5\4\7\2\9\1\8\3\6\5\4\7\eeee
\1\8\3\6\5\4\7\2\9\1\8\3\6\5\4\7\eeee
\1\8\3\6\5\4\7\2\9\1\8\3\6\5\4\7\eeee
\1\8\3\6\5\4\7\2\9\1\8\3\6\5\4\7\eeee
\1\8\3\6\5\4\7\2\9\1\8\3\6\5\4\7\eeee
\1\8\3\6\5\4\7\2\9\1\8\3\6\5\4\7\eeee
\1\8\3\6\5\4\7\2\9\1\8\3\6\5\4\7\eeee
\1\8\3\6\5\4\7\2\9\1\8\3\6\5\4\7\eeee
\1\8\3\6\5\4\7\2\9\1\8\3\6\5\4\7\eeee
\1\8\3\6\5\4\7\2\9\1\8\3\6\5\4\7\eeee
} 

\baaa
9-194
\eaaa
\bbbb
2&0&0&0&0&0&0&1&1\\
0&2&0&0&0&0&1&0&1\\
0&0&2&0&0&1&0&1&0\\
0&0&0&2&1&0&1&0&0\\
0&0&0&1&3&0&0&0&0\\
0&0&1&0&0&3&0&0&0\\
0&1&0&1&0&0&2&0&0\\
1&0&1&0&0&0&0&2&0\\
1&1&0&0&0&0&0&0&2\\
\ebbb
\parbox{7cm}{ 
\1\8\3\6\6\3\8\1\9\2\7\4\5\5\4\7\eeee
\1\8\3\6\6\3\8\1\9\2\7\4\5\5\4\7\eeee
\1\8\3\6\6\3\8\1\9\2\7\4\5\5\4\7\eeee
\1\8\3\6\6\3\8\1\9\2\7\4\5\5\4\7\eeee
\1\8\3\6\6\3\8\1\9\2\7\4\5\5\4\7\eeee
\1\8\3\6\6\3\8\1\9\2\7\4\5\5\4\7\eeee
\1\8\3\6\6\3\8\1\9\2\7\4\5\5\4\7\eeee
\1\8\3\6\6\3\8\1\9\2\7\4\5\5\4\7\eeee
\1\8\3\6\6\3\8\1\9\2\7\4\5\5\4\7\eeee
\1\8\3\6\6\3\8\1\9\2\7\4\5\5\4\7\eeee
\1\8\3\6\6\3\8\1\9\2\7\4\5\5\4\7\eeee
\1\8\3\6\6\3\8\1\9\2\7\4\5\5\4\7\eeee
} 


\bibliographystyle{plain}
\bibliography{../k}

\begin{thebibliography}{1}

\bibitem{Puz2004}
С.~А. Пузынина.
\newblock Периодичность совершенных раскрасок бесконечной прямоугольной
  решетки.
\newblock {\em Дискрет. анализ и исслед. операций. Сер. 1}, 11(1):79--92, 2004.

\bibitem{Puz2005}
С.~А. Пузынина.
\newblock Совершенные раскраски вершин графа $G(Z^2)$ в три цвета.
\newblock {\em Дискрет. анализ и исслед. операций. Сер. 2}, 12(1):37--54, 2005.

\bibitem{Axe2003}
M.~A. Axenovich.
\newblock On multiple coverings of the infinite rectangular grid with
balls of constant radius.
\newblock {\em Discrete Math.}, 268(1-3):37--54, 2003.

\end{thebibliography}
\end{document}